%% file: 1601.03259.English.tex
\documentclass{amsbook}

\input{1601.03259}

\begin{document}
\title{Introduction into Calculus over  Banach Algebra}
\input cyracc.def
\font\tencyr=wncyr10
\def\cyr{\tencyr\cyracc}
\ShowEq{contents}
\end{document}

%% file: 1601.03259.tex
\def\BookNumber{1601.03259}
\def\PrintBook{}
\input{Commands}
\usepackage{chngcntr}
\usepackage{array,arydshln}
\xRefDef{0906.0135}
\xRefDef{0912.3315}
\xRefDef{1003.1544}
\xRefDef{1305.4547}
\xRefDef{1302.7204}
\xRefDef{1310.5591}
\xRefDef{1502.04063}

\DefEq
{
\input{\FilePrefix Calculus.16.Absract.\TheLanguage}

\maketitle
\tableofcontents

\input{\FilePrefix Preface.\BookNumber.\TheLanguage}
\input{\FilePrefix Preface.1610.309618526.\TheLanguage}
\input{\FilePrefix Preliminary.\TheLanguage}
\input{\FilePrefix Derivative.16.\TheLanguage}
\input{\FilePrefix Second.Derivative.16.\TheLanguage}
\input{\FilePrefix Indefinite.Integral.\TheLanguage}

\input{\FilePrefix Two.Integrals.\TheLanguage}
\input{\FilePrefix Diff.Form.1.\TheLanguage}

\input{\FilePrefix Diff.Form.2.\TheLanguage}
\newcounter{Chapter}
\setcounter{Chapter}{\thechapter}
\appendix
\input{\FilePrefix Appendix.16.Calculate.\TheLanguage}
\input{\FilePrefix Summary.Calculus.\TheLanguage}

\input{\FilePrefix Preliminary.Omega.\TheLanguage}
\input{\FilePrefix Biblio.\TheLanguage}
\input{\FilePrefix Index.\TheLanguage}
\input{\FilePrefix Symbol.\TheLanguage}
}
{contents}

%% file: Commands.tex
\def\Defined{}
\scrollmode
\ifx\FilePrefix\undefined
\newcommand{\FilePrefix}{}
\fi
\ifx\UseRussian\Defined
\usepackage{cmap}
\usepackage[T2A,T2B]{fontenc}
\usepackage[cp1251]{inputenc}
\fi
\ifx\GJSFRA\Defined
\usepackage[top=1.905cm,bottom=1.905cm,inner=1.65cm,outer=1.65cm]{geometry}
\fi
\ifx\Presentation\Defined
\usepackage[margin=1cm]{geometry}
\fi

\ifx\setCACAA\Defined
\usepackage{times,mathptm}
\usepackage{amsmath,amsfonts,amssymb}
\usepackage{graphicx,graphics,epsfig}
\usepackage{fancyhdr,fancybox}
\usepackage{verbatim}
\usepackage{setspace}
\fi

\ifx\CreateSpace\Defined
\usepackage[top=1.905cm,bottom=1.905cm,inner=1.905cm,outer=1.27cm]{geometry}
\def\Publisher{CreateSpace Independent Publishing Platform}
\fi
\ifx\PublishBook\Defined
\def\PrintPaper{}
\usepackage{setspace}
\ifx\UseRussian\undefined
\usepackage{pslatex}
\fi
\usepackage[margin=2cm]{geometry}
\fi
\usepackage{graphicx}
\usepackage{footmisc}
\usepackage[all]{xy}
\usepackage{color}
\ifx\PrintPaper\undefined
\definecolor{CoverColor}{rgb}{.82,.7,.55}
\definecolor{UrlColor}{rgb}{.9,0,.3}
\definecolor{SymbColor}{rgb}{.4,0,.9}
\definecolor{IndexColor}{rgb}{1,.3,.6}
\newcommand\BlueText[1]{\textcolor{blue}{#1}}
\newcommand\RedText[1]{\textcolor{red}{#1}}
\else
\definecolor{UrlColor}{rgb}{.1,.1,.1}
\definecolor{SymbColor}{rgb}{.1,.1,.1}
\definecolor{IndexColor}{rgb}{.1,.1,.1}
\newcommand\BlueText[1]{#1}
\newcommand\RedText[1]{#1}
\fi

\usepackage{chngcntr}
\usepackage{xr-hyper}
\usepackage[unicode]{hyperref}
\hypersetup{pdfdisplaydoctitle=true}
\hypersetup{colorlinks}
\hypersetup{citecolor=UrlColor}
\hypersetup{urlcolor=UrlColor}
\hypersetup{linkcolor=UrlColor}
\hypersetup{pdffitwindow=true}
\hypersetup{pdfnewwindow=true}
\hypersetup{pdfstartview={FitH}}
\ifx\UseRussian\Defined
\usepackage[english,russian]{babel}
\fi

\newcounter{Index}
\newcounter{Symbol}
\newcounter{Symbols}

\def\hyph{\penalty0\hskip0pt\relax-\penalty0\hskip0pt\relax}
\def\Hyph{-\penalty0\hskip0pt\relax}%

\def\ValueOff{off}%
\def\ValueOn{on}%
\def\Items#1{\ItemList#1,LastItem,}%
\def\LastItem{LastItem}%
\def\ItemList#1,{\def\ViewBook{#1}%
\ifx\ViewBook\LastItem%
\else%
\ifx\ViewBook\BookNumber%
\def\Semafor{on}%
\fi%
\expandafter\ItemList%
\fi%
}%

\newcommand{\ePrints}[1]%
{%
\def\Semafor{off}%
\Items{#1}%
}%

\newcommand{\Basis}[1]{\overline{\overline{#1}}{}}
\newcommand{\Vector}[1]{\overline{#1}{}}
\ifx\PrintPaper\undefined
\newcommand{\gi}[1]{\boldsymbol{\textcolor{IndexColor}{#1}}}
\else
\newcommand{\gi}[1]{\boldsymbol{#1}}
\fi
\newcommand\gii{\gi i}
\newcommand\giI{\gi I}
\newcommand\gij{\gi j}
\newcommand\giJ{\gi J}
\newcommand\gik{\gi k}
\newcommand\gil{\gi l}
\newcommand\gin{\gi n}
\newcommand\gim{\gi m}
\newcommand\giA{\gi 1}

\newcommand{\VX}[1]{\Vector{#1}_{[1]}}
\makeatletter
\newcommand{\NameDef}[1]{%
\expandafter\gdef\csname #1\endcsname%
}%
\newcommand{\xNameDef}[1]{%
\expandafter\xdef\csname #1\endcsname%
}%
\newcommand{\ShowSymbol}[2]{%
\@nameuse{ViewSymbol#1,,,#2}%
}%
\newcommand{\symb}[4][]{%
\def\TempA{}%
\def\TempB{#1}%
\ifx\TempA\TempB%
\def\ThisSymbol{#3}%
\else%
\edef\ThisSymbol{#3(#1)}%
\fi%
\@ifundefined{ViewSymbol\ThisSymbol}{%
\addtocounter{Symbols}{1}%
\edef\SymbolId{\arabic{Symbols}}%
\xNameDef{ViewSymbol\ThisSymbol}{\SymbolId}%
\NameDef{ViewSymbol\ThisSymbol:::\SymbolId}{#2}%
\@namedef{RefSymbol}{:}%
}{%
\edef\Symbols{\@nameuse{ViewSymbol\ThisSymbol}}%
\def\aSymbolId{0}%
\@for\Symbol:=\Symbols\do{%
\protected@edef\TempA{#2}%
\protected@edef\TempB{\@nameuse{ViewSymbol\ThisSymbol:::\Symbol}}%
\ifx\TempA\TempB%
\edef\aSymbolId{\Symbol}%
\fi%
}%
\def\Zero{0}%
\ifx\aSymbolId\Zero%
\addtocounter{Symbols}{1}%
\edef\SymbolIds{\@nameuse{ViewSymbol\ThisSymbol},\arabic{Symbols}}%
\xNameDef{ViewSymbol\ThisSymbol}{\SymbolIds}%
\edef\SymbolId{\arabic{Symbols}}%
\NameDef{ViewSymbol\ThisSymbol:::\SymbolId}{#2}%
\else%
\def\SymbolId{\aSymbolId}%
\fi%
\addtocounter{Symbol}{1}%
\@namedef{RefSymbol}{\arabic{Symbol}}%
}%
\@namedef{LabelSymbol}{\label{symbol: \ThisSymbol:\@nameuse{RefSymbol}}}%
\edef\RefIds{RefSymbol\ThisSymbol===\SymbolId}%
\@ifundefined{\RefIds}{%
\xNameDef{\RefIds}{\@nameuse{RefSymbol}}%
}{%
\xNameDef{\RefIds}{\@nameuse{\RefIds},\@nameuse{RefSymbol}}%
}%
\NameDef{ViewSymbol#3,,,#4}{\textcolor{SymbColor}{#2}}%
\def\Temp{#4}%
\def\One{1}%
\def\Two{2}%
\ifx\Temp\One%
$\@nameuse{ViewSymbol#3,,,#4}$%
\fi%
\ifx\Temp\Two%
\[\@nameuse{ViewSymbol#3,,,#4}\]%
\fi%
\@nameuse{LabelSymbol}%
}%
\newcommand\AddEq[3][0]%
{%
\@ifundefined{ViewEq #1}{%
\expandafter\newcommand\csname ViewEq #2\endcsname[#1]{#3}%
}{%
\errmessage {second entry of AddEq: #2}%
}%
}%
\newcommand\DefEq[2]{%
\@ifundefined{ViewEq #2}{%
\NameDef{ViewEq #2}{#1}%
}{%
\errmessage {second entry of DefEq: #2}%
}%
}%
\newcommand{\DefEquation}[2]{%
\AddEq{#2}%
{%
\begin{equation}%
#1%
\EqLabel{#2}%
\end{equation}%
}%
}%
\newcommand{\AddEquation}[2]{%
\AddEq{#1}{\begin{equation}#2\EqLabel{#1}\end{equation}}%
}%
\def\ViewParm#1{\protect\getParm#1,endParm,}%
\def\endParm{endParm}%
\def\getParm#1,{\def\temp{#1}%
\ifx\temp\endParm%
\else%
\ShowEq{#1}%
\expandafter\getParm%
\fi%
}%
\newcommand{\DefTheorem}[2]{%
\AddEq{theorem: #1}
{
\begin{theorem}
\label{theorem: #1}
#2
\end{theorem}
}
}
\newcommand{\ShowTheorem}[1]{%
\ShowEq{theorem: #1}%
}%
\newcommand{\DefProof}[2]{%
\AddEq{proof: #1}
{
\begin{proof}
#2
\end{proof}%
}
}
\newcommand{\ShowProof}[1]{\ShowEq{proof: #1}}%
\newcommand{\DefDefinition}[2]{%
\AddEq{definition: #1}
{
\begin{definition}
\labelDefinition{#1}
{\it
#2
}
\qed
\end{definition}
}
}%
\newcommand\ShowDefinition[1]{\ShowEq{definition: #1}}
\newcommand{\DefExample}[2]{%
\AddEq{example: #1}
{
\begin{example}
\labelExample{#1}
#2
\qed
\end{example}
}
}%
\newcommand{\ShowExample}[1]{%
\ShowEq{example: #1}%
}%
\newcommand\DefRemark[2]{%
\AddEq{remark: #1}
{
\begin{remark}
\labelRemark{#1}
#2
\qed
\end{remark}
}
}%
\newcommand{\ShowRemark}[1]{%
\ShowEq{remark: #1}%
}%
\newcommand{\EqParm}[2]{%
\ViewParm{#2}\ShowEq{#1}%
}%
\newcommand{\EquationParm}[2]{%
\@ifundefined{ViewEq #1[#2]}%
{%
\ViewParm{#2}%
\DefEquation{\ShowEq{#1}}{#1[#2]}%
}{}%
\ShowEq{#1[#2]}%
}%
\newcommand\DrawEqParm[3]{%
\ViewParm{#2}%
\@ifundefined{ViewEq #1(#2)}{%
\DefEq%
{%
\ShowEq{#1}%
}{#1(#2)}%
}{%
}%
\DrawEq{#1(#2)}{#3}%
}%
\newcommand\EqRef[2][]%
{%
\def\Semafor{on}%
\def\Temp{}%
\edef\Tempa{#1}%
\ifx\Tempa\Temp%
\def\Semafor{off}%
\fi%
\@ifundefined{xRefDef#1}{%
\def\Semafor{off}
}{}%
\ifx\Semafor\ValueOff%
\eqref{eq: #2}%
\else%
\citeBib{#1}-\eqref{#1-\TheLanguage-eq: #2}%
\fi%
}%
\newcommand\eqRef[3][]{\EqRef[#1]{#2(#3)}}%
\newcommand\EqRefParm[2]{\EqRef{#1[#2]}}%
\newcommand\EqLabel[1]{\label{eq: #1}}%
\newcommand\ShowEq[1]{%
\@ifundefined{ViewEq #1}{%
\message {error: missed ShowEq #1}%
\newline%
\RedText{missed ShowEq #1}%
\newline%
}{%
\csname ViewEq #1\endcsname
}%
}%
\newcommand\DrawEq[3][]{%
\def\Temp{}%
\def\Tempa{#3}%
\ifx\Tempa\Temp%
\[\ShowEq{#2}#1\]
\else%
\def\Temp{-}%
\ifx\Tempa\Temp%
$\ShowEq{#2}#1$
\else%
\@ifundefined{ViewEq #2(#3)#1}{%
\AddEquation{#2(#3)}{\ShowEq{#2}#1}%
  }{%
}%
\ShowEq{#2(#3)}
\fi%
\fi%
}%
\makeatother
\DeclareMathOperator{\Hom}{\mathrm{Hom}} 
\DeclareMathOperator{\End}{\mathrm{End}} 
 
\DeclareMathOperator{\id}{\mathrm{id}} 
 
\newcommand{\subs}{${}_*$\Hyph}
\newcommand{\sups}{${}^*$\Hyph}

\newcommand{\CRstar}{{}^*{}_*}
\newcommand{\RCstar}{{}_*{}^*}
\newcommand{\CRcirc}{{}^{\circ}{}_{\circ}}
\newcommand{\RCcirc}{{}_{\circ}{}^{\circ}}

\newcommand{\RC}{$\RCstar$\Hyph}
\newcommand{\CR}{$\CRstar$\Hyph}
\newcommand{\drc}{$D\RCstar$\Hyph}
\newcommand{\Drc}{$\mathcal D\RCstar$\Hyph}
\newcommand{\dcr}{$D\CRstar$\hyph}
\newcommand{\rcd}{$\RCstar D$\Hyph}
\newcommand{\crd}{$\CRstar D$\Hyph}

\newcommand\sT[1]{$*#1$\Hyph}%
\newcommand\Ts[1]{$#1*$\Hyph}%
\newcommand\sD{$\star D$\Hyph}%
\newcommand\Ds{$D\star$\Hyph}%
\newcommand\VirtFrac{\vphantom{\overset{\rightarrow}{\frac 11}^{\frac 11}}}
\newcommand\VirtVar{\vphantom{\overset{\rightarrow}{1}^1}}
\newcommand\pC[2]{{}_{#1\cdot #2}}%
\newcommand\DcrPartial[1]%
{%
\def\tempa{}%
\def\tempb{#1}%
\ifx\tempa\tempc%
(\partial\CRstar)%
\else%
(\partial_{\gi{#1}}\CRstar)%
\fi%
}%
\newcommand\rcDPartial[1]%
{%
\def\tempa{}%
\def\tempb{#1}%
\ifx\tempa\tempc%
(\RCstar\partial)%
\else%
(\RCstar\partial_{\gi{#1}})%
\fi%
}%
\newcommand\StandPartial[3]%
{%
\frac{d^{\gi{#3}} #1}{d #2}%
}%

\ifx\setCACAA\undefined
\renewcommand{\uppercasenonmath}[1]{}

\makeatletter
\newcommand\@dotsep{4.5}
\def\@tocline#1#2#3#4#5#6#7
{\relax
\par \addpenalty\@secpenalty\addvspace{#2}%
\begingroup 
\@ifempty{#4}{%
\@tempdima\csname r@tocindent\number#1\endcsname\relax
}{%
\@tempdima#4\relax
}%
\parindent\z@ \leftskip#3\relax \advance\leftskip\@tempdima\relax
\rightskip\@pnumwidth plus1em \parfillskip-\@pnumwidth
#5\leavevmode\hskip-\@tempdima #6\relax
\leaders\hbox{$\m@th
\mkern \@dotsep mu\hbox{.}\mkern \@dotsep mu$}\hfill
\hbox to\@pnumwidth{\@tocpagenum{#7}}\par
\nobreak
\endgroup
}
\makeatother 
\fi

\ifx\PrintBook\undefined
\def\Chapter{\section}
\def\Section{\subsection}
\ifx\setCACAA\undefined
\makeatletter
\renewcommand{\@indextitlestyle}{%
\twocolumn[\section*{\indexname}]%
\def\IndexSpace{off}%
}
\makeatother 
\ifx\PrintPaper\undefined
\thanks{\href{mailto:Aleks\_Kleyn@MailAPS.org}{Aleks\_Kleyn@MailAPS.org}}
\ePrints{1102.1776,1201.4158}
\ifx\Semafor\ValueOff
\thanks{\ \ \ \url{http://AleksKleyn.dyndns-home.com:4080/}\ \ \ \ \ \url{http://arxiv.org/a/kleyn\_a\_1}}
\thanks{\ \ \ \url{http://sites.google.com/site/AleksKleyn/}\ \ \ \ \url{http://AleksKleyn.blogspot.com/}}
\fi
\fi
\fi
\else
\def\Chapter{\chapter}
\def\Section{\section}

\makeatletter
\def\@maketitle{%
  \cleardoublepage \thispagestyle{empty}%
  \begingroup \topskip\z@skip
  \null\vfil
  \begingroup
  \LARGE\bfseries \centering
  \openup\medskipamount
  \@title
  \par
  \ifx\subtitle\undefined
  \else
  \centerline{\ }
  \centerline{\emph\subtitle}
  \fi
  \ifx\subtitleA\undefined
  \else
  \centerline{\emph\subtitleA}
  \fi
  \ifx\edition\undefined
  \else
  \centerline{\emph\edition}
  \fi
  \par\vspace{24pt}%
  \def\and{\par\medskip}\centering
  \mdseries\authors\par\bigskip
  \endgroup
  \vfill
\noindent
\href{mailto:Aleks\_Kleyn@MailAPS.org}{Aleks\_Kleyn@MailAPS.org}
\newline
\url{http://AleksKleyn.dyndns-home.com:4080/}
\newline
\url{http://sites.google.com/site/AleksKleyn/}
\newline
\url{http://arxiv.org/a/kleyn\_a\_1}
\newline
\url{http://AleksKleyn.blogspot.com/}
  \newpage\thispagestyle{empty}
  \begin{center}
    \ifx\@empty\@subjclass\else\@setsubjclass\fi
    \ifx\@empty\@keywords\else\@setkeywords\fi
    \ifx\@empty\@translators\else\vfil\@settranslators\fi
    \ifx\@empty\thankses\else\vfil\@setthanks\fi
  \end{center}
  \vfil
  \@setabstract
\vfil
  \def\Temp{0000}
  \ifx\copyrightyear\Temp
  \else
  \begin{center}
\begin{tabular}{@{}c}
Copyright\ \copyright\ \copyrightyear\ \copyrightholder
\\
All rights reserved.
\end{tabular}
  \end{center}
  \fi
  \ifx\Publisher\undefined%
  \else
  \begin{center}
\begin{tabular}{@{}c}
\Publisher
\end{tabular}
  \end{center}
  \fi
  \ifx\ISBN\undefined%
  \else%
 \begin{center}
\begin{tabular}{@{}r@{\ }l}
ISBN:&\ISBN
\\
ISBN-13:&\ISBNa
\end{tabular}
  \end{center}
  \fi%
  \ifx\titleNote\undefined
  \else
  \par\vspace{24pt}%
  \centerline{\mdseries\titleNote}
	  \centerline{\Title}
	  \ifx\Subtitle\undefined
	  \else
	  \centerline{\emph\Subtitle}
	  \fi
	  \ifx\Edition\undefined
	  \else
	  \centerline{\Edition}
	  \fi
	  \centerline{\Authors}
  \fi
  \endgroup}
\def\chapter{%
	\clearpage
  \thispagestyle{plain}\global\@topnum\z@
  \@afterindenttrue \secdef\@chapter\@schapter}
\renewcommand{\@indextitlestyle}{%
\twocolumn[\chapter*{\indexname}]%
\def\IndexSpace{off}%
\let\@secnumber\@empty
\chaptermark{\indexname}%
}
\makeatother 
\email{\href{mailto:Aleks\_Kleyn@MailAPS.org}{Aleks\_Kleyn@MailAPS.org}}
\ePrints{1102.1776,1201.4158}
\ifx\Semafor\ValueOff
\urladdr{\url{http://AleksKleyn.dyndns-home.com:4080/}}
\urladdr{\url{http://sites.google.com/site/alekskleyn/}}
\urladdr{\url{http://arxiv.org/a/kleyn\_a\_1}}
\urladdr{\url{http://AleksKleyn.blogspot.com/}}
\fi
\fi

\ifx\SelectlEnglish\undefined
\ifx\UseRussian\undefined
\def\SelectlEnglish{}
\fi
\fi

\newcommand\arXivOldRef{http://arxiv.org/PS_cache/}
\newcommand\arXivRef{http://arxiv.org/pdf/}
\newcommand\AmazonRef{http://www.amazon.com/s/ref=nb_sb_noss?url=search-alias=aps&field-keywords=aleks+kleyn}
\newcommand\wRefDef[2]
{
\def\Tempa{#1}
\def\Tempb{0405.027}
\ifx\Tempa\Tempb
\def\wRef{\arXivOldRef gr-qc/pdf/0405/0405027v3.pdf}
\fi
\def\Tempb{0405.028}
\ifx\Tempa\Tempb
\def\wRef{\arXivOldRef gr-qc/pdf/0405/0405028v5.pdf}
\fi
\def\Tempb{0412.391}
\ifx\Tempa\Tempb
\def\wRef{\arXivOldRef math/pdf/0412/0412391v4.pdf}
\fi
\def\Tempb{0612.111}
\ifx\Tempa\Tempb
\def\wRef{\arXivOldRef math/pdf/0612/0612111v2.pdf}
\fi
\def\Tempb{0701.238}
\ifx\Tempa\Tempb
\def\wRef{\arXivOldRef math/pdf/0701/0701238v6.pdf}
\fi
\def\Tempb{0702.561}
\ifx\Tempa\Tempb
\def\wRef{\arXivOldRef math/pdf/0702/0702561v3.pdf}
\fi
\def\Tempb{0707.2246}
\ifx\Tempa\Tempb
\def\wRef{\arXivRef 0707.2246v2.pdf}
\fi
\def\Tempb{0803.3276}
\ifx\Tempa\Tempb
\def\wRef{\arXivRef 0803.3276v3.pdf}
\fi
\def\Tempb{0812.4763}
\ifx\Tempa\Tempb
\def\wRef{\arXivRef 0812.4763v7.pdf}
\fi
\def\Tempb{0906.0135}
\ifx\Tempa\Tempb
\def\wRef{\arXivRef 0906.0135v3.pdf}
\fi
 \def\Tempb{0909.0855}
\ifx\Tempa\Tempb
\def\wRef{\arXivRef 0909.0855v5.pdf}
\fi
 \def\Tempb{0912.3315}
\ifx\Tempa\Tempb
\def\wRef{\arXivRef 0912.3315v3.pdf}
\fi
 \def\Tempb{0912.4061}
\ifx\Tempa\Tempb
\def\wRef{\arXivRef 0912.4061v2.pdf}
\fi
 \def\Tempb{1001.4852}
\ifx\Tempa\Tempb
\def\wRef{\arXivRef 1001.4852.pdf}
\fi
 \def\Tempb{1003.3714}
\ifx\Tempa\Tempb
\def\wRef{\arXivRef 1003.3714v2.pdf}
\fi
 \def\Tempb{1003.1544}
\ifx\Tempa\Tempb
\def\wRef{\arXivRef 1003.1544v2.pdf}
\fi
 \def\Tempb{1006.2597}
\ifx\Tempa\Tempb
\def\wRef{\arXivRef 1006.2597v2.pdf}
\fi
 \def\Tempb{1011.3102}
\ifx\Tempa\Tempb
\def\wRef{\arXivRef 1011.3102.pdf}
\fi
 \def\Tempb{1104.5197}
\ifx\Tempa\Tempb
\def\wRef{\arXivRef 1104.5197.pdf}
\fi
 \def\Tempb{1105.4307}
\ifx\Tempa\Tempb
\def\wRef{\arXivRef 1105.4307.pdf}
\fi
 \def\Tempb{1107.1139}
\ifx\Tempa\Tempb
\def\wRef{\arXivRef 1107.1139.pdf}
\fi
 \def\Tempb{1107.5037}
\ifx\Tempa\Tempb
\def\wRef{\arXivRef 1107.5037.pdf}
\fi
 \def\Tempb{1111.6035}
\ifx\Tempa\Tempb
\def\wRef{\arXivRef 1111.6035.pdf}
\fi
 \def\Tempb{1202.6021}
\ifx\Tempa\Tempb
\def\wRef{\arXivRef 1202.6021v2.pdf}
\fi
 \def\Tempb{1211.6965}
\ifx\Tempa\Tempb
\def\wRef{\arXivRef 1211.6965.pdf}
\fi
 \def\Tempb{1302.7204}
\ifx\Tempa\Tempb
\def\wRef{\arXivRef 1302.7204v1.pdf}
\fi
 \def\Tempb{1305.4547}
\ifx\Tempa\Tempb
\def\wRef{\arXivRef 1305.4547.pdf}
\fi
 \def\Tempb{1310.5591}
\ifx\Tempa\Tempb
\def\wRef{\arXivRef 1310.5591.pdf}
\fi
 \def\Tempb{1502.04063}
\ifx\Tempa\Tempb
\def\wRef{\arXivRef 1502.04063v2.pdf}
\fi
 \def\Tempb{1505.03625}
\ifx\Tempa\Tempb
\def\wRef{\arXivRef 1505.03625v1.pdf}
\fi
 \def\Tempb{1601.03259}
\ifx\Tempa\Tempb
\def\wRef{\arXivRef 1601.03259v1.pdf}
\fi
 \def\Tempb{MRepro}
\ifx\Tempa\Tempb
\def\wRef{\arXivRef MReprov1.pdf}
\fi
 \def\Tempb{8433-5163}
\ifx\Tempa\Tempb
\def\wRef{\AmazonRef}
\fi
 \def\Tempb{8443-0072}
\ifx\Tempa\Tempb
\def\wRef{\AmazonRef}
\fi
 \def\Tempb{4776-3181}
\ifx\Tempa\Tempb
\def\wRef{\AmazonRef}
\fi
 \def\Tempb{5059-9176}
\ifx\Tempa\Tempb
\def\wRef{\AmazonRef}
\fi
 \def\Tempb{5114-6019}
\ifx\Tempa\Tempb
\def\wRef{\AmazonRef}
\fi
 \def\Tempb{5148-4632}
\ifx\Tempa\Tempb
\def\wRef{\AmazonRef}
\fi
 \def\Tempb{5410-9916}
\ifx\Tempa\Tempb
\def\wRef{\AmazonRef}
\fi
 \def\Tempb{BRepro}
\ifx\Tempa\Tempb
\def\wRef{\AmazonRef}
\fi
 \def\Tempb{CACAA.01.291}
\ifx\Tempa\Tempb
\def\wRef{http://www.cliffordanalysis.com/}
\fi
 \def\Tempb{CACAA.02.097}
\ifx\Tempa\Tempb
\def\wRef{http://www.cliffordanalysis.com/}
\fi
 \def\Tempb{CACAA.04.001}
\ifx\Tempa\Tempb
\def\wRef{http://www.cliffordanalysis.com/}
\fi
 \def\Tempb{GJSFRA.13.1.39}
\ifx\Tempa\Tempb
\def\wRef{http://www.cliffordanalysis.com/}
\fi
\externaldocument[#1-#2-]{\FilePrefix #1.#2}[\wRef]
}
\newcommand\LanguagePrefix{}%
\makeatletter
\newcommand\StartLabelItem[1][theorem]
{
\counterwithin{enumi}{#1}%
\def\labelenumi{\theenumi:}%
}
\newcommand\StopLabelItem
{
\def\theenumi{\@arabic\c@enumi}
\def\labelenumi{(\theenumi)}
}
\makeatother
\newcommand\labelConvention[1]{\label{convention: #1}}%
\newcommand\labelTheorem[1]{\label{theorem: #1}}%
\newcommand\labelLemma[1]{\label{lemma: #1}}%
\newcommand\labelCorollary[1]{\label{corollary: #1}}%
\newcommand\labelDefinition[1]{\label{definition: #1}}%
\newcommand\labelExample[1]{\label{example: #1}}%
\newcommand\labelRemark[1]{\label{remark: #1}}%
\newcommand\labelSection[1]{\label{section: #1}}%
\newcommand\labelChapter[1]{\label{chapter: #1}}%
\newcommand\labelItem[1]{\label{item: #1}}%

\makeatletter
\newcommand\Ref[2][]
{%
\def\Semafor{on}%
\def\Temp{}%
\edef\Tempa{#1}%
\ifx\Tempa\Temp%
\def\Semafor{off}%
\fi
\@ifundefined{xRefDef#1}{%
\def\Semafor{off}
}{}%
\ifx\Semafor\ValueOff%
\ref{#2}%
\else%
\citeBib{#1}-\ref{#1-\TheLanguage-#2}%
\fi%
}%
\makeatother
\newcommand\RefConvention[2][]{\Ref[#1]{convention: #2}}
\newcommand\RefTheorem[2][]{\Ref[#1]{theorem: #2}}
\newcommand\RefLemma[2][]{\Ref[#1]{lemma: #2}}
\newcommand\RefCorollary[2][]{\Ref[#1]{corollary: #2}}
\newcommand\RefDefinition[2][]{\Ref[#1]{definition: #2}}
\newcommand\RefExample[2][]{\Ref[#1]{example: #2}}

\newcommand\RefSection[2][]{\Ref[#1]{section: #2}}

\newcommand\RefItem[2][]{\Ref[#1]{item: #2}}

\ifx\SelectlEnglish\undefined
\newcommand\TheLanguage{Russian}%
\author{Александр Клейн}
\ifx\setCACAA\undefined
\newtheorem{theorem}{Теорема}[section]
\newtheorem{corollary}[theorem]{Следствие}
\newtheorem{example}[theorem]{Пример}
\newtheorem{definition}[theorem]{Определение}
\newtheorem{remark}[theorem]{Замечание}
\newtheorem{question}[theorem]{Вопрос}
\newtheorem{lemma}[theorem]{Лемма}

\else
\theoremstyle{definition}
\theoremstyle{remark}
\fi
\newtheorem{convention}[theorem]{Соглашение}
\newtheorem{xca}[theorem]{Exercise}
\ifx\PrintBook\undefined
\newcommand{\BibTitle}{%
\section*{Список литературы}%
}
\else
\newcommand{\BibTitle}{%
\chapter*{Список литературы}%
}
\fi
\else
\newcommand\TheLanguage{English}%
\author{Aleks Kleyn}
\ifx\setCACAA\undefined
\newtheorem{theorem}{Theorem}[section]
\newtheorem{corollary}[theorem]{Corollary}
\newtheorem{example}[theorem]{Example}
\newtheorem{definition}[theorem]{Definition}
\newtheorem{remark}[theorem]{Remark}

\newtheorem{lemma}[theorem]{Lemma}
\else
\theoremstyle{definition}
\theoremstyle{remark}
\fi
\newtheorem{convention}[theorem]{Convention}

\ifx\PrintBook\undefined
\newcommand{\BibTitle}{%
\section*{References}%
}
\else
\newcommand{\BibTitle}{%
\chapter*{References}%
}
\fi
\fi
\newcommand\CurrentLanguage{\TheLanguage.}%
\newcommand\xRefDef[1]
{
\wRefDef{#1}{\TheLanguage}
\NameDef{xRefDef#1}{}%
}

\theoremstyle{definition}
\theoremstyle{remark}
\renewenvironment{proof}[1][\proofname]
{\par{\sc #1. }}{\qed}%

\ifx\PrintBook\undefined
\numberwithin{Hfootnote}{section}
\else
\numberwithin{section}{chapter}
\numberwithin{footnote}{chapter}
\numberwithin{Hfootnote}{chapter}
\fi

\ifx\Presentation\undefined
\numberwithin{equation}{section}
\numberwithin{figure}{section}
\numberwithin{table}{section}
\numberwithin{Item}{section}
\fi

\makeatletter
\newcommand\org@maketitle{}
\let\org@maketitle\maketitle
\def\maketitle{%
\hypersetup{pdftitle={\@title}}%
\hypersetup{pdfauthor={\authors}}%
\hypersetup{pdfsubject=\@keywords}%
\ifx\UseRussian\Defined
\pdfbookmark[1]{\@title}{TitleRussian}
\else
\pdfbookmark[1]{\@title}{TitleEnglish}
\fi
\org@maketitle
}
\def\make@stripped@name#1{%
\begingroup
\escapechar\m@ne
\global\let\newname\@empty
\protected@edef\Hy@tempa{\CurrentLanguage #1}%
\edef\@tempb{%
\noexpand\@tfor\noexpand\Hy@tempa:=%
\expandafter\strip@prefix\meaning\Hy@tempa
}%
\@tempb\do{%
\if\Hy@tempa\else
\if\Hy@tempa\else
\xdef\newname{\newname\Hy@tempa}%
\fi
\fi
}%
\endgroup
}%
\newenvironment{enumBib}{%
\BibTitle
\advance\@enumdepth \@ne
\edef\@enumctr{enum\romannumeral\the\@enumdepth}\list
{\csname biblabel\@enumctr\endcsname}{\usecounter
{\@enumctr}\def\makelabel##1{\hss\llap{\upshape##1}}}
}{%
\endlist
}

\makeatletter

\newcommand{\BiblioItem}[2]
{
\def\Semafor{off}
\@ifundefined{\LanguagePrefix ViewCite#1}{}{%
\def\Semafor{on}%
}%
\ifx\Semafor\ValueOff
\@ifundefined{xRefDef#1}{}{%
\def\Semafor{on}%
}%
\fi
\ifx\Semafor\ValueOn
\ifx\IndexState\ValueOff
\begin{enumBib}
\def\IndexState{on}
\fi
\item \label{\LanguagePrefix bibitem: #1}#2%
\fi
}
\makeatother
\newcommand{\OpenBiblio}
{
\def\IndexState{off}
}
\newcommand{\CloseBiblio}
{
\ifx\IndexState\ValueOn
\end{enumBib}
\def\IndexState{off}
\fi
}

\makeatletter
\def\StartCite{[}%
\def\citeBib#1{\protect\showCiteBib#1,endCite,}%
\def\endCite{endCite}%
\def\showCiteBib#1,{\def\temp{#1}%
\ifx\temp\endCite
]%
\def\StartCite{[}%
\else
\StartCite\LanguagePrefix \ref{\LanguagePrefix bibitem: #1}%
\@ifundefined{\LanguagePrefix ViewCite#1}{%
\NameDef{\LanguagePrefix ViewCite#1}{}%
}{%
}%
\def\StartCite{, }%
\expandafter\showCiteBib%
\fi}%
\makeatother

\newcommand{\arp}{\ar @{-->}}

\newcommand\Bundle[1]{{\mathbb #1}}
\newcommand{\bundle}[4]%
{%
\def\tempa{}%
\def\tempb{#3}%
\def\tempc{#1}%
\ifx\tempa\tempb%
\ifx\tempa\tempc%
#2%
\else%
\xymatrix{#2:#1\arp[r]&#4}%
\fi%
\else%
\ifx\tempa\tempc%
#2[#3]%
\else%
\xymatrix{#2[#3]:#1\arp[r]&#4}%
\fi%
\fi%
}%
\makeatletter
\newcommand{\AddIndex}[2]%
{%
\@ifundefined{RefIndex#2}{%
\xNameDef{RefIndex#2}{:}%
\@namedef{LabelIndex}{\label{index: #2::}}%
}{%
\addtocounter{Index}{1}%
\xNameDef{RefIndex#2}{\@nameuse{RefIndex#2},\arabic{Index}}%
\@namedef{LabelIndex}{\label{index: #2:\arabic{Index}}}%
}%
\@nameuse{LabelIndex}%
{\bf #1}%
}%

\newcommand{\Index}[2]%
{%
\@ifundefined{RefIndex#2}{%
\def\Semafor{off}%
}{%
\def\Semafor{on}%
}%
\ifx\Semafor\ValueOn%
\def\tempa{}%
\def\tempb{#2}%
\ifx\IndexState\ValueOff%
\ifx\setCACAA\undefined
\begin{theindex}%
\else
\section*{\indexname}%
\begin{itemize}
\fi
\def\IndexState{on}%
\fi%
\ifx\IndexSpace\ValueOn%
\indexspace%
\def\IndexSpace{off}%
\fi%
\item #1%
\ifx\tempa\tempb%
\else%
\edef\PageRefs{\@nameuse{RefIndex#2}}
\def\Sep{\ }%
\@for\PageRef:=\PageRefs\do{%
\Sep
\pageref{index: #2:\PageRef}%
\def\Sep{,\ }%
}%
\fi%
\fi%
}%

\newcommand{\Symb}[4]%
{%
\def\Semafor{off}%
\@ifundefined{ViewSymbol#2}{%
\@ifundefined{ViewSymbol#2(#3-#4)}{%
}{%
\def\Semafor{on}
\edef\ThisSymbol{#2(#3-#4)}%
}%
}{%
\def\Semafor{on}%
\edef\ThisSymbol{#2}%
}%
\ifx\Semafor\ValueOn%
\ifx\IndexState\ValueOff%
\ifx\setCACAA\undefined
\begin{theindex}%
\else
\section*{\indexname}%
\begin{itemize}
\fi
\def\IndexState{on}%
\fi%
\ifx\IndexSpace\ValueOn%
\indexspace%
\def\IndexSpace{off}%
\fi%
\edef\Symbols{\@nameuse{ViewSymbol\ThisSymbol}}%
\@for\Symbol:=\Symbols\do{%
\edef\Temp{ViewSymbol\ThisSymbol:::\Symbol}%
\item $\displaystyle\textcolor{SymbColor}{\@nameuse{\Temp}}$
\ \ #1
\edef\PageRefs{\@nameuse{RefSymbol\ThisSymbol===\Symbol}}
\def\Sep{}%
\@for\PageRef:=\PageRefs\do{%
\Sep
\pageref{symbol: \ThisSymbol:\PageRef}%
\def\Sep{,\ }%
}%
}%
\fi
}

\newcommand{\Symba}[2]
{
\def\Semafor{off}
\@ifundefined{ViewSymbol#2}{%
}{%
\def\Semafor{on}
}%
\ifx\Semafor\ValueOn
\ifx\IndexState\ValueOff
\begin{theindex}
\def\IndexState{on}
\fi
\ifx\IndexSpace\ValueOn
\indexspace
\def\IndexSpace{off}
\fi
\item $\displaystyle\@nameuse{ViewSymbol#2}$\ \ #1
\edef\PageRefs{\@nameuse{RefSymbol#2}}
\def\Sep{}%
\@for\PageRef:=\PageRefs\do{%
\Sep
\pageref{symbol: #2:\PageRef}%
\def\Sep{,\ }%
}%
\fi
}

\makeatother

\newcommand{\SetIndexSpace}%
{%
\def\IndexSpace{on}%
}%

\newcommand{\OpenIndex}
{
\def\IndexState{off}
}
\newcommand{\CloseIndex}
{
\ifx\IndexState\ValueOn
\ifx\setCACAA\undefined
\end{theindex}
\else
\end{itemize}
\fi
\def\IndexState{off}
\fi
}

\def\LastMemo{LastMemo}%
\def\MemoList#1//{\def\temp{#1}%
\ifx\temp\LastMemo
\else%
\setlength{\parindent}{5mm}
\par
\BlueText{#1}%
\expandafter\MemoList%
\fi%
}     



%% file: Calculus.16.Absract.English.tex

\input{Calculus.16.Absract.Eq}

\begin{abstract}
Let $A$, $B$ be Banach $D$\Hyph algebras.
The map
\ShowEq{f:A->B abs}
is called differentiable
on the set $U\subset A$,
if at every point $x\in U$
the increment of map
$f$ can be represented as
\ShowEq{derivative of map def}
where
\ShowEq{df:A->B abstract}
is linear map and
\ShowEq{map o:A->B}
is such continuous map that
\ShowEq{derivative of map, algebra, 2}
Linear map
$\displaystyle\frac{d f(x)}{d x}$
is called
derivative of map
$f$.

I considered differential forms in Banach Algebra.
Differential form
\ShowEq{omega in A->B}
is defined by map
\ShowEq{g:A->BxB}
If the map $g$,
is derivative of the map
\ShowEq{f:A->B abs},
then the map $f$ is called indefinite integral of the map $g$
\ShowEq{f=int g abs}
Then, for any $A$\Hyph numbers $a$, $b$,
we define definite integral by the equality
\ShowEq{definite integral omega ab}
for any path $\gamma$ from $a$ to $b$.
\end{abstract}

%% file: Calculus.16.Absract.Eq.tex

\DefEq
{
$f:A\rightarrow B$
}
{f:A->B abs}

\DefEq
{
\[
f(x+dx)-f(x)
=\frac{d f(x)}{d x}\circ dx
+o(dx)
\]
}
{derivative of map def}

\AddEq{df:A->B abstract}
{
\[
\frac{d f(x)}{d x}:A\rightarrow B
\]
}

\DefEq
{
$o:A\rightarrow B$
}
{map o:A->B}

\DefEq
{
\[
\lim_{a\rightarrow 0}\frac{\|o(a)\|_B}{\|a\|_A}=0
\]
}
{derivative of map, algebra, 2}

\DefEq
{
\[
f(x)=\int g(x)\circ dx=\int\omega
\]
}
{f=int g abs}

\DefEq
{
$\omega\in\mathcal {LA}(D;A\rightarrow B)$
}
{omega in A->B}

\DefEq
{
$g:A\rightarrow B\otimes B$,
$\omega=g\circ dx$.
}
{g:A->BxB}

\DefEq
{
\[
\int_a^b\omega=\int_{\gamma}\omega
\]
}
{definite integral omega ab}

%% file: Preface.1610.309618526.English.tex
\input{Preface.1610.309618526.Eq}

\ePrints{1610.309618526,CACAA.06}%
\ifx\Semafor\ValueOn%
\Chapter{Preface}
\fi

\Section{Integration in \texorpdfstring{$D$}{D}-algebra}

When professor Svetlin asked me if I can solve differential equations
over quaternion algebra my first answer was that I am not ready.
Diversity of differential equations of first order and
identification of the derivative of any order with function
of real field make this task extremely difficult.

However later I recalled that, in the paper
\citeBib{0812.4763},
I considered the operation inverse to differentiation.
I considered the indefinite integral
\ShowEq{f=indefinite integral g}
as solution of differential equation
\ShowEq{partial f=g}
and I represented this solution as Taylor series.
So, when I was writing the paper
\citeBib{CACAA.05.001},
I decided explicitly use integral notation.

If this had happened in 2012,
I would have stopped at this point and would start to study another subject.
\ePrints{4975-6381}
\ifx\Semafor\ValueOff
But in 2013 I had published the paper
\ePrints{CACAA.06}
\ifx\Semafor\ValueOff
\citeBib{1310.5591}.
\else
\citeBib{CACAA.04.001}.
\fi
\else
But in 2013 I started to write the book
\citeBib{5410-9916}.
\fi
And I had seen a serious problem.

I have two definitions of integral:
namely Lebesgue integral and indefinite integral.
And they have different formats of integrand.
In an indefinite integral,
we consider the representation of $D$\Hyph module
\ShowEq{AoA}
in $D$\Hyph algebra $A$
and integrand is an action of a tensor
\ShowEq{a in AA}
over differential.
In Lebesgue integral,
we consider the representation of algebra with measure
in $D$\Hyph algebra $A$.
What is the relation between
Lebesgue integral and indefinite integral?

At first sight, there is no answer to this question.
In the calculus over real field,
we use the same map in both cases.
And the relation between
Lebesgue integral and indefinite integral
can be expressed by either of following statements.

\begin{itemize}
\item
\ShowEq{int f=int 0x f}
\item
If
\ShowEq{g=int f}
then
\ShowEq{int ab f=g(b)-g(a)}
\end{itemize}

If we consider Lebesgue integral in $D$\Hyph algebra,
what does it mean to get an integral from $a$ to $b$?
The simplest solution is to use a path from $a$ to $b$.
In the paper
\citeBib{CACAA.05.001},
I considered the simplest path
\ShowEq{path abt}
However, the question remained open for arbitrary path.
If Lebesgue integral depends on the path,
then we need to study a dependence of Lebesgue integral from the path
instead of the relation between
Lebesgue integral and indefinite integral.
But this is absolutely different theory.

There exist maps which have
Lebesgue integral dependent on the path.
It follows from the remark
\ref{remark: differential equation 3x^2 does not possess a solution}.
However, I had great interest in the study of maps which have indefinite integral.

When I discussed this question with the professor Sudbery,
I suddenly found simple proof
of the statement that Lebesgue integral
of the derivative of the map does not depend on path (the theorem
\RefTheorem{Lebesgue Integral along Path}).
However I was not satisfied by this statement.
I knew that new adventure is waiting for me.

\Section{Differential Form}

The question about possibility to solve differential equation
\ShowEq{df=g partial}
remained opened.
If solution does not exist,
then we can assume that Lebesgue integral of the map $g$
depends on the path.
In the differential geometry,
the map $f$ is called anholonomic.
Anholonomic maps are important in general relativity.
And I could not ignore this problem.

However if the Lebesgue integral depends on path,
then, of course, the question about integral over a closed circuit arises.
At this point, I realized that the map $g$ is differential form.

The study of differential forms demanded from me to expand the scope of my research.
In the definition
\RefDefinition[1601.03259]{differentiable map, algebra},
I considered differentiable map into Banach $D$\Hyph algebra.
However this definition is still true in case of Banach $D$\Hyph module.
To write the derivative of the map into Banach $D$\Hyph module,
we need a basis of the module and we write the Jacobian matrix of this map.
However, we can write the derivative of the map $f$ into Banach $D$\Hyph algebra $A$ analytically
using tensor
\ShowEq{df/dx in AxA}
I considered relation between two forms of representation of linear map
in the theorem
\RefTheorem{coordinates of map A1 A2, algebra}.

Moreover, we have the opportunity to write analytically a map into Banach $D$\Hyph module
when we consider $D$\Hyph module of polylinear maps
into $D$\Hyph algebra $B$
(the theorem
\RefTheorem{dfa1p/dx=df/dx a1p}).
In this case we can use
product in $D$\Hyph algebra $B$
to represent linear map in $D$\Hyph module
\ShowEq{L(AnB)}

\ePrints{1610.309618526,CACAA.06}%
\ifx\Semafor\ValueOn%
When I started my research in July of 2008,
it was a blind search.
My research was based on the definition of linear map
which I learned in
\citeBib{0701.238}.
I wanted to clarify some details;
and I chose calculus because
this subject was interesting for me,
as well I knew that I need calculus
to study differential geometry.
At that time, I knew only map linear from left
and map linear from right.
Left and right derivatives looked like a layered cake;
the increment of the argument was its stuffing.

Without a doubt, it was Gateaux derivative.
And I used symbol $\partial$
to distinguish this derivative.
But I was not satisfied with this definition;
starting from 2010, I started to consider linear map
of $D$\Hyph algebra $A$ as representation of algebra $A\otimes A$
in algebra $A$. However, I still used old notation.

From the definition
\RefDefinition{linear map from A1 to A2, algebra},
it follows that the definition
\RefDefinition{differentiable map, algebra}
of derivative of map into Banach algebra
is not different from standard definition
of derivative of map into Banach module.
So I decided (starting from this paper)
to use standard notation for derivative.
\fi


\ePrints{1601.03259}%
\ifx\Semafor\ValueOn%
\begin{flushright}
January, 2017
\end{flushright}
\fi%

%% file: Preface.1610.309618526.Eq.tex

\DefEq
{
\[
\frac{d f(x)}{d x}=g(x)
\]
}
{df=g partial}

\AddEq{df/dx in AxA}
{
$\displaystyle\frac{df(x)}{dx}\in A\otimes A$.
}

\DefEq
{
\[
f(x)=\int g(x)\circ dx
\]
}
{f=indefinite integral g}

\DefEq
{
$a\in A\otimes A$
}
{a in AA}

\DefEq
{
$\mathcal L(D;A^n\rightarrow B)$.
}
{L(AnB)}

\DefEq
{
$\displaystyle \int_a^bf(x)dx=g(b)-g(a)$\,.
}
{int ab f=g(b)-g(a)}

\DefEq
{
$\displaystyle g(x)=\int f(x)dx$\,,
}
{g=int f}

\DefEq
{
\[x=a+t(b-a)\]
}
{path abt}

\DefEq
{
$\displaystyle\int f(x)dx=\int_a^x f(x)dx$
}
{int f=int 0x f}

\DefEq
{
$A\otimes A$
}
{AoA}

\DefEq
{
\[
\frac{d f(x)}{d x}\circ dx=g(x)\circ dx
\]
}
{partial f=g}

%% file: Preliminary.English.tex
\def\Preliminary{on}
\input{Preliminary.Eq}

\Chapter{Preliminary Definitions}

\ifx\PrintBook\undefined
This section
\else
This chapter
\fi
contains definitions and theorems
which are necessary for an understanding of the text of this
\ifx\PrintBook\undefined
paper.
\else
book.
\fi
So the reader may read the statements from this
\ifx\PrintBook\undefined
section
\else
chapter
\fi
in process of reading the main text of the
\ifx\PrintBook\undefined
paper.
\else
book.
\fi

\ePrints{4975-6381,1506.00061,MQuater,5148-4632,MRepro,BRepro}
\Items{1601.03259,5410-9916}
\ifx\Semafor\ValueOn
\Section{Universal Algebra}

\ePrints{5410-9916}
\ifx\Semafor\ValueOff
\begin{theorem}
\labelTheorem{maps and kernel equivalence}
Let $N$ be equivalence
on the set $A$.
Consider category $\mathcal A$ whose objects are
maps\,\footnote{
The statement of lemma is
similar to the statement on p. \citeBib{Serge Lang}-119.}
\ShowEq{maps category}
We define morphism $f_1\rightarrow f_2$
to be map $h:S_1\rightarrow S_2$
making following diagram commutative
\ShowEq{maps category, diagram}
The map
\ShowEq{maps category, universal}
is universally repelling in the category $\mathcal A$.\,\footnote{See
definition of universal object of category in definition
on p. \citeBib{Serge Lang}-57.}
\end{theorem}
\begin{proof}
Consider diagram
\ShowEq{maps category, universal, diagram}
\ShowEq{maps category, universal, ker}
From the statement
\EqRef{maps category, universal, ker}
and the equality
\ShowEq{maps category 1}
it follows that
\ShowEq{maps category 2}
Therefore, we can uniquely define the map $h$
using the equality
\ShowEq{maps category, h}
\end{proof}
\fi

\begin{definition}
\labelDefinition{Cartesian power}
{\it
For any sets\,\footnote{
I follow the definition from the example (iv) on the
page \citeBib{Cohn: Universal Algebra}\Hyph 5.
}
$A$, $B$,
\AddIndex{Cartesian power}{Cartesian power}
\ShowEq{Cartesian power}
is the set of maps
\ShowEq{f:A->B}fAB
}
\qed
\end{definition}

\begin{definition}
\labelDefinition{operation on set}
{\it
For any $n\ge 0$, a map\,\footnote{
\ePrints{MVector,5148-4632,MRepro,BRepro}
\ifx\Semafor\ValueOn
Definitions
\RefDefinition{operation on set},
\RefDefinition{set is closed with respect to operation}
follow
\else
Definition
\RefDefinition{operation on set}
follows
\fi
the definition in the example (vi) on the page
page \citeBib{Cohn: Universal Algebra}\Hyph 13.
}
\ShowEq{o:An->A}
is called
\AddIndex{$n$\Hyph ary operation on set}{n-ary operation on set}
$A$ or just
\AddIndex{operation on set}{operation on set}
$A$.
For any
\EqParm{b1n in B}{B=A}
we use either notation
\ShowEq{a1no=oa1n}
to denote image of map $\omega$.
}
\qed
\end{definition}

\begin{remark}
\labelRemark{o in AAn}
According to definitions
\RefDefinition{Cartesian power},
\RefDefinition{operation on set},
$n$\Hyph ari operation
\ShowEq{o in AAn}
\qed
\end{remark}

\begin{definition}
{\it An
\AddIndex{operator domain}{operator domain}
is the set of operators\,\footnote{
I follow the definition (1),
page \citeBib{Cohn: Universal Algebra}\Hyph 48.
}
\ShowEq{operator domain}
with a map
\ShowEq{a:O->N}
If
\ShowEq{o in O},
then $a(\omega)$ is called the
\AddIndex{arity}{arity}
of operator $\omega$. If
\ShowEq{a(o)=n}
then operator $\omega$ is called $n$\Hyph ary.
\ShowEq{set of n-ary operators}
We use notation
\ShowEq{set of n-ary operators =}
for the set of $n$\Hyph ary operators.
}
\qed
\end{definition}

\begin{definition}
\labelDefinition{Omega-algebra}
{\it
Let $A$ be a set. Let $\Omega$ be an operator domain.\,\footnote{
I follow the definition (2),
page \citeBib{Cohn: Universal Algebra}\Hyph 48.
}
The family of maps
\ShowEq{O(n)->AAn}
is called $\Omega$\Hyph algebra structure on $A$.
The set $A$ with $\Omega$\Hyph algebra structure is called
\AddIndex{$\Omega$\Hyph algebra}{Omega-algebra}
\ShowEq{Omega-algebra}
or
\AddIndex{universal algebra}{universal algebra}.
The set $A$ is called
\AddIndex{carrier of $\Omega$\Hyph algebra}{carrier of Omega-algebra}.
}
\qed
\end{definition}

The operator domain $\Omega$ describes a set of $\Omega$\Hyph algebras.
An element of the set $\Omega$ is called operator,
because an operation assumes certain set.
According to the remark
\ref{remark: o in AAn}
and the definition
\RefDefinition{Omega-algebra},
for each operator
\EqParm{omega n ari}{=c}
we match $n$\Hyph ary operation $\omega$ on $A$.

\ePrints{MVector,5148-4632,MRepro,BRepro}
\ifx\Semafor\ValueOn

\begin{theorem}
\labelTheorem{Cartesian power is universal algebra}
Let the set $B$ be $\Omega$\Hyph algebra.
Then the set $B^A$ of maps
\ShowEq{f:A->B}fAB
also is $\Omega$\Hyph algebra.
\end{theorem}
\begin{proof}
Let
\EqParm{omega n ari}{=.}
For maps
\ShowEq{f1n in B**A}
we define the operation $\omega$ by the equality
\DrawEq{f1n omega=}{}
\end{proof}
\fi

\ePrints{MVector,5148-4632,MRepro,BRepro}
\ifx\Semafor\ValueOn
\begin{definition}
\labelDefinition{set is closed with respect to operation}
{\it
Let
\ShowEq{B subset A}.
Since, for any
\EqParm{b1n in B}{B=B}
\ShowEq{b1no in B}
then we say that $B$
\AddIndex{is closed with respect to}{set is closed with respect to operation}
$\omega$ or that $B$
\AddIndex{admits operation}{set admits operation}
$\omega$.
}
\qed
\end{definition}

\begin{definition}
$\Omega$\Hyph algebra $B_{\Omega}$ is
\AddIndex{subalgebra}{subalgebra of Omega-algebra}
of $\Omega$\Hyph algebra $A_{\Omega}$
if following statements are true\,\footnote{
I follow the definition on
page \citeBib{Cohn: Universal Algebra}\Hyph 48.
}
\StartLabelItem
\begin{enumerate}
\item
\ShowEq{B subset A}.
\item
if operator
\ShowEq{o in O}
defines operations $\omega_A$ on $A$ and $\omega_B$ on $B$, then
\ShowEq{oAB=oB}
\end{enumerate}
\qed
\end{definition}
\fi

\begin{definition}
\labelDefinition{homomorphism}
{\it
Let $A$, $B$ be $\Omega$\Hyph algebras and
\EqParm{omega n ari}{=.}
The map\,\footnote{
I follow the definition on
page \citeBib{Cohn: Universal Algebra}\Hyph 49.
}
\ShowEq{f:A->B}fAB
\AddIndex{is compatible with operation}{map is compatible with operation}
$\omega$, if, for all
\EqParm{b1n in B}{B=A}
\ShowEq{afo=aof}
The map $f$ is called
\AddIndex{homomorphism}{homomorphism}
from $\Omega$\Hyph algebra $A$ to $\Omega$\Hyph algebra $B$,
if $f$ is compatible with each
\ShowEq{o in O}.
\ePrints{MRepro,BRepro}%
\ifx\Semafor\ValueOn%
We use notation
\ShowEq{set of homomorphisms}
for the set of homomorphisms
from $\Omega$\Hyph algebra $A$ to $\Omega$\Hyph algebra $B$.
\fi
}
\qed
\end{definition}

\ePrints{MRepro,BRepro}
\ifx\Semafor\ValueOn
\begin{theorem}
\labelTheorem{Hom empty A B=B**A}
Since operator domain is empty,
then a homomorphism
from $\Omega$\Hyph algebra $A$ to $\Omega$\Hyph algebra $B$
is a map
\ShowEq{f:A->B}fAB
Therefore,
\ShowEq{Hom empty A B=B**A}
\end{theorem}
\begin{proof}
The theorem follows from definitions
\RefDefinition{Cartesian power},
\RefDefinition{homomorphism}.
\end{proof}
\fi

\ePrints{MVector,5148-4632,MRepro,BRepro}
\ifx\Semafor\ValueOn
\begin{definition}
\labelDefinition{isomorphism}
{\it
Homomorphism $f$ is called\,\footnote{
I follow the definition on
page \citeBib{Cohn: Universal Algebra}\Hyph 49.
}
\AddIndex{isomorphism}{isomorphism}
between $A$ and $B$, if correspondence $f^{-1}$ is homomorphism.
If there is an isomorphism between $A$ and $B$, then we
say that $A$ and $B$ are isomorphic and write
\ShowEq{isomorphic}
An injective homomorphis is called
\AddIndex{monomorphism}{monomorphism}.
A surjective homomorphis is called
\AddIndex{epimorphism}{epimorphism}.
}
\qed
\end{definition}
\fi

\begin{definition}
\labelDefinition{endomorphism}
{\it
A homomorphism in which source and target are the same algebra is called
\AddIndex{endomorphism}{endomorphism}.
We use notation
\ShowEq{set of endomorphisms}
for the set of endomorphisms
of $\Omega$\Hyph algebra $A$.
\ePrints{4975-6381,1506.00061,MQuater,1601.03259,MVector,5148-4632}%
\ifx\Semafor\ValueOff%
An endomorphism which is also an isomorphism is called
\AddIndex{automorphism}{automorphism}.
\fi
}
\qed
\end{definition}

\ePrints{MRepro,BRepro}
\ifx\Semafor\ValueOn

\begin{theorem}
\labelTheorem{End A=Hom AA}
\ShowEq{End A=Hom AA}
\end{theorem}
\begin{proof}
The theorem follows from the definitions
\RefDefinition{homomorphism},
\RefDefinition{endomorphism}.
\end{proof}

\begin{theorem}
\labelTheorem{End empty A=A**A}
Since operator domain is empty,
then an endomorphism of the set $A$
is a map
\ShowEq{t:A->A}
Therefore,
\ShowEq{End empty A=A**A}
\end{theorem}
\begin{proof}
The theorem follows from the theorems
\RefTheorem{Hom empty A B=B**A},
\RefTheorem{End A=Hom AA}.
\end{proof}
\fi

\ePrints{4975-6381,1506.00061,MQuater,1601.03259,MVector,5148-4632,5410-9916}%
\ifx\Semafor\ValueOff%
\begin{definition}
{\it
If there is a monomorphism from $\Omega$\Hyph algebra $A$
to $\Omega$\Hyph algebra $B$, then we say that
\AddIndex{$A$ can be embeded in $B$}{can be embeded}.
}
\qed
\end{definition}

\begin{definition}
{\it
If there is an epimorphism from $A$ to $B$, then $B$ is called
\AddIndex{homomorphic image}{homomorphic image} of algebra $A$.
}
\qed
\end{definition}
\fi
\fi

\ePrints{MVector,5148-4632,4975-6381,5410-9916}
\ifx\Semafor\ValueOn
\begin{convention}
Element of
$\Omega$\Hyph algebra
$A$ is called
\AddIndex{$A$\Hyph number}{A number}.
For instance, complex number is also called
$C$\Hyph number, and quaternion is called $H$\Hyph number.
\qed
\end{convention}
\fi

\ePrints{MRepro,BRepro}
\ifx\Semafor\ValueOn

\Section{Cartesian Product of Universal Algebras}

\ShowDefinition{product in category}

\begin{example}
Let \(\mathcal S\) be the category of sets.\,\footnote{
See also the example in
\citeBib{Serge Lang},
page 59.
}
According to the definition
\RefDefinition{product in category},
Cartesian product
\ShowEq{Cartesian product of sets}
of family of sets
\ShowEq{Ai iI}
and family of projections on the \(i\)\Hyph th factor
\ShowEq{projection on i factor}
are product in the category \(\mathcal S\).
\qed
\end{example}

\begin{theorem}
\labelTheorem{product exists in category of Omega algebras}
The product exists in the category \(\mathcal A\) of \(\Omega\)\Hyph algebras.
Let \(\Omega\)\Hyph algebra \(A\)
and family of morphisms
\ShowEq{p:A->Ai i in I}
be product in the category \(\mathcal A\).
Then
\StartLabelItem
\begin{enumerate}
\item
The set \(A\) is Cartesian product
of family of sets
\ShowEq{Ai iI}
\item
The homomorphism of \(\Omega\)\Hyph algebra
\ShowEq{projection on i factor}
is projection on \(i\)\Hyph th factor.
\item
We can represent any \(A\)\Hyph number $a$
as tuple
\ShowEq{tuple represent A number}
of \(A_i\)\Hyph numbers.
\labelItem{tuple represent A number}
\item
Let
\ShowEq{omega in Omega}{}{}
be n\Hyph ary operation.
Then operation $\omega$ is defined componentwise
\ShowEq{operation is defined componentwise}
where
\ShowEq{a=ai 1n}.
\labelItem{operation is defined componentwise}
\end{enumerate}
\end{theorem}
\begin{proof}
Let
\ShowEq{Cartesian product of sets}
be Cartesian product
of family of sets
\ShowEq{Ai iI}
and, for each \iI, the map
\ShowEq{projection on i factor}
be projection on the \(i\)\Hyph th factor.
Consider the diagram of morphisms in category of sets $\mathcal S$
\ShowEq{operation is defined componentwise, diagram}
where the map $g_i$ is defined by the equality
\ShowEq{gi()=}
According to the definition
\RefDefinition{product in category},
the map $\omega$ is defined uniquely from the set of diagrams
\EqRef{operation is defined componentwise, diagram}
\ShowEq{omega(ai)=(omega ai)}
The equality
\EqRef{operation is defined componentwise}
follows from the equality
\EqRef{omega(ai)=(omega ai)}.
\ifx\texFuture\Defined

\begin{lemma}
\labelLemma{tuple represent A number}
{\it
Let \(\Omega\)\Hyph algebra \(A\)
and family of morphisms
\ShowEq{projection on i factor}
be product in the category \(\mathcal A\).
We can represent any \(A\)\Hyph number $a$
as tuple
\ShowEq{tuple represent A number}
of \(A_i\)\Hyph numbers.
}
\end{lemma}

{\sc Proof.}
Let
\ShowEq{a ne b}
be \(A\)\Hyph numbers such that for any \iI
\ShowEq{p(a)=p(b)}
Assuming that the unequal \(A\)\Hyph numbers
have different expansion as tuple,
то мы можем определить отношение эквивалентности,
фактор по которому скорее кандидат на произведение.
\hfill\(\odot\)

From the equality
\ShowEq{pj(a)=aj}
it follows that the map \(p_j\) is projection.
\fi
\end{proof}

\begin{definition}
If \(\Omega\)\Hyph algebra \(A\)
and family of morphisms
\ShowEq{p:A->Ai i in I}
is product in the category \(\mathcal A\),
then \(\Omega\)\Hyph algebra \(A\) is called
\AddIndex{direct}{direct product of Omega algebras}
or
\AddIndex{Cartesian product of \(\Omega\)\Hyph algebras}
{Cartesian product of Omega algebras}
\ShowEq{Ai iI}.
\qed
\end{definition}

\begin{theorem}
\labelTheorem{map from product into product}
Let set \(A\) be
Cartesian product of sets
\ShowEq{Ai iI}
and set \(B\) be
Cartesian product of sets
\ShowEq{Bi iI}.
For each \iI, let
\ShowEq{f:A->B i}
be the map from the set $A_i$ into the set $B_i$.
For each \iI, consider commutative diagram
\ShowEq{homomorphism of Cartesian product of Omega algebras diagram}
where maps
\ShowEq{pi p'i}
are projection on the \(i\)\Hyph th factor.
The set of commutative diagrams
\EqRef{homomorphism of Cartesian product of Omega algebras diagram}
uniquely defines map
\ShowEq{f:A->B}fAB
\DrawEq{f:A->B=}{}
\end{theorem}
\begin{proof}
For each \iI, consider commutative diagram
\ShowEq{homomorphism of Cartesian product of Omega algebras}
Let \(a\in A\).
According to the statement
\RefItem{tuple represent A number},
we can represent \(A\)\Hyph number \(a\)
as tuple of \(A_i\)\Hyph numbers
\ShowEq{a=p(a)i}
Let
\ShowEq{b=f(a)}
According to the statement
\RefItem{tuple represent A number},
we can represent \(B\)\Hyph number \(b\)
as tuple of \(B_i\)\Hyph numbers
\ShowEq{b=p(b)i}
From commutativity of diagram $(1)$
and from equalities
\EqRef{b=f(a)},
\EqRef{b=p(b)i},
it follows that
\ShowEq{b=g(a)i}
From commutativity of diagram $(2)$
and from the equality
\EqRef{a=p(a)i},
it follows that
\ShowEq{b=f(a)i}
\end{proof}

\begin{theorem}
\labelTheorem{homomorphism of Cartesian product of Omega algebras}
Let \(\Omega\)\Hyph algebra \(A\) be
Cartesian product of \(\Omega\)\Hyph algebras
\ShowEq{Ai iI}
and \(\Omega\)\Hyph algebra \(B\) be
Cartesian product of \(\Omega\)\Hyph algebras
\ShowEq{Bi iI}.
For each \iI,
let the map
\ShowEq{f:A->B i}
be homomorphism of \(\Omega\)\Hyph algebra.
Then the map
\ShowEq{f:A->B}fAB
defined by the equality
\DrawEq{f:A->B=}{homomorphism}
is homomorphism of \(\Omega\)\Hyph algebra.
\end{theorem}
\begin{proof}
Let
\ShowEq{omega in Omega}{}{}
be n\Hyph ary operation.
Let
\ShowEq{a=ai 1n},
\ShowEq{b=bi 1n}.
From equalities
\EqRef{operation is defined componentwise},
\eqRef{f:A->B=}{homomorphism},
it follows that
\ShowEq{f:A->B omega}
\end{proof}

\begin{definition}
{\it
Equivalence on $\Omega$\Hyph algebra $A$,
which is subalgebra of $\Omega$\Hyph algebra $A^2$,
is called
\AddIndex{congruence}{congruence}\,\footnote{
I follow the definition on
page \citeBib{Cohn: Universal Algebra}\Hyph 57.
}
on $A$.
}
\qed
\end{definition}

\begin{theorem}[first isomorphism theorem]
\labelTheorem{first isomorphism theorem}
Let
\ShowEq{f:A->B}fAB
be homomorphism of $\Omega$\Hyph algebras with kernel $s$.
Then there is decomposition
\ShowEq{decomposition of map f}
\StartLabelItem
\begin{enumerate}
\item
The \AddIndex{kernel of homomorphism}{kernel of homomorphism}
\ShowEq{kernel of homomorphism}
is a congruence on $\Omega$\Hyph algebra $A$.
\labelItem{kernel of homomorphism}
\item
The set
\ShowEq{A/ker f}
is $\Omega$\Hyph algebra.
\labelItem{A/ker f is Omega-algebra}
\item
The map
\ShowEq{p:A->/ker}
is epimorphism and is called
\AddIndex{natural homomorphism}{natural homomorphism}.
\labelItem{natural homomorphism}
\item
The map 
\ShowEq{q:A/ker->f(A)}
is the isomorphism
\labelItem{q:A/ker->f(A) isomorphism}
\item
The map 
\ShowEq{r:f(A)->B}
is the monomorphism
\labelItem{r:f(A)->B monomorphism}
\end{enumerate}
\end{theorem}
\begin{proof}
The statement
\RefItem{kernel of homomorphism}
follows from the proposition II.3.4
(\citeBib{Cohn: Universal Algebra}, page 58).
Statements
\RefItem{A/ker f is Omega-algebra},
\RefItem{natural homomorphism}
follow from the theorem II.3.5
(\citeBib{Cohn: Universal Algebra}, page 58)
and from the following definition.
Statements
\RefItem{q:A/ker->f(A) isomorphism},
\RefItem{r:f(A)->B monomorphism}
follow from the theorem II.3.7
(\citeBib{Cohn: Universal Algebra}, page 60).
\end{proof}
\fi

\ePrints{MRepro,BRepro}
\ifx\Semafor\ValueOn
\Section{Semigroup}

Usually the operation
\ShowEq{omega 2 ari}
is called product
\ShowEq{abo=ab}
or sum
\ShowEq{abo=a+b}

\begin{definition}
{\it
Let $A$ be $\Omega$\Hyph algebra and
\ShowEq{omega 2 ari.}
$A$\Hyph number $e$ is called
\AddIndex{neutral element of operation}{neutral element of operation}
$\omega$, when for any $A$\Hyph number $a$ following equations are true
\ShowEq{left neutral element}
\ShowEq{right neutral element}
}
\qed
\end{definition}

\begin{definition}
{\it
Let $A$ be $\Omega$\Hyph algebra.
The operation
\ShowEq{omega 2 ari}
is called
\AddIndex{associative}{associative operation}
if the following equality is true
\ShowEq{associative operation}
}
\qed
\end{definition}

\begin{definition}
{\it
Let $A$ be $\Omega$\Hyph algebra.
The operation
\ShowEq{omega 2 ari}
is called
\AddIndex{commutative}{commutative operation}
if the following equality is true
\ShowEq{commutative operation}
}
\qed
\end{definition}

\begin{definition}
\labelDefinition{semigroup}
{\it
Let
\ShowEq{Omega=omega}
If the operation
\ShowEq{omega 2 ari}
is associative, then $\Omega$\Hyph algebra is called
\AddIndex{semigroup}{semigroup}.
If the operation in the semigroup is commutative,
then the semigroup is called
\AddIndex{Abelian semigroup}{Abelian semigroup}.
}
\qed
\end{definition}
\fi

\ePrints{4975-6381,1506.00061,MQuater,1601.03259,5148-4632}
\Items{MVector,5410-9916}
\ifx\Semafor\ValueOn
\Section{Representation of Universal Algebra}
\labelSection{Representation of Universal Algebra}

\begin{definition}
\labelDefinition{representation of algebra} 
Let the set $A_2$ be $\Omega_2$\Hyph algebra.
Let
the set of transformations
\ShowEq{End O2}{A_2}{}
be  $\Omega_1$\Hyph algebra.
The homomorphism
\ShowEq{representation of algebra}
of $\Omega_1$\Hyph algebra $A_1$ into
$\Omega_1$\Hyph algebra
\ShowEq{End O2}{A_2}{}
is called
\AddIndex{representation of $\Omega_1$\Hyph algebra}
{representation of algebra}
$A_1$ or
\AddIndex{$A_1$\Hyph representation}{A representation of algebra}
in $\Omega_2$\Hyph algebra $A_2$.
\qed
\end{definition}

We also use notation
\ShowEq{f:A->*B}f{A_1}{A_2}
to denote the representation of $\Omega_1$\Hyph algebra $A_1$
in $\Omega_2$\Hyph algebra $A_2$.

\begin{definition}
\labelDefinition{effective representation of algebra}
Let the map
\ShowEq{representation of algebra}
be an isomorphism of the $\Omega_1$\Hyph algebra $A_1$ into
\ShowEq{End O2}{A_2}.
Then the representation
\ShowEq{f:A->*B}f{A_1}{A_2}
of the $\Omega_1$\Hyph algebra $A_1$ is called
\AddIndex{effective}{effective representation}.\,\footnote{
See similar definition of effective representation of group in
\citeBib{Postnikov: Differential Geometry}, page 16,
\citeBib{Basic Concepts of Differential Geometry}, page 111,
\citeBib{Cohn: Algebra 1}, page 51
(Cohn calls such representation faithful).
}
\qed
\end{definition}

\ePrints{MVector,5148-4632}
\ifx\Semafor\ValueOn

\begin{theorem}
\labelTheorem{representation is effective}
The representation
\ShowEq{f:A->*B}g{A_1}{A_2}
is effective iff the statement
\ShowEq{a1 ne b1}
implies that there exists
\EqParm{a in A}{n=2,=z}
such that\,\footnote{
In case of group, the theorem
\RefTheorem{representation is effective}
has the following form.
{\it
The representation
\ShowEq{f:A->*B}g{A_1}{A_2}
is effective iff, for any $A_1$\Hyph number $a\ne e$,
there exists
\EqParm{a in A}{n=2,=z}
such that
\ShowEq{fam ne m}
}
}
\ShowEq{fam ne fbm}
\end{theorem}
\ProofTheorem{\RefRepresentation}{representation is effective}

\begin{theorem}
\labelTheorem{effective representation of the ring}
Let ring $D$ has unit $e$.
Representation
\ShowEq{f:A->*B}fDA
of the ring $D$
in an Abelian group $A$ is
\AddIndex{effective}{effective representation of ring}
iff
$a=0$ follows from equation $f(a)=0$.
\end{theorem}
\begin{proof}
We define the sum of transformations $f$ and $g$ of an Abelian group
according to rule
\ShowEq{(f+g)(a)=}
Therefore, considering the representation of the ring $D$ in
the Abelian group $A$, we assume
\ShowEq{sum of transformations of Abelian group, 1}

Suppose $a$, $b\in R$
cause the same transformation. Then
\ShowEq{representation of ring, 1}
for any $m\in A$.
From the equation
\EqRef{representation of ring, 1}
it follows that $a-b$ generates zero transformation
\ShowEq{representation of ring, 2}
Element $e+a-b$ generates an identity transformation.
Therefore, the representation $f$ is effective
iff $a=b$.
\end{proof}
\fi

\ePrints{1506.00061}
\ifx\Semafor\ValueOn
\begin{definition}
\labelDefinition{free representation of algebra}
{\it
The representation
\ShowEq{f:A->*B}g{A_1}{A_2}
of the $\Omega_1$\Hyph algebra $A_1$ is called
\AddIndex{free}{free representation},\,\footnote{
See similar definition of free representation of group in
\citeBib{Postnikov: Differential Geometry}, page 16.
}
if the statement
\ShowEq{faa=fba}
for any
\EqParm{a in A}{n=2,=z}
implies that $a=b$.
}
\qed
\end{definition}

\begin{theorem}
\labelTheorem{Free representation is effective}
Free representation is effective.
\end{theorem}
\ProofTheorem{\RefRepresentation}{Free representation is effective}

\begin{remark}
Representation of rotation group in affine space is effective.
However this representation is not free, since origin
is fixed point of every transformation.
\qed
\end{remark}

\begin{definition}
\labelDefinition{transitive representation of algebra}
The representation
\ShowEq{f:A->*B}g{A_1}{A_2}
of $\Omega_1$\Hyph algebra is called
\AddIndex{transitive}{transitive representation of algebra}\,\footnote{
See similar definition of transitive representation of group in
\citeBib{Basic Concepts of Differential Geometry}, page 110,
\citeBib{Cohn: Algebra 1}, page 51.
}
if, for any
\EqParm{ab in A}{A=2,=c}
there exists such $g$ that
\[a=f(g)(b)\]
The representation of $\Omega_1$\Hyph algebra is called
\AddIndex{single transitive}{single transitive representation of algebra}
if it is transitive and free.
\qed
\end{definition}

\begin{theorem}
\labelTheorem{Representation is single transitive iff}
Representation is single transitive iff for any $a, b \in A_2$
exists one and only one $g\in A_1$ such that $a=f(g)(b)$
\end{theorem}
\begin{proof}
The theorem follows from definitions \RefDefinition{free representation of algebra}
and \RefDefinition{transitive representation of algebra}.
\end{proof}

\begin{theorem}
\labelTheorem{single transitive representation generates algebra}
Let
\ShowEq{f:A->*B}g{A_1}{A_2}
be a single transitive representation
of $\Omega_1$\Hyph algebra $A_1$
in $\Omega_2$\Hyph algebra $A_2$.
There is the structure of $\Omega_1$\Hyph algebra on the set $A_2$.
\end{theorem}
\ProofTheorem{\RefRepresentation}{single transitive representation generates algebra}
\fi

\begin{definition}
\labelDefinition{morphism of representations of universal algebra}
Let
\ShowEq{f:A->*B}f{A_1}{A_2}
be representation of $\Omega_1$\Hyph algebra $A_1$
in $\Omega_2$\Hyph algebra $A_2$ and
\ShowEq{f:A->*B}g{B_1}{B_2}
be representation of $\Omega_1$\Hyph algebra $B_1$
in $\Omega_2$\Hyph algebra $B_2$.
For
\EqParm{i=1,2}{=c}
let the map
\ShowEq{ri:A->B}
be homomorphism of $\Omega_i$\Hyph algebra.
The matrix of maps
\ShowEq{map r,R}{r_1}{r_2}{}
such, that
\ShowEq{morphism of representations of universal algebra, definition, 2}
is called
\AddIndex{morphism of representations from $f$ into $g$}
{morphism of representations from f into g}.
We also say that
\AddIndex{morphism of representations of $\Omega_1$\Hyph algebra
in $\Omega_2$\Hyph algebra}
{morphism of representations of Omega1 algebra in Omega2 algebra} is defined.
\qed
\end{definition}

\begin{remark}
We may consider a pair of maps $r_1$, $r_2$ as map
\ShowEq{F:A1+A2->B1+B2}
such that
\ShowEq{F:A1+A2->B1+B2 1}
Therefore, hereinafter the matrix of maps
\ShowEq{map r,R}{r_1}{r_2}{}
also is called map.
\qed
\end{remark}

\ePrints{5148-4632,MVector}
\ifx\Semafor\ValueOn

\begin{remark}
Consider morphism of representations
\ShowEq{r12:A->B}rAB
We denote elements of the set $B_1$ by letter using pattern $b\in B_1$.
However if we want to show that $b$ is image of element
\EqParm{a in A1}{=c}
we use notation $\RedText{r_1(a)}$.
Thus equation
\ShowEq{r1(a)=r1(a)}
means that $\RedText{r_1(a)}$ (in left part of equation)
is image
\EqParm{a in A1}{=z}
(in right part of equation).
Using such considerations, we denote
element of set $B_2$ as $\BlueText{r_2(m)}$.
We will follow this convention when we consider correspondences
between homomorphisms of $\Omega_1$\Hyph algebra
and maps between sets
where we defined corresponding representations.
\qed
\end{remark}

\begin{theorem}
\labelTheorem{Tuple of maps is morphism of representations iff}
Let
\ShowEq{f:A->*B}g{A_1}{A_2}
be representation of $\Omega_1$\Hyph algebra $A_1$
in $\Omega_2$\Hyph algebra $A_2$ and
\ShowEq{f:A->*B}g{B_1}{B_2}
be representation of $\Omega_1$\Hyph algebra $B_1$
in $\Omega_2$\Hyph algebra $B_2$.
The map
\ShowEq{r12:A->B}rAB
is morphism of representations iff
\DrawEq{morphism of representations of universal algebra, 2m}{definition}
\end{theorem}
\ProofTheorem{\RefRepresentation}{Tuple of maps is morphism of representations iff}

\begin{remark}
\labelRemark{morphism of representations of universal algebra}
There are two ways
to interpret
\eqRef{morphism of representations of universal algebra, 2m}{definition}
\begin{itemize}
\item Let transformation $\BlueText{f(a)}$ map $m\in A_2$
into $\BlueText{f(a)}(m)$.
Then transformation
\ShowEq{g(r1(a))}
maps
\ShowEq{r2(m)in B2}
into
\ShowEq{r2(f(a,m))}
\item We represent morphism of representations from $f$ into $g$
using diagram
\DrawEq{morphism of representations of universal algebra, 2m 1}{definition}
From \EqRef{morphism of representations of universal algebra, definition, 2},
it follows that diagram $(1)$ is commutative.
\end{itemize}
We also use diagram
\ShowEq{morphism of representations of universal algebra, definition, 2m 2}
instead of diagram
\eqRef{morphism of representations of universal algebra, 2m 1}{definition}.
\qed
\end{remark}
\fi

\begin{definition}
If representation $f$ and $g$ coincide, then morphism of representations
\ShowEq{map r,R}{r_1}{r_2}{}
is called
\AddIndex{morphism of representation $f$}{morphism of representation f}.
\qed
\end{definition}

\ePrints{aa}
\ifx\Semafor\ValueOn

\begin{theorem}
\labelTheorem{unique morphism of representations of universal algebra}
Let the representation
\ShowEq{f:A->*B}g{A_1}{A_2}
of $\Omega_1$\Hyph algebra $A_1$ be single transitive representation
and the representation
\ShowEq{f:A->*B}g{B_1}{B_2}
of $\Omega_1$\Hyph algebra $B_1$ be single transitive representation.
Given homomorphism of $\Omega_1$\Hyph algebra
\ShowEq{f:A->B}{r_1}{A_1}{B_1}
consider a homomorphism of $\Omega_2$\Hyph algebra
\ShowEq{f:A->B}{r_2}{A_2}{B_2}
such that
\ShowEq{map r,R}{r_1}{r_2}{}
is morphism
of representations from $f$ into $g$.
The map $H$ is unique up to
choice of image $n=H(m)\in N$
of given element $m\in A_2$.
\end{theorem}
\ProofTheorem{\RefRepresentation}{unique morphism of representations of universal algebra}

\begin{definition}
\labelDefinition{transformation coordinated with equivalence}
Let us define equivalence $S$ on the set $A_2$.
Transformation $f$ is called
\AddIndex{coordinated with equivalence}{transformation coordinated with equivalence} $S$,
when
$f(m_1)\equiv f(m_2)(\mathrm{mod}\ S)$
follows from condition
\ShowEq{m1 m2 modS}.
\qed
\end{definition}

\begin{theorem}
\labelTheorem{transformation correlated with equivalence}
Consider equivalence $S$ on set $A_2$.
Consider
$\Omega_1$\Hyph algebra on set
\ShowEq{End O2}{A_2}.
If any transformation  $f\in{}^* A_2$ is coordinated with equivalence $S$,
then we can define the structure of $\Omega_1$\Hyph algebra
on the set ${}^*(A_2/S)$.
\end{theorem}
\ProofTheorem{\RefRepresentation}{transformation correlated with equivalence}

\begin{theorem}
\labelTheorem{decompositions of morphism of representations}
Let
\ShowEq{f:A->*B}g{A_1}{A_2}
be representation of $\Omega_1$\Hyph algebra $A_1$,
\ShowEq{f:A->*B}g{B_1}{B_2}
be representation of $\Omega_1$\Hyph algebra $B_1$.
Let
\ShowEq{r12:A->B}tAB
be morphism of representations from $f$ into $g$.
Then there exist decompositions of $t_1$ and $t_2$,
which we describe using diagram
\ShowEq{decompositions of morphism of representations, diagram}
\StartLabelItem
\begin{enumerate}
\item
The kernel of homomorphism
\ShowEq{kernel of homomorphism i}
is a congruence on $\Omega_i$\Hyph algebra $A_i$,
\EqParm{i=1,2}{=.}
\labelItem{kernel of homomorphism i}
\item
There exists decomposition of homomorphism $t_i$,
\EqParm{i=1,2}{=c}
\ShowEq{morphism of representations of algebra, homomorphism, 1}
\labelItem{exists decomposition of homomorphism}
\item
Maps
\ShowEq{morphism of representations of algebra, p1=}
\ShowEq{morphism of representations of algebra, p2=}
are natural homomorphisms.
\labelItem{p12 are natural homomorphisms}
\item
Maps
\ShowEq{morphism of representations of algebra, q1=}
\ShowEq{morphism of representations of algebra, q2=}
are isomorphisms.
\labelItem{q12 are isomorphisms}
\item
Maps
\ShowEq{morphism of representations of algebra, r1=}
\ShowEq{morphism of representations of algebra, r2=}
are monomorphisms.
\labelItem{r12 are monomorphisms}
\item $F$ is representation of $\Omega_1$\Hyph algebra $A_1/s$ in $A_2/s_2$
\item $G$ is representation of $\Omega_1$\Hyph algebra $t_1A_1$ in $t_2A_2$
\item The map
\ShowEq{map r,R}{p_1}{p_2}{}
is morphism of representations $f$ and $F$
\labelItem{(p12) is morphism of representations}
\item The map
\ShowEq{map r,R}{q_1}{q_2}{}
is isomorphism of representations $F$ and $G$
\labelItem{(q12) is isomorphism of representations}
\item The map
\ShowEq{map r,R}{r_1}{r_2}{}
is morphism of representations $G$ and $g$
\labelItem{(r12) is morphism of representations}
\item There exists decompositions of morphism of representations
\labelItem{exists decompositions of morphism of representations}
\ShowEq{decompositions of morphism of representations}
\end{enumerate}
\end{theorem}
\ProofTheorem{\RefRepresentation}{decompositions of morphism of representations}
\fi

\ShowDefinition{reduced morphism of representations}
\fi

\ePrints{5148-4632,MVector,5410-9916}
\ifx\Semafor\ValueOn
\input{\FilePrefix Preliminary.Omega.\TheLanguage}
\fi

\ePrints{MVector,5148-4632}
\ifx\Semafor\ValueOn

\Section{Basis of Representation of Universal Algebra}

\begin{definition}
Let
\ShowEq{f:A->*B}g{A_1}{A_2}
be representation of $\Omega_1$\Hyph algebra $A_1$
in $\Omega_2$\Hyph algebra $A_2$.
The set
\ShowEq{B2 subset A2}
is called
\AddIndex{stable set of representation}
{stable set of representation} $f$,
if
\ShowEq{fam in B2}
for each
\EqParm{a in A1}{=c}
$m\in B_2$.
\qed
\end{definition}

\begin{theorem}%
\labelTheorem{subrepresentation of representation}
Let
\ShowEq{f:A->*B}g{A_1}{A_2}
be representation of $\Omega_1$\Hyph algebra $A_1$
in $\Omega_2$\Hyph algebra $A_2$.
Let set
\ShowEq{B2 subset A2}
be subalgebra of $\Omega_2$\Hyph algebra $A_2$
and stable set of representation $f$.
Then there exists representation
\ShowEq{representation of algebra A in algebra B}
such that
\ShowEq{fB2(a)=}
Representation $f_{B_2}$ is called
\AddIndex{subrepresentation of representation}
{subrepresentation of representation} $f$.
\end{theorem}
\ProofTheorem{\RefRepresentation}{subrepresentation of representation}

\begin{theorem}%
\labelTheorem{lattice of subrepresentations}
The set\,\footnote{This definition is
similar to definition of the lattice of subalgebras
(\citeBib{Cohn: Universal Algebra}, p. 79, 80)}
\ShowEq{lattice of subrepresentations}
of all subrepresentations of representation $f$ generates
a closure system on $\Omega_2$\Hyph algebra $A_2$
and therefore is a complete lattice.
\end{theorem}
\ProofTheorem{\RefRepresentation}{lattice of subrepresentations}

We denote the corresponding closure operator by
\ShowEq{closure operator, representation}
Thus
\ShowEq{subrepresentation generated by set}
is the intersection of all subalgebras
of $\Omega_2$\Hyph algebra $A_2$
containing $X$ and stable with respect to representation $f$.

\begin{definition}
\labelDefinition{generating set of representation}
\ShowEq{show closure operator, representation}
is called
\AddIndex{subrepresentation generated by set $X$}
{subrepresentation generated by set},
and $X$ is a
\AddIndex{generating set of subrepresentation $J_f(X)$}
{generating set of subrepresentation}.
In particular, a
\AddIndex{generating set of representation $f$}
{generating set of representation}
is a subset $X\subset A_2$ such that
\ShowEq{generating set of representation}
\qed
\end{definition}

\ShowTheorem{structure of subrepresentations}
\ProofTheorem{\RefRepresentation}{structure of subrepresentations}

\begin{definition}
\labelDefinition{basis of representation}
If the set $X\subset A_2$ is generating set of representation
$f$, then any set $Y$, $X\subset Y\subset A_2$
also is generating set of representation $f$.
If there exists minimal set $X$ generating
the representation $f$, then the set $X$ is called
\AddIndex{basis of representation}{basis of representation} $f$.
\qed
\end{definition}

\begin{theorem}
\labelTheorem{X is basis of representation}
The generating set $X$ of representation $f$ is basis
iff for any $m\in X$
the set $X\setminus\{m\}$ is not
generating set of representation $f$.
\end{theorem}
\ProofTheorem{\RefRepresentation}{X is basis of representation}

\begin{remark}
\labelRemark{X is basis of representation}
The proof of the theorem
\RefTheorem{X is basis of representation}
gives us effective method for constructing the basis of the representation $f$.
Choosing an arbitrary generating set, step by step, we
remove from set those elements which have coordinates
relative to other elements of the set.
If the generating set of the representation is infinite,
then this construction may not have the last step.
If the representation has finite generating set,
then we need a finite number of steps to construct a basis of this representation.

As noted by Paul Cohn in
\citeBib{Cohn: Universal Algebra}, p. 82, 83,
the representation may have inequivalent bases.
For instance, the cyclic group of order six has
bases $\{a\}$ and $\{a^2,a^3\}$ which we cannot map
one into another by endomorphism of the representation.
\qed
\end{remark}

\begin{theorem}
\labelTheorem{automorphism uniquely defined by image of basis}
Let $X$ be the basis of the representation $f$.
Let
\ShowEq{R1 X->}
be arbitrary map of the set $X$.
There exists unique endomorphism of representation $f$\,\footnote{This statement is
similar to the theorem \citeBib{Serge Lang}-4.1, p. 135.}
\ShowEq{R:A2->A2}
defined by rule
\ShowEq{Rm=R1m}
\end{theorem}
\ProofTheorem{\RefRepresentation}{automorphism uniquely defined by image of basis}
\fi

\ePrints{MVector,5148-4632}
\ifx\Semafor\ValueOn

\Section{Tower of Representations of Universal Algebras}

\begin{definition}
\labelDefinition{tower of representations}
Consider set of $\Omega_k$\Hyph algebras $A_k$, \Kn1.
Assume $A=(A_1,...,A_n)$.
Assume $f=(f_{1,2},...,f_{n-1,n})$.
Set of representations $f_{k,k+1}$, \Kn1,
of $\Omega_k$\Hyph algebra $A_k$ in
$\Omega_{k+1}$\Hyph algebra $A_{k+1}$
is called
\AddIndex{tower $(A,f)$ of representations of $\Omega$\Hyph algebras}
{tower of representations of algebras}.
\qed
\end{definition}

Consider following diagram for the purposes of illustration of
definition \RefDefinition{tower of representations}
\ShowEq{Tower of Representations}
$f_{i,i+1}$ is representation of $\Omega_i$\Hyph algebra $A_i$
in $\Omega_{i+1}$\Hyph algebra $A_{i+1}$.
$f_{i+1,i+2}$ is representation of $\Omega_{i+1}$\Hyph algebra $A_{i+1}$
in $\Omega_{i+2}$\Hyph algebra $A_{i+2}$.

\begin{definition}
Consider the set of $\Omega_k$\Hyph algebr $A_k$, $B_k$, \Kn1.
The set of maps
$(h_1,...,h_n)$ is called \AddIndex{morphism from tower of representations
$(A,f)$ into tower of representations $(B,g)$}
{morphism from tower of representations into tower of representations},
if for any $k$, $k=1$, ..., $n-1$,
the tuple of maps $(h_k,h_{k+1})$ is
morphism of representations from $f_{k,k+1}$ into $g_{k,k+1}$.
\qed
\end{definition}

For any $k$, $k=1$, ..., $n-1$, we have diagram
\ShowEq{morphism of tower of representations of F algebra, level k, diagram}
Equations
\ShowEq{morphism of tower of representations, level k}
\ShowEq{morphism of tower of representations, levels k k+1}
express commutativity of diagram (1).

\begin{definition}
{\it
Let
\ShowEq{Aij}
be a set of
\ShowEq{Omegaij}%
algebras.
}
\qed
\end{definition}
\fi

\ePrints{MVector,5148-4632,4975-6381,1601.03259,1610.309618526,CACAA.06}
\ifx\Semafor\ValueOn
\Section{Permutation}

\begin{definition}
{\it
An injective map of finite set into itself is called
\AddIndex{permutation}{permutation}.\,\footnote{
You can see definition and properties of permutation in
\citeBib{Kurosh: High Algebra}, pages 27 - 32,
\citeBib{Cohn: Algebra 1}, pages 58, 59.
}
}
\qed
\end{definition}

Usually we write a permutation $\sigma$ as a matrix
\ShowEq{permutation as matrix}
The notation
\EqRef{permutation as matrix}
is equivalent to the statement
\ShowEq{a->sigma a}
So the order of columns in the notation
\EqRef{permutation as matrix}
is not essential.

Since there is an order
on the set
\EqParm{set a1n}{=z}
(for instance, we assume, that $a_i$ precedes $a_j$
when $i<j$),
then we may assume that elements of first row
are written according to the intended order
and we will identify permutation with second row
\ShowEq{permutation as matrix 2}

\begin{definition}
{\it
The map
\ShowEq{sigma->+-}
defined by the equality
\def\permutation{permutation }
\def\even{even}
\def\odd{odd}
\ShowEq{parity of permutation}
is called
\AddIndex{parity of permutation}{parity of permutation}.
}
\qed
\end{definition}
\fi

\ePrints{4975-6381,1506.00061,MQuater,1601.03259,MVector,1610.309618526,CACAA.06}
\ifx\Semafor\ValueOn
\Section{Module over Ring}

\ShowEq{=module}

\ShowDefinition{module over commutative ring}

\ShowTheorem{definition of module over commutative ring}
\ProofTheorem{\RefLinearMap}{definition of module over commutative ring}

\ShowDefinition{coordinates of vector 2016}

\ePrints{MVector}
\ifx\Semafor\ValueOn

\ShowEq{=left}
\ShowTheorem{set of vectors generated by set of vectors}
\ProofTheorem{\RefLinearMap}{set of vectors generated by set of vectors}

\ShowDefinition{vector linearly dependent on vectors}

\ShowEq{remark: scalars and vectors as matrix}

\ShowDefinition{linearly independent vectors}

\ShowTheorem{basis of module}
\ProofTheorem{\RefLinearMap}{basis of module}

\ShowDefinition{free module over ring}

\ShowTheorem{division algebra, basis}
\ProofTheorem{\RefLinearMap}{division algebra, basis}
\fi

\ShowDefinition{linear map from A1 to A2, commutative module}

\ShowTheorem{linear map from A1 to A2, commutative module}
\ProofTheorem{\RefLinearMap}{linear map from A1 to A2, commutative module}
\fi

\ePrints{4975-6381,MQuater,1601.03259}
\ifx\Semafor\ValueOn
\begin{theorem}
\labelTheorem{linear map, 0, D algebra}
Let map
\ShowEq{f:A->B}f{A_1}{A_2}
be linear map of $D$\Hyph module $A_1$ into $D$\Hyph module $A_2$.
Then
\ShowEq{linear map, 0, D algebra}
\end{theorem}
\ProofTheorem{\RefLinearMap}{linear map, 0, D algebra}
\fi

\ePrints{1601.03259,4975-6381}
\ifx\Semafor\ValueOn

\begin{definition}
\labelDefinition{polylinear map of modules}
Let $D$ be the commutative ring.
Let $A_1$, ..., $A_n$, $S$ be $D$\Hyph modules.
We call map
\ShowEq{polylinear map of algebras}
\AddIndex{polylinear map}{polylinear map} of modules
$A_1$, ..., $A_n$
into module $S$,
if
\ShowEq{polylinear map of algebras, 1}
\qed
\end{definition}

\ShowTheorem{sum of polylinear maps, module}
\ProofTheorem{\RefLinearMap}{sum of polylinear maps, module}

\ShowEq{corollary: sum of linear maps, D module}
\fi

\ePrints{1601.03259,4975-6381}
\ifx\Semafor\ValueOn
\begin{definition}
\labelDefinition{tensor product of algebras}
Let $A_1$, ..., $A_n$ be
free modules over commutative ring $D$.\,\footnote{
I give definition
of tensor product of $D$\Hyph modules
following to definition in \citeBib{Serge Lang}, p. 601 - 603.
}
Consider category $\mathcal A_1$ whose objects are
polylinear maps
\ShowEq{polylinear maps category}
where $S_1$, $S_2$ are modules over ring $D$,
We define morphism
\ShowEq{f->g}
to be linear map
\ShowEq{f->g 1}
making following diagram commutative
\ShowEq{polylinear maps category, diagram}
Universal object
\ShowEq{tensor product of algebras}
of category $\mathcal A_1$ is called
\AddIndex{tensor product}{tensor product}
of modules $A_1$, ..., $A_n$.
\qed
\end{definition}

\begin{theorem}
\labelTheorem{Tensor product is distributive over sum}
Let $D$ be the commutative ring.
Let $A_1$, ..., $A_n$ be $D$\Hyph modules.
Tensor product is distributive over sum
\ShowEq{tensors 1, tensor product}
The representation of the ring $D$
in tensor product is defined by equation
\ShowEq{tensors 2, tensor product}
\end{theorem}
\ProofTheorem{\RefLinearMap}{Tensor product is distributive over sum}

\begin{theorem}
\labelTheorem{tensor product and polylinear map}
Let $A_1$, ..., $A_n$ be
modules over commutative ring $D$.
Tensor product
\ShowEq{f:xA->oxA}
exists and unique.
We use notation
\ShowEq{fxa=oxa}
for image of the map $f$.
Let
\ShowEq{map g, algebra, tensor product}
be polylinear map into $D$\Hyph module $V$.
There exists a linear map
\ShowEq{map h, algebra, tensor product}
such that the diagram
\ShowEq{map gh, algebra, tensor product}
is commutative.
The map \(h\) is defined by the equation
\ShowEq{g=h, algebra, tensor product}
\end{theorem}
\begin{proof}
See the proof of theorems
\ShowEq{ref tensor product and polylinear map}
\end{proof}

\begin{convention}
\labelConvention{isomorphic representations S1=S2}
Algebras $S_1$, $S_2$ may be different sets.
However they are indistinguishable for us when we consider them
as isomorphic representations.
In such case, we write the statement $S_1=S_2$.
\qed
\end{convention}

\begin{theorem}
\labelTheorem{tensor product is associative}
\ShowEq{A1xA2xA3}
\end{theorem}
\ProofTheorem{\RefLinearMap}{tensor product is associative}

\begin{definition}
\labelDefinition{tensor power of algebra}
Tensor product
\ShowEq{tensor power of algebra}
is called
\AddIndex{tensor power}{tensor power} of module $A_1$.
\qed
\end{definition}

\begin{theorem}
\labelTheorem{V times->V otimes}
The map
\ShowEq{V times->V otimes}
is polylinear map.
\end{theorem}
\ProofTheorem{\RefLinearMap}{V times->V otimes}
\fi

\ePrints{4975-6381,MQuater,1601.03259}
\ifx\Semafor\ValueOn
\begin{theorem}
\labelTheorem{standard component of tensor, algebra}
Let
\ShowEq{tensor product of algebras, basis i}
be the basis of module $A_i$ over ring $D$.
We can represent any tensor $a\in\Tensor A$ in the following form
\ShowEq{standard component of tensor}
\ShowEq{tensor canonical representation, algebra}
Expression
$\ShowSymbol{standard component of tensor}{}$
is called
\AddIndex{standard component of tensor}{standard component of tensora}.
\end{theorem}
\ProofTheorem{\RefLinearMap}{standard component of tensor, algebra}
\fi

\ePrints{4975-6381,1506.00061,MQuater,1601.03259,MVector,1610.309618526,CACAA.06}
\ifx\Semafor\ValueOn
\Section{Algebra over Commutative Ring}

\ShowDefinition{algebra over ring}

\ePrints{1601.03259}
\ifx\Semafor\ValueOff
\begin{theorem}
\labelTheorem{multiplication in algebra is distributive over addition}
The multiplication in the algebra $A_1$ is distributive over addition.
\end{theorem}
\ProofTheorem{\RefLinearMap}{multiplication in algebra is distributive over addition}
\fi

\ePrints{1506.00061,MQuater,4975-6381,1601.03259,1610.309618526,CACAA.06}
\ifx\Semafor\ValueOn
\begin{convention}
Element of
$D$\Hyph algebra
$A$ is called
\AddIndex{$A$\Hyph number}{A number}.
For instance, complex number is also called
$C$\Hyph number, and quaternion is called $H$\Hyph number.
\qed
\end{convention}
\fi

The multiplication in algebra can be neither commutative
nor associative. Following definitions are based
on definitions given in \citeBib{Richard D. Schafer}, page 13.

\ShowDefinition{commutator of algebra}

\ShowDefinition{associator of algebra}

\ShowDefinition{nucleus of algebra}

\ShowDefinition{center of algebra}

\ShowEq{convention: we use separate color for index of element}

\ShowEq{convention: unit of algebra in basis}

\ePrints{1610.309618526,1601.03259,CACAA.06}
\ifx\Semafor\ValueOff
\begin{theorem}
\labelTheorem{product in algebra}
Let $\Basis e$ be the basis of free algebra $A_1$ over ring $D$.
Let
\ShowEq{a b in basis of algebra}
We can get the product of $a$, $b$ according to rule
\ShowEq{product in algebra}
where
\ShowEq{structural constants of algebra}
are \AddIndex{structural constants}{structural constants}
of algebra $A_1$ over ring $D$.
The product of basis vectors in the algebra $A_1$ is defined according to rule
\ShowEq{product of basis vectors, algebra}
\end{theorem}
\ProofTheorem{\RefLinearMap}{product in algebra}
\fi

\ePrints{4975-6381,1601.03259,1610.309618526,CACAA.06}
\ifx\Semafor\ValueOn
\begin{theorem}
\labelTheorem{Free Algebra over Ring}
\ShowEq{text: Free Algebra over Ring}
\end{theorem}
\begin{proof}
\TheoremFollows
\RefTheorem[\RefLinearMap]{Free Algebra over Ring}
and the corollary
\RefCorollary[\RefLinearMap]{Free Algebra over Ring}.
\end{proof}
\fi

\ShowEq{=DD}

\ShowDefinition{linear map from A1 to A2, algebra}

\ePrints{1601.03259,4975-6381}
\ifx\Semafor\ValueOn

\begin{theorem}
\labelTheorem{module L(A;A) is algebra}
Let $A$ be algebra over commutative ring $D$.
$D$\Hyph module $\mathcal L(D;A;A)$ equiped by product
\ShowEq{module L(A;A) is algebra}
\ShowEq{product of linear map, algebra 1}
is algebra over $D$.
\end{theorem}
\ProofTheorem{\RefLinearMap}{module L(A;A) is algebra}
\fi

\begin{definition}
\labelDefinition{polylinear map of algebras}
Let $A_1$, ..., $A_n$, $S$ be $D$\Hyph algebras.
Polylinear map
\ShowEq{polylinear map of algebras}
of $D$\Hyph modules
$A_1$, ..., $A_n$
into $D$\Hyph module $S$
is called
\AddIndex{polylinear map}{polylinear map}
of $D$\Hyph algebras
$A_1$, ..., $A_n$
into $D$\Hyph algebra $S$.
Let us denote
\ShowEq{set polylinear maps}
set of polylinear maps
of $D$\Hyph algebras
$A_1$, ..., $A_n$
into $D$\Hyph algebra
$S$.
Let us denote
\ShowEq{set polylinear maps An}
set of $n$\hyph linear maps
of $D$\Hyph algebra $A_1$ ($A_1=...=A_n=A_1$)
into $D$\Hyph algebra
$S$.
\qed
\end{definition}

\ShowTheorem{tensor product of D-algebras is D-algebra}
\ProofTheorem{\RefLinearMap}{tensor product of D-algebras is D-algebra}

\ShowTheorem{representation of algebra A2 in LA}
\ProofTheorem{\RefLinearMap}{representation of algebra A2 in LA}

\begin{convention}
I assume sum over index $i$
in expression like
\ShowEq{Sum over repeated index}
\qed
\end{convention}


\ePrints{4975-6381,1506.00061,MQuater,1601.03259}
\ifx\Semafor\ValueOn
\ShowTheorem{standard representation of map A1 A2, associative algebra}
\ProofTheorem{\RefLinearMap}{standard representation of map A1 A2, associative algebra}

\ePrints{4975-6381}
\ifx\Semafor\ValueOff
\ShowDefinition{component of linear map}
\fi

\ShowTheorem{conjugation transformation}
\ProofTheorem{\RefLinearMap}{conjugation transformation}

\ShowTheorem{representation of composition of linear maps}
\ProofTheorem{\RefLinearMap}{representation of composition of linear maps}

\ShowTheorem{representation of composition of linear maps A->A}
\ProofTheorem{\RefLinearMap}{representation of composition of linear maps A->A}

\fi

\ePrints{1601.03259,4975-6381,1610.309618526,CACAA.06}
\ifx\Semafor\ValueOn
\ShowTheorem{coordinates of map A1 A2, algebra}
\ProofTheorem{\RefLinearMap}{coordinates of map A1 A2, algebra}

According to the definition
\ShowEq{Basis F=...}
We also consider the set
\ShowEq{F=...}
having in mind that the equality
\ShowEq{Fi=Fj}
is possible when $i\ne j$.
\fi

\ePrints{4975-6381,1506.00061}
\ifx\Semafor\ValueOn
\begin{convention}
In the equation
\ShowEq{n linear map A LA}
as well as in other expressions of polylinear map,
we have convention that map $f_i$ has variable $x_i$ as argument.
\qed
\end{convention}
\fi

\ePrints{1601.03259,4975-6381,1506.00061,1610.309618526,CACAA.06}
\ifx\Semafor\ValueOn

\ShowTheorem{L(An;B) is free D module}
\ProofTheorem{\RefLinearMap}{L(An;B) is free D module}

\begin{theorem}
\labelTheorem{polylinear map, algebra} 
Let $A$ be associative $D$\Hyph algebra.
Polylinear map
\ShowEq{polylinear map, algebra}
generated by maps
\ShowEq{I1n in L(A;A)}
has form
\ShowEq{polylinear map, algebra, canonical morphism}
where $\sigma_s$ is a transposition of set of variables
\ShowEq{transposition of set of variables, algebra}
\end{theorem}
\ProofTheorem{\RefLinearMap}{polylinear map, algebra}

\ShowTheorem{representation of algebra An in LAnA}
\ProofTheorem{\RefLinearMap}{representation of algebra An in LAnA}

\begin{convention}
Since the tensor
\EqParm{a in Aoxn+1}{=z}
has the expansion
\ShowEq{a=oxi}
then set of permutations
\ShowEq{si in Sn}
and tensor $a$
generate the map
\ShowEq{ox:An->A}
defined by rule
\ShowEq{ox circ =}
\qed
\end{convention}
\fi
\fi

\ePrints{}
\ifx\Semafor\ValueOn
\Section{Module over Algebra}

\ShowEq{=left}

\ShowDefinition{module over algebra}

\ShowTheorem{definition of A module}
\ProofTheorem{\RefLinearMap}{definition of left A module}

\ShowEq{=right}
\ShowDefinition{module over algebra}
\fi

\ePrints{1506.00061}
\ifx\Semafor\ValueOn
\Section{Algebra with Conjugation}

Let $D$ be commutative ring.
Let $A$ be $D$\hyph algebra with unit $e$, $A\ne D$.

Let there exist subalgebra $F$ of algebra $A$ such that
$F\ne A$, $D\subseteq F\subseteq Z(A)$, and algebra $A$ is a free
module over the ring $F$.
Let $\Basis e$ be the basis of free module $A$ over ring $F$.
We assume that $e_{\gi 0}=1$.

\begin{theorem}
\labelTheorem{structural constant of algebra with unit}
Structural constants of
$D$\hyph algebra with unit $e$ satisfy condition
\ShowEq{structural constant re algebra, 1}
\end{theorem}
\begin{proof}
\TheoremFollows
\ShowEq{ref structural constant of algebra with unit}
\end{proof}

Consider maps
\ShowEq{Re Im A->A}
defined by equation
\ShowEq{Re Im A->F 0}
\ShowEq{Re Im A->F 1}
The expression
\ShowEq{Re Im A->F 2}
is called
\AddIndex{scalar of element}{scalar of algebra} $d$.
The expression
\ShowEq{Re Im A->F 3}
is called
\AddIndex{vector of element}{vector of algebra} $d$.

According to
\EqRef{Re Im A->F 1}
\ShowEq{Re Im A->F 4}
We will use notation
\ShowEq{scalar algebra of algebra}
to denote
\AddIndex{scalar algebra of algebra}{scalar algebra of algebra} $A$.

\begin{theorem}
\labelTheorem{Re Im A->F}
The set
\ShowEq{vector module of algebra}
is $(\re A)$\Hyph module
which is called
\AddIndex{vector module of algebra}{vector module of algebra} $A$.
\ShowEq{A=re+im}
\end{theorem}
\begin{proof}
\TheoremFollows
\ShowEq{ref Re Im A->F}
\end{proof}

According to the theorem \RefTheorem{Re Im A->F},
there is unique defined representation
\ShowEq{d Re Im}

\begin{definition}
\labelDefinition{conjugation in algebra}
The map
\ShowEq{conjugation in algebra}
\ShowEq{conjugation in algebra, 0}
is called
\AddIndex{conjugation in algebra}{conjugation in algebra}
provided that this map satisfies
\ShowEq{conjugation in algebra, 1}
$(\re A)$\Hyph algebra $A$
equipped with conjugation is called
\AddIndex{algebra with conjugation}{algebra with conjugation}.
\qed
\end{definition}

\begin{theorem}
\labelTheorem{conjugation in algebra}
The $(\re A)$\Hyph algebra $A$ is
algebra with conjugation
iff structural constants of
$(\re A)$\Hyph algebra $A$ satisfy condition
\ShowEq{conjugation in algebra, 5}
\ShowEq{conjugation in algebra, 4}
\end{theorem}
\begin{proof}
\TheoremFollows
\ShowEq{ref conjugation in algebra}
\end{proof}
\fi

\ePrints{1601.03259,4975-6381,1506.00061}
\ifx\Semafor\ValueOn
\Section{Polynomial over Associative \texorpdfstring{$D$}{D}-Algebra}

Let $D$ be commutative ring and
$A$ be associative $D$\Hyph algebra with unit.
Let $\Basis F$ be basis of algebra
\ShowEq{L(A;B)}DAA.

\begin{theorem}
\labelTheorem{monomial of power k}
Let $p_k(x)$ be
\AddIndex{monomial of power}{monomial of power} $k$
over $D$\Hyph algebra $A$.
Then
\StartLabelItem
\begin{enumerate}
\item
Monomial of power $0$ has form
\ShowEq{p0(x)=a0}
\item
If $k>0$, then
\ShowEq{monomial of power k}
where $a_k\in A$ and $F\in\Basis F$.
\end{enumerate}
\end{theorem}
\begin{proof}
\TheoremFollows
\ShowEq{ref monomial of power k}
\end{proof}

In particular, monomial of power $1$ has form
\ShowEq{p1(x)=aFxa}

\begin{definition}
\labelDefinition{Abelian group of homogeneous polynomials over algebra}
We denote
\ShowEq{module of homogeneous polynomials over algebra}
Abelian group
generated by the set of monomials of power $k$.
Element $p_k(x)$ of Abelian group $A_k[x]$ is called
\AddIndex{homogeneous polynomial of power}
{homogeneous polynomial of power} $k$.
\qed
\end{definition}

\begin{convention}
Let the tensor
\EqParm{a in Aoxn+1}{=.}
Let
\ShowEq{F1n in F}
When
\ShowEq{x1=xn=x}
we assume
\ShowEq{a xn=}
\qed
\end{convention}

\begin{convention}
If we have few tuples of maps $F\in\Basis F$,
then we will use index like $[k]$
to index tuple
\ShowEq{F[k]=...}
\qed
\end{convention}

\begin{theorem}
\labelTheorem{homogeneous polynomial a circ xk}
We can present homogeneous polynomial $p(x)$
in the following form
\ShowEq{p(x)=a circ xk}
\end{theorem}
\begin{proof}
\TheoremFollows
\ShowEq{ref map A(k+1)->pk otimes is polylinear map}
\end{proof}

\begin{definition}
\labelDefinition{Polynomial over algebra}
We denote
\ShowEq{algebra of polynomials over algebra}
direct sum\,\footnote{See the definition
of direct sum of modules in
\citeBib{Serge Lang},
page 128.
On the same page, Lang proves the existence of direct sum of modules.
}
of $A$\Hyph modules $A_n[x]$.
An element $p(x)$ of $A$\Hyph module $A[x]$ is called
\AddIndex{polynomial}{polynomial}
over $D$\Hyph algebra $A$.
\qed
\end{definition}

\ePrints{4975-6381,1601.03259}
\ifx\Semafor\ValueOn
The following definition
\RefDefinition{otimes -}
is based on the definition
\Ref[1302.7204]{definition: otimes -}.

\ShowDefinition{otimes -}

\begin{definition}
\labelDefinition{product of homogeneous polynomials}
Product of homogeneous polynomials
\ShowEq{polynomials p,r}
is defined by the equality\,\footnote{
The definition
\RefDefinition{product of homogeneous polynomials}
is based on the definition
\Ref[1302.7204]{definition: product of homogeneous polynomials}.
}
\ShowEq{pr=p circ r}
\qed
\end{definition}

\ePrints{}
\ifx\Semafor\ValueOn
\begin{theorem}
\labelTheorem{p(x)=a k-1 k}
Let
\ShowEq{p(x)=a circ xk}
be monomial of power $k>1$.
Then the polynomial $p(x)$ can be represented using
one of the following forms
\ShowEq{p(x)=a k-1 k 1}
\ShowEq{p(x)=a k-1 k 2}
where
\ShowEq{p(x)=a k-1 k 3}
\end{theorem}
\begin{proof}
See the proof of the theorem
\RefTheorem[1302.7204]{p(x)=a k-1 k}.
\end{proof}
\fi
\fi

\ePrints{1506.00061}
\ifx\Semafor\ValueOn

\begin{definition}
\labelDefinition{divisor of polynomial}
The polynomial $p(x)$ is called
\AddIndex{divisor of polynomial}{divisor of polynomial} $r(x)$,
if we can represent the polynomial $r(x)$ as
\ShowEq{r=qpq(x)}
\qed
\end{definition}

\begin{theorem}
\labelTheorem{r=+q circ p}
Let
\ShowEq{p=p circ x}
be homogeneous polynomial of power $1$ and $p_1$ be nonsingular tensor.
Let
\ShowEq{r power k}
be polynomial of power $k$.
Then
\ShowEq{r=+q circ p}
\end{theorem}
\begin{proof}
\TheoremFollows
\ShowEq{ref r=+q circ p}
\end{proof}

\begin{theorem}
\labelTheorem{r=+q circ p 1}
Let
\ShowEq{p=p+p circ x}
be polynomial of power $1$ and $p_1$ be nonsingular tensor.
Let
\ShowEq{r power k}
be polynomial of power $k$.
Then
\ShowEq{r=+q circ p 1}
\end{theorem}
\begin{proof}
\TheoremFollows
\ShowEq{ref r=+q circ p 1}
\end{proof}
\fi
\fi

\ePrints{MQuater,1610.309618526,CACAA.06}
\ifx\Semafor\ValueOn
\Section{Derivative of Map of Banach Algebra}

\begin{definition}
Normed $D$\Hyph algebra $A_1$ is called
\AddIndex{Banach $D$\Hyph algebra}{Banach algebra}
if any fundamental sequence of elements
of algebra $A_1$ converges, i.e.
has limit in algebra $A_1$.
\qed
\end{definition}

\begin{definition}
\labelDefinition{continuous map, algebra}
A map
\ShowEq{f:A->B}f{A_1}{A_2}
of Banach $D_1$\Hyph algebra $A_1$ with norm $|x|_1$
into Banach $D_2$\Hyph algebra $A_2$ with norm $|y|_2$
is called \AddIndex{continuous}{continuous map}, if
for every as small as we please $\epsilon>0$
there exist such $\delta>0$, that
\[
|x'-x|_1<\delta
\]
implies
\[
|f(x')-f(x)|_2<\epsilon
\]
\qed
\end{definition}

\ePrints{1610.309618526,CACAA.06}
\ifx\Semafor\ValueOff

\ShowTheorem{set of A->B is D module}
\begin{proof}
\TheoremFollows
ффф
\end{proof}

\else

\ShowDefinition{norm of polylinear map}

\ShowTheorem{|f(a)|<|f||a| 1n}
\ePrints{CACAA.06}
\ifx\Semafor\ValueOff
\ProofTheorem{\RefCalculus}{|f(a)|<|f||a| 1n}
\else
\ShowProof{|f(a)|<|f||a| 1n}
\fi

\ShowTheorem{|on|->0 ona1p->0}
\ePrints{CACAA.06}
\ifx\Semafor\ValueOff
\ProofTheorem{\RefCalculus}{|on|->0 ona1p->0}
\else
\ShowProof{|on|->0 ona1p->0}
\fi
\fi

\ShowDefinition{differentiable map}DAB{algebra}

\begin{definition}
\labelDefinition{differential of map}
{\it
Since, for given $x$, we consider the increment
\eqRef{derivative of map, def}{algebra}
of the map
\ShowEq{f:A->B}fAA
as function of differential $dx$ of variable $x$,
then the linear part of this function
\ShowEq{differential of independent variable}
\ShowEq{differential of map}
\DrawEq{differential of map =}{}
is called
\AddIndex{differential of map}{differential of map}
$f$.
}
\qed
\end{definition}

\ShowEq{remark: differential L(A,A)}

\ePrints{1610.309618526,CACAA.06}
\ifx\Semafor\ValueOff
\ShowTheorem{derivative, representation in algebra}
\begin{proof}
See the proof of the theorem
\ShowEq{ref derivative, representation in algebra}
\end{proof}
\fi

\ShowTheorem{representation of derivative, algebra A->B}
\begin{proof}
\TheoremFollows
\ShowEq{ref differentiable map A->B}
\end{proof}

\ShowDefinition{coordinates of derivative, algebra A->B}

\ePrints{1610.309618526,CACAA.06}
\ifx\Semafor\ValueOn

\ShowTheorem{dfa1p/dx=df/dx a1p}
\ShowProof{dfa1p/dx=df/dx a1p}

\ShowTheorem{bilinear map and differential}

\ShowTheorem{derivative of tensor product}

\fi

\ePrints{1610.309618526,CACAA.06}
\ifx\Semafor\ValueOn
\ShowDefinition{derivative of Second Order, algebra}

\ShowDefinition{derivative of Order n, algebra}
\fi
\fi

\ePrints{1601.03259,4975-6381,MQuater,1610.309618526,CACAA.06}
\ifx\Semafor\ValueOn

\Section{Complex Field}

\begin{theorem}
\labelTheorem{complex field over real field}
Let us consider complex field $C$ as two-dimensional algebra over real field.
Let
\ShowEq{basis of complex field}
be the basis of algebra $C$.
Then in this basis product has form
\ShowEq{product of complex field}
and structural constants have form
\ShowEq{structural constants of complex field}
\end{theorem}
\ProofTheorem{1003.1544}{complex field over real field}

\begin{definition}
{\it
Complex field has following
\AddIndex{maps of conjugation}{map of conjugation}
\ShowEq{Ex= C}
\ShowEq{Ix= C}
}
\qed
\end{definition}

\begin{theorem}
\labelTheorem{linear map of complex field}
A linear map
\ShowEq{f:C->C y=}
of complex field has form
\ShowEq{linear map of complex field, structure}
where $C$\Hyph numbers
\ShowEq{a01}
are defined by the equality
\ShowEq{a01=}
\end{theorem}
\ProofTheorem{1003.1544}{additive map of complex field, structure}

\begin{corollary}
\labelCorollary{basis of L(R;C;C)}
{\it
$C\otimes C$\Hyph module
\ShowEq{L(R;C;C)}
is $C$\Hyph vector space
and has the basis
\EqParm{e=(E,I)}{=.}
}
\qed
\end{corollary}

\begin{theorem}
\labelTheorem{CE is algebra isomorphic to complex field}
The set
\ShowEq{CE set}
is $R$\Hyph algebra isomorphic to complex field.
\end{theorem}
\begin{proof}
The theorem follows from equalities
\ShowEq{aE+bE=(a+b)E}
\ShowEq{aE o bE=(ab)E}
based on the theorem
\RefTheorem{linear map of complex field}.
\end{proof}

\begin{theorem}
\labelTheorem{L(CCC)=CE}
\ShowEq{L(CCC)=CE}
\end{theorem}
\begin{proof}
The theorem follows from the equality
\EqRef{linear map of complex field, structure, 2}
and from commutativity of product of complex numbers.
\end{proof}

\begin{theorem}[the Cauchy\Hyph Riemann equations]
\labelTheorem{Cauchy Riemann equations linear}
Matrix of linear map $f\in CE$
\ShowEq{yi=xj fji}
satisfies relationship
\ShowEq{complex field over real field}
\end{theorem}
\ProofTheorem{1003.1544}{complex field over real field}

\begin{theorem}
\labelTheorem{matrix of linear map f in CI}
Matrix of linear map $f\in CI$
\ShowEq{yi=xj fji}
satisfies relationship\,\footnote{
See also section
\Ref[1003.1544]{section: System of Additive Equations in Complex Field}.
}
\ShowEq{matrix of linear map f in CI}
\end{theorem}
\begin{proof}
The statement follows from equations
\ShowEq{congugate map complex field, 3}
\end{proof}

\ePrints{4975-6381,MQuater}
\ifx\Semafor\ValueOn

\begin{theorem}
The set
\ShowEq{CI set}
is not $R$\Hyph algebra.
\end{theorem}
\begin{proof}
The theorem follows from the equality
\ShowEq{aI o bI=...E}
based on the theorem
\RefTheorem{linear map of complex field}.
\end{proof}

\ePrints{MQuater}
\ifx\Semafor\ValueOn

\begin{definition}
\labelDefinition{differentiable map of complex field}
{\it
The map
\ShowEq{f:C->C}
of complex field
is called
\AddIndex{differentiable}{differentiable map}
on the set $U\subset C$,
if at every point $x\in U$
the increment of map
$f$ can be represented as
\ShowEq{derivative of map}
\DrawEq{derivative of map, def}{}
where
\ShowEq{df:A->B}CC
is linear map of $C$\Hyph vector space $C$ and
\ShowEq{o:A->B}oCC
is such continuous map that
\ShowEq{lim |o|/|a|}CC
Linear map
$\displaystyle\ShowSymbol{derivative of map}{}$
is called
\AddIndex{derivative of map}{derivative of map}
$f$.
}
\qed
\end{definition}

\begin{theorem}
\labelTheorem{dx fx in CE}
Since the map
\ShowEq{f:C->C}
of complex field
is differentiable,
then
\ShowEq{dx fx in CE}
\end{theorem}
\begin{proof}
\TheoremFollows
\RefTheorem{CE is algebra isomorphic to complex field}
and the definition
\RefDefinition{differentiable map of complex field}.
\end{proof}

\begin{definition}
{\it
Differentiable map of complex field is called
\AddIndex{holomorphic}{holomorphic map}.
}
\qed
\end{definition}

\begin{theorem}[the Cauchy\Hyph Riemann equations]
Since matrix
\ShowEq{Jacobian of map of complex variable}
is Jacobian matrix of map of complex variable
\ShowEq{map of complex variable}
over real field,
then
\ShowEq{Cauchy Riemann equations, complex field, 1}
\end{theorem}
\begin{proof}
The theorem follows from theorems
\RefTheorem{Cauchy Riemann equations linear},
\RefTheorem{dx fx in CE}.
\end{proof}

\begin{theorem}[the Cauchy\Hyph Riemann equation]
\labelTheorem{Cauchy Riemann equations, complex field}
Derivative of map of complex variable satisfies
to the Cauchy\Hyph Riemann equation
\DrawEq{Cauchy Riemann equation, complex field}{theorem}
\end{theorem}
\begin{proof}
The equality
\ShowEq{Cauchy Riemann equations, complex field, 2, 1}
follows from equalities
\EqRef{Cauchy Riemann equations, complex field, 1}.
The equality
\eqRef{Cauchy Riemann equation, complex field}{theorem}
follows from the equality
\EqRef{Cauchy Riemann equations, complex field, 2, 1}.
\end{proof}

\begin{definition}
\labelDefinition{projection maps, complex field}
{\it
Following projection maps are defined in complex field
\ShowEq{projection maps, complex field}
}
\qed
\end{definition}

\begin{theorem}
\labelTheorem{Expansion of projection maps relative E,I}
Expansion of projection maps relative basis
\EqParm{e=(E,I)}{=z}
has form
\ShowEq{projection mappings 1, complex field}
\end{theorem}
\begin{proof}
The theorem follows from the theorem
\RefTheorem{linear map of complex field}
and from the definition
\RefDefinition{projection maps, complex field}.
\end{proof}

\begin{theorem}
\labelTheorem{df=...dx= CE}
The derivative of map
\ShowEq{f:C->C}
of complex field
has form\,\footnote{
Theorems in this section
are similar to theorems in \citeBib{Shabat: Complex Analysis}, p. 15 - 19.}
\ShowEq{df=...dx= CE}
\end{theorem}
\begin{proof}
According to the theorem
\RefTheorem{Expansion of projection maps relative E,I},
\ShowEq{dx01=}
follows from the equality
\ShowEq{dx=dx0+dx1}
The equality
\ShowEq{df=...dx}
follows from equalities
\EqRef{dx01=}
and
\ShowEq{df=df/dx0+df/dx1}
The equality
\DrawEq{df=...dx=}{theorem}
follows from the equality
\EqRef{df=...dx}
and from the definition
\RefDefinition{differential of map}.
The equality
\EqRef{df=...dx= CE}
follows from the equality
\eqRef{df=...dx=}{theorem}
and from the theorem
\RefTheorem{Cauchy Riemann equations, complex field}.
\end{proof}
\fi
\fi
\fi

\ePrints{1506.00061,MQuater,1610.309618526,CACAA.06}
\ifx\Semafor\ValueOn

\Section{Quaternion Algebra}

\begin{definition}
Let $R$ be real field.
Extension field
\ShowEq{quaternion algebra}
is called
\AddIndex{the quaternion algebra}{quaternion algebra}
if multiplication in algebra $H$ is defined according to rule
\ShowEq{product of quaternions}
\qed
\end{definition}

Elements of the algebra $H$ have form
\ShowEq{quaternion algebra, element}
Quaternion
\ShowEq{conjugate to the quaternion}
is called conjugate to the quaternion $x$.
We define
\AddIndex{the norm of the quaternion}{norm of quaternion}
$x$ using equation
\ShowEq{norm of quaternion}
From the equality
\EqRef{norm of quaternion},
it follows that inverse element has form
\ShowEq{inverce quaternion}

\ePrints{MQuater}
\ifx\Semafor\ValueOn
\begin{theorem}
\labelTheorem{quaternion conjugation}
Quaternion conjugation satisfies equation
\ShowEq{quaternion conjugation}
\end{theorem}
\begin{proof}
See the proof of the theorem
\ShowEq{ref quaternion conjugation}
\end{proof}
\fi

\ePrints{1610.309618526,CACAA.06}
\ifx\Semafor\ValueOff
\begin{theorem}
\labelTheorem{Quaternion over real field}
Let
\DrawEq{basis of quaternion}{algebra}
be basis of quaternion algebra $H$.
Then in the basis \eqRef{basis of quaternion}{algebra}, structural constants have form
\ShowEq{structural constants, quaternion}
\end{theorem}
\begin{proof}
See the proof of the theorem
\ShowEq{ref Quaternion over real field}
\end{proof}
\fi

\ePrints{MQuater}
\ifx\Semafor\ValueOn
\begin{theorem}
\labelTheorem{Quaternion over real field, matrix}
Consider quaternion algebra $H$ with the basis
\eqRef{basis of quaternion}{algebra}.
Standard components of linear map over field $R$
and coordinates of this map over field $R$
satisfy relationship
\EqParm{quaternion, 3, B->A}{f=f}
\EqParm{quaternion, 3, A->B}{f=f}
where
\ShowEq{quaternion, 3, 1}
\end{theorem}
\begin{proof}
See the proof of the theorem
\ShowEq{ref Quaternion over real field, matrix}
\end{proof}
\fi

\ePrints{MQuater,1610.309618526,CACAA.06}
\ifx\Semafor\ValueOn

\begin{definition}
\labelDefinition{quaternion maps of conjugation}
{\it
Quaternion algebra has following
\AddIndex{maps of conjugation}{map of conjugation}
\ShowEq{Ex= H}
\ShowEq{Ix= H}
\ShowEq{Jx= H}
\ShowEq{Kx= H}
}
\qed
\end{definition}

\ePrints{1610.309618526,CACAA.06}
\ifx\Semafor\ValueOff
\begin{theorem}
\labelTheorem{aE, quaternion, Jacobian matrix}
We can identify the map
\ShowEq{aE, quaternion}
and matrix
\ShowEq{aE, quaternion, Jacobian matrix}
\ShowEq{aE, quaternion, matrix}
\end{theorem}
\begin{proof}
\TheoremFollows
\ShowEq{ref aE, quaternion, Jacobian matrix}
\end{proof}
\fi

\begin{theorem}
\labelTheorem{linear map of quaternion algebra}
A linear map of quaternion algebra
\ShowEq{f:H->H ij}
has form
\ShowEq{linear map of quaternion algebra, structure}
where $H$\Hyph numbers
\ShowEq{a0123}
are defined by the equality
\ShowEq{a0123=}
\end{theorem}
\begin{proof}
\TheoremFollows
\ShowEq{ref expand linear mapping, quaternion}
\end{proof}

\begin{corollary}
\labelCorollary{linear map of quaternion algebra}
{\it
$H\otimes H$\Hyph module
\ShowEq{L(R;H;H)}
is left $H$\Hyph vector space
and has the basis
\EqParm{e=(E,I,J,K)}{=.}
\qed
}
\end{corollary}

\ePrints{1610.309618526,CACAA.06}
\ifx\Semafor\ValueOff
\begin{theorem}
\labelTheorem{HE is algebra isomorphic to quaternion algebra}
The set
\ShowEq{HE set}
is $R$\Hyph algebra isomorphic to quaternion algebra.
\end{theorem}
\begin{proof}
The theorem follows from equalities
\ShowEq{aE+bE=(a+b)E}
\ShowEq{aE o bE=(ab)E}
based on the theorem
\RefTheorem{linear map of quaternion algebra}.
\end{proof}
\fi
\fi

\ePrints{1506.00061}
\ifx\Semafor\ValueOn
\begin{theorem}
\labelTheorem{ax-xa=1 quaternion algebra}
Equation
\ShowEq{ax-xa=1}
in quaternion algebra does not have solutions.
\end{theorem}
\begin{proof}
\TheoremFollows
\ShowEq{ref ax-xa=1 quaternion algebra}
\end{proof}
\fi
\fi

\def\Preliminary{off}
\ShowEq{SetupRef}

%% file: Preliminary.Eq.tex
\ifx\Preliminary\undefined
\def\Preliminary{off}
\fi

\input{Preliminary.Ref}
\input{Preliminary.Parm}

\ePrints{MVector,MRepro,1601.03259,4975-6381,5410-9916,1305.4547,1310.5591}
\Items{1610.309618526,CACAA.06,5059-9176,1702.01,17.Diary}
\ifx\Semafor\ValueOn
\input{\FilePrefix Stmt.Representation.\TheLanguage}
\fi

\ePrints{MVector,1610.309618526,CACAA.06,1601.03259,4975-6381,5410-9916,1502.04063,5114-6019}
\Items{1702.01,MRepro,17.Diary}
\ifx\Semafor\ValueOn
\input{\FilePrefix Stmt.Module.\TheLanguage}
\fi

\ePrints{MRepro,5148-4632,MVector,5410-9916,1601.03259,1610.309618526,CACAA.06,4975-6381,17.Diary}
\ifx\Semafor\ValueOn
\input{\FilePrefix Stmt.Omega.\TheLanguage}
\fi

\ePrints{1610.309618526,CACAA.06,4975-6381,1601.03259,1702.01}
\ifx\Semafor\ValueOn
\input{\FilePrefix Stmt.Derivative.\TheLanguage}
\fi

\def\Ii{i\in I}
\def\iI{$\Ii$}
\def\Times{A_1\times...\times A_n}

\newcommand\Kn[1]{$k=#1$, ..., $n$}
\def\ATwo{A_2\otimes A_2}
\newcommand\BoxB[1]{$#1\otimes #1$\Hyph}

\def\MatrixA
{
\left(
\begin{array}{rrrr}
\FM 00
&
\FM 01
&
\FM 02
&
\FM 03
\\
\FM 11
&
-\FM 01
&
\FM 13
&
-\FM 12
\\
\FM 22
&
-\FM 23
&
-\FM 20
&
\FM 21
\\
\FM 33
&
\FM 32
&
-\FM 31
&
-\FM 30
\end{array}
\right)
}

\def\MatrixB
{
\left(
\begin{array}{rrrr}
\FC 00
&
-\FC 01
&
-\FC 02
&
-\FC 03
\\
\FC 11
&
\FC 10
&
\FC 13
&
-\FC 12
\\
\FC 22
&
-\FC 23
&
\FC 20
&
\FC 21
\\
\FC 33
&
\FC 32
&
-\FC 31
&
\FC 30
\end{array}
\right)
}

\def\MatrixAB
{
\left(
\begin{array}{rrrr}
1&-1&-1&-1
\\
1&-1&1&1
\\
1&1&-1&1
\\
1&1&1&-1
\end{array}
\right)
}

\def\MatrixBA
{
\frac 14
\left(
\begin{array}{rrrr}
1&1&1&1
\\
-1&-1&1&1
\\
-1&1&-1&1
\\
-1&1&1&-1
\end{array}
\right)
}

\newcommand\EL[1]
{
\left(
\begin{array}{rr@{}rr@{}rr@{}r}
#1^{\gi 0}&-&#1^{\giA}&-&#1^{\gi 2}&-&#1^{\gi 3}
\\
#1^{\giA}&&#1^{\gi 0}&-&#1^{\gi 3}&&#1^{\gi 2}
\\
#1^{\gi 2}&&#1^{\gi 3}&&#1^{\gi 0}&-&#1^{\giA}
\\
#1^{\gi 3}&-&#1^{\gi 2}&&#1^{\giA}&&#1^{\gi 0}
\end{array}
\right)
}

\newcommand\FM[2]{f^{\gi #1}_{\gi #2}}%
\newcommand\FC[2]{\Pf^{\gi{#1#2}}}%
\newcommand\TensorBasis[1]
{e_{1\cdot\gi{#1_1}}\otimes...\otimes e_{n\cdot\gi{#1_n}}}
\newcommand\re{\mathrm{Re}\,} 
\newcommand\im{\mathrm{Im}\,} 

\AddEq [4]{L(A;B)}
{
$\mathcal L(#1;#2\rightarrow #3)$#4
}

\DefEq
{
\symb{|\sigma|}{parity of permutation}{}
\[
\ShowSymbol{parity of permutation}{}
:\sigma\in S(n)\rightarrow \{-1,1\}
\]
}
{sigma->+-}

\DefEq
{
$\{a_1,...,a_n\}$\Pt
}
{set a1n}

\DefEq
{
\[
\ShowSymbol{parity of permutation}{}=
\left\{
\begin{matrix}
1&\textrm{\permutation}\sigma\textrm{ \even}
\\
-1&\textrm{\permutation}\sigma\textrm{ \odd}
\end{matrix}
\right.
\]
}
{parity of permutation}

\DefEq
{
\symb{H}{quaternion algebra}{R}
$\ShowSymbol{quaternion algebra}{R}=R(i,j,k)$
}
{quaternion algebra}

\DefEq
{
\begin{matrix}
e_{\gi 0}=1&e_{\giA}=i&
e_{\gi 2}=j&e_{\gi 3}=k
\end{matrix}
}
{basis of quaternion}

\DefEquation
{
\begin{array}{r@{\ }rr@{\ }rr@{\ }rr@{\ }r}
C_{\gi{00}}^{\gi 0}=&1&
C_{\gi{01}}^{\gi 1}=&1&
C_{\gi{02}}^{\gi 2}=&1&
C_{\gi{03}}^{\gi 3}=&1
\\
\VirtVar
C_{\gi{10}}^{\gi 1}=&1&
C_{\gi{11}}^{\gi 0}=&-1&
C_{\gi{12}}^{\gi 3}=&1&
C_{\gi{13}}^{\gi 2}=&-1
\\
\VirtVar
C_{\gi{20}}^{\gi 2}=&1&
C_{\gi{21}}^{\gi 3}=&-1&
C_{\gi{22}}^{\gi 0}=&-1&
C_{\gi{23}}^{\gi 1}=&1
\\
\VirtVar
C_{\gi{30}}^{\gi 3}=&1&
C_{\gi{31}}^{\gi 2}=&1&
C_{\gi{32}}^{\gi 1}=&-1&
C_{\gi{33}}^{\gi 0}=&-1
\end{array}
}
{structural constants, quaternion}

\DefEq
{
\[
\MatrixAB^{-1}=\MatrixBA
\]
}
{quaternion, 3, 1}

\DefEquation
{
\begin{split}
&
\MatrixA
\\
=&
\MatrixAB
\MatrixB
\end{split}
}
{quaternion, 3, B->A}

\DefEquation
{
\begin{split}
&
\MatrixB
\\
=&
\MatrixBA
\MatrixA
\end{split}
}
{quaternion, 3, A->B}

\DefEq
{
$\mathcal L(C;C\rightarrow C)=CE$.
}
{L(CCC)=CE}

\DefEquation
{
|x|^2=x\overline x=(x^{\gi 0})^2+(x^{\giA})^2+(x^{\gi 2})^2+(x^{\gi 3})^2
}
{norm of quaternion}

\DefEquation
{
x^{-1}=|x|^{-2}\overline x
}
{inverce quaternion}

\DefEq
{
\[
x=x^{\gi 0}+x^{\giA}i+x^{\gi 2}j+x^{\gi 3}k
\]
\[
x^{\gi 0}, x^{\giA}, x^{\gi 2}, x^{\gi 3}\in R
\]
}
{quaternion algebra, element}

\DefEquation
{
\begin{array}{c|ccc}
&i&j&k\\
\hline
i&-1&k&-j\\
j&-k&-1&i\\
k&j&-i&-1\\
\end{array}
}
{product of quaternions}

\DefEq
{
$\Basis e_i$
}
{tensor product of algebras, basis i}

\DefEq
{
\symb{a^{\gi{i_1...i_n}}}{standard component of tensor}{}
}
{standard component of tensor}

\DefEquation
{
\circ:(g,f)\in\mathcal L(D;A\rightarrow A)\times\mathcal L(D;A\rightarrow A)
\rightarrow g\circ f\in\mathcal L(D;A\rightarrow A)
}
{module L(A;A) is algebra}

\DefEq
{
\[
\xymatrix{
&A\ar[rd]^g
\\
A\ar[ru]^f\ar[rr]^{g\circ f}&&A
}
\]
}
{product of linear map, algebra 1}

\DefEquation
{
a=
\ShowSymbol{standard component of tensor}{}
\TensorBasis i
}
{tensor canonical representation, algebra}

\AddEq{k,1n}
{
$k$, $k=1$, ..., $n$,
}

\AddEq{n,1.}
{
$n$, $n=1$, ...,
}

\newcommand{\Tensor}[1]{#1_1\otimes...\otimes #1_n}

\AddEq[3]{set f:A->B}
{
\[
\{#1_i:#2\rightarrow #3,i\in I\}
\]
}

\AddEq[1]{set Bi}
{
$\{#1_i,i\in I\}$
}

\AddEq{product in category diagram}
{
\[
\xymatrix{
P\ar[r]^{f_i}&B_i&f_i\circ h=g_i
\\
R\ar[ur]_{g_i}\ar[u]^h
}
\]
}

\AddEq{coproduct in category diagram}
{
\[
\xymatrix{
P\ar[d]_h&B_i\ar[l]_{f_i}\ar[dl]^{g_i}&h\circ f_i=g_i
\\
R
}
\]
}

\AddEq{product in category, 1 n}
{
\symb[Pi-1]{\prod_{i=1}^nB_i}{product in category}{i 1 n}
\symb[B-0]{B_1\times...\times B_n}{product in category}{1 n}
\[
P=\ShowSymbol{product in category}{i 1 n}
=\ShowSymbol{product in category}{1 n}
\]
}

\AddEq{coproduct in category, 1 n}
{
\symb[Pi-1]{\coprod_{i=1}^nB_i}{coproduct in category}{i 1 n}
\symb[B-0]{B_1\coprod...\coprod B_n}{coproduct in category}{1 n}
\[
P=\ShowSymbol{coproduct in category}{i 1 n}
=\ShowSymbol{coproduct in category}{1 n}
\]
}

\DefEq
{
\[A=\prod_{\Ii}A_i\]
}
{Cartesian product of sets}

\DefEq
{
\((A_i,\Ii)\)
}
{Ai iI}

\DefEq
{
\[p_i:A\rightarrow A_i\]
}
{projection on i factor}

\DefEq
{
\[p_i:A\rightarrow A_i\ \ \ \Ii\]
}
{p:A->Ai i in I}

\DefEq
{
\begin{align*}
\re&:A\rightarrow A
\\
\im&:A\rightarrow A
\end{align*}
}
{Re Im A->A}

\DefEq
{
\symb{\re d}{scalar of algebra}{}
\symb{\im d}{vector of algebra}{}
}
{Re Im A->F 0}

\DefEquation
{
\begin{matrix}
\ShowSymbol{scalar of algebra}{}=d^{\gi 0}
&\ShowSymbol{vector of algebra}{}=d-d^{\gi 0}
&d\in D
&d=d^{\gii}e_{\gii}
\end{matrix}
}
{Re Im A->F 1}

\DefEq
{
$\ShowSymbol{scalar of algebra}{}$
}
{Re Im A->F 2}

\DefEq
{
$\ShowSymbol{vector of algebra}{}$
}
{Re Im A->F 3}

\DefEq
{
\[
F=\{d\in A:\re d=d\}
\]
}
{Re Im A->F 4}

\DefEq
{
\symb{\re A}{scalar algebra of algebra}1
}
{scalar algebra of algebra}

\DefEq
{
\symb{\im A}{vector module of algebra}{}
\begin{equation}
\EqLabel{vector module of algebra}
\ShowSymbol{vector module of algebra}{}=\{d\in A:\re d=0\}
\end{equation}
}
{vector module of algebra}

\DefEquation
{
C_{\gi{0k}}^{\gil}=C^{\gil}_{\gi{k0}}=\delta^{\gik}_{\gil}
}
{structural constant re algebra, 1}

\DefEquation
{
d=\re d+\im d
}
{d Re Im}

\DefEq
{
\symb{d^*}{conjugation in algebra}{}
}
{conjugation in algebra}

\DefEquation
{
\ShowSymbol{conjugation in algebra}{}=\re d-\im d
}
{conjugation in algebra, 0}

\DefEquation
{
(cd)^*=d^*\,c^*
}
{conjugation in algebra, 1}

\DefEquation
{
\begin{matrix}
C^{\gi 0}_{\gi{kl}}=C^{\gi 0}_{\gi{lk}}
&
C^{\gi p}_{\gi{kl}}=-C^{\gi p}_{\gi{lk}}
\end{matrix}
}
{conjugation in algebra, 5}

\DefEq
{
\[
\begin{matrix}
\giA\le\gik\le\gin&\giA\le\gil\le\gin&\giA\le\gi p\le\gin
\end{matrix}
\]
}
{conjugation in algebra, 4}

\DefEq
{
$p(x)=p_1\circ x$
}
{p=p circ x}

\DefEquation
{
\begin{split}
r(x)&=r_0+q_{1\cdot 0} p(x)q_{1\cdot 1}
+q_{2\cdot 0}(x) p(x)q_{2\cdot 1}+...
+q_{k\cdot 0}(x) p(x)q_{k\cdot 1}
\\
&=r_0+(q_{1\cdot 0}\otimes q_{1\cdot 1})\circ p(x)
+(q_{2\cdot 0}(x)\otimes q_{2\cdot 1})\circ p(x)
\\&+...
+(q_{k\cdot 0}(x)\otimes q_{k\cdot 1})\circ p(x)
\end{split}
}
{r=+q circ p}

\DefEquation
{
p(x)=p_0+p_1\circ x
}
{p=p+p circ x}

\DefEq
{
\[
r(x)=r_0+r_1\circ x+...+r_k\circ x^k
\]
}
{r power k}

\DefEquation
{
\begin{split}
r(x)&=r_0
-((r_{1\cdot 0\cdot s}\otimes r_{1\cdot 1\cdot s})\circ
p^{-1}_1)\circ p_0
\\
&-(((r_{2\cdot 0\cdot s}\circ x)\otimes r_{2\cdot 1\cdot s})\circ
p^{-1}_1)\circ p_0
\\
&-...-(((r_{k\cdot 0\cdot s}\circ x^{k-1})
\otimes r_{k\cdot 1\cdot s})\circ
p^{-1}_1)\circ p_0
\\
&+((r_{1\cdot 0\cdot s}\otimes r_{1\cdot 1\cdot s})\circ
p^{-1}_1)\circ p(x)
\\
&+(((r_{2\cdot 0\cdot s}\circ x)\otimes r_{2\cdot 1\cdot s})\circ
p^{-1}_1)\circ p(x)
\\
&+...+(((r_{k\cdot 0\cdot s}\circ x^{k-1})
\otimes r_{k\cdot 1\cdot s})\circ
p^{-1}_1)\circ p(x)
\\
&=r_0
-((r_{1\cdot 0\cdot s}\otimes r_{1\cdot 1\cdot s}
+(r_{2\cdot 0\cdot s}\circ x)\otimes r_{2\cdot 1\cdot s}
\\
&+...+(r_{k\cdot 0\cdot s}\circ x^{k-1})
\otimes r_{k\cdot 1\cdot s})\circ
p^{-1}_1)\circ p_0
\\
&+((r_{1\cdot 0\cdot s}\otimes r_{1\cdot 1\cdot s}
+(r_{2\cdot 0\cdot s}\circ x)\otimes r_{2\cdot 1\cdot s}
\\
&+...+(r_{k\cdot 0\cdot s}\circ x^{k-1})
\otimes r_{k\cdot 1\cdot s})\circ
p^{-1}_1)\circ p(x)
\end{split}
}
{r=+q circ p 1}

\DefEq
{
\[
ax-xa=1
\]
}
{ax-xa=1}

\DefEquation
{
r(x)=q_{i\cdot 0}(x)p(x)q_{i\cdot 1}(x)
=(q_{i\cdot 0}(x)\otimes q_{i\cdot 1}(x))\circ p(x)
}
{r=qpq(x)}

\DefEquation
{
a_1...a_n\omega=(p_i(a_1)...p_i(a_n)\omega,\Ii)
}
{omega(ai)=(omega ai)}

\DefEquation
{
\xymatrix{
A\ar[r]^{p_i}&A_i&p_i\circ \omega=g_i
\\
A^n\ar[ur]_{g_i}\ar[u]^{\omega}
}
}
{operation is defined componentwise, diagram}

\DefEquation
{
a_1...a_n\omega=(a_{1i}...a_{ni}\omega,\Ii)
}
{operation is defined componentwise}

\DefEq
{
\((p_i(a),\Ii)\)
}
{tuple represent A number}

\AddEq[2]{omega in Omega}
{
\(\omega\in\Omega_{#1}\)#2
}

\DefEq
{
h=f_1...f_n\omega
}
{h=f1...fn omega}

\DefEq
{
\(f_X(g_{1\cdot n})(g_{2\cdot n})\)
}
{f(g1n)(g2n)}

\DefEq
{
\(a_1=(a_{1i},\Ii)\), ..., \(a_n=(a_{ni},\Ii)\)
}
{a=ai 1n}

\DefEq
{
\((B_i,\Ii)\)
}
{Bi iI}

\DefEq
{
f(a_i,\Ii)=(f_i(a_i),\Ii)
}
{f:A->B=}

\DefEq
{
\[
g_i(a_1, ..., a_n)=p_i(a_1)...p_i(a_n)\omega
\]
}
{gi()=}

\DefEq
{
\[f_i:A_i\rightarrow B_i\]
}
{f:A->B i}

\DefEquation
{
\xymatrix
{
B\ar[r]^{p'_i}&B_i
\\
A\ar[u]^f\ar[r]_{p_i}&A_i\ar[u]_{f_i}
}
}
{homomorphism of Cartesian product of Omega algebras diagram}

\DefEquation
{
\xymatrix
{
B\ar[rrr]^{p'_i}\ar@{}[dr]^(.6){(1)}&&&B_i
\\
&&
\\
A\ar[uu]^f\ar[uurrr]^{g_i}\ar[rrr]_{p_i}&&&A_i\ar[uu]_{f_i}\ar@{}[ul]^(.8){(2)}
}
}
{homomorphism of Cartesian product of Omega algebras}

\DefEquation
{
b=f(a)\in B
}
{b=f(a)}

\DefEquation
{
\begin{matrix}
b=(b_i,\Ii)&b_i=p'_i(b)\in B_i
\end{matrix}
}
{b=p(b)i}

\DefEquation
{
b_i=g_i(b)
}
{b=g(a)i}

\DefEq
{
\[b_i=f_i(a_i)\]
}
{b=f(a)i}

\DefEq
{
\(b_1=(b_{1i},\Ii)\), ..., \(b_n=(b_{ni},\Ii)\)
}
{b=bi 1n}

\DefEq
{
\begin{align*}
f(a_1...a_n\omega)&=f(a_{1i}...a_{ni}\omega,\Ii)
\\&=(f_i(a_{1i}...a_{ni}\omega),\Ii)
\\&=((f_i(a_{1i}))...(f_i(a_{ni})),\Ii)
\\&=(b_{1i}...b_{ni}\omega,\Ii)
\end{align*}
\[
f(a_1)...f(a_n)\omega=b_1...b_n\omega
=(b_{1i}...b_{ni}\omega,\Ii)
\]
}
{f:A->B omega}

\DefEquation
{
\begin{matrix}
a=(a_i,\Ii)&a_i=p_i(a)\in A_i
\end{matrix}
}
{a=p(a)i}

\DefEq
{
\(p_i\), \(p'_i\)
}
{pi p'i}

\DefEq
{
\symb{f(A_1)a_2}{orbit of representation}{}
\[
\ShowSymbol{orbit of representation}{}=\{b_2=f(a_1)a_2:a_1\in A_1\}
\]
}
{orbit of representation}

\DefEq
{
\[f\rightarrow g\]
}
{f->g}

\DefEq
{
\[h:S_1\rightarrow S_2\]
}
{f->g 1}

\DefEq
{
\[
\xymatrix{
&S_1\ar[dd]^h
\\
\Times
\ar[ru]^f\ar[rd]_g
\\
&S_2
}
\]
}
{polylinear maps category, diagram}

\DefEq
{
$a_1\ne b_1$, $a_1$, $b_1\in A_1$,
}
{a1 ne b1}

\DefEq
{
\[
F:A_1\cup A_2\rightarrow B_1\cup B_2
\]
}
{F:A1+A2->B1+B2}

\DefEq
{
\[
F(A_1)=B_1 \ \ \ \ F(A_2)=B_2
\]
}
{F:A1+A2->B1+B2 1}

\AddEq [3]{map r,R}
{
$(#1\ \ #2)$#3
}

\DefEq
{
\[
r_i:A_i\rightarrow B_i
\]
}
{ri:A->B}

\DefEq
{
$i=1$, $2$\Pt
}
{i=1,2}

\DefEq
{
\symb{\mathrm{ker}\,f}{kernel of homomorphism}0
$\ShowSymbol{kernel of homomorphism}0
=f\circ f^{-1}$
}
{kernel of homomorphism}

\DefEq
{
$\mathrm{ker}\,t_i
=t_i\circ t_i^{-1}$
}
{kernel of homomorphism i}

\DefEq
{
$A/\mathrm{ker}\,f$
}
{A/ker f}

\DefEq
{
\[
p:a\in A\rightarrow a^{\mathrm{ker}\,f}\in A/\mathrm{ker}\,f
\]
}
{p:A->/ker}

\DefEq
{
$a$, $b\in \pA$\Pt
}
{ab in A}

\DefEq
{
\[f(a)(m)\ne f(b)(m)\]
}
{fam ne fbm}

\DefEq
{
\[f(a_1)(a_2)\ne a_2\]
}
{fam ne m}

\DefEq
{
$B_2\subset A_2$
}
{B2 subset A2}

\DefEq
{
\[f(a_1)(a_2)=f(b_1)(a_2)\]
}
{faa=fba}

\DefEq
{
$a_{\Pn}\in A_{\Pn}$\Pt
}
{a in A}

\DefEq
{
\symb{A^{\otimes n}}{tensor power of algebra}{}
\[
\begin{matrix}
\ShowSymbol{tensor power of algebra}{}
=\Tensor A
&
A_1=...=A_n=A
\end{matrix}
\]
}
{tensor power of algebra}

\DefEq
{
\[
g:A_1\times...\times A_n\rightarrow V
\]
}
{map g, algebra, tensor product}

\DefEq
{
\[
h:\Tensor A\rightarrow V
\]
}
{map h, algebra, tensor product}

\DefEquation
{
\xymatrix{
&\Tensor A\ar[dd]^h
\\
\Times
\ar[ru]^f\ar[rd]_g
\\
&V
}
}
{map gh, algebra, tensor product}

\DefEquation
{
h\circ(\Tensor a)=g\circ(a_1,...,a_n)
}
{g=h, algebra, tensor product}

\DefEquation
{
\begin{array}{rr}
e_{\gi 0}=1&e_{\giA}=i
\end{array}
}
{basis of complex field}

\DefEquation
{
e_{\giA}^2=-e_{\gi 0}
}
{product of complex field}

\DefEq
{
$\mathcal L(R;C\rightarrow C)$
}
{L(R;C;C)}

\DefEq
{
$\mathcal L(R;H\rightarrow H)$
}
{L(R;H;H)}

\DefEquation
{
E\circ(x_{\gi 0}+x_{\gi 1}i)=x_{\gi 0}+x_{\gi 1}i
}
{Ex= C}

\DefEquation
{
x=E\circ(x_{\gi 0}+x_{\gi 1}i+x_{\gi 2}j+x_{\gi 3}k)
=x_{\gi 0}+x_{\gi 1}i+x_{\gi 2}j+x_{\gi 3}k
}
{Ex= H}

\DefEq
{
\begin{align*}
P^{\gi 0}\circ a&=a^{\gi 0}&P^{\gi 0\cdot}_{}{}^{\gi 0}_{\gi 0}&=1&
R^{\gi 0}\circ a&=a^{\gi 0}&R^{\gi 0\cdot}_{}{}^{\gi 0}_{\gi 0}&=1
\\
P^{\giA}\circ a&=a^{\giA}i&P^{\giA\cdot}_{}{}^{\giA}_{\giA}&=1&
R^{\giA}\circ a&=a^{\giA}&R^{\giA\cdot}_{}{}^{\gi 0}_{\giA}&=1
\end{align*}
}
{projection maps, complex field}

\DefEquation
{
I\circ(x_{\gi 0}+x_{\gi 1}i)=x_{\gi 0}-x_{\gi 1}i
}
{Ix= C}

\DefEquation
{
x^{*_1}=I\circ(x_{\gi 0}+x_{\gi 1}i+x_{\gi 2}j+x_{\gi 3}k)
=x_{\gi 0}-x_{\gi 1}i+x_{\gi 2}j+x_{\gi 3}k
}
{Ix= H}

\DefEquation
{
x^{*_2}=J\circ(x_{\gi 0}+x_{\gi 1}i+x_{\gi 2}j+x_{\gi 3}k)
=x_{\gi 0}+x_{\gi 1}i-x_{\gi 2}j+x_{\gi 3}k
}
{Jx= H}

\DefEquation
{
x^{*_3}=K\circ(x_{\gi 0}+x_{\gi 1}i+x_{\gi 2}j+x_{\gi 3}k)
=x_{\gi 0}+x_{\gi 1}i+x_{\gi 2}j-x_{\gi 3}k
}
{Kx= H}

\DefEq
{
\begin{align}
f&=a_0\circ E+a_1\circ I
\EqLabel{linear map of complex field, structure, 1}
\\
f\circ x&=a_0x+a_1\overline x
\EqLabel{linear map of complex field, structure, 2}
\end{align}
}
{linear map of complex field, structure}

\DefEq
{
\begin{align}
f&=a_0\circ E+a_1\circ I+a_2\circ J+a_3\circ K
\EqLabel{linear map of quaternion algebra, structure, 1}
\\
f\circ x&=a_0x+a_1\circ I\circ x+a_2\circ J\circ x+a_3\circ K\circ x
\EqLabel{linear map of quaternion algebra, structure, 2}
\\
&=a_0x+a_1x^{*_1}+a_2x^{*_2}+a_3x^{*_3}
\end{align}
}
{linear map of quaternion algebra, structure}

\DefEq
{
\[
a_i=a_i^{\gi 0}+a_i^{\gi 1}i+a_i^{\gi 2}j+a_i^{\gi 3}k\ \ \ i=0,1,2,3
\]
}
{a0123}

\DefEq
{
\[
a_i=a_i^{\gi 0}+a_i^{\gi 1}i\ \ \ i=0,1
\]
}
{a01}

\DefEquation
{
\begin{pmatrix}
a_0^{\gi 0}&a_0^{\gi 1}&a_0^{\gi 2}&a_0^{\gi 3}
\\
a_1^{\gi 0}&a_1^{\gi 1}&a_1^{\gi 2}&a_1^{\gi 3}
\\
a_2^{\gi 0}&a_2^{\gi 1}&a_2^{\gi 2}&a_2^{\gi 3}
\\
a_3^{\gi 0}&a_3^{\gi 1}&a_3^{\gi 2}&a_3^{\gi 3}
\end{pmatrix}
=\frac 12
\left(
\begin{array}{r@{\,\ }r@{\,\ }r@{\,\ }r}
-1&1&1&1
\\
1&-1&0&0
\\
1&0&-1&0
\\
1&0&0&-1
\end{array}
\right)
\begin{pmatrix}
f_{\gi 0}^{\gi 0}&f_{\gi 0}^{\gi 1}&f_{\gi 0}^{\gi 2}&f_{\gi 0}^{\gi 3}
\\
f_{\gi 1}^{\gi 1}&-f_{\gi 1}^{\gi 0}&-f_{\gi 1}^{\gi 3}&f_{\gi 1}^{\gi 2}
\\
f_{\gi 2}^{\gi 2}&f_{\gi 2}^{\gi 3}&-f_{\gi 2}^{\gi 0}&-f_{\gi 2}^{\gi 1}
\\
f_{\gi 3}^{\gi 3}&-f_{\gi 3}^{\gi 2}&f_{\gi 3}^{\gi 1}&-f_{\gi 3}^{\gi 0}
\end{pmatrix}
}
{a0123=}

\DefEquation
{
\begin{pmatrix}
a_0^{\gi 0}&a_0^{\gi 1}
\\
a_1^{\gi 0}&a_1^{\gi 1}
\end{pmatrix}
=\frac 12
\left(
\begin{array}{r@{\,\ }r@{\,\ }r@{\,\ }r}
1&1
\\
1&-1
\end{array}
\right)
\begin{pmatrix}
f_{\gi 0}^{\gi 0}&f_{\gi 0}^{\gi 1}
\\
f_{\gi 1}^{\gi 1}&-f_{\gi 1}^{\gi 0}
\end{pmatrix}
}
{a01=}

\DefEq
{
$\Basis e=(E,I)$\Pt
}
{e=(E,I)}

\DefEq
{
$\Basis e=(E,I,J,K)$\Pt
}
{e=(E,I,J,K)}

\DefEq
{
\begin{align*}
P^{\gi 0}&=\frac 12\circ E+\frac 12\circ I&
R^{\gi 0}&=\frac 12\circ E+\frac 12\circ I
\\
P^{\giA}&=\frac 12\circ E-\frac 12\circ I&
R^{\giA}&=-\frac i2\circ E+\frac i2\circ I
\end{align*}
}
{projection mappings 1, complex field}

\DefEq
{
\[f:H\rightarrow H\ \ \ y^{\gii}=f^{\gii}_{\gij}x^{\gij}\]
}
{f:H->H ij}

\DefEq
{
\[f:C\rightarrow C\ \ \ y^{\gii}=f^{\gii}_{\gij}x^{\gij}\]
}
{f:C->C y=}

\DefEq
{
\[f:C\rightarrow C\]
}
{f:C->C}

\DefEq
{
\[(aE)\circ x+(bE)\circ x=ax+bx=(a+b)x=((a+b)E)\circ x\]
}
{aE+bE=(a+b)E}

\DefEq
{
\[(aE)\circ (bE)\circ x=(aE)\circ (bx)=a(bx)=(ab)x=((ab)E)\circ x\]
}
{aE o bE=(ab)E}

\DefEq
{
\[(aI)\circ (bI)\circ x=(aI)\circ (b\overline x)=a\overline{(b\overline x)}
=(a\overline b)x=((a\overline b)E)\circ x\]
}
{aI o bI=...E}

\DefEq
{
\[CE=\{aE:a\in C\}\]
}
{CE set}

\DefEq
{
\[HE=\{aE:a\in H\}\]
}
{HE set}

\DefEq
{
\[CI=\{aI:a\in C\}\]
}
{CI set}

\DefEquation
{
\begin{array}{r@{}rr@{}r}
C_{\gi{00}}^{\gi 0}=&1&C_{\gi{01}}^{\giA}=&1
\\
\VirtVar
C_{\gi{10}}^{\giA}=&1&C_{\gi{11}}^{\gi 0}=&-1
\end{array}
}
{structural constants of complex field}

\DefEquation
{
\begin{split}
\overline x&
=-\frac 12(1\otimes 1+i\otimes i+j\otimes j+k\otimes k)\circ x
\\
&=-\frac 12(x+ixi+jxj+kxk)
\end{split}
}
{quaternion conjugation}

\DefEquation
{
\overline x=x^{\gi 0}-x^{\giA}i-x^{\gi 2}j-x^{\gi 3}k
}
{conjugate to the quaternion}

\DefEquation
{
\begin{split}
f_{\gi 0}^{\gi 0}&=\hphantom{-\,}f_{\giA}^{\giA}
\\
f_{\gi 0}^{\giA}&=-f_{\giA}^{\gi 0}
\end{split}
}
{complex field over real field}

\DefEquation
{
\begin{split}
f_{\gi 0}^{\gi 0}&=-f_{\giA}^{\giA}
\\
f_{\gi 0}^{\giA}&=\hphantom{-\,}f_{\giA}^{\gi 0}
\end{split}
}
{matrix of linear map f in CI}

\DefEq
{
\[
(b_{\gi 0}+b_{\giA}i)I(x_{\gi 0}+x_{\giA}i)
=(b_{\gi 0}+b_{\giA}i)(x_{\gi 0}-x_{\giA}i)
=b_{\gi 0}x_{\gi 0}+b_{\giA}x_{\giA}
+(-b_{\gi 0}x_{\giA}+b_{\giA}x_{\gi 0})i
\]
\[
\begin{pmatrix}
b_{\gi 0}&b_{\giA}
\\
b_{\giA}&-b_{\gi 0}
\end{pmatrix}
\begin{pmatrix}
x_{\gi 0}\\x_{\giA}
\end{pmatrix}
=
\begin{pmatrix}
b_{\gi 0}x_{\gi 0}+b_{\giA}x_{\giA}
\\
b_{\giA}x_{\gi 0}-b_{\gi 0}x_{\giA}
\end{pmatrix}
\]
}
{congugate map complex field, 3}

\DefEq
{
\[
x=x^{\gi 0}+x^{\giA}i\rightarrow
f=f^{\gi 0}(x^{\gi 0},x^{\giA})+f^{\giA}(x^{\gi 0},x^{\giA})i
\]
}
{map of complex variable}

\DefEquation
{
\begin{split}
\frac{\partial f^{\giA}}{\partial x^{\gi 0}}
&=-\frac{\partial f^{\gi 0}}{\partial x^{\giA}}
\\[5pt]
\frac{\partial f^{\gi 0}}{\partial x^{\gi 0}}
&=\hphantom{-\,}\frac{\partial f^{\giA}}{\partial x^{\giA}}
\end{split}
}
{Cauchy Riemann equations, complex field, 1}

\DefEquation
{
\frac{\partial f^{\gi 0}}{\partial x^{\gi 0}}
+i\frac{\partial f^{\giA}}{\partial x^{\gi 0}}
+i\left(\frac{\partial f^{\gi 0}}{\partial x^{\giA}}
+i\frac{\partial f^{\giA}}{\partial x^{\giA}}\right)=
\frac{\partial f^{\gi 0}}{\partial x^{\gi 0}}
+i\frac{\partial f^{\giA}}{\partial x^{\gi 0}}
+i\frac{\partial f^{\gi 0}}{\partial x^{\giA}}
-\frac{\partial f^{\giA}}{\partial x^{\giA}}=
0
}
{Cauchy Riemann equations, complex field, 2, 1}

\DefEq
{
\[
\begin{matrix}
a\circ E:H\rightarrow H&
a=a^{\gi 0}+a^{\gi 1}i+a^{\gi 2}j+a^{\gi 3}k\in H
\end{matrix}
\]
}
{aE, quaternion}

\DefEq
{
\symb{E_{l\cdot a}}{aE, quaternion, Jacobian matrix}{}
}
{aE, quaternion, Jacobian matrix}

\DefEquation
{
\ShowSymbol{aE, quaternion, Jacobian matrix}{}=\EL a
}
{aE, quaternion, matrix}

\DefEquation
{
\frac{\partial f}{\partial x}
=\frac 12\left(\frac{\partial f}{\partial x^{\gi 0}}
-i\frac{\partial f}{\partial x^{\gi 1}}\right)E
}
{df=...dx= CE}

\DefEquation
{
\begin{split}
dx^{\gi 0}&=\frac 12(E+I)\circ dx
\\
dx^{\gi 1}&=-\frac i2(E-I)\circ dx
\end{split}
}
{dx01=}

\DefEquation
{
\begin{split}
df&=\frac 12\frac{\partial f}{\partial x^{\gi 0}}(E+I)\circ dx
-\frac i2\frac{\partial f}{\partial x^{\gi 1}}(E-I)\circ dx
\\
&=\frac 12\left(\left(\frac{\partial f}{\partial x^{\gi 0}}
-i\frac{\partial f}{\partial x^{\gi 1}}\right)E
+\left(\frac{\partial f}{\partial x^{\gi 0}}
+i\frac{\partial f}{\partial x^{\gi 1}}\right)I\right)\circ dx
\end{split}
}
{df=...dx}

\DefEq
{
\[
df=\frac{\partial f}{\partial x^{\gi 0}}dx^{\gi 0}
+\frac{\partial f}{\partial x^{\gi 1}}dx^{\gi 1}
\]
}
{df=df/dx0+df/dx1}

\DefEq
{
\[dx=dx^{\gi 0}+dx^{\gi 1}i\]
}
{dx=dx0+dx1}

\DefEq
{
\[
\begin{pmatrix}
\displaystyle\frac{\partial f^{\gi 0}}{\partial x^{\gi 0}}
&
\displaystyle\frac{\partial f^{\gi 0}}{\partial x^{\giA}} 
\\
\VirtFrac
\displaystyle\frac{\partial f^{\giA}}{\partial x^{\gi 0}}
&
\displaystyle\frac{\partial f^{\giA}}{\partial x^{\giA}}
\end{pmatrix}
\]
}
{Jacobian of map of complex variable}

\DefEq
{
\[
y^{\gi i}=x^{\gi j}f_{\gi j}^{\gi i}
\]
}
{yi=xj fji}

\DefEq
{
\symb{B^A}{Cartesian power}1
}
{Cartesian power}

\DefEq
{
$m_1\equiv m_2(\mathrm{mod} S)$
}
{m1 m2 modS}

\DefEq
{
\[f\circ(a_1,...,a_n)=\Tensor a\]
}
{fxa=oxa}

\DefEq
{
\[
\sigma:a\in A\rightarrow \sigma(a)\in A\ \ \ A=\{a_1,...,a_n\}
\]
}
{a->sigma a}

\DefEquation
{
\sigma=
\begin{pmatrix}
a_1&...&a_n
\\
\sigma(a_1)&...&\sigma(a_n)
\end{pmatrix}
}
{permutation as matrix}

\DefEquation
{
\sigma=
\begin{pmatrix}
\sigma(a_1)&...&\sigma(a_n)
\end{pmatrix}
}
{permutation as matrix 2}

\DefEq
{
\[f:A_1\times...\times A_n\rightarrow\Tensor A\]
}
{f:xA->oxA}

\DefEq
{
\symb{\Tensor A}{tensor product}1
}
{tensor product of algebras}

\DefEq
{
\symb{\Omega(n)}{set of n-ary operators}{}
}
{set of n-ary operators}

\DefEq
{
$\omega\in\Omega(n)$\Pt
}
{omega n ari}

\DefEq
{
$\omega\in\Omega(2)$
}
{omega 2 ari}

\DefEq
{
$\omega\in\Omega(2)$.
}
{omega 2 ari.}

\DefEquation
{
ea\omega=a
}
{left neutral element}

\DefEquation
{
ae\omega=a
}
{right neutral element}

\DefEq
{
\[ab\omega=ba\omega\]
}
{commutative operation}

\DefEq
{
\[
a(bc\omega)\omega=(ab\omega)c\omega
\]
}
{associative operation}

\DefEq
{
\symb{A_{\Omega}}{Omega-algebra}1
}
{Omega-algebra}

\DefEq
{
$\omega\in A^{A^n}$.
}
{o in AAn}

\DefEq
{
\[\Omega(n)\rightarrow A^{A^n}\ \ \ n\in N\]
}
{O(n)->AAn}

\AddEq [1]{B subset A}
{
$B\subseteq A$#1
}

\DefEq
{
\[
\ShowSymbol{set of n-ary operators}{}=
\{\omega\in\Omega:a(\omega)=n\}
\]
}
{set of n-ary operators =}

\DefEq
{
$a(\omega)=n$,
}
{a(o)=n}

\DefEq
{
$\omega\in\Omega$
}
{o in O}

\DefEq
{
\[\omega_A|B=\omega_B\]
}
{oAB=oB}

\DefEq
{
\[a:\Omega\rightarrow N\]
}
{a:O->N}

\DefEq
{
\symb{\Omega}{operator domain}1
}
{operator domain}

\DefEq
{
$\omega(a_1,...,a_n)$, $a_1...a_n\omega$
}
{a1no=oa1n}

\DefEq
{
$b_1...b_n\omega\in B$,
}
{b1no in B}

\DefEq
{
$\b_1$, ..., $\b_n\in\B$,
}
{b1n in B}

\DefEq
{
\[\omega:A^n\rightarrow A\]
}
{o:An->A}

\DefEq
{
\[
\begin{matrix}
f_1:A\rightarrow S_1&\mathrm{ker}\,f_1\supseteq N
\\
f_2:A\rightarrow S_2&\mathrm{ker}\,f_2\supseteq N
\end{matrix}
\]
}
{maps category}

\DefEq
{
\[
\xymatrix{
&S_1\ar[dd]^h
\\
A
\ar[ru]^{f_1}\ar[rd]_{f_2}
\\
&S_2
}
\]
}
{maps category, diagram}

\DefEq
{
\[
\xymatrix
{
&A/N\ar[dd]^h
\\
A\ar[ur]^{j=\mathrm{nat}\,N}\ar[dr]_f
\\
&S
}
\]
}
{maps category, universal, diagram}

\DefEquation
{
\mathrm{ker}\,f\supseteq N
}
{maps category, universal, ker}

\DefEq
{
\[
\mathrm{nat}\,N:A\rightarrow A/N
\]
}
{maps category, universal}

\DefEq
{
\[j(a_1)=j(a_2)\]
}
{maps category 1}

\DefEq
{
\[f(a_1)=f(a_2)\]
}
{maps category 2}

\DefEq
{
\[h(\BlueText{j(b)})=f(b)\]
}
{maps category, h}

\DefEq
{
\[
(f+g)(a)=f(a)+g(a)
\]
}
{(f+g)(a)=}

\DefEq
{
\[
f(a+b)(x)=f(a)(x)+f(b)(x)
\]
}
{sum of transformations of Abelian group, 1}

\ePrints{1502.04063,5114-6019}
\ifx\Semafor\ValueOff
\DefEquation
{
f(a)(m)=f(b)(m)
}
{representation of ring, 1}

\DefEq
{
\[
f(a-b)(m)=0
\]
}
{representation of ring, 2}
\fi

\AddEq{polynomials p,r}
{
$p\circ F_{[1]}\circ x^n$, $r\circ F_{[2]}\circ x^m$
}

\DefEq
{
\[
(p\circ F_{[1]}\circ x^n)(r\circ F_{[2]}\circ x^m)
=(p\underline{\otimes}r)\circ
(F_{[1](1)},...,F_{[1](n)},F_{[2](1)},...,F_{[2](m)})
\circ x^{n+m}
\]
}
{pr=p circ r}

\DefEquation
{
p(x)=
(a_{k\cdot 0}\circ x^{k-1})
((1\otimes a_{k\cdot 1})\circ x)
}
{p(x)=a k-1 k 1}

\DefEquation
{
p(x)=
((a_{k\cdot 0}\circ x^{k-1})
\otimes a_{k\cdot 1})\circ x
}
{p(x)=a k-1 k 2}

\DefEq
{
$A_{ij}$, $i=1$, ..., $n$, $j=1$, ..., $m$,
}
{Aij}

\DefEq
{
$\Omega_{ij}$\Hyph%
}
{Omegaij}

\DefEq
{
$f_1$, ..., $f_n\in B^A$,
}
{f1n in B**A}

\DefEq
{
$\Hom(\emptyset;A\rightarrow B)=B^A$.
}
{Hom empty A B=B**A}

\DefEq
{
$\End(\emptyset;A)=A^A$.
}
{End empty A=A**A}

\DefEq
{
$\End(\Omega;A)=\Hom(\Omega;A\rightarrow A)$
}
{End A=Hom AA}

\DefEq
{
\symb{\Hom(\Omega;A\rightarrow B)}{set of homomorphisms}1
}
{set of homomorphisms}

\DefEq
{
\symb{\End(\Omega;A)}{set of endomorphisms}1
}
{set of endomorphisms}

\DefEq
{
(f_1...f_n\omega)(x)=f_1(x)...f_n(x)\omega
}
{f1n omega=}

\DefEq
{
\[t:A\rightarrow A\]
}
{t:A->A}

\DefEq
{
\[
\xymatrix@C=15pt{
f:\Times\ar[r]&S_1
&
g:\Times\ar[r]&S_2
}
\]
}
{polylinear maps category}

\DefEq
{
\[f:A_1\rightarrow\End(\Omega_2;A_2)\]
}
{representation of algebra}

\DefEq
{
\[a_2'=a_2R(a_1)=a_2a_1\]
}
{effective right-side representation}

\DefEq
{
\[a_2'=L(a_1)a_2=a_1a_2\]
}
{effective left-side representation}

\DefEq
{
\symb{\mathcal B_f}{lattice of subrepresentations}1
}
{lattice of subrepresentations}

\DefEq
{
\[R_1:X\rightarrow X'\]
}
{R1 X->}

\DefEq
{
\[
R\circ m=R_1(m)
\]
}
{Rm=R1m}

\DefEq
{
$f(a)(m)\in B_2$
}
{fam in B2}

\ePrints{1502.04063,5114-6019}
\ifx\Semafor\ValueOff
\DefEq
{
$a\in A_1$\Pt
}
{a in A1}
\fi

\DefEq
{
\[
\xymatrix{
&&A_2/s_2\ar[rrrrr]^{q_2}\ar@{}[drrrrr]|{(5)}&&\ar@{}[dddddll]|{(4)}&
\ar@{}[dddddrr]|{(6)}&&t_2A_2\ar[ddddd]^{r_2}\\
&&&&&&&\\
A_1/s_1\ar[r]^{q_1}\ar@/^2pc/@{=>}[urrr]^F&
t_1A_1\ar[d]^{r_1}\ar@{=>}[urrrrr]^(.4)G&&&
A_2/s_2\ar[r]^{q_2}\ar[lluu]_{F(\RedText{p_1(a)})}&
t_2A_2\ar[d]^{r_2}\ar[rruu]_{G(\RedText{t_1(a)})}\\
A_1\ar[r]_{t_1}\ar[u]^{p_1}\ar@{}[ur]|{(1)}\ar@/_2pc/@{=>}[drrr]^f&
B_1\ar@{=>}[drrrrr]^(.4)g&&&
A_2\ar[r]_{t_2}\ar[u]^{p_2}\ar@{}[ur]|{(2)}\ar[ddll]^{f(a)}&
B_2\ar[ddrr]^{g(\RedText{t_1(a)})}\\
&&&&&&&\\
&&A_2\ar[uuuuu]^{p_2}\ar[rrrrr]_{t_2}\ar@{}[urrrrr]|{(3)}&&&&&B_2
}
\]
}
{decompositions of morphism of representations, diagram}

\DefEquation
{
t_i=r_i\circ q_i\circ p_i
}
{morphism of representations of algebra, homomorphism, 1}

\DefEq
{
\[
\RedText{p_1(a)}=a^{\mathrm{ker}\,t_1}
\]
}
{morphism of representations of algebra, p1=}

\DefEq
{
\[
\BlueText{p_2(a)}=a^{\mathrm{ker}\,t_2}
\]
}
{morphism of representations of algebra, p2=}

\DefEquation
{
q_1(\RedText{p_1(a)})=\RedText{t_1(a)}
}
{morphism of representations of algebra, q1=}

\DefEquation
{
q_2(\BlueText{p_2(a)})=\BlueText{t_2(a)}
}
{morphism of representations of algebra, q2=}

\DefEq
{
\[
r_1:\RedText{t_1(a)}\in f(A_1)\rightarrow \RedText{t_1(a)}\in B_1
\]
}
{morphism of representations of algebra, r1=}

\DefEq
{
\[
r_2:\BlueText{t_2(a)}\in f(A_2)\rightarrow \BlueText{t_2(a)}\in B_2
\]
}
{morphism of representations of algebra, r2=}

\DefEquation
{
\begin{pmatrix}
t_1&t_2
\end{pmatrix}
=
\begin{pmatrix}
r_1&r_2
\end{pmatrix}
\circ
\begin{pmatrix}
q_1&q_2
\end{pmatrix}
\circ
\begin{pmatrix}
p_1&p_2
\end{pmatrix}
}
{decompositions of morphism of representations}

\DefEq
{
\[
\xymatrix
{
A/\mathrm{ker}\,f\ar[r]^q&f(A)\ar[d]_r
\\
A\ar[u]^p\ar[r]^f&B
}
\ \ \ f=p\circ q\circ r
\]
}
{decomposition of map f}

\DefEq
{
\[
q:p(a)\in A/\mathrm{ker}\,f\rightarrow f(a)\in f(A)
\]
}
{q:A/ker->f(A)}

\DefEq
{
\[
r:f(a)\in f(A)\rightarrow f(a)\in B
\]
}
{r:f(A)->B}

\DefEq
{
$r_1=q_1\circ p_1$, $r_2=q_2\circ p_2$.
}
{r=q*p}

\DefEq
{
\[
r_2:A_2\rightarrow A_2
\]
}
{R:A2->A2}

\DefEq
{
\symb{J_f}{closure operator, representation}1.
}
{closure operator, representation}

\AddEq{show closure operator, representation}
{
$\ShowSymbol{subrepresentation generated by set}1$
}

\DefEq
{
\symb{J_f(X)}{subrepresentation generated by set}1
}
{subrepresentation generated by set}

\AddEq{generating set of representation}
{
$J_f(X)=A_2$.
}

\DefEq
{
\[
\xymatrix
{
f_{B_2}:A_1\ar[r]|{*}&B_2
}
\]
}
{representation of algebra A in algebra B}

\DefEq
{
$f_{B_2}(a)=f(a)|_{B_2}$.
}
{fB2(a)=}

\AddEq [3]{r12:A->B}
{
\[
\begin{pmatrix}
#1_1:#2_1\rightarrow #3_1&#1_2:#2_2\rightarrow #3_2
\end{pmatrix}
\]
}

\DefEq
{
$\BlueText{r_2(f(a)(m))}$.
}
{r2(f(a,m))}

\DefEq
{
$\BlueText{r_2(m)}\in B_2$
}
{r2(m)in B2}

\DefEq
{
$g(\RedText{r_1(a)})$
}
{g(r1(a))}

\DefEq
{
\[
\RedText{r_1(a)}=r_1(a)
\]
}
{r1(a)=r1(a)}

\DefEq
{
\BlueText{r_2(f(a)(m))}=g(\RedText{r_1(a)})(\BlueText{r_2(m)})
}
{morphism of representations of universal algebra, 2m}

\DefEq
{
\xymatrix{
&A_2\ar[dd]_(.3){f(a)}\ar[rr]^{r_2}&&B_2\ar[dd]^(.3){g(\RedText{r_1(a)})}\\
&\ar @{}[rr]|{(1)}&&\\
&A_2\ar[rr]^{r_2}&&B_2\\
A_1\ar[rr]^{r_1}\ar@{=>}[uur]^(.3)f&&B_1\ar@{=>}[uur]^(.3)g
}
}
{morphism of representations of universal algebra, 2m 1}

\DefEquation
{
\xymatrix{
A_2\ar[rr]^{r_2}&&B_2\\
A_1\ar[rr]^{r_1}\ar[u]^f|{*}&&B_1\ar[u]^g|{*}
}
}
{morphism of representations of universal algebra, definition, 2m 2}

\DefEq
{
\[
a\pC i0xa\pC i1
\]
}
{Sum over repeated index}

\DefEq
{
$\Omega=\{\omega\}$.
}
{Omega=omega}

\DefEq
{
\[ab\omega=ab\]
}
{abo=ab}

\DefEq
{
\[ab\omega=a+b\]
}
{abo=a+b}

\DefEq
{
\[f:A_1\times...\times A_n\rightarrow S\]
}
{polylinear map of algebras}

\DefEq
{
\symb{\mathcal L(D;A_1\times...\times A_n\rightarrow S)}{set polylinear maps}1
}
{set polylinear maps}

\DefEquation
{
(A_1\otimes A_2)\otimes A_3=A_1\otimes(A_2\otimes A_3)
=A_1\otimes A_2\otimes A_3
}
{A1xA2xA3}

\DefEq
{
\[
(v_1, ..., v_n)\in V_1\times...\times V_n
\rightarrow v_1\otimes...\otimes v_n\in V_1\otimes...\otimes V_n
\]
}
{V times->V otimes}

\DefEquation
{
\begin{split}
&\,a_1\otimes...\otimes(a_i+b_i)\otimes...\otimes a_n
\\
=&\,a_1\otimes...\otimes a_i\otimes...\otimes a_n+
a_1\otimes...\otimes b_i\otimes...\otimes a_n
\\
&\,a_i, b_i\in A_i
\end{split}
}
{tensors 1, tensor product}

\DefEquation
{
\begin{split}
&a_1\otimes...\otimes(ca_i)\otimes...\otimes a_n=
c(a_1\otimes...\otimes a_i\otimes...\otimes a_n)
\\
&a_i\in A_i\ \ \ c\in D
\end{split}
}
{tensors 2, tensor product}

\DefEq
{
\[f\circ 0=0\]
}
{linear map, 0, D algebra}

\AddEq [3]{f:A->B}
{
\[#1:#2\rightarrow #3\]
}

\DefEquation
{
f:A^n\rightarrow A,
a=f\circ(a_1,...,a_n)
}
{polylinear map, algebra}

\DefEquation
{
a=f\pC{s}{0}^n\ \sigma_s(I_{(s\cdot 1)}\circ a_1)
\ f\pC{s}{1}^n\ ...\ \sigma_s(I_{(s\cdot n)}\circ a_n)\ f\pC{s}{n}^n
}
{polylinear map, algebra, canonical morphism}

\DefEq
{
$\{a_1,...,a_n\}$
\[
\sigma_s=
\begin{pmatrix}
a_1&...&a_n
\\
\sigma_s(a_1)&...&\sigma_s(a_n)
\end{pmatrix}
\]
}
{transposition of set of variables, algebra}

\DefEq
{
$I_{(s\cdot 1)}$, ..., $I_{(s\cdot n)}\in\mathcal L(D;A\rightarrow A)$
}
{I1n in L(A;A)}

\DefEq
{
\labelItem{monomial of power 0}
$p_0(x)=a_0$, $a_0\in A$.
}
{p0(x)=a0}

\DefEq
{
\labelItem{monomial of power k}
\[
p_k(x)=p_{k-1}(x)(F\circ x)a_k
\]
}
{monomial of power k}

\AddEq{p1(x)=aFxa}
{
$p_1(x)=a_0(F\circ x)a_1$.
}

\AddEq{Basis F=...}
{
$\Basis F=(F_1,...,F_n)$.
}

\AddEq{Fi=Fj}
{
$F_{(i)}=F_{(j)}$
}

\AddEq{F=...}
{
$F=(F_{(1)},...,F_{(k)})$
}

\DefEq
{
\symb{A_k[x]}{module of homogeneous polynomials over algebra}1
}
{module of homogeneous polynomials over algebra}

\DefEq
{
$a\in A^{\otimes (n+1)}$\Pt
}
{a in Aoxn+1}

\DefEq
{
$x_1=...=x_n=x$,
}
{x1=xn=x}

\DefEq
{
\[
a\circ F\circ x^n=a\circ(F_{(1)},...,F_{(n)})\circ(x_1\otimes...\otimes x_n)
\]
}
{a xn=}

\AddEq{F1n in F}
{
$F_{(1)}$, ..., $F_{(n)}\in\Basis F$.
}

\AddEq{F[k]=...}
{
\[
F_{[k]}=(F_{[k](1)},...,F_{[k](n)})
\]
}

\DefEq
{
\[
a=a_{i\cdot 0}\otimes a_{i\cdot 1}\otimes...\otimes a_{i\cdot n}
\ \ \ i\in I
\]
}
{a=oxi}

\DefEq
{
$\sigma=\{\sigma_i\in S(n):i\in I\}$
}
{si in Sn}

\DefEq
{
\[
(a,\sigma):A^{\times n}\rightarrow A
\]
}
{ox:An->A}

\DefEq
{
\begin{align*}
(a,\sigma)\circ (b_1,...,b_n)&=
(a_{i\cdot 0}\otimes a_{i\cdot 1}\otimes...\otimes a_{i\cdot n},\sigma_i)\circ (b_1,...,b_n)
\\&=
a_{i\cdot 0}\sigma_i(b_1)a_{i\cdot 1}...\sigma_i(b_n)a_{i\cdot n}
\end{align*}
}
{ox circ =}

\AddEq{p(x)=a circ xk}
{
\[
\begin{matrix}
p(x)=a_{[s]}\circ F_{[s]}\circ x^k
&a_{[s]}\in A^{\otimes(k+1)}
\end{matrix}
\]
}

\DefEq
{
\symb{A[x]}{algebra of polynomials over algebra}{}
\[
\ShowSymbol{algebra of polynomials over algebra}{}
=\bigoplus_{n=0}^{\infty}A_n[x]
\]
}
{algebra of polynomials over algebra}

\DefEquation
{
}
{derivative in L = +o}

\AddEq [2]{a1n}
{
$a_1$, ..., $a_{#1}$#2
}

\DefEq
{
$\epsilon\in R$, $\epsilon>0$,
}
{epsilon in R}

\DefEq
{
$a_0\in U$.
}
{a0 in U}

\DefEq
{
\[
f:U\rightarrow \mathcal L(D;A^p\rightarrow B)
\]
}
{U->L(Ap,B)}

\DefEq
{
$\partial_x f(x)\in CE$.
}
{dx fx in CE}

\DefEquation
{
\begin{array}{r@{\ }l}
&((a_0,...,a_n,\sigma)\circ(f_1,...,f_n))\circ(x_1,...,x_n)
\\=&
(a_0\sigma(f_1)a_1...a_{n-1}\sigma(f_n)a_n)\circ(x_1,...,x_n)
\\=&
a_0\sigma(f_1\circ x_1)a_1...a_{n-1}\sigma(f_n\circ x_n)a_n
\end{array}
}
{n linear map A LA}

\DefEq
{
\symb{\mathcal L(D;A^n\rightarrow S)}{set polylinear maps An}1
}
{set polylinear maps An}

\DefEq
{
\begin{align*}
f\circ(
a_1, ...,
a_i+ b_i, ...,
a_n)
&=
f\circ(
a_1, ...,
a_i, ...,
a_n)
+
f\circ(
a_1, ...,
b_i, ...,
a_n)
\\
f\circ(
a_1, ...,
pa_i, ...,
a_n)
&=
pf\circ(
a_1, ...,
a_i, ...,
a_n)
\end{align*}
\[
\begin{matrix}
1\le i\le n
&
a_i, b_i \in A_i
&
p\in D
\end{matrix}
\]
}
{polylinear map of algebras, 1}

\DefEquation
{
r_2\circ\BlueText{f(a)}=g(\RedText{r_1(a)})\circ r_2
}
{morphism of representations of universal algebra, definition, 2}

\DefEq
{
\symb{A\cong B}{isomorphic}1.
}
{isomorphic}

\DefEquation
{
f(a_1)...f(a_n)\omega=f(a_1...a_n\omega)
}
{afo=aof}

\DefEquation
{
\begin{array}{r@{\,}l@{\ \ \ }r@{\,}l@{\ \ \ }r@{\,}l}
a_k&=a_{k\cdot 0}\underline{\otimes}(1\otimes a_{k\cdot 1})
&a_{k\cdot 0}&\in A^{\otimes k}&a_{k\cdot 1}&\in A
\end{array}
}
{p(x)=a k-1 k 3}

\DefEq
{
\[M=\max(M_1,...,M_n)\]
}
{M M1 Mn}

\DefEq
{
\[N=\max(N_1,N_2)\]
}
{N N1 N2}

\DefEq
{
f_i(x)=\lim_{m\rightarrow\infty}f_{i\cdot m}(x)
}
{fi(x)=lim}

\DefEq
{
\delta_1\in R,\ \delta_1>0
}
{delta1 in R}

\DefEq
{
\(m>M_i\)
}
{m>Mi}

\DefEq
{
\[h_m=f_{1\cdot m}...f_{n\cdot m}\omega\]
}
{hm=f1m...fnm omega}

\DefEq
{
\[
\epsilon_1(\delta)<\epsilon
\]
}
{e(d)<e}

\DefEq
{
\[\delta_1=0\Rightarrow\epsilon_1=0\]
}
{d1=0=>e1=0}

\DefEq
{
F_i\ge 0
}
{F>0}

\DefEq
{
F_i=\sup\|f_i(x)\|
}
{F=sup|f|}

%% file: Preliminary.Ref.tex

\ifx\UseRussian\Defined
\def\TheoremFollows{Теорема является следствием теоремы\ }
\else
\def\TheoremFollows{The theorem follows from the theorem\ }
\fi

\ePrints{4975-6381}%
\ifx\Semafor\ValueOn%
\def\RefLinearMap{5114-6019}%
\else
\ePrints{1601.03259,1506.00061,MQuater,1610.309618526,CACAA.06}%
\ifx\Semafor\ValueOn%
\def\RefLinearMap{1502.04063}%
\else
\def\RefLinearMap{MVector}%
\fi
\fi

\ePrints{1610.309618526}
\ifx\Semafor\ValueOn
\def\RefCalculus{1601.03259}%
\else
\def\RefCalculus{}%
\fi

\AddEq{SetupRefRepresentation}
{
\ifx\Preliminary\ValueOn%
\ePrints{1506.00061,1601.03259,MQuater}%
\ifx\Semafor\ValueOn%
\def\RefRepresentation{1502.04063}%
\fi
\ePrints{5148-4632}
\ifx\Semafor\ValueOn
\def\RefRepresentation{BRepro}%
\fi
\ePrints{4975-6381}
\ifx\Semafor\ValueOn
\def\RefRepresentation{5114-6019}%
\fi
\else
\def\RefRepresentation{}%
\ePrints{1102.5168,5114-6019}
\ifx\Semafor\ValueOn
\def\RefRepresentation{0912.3315}%
\fi
\ePrints{MVector}
\ifx\Semafor\ValueOn
\def\RefRepresentation{0912.3315}%
\fi
\fi
}

\AddEq{SetupRefPolymorphism}
{
\ifx\Preliminary\ValueOn%
\ePrints{5148-4632}%
\ifx\Semafor\ValueOn%
\def\RefPolymorphism{5114-6019}%
\fi
\else
\def\RefPolymorphism{}%
\ePrints{MVector}%
\ifx\Semafor\ValueOn
\def\RefPolymorphism{1502.04063}%
\fi
\fi
}

\AddEq{SetupRefOmegaNorm}
{
\ePrints{1310.5591}%
\ifx\Semafor\ValueOff%
\ifx\Preliminary\ValueOn%
\def\RefTheoremOmegaNorm{1305.4547}
\ePrints{5410-9916,4975-6381}%
\ifx\Semafor\ValueOn%
\def\RefTheoremOmegaNorm{5059-9176}
\fi
\else
\def\RefTheoremOmegaNorm{}
\ePrints{CACAA.06}%
\ifx\Semafor\ValueOn%
\def\RefTheoremOmegaNorm{CACAA.04.001}
\fi
\fi
\else%
\def\RefTheoremOmegaNorm{1305.4547}
\fi%
}

\AddEq{SetupRefMeasure}
{
\ifx\Preliminary\ValueOn%
\ePrints{1601.03259,1610.309618526,1601.03259}
\ifx\Semafor\ValueOn
\def\RefMeasure{1310.5591}
\else
\ePrints{4975-6381}
\ifx\Semafor\ValueOn
\def\RefMeasure{5410-9916}
\else
\def\RefMeasure{CACAA.04.001}
\fi
\fi
\else
\def\RefMeasure{}
\fi
}

\AddEq{SetupRef}
{
\ShowEq{SetupRefRepresentation}
\ShowEq{SetupRefOmegaNorm}
\ShowEq{SetupRefMeasure}
\ShowEq{SetupRefPolymorphism}
}

\ShowEq{SetupRef}

\newcommand\ProofTheorem[2]
{
\begin{proof}
\TheoremFollows
\RefTheorem[#1]{#2}.
\end{proof}
}%


\DefEq
{
\ePrints{1310.5591,5059-9176}%
\ifx\Semafor\ValueOff%
\ePrints{5410-9916}%
\ifx\Semafor\ValueOn%
\RefDefinition{representation of algebra}\Pt
\else
\RefDefinition{left-side representation of algebra}\Pt
\fi
\else
\RefDefinition[0912.3315]{left-side representation of algebra}\Pt
\fi
}
{ref definition: left-side representation of algebra}

\DefEq
{
\ePrints{1310.5591}
\ifx\Semafor\ValueOff
\RefDefinition{morphism of representations of universal algebra}.
\else
\RefDefinition[0912.3315]{morphism of representations of universal algebra}.
\fi
}
{ref definition: morphism of representations of universal algebra}

\DefEq
{
\ePrints{1310.5591}%
\ifx\Semafor\ValueOff%
\RefDefinition{closed ball}
\else%
\RefDefinition[1305.4547]{closed ball}
\fi%
}
{ref definition: closed ball}

\DefEq
{
\ePrints{1310.5591}%
\ifx\Semafor\ValueOff%
\RefDefinition{open ball}\Pt
\else%
\RefDefinition[1305.4547]{open ball}\Pt
\fi%
}
{ref definition: open ball}

\DefEq
{
\ePrints{1310.5591}%
\ifx\Semafor\ValueOff%
\RefDefinition{open set},
\else%
\RefDefinition[1305.4547]{open set},
\fi%
}
{ref definition: open set}

\DefEq
{
\ePrints{1310.5591}%
\ifx\Semafor\ValueOff%
\RefDefinition{fundamental sequence}\Pt
\else%
\RefDefinition[1305.4547]{fundamental sequence}\Pt
\fi%
}
{ref definition: fundamental sequence}

\DefEq
{
\ePrints{1310.5591}%
\ifx\Semafor\ValueOff%
\RefDefinition{sequence converges uniformly}\Pt
\else%
\RefDefinition[1305.4547]{sequence converges uniformly}\Pt
\fi%
}
{ref definition: sequence converges uniformly}

\DefEq
{
\ePrints{1310.5591}%
\ifx\Semafor\ValueOff%
\RefDefinition{compact set},
\else%
\RefDefinition[1305.4547]{compact set},
\fi%
}
{ref definition: compact set}

\DefEq
{
\ePrints{1310.5591}%
\ifx\Semafor\ValueOff%
\RefDefinition{Omega group},
\else%
\RefDefinition[1305.4547]{Omega group},
\fi%
}
{ref definition: Omega group}

\DefEq
{
\ePrints{1310.5591}%
\ifx\Semafor\ValueOff%
\RefItem{|a|>=0}
\else%
\RefItem[1305.4547]{|a|>=0}
\fi%
}
{ref item |a|>=0}

\DefEq
{
\ePrints{1310.5591}%
\ifx\Semafor\ValueOff%
\RefItem{|a|=0}.
\else%
\RefItem[1305.4547]{|a|=0}.
\fi%
}
{ref item |a|=0}

\DefEq
{
\ePrints{1310.5591}%
\ifx\Semafor\ValueOff%
\RefItem{|a+b|<=|a|+|b|}\Pt
\else%
\RefItem[1305.4547]{|a+b|<=|a|+|b|}\Pt
\fi%
}
{ref item |a+b|<=|a|+|b|}

\DefEq
{
\ePrints{1310.5591}%
\ifx\Semafor\ValueOff%
\EqRef{|a omega|<|omega||a|1n}.
\else%
\EqRef[1305.4547]{|a omega|<|omega||a|1n}.
\fi%
}
{ref EqRef |a omega|<|omega||a|1n}

\DefEq
{
\ePrints{1310.5591}%
\ifx\Semafor\ValueOff%
\EqRef{|a-b|>|a|-|b|},
\else%
\EqRef[1305.4547]{|a-b|>|a|-|b|},
\fi%
}
{ref EqRef |a-b|>|a|-|b|}

\DefEq
{
\ePrints{1310.5591}%
\ifx\Semafor\ValueOff%
\EqRef{|fab|<|f||a||b|}.
\else%
\EqRef[1305.4547]{|fab|<|f||a||b|}.
\fi%
}
{ref |fab|<|f||a||b|}

\AddEq{ref tensor product and polylinear map}
{
\RefTheorem[\RefLinearMap]{there exists tensor product of modules},
\RefTheorem[\RefLinearMap]{tensor product and polylinear map}.
}

\DefEq
{
\ePrints{1502.04063,5114-6019}
\ifx\Semafor\ValueOff
\RefDefinition{reduced morphism of representations},
\else
\Ref[0912.3315]{remark: reduced morphism of representations},
\fi
}
{ref reduced morphism of representations}

\DefEq
{
\RefTheorem[1202.6021]{expand linear mapping, quaternion}.
}
{ref expand linear mapping, quaternion}

\DefEq
{
\RefTheorem[0912.4061]{ax-xa=1 quaternion algebra}.
}
{ref ax-xa=1 quaternion algebra}

\DefEq
{
\RefTheorem[1202.6021]{aE, quaternion, Jacobian matrix}.
}
{ref aE, quaternion, Jacobian matrix}

\DefEq
{
\RefTheorem[1003.1544]{Quaternion over real field}.
}
{ref Quaternion over real field}

\DefEq
{
\RefTheorem[1003.1544]{quaternion conjugation}.
}
{ref quaternion conjugation}

\DefEq
{
\RefTheorem[1003.1544]{Quaternion over real field, matrix}.
}
{ref Quaternion over real field, matrix}

\DefEq
{
\RefTheorem[1601.03259]{representation of derivative, algebra A->B}.
}
{ref differentiable map A->B}

\DefEq
{
\RefTheorem[1601.03259]{map is continuous, derivative}.
}
{ref map is continuous, derivative}

\DefEq
{
\RefTheorem[1302.7204]{monomial of power k}.
}
{ref monomial of power k}

\DefEq
{
\RefTheorem[1302.7204]{map A(k+1)->pk otimes is polylinear map}.
}
{ref map A(k+1)->pk otimes is polylinear map}

\DefEq
{
\RefTheorem[1302.7204]{r=+q circ p}.
}
{ref r=+q circ p}

\DefEq
{
\RefTheorem[1302.7204]{r=+q circ p 1}.
}
{ref r=+q circ p 1}

\DefEq
{
\RefTheorem[1105.4307]{Re Im A->F}.
}
{ref Re Im A->F}

\DefEq
{
\RefTheorem[1105.4307]{conjugation in algebra}.
}
{ref conjugation in algebra}

\DefEq
{
\RefTheorem[1105.4307]{structural constant of algebra with unit}.
}
{ref structural constant of algebra with unit}

\DefEq
{
\ePrints{MQuater}%
\ifx\Semafor\ValueOn%
\RefTheorem[1601.03259]{derivative, representation in algebra}.
\else%
\EqRef{derivative, linear path}
\fi%
}
{ref derivative, representation in algebra}

%% file: Preliminary.Parm.tex

\AddEq{=D1D2}
{
\def\DFDT{D1 D2 }%
\def\DF{1}%
\def\DT{2}%
\def\VF{1}%
\def\VT{2}%
\def\MapE{hf}%
}

\AddEq{=DD}
{
\def\DFDT{D }%
\def\DF{}%
\def\DT{}%
\def\VF{1}%
\def\VT{2}%
\def\MapE{f}%
}

\AddEq{=left}
{
\def\SideWS{left }%
\def\SideNS{left}%
\def\HSide{\Hyph side }%
\def\SideA{A*}%
\def\CBase{D}%
\def\Base{A}%
\def\Module{V}%
\def\BaseRings{of commutative ring $D$ and $D$\Hyph algebra $A$ }
\def\Algebra{associative division $D$\Hyph algebra}%
\def\algebra{associative $D$\Hyph algebra}%
\def\algebraa{$D$\Hyph algebra }%
\def\ATransf{g_{3,4}}%
\def\DTransf{g_{1,4}}%
\def\DArg{d}%
\ifx\UseRussian\Defined%
\def\HSide{востороннее }%
\def\BaseRings{коммутативного кольца $D$ и $D$\Hyph алгебры $A$ }
\def\Algebra{ассоциативная $D$\Hyph алгебра с делением}%
\def\algebra{ассоциативная $D$\Hyph алгебра}%
\def\algebraa{$D$\Hyph алгебре }%
\fi%
}

\AddEq{=right}
{
\def\SideWS{right }%
\def\SideNS{right}%
\def\HSide{\Hyph side }%
\def\SideA{*A}%
\def\CBase{D}%
\def\Base{A}%
\def\Module{V}%
\def\BaseRings{of commutative ring $D$ and $D$\Hyph algebra $A$ }
\def\Algebra{associative division $D$\Hyph algebra}%
\def\algebra{associative $D$\Hyph algebra}%
\def\algebraa{$D$\Hyph algebra }%
\def\ATransf{g_{3,4}}%
\def\DTransf{g_{1,4}}%
\def\DArg{d}%
\ifx\UseRussian\Defined%
\def\HSide{востороннее }%
\def\BaseRings{коммутативного кольца $D$ и $D$\Hyph алгебры $A$ }
\def\Algebra{ассоциативная $D$\Hyph алгебра с делением}%
\def\algebra{ассоциативная $D$\Hyph алгебра}%
\def\algebraa{$D$\Hyph алгебре }%
\fi%
}

\AddEq{=module}
{
\def\SideWS{}%
\def\SideNS{}%
\def\HSide{}%
\def\SideA{D}%
\def\CBase{Z}%
\def\Base{D}%
\def\Module{A}%
\def\BaseRings{of ring of rational integers $Z$ and commutative ring $D$ }
\def\Algebra{field}%
\def\algebra{ring}%
\def\algebraa{ring }%
\def\Algebrab{module}%
\def\Algebrac{module}%
\def\ATransf{g_2}%
\def\DTransf{g_1}%
\def\DArg{n}%
\ifx\UseRussian\Defined%
\def\BaseRings{кольца целых чисел $Z$ и коммутативного кольца $D$ }
\def\Algebra{поле}%
\def\algebra{кольцо}%
\def\algebraa{кольце }%
\def\Algebrab{модуль}%
\fi%
}

\AddEq{=Omega}
{
\def\SideA{\Omega}%
\def\Algebrab{group}%
\def\Algebrac{Omega group}%
\def\Module{A}%
\ifx\UseRussian\Defined%
\def\Algebrab{группа}%
\fi%
}

\DefEq%
{%
\def\Pt{.}%
}%
{=.}%

\DefEq%
{%
\def\Pt{;}%
}%
{=.c}%

\DefEq%
{%
\def\Pt{,}%
}%
{=c}%

\DefEq%
{%
\def\Pt{}%
}%
{=z}%

\DefEq
{
\def\PA{}%
}
{D= F= A=}

\DefEq%
{%
\def\b{a}%
\def\B{A}%
}%
{B=A}%

\DefEq%
{%
\def\b{b}%
\def\B{B}%
}%
{B=B}%

\DefEq
{
\def\PX{X}
}
{X}

\DefEq%
{%
\def\Pf{f}%
}%
{f=f}

\DefEq%
{%
\def\Pf{g}%
}%
{f=g}%

\DefEq%
{%
\def\Pf{m}%
}%
{f=m}%

\DefEq
{
\def\Pf{I}%
}
{f=I}

\DefEq
{
\def\Pf{P^{\gi 0\cdot}_{}{}}%
}
{f=P0}

\DefEq
{
\def\Pf{P^{\gi 1\cdot}_{}{}}%
}
{f=P1}

\DefEq
{
\def\Pf{P^{\gi 2\cdot}_{}{}}%
}
{f=P2}

\DefEq
{
\def\Pf{P^{\gi 3\cdot}_{}{}}%
}
{f=P3}

\DefEq
{
\def\pD{D}%
\def\pA{A}%
\def\pB{A}%
}
{A=A}

\DefEq%
{%
\def\pD{D}%
\def\pA{A}%
\def\pB{B}%
}%
{A=AB}

\DefEq
{
\def\pD{D}%
\def\pA{B}%
\def\pB{C}%
}
{A=BC}

\DefEq%
{%
\def\pD{R}%
\def\pA{C}%
\def\pB{C}%
}%
{A=C}%

\DefEq%
{%
\def\pD{C}%
\def\pA{C}%
\def\pB{C}%
}%
{A=CC}%

\DefEq
{
\def\pD{D}%
\def\pA{A_1}%
\def\pB{A_2}%
}
{A=12}

\DefEq
{
\def\pD{D}%
\def\pA{A_2}%
\def\pB{A_2}%
}
{A=2}

\DefEq
{
\def\pD{D}%
\def\pA{A_1\rightarrow A_2}%
\def\pB{A_3}%
}
{A=123}

\DefEq
{
\def\pD{D}%
\def\pA{A_1}%
\def\pB{A_3}%
}
{A=13}

\DefEq
{
\def\pD{D}%
\def\pA{A}%
\def\pB{C}%
}
{A=AC}

\DefEq
{%
\def\Pn{0}%
}%
{n=0}

\DefEq
{%
\def\Pn{n}%
}%
{n=n}

\DefEq
{%
\def\Pn{k}%
}%
{n=k}

\DefEq
{%
\def\Pn{1}%
\def\Pp{p}%
}%
{n=1}

\DefEq
{%
\def\Pn{2}%
\def\Pp{q}%
}%
{n=2}

\DefEq
{%
\def\Pn{3}%
\def\Pp{r}%
}%
{n=3}

\DefEq
{%
}%
{m=1}

\DefEq
{%
}%
{m=2}

%% file: Stmt.Representation.English.tex
\input{Stmt.Representation.Eq}

\DefDefinition{morphism of representations of universal algebra}
{
Let
\ShowEq{f:A->*B}f{A_1}{A_2}
be representation of $\Omega_1$\Hyph algebra $A_1$
in $\Omega_2$\Hyph algebra $A_2$ and
\ShowEq{f:A->*B}g{B_1}{B_2}
be representation of $\Omega_1$\Hyph algebra $B_1$
in $\Omega_2$\Hyph algebra $B_2$.
For
\EqParm{i=1,2}{=c}
let the map
\ShowEq{ri:A->B}
be homomorphism of $\Omega_i$\Hyph algebra.
The matrix of maps
\ShowEq{map r,R}{r_1}{r_2}{}
such, that
\ShowEq{morphism of representations of universal algebra, definition, 2}
is called
\AddIndex{morphism of representations from $f$ into $g$}
{morphism of representations from f into g}.
We also say that
\AddIndex{morphism of representations of $\Omega_1$\Hyph algebra
in $\Omega_2$\Hyph algebra}
{morphism of representations of Omega1 algebra in Omega2 algebra} is defined.
}

\DefTheorem{reduced polymorphism of representations}
{
Let the map $r_2$ be reduced polymorphism of
effective representations $f_1$, ..., $f_n$ into effective representation $f$.
\begin{itemize}
\item
For any
\ShowEq{k,1n}
the map $r_2$ satisfies to the equality
\ShowEq{reduced polymorphism of representation, ak}
\item
For any
\ShowEq{kl,1n}
the map $r_2$ satisfies to the equality
\ShowEq{reduced polymorphism of representation, akl}
\item
Let $\omega_2\in\Omega_2(p)$.
For any
\ShowEq{k,1n}
the map $r_2$ satisfies to the equality
\ShowEq{reduced polymorphism of representation, omega2}
\end{itemize}
}

\DefEq
{
\begin{remark}
\labelRemark{notation for effective \SideNS-side representation of algebra}
If the \SideNS\Hyph side representation of multiplicative $\Omega$\Hyph group is effective, then we identify
an element of multiplicative $\Omega$\Hyph group and its image and write \SideNS\Hyph side transformation
caused by $A_1$\Hyph number $a_1$
as
\ShowEq{effective \SideNS-side representation}
The map
\ShowEq{A12->A2 \SideNS}
generated by \SideNS\Hyph side effective representation $f$ is called
\AddIndex{\SideNS\Hyph side product}{\SideNS-side product}
of $A_2$\Hyph number over $A_1$\Hyph number.
If representations
\ShowEq{f:A->*B}fAM
\ShowEq{f:A->*B}gBN
are effective, then we can write the equality
\eqRef{morphism of representations of universal algebra, 2m}{definition}
in the following form
\ShowEq{morphism of \SideWS representation of Omega group}
Orbit of representation $f$ has form
\ShowEq{orbit of effective \SideNS-side representation}
We will use notation
\ShowEq{space of orbits of effective \SideNS-side representation}
for the space of orbits of \SideNS\Hyph side effective $A_1$\Hyph representation.
\qed
\end{remark}
}
{remark: notation for effective representation}

\DefDefinition{product in category}
{
Let $\mathcal A$ be a category.
Let
\ShowEq{set Bi}B
be the set of objects of $\mathcal A$.
Object
\ShowEq{product in category}
and set of morphisms
\ShowEq{set f:A->B}fP{B_i}
is called a
\AddIndex{product of set of objects
\ShowEq{set Bi}B
in category $\mathcal A$}
{product in category}\,\footnote{
I made definition according to \citeBib{Serge Lang}, page 58.}
if for any object $R$ and set of morphisms
\ShowEq{set f:A->B}gR{B_i}
there exists a unique morphism
\ShowEq{f:A->B}hRP
such that diagram
\ShowEq{product in category diagram}
is commutative for all $i\in I$.

If $|I|=n$, then we also will use notation
\ShowEq{product in category, 1 n}
for product of set of objects
$\{B_i,\Ii\}$ in $\mathcal A$.
}

\DefDefinition{coproduct in category}
{
Let $\mathcal A$ be a category.
Let
\ShowEq{set Bi}B
be the set of objects of $\mathcal A$.
Object
\ShowEq{coproduct in category}
and set of morphisms
\ShowEq{set f:A->B}f{B_i}P
is called a
\AddIndex{coproduct of set of objects
\ShowEq{set Bi}B
in category $\mathcal A$}
{coproduct in category}\,\footnote{
I made definition according to \citeBib{Serge Lang}, page 59.}
if for any object $R$ and set of morphisms
\ShowEq{set f:A->B}g{B_i}R
there exists a unique morphism
\ShowEq{f:A->B}hPR
such that diagram
\ShowEq{coproduct in category diagram}
is commutative for all $i\in I$.

If $|I|=n$, then we also will use notation
\ShowEq{coproduct in category, 1 n}
for coproduct of set of objects
\ShowEq{set Bi}B
in $\mathcal A$.
}

\DefTheorem{structure of subrepresentations}
{
Let\,\footnote{The statement of theorem is similar to the
statement of theorem 5.1, \citeBib{Cohn: Universal Algebra}, page 79.}
\ShowEq{f:A->*B}g{A_1}{A_2}
be representation of $\Omega_1$\Hyph algebra $A_1$
in $\Omega_2$\Hyph algebra $A_2$.
Let $X\subset A_2$.
Define a subset $X_k\subset A_2$ by induction on $k$.
\ShowEq{structure of subrepresentations}
Then
\ShowEq{structure of subrepresentations, 1}
}

\DefDefinition{direct sum of Abelian groups}
{
Coproduct in category of Abelian groups $Ab$ is called
\AddIndex{direct sum}{direct sum}.\,\footnote{
See also definition in \citeBib{Serge Lang}, pages 36, 37.}
We will use notation
\ShowEq{direct sum of Abelian groups}
for direct sum of Abelian groups $A$ and $B$.
}

\DefTheorem{direct sum of Abelian groups}
{
Let
\ShowEq{set Bi}A
be set of Abelian groups.
Let
\ShowEq{A in xAi}
be such set that
\ShowEq{(xi)in A},
if
$x_i\ne 0$
for finite number of indices $i$.
Then\,\footnote{
See also proposition \citeBib{Serge Lang}-7.1, page 37.}
\ShowEq{A=o+Ai}
}

\DefProof{direct sum of Abelian groups}
{
According to construction,
$A$ is subgroup of Abelian group $\prod A_i$.
The map
\ShowEq{f:A->B}{\lambda_j}{A_j}A
defined by the equality
\ShowEq{lj(x)=0x0}
is an injective homomorphism.

Let
\ShowEq{set f:A->B}f{A_i}B
be set of homomorphisms into Abelian group $B$.
We define the map
\ShowEq{f:A->B}fAB
by the equality
\ShowEq{fxi,i=sum fxi}
The sum in the right side of the equality
\EqRef{fxi,i=sum fxi}
is finite, since all summands, except for a finite number, equal $0$.
From the equality
\ShowEq{fi(x+y)=}
and the equality
\EqRef{fxi,i=sum fxi},
it follows that
\ShowEq{f(x+y)i=}
Therefore, the map $f$ is homomorphism of Abelian group.
The equality
\ShowEq{flj=fj}
follows from equalities
\EqRef{lj(x)=0x0},
\EqRef{fxi,i=sum fxi}.
Since the map $\lambda_i$ is injective,
then the map $f$ is unique.
Therefore, the theorem folows from definitions
\RefDefinition{coproduct in category},
\RefDefinition{direct sum of Abelian groups}.
}

\DefEq
{
Let an endomorphism $f$ act on $A_2$\Hyph number $a$ on the \SideNS
\ShowEq{f(a) \SideNS}
Let
\ShowEq{End O2}{A_2}{}
be a multiplicative $\Omega$\Hyph group with product
\ShowEq{fga= \SideNS}
and $\delta$ be unit of group
\ShowEq{End O2}{A_2}.
We consider the equality
\EqRef{fga= \SideNS}
as
\AddIndex{associative law}{associative law}.
This allows writing of equality
\EqRef{fga= \SideNS}
without using of brackets
\ShowEq{fga= 1 \SideNS}

\begin{definition}
\labelDefinition{\SideNS-side representation of group}
Let $A_1$ be multiplicative $\Omega$\Hyph group.
We call a homomorphism of multiplicative $\Omega$\Hyph group
\DrawEq{representation of group, map}{\SideNS-side}
\AddIndex{\SideNS\Hyph side representation}{\SideNS-side representation}
of multiplicative $\Omega$\Hyph group $A_1$
or
\AddIndex{\SideNS\Hyph side $A_1$\Hyph representation}{\SideNS-side A representation}
in $\Omega_2$\Hyph algebra $A_2$ if the map $f$ holds
\DrawEq{\SideNS-side representation of group}{def}
\qed
\end{definition}
}
{definition: side representation of group}

\DefDefinition{reduced morphism of representations}
{
Let
\ShowEq{f:A->*B}f{A_1}{A_2}
be representation of $\Omega_1$\Hyph algebra $A_1$
in $\Omega_2$\Hyph algebra $A_2$ and
\ShowEq{f:A->*B}g{A_1}{B_2}
be representation of $\Omega_1$\Hyph algebra $A_1$
in $\Omega_2$\Hyph algebra $B_2$.
Let
\ShowEq{id:A1->A1 A2->B2}
be morphism of representations.
In this case we identify morphism
\ShowEq{map r,R}{\id}{r_2}{}
of representations of $\Omega_1$\Hyph algebra and corresponding homomorphism $r_2$ of $\Omega_2$\Hyph algebra
and the homomorphism $r_2$ is called
\AddIndex{reduced morphism of representations}
{reduced morphism of representations}.
We will use diagram
\ShowEq{morphism id,R of representations}
to represent reduced morphism $r_2$ of representations of
$\Omega_1$\Hyph algebra.
From diagram it follows
\ShowEq{morphism of representations of universal algebra}
We also use diagram
\ShowEq{morphism id,R of representations 2}
instead of diagram
\EqRef{morphism id,R of representations}.
}

\DefEq
{
\begin{definition}
\labelDefinition{tensor product of representations}
Let $A$ be Abelian multiplicative $\Omega_1$\Hyph group.
Let
\EqParm{A1n}{=z}
be $\Omega_2$\Hyph algebras.\,\footnote{
I give definition
of tensor product of representations of universal algebra
following to definition in \citeBib{Serge Lang}, p. 601 - 603.
}
Let, for any
\ShowEq{k,1n}
\ShowEq{representation A Ak 1}
be effective representation of multiplicative $\Omega_1$\Hyph group $A$
in $\Omega_2$\Hyph algebra $A_k$.
Consider category $\mathcal A$ whose objects are
reduced polymorphisms of representations $f_1$, ..., $f_n$
\ShowEq{polymorphisms category}
where $S_1$, $S_2$ are $\Omega_2$\Hyph algebras and
\ShowEq{representation of algebra in S1 S2}
are effective representations of multiplicative $\Omega_1$\Hyph group $A$.
We define morphism
\ShowEq{r1->r2}
to be reduced morphism of representations $h:S_1\rightarrow S_2$
making following diagram commutative
\ShowEq{polymorphisms category, diagram}
Universal object
\ShowEq{tensor product of representations}
of category $\mathcal A$ is called
\AddIndex{tensor product}{tensor product}
of representations
\EqParm{A1n}{=.}
\qed
\end{definition}
}
{definition: tensor product of representations}

\DefEq
{
\begin{theorem}
\labelTheorem{tensor product of representations}
Since there exists polymorphism of representations,
then there exists tensor product of representations.
\end{theorem}
}
{theorem: tensor product of representations}

\DefEq
{
\begin{theorem}
\labelTheorem{representation, tensor product}
Let
\ShowEq{equivalence, 1, representation, tensor product}
Tensor product is distributive over operation $\omega$
\ShowEq{tensors 1, representation, tensor product}
The representation of multiplicative $\Omega_1$\Hyph group $A$
in tensor product is defined by equality
\ShowEq{tensors 2, representation, tensor product}
\end{theorem}
}
{theorem: representation, tensor product}

\DefEq
{
\begin{theorem}
\labelTheorem{tensor product and polymorphism}
Let $B_1$, ..., $B_n$ be
$\Omega_2$\Hyph algebras.
Let
\ShowEq{map f, 1, representation, tensor product}
be reduced polymorphism defined by equality
\ShowEq{map f, representation, tensor product}
Let
\ShowEq{map g, representation, tensor product}
be reduced polymorphism into $\Omega$\Hyph algebra $V$.
There exists morphism of representations
\ShowEq{map h, representation, tensor product}
such that the diagram
\ShowEq{map gh, representation, tensor product}
is commutative.
\end{theorem}
}
{theorem: tensor product and polymorphism}

\DefEq
{
\begin{theorem}
\labelTheorem{B times->B otimes}
The map
\ShowEq{B times->B otimes}
is polymorphism.
\end{theorem}
}
{theorem: B times->B otimes}

\DefTheorem{|a-b|>|a|-|b|}
{
Let $A$ be normed $\Omega$\Hyph group.
Then
\ShowEq{|a-b|>|a|-|b|}
}

\DefDefinition{continuous map, Omega group}
{
A map
\ShowEq{f:A->B}f{A_1}{A_2}
of normed $\Omega_1$\Hyph group $A_1$ with norm $\|x\|_1$
into normed $\Omega_2$\Hyph group $A_2$ with norm $\|y\|_2$
is called \AddIndex{continuous}{continuous map}, if
for every as small as we please $\epsilon>0$
there exist such $\delta>0$, that
\ShowEq{|x'-x|<delta}
implies
\ShowEq{|f(x)-f(x')|<epsilon}
}

\DefTheorem{continuous map, open set}
{
A map
\ShowEq{f:A->B}f{A_1}{A_2}
of normed $\Omega_1$\Hyph group $A_1$ with norm $\|x\|_1$
into normed $\Omega_2$\Hyph group $A_2$ with norm $\|y\|_2$
is continuous, iff
preimage of an open set
is the open set.
}

\DefTheorem{image of interval is interval, real field}
{
Let
\ShowEq{f:A->B}fRR
be continuous map of real field.
Then image of interval is interval.
}

\DefEq
{
\begin{definition}
\labelDefinition{norm of operation}
Let $A$
be normed $\Omega$\Hyph group.
For $n$\Hyph ary operation $\omega$, the value
\ShowEq{norm of operation}
\ShowEq{norm of operation, definition}
is called
\AddIndex{norm of operation}{norm of operation} $\omega$.
\qed
\end{definition}
}
{definition: norm of operation}

\DefTheorem{|fx|<|f||x|1n}
{
Let $A$ be normed $\Omega$\Hyph group.
For $n$\Hyph ary operation $\omega$,
\ShowEq{|a omega|<|omega||a|1n}
}

\DefEq
{
\begin{definition}
\labelDefinition{norm of representation}
Let
\ShowEq{f:A->*B}g{A_1}{A_2}
be representation
\ePrints{1305.4547}
\ifx\Semafor\ValueOn
\footnote{
See the definition
\EqParm{ref definition: left-side representation of algebra}{=z}
of the representation of universal algebra.
According to the definition
\Ref[1211.6965]{definition: module over ring},
module is representation of ring in Abelian group.
Since ring and Abelian group are $\Omega$\Hyph groups,
then module is representation of $\Omega$\Hyph group.
}
\fi
of $\Omega_1$\Hyph group $A_1$ with norm $\|x\|_1$
in $\Omega_2$\Hyph group $A_2$ with norm $\|x\|_2$.
The value
\ShowEq{norm of representation}
\ShowEq{norm of representation, definition}
is called
\AddIndex{norm of representation}{norm of representation} $f$.
\qed
\end{definition}
}
{definition: norm of representation}

\DefTheorem{|fab|<|f||a||b|}
{
Let
\ShowEq{f:A->*B}g{A_1}{A_2}
be representation
of $\Omega_1$\Hyph group $A_1$ with norm $\|x\|_1$
in $\Omega_2$\Hyph group $A_2$ with norm $\|x\|_2$.
Then
\ShowEq{|fab|<|f||a||b|}
}

\DefEq
{
\begin{definition}
\labelDefinition{open set}
Let $A$ be normed $\Omega$\Hyph group.
A set
$U\subset A$
is called
\AddIndex{open}{open set},\,\footnote{
In topology,
we usually define an open set before we
define base of topology.
In the case of a metric or normed space,
it is more convenient to define an open set,
based on the definition of base of topology.
In such case, the definition is based on one of the properties of
base of topology. An immediate proof
allows us to see that defined such
an open set satisfies the basic properties.
}
if, for any $A$\Hyph number
\EqParm{a in U}{=c}
there exists
\ShowEq{epsilon in R}
such that
\ShowEq{B(a) subset U}
\qed
\end{definition}
}
{definition: open set}

\DefEq
{
\begin{definition}
\labelDefinition{compact set}
A set $T$ of topological space is called
\AddIndex{compact}{compact set},
if every open cover of $T$ has
finite subcover.\,\footnote{
See also the definition
\citeBib{Kolmogorov Fomin}-1, page 92.
}
\qed
\end{definition}
}
{definition: compact set}

\DefTheorem{c1+c2 in B}
{
Let $A$ be normed $\Omega$\Hyph group.
For
\ShowEq{c1 c2 in A}
let
\ShowEq{c1 c2 in B}
Then
\ShowEq{c1+c2 in B}
}

\AddEq{definition: limit of sequence}
{
\begin{definition}
\labelDefinition{limit of sequence, normed \Algebrac}
Let $\Module$
be normed $\SideA$\Hyph \Algebrab.
$\Module$\Hyph number $a$ is called 
\AddIndex{limit of a sequence}{limit of sequence}
\ShowEq{[an] in A}
\ShowEq{limit of sequence}
if for any
\ShowEq{epsilon in R}
there exists positive integer $n_0$ depending on $\epsilon$ and such, that
\DrawEq{|an-a|<epsilon}{-}
for every $n>n_0$.
We also say that
\AddIndex{sequence $a_n$ converges}{sequence converges}
to $a$.
\qed
\end{definition}
}

\AddEq{theorem: limit of sequence}
{
\begin{theorem}
\labelTheorem{limit of sequence, normed \Algebrac}
Let $\Module$
be normed $\SideA$\Hyph \Algebrab.
$\Module$\Hyph number $a$ is
limit of a sequence
\ShowEq{[an] in A}
\DrawEq{a=lim an}{}
if for any
\ShowEq{epsilon in R}
there exists positive integer $n_0$ depending on $\epsilon$ and such, that
\DrawEq{an in B(a)}{}
for every $n>n_0$.
\end{theorem}
}

\DefProof{limit of sequence}
{
The theorem follows from definitions
\RefDefinition{open ball},
\RefDefinition{limit of sequence, normed \Algebrac}.
}

\AddEq{definition: fundamental sequence}
{
\begin{definition}
\labelDefinition{fundamental sequence, normed \Algebrac}
Let $\Module$
be normed $\SideA$\Hyph \Algebrab.
The sequence
\ShowEq{[an] in A}
is called 
\AddIndex{fundamental}{fundamental sequence}
or \AddIndex{Cauchy sequence}{Cauchy sequence},
if for every
\ShowEq{epsilon in R}
there exists positive integer $n_0$ depending on $\epsilon$ and such, that
\DrawEq{|ap-aq|<epsilon}{-}
for every $p$, $q>n_0$.
\qed
\end{definition}
}

\AddEq{theorem: fundamental sequence}
{
\begin{theorem}
\labelTheorem{fundamental sequence, normed \Algebrac}
Let $\Module$
be normed $\SideA$\Hyph \Algebrab.
The sequence
\ShowEq{[an] in A}
is fundamental sequence,
iff for every
\ShowEq{epsilon in R}
there exists positive integer $n_0$ depending on $\epsilon$ and such, that
\ShowEq{aq in B(ap)}
for every $p$, $q>n_0$.
\end{theorem}
}

\DefProof{fundamental sequence}
{
The theorem follows from definitions
\RefDefinition{open ball},
\RefDefinition{fundamental sequence, normed \Algebrac}.
}

\DefTheorem{norm of compact set is bounded}
{
Let $C$ be compact set of
normed $\Omega$\Hyph group $A$.
Then the norm
\ShowEq{|x in C|}
is bounded from both sides.
}

\DefTheorem{lim a=lim b}
{
Let $A$ be normed $\Omega$\Hyph group.
Let $a_n$, $b_n$, $n=1$, ..., be fundamental sequences.
Let
\DrawEq{lim a-b=0}{limit}
If the sequence $a_n$ converges, then
the sequence $b_n$ converges and
\ShowEq{lim a=lim b}
}

\DefEq
{
\begin{definition}
\labelDefinition{complete Omega group}
Normed $\Omega$\Hyph group $A$ is called
\AddIndex{complete}{complete Omega group}
if any fundamental sequence of elements
of $\Omega$\Hyph group $A$ converges, i.e.
has limit in $\Omega$\Hyph group $A$.
\qed
\end{definition}
}
{definition: complete Omega group}

\DefEq
{
\begin{definition}
\labelDefinition{limit of sequence, map to Omega group}
Let
\ShowEq{fn M(X,A)}
be sequence of maps into normed
$\Omega$\Hyph group $A$.
The map
\ShowEq{f M(X,A)}
is called
\AddIndex{limit of sequence}{limit of sequence}
$f_n$, if for any $x\in X$
\DrawEq{f(x)=lim}{}
We also say that
\AddIndex{sequence $f_n$ converges}{sequence converges}
to the map $f$.
\qed
\end{definition}
}
{definition: limit of sequence, map to Omega group}

\DefEq
{
\begin{definition}
\labelDefinition{sequence converges uniformly}
Let
\ShowEq{fn M(X,A)}
be sequence of maps into normed
$\Omega$\Hyph group $A$.
\AddIndex{Sequence $f_n$ converges uniformly}
{sequence converges uniformly}
to the map $f$,
if for any
\ShowEq{epsilon in R}
there exists $N$ such that
\ShowEq{fn(x) - f(x)}
for any $n>N$.
\qed
\end{definition}
}
{definition: sequence converges uniformly}

\DefEq
{
\begin{theorem}
\labelTheorem{sequence converges uniformly, fn-fm}
Sequence of maps
\ShowEq{fn M(X,A)}
into normed
$\Omega$\Hyph group $A$ converges uniformly
to the map $f$,
if for any
\ShowEq{epsilon in R}
there exists $N$ such that
\ShowEq{|fn(x)-fm(x)|<e}
for any \(n\), \(m>N\).
\end{theorem}
}
{theorem: sequence converges uniformly, fn-fm}

\DefTheorem{h=f+g, converges uniformly}
{
Let sequence of maps
\ShowEq{fn M(X,A)}
into complete $\Omega$\Hyph group $A$
converge uniformly
to the map $f$.
Let sequence of maps
\ShowEq{gn M(X,A)}
into complete $\Omega$\Hyph group $A$
converge uniformly
to the map $g$.
Then sequence of maps
\ShowEq{hn=fn+gn}
into complete $\Omega$\Hyph group $A$
converges uniformly
to the map
\DrawEq{h=f+g}{converges uniformly}
}

\DefTheorem{representation f generates representation fX}
{
The representation
\ShowEq{f:A->*B}g{A_1}{A_2}
of $\Omega_1$\Hyph group $A_1$ with norm $\|x\|_1$
in $\Omega_2$\Hyph group $A_2$ with norm $\|x\|_2$
generates representation
\ShowEq{f:A->*B}{f_X}{M(X,A_1)}{M(X,A_2)}
of $\Omega_1$\Hyph group
\ShowEq{M(X,A1)}
in $\Omega_2$\Hyph group
\ShowEq{M(X,A2)}
where (
\ShowEq{g1 in M(X,A1)},
\ShowEq{g2 in M(X,A2)}
)
\ShowEq{f(x)g(x):X->A2}
}

\DefTheorem{fX(g1)(g2), converges uniformly}
{
Let
\ShowEq{f:A->*B}g{A_1}{A_2}
be representation
of complete $\Omega_1$\Hyph group $A_1$ with norm $\|x\|_1$
in complete $\Omega_2$\Hyph group $A_2$ with norm $\|x\|_2$.
Let sequence of maps
\ShowEq{g1n M(X,A1)}
converge uniformly
to the map $g_1$.
Let sequence of maps
\ShowEq{g2n M(X,A2)}
converge uniformly
to the map $g_2$.
Let range of the map $g_i$, $i=1$, $2$,
be compact set.
Then the sequence of maps
\ShowEq{f(g1n)(g2n)}
converge uniformly
to the map
\ShowEq{fX(g1)(g2)}.
}

\DefTheorem{set of maps to Omega group}
{
Let
\ShowEq{set of maps to Omega group}
be set of maps of a set $X$ to $\Omega$\Hyph group $A$.
We can define the structure of $\Omega$\Hyph group on the set
\EqParm{M(X,A)}{=.}
}

\DefEq
{
Since $X$ is an arbitrary set, we cannot define
norm in $\Omega$\Hyph group
\EqParm{M(X,A)}{=.}
However we can define the convergence of the sequence in
\EqParm{M(X,A)}{=.c}
therefore, we can define a topology in
\EqParm{M(X,A)}{=.}
}
{remark: set of maps to Omega group}

\DefDefinition{closed ball}
{
Let $\Module$
be normed $\SideA$\Hyph \Algebrab.
Let $a\in \Module$.
The set
\ShowEq{closed ball}
is called
\AddIndex{closed ball}{closed ball}
with center at $a$.
}

\DefDefinition{open ball}
{
Let $\Module$
be normed $\SideA$\Hyph \Algebrab.
Let $a\in \Module$.
The set
\ShowEq{open ball}
is called
\AddIndex{open ball}{open ball}
with center at $a$.
}

\DefDefinition{reduced polymorphism of representations}
{
Let
\ShowEq{set of universal algebras 1}
be universal algebras.
Let, for any
\ShowEq{k,1n}
\ShowEq{f:A->*B}{f_k}A{B_k}
be effective representation of $\Omega_1$\Hyph algebra $A$
in $\Omega_2$\Hyph algebra $B_k$.
Let
\ShowEq{f:A->*B}fAB
be effective representation of $\Omega_1$\Hyph algebra $A$
in $\Omega_2$\Hyph algebra $B$.
The map
\ShowEq{reduced polymorphism of representation}
is called
\AddIndex{reduced polymorphism of representations}
{reduced polymorphism of representations}
$f_1$, ..., $f_n$ into representation $f$,
if, for any
\ShowEq{k,1n}
provided that all variables except the variable $x_k\in B_k$
have given value, the map $r_2$
is a reduced morphism of representation $f_k$ into representation $f$.

If $f_1=...=f_n$, then we say that the map $r_2$
is reduced polymorphism of representation $f_1$ into representation $f$.

If $f_1=...=f_n=f$, then we say that the map $r_2$
is reduced polymorphism of representation $f$.
}

\DefEq
{
\begin{convention}
We will use Einstein summation convention.
When an index is present in
an expression twice (one above and one below) and a set of index is known,
we have the sum with respect to repeated index.
In this case we assume that we know the set
of summation index and do not use summation symbol
\ShowEq{av=sum av}av
If needed to clearly
show set of index, I will do it.
\qed
\end{convention}
}
{convention: Einstein summation}

\AddEq [1]{step of proof: action in Abelian group}
{
From the equality
\EqRef{0g=0},
it follows that the equality
\EqRef{#1}
is true when $n=0$.
\def\TempA{(nm)a=n(ma)}
\def\TempB{#1}
\begin{itemize}
\item
Let the equality
\EqRef{#1}
is true when $n=k\ge 0$.
Then
\ShowEq{#1,n=k}
The equality
\ShowEq{#1,n=k+1}
follows from
\ifx\TempA\TempB
equalities
\EqRef{(n+1)g=},
\EqRef{(m+n)a=ma+na}.
\else
the equality
\EqRef{(n+1)g=}.
\fi
Therefore, the equality
\EqRef{#1}
is true when $n=k+1$.
According to mathematical induction,
the equality
\EqRef{#1}
is true for any $n\ge 0$.
\item
Let the equality
\EqRef{(m+n)a=ma+na}
is true when $n=k\le 0$.
Then
\ShowEq{#1,n=k}
The equality
\ShowEq{#1,n=k-1}
follows from
\ifx\TempA\TempB
equalities
\EqRef{(n-1)g=},
\EqRef{(m-n)a=ma-na}.
\else
the equality
\EqRef{(n-1)g=}.
\fi
Therefore, the equality
\EqRef{#1}
is true when $n=k-1$.
According to mathematical induction,
the equality
\EqRef{#1}
is true for any $n\le 0$.
\item
Therefore, the equality
\EqRef{#1}
is true for any $n\in Z$.
\end{itemize}
}

\DefDefinition{category of representations A1(mA2)}
{
Let $A_1$ be $\Omega_1$\Hyph algebra.
Let $\mathcal A_2$ be category of $\Omega_2$\Hyph algebras.
We define \AddIndex{category
\ShowEq{A1(mA2) symb}
of representations}
{category of representations}
of $\Omega_1$\Hyph algebra $A_1$ in $\Omega_2$\Hyph algebra.
Representations of $\Omega_1$\Hyph algebra $A_1$ in $\Omega_2$\Hyph algebra
are objects of this category.
Reduced morphisms of corresponding representations are morphisms of this category.
}

\DefEq
{
\begin{definition}
\labelDefinition{matrix operations}
{\it
We define following operations on the set of matrices
over $\Omega$\Hyph ring $A$.
\StartLabelItem
\begin{enumerate}
\item
Let the nubmer of columns of the matrix $a$ equal the number of rows of the matrix $b$.
\AddIndex{\RC product}{rc-product}
of matrices $a$ and $b$
has form
\ShowEq{rc-product}
\ShowEq{rc-product of matrices}
and can be expressed as product of a row of matrix $a$
over a column of matrix $b$.\,\footnote{
To draw symbol of \RC product, we put two symbols of product in the place
of index which participate in sum.
}
\item
Let the nubmer of rows of the matrix $a$ equal the number of columns of the matrix $b$.
\AddIndex{\CR product}{cr-product}
of matrices $a$ and $b$
has form
\ShowEq{cr-product}
\ShowEq{cr-product of matrices}
and can be expressed as product of a column of matrix $a$
over a row of matrix $b$.\,\footnote{
To draw symbol of \CR product, we put two symbols of product in the place
of index which participate in sum.
}
\item
The sum of matrices $a$ and $b$ is defined by the equality
\ShowEq{(a+b)=}
\item
The transpose $a^T$ of the matrix $a$
exchanges rows and columns
\ShowEq{transpose of matrix, 1}
\end{enumerate}
}
\qed
\end{definition}
}
{definition: matrix operations}

\DefEq
{
\begin{remark}
\labelRemark{pronunciation of product}
We will use symbol \RC\  or \CR\  in name of properties of each product
and in the notation.
We can read symbols
$\RCstar$ and $\CRstar$ as $rc$\hyph product and $cr$\hyph product.
This rule we extend to following terminology.
\qed
\end{remark}
}
{remark: pronunciation of product}

\DefEq
{
\begin{definition}
\labelDefinition{biring}
The set $A$ is a \AddIndex{biring}{biring}
if we defined on $A$ an unary operation, say transpose,
and three binary operations,
say \RC product, \CR product and sum, such that
\begin{itemize}
\item \RC product and sum define structure of ring on $A$
\item \CR product and sum define structure of ring on $A$
\item both products have common identity $\delta$
\item products satisfy equation
\ShowEq{rcstar transpose}
\item transpose of identity is identity
\ShowEq{transpose of identity}
\item double transpose is original element
\ShowEq{double transpose}
\end{itemize}
\qed
\end{definition} 
}
{definition: biring}

\DefEq
{
\begin{theorem}
[\AddIndex{duality principle for biring}
{duality principle for biring}]
\labelTheorem{duality principle for biring}
Let $\mathcal{A}$ be true statement about
biring $A$.
If we exchange the same time
\begin{itemize}
\item $a\in A$ and $a^T$
\item \RC product and \CR product
\end{itemize}
then we soon get true statement.
\end{theorem}
}
{theorem: duality principle for biring}

\DefEq
{
\begin{theorem}
[\AddIndex{duality principle for biring of matrices}
{duality principle for biring of matrices}]
\labelTheorem{duality principle for biring of matrices}
Let $A$ be biring of matrices.
Let $\mathcal{A}$ be true statement about matrices.
If we exchange the same time
\begin{itemize}
\item rows and columns of all matrices
\item \RC product and \CR product
\end{itemize}
then we soon get true statement.
\end{theorem}
\begin{proof}
This is the immediate consequence
of the theorem \RefTheorem{duality principle for biring}.
\end{proof}
}
{theorem: duality principle for biring of matrices}

\DefEq
{
\begin{remark}
\labelRemark{reducible biring}
If product in $\Omega$\Hyph $A$ ring is commutative, then
\ShowEq{reducibility of products}
\AddIndex{Reducible biring}{reducible biring} is
the biring which holds
\AddIndex{condition of reducibility of products}
{condition of reducibility of products}
\EqRef{reducibility of products}.
So, in reducible biring, it is enough to consider only
\RC product. However in case when
the order of factors is essential
we will use \CR product also.
\qed
\end{remark}
}
{remark: reducible biring}

\DefEq
{
Since left\Hyph side representation and right\Hyph side representation
are based on homomorphism of $\Omega$\Hyph group,
then the following statement is true.

\begin{theorem}[duality principle for representation of multiplicative $\Omega$\Hyph group]
\labelTheorem{duality principle, algebra representation}
Any statement which holds for
left\Hyph side representation of multiplicative $\Omega$\Hyph group $A_1$
holds also for right\Hyph side representation of multiplicative $\Omega$\Hyph group $A_1$, if
\begin{itemize}
\item
we will use right\Hyph side product over $A_1$\Hyph number $a_1$ instead of
left\Hyph side product over $A_1$\Hyph number $a_1$;
\item
we will use right shift $R(a_1)$ of multiplicative $\Omega$\Hyph group $A_1$
instead of left shift $L(a_1)$.
\end{itemize}
\qed
\end{theorem}
}
{theorem: duality principle, algebra representation}

\DefEq
{
\begin{theorem}
\labelTheorem{left right shift}
The product in multiplicative $\Omega$\Hyph group $A_1$ determines two different representations
of multiplicative group $A_1$
in $\Omega$\Hyph algebra $A_1$:
\begin{itemize}
\EqParm{theorem: shift}{=left}
\EqParm{theorem: shift}{=right}
\end{itemize}
\end{theorem}
}
{theorem: left right shift}

\DefEq
{
\item
The \AddIndex{\SideWS shift}{\SideWS shift}
\ShowEq{\SideWS shift =}
\ShowEq{\SideWS shift}
is \SideNS\Hyph side representation
\ShowEq{\SideWS shift, product}
}
{theorem: shift}

\DefEq
{
If multiplicative $\Omega$\Hyph group $A_1$ is Abelian,
then there is no difference between left\Hyph side and right\Hyph side representations.

\begin{definition}
\labelDefinition{representation of Abelian group}
Let $A_1$ be Abelian multiplicative $\Omega$\Hyph group.
We call a homomorphism of multiplicative $\Omega$\Hyph group
\DrawEq{representation of group, map}{Abelian group}
\AddIndex{representation}{representation of algebra}
of multiplicative $\Omega$\Hyph group $A_1$
or
\AddIndex{$A_1$\Hyph representation}{A representation of algebra}
in $\Omega_2$\Hyph algebra $A_2$ if the map $f$ holds
\DrawEq{left-side representation of group}{Abelian}
\qed
\end{definition}

Usually we identify a representation of the Abelian multiplicative $\Omega$\Hyph group $A_1$
and a left\Hyph side representation of the multiplicative $\Omega$\Hyph group $A_1$.
However, if it is necessary for us,
we identify a representation of the Abelian multiplicative $\Omega$\Hyph group $A_1$
and a right\Hyph side representation of the multiplicative $\Omega$\Hyph group $A_1$.
}
{definition: representation of Abelian group}

\DefDefinition{algebra of sets}
{
A nonempty system of sets $\mathcal R$ is called
\AddIndex{ring of sets}{ring of sets},\,\footnote{
See also the definition
\citeBib{Kolmogorov Fomin}\Hyph 1,  page 31.}
if condition
\ShowEq{A, B in R}
imply
\ShowEq{A, B in R 1}
A set
\ShowEq{E in R}
is called
\AddIndex{unit of ring of sets}{unit of ring of sets}
if
\ShowEq{A cap E=A}
A ring of sets with unit is called
\AddIndex{algebra of sets}{algebra of sets}.
}

\DefEq
{
\begin{remark}
\labelRemark{A cup/minus B}
For any
\ShowEq{A cup/minus B}
Therefore, if
\ShowEq{A, B in R},
then
\ShowEq{A cup/minus B R}
\qed
\end{remark}
}
{remark: A cup/minus B}

\DefDefinition{sigma algebra of sets}
{
The ring of sets $\mathcal R$ is called
\AddIndex{\(\sigma\)\Hyph ring of sets}{sigma ring of sets},\,\footnote{
See similar definition in
\citeBib{Kolmogorov Fomin}, page 35, definition 3.}
if condition
\ShowEq{Ai in R}
imply
\ShowEq{Ai in R 1}
\(\sigma\)\Hyph Ring of sets with unit is called
\AddIndex{\(\sigma\)\Hyph algebra of sets}{sigma algebra of sets}.
}

\DefDefinition{Borel algebra}
{
Minimal \(\sigma\)\Hyph algebra
\ShowEq{Borel algebra}
generated by the set of all open balls of normed $\Omega$\Hyph group $A$
is called
\AddIndex{Borel algebra}{Borel algebra}.\,\footnote{
See remark in
\citeBib{Kolmogorov Fomin}, p. 36.
According to the remark
\ref{remark: A cup/minus B},
the set of closed balls also generates
Borel algebra.
}
A set
belonging to Borel algebra
is called
\AddIndex{Borel set}{Borel set}
or
\AddIndex{$B$\Hyph set}{B set}.
}

\DefDefinition{C_X C_Y measurable}
{
Let $\mathcal C_X$ be
\(\sigma\)\Hyph algebra of sets of set $X$.
Let $\mathcal C_Y$ be
\(\sigma\)\Hyph algebra of sets of set $Y$.
The map
\ShowEq{f:A->B}fXY
is called
$(\mathcal C_X,\mathcal C_Y)$\Hyph measurable\,\footnote{
See similar definition in
\citeBib{Kolmogorov Fomin}, page 284.}
if for any set
\ShowEq{CX CY measurable map}
}

\DefEq
{
\begin{example}
\labelExample{measurable map}
Let $\mu$ be a \(\sigma\)\Hyph additive measure defined on the set $X$.
Let
\ShowEq{algebra of sets measurable with respect to measure}
be \(\sigma\)\Hyph algebra of sets measurable with respect to measure $\mu$.
Let
$\mathcal B(A)$
be Borel algebra of
normed $\Omega$\Hyph group $A$.
The map
\ShowEq{f:A->B}fXA
is called
\AddIndex{$\mu$\Hyph measurable}{mu measurable map}\,\footnote{
See similar definition in
\citeBib{Kolmogorov Fomin},  pp. 284, 285, definition 1.
If the measure $\mu$ is defined on the set $X$ by context,
then we also call the map
\ShowEq{f:A->B}fXA
\AddIndex{measurable}{measurable map}.
}
if for any set
\ShowEq{measurable map}
\qed
\end{example}
}
{example: measurable map}

\DefDefinition{simple map}
{
Let $\mu$ be a \(\sigma\)\Hyph additive measure defined on the set $X$.
The map
\ShowEq{f:A->B}fXA
into normed $\Omega$\Hyph group $A$ is called
\AddIndex{simple map}{simple map},
if this map is $\mu$\Hyph measurable and
its range is finite or countable set.
}

\DefTheorem{measurable simple map}
{
Let range
\ShowEq{range f}
of the map
\ShowEq{f:A->B}fXA
be finite or countable set.
The map $f$ is $\mu$\Hyph measurable iff
all sets
\DrawEq{Fn=}{}
are $\mu$\Hyph measurable.\,\footnote{
See similar theorem in
\citeBib{Kolmogorov Fomin},  page 286, theorem 4.}
}

\DefDefinition{measure}
{
Let \(\mathcal C_m\) be semiring of sets.\,\footnote{
See also the definition
\citeBib{Kolmogorov Fomin}\Hyph 1 on page 270.}
The map
\ShowEq{f:A->B}f{\mathcal C_m}R
is called
\AddIndex{measure}{measure},
if
\StartLabelItem
\begin{enumerate}
\item
\ShowEq{m(A)>0}
\item
The map \(m\) is additive map.
If a set
\ShowEq{A in Cm}
has finite expansion
\DrawEq{A=+Ai Cm}{}
where
\ShowEq{Ai*Aj=0},
then
\ShowEq{m A=+m Ai}
\end{enumerate}
}

\DefDefinition{complete measure}
{
A measure \(\mu\) is called
\AddIndex{complete measure}{complete measure},
if conditions
\ShowEq{B in A, mu(A)=0}
imply
that \(B\) is measurable set.
}

\DefDefinition{sigma-additive measure}
{
Let $\mathcal C_{\mu}$ be
\(\sigma\)\Hyph algebra of sets of set $F$.\,\footnote{
See similar definitions in
\citeBib{Kolmogorov Fomin}, definition 1 on page 270
and definition 2 on page 272.}
The map
\ShowEq{f:A->B}{\mu}{\mathcal C_{\mu}}R
into real field $R$ is called
\AddIndex{\(\sigma\)\Hyph additive measure}{sigma-additive measure},
if, for any set
\ShowEq{X in Cmu},
following conditions are true.
\StartLabelItem
\begin{enumerate}
\item
\ShowEq{muX>0}
\item
Let
\DrawEq{X=+Xi}{}
be finite or countable union of sets
\ShowEq{Xn in Cmu}
Then
\DrawEqParm{mu(X=+Xi)}X{}
where series on the right converges absolutly.
\labelItem{mu(X=+Xi)}
\end{enumerate}
}

\DefEq
{
Let $\mu$ be a \(\sigma\)\Hyph additive measure defined on the set $X$.
Let effective representation of real field \(R\)
in complete Abelian \(\Omega\)\Hyph group \(A\) be defined.\,\footnote{In
other words,
\(\Omega\)\Hyph group \(A\) is
\(R\)\Hyph vector space.}
}
{remark - measure - effective representation of real field}

\DefDefinition{series converges normally}
{
Let \(a_i\) be a sequence of \(A\)\Hyph numbers.
If
\ShowEq{sum |ai|<infty}
then we say that the {\bf series}
\ShowEq{sum ai}
\AddIndex{converges normally}
{series converges normally}.\,\footnote{See also
the definition of normal convergence of the series
on page \citeBib{Cartan differential form}-12.}
}

\DefDefinition{Integral of Map over Set of Finite Measure}
{
\(\mu\)\Hyph measurable map
\ShowEq{f:A->B}fXA
is called
\AddIndex{integrable map}{integrable map}
over the set \(X\),\,\footnote{
See also the definition in
\citeBib{Kolmogorov Fomin}, page 296.
}
if there exists a sequence
of simple integrable over the set \(X\) maps
\ShowEq{f:A->B}{f_n}XA
converging uniformly to \(f\).
Since the map \(f\) is integrable map,
then the limit
\ShowEq{integral of map}
\ShowEq{integral of measurable map}
is called
\AddIndex{Lebesgue integral}{Lebesgue integral}
of map \(f\) over the set \(X\).
}

\DefDefinition{integrable map}
{
For simple map
\ShowEq{f:A->B}fXA
consider series
\DrawEq{sum mu(F) f}{integral}
where
\begin{itemize}
\item
The set
\ShowEq{f1...}
is domain of the map \(f\)
\item
Since \(n\ne m\), then
\ShowEq{fn ne fm}
\item
\DrawEq{Fn=}{-}
\end{itemize}
Simple map
\ShowEq{f:A->B}fXA
is called
\AddIndex{integrable map}{integrable map}
over the set \(X\)
if series
\eqRef{sum mu(F) f}{integral}
converges normally.\,\footnote{
See similar definition in
\citeBib{Kolmogorov Fomin}, p. 294.}
Since the map \(f\) is integrable map,
then sum of series
\eqRef{sum mu(F) f}{integral}
is called
\AddIndex{Lebesgue integral}{Lebesgue integral}
of map \(f\) over the set \(X\)
\ShowEq{integral of map}
\ShowEq{integral of simple map}
}

\DefTheorem{|int f|<int|f|}
{
Let
\ShowEq{f:A->B}fXA
be measurable map.
Integral
\ShowEq{int f X}
exists
iff
integral
\ShowEq{int |f| X}
exists.
Then
\DrawEq{|int f|<int|f|}{measurable}
}

\DefTheorem{int |f|<M mu}
{
Let
\ShowEq{f:A->B}fXA
be $\mu$\Hyph measurable map such that
\DrawEq{|f(x)|<M}{}
Since the measure of the set \(X\) is finite, then
\DrawEq{int |f|<M mu}{measurable}
}

\DefTheorem{int f+g X}
{
Let
\ShowEq{f:A->B}fXA
\ShowEq{f:A->B}gXA
be $\mu$\Hyph measurable maps
with compact range.
Since there exist integrals
\ShowEq{int f X}
\ShowEq{int g X}
then there exists integral
\ShowEq{int f+g X}
and
\DrawEq{int f+g X =}{measurable}
}

%% file: Stmt.Representation.Eq.tex

\DefEq
{
\begin{theorem}
\labelTheorem{rcstar transpose}
\ShowEq{rcstar transpose, 0}
\end{theorem}
}
{theorem: rcstar transpose}

\AddEq{direct sum of Abelian groups}
{
\symb{A\oplus B}{direct sum}1
}

\AddEquation{reduced polymorphism of representation, akl}
{
\begin{array}{r@{\,}l}
&r_2(m_1,...,\BlueText{f_k(a)(m_k)},...,m_l,...,m_n)
\\
=&r_2(m_1,...,m_k,...,\BlueText{f_l(a)(m_l)},...,m_n)
\end{array}
}

\AddEquation{reduced polymorphism of representation, omega2}
{
\begin{array}{r@{\,}l}
&\,r_2(m_1,...,m_{k\cdot 1}...m_{k\cdot p}\omega_2,...,m_n)
\\
=&\,
r_2(m_1,...,m_{k\cdot 1},...,m_n)
...
r_2(m_1,...,m_{k\cdot p},...,m_n)
\omega_2
\end{array}
}

\AddEq{kl,1n}
{
$k$, $l$, $k=1$, ..., $n$, $l=1$, ..., $n$,
}

\AddEq{reduced morphism representation =}
{
r_2(f(a)(m))=g(a)(r_2(m))
}

\AddEquation{morphism of representations of universal algebra, definition, 2m-1}
{
\BlueText{r_2(f(a')(r_2^{-1}(m')))}
=f(\RedText{a'})(\BlueText{m'})
}

\AddEq{A in xAi}
{
\[A\subseteq\prod_{i\in I}A_i\]
}

\AddEquation{A=o+Ai}
{
A=\bigoplus_{i\in I}A_i
}

\AddEq[2]{End O2}
{
$\End(\Omega_2,#1)$#2
}

\AddEquation{fxi,i=sum fxi}
{
f((x_i,i\in I))=\sum_{i\in I}f_i(x_i)
}

\AddEq{fi(x+y)=}
{
\[
f_i(x_i+y_i)=f_i(x_i)+f_i(y_i)
\]
}

\AddEq{f(x+y)i=}
{
\begin{align*}
f((x_i+y_i,i\in I))&=
\sum_{i\in I}f_i(x_i+y_i)=
\sum_{i\in I}(f_i(x_i)+f_i(y_i))
\\&=\sum_{i\in I}f_i(x_i)+\sum_{i\in I}f_i(y_i)
\\&=f((x_i,i\in I))+f((y_i,i\in I))
\end{align*}
}

\AddEq{structure of subrepresentations}
{
\StartLabelItem
\begin{enumerate}
\item
$X_0=X$
\labelItem{X0=X}
\item
$x\in X_k=>x\in X_{k+1}$
\labelItem{Xk in Xk+1}
\item
$x_1\in X_k$, ..., $x_n\in X_k$, $\omega\in\Omega_2(n)
=>x_1...x_n\omega\in X_{k+1}$
\labelItem{x1n omega in Xk+1}
\item
$x\in X_k$, $a\in A=>f(a)(x)\in X_{k+1}$
\labelItem{ax in Xk+1}
\end{enumerate}
}

\AddEquation{structure of subrepresentations, 1}
{
\bigcup_{k=0}^{\infty}X_k=J_f(X)
}

\AddEquation{lj(x)=0x0}
{
\lambda_j(x)=(\delta^i_jx,i\in I)
}

\AddEq{flj=fj}
{
\[
f\circ\lambda_j(x)=\sum_{i\in I}f_i(\delta_j^ix)=f_j(x)
\]
}

\AddEq[1]{(xi)in A}
{
$(x_i,i\in I)\in A$#1
}

\DefEq
{
$\|x\|$, $x\in C$,
}
{|x in C|}

\DefEquation
{
\|a-b\|\ge|\,\|a\|-\|b\|\,|
}
{|a-b|>|a|-|b|}

\DefEq
{
\[
\|x'-x\|_1<\delta
\]
}
{|x'-x|<delta}

\DefEq
{
\[
\|f(x')-f(x)\|_2<\epsilon
\]
}
{|f(x)-f(x')|<epsilon}

\AddEq[1]{End A(O2) B}
{
$\End(A(\Omega_2);B)$#1
}

\AddEq{open ball}
{
\symb{B_o(a,\rho)}{open ball}{}
\[
\ShowSymbol{open ball}{}=
\{
b\in A:\|b-a\|<\rho
\}
\]
}

\AddEq{closed ball}
{
\symb{B_c(a,\rho)}{closed ball}{}
\[
\ShowSymbol{closed ball}{}=
\{
b\in A:\|b-a\|\le\rho
\}
\]
}

\DefEq
{
$f_n\in M(X,A)$, $n=1$, ...,
}
{fn M(X,A)}

\DefEquation
{
\|f_n(x)-f_m(x)\|<\epsilon
}
{|fn(x)-fm(x)|<e}

\DefEq
{
\|a_p-a_q\|<\epsilon
}
{|ap-aq|<epsilon}

\AddEq{[an] in A}
{
$\{a_n\}$, $a_n\in \Module$,
}

\AddEq{aq in B(ap)}
{
\[a_q\in B_o(a_p,\epsilon)\]
}

\AddEq{an in B(a)}
{
a_n\in B_o(a,\epsilon)
}

\AddEq{limit of sequence}
{
\symb{\lim_{n\rightarrow\infty}a_n}{limit of sequence}{}
\[
a=\ShowSymbol{limit of sequence}{}
\]
}

\AddEq{a=lim an}
{
a=\lim_{n\rightarrow\infty}a_n
}

\DefEq
{
\|a_n-a\|<\epsilon
}
{|an-a|<epsilon}

\DefEquation
{
\lim_{n\rightarrow\infty}a_n=\lim_{n\rightarrow\infty}b_n
}
{lim a=lim b}

\DefEq
{
\lim_{n\rightarrow\infty}(a_n-b_n)=0
}
{lim a-b=0}

\DefEq
{
$g_n\in M(X,A)$, $n=1$, ...,
}
{gn M(X,A)}

\DefEq
{
\[h_n=f_n+g_n\]
}
{hn=fn+gn}

\DefEq
{
h=f+g
}
{h=f+g}

\DefEq
{
$c_1$, $c_2\in A$,
}
{c1 c2 in A}

\DefEq
{
$c_1\in B_c(a_1,R_1)$, $c_2\in B_c(a_2,R_2)$.
}
{c1 c2 in B}

\DefEq
{
$g_{1\cdot n}\in M(X,A_1)$, $n=1$, ...,
}
{g1n M(X,A1)}

\DefEq
{
$g_{2\cdot n}\in M(X,A_2)$, $n=1$, ...,
}
{g2n M(X,A2)}

\DefEq
{
$f_X(g_1)(g_2)$
}
{fX(g1)(g2)}

\DefEq
{
$M(X,A_1)$
}
{M(X,A1)}

\DefEq
{
$M(X,A_2)$
}
{M(X,A2)}

\DefEq
{
$g_1\in M(X,A_1)$
}
{g1 in M(X,A1)}

\DefEq
{
$g_2\in M(X,A_2)$
}
{g2 in M(X,A2)}

\DefEquation
{
\begin{split}
f_X(g_1)(g_2)&:X->A_2
\\
(f_X(g_1)(g_2))(x)&=f(g_1(x))(g_2(x))
\end{split}
}
{f(x)g(x):X->A2}

\DefEq
{
\symb{\|\omega\|}{norm of operation}{}
}
{norm of operation}

\DefEquation
{
\ShowSymbol{norm of operation}{}=
\text{sup}\frac{\| a_1...a_n\omega\|}{\|a_1\|...\|a_n\|}
}
{norm of operation, definition}

\DefEquation
{
\|a_1...a_n\omega\|\le\|\omega\|\|a_1\|...\|a_n\|
}
{|a omega|<|omega||a|1n}

\DefEq
{
\symb{\|f\|}{norm of representation}{}
}
{norm of representation}

\DefEquation
{
\ShowSymbol{norm of representation}{}=
\text{sup}\frac{\|f(a_1)(a_2)\|_2}{\|a_1\|_1\|a_2\|_2}
}
{norm of representation, definition}

\DefEquation
{
\|f(a_1)(a_2)\|_2\le\|f\|\|a_1\|_1\|a_2\|_2
}
{|fab|<|f||a||b|}

\DefEquation
{
c_1+c_2\in B_c(a_1+a_2,R_1+R_2)
}
{c1+c2 in B}

\DefEq
{
$f\in M(X,A)$
}
{f M(X,A)}

\DefEq
{
f(x)=\lim_{n\rightarrow\infty}f_n(x)
}
{f(x)=lim}

\DefEq
{
\[\|f_n(x)-f(x)\|<\epsilon\]
}
{fn(x) - f(x)}

\DefEq
{
$a \in U$\Pt
}
{a in U}

\AddEq{sum mu(F) f}
{
\displaystyle\sum_n\mu(F_n)f_n
}

\AddEq{f1...}
{
\(\{f_1,f_2,...\}\)
}

\AddEq{fn ne fm}
{
\(f_n\ne f_m\)
}

\AddEq{Fn=}
{
F_n=\{x\in X:f(x)=f_n\}
}

\AddEq{mu(X=+Xi)}
{
\mu(\PX)=\sum_i\mu(\PX_i)
}

\AddEq{m A=+m Ai}
{
\[
m(A)=\sum_{i=1}^nm(A_i)
\]
\labelItem{m A=+m Ai}
}

\DefEq
{
\(A\), \(B\in\mathcal R\)
}
{A, B in R}

\DefEq
{
$A\cup B\in\mathcal R$, $A\setminus B\in\mathcal R$.
}
{A cup/minus B R}

\DefEq
{
\symb{\mathcal B(A)}{Borel algebra}1
}
{Borel algebra}

\DefEq
{
$C\in\mathcal C_Y$
\[f^{-1}(C)\in\mathcal C_X\]
}
{CX CY measurable map}

\DefEq
{
$C\in\mathcal B(A)$
\[f^{-1}(C)\in\mathcal C_{\mu}\]
}
{measurable map}

\AddEq{range f}
{
\(y_1\), \(y_2\), ...
}

\DefEq
{
\(m(A)\ge 0\)
\labelItem{m(A)>0}
}
{m(A)>0}

\DefEq
{
\(A\in\mathcal C_m\)
}
{A in Cm}

\DefEq
{
A=\bigcup_{i=1}^nA_i\ \ \ A_i\in\mathcal C_m
}
{A=+Ai Cm}

\DefEq
{
\symb{\mathcal C_{\mu}}
{algebra of sets measurable with respect to measure}1
}
{algebra of sets measurable with respect to measure}

\DefEq
{
\(A\triangle B\), \(A\cap B\in\mathcal R\).
}
{A, B in R 1}

\DefEq
{
$E\in\mathcal R$
}
{E in R}

\DefEq
{
\[A\cap E=A\]
}
{A cap E=A}

\DefEq
{
$A_i\in\mathcal R$, $i=1$, ..., $n$, ...,
}
{Ai in R}

\DefEq
{
$A$, $B$
\begin{align*}
A\cup B&=(A\Delta B)\Delta (A\cap B)
\\
A\setminus B&=A\Delta (A\cap B)
\end{align*}
}
{A cup/minus B}

\DefEq
{
\[\bigcup_nA_n\in\mathcal R\]
}
{Ai in R 1}

\DefEq
{
$B_o(a,\epsilon)\subset U$.
}
{B(a) subset U}

\DefEq
{
$M(X,A)$\Pt
}
{M(X,A)}

\DefEq
{
\symb{M(X,A)}{set of maps to Omega group}1
}
{set of maps to Omega group}

\DefEquation
{
(A\RCstar B)^T=A^T\CRstar B^T
}
{rcstar transpose, 0}

\DefEquation
{
(a\RCstar b)^T=a^T\CRstar b^T
}
{rcstar transpose}

\DefEquation
{
\delta^T=\delta
}
{transpose of identity}

\DefEquation
{
(a^T)^T=a
}
{double transpose}

\DefEq
{
$A_1$, ..., $A_n$\Pt
}
{A1n}

\DefEq
{
\[
\xymatrix
{
f_k:A\ar[r]|{*}&A_k
}
\]
}
{representation A Ak 1}

\DefEq
{
\[
\xymatrix{
r_1:\Times\ar[r]&S_1
&
r_2:\Times\ar[r]&S_2
}
\]
}
{polymorphisms category}

\DefEq
{
\[
\begin{matrix}
\xymatrix
{
g_1:A\ar[r]|{*}&S_1
}
&
\xymatrix
{
g_2:A\ar[r]|{*}&S_2
}
\end{matrix}
\]
}
{representation of algebra in S1 S2}

\DefEq
{
$r_1\rightarrow r_2$
}
{r1->r2}

\DefEq
{
\[
\xymatrix{
&S_1\ar[dd]^h
\\
\Times
\ar[ru]^{r_1}\ar[rd]_{r_2}
\\
&S_2
}
\]
}
{polymorphisms category, diagram}

\DefEq
{
\symb{\Tensor B}{tensor product}1
}
{tensor product of representations}

\DefEq
{
\[
\begin{matrix}
b_k\in B_k&k=1,...,n&
b_{i\cdot 1},...,b_{i\cdot p}\in B_i&\omega\in\Omega_2(p)&a\in A
\end{matrix}
\]
}
{equivalence, 1, representation, tensor product}

\DefEquation
{
\begin{split}
&b_1\otimes...\otimes(b_{i\cdot 1}...b_{i\cdot p}\omega)\otimes...\otimes b_n
\\
=&(b_1\otimes...\otimes b_{i\cdot 1}\otimes...\otimes b_n)...
(b_1\otimes...\otimes b_{i\cdot p}\otimes...\otimes b_n)
\omega
\end{split}
}
{tensors 1, representation, tensor product}

\DefEquation
{
b_1\otimes...\otimes(f_i(a)\circ b_i)\otimes...\otimes b_n=
f(a)\circ(b_1\otimes...\otimes b_i\otimes...\otimes b_n)
}
{tensors 2, representation, tensor product}

\DefEq
{
\[
f:\Times\rightarrow\Tensor B
\]
}
{map f, 1, representation, tensor product}

\DefEquation
{
f\circ(b_1,...,b_n)=\Tensor b
}
{map f, representation, tensor product}

\DefEq
{
\[
g:\Times\rightarrow V
\]
}
{map g, representation, tensor product}

\DefEq
{
\[
h:\Tensor B\rightarrow V
\]
}
{map h, representation, tensor product}

\AddEquation{Tower of Representations}
{
\xymatrix{
\\
A_{i+2}\ar@/_3pc/[rrrr]^{f_{i+1,i+2}(a_{i+1})}="f3a2"
\ar@/^3pc/[rrrr]^{f_{i+1,i+2}(\RedText{f_{i,i+1}(a_i)(a_{i+1})})}="f3b2"&&&&A_{i+2}
\\
\\
&A_{i+1}\ar[rr]^{f_{i,i+1}(a_i)}\ar @{}[rl]="a2"&&A_{i+1}\ar @{}[rl]="b2"
\\
\\
&&A_i\ar @{}[rl]="a1"\ar@{=>}[uu]^{f_{i,i+1}}&&
\ar @{=>}^{f_{i+1,i+2}} "a2";"f3a2"
\ar @{=>}@/_3pc/^{f_{i+1,i+2}} "b2";"f3b2"
\ar @{=>} "f3a2";"f3b2"_{f_{i,i+2}(a_i)}="a2b2"
\ar @{=>}@/^6pc/ "a1";"a2b2"^{f_{i,i+2}}
}
}

\AddEquation{morphism of tower of representations, level k}
{
h_{k+1}\circ\BlueText{f_{k,k+1}(a_k)}
=g_{k,k+1}(\RedText{h_k(a_k)})\circ h_{k+1}
}

\AddEquation{morphism of tower of representations, levels k k+1}
{
\BlueText{h_{k+1}(f_{k,k+1}(a_k)(a_{k+1}))}
=g_{k,k+1}(\RedText{h_k(a_k)})(\BlueText{h_{k+1}(a_{k+1})})
}

\AddEquation{morphism of tower of representations of F algebra, level k, diagram}
{
\xymatrix{
&A_{k+1}\ar[dd]_(.3){f_{k,k+1}(a_k)}
\ar[rr]^{h_{k+1}}
&&
B_{k+1}\ar[dd]^(.3){g_{k,k+1}(\RedText{h_k(a_k)})}\\
&\ar @{}[rr]|{(1)}&&\\
&A_{k+1}\ar[rr]^{h_{k+1}}&&B_{k+1}\\
A_k\ar[rr]^{h_k}\ar@{=>}[uur]^(.3){f_{k,k+1}}
&&
B_k\ar@{=>}[uur]^(.3){g_{k,k+1}}
}
}

\DefEq
{
\[
\xymatrix{
&\Tensor B\ar[dd]^h
\\
\Times
\ar[ru]^f\ar[rd]_g
\\
&V
}
\]
}
{map gh, representation, tensor product}

\DefEq
{
$A$, $B_1$, ..., $B_n$, $B$
}
{set of universal algebras 1}

\DefEq
{
\[
r_2:B_1\times...\times B_n\rightarrow B
\]
}
{reduced polymorphism of representation}

\DefEq
{
\symb{a\RCstar b}{rc-product}{}
}
{rc-product}

\DefEquation
{
\left\{\begin{array}{r@{\,=\,}l}
\ShowSymbol{rc-product}{}&
\begin{pmatrix}
a^i_kb^k_j
\end{pmatrix}
\\
(\ShowSymbol{rc-product}{}){}^i_j&
{}a^i_kb^k_j
\end{array}\right.
}
{rc-product of matrices}

\DefEq
{
\symb{a\CRstar b}{cr-product}{}
}
{cr-product}

\DefEquation
{
\left\{\begin{array}{r@{\,=\,}l}
\ShowSymbol{cr-product}{}&
\begin{pmatrix}
a^k_ib^j_k
\end{pmatrix}
\\
(\ShowSymbol{cr-product}{}){}^i_j&
{}a^k_ib^j_k
\end{array}\right.
}
{cr-product of matrices}

\DefEq
{
\[(a+b)^i_j=a^i_j+b^i_j\]
}
{(a+b)=}

\DefEquation
{
(a^T)^i_j=a^j_i
}
{transpose of matrix, 1}

\DefEquation
{
a\RCstar b
=(a_i^kb_k^j)
=(b_k^ja_i^k)
=b\CRstar a
}
{reducibility of products}

\AddEq{(m+n)a=ma+na,n=k}
{
\[(m+k)a=ma+ka\]
}

\AddEq{(m+n)a=ma+na,n=k+1}
{
\begin{align*}
(m+(k+1))a&=((m+k)+1)a=(m+k)a+a=ma+ka+a\\&=ma+(k+1)a
\end{align*}
}

\AddEq{(m+n)a=ma+na,n=k-1}
{
\begin{align*}
(m+(k-1))a&=((m+k)-1)a=(m+k)a-a=ma+ka-a\\&=ma+(k-1)a
\end{align*}
}

\AddEq{n(a+b)=na+nb,n=k}
{
\[k(a+b)=ka+kb\]
}

\AddEq{n(a+b)=na+nb,n=k+1}
{
\begin{align*}
(k+1)(a+b)&=k(a+b)+a+b=ka+kb+a+b\\&=ka+a+kb+b\\&=(k+1)a+(k+1)b
\end{align*}
}

\AddEq{n(a+b)=na+nb,n=k-1}
{
\begin{align*}
(k-1)(a+b)&=k(a+b)-(a+b)=ka+kb-a-b\\&=ka-a+kb-b\\&=(k-1)a+(k-1)b
\end{align*}
}

\AddEq{(nm)a=n(ma),n=k}
{
\[(km)a=k(ma)\]
}

\AddEq{(nm)a=n(ma),n=k+1}
{
\begin{align*}
((k+1)m)a&=(km+m)a=(km)a+ma=k(ma)+ma\\&=(k+1)(ma)
\end{align*}
}

\AddEq{(nm)a=n(ma),n=k-1}
{
\begin{align*}
((k-1)m)a&=(km-m)a=(km)a-ma=k(ma)-ma\\&=(k-1)(ma)
\end{align*}
}

\AddEq{f(na)=nf(a),n=k}
{
\[f(ka)=kf(a)\]
}

\AddEq{f(na)=nf(a),n=k+1}
{
\begin{align*}
f((k+1)a)&=f(ka+a)=f(ka)+f(a)=kf(a)+f(a)\\&=(k+1)f(a)
\end{align*}
}

\AddEq{f(na)=nf(a),n=k-1}
{
\begin{align*}
f((k-1)a)&=f(ka-a)=f(ka)-f(a)=kf(a)-f(a)\\&=(k-1)f(a)
\end{align*}
}

\AddEq{A1(mA2) symb}
{%
\symb{A_1(\mathcal A_2)}{category of representations}1
}

\AddEq{A1(mA2)}
{
$A_1(\mathcal A_2)$
}

\AddEq{product in category}
{
\symb[Pi-1]{\prod_{i\in I}B_i}{product in category}{}
\[P=\ShowSymbol{product in category}{}\]
}

\AddEq{coproduct in category}
{
\symb[Pi-1]{\coprod_{i\in I}B_i}{coproduct in category}{}
\[P=\ShowSymbol{coproduct in category}{}\]
}

\AddEq [2]{av=sum av}
{
\[#1^i#2_i=\sum_{i\in I}#1^i#2_i\]
}

\DefEq
{
\[
(a_1,a_2)\in A_1\times A_2\rightarrow a_1a_2\in A_2
\]
}
{A12->A2 left}

\DefEq
{
\[f(a)=fa\]
}
{f(a) left}

\DefEq
{
\[f(a)=af\]
}
{f(a) right}

\DefEquation
{
(fg)a=f(ga)
}
{fga= left}

\DefEquation
{
a(gf)=(ag)f
}
{fga= right}

\DefEq
{
\[fga=f(ga)=(fg)a\]
}
{fga= 1 left}

\DefEq
{
\[agf=(ag)f=a(gf)\]
}
{fga= 1 right}

\DefEquation
{
R(bc)=R(b)R(c)
}
{right shift, product}

\DefEquation
{
L(cb)=L(c)L(b)
}
{left shift, product}

\DefEq
{
\symb{R(b)}{right shift}{}
}
{right shift =}

\DefEquation
{
\begin{split}
a'&=a\ShowSymbol{right shift}{}=ab
\\
a'&=\ShowSymbol{right shift}{}(a)=ab
\end{split}
}
{right shift}

\DefEq
{
\symb{L(b)}{left shift}{}
}
{left shift =}

\DefEquation
{
\begin{split}
a'&=\ShowSymbol{left shift}{}a=ba
\\
a'&=\ShowSymbol{left shift}{}(a)=ba
\end{split}
}
{left shift}

\DefEq
{
\[
(x_1, ..., x_n)\in B_1\times...\times B_n
\rightarrow x_1\otimes...\otimes x_n\in B_1\otimes...\otimes B_n
\]
}
{B times->B otimes}

\DefEq
{
f:A_1\rightarrow\End(\Omega_2,A_2)
}
{representation of group, map}

\DefEq
{
\[
\begin{pmatrix}
\id:A_1\rightarrow A_1&r_2:A_2\rightarrow B_2
\end{pmatrix}
\]
}
{id:A1->A1 A2->B2}

\DefEquation
{
\xymatrix{
&A_2\ar[dd]_(.3){f(a)}\ar[rrr]^{r_2}&&&B_2\ar[dd]^(.3){g(a)}\\
&&&&\\
&A_2\ar[rrr]_(.65){r_2}&&&B_2\\
A_1\ar@{=>}[uur]^(.3)f\ar@{=>}[uurrrr]^(.7)g
}
}
{morphism id,R of representations}

\DefEquation
{
r_2\circ f(a)=g(a)\circ r_2
}
{morphism of representations of universal algebra}

\DefEq
{
\[
\xymatrix{
A_2\ar[rr]^{r_2}&&B_2\\
&A_1\ar[lu]|{*}^f\ar[ru]|{*}_g
}
\]
}
{morphism id,R of representations 2}

\AddEq[3]{f:A->*B}
{
\[
\xymatrix
{
#1:#2\ar[r]|{*}&#3
}
\]
}

\DefEq
{
a_2f(a_1b_1)=a_2f(a_1)f(b_1)
}
{right-side representation of group}

\DefEq
{
f(a_1b_1)a_2=f(a_1)f(b_1)a_2
}
{left-side representation of group}

\DefEq
{
\[
(a_2,a_1)\in A_2\times A_1\rightarrow a_2a_1\in A_2
\]
}
{A12->A2 right}

\DefEq
{
\[
\BlueText{r_2(a_1a_2)}=\RedText{r_1(a_1)}\BlueText{r_2(a_2)}
\]
}
{morphism of left representation of Omega group}

\DefEq
{
\[
\BlueText{r_2(a_2a_1)}=\BlueText{r_2(a_2)}\RedText{r_1(a_1)}
\]
}
{morphism of right representation of Omega group}

\DefEq
{
\symb{A_1a_2}{orbit of effective left-side representation}{}
\[
\ShowSymbol{orbit of effective left-side representation}{}
=L(A_1) a_2
\]
}
{orbit of effective left-side representation}

\DefEq
{
\symb{a_2A_1}{orbit of effective right-side representation}{}
\[
\ShowSymbol{orbit of effective right-side representation}{}
=a_2R(A_1)
\]
}
{orbit of effective right-side representation}

\DefEq
{
\symb{A_2/L(A_1)}{space of orbits of effective left-side representation}1
}
{space of orbits of effective left-side representation}

\DefEq
{
\symb{A_2/R(A_1)}{space of orbits of effective right-side representation}1
}
{space of orbits of effective right-side representation}

\DefEq
{
\begin{itemize}
\item
\EqParm{definition: side representation of group}{=left}
\item
\EqParm{definition: side representation of group}{=right}
\item
\ShowDefinition{representation of Abelian group}
\end{itemize}
}
{definition: representation of group}

\DefEq
{
\[\sum_{i=1}^{\infty}\|a^i\|<\infty\]
}
{sum |ai|<infty}

\DefEq
{
\[\sum_{i=1}^{\infty}a^i\]
}
{sum ai}

\AddEq{integral of map}
{
\symb{\int_X d\mu(x)f(x)}{Lebesgue integral}{}
}

\AddEquation{integral of simple map}
{
\ShowSymbol{Lebesgue integral}{}=
\sum_n \mu(F_n)f_n
}

\AddEquation{integral of measurable map}
{
\ShowSymbol{Lebesgue integral}{}=
\lim_{n\rightarrow\infty}\int_X d\mu(x)f_n(x)
}

\DefEq
{
\[\int_Xd\mu(x)f(x)\]
}
{int f X}

\DefEq
{
\[\int_Xd\mu(x)g(x)\]
}
{int g X}

\DefEq
{
\[\int_Xd\mu(x)(f(x)+g(x))\]
}
{int f+g X}

\DefEq
{
\int_Xd\mu(x)(f(x)+g(x))=
\int_Xd\mu(x)f(x)+
\int_Xd\mu(x)g(x)
}
{int f+g X =}

\DefEq
{
\[\int_Xd\mu(x)\|f(x)\|\]
}
{int |f| X}

\DefEq
{
\left\|\int_Xd\mu(x)f(x)\right\|\le\int_Xd\mu(x)\|f(x)\|
}
{|int f|<int|f|}

\DefEq
{
\|f(x)\|\le M
}
{|f(x)|<M}

\DefEq
{
\int_Xd\mu(x) \|f(x)\|\le M \mu(X)
}
{int |f|<M mu}

%% file: Stmt.Module.English.tex
\input{Stmt.Module.Eq}

\DefEq
{
\begin{definition}
\labelDefinition{module over ring}
Let commutative ring $D$ has unit $1$.
Effective representation
\DrawEq{D->*V}{def}
of ring $D$
in an Abelian group $V$
is called
\AddIndex{module over ring}{module over ring} $D$
or
\AddIndex{$D$\Hyph module}{D module}.
Effective representation
\eqRef{D->*V}{def}
of commutative ring $D$ in an Abelian group $V$
is called
\AddIndex{vector space over field}{vector space over field} $D$
or
\AddIndex{$D$\Hyph vector space}{D vector space}.
\qed
\end{definition}
}
{... definition: module over ring}

\DefTheorem{tensor product of D-algebras is D-algebra}
{
Let
\ShowEq{A1...n}
be $D$\Hyph algebras.
Tensor product
$\Tensor A$
of $D$\Hyph modules
\ShowEq{A1...n}
is $D$\Hyph algebra,
if we define product by the equality
\ShowEq{xA1n*xA1n->oxA1n=}
}

\DefTheorem{representation of algebra A2 in LA}
{
Let $A$ be $D$\Hyph algebra.
Let product in algebra
\ShowEq{AoxA}
be defined according to rule
\DrawEq{product in algebra AA}{}
A representation
\ShowEq{h:AoxA->L(A)}
of $D$\Hyph algebra
\ShowEq{AoxA}
in module
\ShowEq{L(A;B)}DAA{}
defined by the equality
\ShowEq{representation AA in LA}
allows us to identify tensor
\ShowEq{d in AxoA}
and linear map
\ShowEq{product in algebra AA 2}
where
\ShowEq{product in algebra AA 3}
is identity map.
Linear map
\ShowEq{a ox b}
has form
\ShowEq{a ox b c=}
}

\DefDefinition{module over commutative ring}
{
Effective representation of commutative ring $D$
in an Abelian group $V$
\DrawEq{D->*V}{module}
is called
\AddIndex{module over ring}{module over ring} $D$
or
\AddIndex{$D$\Hyph module}{D module}.
$V$\Hyph number is called
\AddIndex{vector}{vector}.
\ePrints{1610.309618526,CACAA.06}
\ifx\Semafor\ValueOn
If $D$ is a field, then the Abelian group $V$
is called
\AddIndex{vector space}{vector space} over field $D$
or
\AddIndex{$D$\Hyph vector space}{D vector space}.
\fi
}

\DefTheorem{definition of module over commutative ring}
{
Following conditions hold for $D$\Hyph module:
\StartLabelItem
\begin{enumerate}
\item
\AddIndex{associative law}{associative law}
\DrawEq{associative law, module}2
\item 
\AddIndex{distributive law}{distributive law}
\ShowEq{distributive law, module}
\item
\AddIndex{unitarity law}{unitarity law}
\ShowEq{unitarity law, D-module}
\end{enumerate}
for any
\ShowEq{c,d in Z, a,b in D, v,w in V}
}

\AddEq{definition: module over algebra}
{
\begin{definition}
\labelDefinition{\SideWS module over algebra}
Effective \SideNS\Hyph side representation
\DrawEq{\SideWS A->*V}{def}
of associative $D$\Hyph alebra $A$
in $D$\Hyph module $V$
is called
\AddIndex{\SideWS module}{\SideWS module} over $D$\Hyph alebra $A$.
We will also say that $D$\Hyph module $V$ is
\AddIndex{\SideWS $A$\Hyph module}{\SideWS A module}
or
\AddIndex{$\SideA$\Hyph module}{\SideA-module}.
$V$\Hyph number is called
\AddIndex{vector}{vector}.
If $A$ is division algebra, then \SideWS $A$\Hyph module
is called
\AddIndex{\SideWS vector space}{\SideWS vector space}
over $D$\Hyph alebra $A$.
We will also say that $D$\Hyph module $V$ is
\AddIndex{\SideWS $A$\Hyph vector space}{\SideWS A vector space}
or
\AddIndex{$\SideA$\Hyph vector space}{\SideA-vector space}.
\qed
\end{definition}
}

\AddEq{theorem: module over algebra}
{
\begin{theorem}
\labelTheorem{\SideWS module over algebra}
The following diagram of representations describes \SideWS $\Base$\Hyph module $V$
\ShowEq{diagram of representations, \SideWS module}
The diagram of representations
\EqRef{diagram of representations, \SideWS module}
holds
\AddIndex{commutativity of representations}{commutativity of representations}
\BaseRings
in абелевой группе $V$
\ShowEq{\SideWS module, a d v}
\end{theorem}
}

\DefProof{module over algebra}
{
The diagram of representations
\EqRef{diagram of representations, \SideWS module}
follows from the definition
\def\Temp{}
\ifx\SideNS\Temp
\RefDefinition{module over commutative ring}
and from the theorem
\RefTheorem{action of ring of rational integers in Abelian group}.
\else
\RefDefinition{\SideWS module over algebra}
and the corollary
\RefCorollary{Free Algebra over Ring}.
\fi
Since \SideNS\HSide transformation $\ATransf(a)$
is endomorphism
of $\CBase$\Hyph module $V$,
we obtain the equality
\EqRef{\SideWS module, a d v}.
}

\AddEq{theorem: A module -> alebra is associative}
{
\begin{theorem}
\labelTheorem{\SideWS A module -> alebra is associative}
Let $g$ be effective left\Hyph side representation of $D$\Hyph algebra $A$
in $D$\Hyph module $V$.
Then $D$\Hyph algebra $A$ is associative.
\end{theorem}
}

\DefProof{A module -> alebra is associative}
{
Let
\ShowEq{abc in A, v in v}
Since \SideNS\Hyph side representation
$g$ is \SideNS\Hyph side representation
of the multiplicative group
of $D$\Hyph algebra $A$,
we obtain the equality
\DrawEq{\SideWS module, associative law}0
The equality
\ShowEq{a(b(cv))=(a(bc))v \SideNS}
follows from the equality
\eqRef{\SideWS module, associative law}0.
Since
\ShowEq{cv \SideNS}
the equality
\ShowEq{a(b(cv))=((ab)c)v \SideNS}
follows from the equality
\eqRef{\SideWS module, associative law}0.
The equality
\ShowEq{(a(bc))v=((ab)c)v \SideNS}
follows from equalities
\EqRef{a(b(cv))=(a(bc))v \SideNS},
\EqRef{(a(bc))v=((ab)c)v \SideNS}.
Since $v$ is any vector of $A$\Hyph module $V$,
the equality
\ShowEq{a(bc)=(ab)c \SideNS}
follows from the equality
\EqRef{(a(bc))v=((ab)c)v \SideNS}.
Therefore, $D$\Hyph algebra $A$ is associative.
}

\AddEq{theorem: definition of A module}
{
\begin{theorem}
\labelTheorem{definition of \SideWS module}
Let $V$ be \SideWS $\Base$\Hyph module.
For any vector $v\in V$,
vector generated by the diagram of representations
\EqRef{diagram of representations, \SideWS module}
has the following form
\ShowEq{\SideWS module, (a+n)v=av+nv}
\StartLabelItem
\begin{enumerate}
\item
The set of maps
\ShowEq{\SideWS module, a+n:V->V}
generates\,\footnote{
See similar statement on the page
\citeBib{McCrimmon: Jordan Algebras}-52.
}
\algebraa $\Base_{(1)}$
where the sum is defined by the equality
\DrawEq{(a+n)+(b+m)=}{\SideWS module}
and the product is defined by the equality
\DrawEq{(a+n)(b+m)=}{\SideWS module}
The \algebraa $\Base_{(1)}$ is called
\AddIndex{unital extension}{unital extension}
of the \algebraa $\Base$.
\item
The \algebraa $\Base$ is \SideWS ideal of \algebraa $\Base_{(1)}$.
\labelItem{Algebra is \SideWS ideal of algebra (1)}
\item
The set of transormations
\EqRef{\SideWS module, (a+n)v=av+nv}
is \SideNS\HSide representation of \algebraa $\Base_{(1)}$ in Abelian group $V$.
\labelItem{\SideWS representation of D1 in Abelian group}
\end{enumerate}
Following conditions hold for \SideWS $\Base$\Hyph module $V$:
\begin{itemize}
\item 
\AddIndex{associative law}{associative law}
\DrawEq{associative law, \SideWS module}1
\item 
\AddIndex{distributive law}{distributive law}
\ShowEq{distributive law, \SideWS module}
\item
\AddIndex{unitarity law}{unitarity law}
\ShowEq{unitarity law, \SideWS \Base-module}
\end{itemize}
for any
\ShowEq{p,q in \Base 1, v,w in V}
\end{theorem}
}

\DefProof{definition of A module}
{
Let $v\in V$.

\begin{lemma}
\labelLemma{\SideWS module, map a+n:V->V is endomorphism of Abelian group}
{\it
Let
\ShowEq{module, d in D}
\ShowEq{module, a in A}
The map
\EqRef{\SideWS module, a+n:V->V}
is endomorphism of Abelian group $V$.
}
\end{lemma}

{\sc Proof.}
Statements
\ShowEq{\SideWS module, dv in V}
\ShowEq{\SideWS module, av in V}
follow from the theorems
\RefTheorem[\RefRepresentation]{structure of subrepresentations},
\RefTheorem{\SideWS module over algebra}.
Since $V$ is Abelian group, then
\ShowEq{\SideWS module, dv+av in V}
Therefore,
for any $\CBase$\Hyph number $\DArg$
and for any $\Base$\Hyph number $a$,
we defined the map
\EqRef{\SideWS module, a+n:V->V}.
Since transformation $\DTransf(\DArg)$
and \SideNS\HSide transformation $\ATransf(a)$
are endomorphisms of Abelian group $V$,
then the map
\EqRef{\SideWS module, a+n:V->V}
is endomorphism of Abelian group $V$.
\hfill\(\odot\)

Let $\Base_{(1)}$ be the set of maps
\EqRef{\SideWS module, a+n:V->V}.
The equality
\EqRef{distributive law, \SideWS module, 1}
follows from the lemma
\RefLemma{\SideWS module, map a+n:V->V is endomorphism of Abelian group}.

Let
\ShowEq{p,q in D1}
According to the statement
\RefItem{representation of D1 in Abelian group},
we define the sum of $\Base_{(1)}$\Hyph numbers $p$ and $q$ by the equality
\EqRef{distributive law, \SideWS module, 2}.
The equality
\ShowEq{\SideWS module, (a+n)+(b+m)=, 1}
follows from the equality
\EqRef{distributive law, \SideWS module, 2}.
Since representation
$\DTransf$ is homomorphism of the aditive group of ring $\CBase$,
we obtain the equality
\ShowEq{distributive law, \CBase, \SideWS module, 2}
Since \SideNS\HSide representation
$\ATransf$ is homomorphism of the aditive group of \algebraa $\Base$,
we obtain the equality
\ShowEq{distributive law, \Base, \SideWS module, 2}
Since $V$ is Abelian group, then the equality
\ShowEq{\SideWS module, (a+n)+(b+m)=}
follows from equalities
\EqRef{\SideWS module, (a+n)+(b+m)=, 1},
\EqRef{distributive law, \CBase, \SideWS module, 2},
\EqRef{distributive law, \Base, \SideWS module, 2}.
From the equality
\EqRef{\SideWS module, (a+n)+(b+m)=},
it follows that the definition
\eqRef{(a+n)+(b+m)=}{\SideWS module}
of sum on the set $\Base_{(1)}$ does not depend on vector $v$.

Equalities
\eqRef{associative law, \SideWS module}1,
\EqRef{unitarity law, \SideWS \Base-module}
follow from the statement
\RefItem{\SideWS representation of D1 in Abelian group}.
Let
\ShowEq{p,q in D1}
\def\Temp{}
\ifx\SideNS\Temp
Since representation $\DTransf$ is representation
of the multiplicative group of ring $\CBase$,
we obtain the equality
\DrawEq{associative law, \CBase, \SideWS module}1
Since representation $g_2$ is representation
of the multiplicative group of ring $D$,
we obtain the equality
\DrawEq{\SideWS module, associative law}1
Since the ring $D$
is Abelian group,
we obtain the equality
\DrawEq{associative law, \CBase\Base, \SideWS module}1
The equality
\ShowEq{module, (a+n)(b+m)=}
follows from equalities
\ShowEq{ref module, (a+n)(b+m)=}
The equality
\eqRef{(a+n)(b+m)=}{\SideWS module}
follows from the equality
\EqRef{module, (a+n)(b+m)=}.
\else%
Так как произведение в $D$\Hyph алгебре $A$ может быть неассоциативным,
то опираясь на теорему
мы рассмотрим произведение
of $\Base_{(1)}$\Hyph numbers $p$ and $q$
как билинейное отображение
\ShowEq{f:D1xD1->D1}
such that following equalities are true
\DrawEq{fab=ab}{\SideWS module}
\DrawEq{f1p=p}{\SideWS module}
The equality
\DrawEq{pq=fpq}{\SideWS module}
follows from equalities
\eqRef{fab=ab}{\SideWS module},
\eqRef{f1p=p}{\SideWS module}.
The equality
\eqRef{(a+n)(b+m)=}{\SideWS module}
follows from the equality
\eqRef{pq=fpq}{\SideWS module}.
\fi%

The statement
\RefItem{Algebra is \SideWS ideal of algebra (1)}
follows from the equality
\eqRef{(a+n)(b+m)=}{\SideWS module}.
}

\AddEq{definition: submodule}
{
\begin{definition}
\labelDefinition{submodule, \SideWS module}
Subrepresentation of \SideWS $D$\Hyph module $V$ is called
\AddIndex{submodule}{submodule}
of $D$\Hyph module $V$.
\qed
\end{definition}
}

\AddEq{theorem: submodule}
{
\begin{theorem}
\labelTheorem{submodule, \SideWS module}
Let
\ShowEq{vi V}{}
be set of vectors of \SideWS $\Base$\Hyph module $V$.
If vectors
\ShowEq{set vi}
belongs submodule $V'$ of \SideWS $\Base$\Hyph module $V$,
then linear combination of vectors
\ShowEq{set vi}
belongs submodule $V'$.
\end{theorem}
}

\DefProof{submodule}
{
The theorem follows from
the theorem
\RefTheorem{set of vectors generated by set of vectors, \SideWS module}
and definitions
\RefDefinition{\SideWS linear combination of vectors},
\RefDefinition{submodule, \SideWS module}.
}

\AddEq{theorem: set of vectors generated by set of vectors}
{
\begin{theorem}
\labelTheorem{set of vectors generated by set of vectors, \SideWS module}
Let $V$ be \SideWS $\Base$\Hyph module.
The set of vectors generated by the set of vectors
\ShowEq{vi V}{}
has form\,\footnote{
For a set $A$,
we denote by $|A|$ the cardinal number of the set $A$.
The notation $|A|<\infty$ means that the set $A$ is finite.
}
\ShowEq{w=sum vi, \SideWS module}
\end{theorem}
}

\DefProof{set of vectors generated by set of vectors}
{
Acording to the theorem
\RefTheorem[\RefRepresentation]{structure of subrepresentations},
to prove the theorem, we need to prove following statements:
\ShowEq{vectors generated by set of vectors}

\begin{itemize}
\item
For any
\ShowEq{vk in v}
let
\ShowEq{ci=dik}
Then
\ShowEq{vk=sum vi, \SideWS module}
The statement
\RefItem{\SideWS module, vk in Jv}
follows from
\EqRef{w=sum vi, \SideWS module},
\EqRef{vk=sum vi, \SideWS module}.
\item
The statement
\RefItem{\SideWS module, cvk in Jv}
follow from the theorems
\RefTheorem[\RefRepresentation]{structure of subrepresentations},
\RefTheorem{definition of \SideWS module}
and from the statement
\RefItem{\SideWS module, vk in Jv}.
\item
Since $V$ is Abelian group,
then the statement
\RefItem{\SideWS module, sum cvk in Jv}
follows from the statement
\RefItem{\SideWS module, cvk in Jv}
and from theorems
\RefTheorem[\RefRepresentation]{structure of subrepresentations},
\RefTheorem{structure of Abelian group}.
\item
Let
\ShowEq{w12 in Jv}wv
Since $V$ is Abelian group,
then, according to the statement
\RefItem[\RefRepresentation]{x1n omega in Xk+1},
\ShowEq{w1+w2 in Jv}wv
According to the equality
\EqRef{w=sum vi, \SideWS module},
there exist $\Base_{(1)}$\Hyph numbers
\ShowEq{gi12}w
such that
\DrawEq[w]{w12= \SideNS}{module}
where sets
\DrawEq[w]{|ci12 ne 0|}{\SideWS module}
are finite.
Since $V$ is Abelian group,
then from the equality
\eqRef{w12= \SideNS}{module}
it follows that
\DrawEq[w]{w1+w2= \SideNS}{module}
The equality
\DrawEq[w]{w1+w2= 1 \SideNS}{module}
follows from equalities
\EqRef{distributive law, \SideWS module, 2},
\eqRef{w1+w2= \SideNS}{module}.
From the equality
\eqRef{|ci12 ne 0|}{\SideWS module},
it follows that
the set
\ShowEq{|gi 1+2 ne 0|}w
is finite.
From the equality
\eqRef{w1+w2= 1 \SideNS}{module},
it follows that
\ShowEq{w1+w2 in Jv}wv
\item
Let
\ShowEq{w in Jv}
According to the statement
\RefItem[\RefRepresentation]{ax in Xk+1},
for any $\Base_{(1)}$\Hyph number $a$,
\ShowEq{aw in Xk \SideWS module}
According to the equality
\EqRef{w=sum vi, \SideWS module},
there exist $\Base_{(1)}$\Hyph numbers
\ShowEq{wi}
such that
\ShowEq{w= \SideWS module}
where
\DrawEq{|wi ne 0|}{\SideWS module}
From the equality
\EqRef{w= \SideWS module}
it follows that
\ShowEq{aw= \SideWS module}
From the statement
\eqRef{|wi ne 0|}{\SideWS module},
it follows that
the set
\ShowEq{\SideWS module, |awi ne 0|}
is finite.
From the equality
\EqRef{aw= \SideWS module},
it follows that
\ShowEq{aw in Jv \SideWS module}
\end{itemize}
}

\AddEq{definition: linear combination of vectors}
{
\begin{definition}
\labelDefinition{\SideWS linear combination of vectors}
Let
\ShowEq{vi V}{}
be set of vectors.
\ShowEq{\SideWS linear combination}%
The expression
\ShowEq{A linear combination =}
is called
\AddIndex{linear combination}{linear combination} of vectors
\ShowEq{vi}
A vector
\ShowEq{w=wi vi\SideNS}
is called
\AddIndex{linearly dependent}{linearly dependent}
on vectors
\ShowEq{vi}
\qed
\end{definition}
}

\AddEq{definition: generating set of module}
{
\begin{definition}
\labelDefinition{generating set of \SideWS module}
$J(v)$
is called
\AddIndex{submodule generated by set}
{submodule generated by set} $v$,
and $v$ is a
\AddIndex{generating set}{generating set}
of submodule $J(v)$.
In particular, a
\AddIndex{generating set}{generating set}
of \SideWS $D$\Hyph module $V$
is a subset $X\subset V$ such that
\ShowEq{generating set of module}
\qed
\end{definition}
}

\AddEq{remark: generating set of module}
{
The following definition follows from the theorems
\RefTheorem{set of vectors generated by set of vectors, \SideWS module},
\RefTheorem[\RefRepresentation]{structure of subrepresentations}
and from the definition
\RefDefinition[\RefRepresentation]{generating set of representation}.
}

\AddEq{remark: basis of module}
{
The following definition follows from the theorems
\RefTheorem{set of vectors generated by set of vectors, \SideWS module},
\RefTheorem[\RefRepresentation]{structure of subrepresentations}
and from the definition
\RefDefinition[\RefRepresentation]{basis of representation}.
}

\AddEq{definition: basis of module}
{
\begin{definition}
\labelDefinition{basis of \SideWS module}
If the set $X\subset V$ is generating set of \SideWS $D$\Hyph module
$V$, then any set $Y$, $X\subset Y\subset V$
also is generating set of \SideWS $D$\Hyph module $V$.
If there exists minimal set $X$ generating
the \SideWS $D$\Hyph module $V$, then the set $X$ is called
\AddIndex{basis of \SideWS $D$\Hyph module}{basis, module} $V$.
\qed
\end{definition}
}

\AddEq{theorem: linearly depends on rest of vectors}
{
\begin{theorem}
\labelTheorem{linearly depends on rest of vectors, \SideWS module}
Let $\Base$ be \Algebra.
Since the equation
\ShowEq{\SideWS wi vi=0}
implies existence of index
\ShowEq{i=j}
such that
\ShowEq{wj ne 0},
then the vector $v_{\gij}$
linearly depends on rest of vectors $v$.
\end{theorem}
}

\AddEq{proof: linearly depends on rest of vectors}
{
\begin{proof}
The theorem follows from the equality
\ShowEq{\SideWS vj=sum vi}v
and from the definition
\RefDefinition{\SideWS linear combination of vectors}.
\end{proof}
}

\AddEq{remark: 0=0vi}
{
It is evident that for any set of vectors $v_{\gii}$
\ShowEq{\SideWS 0=0vi}
}

\AddEq{remark: scalars and vectors as matrix}
{
We represent the set of $\Base_{(1)}$\Hyph numbers
\ShowEq{wi}
as matrix
\ShowEq{w=wi}
We represent the set of vectors
\ShowEq{set vi}
as matrix
\ShowEq{v=vi}
Then we can represent linear combination of vectors
\ShowEq{w=wi vi\SideNS}
as
\ShowEq{w=w cr v}
}

\AddEq{definition: linearly independent vectors}
{
\begin{definition}
\labelDefinition{\SideWS linearly independent vectors}
{\it
The set of vectors\,\footnote{
I follow to the definition in
\citeBib{Serge Lang}, page 130.}
\ShowEq{set vi}
of \SideWS $\Base$\Hyph module $V$ is
\AddIndex{linearly independent}{linearly independent set}
if $w=0$ follows from the equation
\ShowEq{\SideWS wi vi=0}
Otherwise the set of vectors
\ShowEq{set vi}
is \AddIndex{linearly dependent}{linearly dependent set}.
}
\qed
\end{definition}
}

\AddEq{theorem: basis of module}
{
\begin{theorem}
\labelTheorem{basis of \SideWS module}
The set of vectors
\ShowEq{basis for module over algebra}
is basis of \SideWS $\Base$\Hyph module
$V$, if following statements are true.
\StartLabelItem
\begin{enumerate}
\item
\labelItem{vector is linear combination of set, \SideWS module}
Arbitrary vector $v\in V$
is linear combination of
vectors of the set $\Basis e$.
\item
\labelItem{cannot be represented as a linear combination, \SideWS module}
Vector $e_{\gii}$
cannot be represented as a linear combination
of the remaining vectors of the set $\Basis e$.
\end{enumerate}
\end{theorem}
}

\DefProof{basis of module}
{
According to the statement
\RefItem{vector is linear combination of set, \SideWS module},
the theorem
\RefTheorem{set of vectors generated by set of vectors, \SideWS module}
and the definition
\RefDefinition{\SideWS linear combination of vectors},
the set $\Basis e$ generates \SideWS $\Base$\Hyph module $V$
(the definition
\RefDefinition{generating set of \SideWS module}).
According to the statement
\RefItem{cannot be represented as a linear combination, \SideWS module},
the set $\Basis e$ is minimal set
generating \SideWS $\Base$\Hyph module $V$.
According to the definitions
\RefDefinition{basis of \SideWS module},
the set $\Basis e$ is a basis of \SideWS $\Base$\Hyph module $V$.
}

\AddEq{theorem: division algebra, basis}
{
\begin{theorem}
\labelTheorem{division algebra, \SideWS basis}
Let $\Base$ be \Algebra.
The set of vectors
\ShowEq{basis, module}
is a
\AddIndex{basis of \SideWS $\Base$\Hyph vector space}{basis, vector space} $V$
if vectors $e_{\gii}$ are linearly independent and any vector $v\in V$
linearly depends on vectors $e_{\gii}$.
\end{theorem}
}

\DefProof{division algebra, basis}
{
Let the set of vectors
\ShowEq{ei, i in I}
be linear dependent. Then the equation
\DrawEq[w]{c*e=0, \SideWS module}{}
implies existence of index $\gii=\gij$ such that
\ShowEq{wj ne 0}.
According to the theorem
\RefTheorem{linearly depends on rest of vectors, \SideWS module},
the vector $e_{\gij}$
linearly depends on rest of vectors of the set $\Basis e$.
According to the definition
\RefDefinition{basis of \SideWS module},
the set of vectors
\ShowEq{ei, i in I}
is not a basis for \SideWS $\Base$\Hyph vector space $V$.

Therefore, if the set of vectors
\ShowEq{ei, i in I}
is a basis, then these vectors
are linearly independent.
Since an arbitrary vector $v\in V$
is linear combination of vectors
\ShowEq{ei, i in I},
then the set of vectors $v$,
\ShowEq{ei, i in I}
is not linearly independent.
}

\AddEq{convention: we use separate color for index of element}
{
\begin{convention}
\labelConvention{we use separate color for index of element}
Let $A$ be free algebra
with finite or countable basis.
Considering expansion of element of algebra $A$ relative basis $\Basis e$
we use the same root letter to denote this element and its coordinates.
In expression $a^2$, it is not clear whether this is component
of expansion of element
$a$ relative basis, or this is operation $a^2=aa$.
To make text clearer we use separate color for index of element
of algebra. For instance,
\ShowEq{Expansion relative basis in algebra}
\qed
\end{convention}
}

\AddEq{convention: unit of algebra in basis}
{
\begin{convention}
Let $\Basis e$ be the basis of free algebra $A$ over ring $D$.
If algebra $A$ has unit,
then we assume that $e_{\gi 0}$ is the unit of algebra $A$.
\qed
\end{convention}
}

\DefDefinition{linear map from D1 A1 to D2 A2, module}
{
Morphism of representations
\ShowEq{h:D1->D2 f:A1->A2}h\Base fV
of $D_1$\Hyph module $A_1$
into $D_2$\Hyph module $A_2$
is called
\AddIndex{linear map}{linear map}
of $D_1$\Hyph module $A_1$ into $D_2$\Hyph module $A_2$.
Let us denote
\ShowEq{set linear maps, D12 module}
set of linear maps
of $D_1$\Hyph module $A_1$ into $D_2$\Hyph module $A_2$.
}

\ifx\texFuture\Defined
\AddEq{linear map from D1 A1 to D2 A2, module}
{
\begin{definition}
\labelDefinition{linear map from D1 A1 to D2 A2, \SideNS module}
{\it
Morphism of representations
\ShowEq{h:D1->D2 f:A1->A2}h\Base fV
of \SideWS $D_1$\Hyph module $A_1$
into \SideWS $D_2$\Hyph module $A_2$
is called
\AddIndex{linear map}{linear map}
of \SideWS $D_1$\Hyph module $A_1$ into \SideWS $D_2$\Hyph module $A_2$.

Let us denote
\ShowEq{set linear maps, D12 module}
set of linear maps
of \SideWS $D_1$\Hyph module $A_1$ into \SideWS $D_2$\Hyph module $A_2$.
}
\qed
\end{definition}
}
\fi

\DefEq
{
\begin{definition}
\labelDefinition{linear map from A1 to A2, commutative module}
Let $A_1$, $A_2$ be $D$\Hyph modules.
Morphism of representations
\ShowEq{f:A->B}f{A_1}{A_2}
of $D$\Hyph module $A_1$
into $D$\Hyph module $A_2$
is called
\AddIndex{linear map}{linear map}
of $D$\Hyph module $A_1$ into $D$\Hyph module $A_2$.
Let us denote
\ShowEq{set linear maps, module}
set of linear maps
of $D$\Hyph module $A_1$ into $D$\Hyph module $A_2$.
\qed
\end{definition}
}
{definition: linear map from A1 to A2, commutative module}

\DefEq
{
\begin{theorem}
\labelTheorem{linear map from A1 to A2, commutative module}
Linear map
\ShowEq{f:A->B}f{A_1}{A_2}
of $D$\Hyph module $A_1$
into $D$\Hyph module $A_2$
satisfies to equalities\,\footnote{In some books
(for instance, \citeBib{Serge Lang}, p. 119) the theorem \RefTheorem{linear map from A1 to A2, commutative module}
is considered as a definition.}
\ShowEq{linear map from A1 to A2, 1}
\end{theorem}
}
{theorem: linear map from A1 to A2, commutative module}

\DefEq
{
\begin{theorem}
\labelTheorem{set of A->B is D module}
Let $A$ be Banach $D$\Hyph module with norm $|x|_A$.
Let $B$ be Banach $D$\Hyph module with norm $|y|_B$.
\StartLabelItem
\begin{enumerate}
\item
The set
$B^A$
of maps
\ShowEq{f:A->B}fAB
is $D$\Hyph module.
\labelItem{set of A->B is D module}
\item
The map
\ShowEq{f in BA->|f|}
defined by the equality
\ShowEq{norm of map}
\ShowEq{norm of map, algebra}
is the norm in $D$\Hyph module $B^A$
and the value
\ShowEq{show|f|}
is called
\AddIndex{norm of map $f$}{norm of map}.
\labelItem{norm of map}
\end{enumerate}
\end{theorem}
}
{theorem: set of A->B is D module}

\DefDefinition{coordinates of vector 2016}
{
Let $\Basis e$ be the basis of \SideWS $\Base$\Hyph module $V$
and vector
\ShowEq{vv in V}
has expansion
\DrawEq{vv=ve \SideWS module}{}
with respect to the basis $\Basis e$.
$\Base$\Hyph numbers $v^{\gii}$ are called
\AddIndex{coordinates}{coordinates}
of vector $\Vector v$ with respect to the basis $\Basis e$.
\ePrints{1610.309618526,CACAA.06}
\ifx\Semafor\ValueOn
Matrix of $\Base$\Hyph numbers
\ShowEq{coordinate matrix of vector}
is called
\AddIndex{coordinate matrix of vector}{coordinate matrix of vector}
$\Vector v$ in basis $\Basis e$.
\fi
}

\AddEq{definition: coordinates of vector}
{
\begin{definition}
\labelDefinition{coordinates of vector, \SideWS module}
{\it
Let $\Basis e$ be the basis of \SideWS $\Base$\Hyph module $V$
and vector
\ShowEq{vv in V}
has expansion
\DrawEq{vv=ve \SideWS module}{}
with respect to the basis $\Basis e$.
$\Base_{(1)}$\Hyph numbers $v^{\gii}$ are called
\AddIndex{coordinates}{coordinates}
of vector $\Vector v$ with respect to the basis $\Basis e$.
Matrix of $\Base_{(1)}$\Hyph numbers
\ShowEq{coordinate matrix of vector}
is called
\AddIndex{coordinate matrix of vector}{coordinate matrix of vector}
$\Vector v$ in basis $\Basis e$.
}
\qed
\end{definition}
}

\AddEq{theorem: linear dependence between vectors of basis}
{
\begin{theorem}
\labelTheorem{linear dependence between vectors of basis, \SideWS module}
Let $\Base$ be \algebra.
Let $\Basis e$ be basis of \SideWS $\Base$\Hyph module $V$.
Let
\DrawEq[w]{c*e=0, \SideWS module}{1}
be linear dependence of vectors of the basis $\Basis e$.
Then
\StartLabelItem
\begin{enumerate}
\item
$\Base_{(1)}$\Hyph number
\ShowEq{wi i in I}
does not have inverse element
in \algebraa $\Base_{(1)}$.
\item
The set $\Base'$ of matrices
\ShowEq{w=wi i in I}
generates \SideWS $\Base$\Hyph module.
\end{enumerate}
\end{theorem}
}

\DefProof{linear dependence between vectors of basis}
{
Let there exist matrix
\ShowEq{w=wi i in I}
such that the equality
\eqRef{c*e=0, \SideWS module}{1}
is true and there exist index
\ShowEq{i=j}
such that
\ShowEq{wj ne 0}.
If we assume that $\Base_{(1)}$\Hyph number $c^{\gij}$
has inverse one, then the equality
\ShowEq{\SideWS vj=sum vi}e
follows from the equality
\eqRef{c*e=0, \SideWS module}{1}.
Therefore, the vector $e_{\gij}$
is linear combination of other vectors of the set $\Basis e$
and the set $\Basis e$ is not basis.
Therefore, our assumption is false,
and $\Base_{(1)}$\Hyph number $c^{\gij}$ does not have inverse.

Let matrices
\ShowEq{b in D'}b,
\ShowEq{b in D'}c.
From equalities
\DrawEq[b]{c*e=0, \SideWS module}{}
\DrawEq[c]{c*e=0, \SideWS module}{}
it follows that
\ShowEq{(b+c)*e, \SideWS module}
Therefore, the set $\Base'$ is Abelian group.

Let matrix
\ShowEq{b in D'}c{}
and $a\in\Base$.
From the equality
\DrawEq[w]{c*e=0, \SideWS module}{}
it follows that
\ShowEq{(ac)*e, \SideWS module}
Therefore, Abelian group $\Base'$ is \SideWS $\Base$\Hyph module.
}

\AddEq{theorem: coordinates of vector with linear dependence}
{
\begin{theorem}
\labelTheorem{coordinates of vector with linear dependence, \SideWS module}
Let \SideWS $\Base$\Hyph module $V$
have the basis $\Basis e$ such that in the equality
\DrawEq[w]{c*e=0, \SideWS module}{2}
there exists of index
\ShowEq{i=j}
such that
\ShowEq{wj ne 0}.
Then
\StartLabelItem
\begin{enumerate}
\item
The matrix
\ShowEq{w=wi i in I}
determines coordinates of vector $0\in V$ with respect to basis $\Basis e$.
\labelItem{coordinates of vector 0 with linear dependence, \SideWS module}
\item
Coordinates of vector $\Vector v$ with respect to basis $\Basis e$
are uniquely determined up to a choice of coordinates of vector $0\in V$.
\labelItem{coordinates of vector with linear dependence, \SideWS module}
\end{enumerate}
\end{theorem}
}

\DefProof{coordinates of vector with linear dependence}
{
The statement
\RefItem{coordinates of vector 0 with linear dependence, \SideWS module}
follows from the equality
\eqRef{c*e=0, \SideWS module}{2}
and from the definition
\RefDefinition{coordinates of vector, \SideWS module}.

Let vector $\Vector v$ have expansion
\DrawEq{vv=ve \SideWS module}{2}
with respect to basis $\Basis e$.
The equality
\ShowEq{v=v+0, \SideWS module}
follows from equalities
\eqRef{c*e=0, \SideWS module}{2},
\eqRef{vv=ve \SideWS module}{2}.
The statement
\RefItem{coordinates of vector with linear dependence, \SideWS module}
follows from equalities
\eqRef{vv=ve \SideWS module}{2},
\EqRef{v=v+0, \SideWS module}
and from the definition
\RefDefinition{coordinates of vector, \SideWS module}.
}

\AddEq{definition: free module over ring}
{
\begin{definition}
\labelDefinition{free \SideWS module}
The \SideWS $\Base$\Hyph module $V$ is
\AddIndex{free \SideWS $\Base$\Hyph module}{free module},\,\footnote{
I follow to the
definition in \citeBib{Serge Lang}, page 135.}
if \SideWS $\Base$\Hyph module $V$ has basis
and vectors of the basis are linearly independent.
\qed
\end{definition}
}

\AddEq{theorem: coordinates of vector of free module}
{
\begin{theorem}
\labelTheorem{coordinates of vector of free \SideWS module}
Coordinates of vector $v\in V$ relative to basis $\Basis e$
of \SideWS $\Base$\Hyph module $V$
are uniquely defined.
\end{theorem}
}

\DefProof{coordinates of vector of free module}
{
The theorem follows from the theorem
\RefTheorem{coordinates of vector with linear dependence, \SideWS module}
and from definitions
\RefDefinition{\SideWS linearly independent vectors},
\RefDefinition{free \SideWS module}.
}

\DefTheorem{standard representation of map A1 A2, associative algebra}
{
Let $A_1$ be free $D$\Hyph module.
Let $A_2$ be free finite dimensional associative $D$\Hyph algebra.
Let $\Basis e$ be basis of $D$ module $A_2$.
Let $\Basis F$
be the basis of left \BoxB{A_2}module
\ShowEq{L(A;B)}D{A_1}{A_2}.\,\footnote{
If $D$\Hyph module $A_1$ or $D$\Hyph module $A_2$
is not free $D$\Hyph nodule,
then we may consider the set
\ShowEq{Ik 1n}
of linear independent linear maps. The theorem is true for any linear map
\ShowEq{f:A->B}f{A_1}{A_2}
generated by the set of linear maps $\Basis F$.
}
\StartLabelItem
\begin{enumerate}
\item
The map
\ShowEq{f:A->B}f{A_1}{A_2}
has the following expansion
\labelItem{map f generated by basis I}
\DrawEq{map f generated by basis I}{expansion}
where
\ShowEq{fk= in A2xA2}
\item
The map $f$ has the standard representation
\labelItem{standard representation of map A1 A2, associative algebra}
\ShowEq{standard representation of map A1 A2, associative algebra}
\end{enumerate}
}

\DefDefinition{component of linear map}
{
Expression
\ShowEq{component of linear map}
in equality
\eqRef{map f generated by basis I}{expansion}
is called
\AddIndex{component of linear map}
{component of linear map} $f$.
Expression
\ShowEq{standard component of linear map}
in the equality
\EqRef{standard representation of map A1 A2, associative algebra}
is called
\AddIndex{standard component of linear map}
{standard component of linear map} $f$.
}

\DefEq
{
\begin{proof}
Since $\Basis F$ is the basis of left \BoxB{A_2}module
\ShowEq{L(A;B)}D{A_1}{A_2},
then according to the definition
\ShowEq{ref definition: basis of representation}
and the theorem
\RefTheorem{set of vectors generated by set of vectors, \SideWS module},
there exists expansion
\ShowEq{expansion of linear map with respect to basis}
of the linear map $f$ with respect to the basis $\Basis I$.
According to the definition
\ShowEq{ref map j, representation, tensor product}
\ShowEq{f=fkxfk}
The equality
\eqRef{map f generated by basis I}{expansion}
follows from equalities
\EqRef{expansion of linear map with respect to basis},
\EqRef{f=fkxfk}.
According to theorem
\RefTheorem{standard component of tensor, algebra},
the standard representation of the tensor $f^k$ has form
\ShowEq{standard representation of map A1 A2, 3, associative algebra}
The equation
\EqRef{standard representation of map A1 A2, associative algebra}
follows from equations
\eqRef{map f generated by basis I}{expansion},
\EqRef{standard representation of map A1 A2, 3, associative algebra}.
\end{proof}
}
{proof: standard representation of map A1 A2, associative algebra}

\DefTheorem{Free Algebra over Ring}
{
The representation
\ShowEq{endomorphism of module from product}
of $D$\Hyph module $A$ in $D$\Hyph module $A$
is equivalent to structure of $D$\Hyph algebra $A$.
}

\DefEq
{
Let $D$ be commutative ring and $A$ be Abelian group.
The diagram of representations
\DrawEq[{}{}{}g]{diagram of representations of D algebra}{}
generates the structure of $D$\Hyph algebra $A$.
}
{text: Free Algebra over Ring}

\AddEq [2]{definition: linear map of A module}
{
\begin{definition}
\labelDefinition{linear map of #1 A module}
{\it
Let $V$ and
$W$ be #2 $A$\Hyph modules.
We call a map
\ShowEq{f:A->B}fVW
\AddIndex{linear map}{linear map}
if the map $f$
is linear map
of $D$\Hyph module $V$ into $D$\Hyph module $W$.
Let us denote
\ShowEq{set linear maps, D vector space}
set of linear maps
of #2 $A$\Hyph module $V$ into #2 $A$\Hyph module $W$.
}
\qed
\end{definition}
}

\DefTheorem{linear map of A module}
{
Linear map
\ShowEq{f:A->B}fVW
relative to basis \eV V
of left $A$\Hyph module $V$ and
basis \eV W
of left $A$\Hyph module $W$
has form
\ShowEq{linear map in A module}
where $f_{\gi i}^{\gi j}\circ v^{\gi i}$ linearly depends on
one variable $v^{\gi i}$ and does not depends on the rest of
coordinates of vector
\ShowEq{v=veV}
}

\DefProof{linear map of A module}
{
According to definition
\RefDefinition{linear map of A module},
\ShowEq{vector expansion, linear map, D vector space, 1}
For any given $\gi i$, vector
\ShowEq{vector expansion, linear map, D vector space, 2}
has only expansion
\ShowEq{vector expansion, linear map, D vector space, 3}
relative to basis \eV W.
From equalities
\begin{align*}
A_{\gi i}^{\gi j}(v^{\gi i}_1+v^{\gi i}_2)&=
A_{\gi i}^{\gi j}(v^{\gi i}_1)+A_{\gi i}^{\gi j}(v^{\gi i}_2)
\\
A_{\gi i}^{\gi j}(pv^{\gi i})&=
pA_{\gi i}^{\gi j}(v^{\gi i})\ \ \ p\in F
\end{align*}
it follows that map
\[A_{\gi i}^{\gi j}:D\rightarrow D\]
is linear map.
Therefore, we can use notation
\ShowEq{Aij vj}
Equation
\EqRef{linear map in vector space}
follows from equalities
\EqRef{vector expansion, linear map, D vector space, 1},
\EqRef{vector expansion, linear map, D vector space, 3},
\EqRef{Aij vj}.
}

\DefTheorem{representation of algebra An in LAnA}
{
Consider $D$\Hyph algebra $A$.
A representation
\ShowEq{representation An in LAnA}
of algebra $\AOn$
in module $\LAnA$
defined by the equation
\DrawEq{representation An in LAnA, 1}{}
allows us to identify tensor
\ShowEq{product in algebra An 1}
and transposition $\sigma\in S^n$
with map \ShowEq{product in algebra An 2} where
\ShowEq{product in algebra AA 3}
is identity map.
}

\AddEq{definition: linear map from A1 to A2, algebra}
{
\begin{definition}
\labelDefinition{linear map from A1 to A2, algebra}
Let $A_1$ and
$A_2$ be algebras over commutative ring $D$.
The linear map
of the $D_{\DF}$\Hyph module $A_\VF$
into the $D_{\DT}$\Hyph module $A_\VT$
is called
\AddIndex{linear map}{linear map}
of $D_{\DF}$\Hyph algebra $A_\VF$ into $D_{\DT}$\Hyph algebra $A_\VT$.

Let us denote
\ShowEq{set linear maps, module}
set of linear maps
of $D$\Hyph algebra
$A_\VF$
into $D$\Hyph algebra
$A_\VT$.
\qed
\end{definition}
}

\DefTheorem{linear map times constant, algebra}
{
Let map
\ShowEq{f:A->B}f{A_1}{A_2}
be linear map of $D$\Hyph algebra $A_1$ into $D$\Hyph algebra $A_2$.
Then maps
\ShowEq{linear map times constant, algebra}
defined by equalities
\ShowEq{linear map times constant, 0, algebra}
are linear.
}

\DefTheorem{linear map AA LAA}
{
\ePrints{1502.04063}
\ifx\Semafor\ValueOn
Consider $D$\Hyph algebras $A_1$ and $A_2$.
For given map
\ShowEq{f in L(A->B)}{A_1}{A_2},
\else
For given map
\ShowEq{f in L(A->B)}{A_1}{A_2}{}
of $D$\Hyph algebra $A_1$ into $D$\Hyph algebra $A_2$,
\fi
there exists linear map
\ShowEq{linear map AA LAA}
defined by the equality
\ShowEq{linear map AA LAA, 1}
}

\DefEq
{
\begin{theorem}
\labelTheorem{conjugation transformation}
Let $A_1$ be free $D$\Hyph module.
Let $A_2$ be free associative $D$\Hyph algebra.
Let $\Basis F$ be the basis of left \BoxB{A_2}module
\ShowEq{L(A;B)}D{A_1}{A_2}.
For any map
\ShowEq{Ik in I}
there exists set of linear maps
\ShowEq{conjugation transformation}
\ShowEq{conjugation transformation:}
of $D$\Hyph module
$A_1\otimes A_1$
into $D$\Hyph module
$A_2\otimes A_2$
such that
\ShowEq{conjugation transformation =}
The map
\ShowEq{show conjugation transformation}
is called
\AddIndex{conjugation transformation}{conjugation transformation}.
\end{theorem}
}
{theorem: conjugation transformation}

\DefEq
{
\begin{proof}
According to the theorem
\RefTheorem{product of linear map, algebra},
for any tensor
$a\in A_1\otimes A_1$,
the map
\ShowEq{x->Ik a x}
is linear.
According to the statement
\RefItem{map f generated by basis I},
there exists expansion
\ShowEq{expansion Ik a x}
Let
\ShowEq{b=Ikl a}
The equality
\EqRef{conjugation transformation =}
follows from equalities
\EqRef{expansion Ik a x},
\EqRef{b=Ikl a}.
From equalities
\ShowEq{Ikl a1+a2}
\ShowEq{Ikl d a}
it follows that the map $I_k^l$ is linear map.
\end{proof}
}
{proof: conjugation transformation}

\DefEq
{
\begin{theorem}
\labelTheorem{representation of composition of linear maps}
Let $A_1$ be free $D$\Hyph module.
Let $A_2$, $A_2$ be free associative $D$\Hyph algebras.
Let $\Basis F$ be the basis of left \BoxB{A_2}module
\ShowEq{L(A;B)}D{A_1}{A_2}.
Let $\Basis G$ be the basis of left \BoxB{A_3}module
\ShowEq{L(A;B)}D{A_2}{A_3}.
\StartLabelItem
\begin{enumerate}
\item
The set of maps
\labelItem{basis of composition of linear maps}
\DrawEq{JlIk}{linear map}
is the basis of left \BoxB{A_3}module
\ShowEq{L(A;B)}D{A_1\rightarrow A_2}{A_3}.
\item
Let
\ShowEq{expansion of f with respect to basis I}
be expansion of linear map
\ShowEq{f:A->B}f{A_1}{A_2}
with respect to the basis $\Basis I$.
Let
\ShowEq{expansion of g with respect to basis J}
be expansion of linear map
\ShowEq{f:A->B}g{A_2}{A_3}
with respect to the basis $\Basis J$.
Then linear map
\DrawEq{h=g o f}{123}
has expansion
\labelItem{hlk=glfk}
\ShowEq{expansion of h with respect to basis K}
with respect to the basis $\Basis K$ where
\ShowEq{hlk=glfk}
\end{enumerate}
\end{theorem}
}
{theorem: representation of composition of linear maps}

\DefEq
{
\begin{proof}
The equality
\ShowEq{h(a)1}
follows from equalities
\EqRef{expansion of f with respect to basis I},
\EqRef{expansion of g with respect to basis J},
\eqRef{h=g o f}{123}.
The equality
\ShowEq{h(a)2}
follows from equalities
\eqRef{JlIk}{linear map},
\EqRef{h(a)1}
and from the theorem
\RefTheorem{conjugation transformation}.
From the equality
\EqRef{h(a)2}
it follows that set of maps $\Basis K$ generates
left \BoxB{A_3}module\newline
\ShowEq{L(A;B)}D{A_1\rightarrow A_2}{A_3}.
From the equality
\ShowEq{aK=aJI}
it follows that
\ShowEq{aJ=0}
and, therefore, $a^{lk}=0$.
Therefore, the set $\Basis K$ is the basis of
left \BoxB{A_3}module
\ShowEq{L(A;B)}D{A_1\rightarrow A_2}{A_3}.
\end{proof}
}
{proof: representation of composition of linear maps}

\DefEq
{
\begin{theorem}
\labelTheorem{representation of composition of linear maps A->A}
Let $A$ be free associative $D$\Hyph algebra.
Let left \BoxB{A}module
\ShowEq{L(A;B)}DAA{}
is generated by the identity map $F_0=\delta$.
Let
\ShowEq{expansion of f A->A}
be expansion of linear map
\ShowEq{f:A->B}fAA
Let
\ShowEq{expansion of g A->A}
be expansion of linear map
\ShowEq{f:A->B}gAA
Then linear map
\DrawEq{h=g o f}{A}
has expansion
\ShowEq{expansion of h A->A}
where
\ShowEq{hlk=glfk 01}
\end{theorem}
}
{theorem: representation of composition of linear maps A->A}

\DefEq
{
\begin{proof}
The equality
\ShowEq{h(a)}
follows from equalities
\EqRef{expansion of f A->A},
\EqRef{expansion of g A->A},
\eqRef{h=g o f}{A}.
The equality
\EqRef{hlk=glfk 01}
follows from the equality
\EqRef{h(a)}.
\end{proof}
}
{proof: representation of composition of linear maps A->A}

\DefTheorem{linear map in L(A,A), associative algebra}
{
Let $A_1$ be $D$\Hyph algebra.
Let $A_2$ be free finite dimensional associative $D$\Hyph algebra.
Let $\Basis e$ be basis of $D$ module $A_2$.
Left \BoxB{A_2}module
\ShowEq{L(A;B)}D{A_1}{A_2}{}
has finite
\AddIndex{basis}{basis of algebra L(A,A)} $\Basis I$.
\StartLabelItem
\begin{enumerate}
\item
The linear map
\ShowEq{f in L(A->B)}{A_1}{A_2}{}
has form
\labelItem{f in L(A,A), 1, associative algebra}
\ShowEq{f in L(A,A), 1, associative algebra}
\item
Its standard representation has form
\labelItem{f in L(A,A), 2, associative algebra}
\ShowEq{f in L(A,A), 2, associative algebra}
\end{enumerate}
}

\DefEq
{
\begin{theorem}
\labelTheorem{h generated by f, associative algebra}
Consider $D$\Hyph algebra $A_1$
and associative $D$\Hyph algebra $A_2$.
Consider the representation of algebra $\ATwo$
in the module $\mathcal L(D;A_1;A_2)$.
The map
\ShowEq{h:A1->A2}
generated by the map
\ShowEq{f:A->B}f{A_1}{A_2}
has form
\ShowEq{h generated by f, associative algebra}
\end{theorem}
}
{theorem: h generated by f, associative algebra}

\DefEq
{
\begin{theorem}
\labelTheorem{coordinates of map A1 A2, algebra}
Let $\Basis e_1$ be basis of the free finite dimensional
$D$\Hyph module $A_1$.
Let $\Basis e_2$ be basis of the free finite dimensional associative
$D$\Hyph algebra $A_2$.
Let
\ShowEq{structural constants, algebra}
be structural constants of algebra $A_2$.
Let $\Basis F$ be the basis
of left \BoxB{A_2}module
\ShowEq{L(A;B)}D{A_1}{A_2}{}
and
\ShowEq{coordinates of map Ik}
be coordinates of map $F_k$ with respect to bases $\Basis e_1$ and $\Basis e_2$.
Coordinates
\ShowEq{Coordinates of map f}
of the map
\ShowEq{f in L(A->B)}{A_1}{A_2}{}
and its standard components
\ShowEq{standard components of map f}
are connected by the equation
\ShowEq{coordinates of map A1 A2, 2, associative algebra}
\end{theorem}
}
{theorem: coordinates of map A1 A2, algebra}

\DefEq
{
\begin{proof}
Relative to bases
$\Basis e_1$ and $\Basis e_2$, linear maps $f$ and $I_k$ have form
\ShowEq{coordinates of map f, associative algebra}
\ShowEq{coordinates of map Ik, associative algebra}
The equality
\ShowEq{coordinates of map A1 A2, 3, associative algebra}
follows from equalities
\EqRef{standard representation of map A1 A2, associative algebra},
\EqRef{coordinates of map f, associative algebra},
\EqRef{coordinates of map Ik, associative algebra}.
Since vectors $\Vector e_{2\cdot\gik}$
are linear independent and $x^{\gi i}$ are arbitrary,
then the equation
\EqRef{coordinates of map A1 A2, 2, associative algebra}
follows from the equation
\EqRef{coordinates of map A1 A2, 3, associative algebra}.
\end{proof}
}
{proof: coordinates of map A1 A2, algebra}

\DefEq
{
\begin{theorem}
\labelTheorem{sum of linear maps, D module}
Let $A_1$, $A_2$ be $D$\Hyph modules.
The map
\ShowEq{sum of maps,,D module}
\ShowEq{sum of maps, 1, D module}
defined by equation
\ShowEq{sum of maps, D module}
is called
\AddIndex{sum of maps}{sum of maps}
$f$ and $g$
and is linear map.
The set
$\mathcal L(D;A_1;A_2)$
is an Abelian group
relative sum of maps.
\end{theorem}
}
{theorem: sum of linear maps, D module}

\DefEq
{
\begin{theorem}
\labelTheorem{sum of polylinear maps, module}
Let $D$ be the commutative ring.
Let $A_1$, ..., $A_n$, $S$ be $D$\Hyph modules.
The map
\ShowEq{sum of maps,,polylinear}
\ShowEq{sum of maps, 1, polylinear}
defined by the equality
\ShowEq{sum of  maps, polylinear}
is called
\AddIndex{sum of polylinear maps}{sum of maps}
$f$ and $g$
and is polylinear map.
The set
\ShowEq{module of polylinear maps}
is an Abelian group
relative sum of maps.
\end{theorem}
}
{theorem: sum of polylinear maps, module}

\AddEq{corollary: sum of linear maps, D module}
{
\begin{corollary}
\labelCorollary{sum of linear maps, D module}
Let $A_1$, $A_2$ be $D$\Hyph modules.
The map
\ShowEq{sum of maps,,D module}
\ShowEq{sum of maps, 1, D module}
defined by equation
\ShowEq{sum of maps, D module}
is called
\AddIndex{sum of maps}{sum of maps}
$f$ and $g$
and is linear map.
The set
$\mathcal L(D;A_1;A_2)$
is an Abelian group
relative sum of maps.
\qed
\end{corollary}
}

\AddEq{definition: linear map A module}
{
\begin{definition}
\labelDefinition{linear map \SideWS A module}
Let
\ShowEq{Ai, i=}A
be algebra over commutative ring $D_i$.
Let
\ShowEq{Ai, i=}V
be
\ShowEq{\SideWS column module}{A_i}
module.
Morphism of diagram of representations
\ShowEq{A module diagram for linear}1
into diagram of representations
\ShowEq{A module diagram for linear}2
is called
\AddIndex{linear map}{linear map}
of
\ShowEq{\SideWS column module}{A_1}
module $V_1$ into
\ShowEq{\SideWS column module}{A_2}
module $V_2$.
Let us denote
\ShowEq{set linear maps, \SideWS A12 module}
set of linear maps
of
\ShowEq{\SideWS column module}{A_1}
module $V_1$ into
\ShowEq{\SideWS column module}{A_2}
module $V_2$.
\qed
\end{definition}
}

\AddEq{definition: homomorphism A module}
{
\begin{definition}
\labelDefinition{homomorphism \SideWS A module}
Let
\ShowEq{Ai, i=}A
be algebra over commutative ring $D_i$.
Let
\ShowEq{Ai, i=}V
be
\ShowEq{\SideWS column module}{A_i}
module.
Morphism of diagram of representations
\ShowEq{A module diagram for homomorphism}1
into diagram of representations
\ShowEq{A module diagram for homomorphism}2
is called
\AddIndex{homomorphism}{homomorphism}
of
\ShowEq{\SideWS column module}{A_1}
module $V_1$ into
\ShowEq{\SideWS column module}{A_2}
module $V_2$.
Let us denote
\ShowEq{set homomorphisms, \SideWS A12 module}
set of homomorphisms
of
\ShowEq{\SideWS column module}{A_1}
module $V_1$ into
\ShowEq{\SideWS column module}{A_2}
module $V_2$.
\qed
\end{definition}
}

\DefTheorem{L(An;B) is free D module}
{
Let $A_1$, ..., $A_n$, $B$ be free modules over commutative ring $D$.
$D$\Hyph module
\ShowEq{L(A;B)}D{A_1\times...\times A_n}B{}
is free $D$\Hyph module.
}

\DefEq
{
\begin{definition}
\labelDefinition{linear map from A1 to A2, module}
Reduced morphism of representations
\ShowEq{f:A->B}f{A_1}{A_2}
of $D$\Hyph module $A_1$
into $D$\Hyph module $A_2$
is called
\AddIndex{linear map}{linear map}
of \SideWS $D$\Hyph module $A_1$ into \SideWS $D$\Hyph module $A_2$.

Let us denote
\ShowEq{set linear maps, module}
set of linear maps
of \SideWS $D$\Hyph module $A_1$ into \SideWS $D$\Hyph module $A_2$.
\qed
\end{definition}
}
{definition: linear map from A1 to A2, module}

\AddEq{theorem: linear map of module}
{
\begin{theorem}
\labelTheorem{linear map from \Base\DF \Module\DF to \Base\DT \Module\DT, \SideWS module}
Linear map
\ShowEq{show linear map \MapE}f
of \SideWS $\Base_\DF$\Hyph module $\Module_\DF$
into \SideWS $\Base_\DT$\Hyph module $\Module_\DT$
satisfies to equations\,\footnote{
In some books
(for instance, on page \citeBib{Serge Lang}\Hyph 119) the theorem
\RefTheorem{linear map from \Base\DF \Module\DF to \Base\DT \Module\DT, \SideWS module}
is considered as a definition.}
\ShowEq{property linear map \MapE}
\end{theorem}
}

\DefProof{linear map of module}
{
From definitions
\RefDefinition[\RefRepresentation]{morphism of representations of universal algebra},
\RefDefinition{linear map from D1 A1 to D2 A2, module},
it follows that
the map $h$ is a homomorphism of the ring $D_1$
into the ring $D_2$ (the equalities
\eqRef{h(d1+d2)=...}{\SideWS module},
\eqRef{h(d1d2)=...}{\SideWS module})
and the map $f$ is a homomorphism of the Abelian group $A_1$
into the Abelian group $A_2$ (the equality
\eqRef{f(a+b)=...}{D1 D2 \SideWS module}).
The equality
\EqRef{\DFDT, f(da)=... \SideWS module}
follows from the equality
\EqRef[\RefRepresentation]{morphism of representations of F algebra, definition, 2m}.
}

\AddEq{theorem: linear map of module, coordinates}
{
\begin{theorem}
\labelTheorem{linear map of \SideWS \DFDT module, coordinates}
Let
\ShowEq{basis, module V1}
be a basis of free \SideWS $\Base_\DF$\Hyph module $\Module_1$
and
\ShowEq{basis, module V2}
be a basis of free \SideWS $\Base_\DT$\Hyph module $\Module_2$.
Then linear map
\ShowEq{show linear map \MapE}{\Vector f}
has presentation
\ShowEq{r2:A1->A2, \SideWS \DFDT module}
relative to selected bases. Here
\begin{itemize}
\item $a$ is coordinate matrix of vector
$\Vector a$ relative the basis $\Basis e_1$
\DrawEq [a1]{va=ae1, \SideWS module}{a \DFDT}
\def\Temp{D1 D2 }
\ifx\DFDT\Temp
\item
\ShowEq{h(a)=...}
is a matrix of $D_2$\Hyph numbers.
\fi
\item $b$ is coordinate matrix of vector
\DrawEq{vb=f(va)}{\DFDT\SideWS module}
relative the basis $\Basis e_2$
\DrawEq[b2]{va=ae1, \SideWS module}{b \DFDT}
\item $f$ is coordinate matrix of set of vectors
\ShowEq{Vector f(e1) \SideWS module}
relative the basis $\Basis e_2$.
The matrix $f$ is
called \AddIndex{matrix of linear map}{matrix of linear map}
relative bases $\Basis e_1$ and $\Basis e_2$.
\end{itemize}
\end{theorem}
}

\DefProof{linear map of module, coordinates}
{
Since
\ShowEq{show linear map \MapE}{\Vector f}
is a linear map, from
\EqRef{\DFDT, f(da)=... \SideWS module},
\eqRef{va=ae1, \SideWS module}{a \DFDT},
\eqRef{vb=f(va)}{\DFDT\SideWS module},
it follows that
\ShowEq{\DFDT vb=h(e1)a \SideNS}
\ShowEq{f(e1) \SideNS}
is a vector of free \SideWS $\Base_\DT$\Hyph module $V_2$ and has
expansion
\DrawEq{f(e1)= \SideNS}{\DFDT}
relative to basis $\Basis e_2$.
Combining \EqRef{\DFDT vb=h(e1)a \SideNS}
and \eqRef{f(e1)= \SideNS}{\DFDT}, we get
\ShowEq{\DFDT vb=a f e \SideNS}
\EqRef{r2:A1->A2, \SideWS \DFDT module}
follows from comparison of
\eqRef{va=ae1, \SideWS module}{b \DFDT}
and \EqRef{\DFDT vb=a f e \SideNS} and
theorem \RefTheorem{coordinates of vector of free \SideWS module}.
}

\DefDefinition{algebra over ring}
{
Let $D$ be commutative ring.
$D$\Hyph module $A$ is called
\AddIndex{algebra over ring}{algebra over ring} $D$
or
\AddIndex{$D$\Hyph algebra}{D algebra},
if we defined product\,\footnote{
I follow the definition
given in
\citeBib{Richard D. Schafer}, page 1,
\citeBib{0105.155}, page 4.  The statement which
is true for any $D$\Hyph module,
is true also for $D$\Hyph algebra.} in $A$
\DrawEq{product in D algebra}{definition}
where $C$ is bilinear map
\ShowEq{product in algebra, definition 1}
If $A$ is free
$D$\Hyph module, then $A$ is called
\AddIndex{free algebra over ring}{free algebra over ring} $D$.
}

\DefDefinition{commutator of algebra}
{
The \AddIndex{commutator}{commutator of algebra}
\ShowEq{commutator of algebra}
measures commutativity in $D$\Hyph algebra $A$.
$D$\Hyph algebra $A$ is called
\AddIndex{commutative}{commutative D algebra},
if
\ShowEq{commutative D algebra}
}

\DefDefinition{nucleus of algebra}
{
The set\,\footnote{
The definition is based on
the similar definition in
\citeBib{Richard D. Schafer}, p. 13}
\ShowEq{nucleus of algebra}
is called the
\AddIndex{nucleus of an $D$\Hyph algebra $A$}{nucleus of algebra}.
}

\DefDefinition{associator of algebra}
{
The \AddIndex{associator}{associator of algebra}
\ShowEq{associator of algebra}
measures associativity in $D$\Hyph algebra $A$.
$D$\Hyph algebra $A$ is called
\AddIndex{associative}{associative D algebra},
if
\ShowEq{associative D algebra}
}

\DefEq
{
\begin{definition}
\labelDefinition{center of algebra}
The set\,\footnote{
The definition is based on
the similar definition in
\citeBib{Richard D. Schafer}, page 14}
\ShowEq{center of algebra}
is called the
\AddIndex{center of an $D$\Hyph algebra $A_1$}{center of algebra}.
\qed
\end{definition}
}
{definition: center of algebra}

\DefEq
{
\begin{definition}
\labelDefinition{otimes -}
Bilinear map
\ShowEq{otimes -}
is defined by the equality
\ShowEq{otimes -, 1}
\qed
\end{definition}
}
{definition: otimes -}

%% file: Stmt.Module.Eq.tex
\input{Stmt.Module.Ref}

\newcommand\eV[1]{$\Basis e_{#1}$}
\newcommand\EV[2]{e_{#1\cdot\gi{#2}}}
\def\AOn{A^{\otimes n+1}}
\def\LAnA{\mathcal L(D;A^n\rightarrow A)}

\AddEq{show definition of A module}
{
\ShowDefinition{module over algebra}

\ShowTheorem{module over algebra}
\ShowProof{module over algebra}

\ShowTheorem{A module -> alebra is associative}
\ShowProof{A module -> alebra is associative}

\ShowTheorem{definition of A module}
\ShowProof{definition of A module}
}

\AddEq{show linear combination of vectors}
{
\ShowTheorem{set of vectors generated by set of vectors}
\ShowProof{set of vectors generated by set of vectors}

\ShowDefinition{linear combination of vectors}

\ShowEq{remark: scalars and vectors as matrix}

\ShowTheorem{linearly depends on rest of vectors}
\ShowProof{linearly depends on rest of vectors}

\ShowRemark{0=0vi}

\ShowDefinition{linearly independent vectors}
}

\AddEq{show basis of module}
{
\ShowRemark{generating set of module}

\ShowDefinition{generating set of module}

\ShowRemark{basis of module}

\ShowDefinition{basis of module}

\ShowTheorem{basis of module}
\ShowProof{basis of module}

\ShowTheorem{division algebra, basis}
\ShowProof{division algebra, basis}

\ShowDefinition{coordinates of vector}

\ShowTheorem{linear dependence between vectors of basis}
\ShowProof{linear dependence between vectors of basis}

\ShowTheorem{coordinates of vector with linear dependence}
\ShowProof{coordinates of vector with linear dependence}

\ShowDefinition{free module over ring}

\ShowTheorem{coordinates of vector of free module}
\ShowProof{coordinates of vector of free module}
}

\DefEq
{
\begin{matrix}
\xymatrix
{
f:D\ar[r]|{*}&V
}
&
f(d):v\rightarrow d\,v
\end{matrix}
}
{D->*V}

\AddEq{v=veV}
{
$v=v^{\gi i}\EV V i$.
}

\AddEquation{linear map in A module}
{
f\circ(v\CRstar e_V )=(f_{\gi i}^{\gi j}\circ v^{\gi i})\ \EV W j
}

\AddEq [1]{Ai, i=}
{
$#1_i$, $i=1$, $2$,
}

\AddEquation{module, (a+n)v=av+nv}
{
(a+n)v=av+nv\ \ \ a\in D\ \ \ n\in Z
}

\AddEquation{left module, (a+n)v=av+nv}
{
(a+n)v=av+nv\ \ \ a\in A\ \ \ n\in D
}

\AddEquation{right module, (a+n)v=av+nv}
{
v(a+n)=va+vn\ \ \ a\in A\ \ \ n\in D
}

\AddEq{(a+n)+(b+m)=}
{
(a+n)+(b+m)=(a+b)+(n+m)
}

\AddEquation{module, (a+n)+(b+m)=}
{
\begin{split}
((a+n)+(b+m))v&=av+nv+bv+mv=av+bv+nv+mv\\
&=(a+b)v+(n+m)v=((a+b)+(n+m))v
\end{split}
}

\AddEquation{left module, (a+n)+(b+m)=}
{
\begin{split}
((a+n)+(b+m))v&=av+nv+bv+mv=av+bv+nv+mv\\
&=(a+b)v+(n+m)v=((a+b)+(n+m))v
\end{split}
}

\AddEquation{right module, (a+n)+(b+m)=}
{
\begin{split}
v((a+n)+(b+m))&=va+vn+vb+vm=va+vb+vn+vm\\
&=v(a+b)+v(n+m)=v((a+b)+(n+m))
\end{split}
}

\AddEquation{module, (a+n)+(b+m)=, 1}
{
((a+n)+(b+m))v=(a+n)v+(b+m)v
}

\AddEquation{left module, (a+n)+(b+m)=, 1}
{
((a+n)+(b+m))v=(a+n)v+(b+m)v
}

\AddEquation{right module, (a+n)+(b+m)=, 1}
{
v((a+n)+(b+m))=v(a+n)+v(b+m)
}

\AddEq{(a+n)(b+m)=}
{
(a+n)(b+m)=(ab+ma+nb)+(nm)
}

\AddEq [1]{A module diagram for linear}
{
\[
\xymatrix
{
A_{#1}&&V_{#1}\\
&D_{#1}\ar[lu]|{*}\ar[ru]|{*}
}
\]
}

\AddEq [1]{A module diagram for homomorphism}
{
\[
\xymatrix
{
A_{#1}\ar[r]|{*}&A_{#1}\ar[r]|{*}&V_{#1}\\
&D_{#1}\ar[lu]|{*}\ar[u]|{*}\ar[ru]|{*}
}
\]
}

\DefEquation
{
f^k=f^{k\cdot\gi{ij}}e_{\gii}\otimes e_{\gij}
}
{standard representation of map A1 A2, 3, associative algebra}

\DefEquation
{
f=
f^{k\cdot\gi{ij}}(e_{\gii}\otimes e_{\gij})\circ F_k
=
f^{k\cdot\gi{ij}}e_{\gii}F_ke_{\gij}
}
{standard representation of map A1 A2, associative algebra}

\AddEq{corollary: Free Algebra over Ring}
{
\begin{corollary}
\labelCorollary{Free Algebra over Ring}
\ShowEq{text: Free Algebra over Ring}
\qed
\end{corollary}
}

\DefEquation
{
\begin{split}
f^{\gik}_{\gil}x^{\gil}e_{2\cdot\gik}
&=
f^{k\cdot\gi{ij}} e_{2\cdot\gii}
I_{k\cdot}^{}{}^{\gim}_{\gil}x^{\gil}e_{2\cdot\gim}
\Vector e_{2\cdot\gij}
\\
&=
f^{k\cdot\gi{ij}}
I_{k\cdot}^{}{}^{\gim}_{\gil}x^{\gil}
C^{\gi p}_{\gi{im}}
C^{\gik}_{\gi{pj}}
e_{2\cdot\gik}
\end{split}
}
{coordinates of map A1 A2, 3, associative algebra}

\DefEq
{
f=f^k\circ F_k
}
{map f generated by basis I}

\AddEquation{fk= in A2xA2}
{
f^k=f^k_{s_k\cdot 0}\otimes f^k_{s_k\cdot 1}\ \ \ \ f^k\in A_2\otimes A_2
}

\DefEquation
{
f=f_{s\cdot 0}\otimes f_{s\cdot 1}
}
{expansion of f A->A}

\DefEquation
{
f=f^k\circ F_k
}
{expansion of f with respect to basis I}

\DefEquation
{
g=g_{t\cdot 0}\otimes g_{t\cdot 1}
}
{expansion of g A->A}

\DefEquation
{
g=g^l\circ G_l
}
{expansion of g with respect to basis J}

\DefEquation
{
h=h_{ts\cdot 0}\otimes h_{ts\cdot 1}
}
{expansion of h A->A}

\DefEquation
{
h=h^{lk}\circ K_{lk}
}
{expansion of h with respect to basis K}

\DefEquation
{
h\circ a=g\circ f\circ a
=g^l\circ J_l\circ f^k\circ I_k\circ a 
}
{h(a)1}

\DefEquation
{
\begin{split}
h\circ a=g\circ f\circ a
&=g^l\circ(J^m_l\circ f^k)\circ
J_m\circ I_k\circ a
\\&=g^l\circ(J^m_l\circ f^k)\circ K_{mk}\circ a
\end{split}
}
{h(a)2}

\DefEquation
{
\begin{split}
h_{ts\cdot 0}&=g_{t\cdot 0}f_{s\cdot 0}
\\
h_{ts\cdot 1}&=f_{s\cdot 1}g_{t\cdot 1}
\end{split}
}
{hlk=glfk 01}

\DefEq
{
\[a^{lk}K_{lk}=(a^{lk}\circ J_l)\circ I_k=0\]
}
{aK=aJI}

\DefEquation
{
\begin{split}
h\circ a&=g\circ f\circ a
\\&=(g_{t\cdot 0}\otimes g_{t\cdot 1})\circ(f_{s\cdot 0}\otimes f_{s\cdot 1})\circ a
\\&=(g_{t\cdot 0}\otimes g_{t\cdot 1})\circ(f_{s\cdot 0}a f_{s\cdot 1})
\\&=g_{t\cdot 0}f_{s\cdot 0}a f_{s\cdot 1} g_{t\cdot 1}
\end{split}
}
{h(a)}

\DefEq
{
\[a^{lk}\circ J_l=0\]
}
{aJ=0}

\DefEquation
{
h^{lk}=g^l\circ (G^k_m\circ f^m)
}
{hlk=glfk}

\DefEq
{
\Basis H=\{H_{lk}:H_{lk}=G_l\circ F_k, G_l\in\Basis G, F_k\in\Basis F\}
}
{JlIk}

\DefEq
{
h=g\circ f
}
{h=g o f}

\DefEquation
{
\xymatrix
{
h:A\otimes A\ar[r]|{*}&\mathcal L(D;A\rightarrow A)
}
\ \ \ h(p):f\rightarrow p\circ f
}
{h:AoxA->L(A)}

\DefEq
{
$A\otimes A$
}
{AoxA}

\ePrints{1502.04063,5114-6019}
\ifx\Semafor\ValueOff
\DefEq
{
\[
\begin{matrix}
(a\otimes b)\circ f=afb
&a,b\in A&f\in\mathcal L(D;A\rightarrow A)
\end{matrix}
\]
}
{representation AA in LA}

\AddEq{component of linear map}
{
\symb{f^k_{s_k\cdot p}}{component of linear map, division ring}1,
$p=0$, $1$,
}

\AddEq{standard component of linear map}
{
\symb{f^{k\cdot\gi{ij}}}{standard component of linear map}1
}
\fi

\AddEq{product in algebra AA}
{
(p_0\otimes p_1)\circ(q_0\otimes q_1)=(p_0q_0)\otimes(q_1p_1)
}

\AddEq{xA1n*xA1n->oxA1n=}
{
\[
(a_1,...,a_n)*(b_1,...,b_n)=(a_1b_1)\otimes...\otimes(a_nb_n)
\]
}

\DefEquation
{
(a\otimes b)\circ c=acb
}
{a ox b c=}

\AddEq{a ox b}
{
$(a\otimes b)\circ\delta$
}

\DefEq
{
$d\in A\otimes A$
}
{d in AxoA}

\AddEq{product in algebra AA 2}
{
$d\circ\delta\in\mathcal L(D;A\rightarrow A)$%
}

\DefEq
{
$A_1$, ..., $A_n$
}
{A1...n}

\AddEq{right A->*V}
{
\begin{matrix}
\xymatrix
{
f:A\ar[r]|{*}&V
}
&
f(a):v\in V\rightarrow va\in V
&a\in A
\end{matrix}
}

\AddEquation{diagram of representations, left module}
{
\begin{matrix}
\vcenter
{
\xymatrix
{
D\ar[r]|(.4){*}^{g_{1,2}}&
A\ar[r]|{*}^{g_{2,3}}&
A\ar[r]|{*}^{g_{3,4}}&
V
\\
&&D\ar[u]|{*}_{g_{1,2}}
&D\ar[u]|{*}_{g_{1,4}}
}
}
&
\begin{array}{r@{\,}l}
g_{1,2}(d):a&\rightarrow d\,a
\\
g_{2,3}(v):w&\rightarrow
C(w, v)
\\
C&\in\mathcal L(A^2\rightarrow A)
\\
g_{3,4}(a):v&\rightarrow a\,v
\\
g_{1,4}(d):v&\rightarrow d\,v
\end{array}
\end{matrix}
}

\AddEquation{diagram of representations, right module}
{
\begin{matrix}
\vcenter
{
\xymatrix
{
D\ar[r]|(.4){*}^{g_{1,2}}&
A\ar[r]|{*}^{g_{2,3}}&
A\ar[r]|{*}^{g_{3,4}}&
V
\\
&&D\ar[u]|{*}_{g_{1,2}}
&D\ar[u]|{*}_{g_{1,4}}
}
}
&
\begin{array}{r@{\,}l}
g_{1,2}(d):a&\rightarrow d\,a
\\
g_{2,3}(v):w&\rightarrow
C(w, v)
\\
C&\in\mathcal L(A^2\rightarrow A)
\\
g_{3,4}(a):v&\rightarrow v\,a
\\
g_{1,4}(d):v&\rightarrow d\,v
\end{array}
\end{matrix}
}

\AddEquation{diagram of representations, module}
{
\xymatrix
{
D\ar[r]|{*}^{g_2}&V\\
&Z\ar[u]|{*}_{g_1}
}
}

\AddEq{left A->*V}
{
\begin{matrix}
\xymatrix
{
f:A\ar[r]|{*}&V
}
&
f(a):v\in V\rightarrow av\in V
&a\in A
\end{matrix}
}

\AddEq{associative law, module}
{
(pq)v=p(qv)
}

\AddEq{module, associative law}
{
(ab)v=a(bv)
}

\AddEquation{module, (a+n)(b+m)=}
{
\begin{split}
((a+n)(b+m))v&=(a+n)((b+m)v)=(a+n)(bv+mv)\\
=&a(bv+mv)+n(bv+mv)\\
=&a(bv)+a(mv)+n(bv)+n(mv)\\
=&(ab)v+m(av)++n(bv)+(nm)v\\
=&(ab)v+(ma)v++(nb)v+(nm)v\\
=&((ab+ma+nb)+nm)v
\end{split}
}

\AddEq{fab=ab}
{
f(a,b)=ab\ \ \ a,b\in\Base
}

\AddEq{f1p=p}
{
f(1,p)=f(p,1)=p\ \ \ p\in\Base_{(1)}\ \ \ 1\in\CBase_{(1)}
}

\AddEq{pq=fpq}
{
\begin{split}
(a+n)(b+m)&=f(a+n,b+m)\\
&=f(a,b)+f(a,m)+f(n,b)+f(n,m)\\
&=f(a,b)+mf(a,1)+nf(1,b)+nf(1,m)\\
&=ab+ma+nb+nm
\end{split}
}

\AddEq{f:D1xD1->D1}
{
\[
f:\Base_{(1)}\times\Base_{(1)}\rightarrow\Base_{(1)}
\]
}

\AddEq{distributive law, module}
{
\begin{align}
\EqLabel{distributive law, module, 1}
p(v+w)&=pv+pw\\
\EqLabel{distributive law, module, 2}
(p+q)v&=pv+qv
\end{align}
}

\AddEq{p,q in D1}
{
$p=a+n\in\Base_{(1)}$, $q=b+m\in\Base_{(1)}$.
}

\AddEq [4]{h:D1->D2 f:A1->A2}
{
\[
\begin{pmatrix}
#1:#2_1\rightarrow #2_2&
#3:#4_1\rightarrow #4_2
\end{pmatrix}
\]
}

\AddEq{set linear maps, D12 module}
{
\symb{\mathcal L(D_1\rightarrow D_2;A_1\rightarrow A_2)}{set linear maps}1
}

\AddEq{set linear maps, left A12 module}
{
\symb{\mathcal L(D_1\rightarrow D_2;A_1\CRstar\rightarrow A_2\CRstar;V_1\rightarrow V_2)}{set linear maps}1
}

\AddEq{set homomorphisms, left A12 module}
{
\symb{\Hom(D_1\rightarrow D_2;A_1\CRstar\rightarrow A_2\CRstar;V_1\rightarrow V_2)}{set homomorphisms}1
}

\AddEq{set linear maps, right A12 module}
{
\symb{\mathcal L(D_1\rightarrow D_2;\RCstar A_1\rightarrow\RCstar A_2;V_1\rightarrow V_2)}{set linear maps}1
}

\AddEq{set homomorphisms, right A12 module}
{
\symb{\Hom(D_1\rightarrow D_2;\RCstar A_1\rightarrow\RCstar A_2;V_1\rightarrow V_2)}{set homomorphisms}1
}

\AddEq [1]{left column module}
{
$#1\CRstar$\Hyph
}

\AddEq [1]{right column module}
{
$\RCstar #1$\Hyph
}

\AddEquation{unitarity law, D-module}
{
1v=v
}

\DefEq
{
\[
a=a^{\gi i}e_{\gi i}
\]
}
{Expansion relative basis in algebra}

\AddEq{associative law, left module}
{
(pq)v=p(qv)
}

\AddEq{left module, associative law}
{
(ab)v=a(bv)
}

\AddEq{associative law, Z, module}
{
(mn)v=m(nv)
}

\AddEq{associative law, ZD, module}
{
(md)v=m(dv)
}

\AddEq{associative law, D, left module}
{
(cd)v=c(dv)
}

\AddEq{associative law, DA, left module}
{
(da)v=d(av)
}

\AddEq{distributive law, left module}
{
\begin{align}
\EqLabel{distributive law, left module, 1}
p(v+w)&=pv+pw\\
\EqLabel{distributive law, left module, 2}
(p+q)v&=pv+qv
\end{align}
}

\AddEquation{distributive law, Z, module, 2}
{
(n+m)v=nv+mv
}

\AddEquation{distributive law, D, module, 2}
{
(a+b)v=av+bv
}

\AddEquation{distributive law, D, left module, 2}
{
(n+m)v=cn+dm
}

\AddEquation{distributive law, A, left module, 2}
{
(a+b)v=av+bv
}

\AddEquation{distributive law, D, right module, 2}
{
v(n+m)=vn+vm
}

\AddEquation{distributive law, A, right module, 2}
{
v(a+b)=va+vb
}

\AddEq{distributive law, right module}
{
\begin{align}
\EqLabel{distributive law, right module, 1}
(v+w)p&=vp+wp\\
\EqLabel{distributive law, right module, 2}
v(p+q)&=vp+vq
\end{align}
}

\AddEq{aw in Xk module}
{
$aw\in X_k$.
}

\AddEq{aw in Xk left module}
{
$aw\in X_k$.
}

\AddEq{aw in Xk right module}
{
$wa\in X_k$.
}

\AddEq[2]{w1+w2 in Jv}
{
$#1_1+#1_2\in J(#2)$.
}

\AddEquation{unitarity law, left A-module}
{
1v=v
}

\AddEquation{left module, a d v}
{
a(dv)=d(av)
}

\AddEquation{module, a d v}
{
a(nv)=n(av)
}

\AddEq{abc in A, v in v}
{
$a$, $b$, $c\in A$, $v\in V$.
}

\AddEquation{a(b(cv))=(a(bc))v left}
{
a(b(cv))=a((bc)v)=(a(bc))v
}

\AddEquation{a(b(cv))=(a(bc))v right}
{
((vc)b)a=(v(cb))a=v((cb)a)
}

\AddEquation{a(b(cv))=((ab)c)v left}
{
a(b(cv))=(ab)(cv)=((ab)c)v
}

\AddEquation{a(b(cv))=((ab)c)v right}
{
((vc)b)a=(vc)(ba)=v(c(ba))
}

\AddEquation{(a(bc))v=((ab)c)v left}
{
(a(bc))v=((ab)c)v
}

\AddEquation{(a(bc))v=((ab)c)v right}
{
v((cb)a)=v(c(ba))
}

\AddEquation{a(bc)=(ab)c left}
{
a(bc)=(ab)c
}

\AddEquation{a(bc)=(ab)c right}
{
(cb)a=c(ba)
}

\AddEq{cv left}
{
$cv\in A$,
}

\AddEq{cv right}
{
$vc\in A$,
}

\AddEq{associative law, right module}
{
v(pq)=(vp)q
}

\AddEq{right module, associative law}
{
v(ab)=(va)b
}

\AddEq{associative law, D, right module}
{
v(dc)=(vd)c
}

\AddEq{associative law, DA, right module}
{
v(ad)=(va)d
}

\AddEquation{unitarity law, right A-module}
{
v1=v
}

\AddEquation{right module, a d v}
{
(vd)a=(va)d
}

\AddEq{c,d in Z, a,b in D, v,w in V}
{
$m$, $n\in Z$, $a$, $b \in D$, $v$, $w \in V$.
}

\AddEq{p,q in D1, v,w in V}
{
$p$, $q \in D_{(1)}$, $v$, $w \in V$.
}

\AddEq{p,q in A1, v,w in V}
{
$p$, $q \in A_{(1)}$, $v$, $w \in V$.
}

\AddEq [1]{vi V}
{
$v=(v_{\gii}\in V,\gii\in\giI)$#1
}

\DefEq
{
$\Vector v\in V$
}
{vv in V}

\AddEq{vv=ve left module}
{
\Vector v=v\CRstar e=v^{\gii}e_{\gii}
}

\AddEq{vv=ve right module}
{
\Vector v=e\RCstar v=e_{\gii}v^{\gii}
}

\AddEq{vv=ve module}
{
\Vector v=v\CRstar e=v^{\gii}e_{\gii}
}

\AddEq{coordinate matrix of vector}
{
$v=(v^{\gii},\Ii)$
}

\AddEquation{w= module}
{
w=\sum_{\gii\in\giI}w^{\gii}v_{\gii}
}

\AddEquation{w= left module}
{
w=\sum_{\gii\in\giI}w^{\gii}v_{\gii}
}

\AddEquation{w= right module}
{
w=\sum_{\gii\in\giI}v_{\gii}w^{\gii}
}

\AddEquation{aw= module}
{
aw=a\sum_{\gii\in\giI}w^{\gii}v_{\gii}
=\sum_{\gii\in\giI}a(w^{\gii}v_{\gii})
=\sum_{\gii\in\giI}(aw^{\gii})v_{\gii}
}

\AddEquation{aw= left module}
{
aw=a\sum_{\gii\in\giI}w^{\gii}v_{\gii}
=\sum_{\gii\in\giI}a(w^{\gii}v_{\gii})
=\sum_{\gii\in\giI}(aw^{\gii})v_{\gii}
}

\AddEquation{aw= right module}
{
wa=\left(\sum_{\gii\in\giI}v_{\gii}w^{\gii}\right)a
=\sum_{\gii\in\giI}(v_{\gii}w^{\gii})a
=\sum_{\gii\in\giI}(v_{\gii}w^{\gii}a)
}

\AddEq{aw in Jv module}
{
$aw\in J(v)$.
}

\AddEq{aw in Jv left module}
{
$aw\in J(v)$.
}

\AddEq{aw in Jv right module}
{
$wa\in J(v)$.
}

\DefEq
{
$w^{\gii}$, $\gii\in\giI$,
}
{wi}

\AddEq{basis for module over algebra}
{
\symb{\Basis{e}=(e_{\gii},\gii\in\giI)}{basis for module}1
}

\AddEq{set vi}
{%
$v_{\gii}$, $\gii\in\giI$,
}

\AddEq{w=wi}
{
\[
w=
\begin{pmatrix}
w^{\gi 1}\\...\\w^{\gi n}
\end{pmatrix}
\]
}

\DefEquation
{
\xymatrix
{
g_{2,3}:A\ar[r]|{*}&A
}
}
{endomorphism of module from product}

\DefEquation
{
l\circ v:w\in A\rightarrow v\,w\in A
}
{l(v):w->vw}

\DefEq
{
\[
\xymatrix
{
h:\AOn\times S_n\ar[r]|{*}&\LAnA
}
\]
}
{representation An in LAnA}

\DefEq
{
\begin{matrix}
(a_0\otimes...\otimes a_n,\sigma)\circ(f_1\otimes...\otimes f_n)
=a_0\sigma(f_1)a_1...a_{n-1}\sigma(f_n)a_n
\\
\begin{matrix}
a_0,...,a_n\in A&\sigma\in S_n&f_1,...,f_n\in\mathcal L(D;A_n\rightarrow A)
\end{matrix}
\end{matrix}
}
{representation An in LAnA, 1}

\DefEq
{
$d\in\AOn$
}
{product in algebra An 1}

\DefEquation
{
\begin{matrix}
(d,\sigma)\circ(f_1,...,f_n)&f_i=\delta\in\mathcal L(D;A\rightarrow A)
\end{matrix}
}
{product in algebra An 2}

\DefEq
{
$\delta\in\mathcal L(D;A\rightarrow A)$
}
{product in algebra AA 3}

\DefEq
{
\[
v=
\begin{pmatrix}
v_{\gi 1}&...&v_{\gi n}
\end{pmatrix}
\]
}
{v=vi}

\AddEq{wi vi=0}
{
\[w^{\gii}v_{\gii}=0\]
}

\AddEq{left wi vi=0}
{
\[w^{\gii}v_{\gii}=0\]
}

\AddEq{right wi vi=0}
{
\[v_{\gii}w^{\gii}=0\]
}

\AddEq{0=0vi}
{
\[
w^{\gii}=0\Rightarrow w\CRstar v=0
\]
}

\AddEq{left 0=0vi}
{
\[
w^{\gii}=0\Rightarrow w\CRstar v=0
\]
}

\AddEq{right 0=0vi}
{
\[
w^{\gii}=0\Rightarrow v\RCstar w=0
\]
}

\AddEq[1]{vj=sum vi}
{
\[
#1_{\gij}=\sum_{\gii\in \giI\setminus\{\gij\}}
\frac{w^{\gii}}{w^{\gij}}#1_{\gii}
\]
}

\AddEq[1]{left vj=sum vi}
{
\[
#1_{\gij}=\sum_{\gii\in \giI\setminus\{\gij\}}
(w^{\gij})^{-1}w^{\gii}#1_{\gii}
\]
}

\AddEq[2]{b in D'}
{
$#1=(#1^{\gii},\gii\in\giI)\in\Base'$#2
}

\AddEq[1]{right vj=sum vi}
{
\[
#1_{\gij}=\sum_{\gii\in \giI\setminus\{\gij\}}
#1_{\gii}w^{\gii}(w^{\gij})^{-1}
\]
}

\AddEq[1]{wj ne 0}
{
$w^{\gij}\ne 0$#1
}

\AddEq[2]{w12 in Jv}
{
$#1_1$, $#1_2\in J(#2)$.
}

\AddEquation{w=sum vi, module}
{
J(v)=\left\{w:w=\sum_{\gii\in \giI}c^{\gii}v_{\gii}, c^{\gii}\in D_{(1)},
|\{\gii:c^{\gii}\ne 0\}|<\infty\right\}
}

\AddEq{generating set of module}
{
$J(X)=V$.
}

\AddEquation{w=sum vi, left module}
{
J(v)=\left\{w:w=\sum_{\gii\in \giI}c^{\gii}v_{\gii}, c^{\gii}\in A_{(1)},
|\{\gii:c^{\gii}\ne 0\}|<\infty\right\}
}

\AddEquation{w=sum vi, right module}
{
J(v)=\left\{w:w=\sum_{\gii\in \giI}v_{\gii}c^{\gii}, c^{\gii}\in A_{(1)},
|\{\gii:c^{\gii}\ne 0\}|<\infty\right\}
}

\AddEq{vectors generated by set of vectors}
{
\StartLabelItem
\begin{enumerate}
\item
\ShowEq{vk in Jv}
\labelItem{\SideWS module, vk in Jv}
\item
\ShowEq{\SideWS module, cvk in Jv}
\item
\ShowEq{\SideWS module, sum cvk in Jv}
\item
\ShowEq{v,w in Jv=>v+w in Jv}
\item
\ShowEq{v in Jv=>av in Jv}
\end{enumerate}
}

\AddEq{vk in Jv}
{%
$v_{\gik}\in J(v)$
}

\AddEq{module, d in D}
{
$\DArg\in \CBase$,
}

\AddEq{module, dv in V}
{%
$\DArg v\in V$,
}

\AddEq{left module, dv in V}
{%
$\DArg v\in V$,
}

\AddEq{right module, dv in V}
{%
$v\DArg \in V$,
}

\AddEq{module, a in A}
{
$a\in \Base$.
}

\AddEq{module, av in V}
{%
$av\in V$
}

\AddEq{left module, av in V}
{%
$av\in V$
}

\AddEq{right module, av in V}
{%
$va\in V$
}

\AddEq{module, dv+av in V}
{%
\[\DArg v+av\in V\ \ \ \DArg\in \CBase\ \ \ a\in \Base\]
}

\AddEq{left module, dv+av in V}
{%
\[\DArg v+av\in V\ \ \ \DArg\in \CBase\ \ \ a\in \Base\]
}

\AddEq{right module, dv+av in V}
{%
\[v\DArg +va\in V\ \ \ \DArg\in \CBase\ \ \ a\in \Base\]
}

\AddEq{module, cvk in Jv}
{%
$c^{\gik}v_{\gik}\in J(v)$, $c^{\gik}\in \Base_{(1)}$,
$\gik\in\giI$
\labelItem{module, cvk in Jv}
}

\AddEquation{module, a+n:V->V}
{
a+n:v\in V\rightarrow (a+n)v\in V
}

\AddEquation{left module, a+n:V->V}
{
a+n:v\in V\rightarrow (a+n)v\in V
}

\AddEquation{right module, a+n:V->V}
{
a+n:v\in V\rightarrow v(a+n)\in V
}

\AddEq{left module, cvk in Jv}
{%
$c^{\gik}v_{\gik}\in J(v)$, $c^{\gik}\in \Base_{(1)}$,
$\gik\in\giI$
\labelItem{left module, cvk in Jv}
}

\AddEq{right module, cvk in Jv}
{%
$v_{\gik}c^{\gik}\in J(v)$, $c^{\gik}\in \Base_{(1)}$,
$\gik\in\giI$
\labelItem{right module, cvk in Jv}
}

\AddEq{module, sum cvk in Jv}
{%
$\displaystyle\sum_{\gik\in\giI}c^{\gik}v_{\gik}\in J(v)$,
$c^{\gik}\in\Base_{(1)}$, $|\{\gii:c^{\gii}\ne 0\}|<\infty$
\labelItem{module, sum cvk in Jv}
}

\AddEq{left module, sum cvk in Jv}
{%
$\displaystyle\sum_{\gik\in\giI}c^{\gik}v_{\gik}\in J(v)$,
$c^{\gik}\in\Base_{(1)}$, $|\{\gii:c^{\gii}\ne 0\}|<\infty$
\labelItem{left module, sum cvk in Jv}
}

\AddEq{right module, sum cvk in Jv}
{%
$\displaystyle\sum_{\gik\in\giI}v_{\gik}c^{\gik}\in J(v)$,
$c^{\gik}\in\Base_{(1)}$, $|\{\gii:c^{\gii}\ne 0\}|<\infty$
\labelItem{right module, sum cvk in Jv}
}

\AddEq[1]{|ci12 ne 0|}
{
H_1=\{\gii\in\giI:#1_1^{\gii}\ne 0\}
\ \ \ \,
H_2=\{\gii\in\giI:#1_2^{\gii}\ne 0\}
}

\AddEq{|wi ne 0|}
{
|\{\gii\in\giI:w^{\gii}\ne 0\}|<\infty
}

\AddEq{module, |awi ne 0|}
{
$\{\gii\in\giI:aw^{\gii}\ne 0\}$
}

\AddEq{left module, |awi ne 0|}
{
$\{\gii\in\giI:aw^{\gii}\ne 0\}$
}

\AddEq{right module, |awi ne 0|}
{
$\{\gii\in\giI:w^{\gii}a\ne 0\}$
}

\AddEq{v,w in Jv=>v+w in Jv}
{
$v$, $w\in J(v)\Rightarrow v+w\in J(v)$
\labelItem{\SideWS module, v,w in Jv=>v+w in Jv}
}

\AddEq{v in Jv=>av in Jv}
{
$a\in\Base_{(1)}$,
$v\in J(v)\Rightarrow av\in J(v)$
\labelItem{\SideWS module, v in Jv=>av in Jv}
}

\AddEq[1]{gi12}
{
$#1^{\gii}_1$, $#1^{\gii}_2$, $\gii\in\giI$,
}

\AddEq[1]{w1+w2= 1 }
{
#1_1+#1_2=\sum_{\gii\in\giI}(#1^{\gii}_1+#1^{\gii}_2)v_{\gii}
}

\AddEq[1]{w1+w2= 1 left}
{
#1_1+#1_2=\sum_{\gii\in\giI}(#1^{\gii}_1+#1^{\gii}_2)v_{\gii}
}

\AddEq[1]{w1+w2= 1 right}
{
#1_1+#1_2=\sum_{\gii\in\giI}v_{\gii}(#1^{\gii}_1+#1^{\gii}_2)
}

\AddEq[1]{|gi 1+2 ne 0|}
{
\[
\{\gii\in\giI:#1_1^{\gii}+#1_2^{\gii}\ne 0\}\subseteq H_1\cup H_2
\]
}

\AddEq[1]{w1+w2= }
{
#1_1+#1_2=\sum_{\gii\in\giI}#1^{\gii}_1v_{\gii}+\sum_{\gii\in\giI}#1^{\gii}_2v_{\gii}
=\sum_{\gii\in\giI}(#1^{\gii}_1v_{\gii}+#1^{\gii}_2v_{\gii})
}

\AddEq[1]{w1+w2= left}
{
#1_1+#1_2=\sum_{\gii\in\giI}#1^{\gii}_1v_{\gii}+\sum_{\gii\in\giI}#1^{\gii}_2v_{\gii}
=\sum_{\gii\in\giI}(#1^{\gii}_1v_{\gii}+#1^{\gii}_2v_{\gii})
}

\AddEq[1]{w1+w2= right}
{
#1_1+#1_2=\sum_{\gii\in\giI}v_{\gii}#1^{\gii}_1+\sum_{\gii\in\giI}v_{\gii}#1^{\gii}_2
=\sum_{\gii\in\giI}(v_{\gii}#1^{\gii}_1+v_{\gii}#1^{\gii}_2)
}

\AddEq[1]{w12= }
{
#1_1=\sum_{\gii\in\giI}#1^{\gii}_1v_{\gii}
\ \ \ \,
#1_2=\sum_{\gii\in\giI}#1^{\gii}_1v_{\gii}
}

\AddEq[1]{w12= left}
{
#1_1=\sum_{\gii\in\giI}#1^{\gii}_1v_{\gii}
\ \ \ \,
#1_2=\sum_{\gii\in\giI}#1^{\gii}_2v_{\gii}
}

\AddEq[1]{w12= right}
{
#1_1=\sum_{\gii\in\giI}v_{\gii}#1^{\gii}_1
\ \ \ \,
#1_2=\sum_{\gii\in\giI}v_{\gii}#1^{\gii}_2
}

\AddEquation{vk=sum vi, module}
{
v_{\gik}=\sum_{\gii\in\giI}c^{\gii}v_{\gii}
}

\AddEquation{vk=sum vi, left module}
{
v_{\gik}=\sum_{\gii\in\giI}c^{\gii}v_{\gii}
}

\AddEquation{vk=sum vi, right module}
{
v_{\gik}=\sum_{\gii\in\giI}v_{\gii}c^{\gii}
}

\AddEq{ci=dik}
{
$c^{\gii}=\delta^{\gii}_{\gik}\in\Base_{(1)}$.
}

\AddEq{vk in v}
{
$v_{\gik}\in v$,
}

\DefEq
{
$v_{\gii}$.
}
{vi}

\AddEq{w=wi vi}
{
$\Vector w=w^{\gii}v_{\gii}$
}

\AddEq{w=wi vileft}
{
$\Vector w=w^{\gii}v_{\gii}$
}

\AddEq{w=wi viright}
{
$\Vector w=v_{\gii}w^{\gii}$
}

\AddEq{w=w cr v}
{
\[\Vector w=\ShowSymbol{linear combination}{\SideWS cr}\]
}

\AddEq{g number linear combination}
{%
\symb{g^is_i}{linear combination}{}%
}

\AddEq{linear combination}
{%
\symb{w^{\gii}v_{\gii}}{linear combination}{i}%
\symb{w\CRstar v}{linear combination}{cr}%
}

\AddEq{left linear combination}
{%
\symb{w^{\gii}v_{\gii}}{linear combination}{left i}%
\symb{w\CRstar v}{linear combination}{left cr}%
}

\AddEq{right linear combination}
{%
\symb{v_{\gii}w^{\gii}}{linear combination}{right i}%
\symb{v\RCstar w}{linear combination}{right cr}%
}

\AddEq{A linear combination =}
{
$\ShowSymbol{linear combination}{\SideWS i}$
}

\DefEq
{
$e_{\gii}$, \iI,
}
{ei, i in I}

\AddEq[1]{c*e=0, module}
{
#1^{\gii}e_{\gii}=0
}

\AddEq[1]{c*e=0, left module}
{
#1^{\gii}e_{\gii}=0
}

\AddEq[1]{c*e=0, right module}
{
e_{\gii}#1^{\gii}=0
}

\AddEq{(b+c)*e, module}
{
\[
(b^{\gii}+c^{\gii})e_{\gii}=0
\]
}

\AddEq{(b+c)*e, left module}
{
\[
(b^{\gii}+c^{\gii})e_{\gii}=0
\]
}

\AddEq{(b+c)*e, right module}
{
\[
e_{\gii}(b^{\gii}+c^{\gii})=0
\]
}

\AddEq{(ac)*e, module}
{
\[
(ac^{\gii})e_{\gii}=0
\]
}

\AddEq{(ac)*e, left module}
{
\[
(ac^{\gii})e_{\gii}=0
\]
}

\AddEq{(ac)*e, right module}
{
\[
e_{\gii}(c^{\gii}a)=0
\]
}

\DefEq
{
\symb{f+g}{sum of maps}{D module}
}
{sum of maps,,D module}

\DefEquation
{
\begin{matrix}
\ShowSymbol{sum of maps}{D module}:A_1\rightarrow A_2
&f,g\in\mathcal L(D;A_1\rightarrow A_2)
\end{matrix}
}
{sum of maps, 1, D module}

\DefEquation
{
(\ShowSymbol{sum of maps}{D module})\circ x=f\circ x+g\circ x
}
{sum of maps, D module}

\AddEq [3]{f in L(A->B)}
{
$f\in\mathcal L(D;#1\rightarrow #2)$#3
}

\DefEq
{
$F_k\in\Basis F$,
}
{Ik in I}

\DefEquation
{
x\rightarrow I_k\circ a\circ x
}
{x->Ik a x}

\DefEquation
{
b^l=I^l_k\circ a
}
{b=Ikl a}

\DefEq
{
\begin{align*}
(I_k^l\circ(a_1+a_2))\circ I_l\circ x
&=I_k\circ(a_1+a_2)\circ x
\\&=I_k\circ a_1\circ x+I_k\circ a_2\circ x
\\&=(I_k^l\circ a_1)\circ I_l\circ x+(I_k^l\circ a_1)\circ I_l\circ x
\end{align*}
}
{Ikl a1+a2}

\DefEq
{
\begin{align*}
(I_k^l\circ(da))\circ I_l\circ x
&=I_k\circ(da)\circ x
=I_k\circ(d(a\circ x))
\\&=d(I_k\circ a\circ x)
=d((I_k^l\circ a)\circ I_l\circ x)
\\&=(d(I_k^l\circ a))\circ I_l\circ x
\end{align*}
}
{Ikl d a}

\DefEquation
{
I_k\circ a\circ x=b^l\circ I_l\circ x\ \ \ \ b^l\in A_2\times A_2
}
{expansion Ik a x}

\DefEq
{
\symb{F_k^l}{conjugation transformation}{}
}
{conjugation transformation}

\DefEq
{
\[
\ShowSymbol{conjugation transformation}{}
:A_1\otimes A_1\rightarrow A_2\otimes A_2
\]
}
{conjugation transformation:}

\DefEquation
{
F_k\circ a\circ x=(\ShowSymbol{conjugation transformation}{}\circ a)\circ F_l\circ x
}
{conjugation transformation =}

\DefEq
{
$\ShowSymbol{conjugation transformation}{}$
}
{show conjugation transformation}

\DefEq
{
$C_{\gi{kl}}^{\gi p}$
}
{structural constants, algebra}

\AddEq{wi i in I}
{
$w^{\gii}$, $\gii\in\giI$,
}

\AddEq{w=wi i in I}
{
$w=(w^{\gii},\gii\in\giI)$
}

\AddEq{i=j}
{
$\gii=\gij$
}

\DefEq
{
$f^{\gik}_{\gil}$
}
{Coordinates of map f}

\DefEq
{
$f^{k\cdot\gi{ij}}$
}
{standard components of map f}

\DefEq
{
$F_{k\cdot}^{}{}_{\gii}^{\gij}$
}
{coordinates of map Ik}

\DefEquation
{
f^{\gik}_{\gil}=f^{k\cdot\gi{ij}}
F_{k\cdot}^{}{}^{\gim}_{\gil}
C^{\gi p}_{\gi{im}}C^{\gik}_{\gi{pj}}
}
{coordinates of map A1 A2, 2, associative algebra}

\DefEquation
{
F_k\circ x=F_{k\cdot}^{}{}^{\gii}_{\gij}x^{\gij}e_{2\cdot\gii}
}
{coordinates of map Ik, associative algebra}

\DefEquation
{
f\circ x=f^{\gii}_{\gij}x^{\gij}e_{2\cdot\gii}
}
{coordinates of map f, associative algebra}

\DefEq
{
\[
\Basis F=\{F_k\in\mathcal L(D;A_1\rightarrow A_2): k=1, ..., n\}
\]
}
{Ik 1n}

\DefEquation
{
h=
(
a\pC{s}{0}
\otimes
a\pC{s}{1}
)
\circ f
=a\pC{s}{0}
f
a\pC{s}{1}
}
{h generated by f, associative algebra}

\DefEq
{
\[
h:A_1\rightarrow A_2
\]
}
{h:A1->A2}

\DefEquation
{
f=
(
a_{k\cdot s_k\cdot 0}
\otimes
a_{k\cdot s_k\cdot 1}
)
\circ I_k
=
\sum_{\gik}
a_{k\cdot s_k\cdot 0}
I_k
a_{k\cdot s_k\cdot 1}
}
{f in L(A,A), 1, associative algebra}

\DefEquation
{
f=
a^{k\cdot \gi{ij}}(e_{\gii}\otimes e_{\gij})\circ I_k
=
a^{k\cdot \gi{ij}}e_{\gii}I_ke_{\gij}
}
{f in L(A,A), 2, associative algebra}

\DefEq
{
\[
f\in B^A\rightarrow \|f\|\in R
\]
}
{f in BA->|f|}

\DefEq
{
\symb{\|f\|}{norm of map}{}
}
{norm of map}

\DefEquation
{
\ShowSymbol{norm of map}{}=
\text{sup}\frac{\|f(x)\|_B}{\|x\|_A}
}
{norm of map, algebra}

\DefEq
{
$\ShowSymbol{norm of map}{}$
}
{show|f|}

\DefEq
{
\symb{f+g}{sum of maps}{polylinear}
}
{sum of maps,,polylinear}

\DefEquation
{
\begin{matrix}
\ShowSymbol{sum of maps}{polylinear}:A_1\times...\times A_n\rightarrow S
&f,g\in\mathcal L(D;A_1\times...\times A_n\rightarrow S)
\end{matrix}
}
{sum of maps, 1, polylinear}

\DefEquation
{
f^k=f^k_{s_k\cdot 0}\otimes f^k_{s_k\cdot 1}
}
{f=fkxfk}

\DefEquation
{
f=f^k\circ I_k\ \ \ f^k\in\ATwo
}
{expansion of linear map with respect to basis}

\DefEquation
{
(\ShowSymbol{sum of maps}{polylinear})\circ (a_1,...,a_n)
=f\circ (a_1,...,a_n)+g\circ (a_1,...,a_n)
}
{sum of  maps, polylinear}

\DefEq
{
$\mathcal L(D;A_1\times...\times A_n\rightarrow S)$
}
{module of polylinear maps}

\DefEq
{
\symb{\mathcal L(D;A_1\rightarrow A_2)}{set linear maps}1
}
{set linear maps, module}

\AddEq{linear map times constant, algebra}
{
$af$, $fb$, $a$, $b\in A_2$,
}

\AddEq{linear map times constant, 0, algebra}
{
\begin{align*}
(af)\circ x&=a(f\circ x)
\\
(fb)\circ x&=(f\circ x)b
\end{align*}
}

\AddEq{linear map AA LAA}
{
\[
h:\ATwo\rightarrow \mathcal L(D;A_1\rightarrow A_2)
\]
}

\DefEquation
{
(a\otimes b)\circ f=afb
}
{linear map AA LAA, 1}

\DefEq
{
\begin{align}
f\circ(a+b)&=f\circ a+f\circ b
\EqLabel{linear map from A1 to A2, 1 1}
\\
f\circ(pa)&=p(f\circ a)
\EqLabel{linear map from A1 to A2, 1 2}
\end{align}
\[
\begin{matrix}
a,b\in A_1
&
p\in D
\end{matrix}
\]
}
{linear map from A1 to A2, 1}

\DefEq
{
\symb{[a,b]}{commutator of algebra}{}
\[
\ShowSymbol{commutator of algebra}{}=ab-ba
\]
}
{commutator of algebra}

\DefEq
{
\[
[a,b]=0
\]
}
{commutative D algebra}

\DefEq
{
\symb{(a,b,c)}{associator of algebra}{}
\begin{equation}
\ShowSymbol{associator of algebra}{}=(ab)c-a(bc)
\EqLabel{associator of algebra}
\end{equation}
}
{associator of algebra}

\DefEq
{
\[
(a,b,c)=0
\]
}
{associative D algebra}

\AddEq{nucleus of algebra}
{
\symb{N(A)}{nucleus of algebra}{}
\[
\ShowSymbol{nucleus of algebra}{}=
\{
a\in A:
\forall b, c\in A,
(a,b,c)=(b,a,c)=(b,c,a)=0
\}
\]
}

\AddEq{center of algebra}
{
\symb{Z(A)}{center of algebra}{}
\[
\ShowSymbol{center of algebra}{}=
\{
a\in A:
a\in N(A),
\forall b\in A,
ab=ba
\}
\]
}

\AddEq{f(a+b)=...}
{
f\circ(a+b)=f\circ a+f\circ b
}

\DefEquation
{
\Vector r_2(da)=d\Vector r_2(a)
}
{D r2(da)=... left}

\DefEq
{
\[\Vector r_2:A_1\rightarrow A_2\]
}
{r2:A1->A2}

\AddEq{D ab in A, d in D}
{
\[
\begin{matrix}
a,b\in A_1
&
d\in D
\end{matrix}
\]
}

\AddEq{D1 D2 ab in A, d in D}
{
\[
\begin{matrix}
a,b\in A_1
&d,d_1,d_2\in D_1
\end{matrix}
\]
}

\AddEq{property linear map hf}
{
\DrawEq{h(d1+d2)=...}{\SideWS module}
\DrawEq{h(d1d2)=...}{\SideWS module}
\DrawEq{f(a+b)=...}{D1 D2 \SideWS module}
\ShowEq{\DFDT, f(da)=... \SideWS module}
\ShowEq{D1 D2 ab in A, d in D}
}

\AddEq[1]{show linear map hf}
{
\ShowEq{h:D1->D2 f:A1->A2}h \Base {#1} \Module
}

\AddEq[1]{show linear map f}
{
\ShowEq{f:A->B}{#1}{\Module_1}{\Module_2}
}

\AddEq{w in Jv}
{
$w\in J(v)$.
}

\AddEquation{D1 D2 , f(da)=... module}
{
f\circ(da)=(h\circ d)(f\circ a)
}

\AddEquation{D1 D2 , f(da)=... left module}
{
f\circ(da)=(h\circ d)(f\circ a)
}

\AddEquation{D1 D2 , f(da)=... right module}
{
f\circ(ad)=(f\circ a)(h\circ d)
}

\AddEquation{D , f(da)=... module}
{
f\circ(da)=(h\circ d)(f\circ a)
}

\AddEquation{D , f(da)=... left module}
{
f\circ(da)=(h\circ d)(f\circ a)
}

\AddEquation{D , f(da)=... right module}
{
f\circ(ad)=(f\circ a)(h\circ d)
}

\AddEq{h(d1+d2)=...}
{
h\circ(d_1+d_2)=h\circ d_1+h\circ d_2
}

\AddEq{h(d1d2)=...}
{
h\circ(d_1d_2)=(h\circ d_1)(h\circ d_2)
}

\AddEq{property linear map f}
{
\DrawEq{f(a+b)=...}{D }
\ShowEq{D r2(da)=... left}
\ShowEq{D ab in A, d in D}
}

\DefEq
{
\symb{\Basis{e}=(e_{\gii},\Ii)}{basis, module}1
}
{basis, module}

\DefEq
{
\[\Basis e_1=(e_{1\cdot\gii},\gii\in\giI)\]
}
{basis, module V1}

\DefEq
{
\[\Basis e_2=(e_{2\cdot\gij},\gij\in\giJ)\]
}
{basis, module V2}

\AddEquation{r2:A1->A2, left D1 D2 module}
{
b=(h\circ a)\CRstar f
}

\AddEquation{r2:A1->A2, D1 D2 module}
{
b=(h\circ a)\CRstar f
}

\AddEquation{r2:A1->A2, right D1 D2 module}
{
b=f\RCstar (h\circ a)
}

\AddEq{r2:A1->A2, left D module}
{
b=a\CRstar f
}

\AddEq{r2:A1->A2, D module}
{
b=a\CRstar f
}

\AddEq{r2:A1->A2, right D module}
{
b=f\RCstar a
}

\AddEq [2]{va=ae1, left module}
{
\Vector #1=#1\CRstar e_{#2}
}

\AddEq [2]{va=ae1, module}
{
\Vector #1=#1\CRstar e_{#2}
}

\AddEq [2]{va=ae1, right module}
{
\Vector #1=e_{#2}\RCstar #1
}

\AddEq{product in D algebra}
{
v\,w=C\circ(v,w)
}

\AddEq{product in algebra, definition 1}
{
\[
C:A\times A\rightarrow A
\]
}

\DefEq
{
\[
\underline{\otimes}:A^{\otimes n}\times A^{\otimes m}
\rightarrow A^{\otimes(n+m-1)}
\]
}
{otimes -}

\AddEq{a b in basis of algebra}
{
\[
\begin{matrix}
a=a^{\gii}e_{\gii}&b=b^{\gii}e_{\gii}&a, b\in A
\end{matrix}
\]
}

\AddEquation{product in algebra}
{
(ab)^{\gik}=C^{\gik}_{\gi{ij}}a^{\gii}b^{\gij}
}

\AddEquation{product of basis vectors, algebra}
{
e_{\gii} e_{\gij}=
C^{\gik}_{\gi{ij}}e_{\gik}
}

\AddEq[4]{diagram of representations of D algebra}
{
\begin{matrix}
\vcenter
{
\xymatrix
{
D_{#1}\ar[r]|(.4){*}^{#4_{#2 1,2}}&
A_{#3}\ar[r]|{*}^{#4_{#2 2,3}}&
A_{#3}
\\
&&D_{#1}\ar[u]|{*}_{#4_{#2 1,2}}
}
}
&
\begin{array}{r@{\,}l}
#4_{#2 1,2}(d):v&\rightarrow d\,v
\\
#4_{#2 2,3}(v):w&\rightarrow
C_{#3}\circ(v, w)
\\
C_{#3}&\in\mathcal L(D_{#1};A_{#3}^2\rightarrow A_{#3})
\end{array}
\end{matrix}
}

\AddEq{structural constants of algebra}
{
\symb{C^{\gik}_{\gi{ij}}}{structural constants}1
}

\DefEquation
{
(a_1\otimes...\otimes a_n)\underline{\otimes}(b_1\otimes...\otimes b_n)
=
a_1\otimes...\otimes a_{n-1}\otimes a_nb_1\otimes b_2\otimes...\otimes b_n
}
{otimes -, 1}

\AddEq{vb=f(va)}
{
\Vector b=\Vector f\circ\Vector a
}

\AddEquation{D1 D2 vb=h(e1)a }
{
\Vector b=\Vector f\circ\Vector a=\Vector f\circ(a\CRstar e_1)
= (h\circ a)\CRstar(\Vector f\circ e_1)
}

\AddEquation{D1 D2 vb=h(e1)a left}
{
\Vector b=\Vector f\circ\Vector a=\Vector f\circ(a\CRstar e_1)
= (h\circ a)\CRstar(\Vector f\circ e_1)
}

\AddEquation{D1 D2 vb=h(e1)a right}
{
\Vector b=\Vector a\diamond\Vector f=(e_1\RCstar a)\diamond\Vector f
= (e_1\diamond\Vector f)\RCstar(h\circ a)
}

\DefEq
{
}
{}

\AddEq{f(e1) }
{
$\Vector f\circ e_{1\cdot\gii}$
}

\AddEq{f(e1) left}
{
$\Vector f\circ e_{1\cdot\gii}$
}

\AddEq{f(e1) right}
{
$e_{1\cdot\gii}\diamond\Vector f$
}

\AddEq{f(e1)= }
{
\Vector f\circ e_{1\cdot\gii}
=f_{\gii}\CRstar e_2
=f_{\gii}^{\gij}e_{2\cdot\gij}
}

\AddEq{f(e1)= left}
{
\Vector f\circ e_{1\cdot\gii}
=f_{\gii}\CRstar e_2
=f_{\gii}^{\gij}e_{2\cdot\gij}
}

\AddEq{f(e1)= right}
{
e_{1\cdot\gii}\diamond\Vector f
=e_2\RCstar f_{\gii}
=e_{2\cdot\gij}f_{\gii}^{\gij}
}

\AddEquation{D1 D2 vb=a f e left}
{
\Vector b=(h\circ a)\CRstar f\CRstar e_2
}

\AddEquation{D1 D2 vb=a f e }
{
\Vector b=(h\circ a)\CRstar f\CRstar e_2
}

\AddEquation{D1 D2 vb=a f e right}
{
\Vector b=e_2\RCstar f\RCstar (h\circ a)
}

\DefEq
{
$h\circ a=(h\circ a_{\gii},\gii\in\giI)$
}
{h(a)=...}

\AddEq{Vector f(e1) module}
{
$(\Vector f\circ e_{1\cdot\gii},\gii\in\giI)$
}

\AddEq{Vector f(e1) left module}
{
$(\Vector f\circ e_{1\cdot\gii},\gii\in\giI)$
}

\AddEq{Vector f(e1) right module}
{
$(e_{1\cdot\gii}\diamond\Vector f,\gii\in\giI)$
}

%% file: Stmt.Module.Ref.tex

\AddEq{ref definition: basis of representation}
{
\ePrints{1502.04063,5114-6019,MVector}
\ifx\Semafor\ValueOn
\RefDefinition[0912.3315]{basis of representation}
\else
\RefDefinition{basis of representation}
\fi
}

\AddEq{ref module, (a+n)(b+m)=}
{
\EqRef{\SideWS module, a d v},
\EqRef{\SideWS module, (a+n)v=av+nv},
\eqRef{associative law, \SideWS module}1,
\eqRef{associative law, \CBase, \SideWS module}1,
\eqRef{\SideWS module, associative law}1,
\eqRef{associative law, \CBase\Base, \SideWS module}1.
}

\AddEq{ref map j, representation, tensor product}
{
\ePrints{MVector}
\ifx\Semafor\ValueOn
\EqRef[\RefPolymorphism]{map j, representation, tensor product},
\else
\EqRef{map j, representation, tensor product},
\fi
}

%% file: Stmt.Omega.English.tex
\input{Stmt.Omega.Eq}

\DefDefinition{polyadditive map}
{
A map
\ShowEq{f:A->B}g{A^n}A
is called
\AddIndex{polyadditive map}{polyadditive map}
if for any
\ShowEq{fi(a+b)=}
}

\DefDefinition{semiring of sets}
{
A nonempty system of sets $\mathcal S$ is called
\AddIndex{semiring of sets}{semiring of sets},\,\footnote{
See also the definition
\citeBib{Kolmogorov Fomin}\Hyph 2,  page 32.}
if
\ifx\setCACAA\undefined
\StartLabelItem
\else
\StartLabelItem[definition]
\fi
\begin{enumerate}
\item
\ShowEq{empty in S}
\item
If
\ShowEq{A, B in S}
then
\ShowEq{AB in S}
\item
If
\ShowEq{A, A1 in S},
then the set \(A\) can be represented as
\DrawEq{A=+Ai}{def}
where
\ShowEq{Ai*Aj=0}
\labelItem{A=+Ai}
\end{enumerate}
The representation
\eqRef{A=+Ai}{def}
of the set \(A\) is called
\AddIndex{finite expansion of set}{finite expansion of set}
\(A\).
}

\DefDefinition{norm on Omega group}
{
\AddIndex{Norm on $\Omega$\Hyph group}
{norm on Omega group} $A$\,\footnote{
I made definition according\,to definition
from\,\citeBib{Bourbaki: General Topology: Chapter 5 - 10},
IX,\,\S 3.2 and definition \citeBib{Arnautov Glavatsky Mikhalev}-1.1.12,
p. 23.} is a map
\ShowEq{d->|d|}
which satisfies the following axioms
\ifx\setCACAA\undefined
\StartLabelItem
\else
\StartLabelItem[definition]
\fi
\begin{enumerate}
\ShowEq{|a|>=0}
\ShowEq{|a|=0}
if, and only if, $a=0$
\ShowEq{|a+b|<=|a|+|b|}
\ShowEq{|-a|=|a|}
\end{enumerate}

The $\Omega$\Hyph group $A$,
endowed with the structure defined by a given norm on $A$, is called
\AddIndex{normed $\Omega$\Hyph group}{normed Omega group}.
}

\DefTheorem{h=f1...fn omega, converges uniformly}
{
Let $A$ be complete $\Omega$\Hyph group.
Let
\ShowEq{omega in Omega}{}{}
be $n$\Hyph ari operation.
Let sequence of maps
\ShowEq{fim M(X,A)}
into complete $\Omega$\Hyph group $A$
converge uniformly
to the map $f_i$.
Let range of the map
$f_i$
be compact set.
Then sequence of maps
\ShowEq{hm=f1m...fnm omega}
into complete $\Omega$\Hyph group $A$
converges uniformly
to the map
\DrawEq{h=f1...fn omega}{converges uniformly}
}

\DefTheorem{|int h|<int|omega||f1n|}
{
Let
\ShowEq{omega in Omega}{}{}
be $n$\Hyph ary operation in
Abelian \(\Omega\)\Hyph group \(A\).
Let
\ShowEq{fi:X->A}
be $\mu$\Hyph measurable map
with compact range.
Since map \(f_i\), \(i=1\), ..., \(n\),
is integrable map, then map
\DrawEq{h=f1...fn omega}{}
is integrable map and
\DrawEq{|int h|<int|omega||f1n|}{measurable}
}

\DefTheorem{h=fX g1 g2 representation}
{
Let $\mu$ be a \(\sigma\)\Hyph additive measure defined on the set $X$.
Let
\ShowEq{f:A->*B}g{A_1}{A_2}
be representation
of $\Omega_1$\Hyph group $A_1$ with norm $\|x\|_1$
in $\Omega_2$\Hyph group $A_2$ with norm $\|x\|_2$.
Let
\ShowEq{gi:X->A 12}
be integrable map
with compact range.
Then map
\ShowEq{h=fX g1 g2}
is integrable map and
\DrawEq{|int h|<int|g1||g2|}{measurable}
}

\DefTheorem{polymorphism of representations}
{
Let the map
\ShowEq{map r1n,R}
be polymorphism of
representations $f_1$, ..., $f_n$ into representation $f$.
For any
\ShowEq{k,1n}
the map
\ShowEq{map r,R}{r_{1k}}{r_2}{}
satisfies to the equality
\ShowEq{polymorphism of representation, ak}
Let $\omega_1\in\Omega_1(p)$.
For any
\ShowEq{k,1n}
the map $r_{1k}$ satisfies to the equality
\ShowEq{polymorphism of representation, omega1}
Let $\omega_2\in\Omega_2(p)$.
For any
\ShowEq{k,1n}
the map $r_2$ satisfies to the equality
\ShowEq{polymorphism of representation, omega2}
}

\DefDefinition{polymorphism of representations}
{
Let
\ShowEq{A1n A}
be $\Omega_1$\Hyph algebras.
Let
\ShowEq{B1n}
be $\Omega_2$\Hyph algebras.
Let, for any
\ShowEq{k,1n}
\ShowEq{f:A->*B}{f_k}{A_k}{B_k}
be representation of $\Omega_1$\Hyph algebra $A_k$
in $\Omega_2$\Hyph algebra $B_k$.
Let
\ShowEq{f:A->*B}fAB
be representation of $\Omega_1$\Hyph algebra $A$
in $\Omega_2$\Hyph algebra $B$.
The map
\ShowEq{polymorphism of representation}
is called
\AddIndex{polymorphism of representations}
{polymorphism of representations}
$f_1$, ..., $f_n$ into representation $f$,
if, for any
\ShowEq{k,1n}
provided that all variables except variables
\ShowEq{ak bk in}
have given value, the map
\ShowEq{map r,R}{r_{1k}}{r_2}{}
is a morphism of representation $f_k$ into representation $f$.

If $f_1=...=f_n$, then we say that the map
\ShowEq{map r1n,R}
is polymorphism of representation $f_1$ into representation $f$.

If $f_1=...=f_n=f$, then we say that the map
\ShowEq{map r1n,R}
is polymorphism of representation $f$.
}

\DefDefinition{Omega group}
{
Let sum which is not necessarily commutative
be defined
in $\Omega_1$\Hyph algebra $A$.
We use the symbol $+$ to denote sum.
Let
\ShowEq{Omega=...-+}
If $\Omega_1$\Hyph algebra $A$ is group
relative to sum
and any operation $\omega\in\Omega$
is polyadditive map,
then $\Omega_1$\Hyph algebra $A$ is called
\AddIndex{$\Omega$\Hyph group}{Omega group}.
If $\Omega$\Hyph group $A$ is associative group
relative to sum,
then $\Omega_1$\Hyph algebra\,$A$\,is\,called
\AddIndex{associative $\Omega$\Hyph group}{associative Omega group}.
If $\Omega$\Hyph group $A$ is Abelian group
relative\,to sum,
then $\Omega_1$\Hyph algebra\,$A$ is\,called
\AddIndex{Abelian $\Omega$\Hyph group}{Abelian Omega group}.%
}

\DefTheorem{operation is distributive over addition}
{
Let
\EqParm{omega n ari}{=z}
be polyadditive map.
The operation $\omega$
\AddIndex{is distributive}{distributive law}
over addition
\ShowEq{omega is distributive over addition}
}

\DefDefinition{multiplicative Omega group}
{
Let product
\DrawEq{c=ab}{}
be operation of $\Omega_1$\Hyph algebra $A$.
Let
\ShowEq{Omega=...-o}
If $\Omega_1$\Hyph algebra $A$ is group
relative to product
and, for any operation
\EqParm{omega n ari}{=c}
the product is distributive over the operation $\omega$
\DrawEq{a omega=omega a left}{}
\DrawEq{a omega=omega a right}{}
then $\Omega_1$\Hyph algebra $A$ is called
\AddIndex{multiplicative $\Omega$\Hyph group}
{multiplicative Omega group}.
}

\DefDefinition{Abelian multiplicative Omega group}
{
If
\ShowEq{Abelian multiplicative Omega group}
then multiplicative $\Omega$\Hyph group is called
\AddIndex{Abelian}{Abelian multiplicative Omega group}.
}

\DefDefinition{Omega ring}
{
Let sum
\ShowEq{c1=a+b}
which is not necessarily commutative
and product
\DrawEq{c=ab}{}
be operations of $\Omega_1$\Hyph algebra $A$.
Let
\ShowEq{Omega=...-o+}
If $\Omega_1$\Hyph algebra $A$ is
$\Omega\cup\{*\}$\Hyph group
and
multiplicative $\Omega\cup\{+\}$\Hyph group,
then $\Omega_1$\Hyph algebra $A$ is called
\AddIndex{$\Omega$\Hyph ring}{Omega ring}.
}

\DefTheorem{product in Omega ring is distributive over addition}
{
The product in $\Omega$\Hyph ring
\AddIndex{is distributive}{distributive law}
over addition
\ShowEq{a*(b1+b2)=}
}

\DefProof{product in Omega ring is distributive over addition}
{
The theorem follows from the definitions
\RefDefinition{Omega group},
\RefDefinition{multiplicative Omega group},
\RefDefinition{Omega ring}.
}

\DefDefinition{matrix over Omega ring}
{
Let $A$ be $\Omega$\Hyph ring. The
\AddIndex{matrix}{matrix}
over $\Omega$\Hyph ring $A$ is a table of $A$\Hyph numbers $a^i_j$
where the index $i$ is the number of row
and the index $j$ is the number of column.
}

%% file: Stmt.Omega.Eq.tex

\AddEq{f(a+b)=}
{
\[
f(a+b)=f(a)+f(b)
\]
}

\AddEq{fi(a+b)=}
{
$i$, $i=1$, ..., $n$,
\[
f(a_1,...,a_i+b_i,...,a_n)=f(a_1,...,a_i,...,a_n)+f(a_1,...,b_i,...,a_n)
\]
}

\AddEq{polymorphism of representation}
{
\[
\begin{pmatrix}
r_{1k}:A_k\rightarrow A
&
k=1,...,n
&
r_2:B_1\times...\times B_n\rightarrow B
\end{pmatrix}
\]
}

\AddEq{map r1n,R}
{
$((r_{1,k}, k=1,...,n)\ \ r_2)$
}

\AddEquation{polymorphism of representation, ak}
{
r_2(m_1,...,\BlueText{f_k(a_k)(m_k)},...,m_n)
=f(\RedText{r_{1k}(a_k)})(\BlueText{r_2(m_1,...,m_n)})
}

\AddEquation{reduced polymorphism of representation, ak}
{
r_2(m_1,...,\BlueText{f_k(a)(m_k)},...,m_n)
=f(a)(\BlueText{r_2(m_1,...,m_n)})
}

\AddEquation{polymorphism of representation, omega1}
{
r_{1k}(a_{k\cdot 1}...a_{k\cdot p}\omega_1)
=
r_{1k}(a_{k\cdot 1})
...
r_{1k}(a_{k\cdot p})\omega_1
}

\AddEquation{polymorphism of representation, omega2}
{
\begin{array}{r@{\,}l}
&\,r_2(m_1,...,m_{k\cdot 1}...m_{k\cdot p}\omega_2,...,m_n)
\\
=&\,
r_2(m_1,...,m_{k\cdot 1},...,m_n)
...
r_2(m_1,...,m_{k\cdot p},...,m_n)
\omega_2
\end{array}
}

\AddEq{gi:X->A 12}
{
\[
\begin{matrix}
g_i:X\rightarrow A_i&i=1,2
\end{matrix}
\]
}

\AddEq{Ai*Aj=0}
{
\(i\ne j=>A_i\cap A_j=\emptyset\)
}

\AddEq{A=+Ai}
{
A=\bigcup_{i=1}^nA_i\ \ \ A_i\in\mathcal S
}

\AddEq{A, A1 in S}
{
\(A\), \(A_1\in\mathcal S\),
\(A_1\subset A\)%
}

\AddEq{AB in S}
{
\(A\cap B\in\mathcal S\)
\labelItem{AB in S}
}

\AddEq{A, B in S}
{
\(A\), \(B\in\mathcal S\),
}

\AddEq{empty in S}
{
\(\emptyset\in\mathcal S\)
\labelItem{empty in S}
}

\AddEq{|-a|=|a|}
{
\item $\|-a\|=\|a\|$
\labelItem{|-a|=|a|}
}

\AddEq{|a+b|<=|a|+|b|}
{
\item $\|a+b\|\le\|a\|+\|b\|$
\labelItem{|a+b|<=|a|+|b|}
}

\AddEq{d->|d|}
{
\[d\in A\rightarrow \|d\|\in R\]
}

\AddEq{|a|=0}
{
\item $\|a\|=0$
\labelItem{|a|=0}
}

\AddEq{|a|>=0}
{
\item $\|a\|\ge 0$
\labelItem{|a|>=0}
}

\AddEq{fim M(X,A)}
{
$f_{i\cdot m}\in M(X,A)$, $i=1$, ..., $n$, $m=1$, ...,
}

\AddEq{|int h|<int|omega||f1n|}
{
\left\|\int_Xd\mu(x) h(x)\right\|\le
\int_Xd\mu(x)(\|\omega\|\|f_1(x)\|...\|f_n(x)\|)
}

\AddEq{fi:X->A}
{
\[
\begin{matrix}
f_i:X\rightarrow A&i=1,...,n
\end{matrix}
\]
}

\AddEq{|int h|<int|g1||g2|}
{
\left\|\int_Xd\mu(x) h(x)\right\|_2\le
\int_Xd\mu(x)(\|f\|\|g_1(x)\|_1\|g_2(x)\|_2)
}

\AddEq{h=fX g1 g2}
{
\[h=f_X(g_1)(g_2)\]
}

\AddEq{ak bk in}
{
$a_k\in A_k$, $b_k\in B_k$
}

\AddEq{B1n}
{
$B_1$, ..., $B_n$, $B$
}

\AddEq{A1n A}
{
$A_1$, ..., $A_n$, $A$
}

\AddEq{omega is distributive over addition}
{
\[
a_1...(a_i+b_i)...a_n\omega=a_1...a_i...a_n\omega+a_1...b_i...a_n\omega
\ \ \ i=1, ..., n
\]
}

\AddEq{Omega=...-o+}
{
$\Omega=\Omega_1\setminus\{+,*\}$.
}

\AddEq{c=ab}
{
c_1=a_1* b_1
}

\AddEq{a*(b1+b2)=}
{
\begin{align*}
a*(b_1+b_2)&=a*b_1+a*b_2\\
(b_1+b_2)*a&=b_1*a+b_2*a
\end{align*}
}

\AddEq{f(ab)=f(a)f(b)}
{
f(a* b)=f(a)\circ f(b)
}

\AddEq{c1=a+b}
{
\[c_1=a_1+ b_1\]
}

\AddEquation{Abelian multiplicative Omega group}
{
a* b=b* a
}

\AddEq{a * b + =}
{
a_{11}a_{21}+a_{12}a_{22}
=
(a_{11}+a_{12})(a_{21}+a_{22})
}

\DefEq
{
\begin{align*}
(a_{11}+a_{12})(a_{21}+a_{22})
&=(a_{11}+a_{12})a_{21}+(a_{11}+a_{12})a_{22}
\\
&=a_{11}a_{21}+a_{12}a_{21}
+a_{11}a_{22}+a_{12}a_{22}
\end{align*}
}
{a+a a+a}

\AddEq{a omega=omega a right}
{
(b_1...b_n\omega)* a=(b_1* a)...(b_n* a)\omega
}

\AddEq{a omega=omega a left}
{
a*(b_1...b_n\omega)=(a* b_1)...(a* b_n)\omega
}

\AddEq{Omega=...-o}
{
$\Omega=\Omega_1\setminus\{*\}$.
}

\AddEq{Omega=...-+}
{
\[\Omega=\Omega_1\setminus\{+\}\]
}

%% file: Stmt.Derivative.English.tex
\input{Stmt.Derivative.Eq}

\DefTheorem{bilinear map and differential}
{
Let $A$ be Banach $D$\Hyph module.
Let $B_1$, $B_2$, $B$ be Banach $D$\Hyph algebras.
Let
\ShowEq{h:BxB->B}
be continuous bilinear map.
Let $f$, $g$ be differentiable maps
\ShowEq{f,g:A->B12}
The map
\ShowEq{h(f,g):A->B}
is differentiable and
the derivative satisfies to relationship
\ShowEq{bilinear map and derivative, 3}
\ShowEq{bilinear map and derivative, 5}
}

\DefTheorem{derivative of product, algebra}
{
Let $A$ be Banach $D$\Hyph module.
Let $B$ be Banach $D$\Hyph algebra.
Let $f$, $g$ be differentiable maps
\ShowEq{f,g:A->B}
The derivative satisfies to relationship
\ShowEq{differential of product, algebra}
\ShowEq{derivative of product, algebra}
}

\DefDefinition{norm of polylinear map}
{
Let $A$ be Banach $D$\Hyph module with norm $\|x\|_A$.
Let $B$ be Banach $D$\Hyph module with norm $\|x\|_B$.
For map
\ShowEq{f:An->B}
the value
\ShowEq{norm of polymap}
\ShowEq{norm of map, module}
is called
\AddIndex{norm of map}{norm of map} $f$.
}

\DefTheorem{|f(a)|<|f||a| 1n}
{
For map
\ShowEq{f:An->B}
of Banach $D$\Hyph module $A$ with norm $\|x\|_A$
to Banach $D$\Hyph module $B$ with norm $\|x\|_B$
\ShowEq{|f(a)|<|f||a| 1n}
}

\DefProof{|f(a)|<|f||a| 1n}
{
The inequality
\ShowEq{|f(a)|/|a|<max}
follows from the equality
\EqRef{norm of map, module}.
The inequality
\EqRef{|f(a)|<|f||a| 1n}
follows from the inequality
\EqRef{|f(a)|/|a|<max}.
}

\DefTheorem{|on|->0 ona1p->0}
{
Let
\ShowEq{f:A->B}{o_n}{A^n}B
be sequence of maps of Banach $D$\Hyph module $A$
into Banach $D$\Hyph module $B$ such that
\ShowEq{|on|->0}
Then, for any $B$\Hyph numbers
\ShowEq{a1n}p,
\ShowEq{ona1p->0}
}

\DefProof{|on|->0 ona1p->0}
{
From the theorem
\RefTheorem{|f(a)|<|f||a| 1n},
it follows that
\ShowEq{|oa1p|<}
The statement
\EqRef{ona1p->0}
follows from statements
\EqRef{|on|->0},
\EqRef{|oa1p|<}.
}

\DefTheorem{derivative x2}
{
Let $D$ be the complete commutative ring of characteristic $0$.
Let $A$ be associative Banach $D$\Hyph algebra.
Then\,\footnote{
The statement of the theorem is similar to example VIII,
\citeBib{Hamilton Elements of Quaternions 1}, p. 451.
If product is commutative, then the equality
\EqRef{derivative x2, algebra}
gets form
\ShowEq{derivative x2, field}}
\ShowEq{derivative x2, algebra}
\ShowEq{derivative x2 differential, algebra}
\ShowEq{derivative x2 component, algebra}
}

\DefProof{derivative x2}
{
Consider increment of map $f(x)=x^2$.
\ShowEq{derivative x2, algebra, 1}
The equality
\EqRef{derivative x2 differential, algebra}
follows from the equality
\EqRef{derivative x2, algebra, 1}
and the definition
\RefDefinition{differentiable map, algebra}.
The equality
\EqRef{derivative x2, algebra}
follows from equalities
\EqRef{a ox b c=},
\EqRef{derivative x2 differential, algebra}
and the definition
\RefDefinition{differential of map}.
The equality
\EqRef{derivative x2 component, algebra}
follows from the equality
\EqRef{derivative x2, algebra}.
}

\DefExample{differential equation y=xx, real number}
{
We consider
\AddIndex{method of successive differentiation}
{method of successive differentiation}
to solve differential equation
\ShowEq{differential equation y=xx, real number}
\DrawEq{differential equation, initial}{y=xx, real number}
over real field.
Differentiating one after another equation
\EqRef{differential equation y=xx, real number},
we get the chain of equations
\ShowEq{differential equation y=xx, 2...n, real number}
The expansion into Taylor series
\DrawEq{y=x3+C}{}
follows from equations
\ShowEq{differential equation y=xx, ref, real number}
}

\DefTheorem{int=x3 algebra}
{
Let $A$ be Banach algebra over commutative ring $D$.
\DrawEq{int=x3}{integral}
where $C$ is any $A$\Hyph number.
}

\DefProof{int=x3 algebra}
{
According to the definition
\RefDefinition{indefinite integral},
the map $y$ is integral
\eqRef{int=x3}{integral},
when the map $y$ satisfies to differential equation
\DrawEq{differential equation y=xx, 1, algebra}{integral}
and initial condition
\DrawEq{differential equation, initial}{y=xx, algebra}
We use the method of successive differentiation
to solve the differential equation
\eqRef{differential equation y=xx, 1, algebra}{integral}.
Successively differentiating equation
\eqRef{differential equation y=xx, 1, algebra}{integral},
we get the chain of equations
\ShowEq{differential equation y=xx, 2, algebra}
\ShowEq{differential equation y=xx, 3, algebra}
\ShowEq{differential equation y=xx, 4, algebra}
The expansion into Taylor series
\DrawEq{y=x3+C}{algebra}
follows from equations
\ShowEq{differential equation y=xx, ref, algebra}
The equality
\eqRef{int=x3}{integral}
follows from
\eqRef{differential equation y=xx, 1, algebra}{integral},
\eqRef{differential equation, initial}{y=xx, algebra},
\eqRef{y=x3+C}{algebra}.
}

\DefRemark{int=x3 algebra}
{
According to the definition
\EqRef{a ox b c=},
we can present integral
\eqRef{int=x3}{integral}
following way
\ShowEq{int=x3 1}
}

\DefRemark{notation int=x3 algebra}
{
In the proof of the theorem
\RefTheorem{int=x3 algebra},
I use notation like following
\ShowEq{notation theorem int=x3 algebra}
I will write following equalities
to show how derivative works. 
\ShowEq{differential equation y=xx, algebra, 1}
}

\DefRemark{differential equation 3x^2 does not possess a solution}
{
Differential equation
\DrawEq{differential equation y=xx, 1, a, algebra}{integral}
\DrawEq{differential equation, initial}{}
also leads to answer $y=x^3$. It is evident that
map $y=x^3$ does not satisfies differential equation
\eqRef{differential equation y=xx, 1, a, algebra}{integral}.
This means that differential equation
\eqRef{differential equation y=xx, 1, a, algebra}{integral}
does not possess a solution.

I advise you to pay attention that
second derivative is
not symmetric polynomial (see Tailor expansion).
}

\DefTheorem{derivative of the sum}
{
Let $A$ be Banach $D$\Hyph module.
Let $B$ be Banach $D$\Hyph algebra.
Let $f$, $g$ be differentiable maps
\ShowEq{f,g:A->B}
The map
\ShowEq{derivative of the sum, 2}
is differentiable and
the derivative satisfies to relationship
\ShowEq{derivative of the sum, 3}
}

\AddEq [4]{definition: differentiable map}
{
\begin{definition}
\labelDefinition{differentiable map, #4}
{\it
Let $#2$ be Banach $#1$\Hyph module with norm $\|a\|_#2$.
Let $#3$ be Banach $#1$\Hyph module with norm $\|a\|_#3$.
The map
\ShowEq{f:A->B}f{#2}{#3}
is called
\AddIndex{differentiable}{differentiable map}
on the set $U\subset #2$,
if at every point $x\in U$
the increment of map $f$ can be represented as
\ShowEq{derivative of map}
\DrawEq{derivative of map, def}{#4}
where
\ShowEq{df:A->B}{#2}{#3}
is linear map of $#1$\Hyph module $#2$ into $#1$\Hyph module $#3$ and
\ShowEq{f:A->B}o{#2}{#3}
is such continuous map that
\ShowEq{lim |o|/|a|}{#2}{#3}
Linear map
\ShowEq{show derivative of map}
is called
\AddIndex{derivative of map}{derivative of map}
$f$.
}
\qed
\end{definition}
}

\AddEq [4]{remark: derivative as differential form}
{
\begin{remark}
\labelRemark{derivative as differential form, #4}
According to definition
\RefDefinition{differentiable map, #4}
for given $x$, the derivative
\ShowEq{df/dx in L(A->B)}{#1}{#2}{#3}
Therefore, the derivative of map $f$ is
differential $#3$\Hyph valued $1$\Hyph form
\ShowEq{df/dx:A->L(A->B)}{#1}{#2}{#3}
\qed
\end{remark}
}

\AddEq [4]{remark: differential of map}
{
\begin{remark}
\labelRemark{differential of map, #4}
In calculus, the map
\ShowEq{differential of independent variable}
\ShowEq{dx=delta}{#1}{#2}
is called
\AddIndex{differential of independent variable}
{differential of independent variable}.
Differential $#3$\Hyph valued $1$\Hyph form
\ShowEq{df/dx:A->L(A->B)}{#1}{#2}{#3}
which is writen as
\ShowEq{differential of independent variable}
\ShowEq{differential of map}
\DrawEq{differential of map =}{#4}
is called
\AddIndex{differential of map}{differential of map}
$f$.
However, the equality
\eqRef{differential of map =}{#4}
has a different interpretation.
Namely, we consider the increment $df$ of the map $f$
as function
\eqRef{differential of map =}{#4}
of the increment $dx$ of independent variable.
\qed
\end{remark}
}

\DefRemark{differential L(A,A)}
{
According to definition
\RefDefinition{differentiable map, algebra},
the derivative of the map $f$ is the map
\EqParm{x->df/dx in L(A,B)}{A=AB}
Expressions $d_x f(x)$ and
\ShowEq{df/dx}
are different notations for the same map.
}

\DefEq
{
\begin{theorem}
\labelTheorem{representation of derivative, algebra A->B}
Let $A$ be free Banach $D$\Hyph module.
Let $B$ be free Banach $D$\Hyph algebra.
Let $\Basis F$ be the basis of left \BoxB{B}module
\ShowEq{L(A;B)}DAB.
It is possible to represent the derivative of the map
\ShowEq{f:A->B}fAB
as
\ShowEq{coordinates of derivative}
\ShowEq{df/dx=dkf/dx Ik}
\end{theorem}
}
{theorem: representation of derivative, algebra A->B}

\DefEq
{
\begin{definition}
\labelDefinition{coordinates of derivative, algebra A->B}
{\it
The expression
\ShowEq{show coordinates of derivative}
is called
\AddIndex{coordinates}{coordinates}
of derivative
\ShowEq{df/dx}
with respect to the basis $\Basis F$.
Expression
\ShowEq{component of derivative A->B}
\ShowEq{component of derivative A->B =}
is called
\AddIndex{component of derivative}
{component of derivative} of map $f(x)$.
}
\qed
\end{definition}
}
{definition: coordinates of derivative, algebra A->B}

\DefEq
{
\begin{theorem}
\labelTheorem{representation of differential, algebra A->B}
Let $A$ be free Banach $D$\Hyph module.
Let $B$ be free Banach $D$\Hyph algebra.
Let $\Basis F$ be the basis of left \BoxB{B}module
\ShowEq{L(A;B)}DAB.
It is possible to represent the differential of the map
\ShowEq{f:A->B}fAB
as
\ShowEq{df=sum dsf/dx dx A->B}
\end{theorem}
}
{theorem: representation of differential, algebra A->B}

\DefTheorem{dfa1p/dx=df/dx a1p}
{
Let $B$ be Banach module over commutative ring $D$.
Let $U$ be open set of Banach $D$\Hyph module $A$.
Let
\ShowEq{U->L(Ap,B)}
be differentiable map.
Then
\ShowEq{dfa1p/dx=df/dx a1p}
}

\AddEq{proof: dfa1p/dx=df/dx a1p}
{
\begin{proof}
Let
\ShowEq{a0 in U}
According to the definition
\RefDefinition{differentiable map, algebra},
the increment of map
$f$ can be represented as
\ShowEq{derivative of map, U->L(B)}
where
\ShowEq{derivative in L(L)}
and
\ShowEq{o:A->LB}
is such continuous map that
\ShowEq{o:A->LB, lim o}
Since
\ShowEq{f(x) in LB}
then, for any
\ShowEq{a1n}p,
the equality
\ShowEq{derivative of map, U->B}
follows from the equality
\EqRef{derivative of map, U->L(B)}.
According to the theorem
\RefTheorem{|on|->0 ona1p->0},
the equality
\ShowEq{o:A->LB, lim o a1p}
follows from the equality
\EqRef{o:A->LB, lim o}.
Considering the expression
\ShowEq{f(x) o a1p}
as map
\ShowEq{x->f(x) o a1p}
we get that the increment of this map
can be represented as
\ShowEq{derivative of map, U->B 1}
where
\ShowEq{derivative in L}
and
\ShowEq{f:A->B}{o_1}AB
is such continuous map that
\ShowEq{o1:A->B, lim o}
The equality
\ShowEq{derivative of map, U->B = +o}
follows from the equalities
\EqRef{derivative of map, U->B},
\EqRef{derivative of map, U->B 1}.
The equality
\ShowEq{lim o a1p-o1}
follows from the equalities
\EqRef{o:A->LB, lim o a1p},
\EqRef{o1:A->B, lim o}.
The equality
\EqRef{dfa1p/dx=df/dx a1p}
follows from the equalities
\EqRef{derivative of map, U->B = +o},
\EqRef{lim o a1p-o1}.
\end{proof}
}

\DefEq
{
\begin{theorem}
\labelTheorem{derivative, representation in algebra}
Definitions of the derivative
\eqRef{derivative of map, def}{algebra}
is equivalent to the definition
\ShowEq{derivative, linear path}
\end{theorem}
}
{theorem: derivative, representation in algebra}

\DefDefinition{derivative of Second Order, algebra}
{
Polylinear map
\ShowEq{derivative of Second Order}
\ShowEq{derivative of Second Order, algebra}
is called
\AddIndex{derivative of second order}
{derivative of Second Order} of map $f$.
}

\DefTheorem{Differential of Second Order, algebra, representation}
{
Let $A$ be free Banach $D$\Hyph module.
Let $B$ be free associative Banach $D$\Hyph algebra.
Let $\Basis F$ be the basis of left \BoxB{B}module
\ShowEq{L(A;B)}DAB.
It is possible to represent
the derivative of second order of the map $f$
as
\ShowEq{Differential of Second Order, algebra, representation}
Expression
\ShowEq{component of derivative of Second Order, algebra}
is called
\AddIndex{component of derivative of second order}
{component of derivative of Second Order} of map $f(x)$.
}

\DefDefinition{derivative of Order n, algebra}
{
By induction, assuming that we defined the derivative
\ShowEq{partial n-1 f}
of order $n-1$, we define
\ShowEq{derivative of Order n}
\ShowEq{derivative of Order n, algebra}
\AddIndex{derivative of order $n$}
{derivative of Order n} of map $f$.

We also assume
\ShowEq{d0f(x)}
}

\DefTheorem{n derivative equal 0, algebra}
{
If function $f(x)$ holds
\ShowEq{n derivatives of function, algebra}
then for $t\rightarrow 0$ expression $f(x+th)$ is infinitesimal of order
higher than $n$ with respect to $t$
\[
f(x_0+th)=o(t^n)
\]
}

\AddEq{Taylor polynomial}
{
Let us form polynomial
\ShowEq{Taylor polynomial, f(x), algebra}
According to theorem \RefTheorem{n derivative equal 0, algebra}
\[
f(x_0+t(x-x_0))-p(x_0+t(x-x_0))=o(t^n)
\]
Therefore, polynomial $p(x)$ is good approximation of map $f(x)$.
}

\AddEq{Taylor series}
{
If the map $f(x)$ has the derivative of any order,
then passing to the limit
$n\rightarrow\infty$, we get expansion into series
\ShowEq{Taylor series, f(x), algebra}
which is called \AddIndex{Taylor series}{Taylor series, algebra}.
}

\DefDefinition{indefinite integral}
{
Let $A$ be Banach $D$\Hyph module.
Let $B$ be Banach $D$\Hyph algebra.
The map
\ShowEq{g:A->BoxB}
is called
\AddIndex{integrable}{integrable map},
if there exists a map
\ShowEq{f:A->B}fAB
such that
\ShowEq{df=g}
Then we use notation
\ShowEq{indefinite integral}
\ShowEq{f=int g}
and the map $f$ is called
\AddIndex{indefinite integral}{indefinite integral}
of the map $g$.
}

\DefTheorem{derivative of tensor product}
{
Let $A$ be Banach $D$\Hyph module.
Let $B$, $C$ be Banach $D$\Hyph algebras.
Let $f$, $g$ be differentiable maps
\ShowEq{f,g:A->BC}
The derivative satisfies to relationship
\ShowEq{differential of tensor product, algebra}
\ShowEq{derivative of tensor product, algebra}
}

\DefTheorem{map is continuous, derivative}
{
Let $A$ be Banach $D$\Hyph module with norm $\|a\|_A$.
Let $B$ be Banach $D$\Hyph algebra with norm $\|b\|_B$.
If the derivative
\ShowEq{df/dx}
of the map
\ShowEq{f:A->B}fAB
exists in point $x$ and has finite norm,
then map $f$ is continuous at point $x$.
}

\DefTheorem{composite map, derivative, D algebra}
{
Let $A$ be Banach $D$\Hyph module with norm $\|a\|_A$.
Let $B$ be Banach $D$\Hyph module with norm $\|b\|_B$.
Let $C$ be Banach $D$\Hyph module with norm $\|c\|_C$.
Let map
\ShowEq{f:A->B}fAB
be differentiable at point
$x$
and norm of the derivative of map $f$
be finite
\ShowEq{composite map, norm f, D algebra}
Let map
\ShowEq{f:A->B}gBC
be differentiable at point
\ShowEq{composite map, y fx, D algebra}
and norm of the derivative of map $g$
be finite
\ShowEq{composite map, norm g, D algebra}
The map\,\footnote{
The notation
\ShowEq{dg/df}
means expression
\ShowEq{dg/df=}
Similar remark is true for components of derivative.
}
\ShowEq{composite map, gfx, D algebra}
is differentiable at point
\labelItem{composite map is differentiable}
$x$
\ShowEq{composite map, derivative, D algebra}
}

%% file: Stmt.Derivative.Eq.tex

\DefEq
{
\[
h:B_1\times B_2\rightarrow B
\]
}
{h:BxB->B}

\AddEquation{|f(a)|/|a|<max}
{
\frac{\|f(a_1,...,a_n)\|_B}{\|a_1\|_A...\|a_n\|_A}
\le\max
\frac{\|f(a_1,...,a_n)\|_B}{\|a_1\|_A...\|a_n\|_A}
=\|f\|
}

\AddEquation{|oa1p|<}
{
0\le\|o_n(a_1,...,a_p)\|_B\le\|o_n\|\|a_1\|_A...\|a_p\|_A
}

\DefEq
{
\[
\begin{matrix}
f:A\rightarrow B&g:A\rightarrow B
\end{matrix}
\]
}
{f,g:A->B}

\AddEq{derivative of the sum, 2}
{
\[
f+g:A\rightarrow B
\]
}

\AddEquation{derivative of the sum, 3}
{
\frac{d(f+g)}{dx}
=\frac{df}{dx}
+\frac{dg}{dx}
}

\DefEq
{
\[
\begin{matrix}
f:A\rightarrow B_1&g:A\rightarrow B_2
\end{matrix}
\]
}
{f,g:A->B12}

\DefEq
{
\[
f:A^n\rightarrow B
\]
}
{f:An->B}

\DefEq
{
\symb{\|f\|}{norm of map}{}
}
{norm of polymap}

\DefEquation
{
\ShowSymbol{norm of map}{}=
\text{sup}\frac{\|f(a_1,...,a_n)\|_B}{\|a_1\|_A...\|a_n\|_A}
}
{norm of map, module}

\DefEq
{
\(B\subset A\), \(\mu(A)=0\)
}
{B in A, mu(A)=0}

\DefEq
{
\(X\in\mathcal C_{\mu}\)
}
{X in Cmu}

\DefEq
{
\(\mu(X)\ge 0\)
\labelItem{muX>0}
}
{muX>0}

\DefEq
{
X=\bigcup_iX_i
\ \ \ i\ne j=>X_i\cap X_j=\emptyset
}
{X=+Xi}

\DefEq
{
\(X_n\in\mathcal C_{\mu}\).
}
{Xn in Cmu}

\AddEquation{composite map, y fx, D algebra}
{
y=f(x)
}

\AddEquation{composite map, norm f, D algebra}
{
\left\|\frac{df(x)}{dx}\right\|=F\le\infty
}

\AddEquation{composite map, norm g, D algebra}
{
\left\|\frac{dg(y)}{dy}\right\|=G\le\infty
}

\AddEq{composite map, gfx, D algebra}
{
\[
(g\circ f)(x)=g(f(x))
\]
}

\AddEquation{composite map, derivative, D algebra}
{
\left\{
\begin{aligned}
\frac{d(g\circ f)(x)}{dx}
&=\frac{dg(f(x))}{df(x)}\circ\frac{df(x)}{dx}
\\
\frac{d(g\circ f)(x)}{dx}\circ a
&=\frac{dg(f(x))}{df(x)}\circ\frac{df(x)}{dx}\circ a
\end{aligned}
\right.
}

\AddEq{dg/df}
{
$\displaystyle\frac{dg(f(x))}{df(x)}$
}

\AddEq{g:A->BoxB}
{
\[
g:A\rightarrow B\otimes B
\]
}

\AddEquation{differential equation y=xx, 2...n, real number}
{
\left\{
\begin{aligned}
y''&=6x
\\
y'''&=6
\\
y^{(n)}&=0\ \ \ n>3
\end{aligned}
\right.
}

\AddEq{int=x3}
{
\int(1\otimes x^2+x\otimes x+x^2\otimes 1)\circ dx=x^3+C
}

\AddEq{differential equation y=xx, 1, algebra}
{
\frac {dy}{dx}=1\otimes x^2+x\otimes x+x^2\otimes 1
}

\AddEq{differential equation y=xx, ref, real number}
{
\EqRef{differential equation y=xx, real number},
\eqRef{differential equation, initial}{y=xx, real number},
\EqRef{differential equation y=xx, 2...n, real number}.
}

\AddEq{y=x3+C}
{
y=x^3+C
}

\AddEq{differential equation y=xx, ref, algebra}
{
\eqRef{differential equation y=xx, 1, algebra}{integral},
\eqRef{differential equation, initial}{y=xx, algebra},
\EqRef{differential equation y=xx, 2, algebra},
\EqRef{differential equation y=xx, 3, algebra},
\EqRef{differential equation y=xx, 4, algebra}.
}

\AddEquation{int=x3 1}
{
\int dx\, x^2+x\,dx\, x+x^2\, dx=x^3+C
}

\AddEq{differential equation y=xx, 1, a, algebra}
{
\frac{dy}{dx}=3\otimes x^2
}

\AddEq{differential equation y=xx, algebra, 1}
{
\begin{align*}
\frac{dy}{dx}\circ h&=hx^2+xhx+x^2h
\\
\frac{d^2y}{dx^2}\circ(h_1;h_2)
&=h_1h_2x+h_1xh_2+h_2h_1x
\\
&+xh_1h_2+h_2xh_1+xh_2h_1
\\
\frac{d^3y}{dx^3}\circ(h_1;h_2;h_3)
&=h_1h_2h_3+h_1h_3h_2+h_2h_1h_3
\\
&+h_3h_1h_2+h_2h_3h_1+h_3h_2h_1
\end{align*}
}

\AddEq{notation theorem int=x3 algebra}
{
\[
(a_1\otimes_1 b_1\otimes_2 c_1+a_2\otimes_2 b_2\otimes_1 c_2)\circ(x_1,x_2)=
a_1x_1 b_1x_2 c_1+a_2x_2 b_2x_1 c_2
\]
}

\AddEq{differential equation, initial}
{
\begin{matrix}
x_0=0&y_0=C
\end{matrix}
}

\AddEquation{differential equation y=xx, 2, algebra}
{
\begin{split}
\frac{d^2y}{dx^2}
&=1\otimes_1 1\otimes_2x+1\otimes_1x\otimes_21
+1\otimes_21\otimes_1x
\\
&+x\otimes_11\otimes_21+1\otimes_2x\otimes_11
+x\otimes_21\otimes_11
\end{split}
}

\AddEquation{differential equation y=xx, 3, algebra}
{
\begin{split}
\frac{d^3y}{dx^3}
&=1\otimes_11\otimes_21\otimes_31+1\otimes_11\otimes_31\otimes_21
+1\otimes_21\otimes_11\otimes_31
\\
&+1\otimes_31\otimes_11\otimes_21+1\otimes_21\otimes_31\otimes_11
+1\otimes_31\otimes_21\otimes_11
\end{split}
}

\AddEquation{differential equation y=xx, 4, algebra}
{
\frac{d^ny}{dx^n}=0\ \ \ n>3
}

\AddEquation{differential equation y=xx, real number}
{
y'=3x^2
}

\AddEq{df=g}
{
\[
\frac{d f(x)}{d x}=g(x)
\]
}

\AddEq{f=int g}
{
\[
f(x)=
\ShowSymbol{indefinite integral}{}
\]
}

\AddEq{indefinite integral}
{
\symb{\int g(x)\circ dx}{indefinite integral}{}
}

\DefEquation
{
\|f(a_1,...,a_n)\|_B\le\|f\|\|a_1\|_A...\|a_n\|_A
}
{|f(a)|<|f||a| 1n}

\DefEq
{
\[
o_n:A\rightarrow B
\]
}
{on:A->B}

\DefEquation
{
\lim_{n\rightarrow \infty}\|o_n\|=0
}
{|on|->0}

\DefEquation
{
\lim_{n\rightarrow \infty}o_n(a_1,...,a_p)=0
}
{ona1p->0}

\DefEq
{
\[
\begin{matrix}
f:A\rightarrow B&g:A\rightarrow C
\end{matrix}
\]
}
{f,g:A->BC}

\AddEq{differential of tensor product, algebra}
{
\[
\frac{d f(x)\otimes g(x)}{d x}\circ a
=\left(\frac{d f(x)}{d x}\circ a\right)\otimes g(x)
+f(x)\otimes \left(\frac{d g(x)}{d x}\circ a\right)
\]
}

\DefEq
{
\symb{\frac{d f(x)}{d x}}{derivative of map}{}
\symb{d_x f(x)}{derivative of map inline}{}
}
{derivative of map}

\AddEq{derivative of map, def}
{
f(x+h)-f(x)
=\ShowSymbol{derivative of map inline}{}\circ h
+o(h)
=\ShowSymbol{derivative of map}{}\circ h
+o(h)
}

\AddEq [3]{df/dx:A->L(A->B)}
{
\[
\frac{df}{dx}:x\in #2\rightarrow
\frac{df(x)}{dx}\in\mathcal L(#1;#2\rightarrow #3)
\]
}

\AddEq [2]{df:A->B}
{
\[
\ShowSymbol{derivative of map}{}:#1\rightarrow #2
\]
}

\DefEq
{
\[
x\in \pA\rightarrow \frac{d f(x)}{d x}\in\mathcal L(D;\pA\rightarrow \pB)
\]
}
{x->df/dx in L(A,B)}

\DefEquation
{
\begin{split}
df=\frac{d f(x)}{d x}\circ dx
&=
\left(
\frac{d^k_{s_k\cdot 0} f(x)}{d x}
\otimes
\frac{d^k_{s_k\cdot 1} f(x)}{d x}
\right)
\circ F_k\circ dx
\\
&=\frac{d^k_{s_k\cdot 0} f(x)}{d x}
(F_k\circ dx)
\frac{d^k_{s_k\cdot 1} f(x)}{d x}
\end{split}
}
{df=sum dsf/dx dx A->B}

\DefEquation
{
\frac{df(x)}{d x}
=
\ShowSymbol{coordinates of derivative}{}
\circ F_k
}
{df/dx=dkf/dx Ik}

\DefEq
{
\[
\ShowSymbol{coordinates of derivative}{}=
\frac{d^k_{s_k\cdot 0} f(x)}{d x}
\otimes
\frac{d^k_{s_k\cdot 1} f(x)}{d x}
\in B\otimes B
\]
}
{show coordinates of derivative}

\DefEq
{
\symb{\frac{d_{s\cdot p} f(x)}{d x}}{component of derivative}{}
}
{component of derivative}

\DefEq
{
\symb{\frac{d^k f(x)}{d x}}{coordinates of derivative}{}
}
{coordinates of derivative}

\DefEq
{
\symb{\frac{d^k_{s_k\cdot p} f(x)}{d x}}{component of derivative}{A->B}
}
{component of derivative A->B}

\DefEquation
{
\frac{d f(x)}{d x}\circ a=
\lim_{t\rightarrow 0,\ t\in R}(t^{-1}(f(x+ta)-f(x)))
}
{derivative, linear path}

\DefEq
{
$\displaystyle\ShowSymbol{component of derivative}{A->B}$,
$p=0$, $1$,
}
{component of derivative A->B =}

\DefEq
{
$\displaystyle\ShowSymbol{component of derivative}{}$,
$p=0$, $1$,
}
{component of derivative =}

\AddEq [2]{lim |o|/|a|}
{
\[
\lim_{a\rightarrow 0}\frac{\|o(a)\|_{#2}}{\|a\|_{#1}}=0
\]
}

\DefEquation
{
\frac{d f(x)\circ (a_1,...,a_p)}{d x}\circ a_0=
\left(
\frac{d f(x)}{d x}\circ a_0
\right)\circ (a_1,...,a_p)
}
{dfa1p/dx=df/dx a1p}

\DefEquation
{
\begin{split}
df=\frac{d f(x)}{d x}\circ dx
&=
\left(
\frac{d_{s\cdot 0} f(x)}{d x}
\otimes
\frac{d_{s\cdot 1} f(x)}{d x}
\right)
\circ dx
\\
&=\frac{d_{s\cdot 0} f(x)}{d x}
dx
\frac{d_{s\cdot 1} f(x)}{d x}
\end{split}
}
{df=sum dsf/dx dx}

\DefEquation
{
f(x+a_0)-f(x)
=\frac{d f(x)}{d x}\circ a_0
+o(a_0)
}
{derivative of map, U->L(B)}

\AddEq{derivative in L(L)}
{
\[
\frac{d f(x)}{d x}\in\mathcal L(D;A\rightarrow \mathcal L(D;A^p\rightarrow B))
\]
}

\AddEq [3]{df/dx in L(A->B)}
{
\[
\frac{d f(x)}{d x}\in\mathcal L(#1;#2\rightarrow #3)
\]
}

\AddEq{o:A->LB}
{
\[o:A\rightarrow\mathcal L(D;A^p\rightarrow B)\]
}

\DefEquation
{
\lim_{a_0\rightarrow 0}\frac{\|o(a_0)\|}{\|a_0\|_A}=0
}
{o:A->LB, lim o}

\AddEq{f(x) in LB}
{
$f(x)\in\mathcal L(D;A^p\rightarrow B)$,
}

\DefEquation
{
\begin{split}
&\,f(x+a_0)\circ(a_1,...,a_p)
-f(x)\circ(a_1,...,a_p)
\\=&\,\left(\frac{d f(x)}{d x}\circ a_0\right)\circ(a_1,...,a_p)
+o(a_0)\circ(a_1,...,a_p)
\end{split}
}
{derivative of map, U->B}

\DefEquation
{
\lim_{a_0\rightarrow 0}\frac{\|o(a_0)\circ(a_1,...,a_p)\|_B}{\|a_0\|_A}=0
}
{o:A->LB, lim o a1p}

\AddEq{f(x) o a1p}
{
\[f(x)\circ(a_1,...,a_p)\]
}

\AddEq{x->f(x) o a1p}
{
\[
x\in U\rightarrow f(x)\circ(a_1,...,a_p)\in B
\]
}

\DefEquation
{
\begin{split}
&\,f(x+a_0)\circ(a_1,...,a_p)
-f(x)\circ(a_1,...,a_p)
\\=&\,\left(\frac{d f(x)\circ(a_1,...,a_p)}{d x}\right)\circ a_0
+o_1(a_0)
\end{split}
}
{derivative of map, U->B 1}

\AddEq{derivative in L}
{
\[
\frac{d f(x)\circ(a_1,...,a_p)}{d x}\in\mathcal L(D;A\rightarrow B)
\]
}

\DefEquation
{
\lim_{a_0\rightarrow 0}\frac{\|o_1(a_0)\|_B}{\|a_0\|_A}=0
}
{o1:A->B, lim o}

\DefEquation
{
\begin{split}
&\,\left(\frac{d f(x)\circ(a_1,...,a_p)}{d x}\right)\circ a_0
\\=&\,\left(\frac{d f(x)}{d x}\circ a_0\right)\circ(a_1,...,a_p)
+o(a_0)\circ(a_1,...,a_p)
-o_1(a_0)
\end{split}
}
{derivative of map, U->B = +o}

\DefEquation
{
\lim_{a_0\rightarrow 0}\frac{\|o(a_0)\circ(a_1,...,a_p)\|_B-\|o_1(a_0)\|_B}{\|a_0\|_A}=0
}
{lim o a1p-o1}

\AddEq [3]{o:A->B}
{
$#1:#2\rightarrow #3$
}

\AddEq{differential of independent variable}
{
\symb{dx}{differential of independent variable}{}
}

\AddEq [2]{dx=delta}
{
\[\ShowSymbol{differential of independent variable}{}
=\delta\in\mathcal L(#1;#2\rightarrow #2)\]
}

\AddEq{differential of map}
{
\symb{df}{differential of map}{}
}

\AddEq{differential of map =}
{
\ShowSymbol{differential of map}{}
=\ShowSymbol{derivative of map}{}
\circ\ShowSymbol{differential of independent variable}{}
}

\DefEq
{
$\displaystyle\frac{d f(x)}{d x}$
}
{df/dx}

\AddEq{show derivative of map}
{
$\displaystyle\ShowSymbol{derivative of map}{}$
}

\AddEq{derivative of tensor product, algebra}
{
\[
\frac{d f(x)\otimes g(x)}{d x}
=\frac{d f(x)}{d x}\otimes g(x)
+f(x)\otimes \frac{d g(x)}{d x}
\]
}

\AddEquation{bilinear map and derivative, 3}
{
\frac{dh(f(x),g(x))}{dx}\circ a
=h\left(\frac{d f(x)}{dx}\circ dx,g(x)\right)
+h\left(f(x),\frac{d g(x)}{dx}\circ dx\right)
}

\AddEquation{bilinear map and derivative, 5}
{
\frac{d\,h(f(x),g(x))}{dx}
=h\left( \frac{df(x)}{dx},g(x)\right)
+h\left(f(x), \frac{dg(x)}{dx}\right)
}

\AddEq{partial n-1 f}
{
$\displaystyle\frac{d^{n-1}f(x)}{dx^{n-1}}$
}

\AddEq{derivative of Order n}
{
\symb{\frac{d^n f(x)}{d x^n}}{derivative of Order n}{}
\symb{d^n_{x^n} f(x)}{derivative of Order n inline}{}
}

\AddEquation{derivative of Order n, algebra}
{
\begin{split}
&\,\ShowSymbol{derivative of Order n}{}
\circ(a_1;...;a_n)
=\ShowSymbol{derivative of Order n inline}{}\circ(a_1;...;a_n)
\\
=&\,\frac d{dx}
\left(\frac{d^{n-1}f(x)}{dx^{n-1}}
\circ(a_1;...;a_{n-1})\right)\circ a_n
\end{split}
}

\DefEquation
{
\frac{d x^2}{dx}=x\otimes 1+1\otimes x
}
{derivative x2, algebra}

\DefEquation
{
dx^2=x\,dx+dx\,x
}
{derivative x2 differential, algebra}

\DefEquation
{
\left\{
\begin{array}{cc}
\displaystyle\VirtFrac
\frac{d\pC{1}{0} x^2}{dx}=x
&
\displaystyle\frac{d\pC{1}{1} x^2}{dx}=1
\\
\displaystyle\VirtFrac
\frac{d\pC{2}{0} x^2}{d x}=1
&
\displaystyle\frac{d\pC{2}{1} x^2}{dx}=x
\end{array}
\right.
}
{derivative x2 component, algebra}

\DefEquation
{
f(x+h)-f(x)
=(x+h)^2-x^2
=xh+hx+h^2
=xh+hx+o(h)
}
{derivative x2, algebra, 1}

\DefEq
{
\begin{align*}
d x^2\circ dx&=2x\ dx
\\
\frac{dx^2}{dx}&=2x
\end{align*}
}
{derivative x2, field}

\AddEq{d0f(x)}
{
$\displaystyle\frac{d^0 f(x)}{dx^0}=f(x)$.
}

\AddEq{derivative of Second Order}
{
\symb{\frac{d^2 f(x)}{d x^2}}{derivative of Second Order}{}
\symb{d^2_{x^2} f(x)}{derivative of Second Order inline}{}
}

\AddEq{n derivatives of function, algebra}
{
\[
f(x_0)=\frac{df(x_0)}{dx}\circ h=...
=\frac{d^n f(x_0)}{dx^n}\circ h^n=0
\]
}

\AddEq{Taylor polynomial, f(x), algebra}
{
\[
p(x)=f(x_0)+\frac{1}{1!}\frac{df(x_0)}{dx}\circ(x-x_0)
+...+\frac{1}{n!}\frac{d^nf(x_0)}{dx^n}\circ(x-x_0)^n
\]
}

\AddEq{Taylor series, f(x), algebra}
{
\[
f(x)=\sum_{n=0}^{\infty}(n!)^{-1}\frac{d^n f(x_0)}{dx^n}\circ(x-x_0)^n
\]
}

\def\DfTwo{\frac d{dx}\left(\frac{d f(x)}{dx}\circ a_1\right)\circ a_2}

\AddEquation{derivative of Second Order, algebra}
{
\ShowSymbol{derivative of Second Order}{}\circ(a_1;a_2)
=\ShowSymbol{derivative of Second Order inline}{}\circ(a_1;a_2)
=\DfTwo
}

\AddEq{component of derivative of Second Order, algebra}
{
\symb{\frac{d\pC{s}{p}^2 f(x)}{d x^2}}
{component of derivative of Second Order}{}
\[
\begin{matrix}
\displaystyle
\ShowSymbol{component of derivative of Second Order}{}
&p=0, 1, 2
\end{matrix}
\]
}

\AddEq{Differential of Second Order, algebra, representation}
{
\[
\begin{split}
\frac{d^2 f(x)}{dx^2}\circ(a_1; a_2)
&=
\left(
\frac{d\pC{s}{0}^2 f(x)}{dx^2}
\otimes
\frac{d\pC{s}{1}^2 f(x)}{dx^2}
\otimes
\frac{d\pC{s}{2}^2 f(x)}{dx^2}
\right)\circ(F_{1\cdot s},F_{2\cdot s})\circ\sigma_s
\circ(a_1; a_2)
\\
&\VirtFrac
=\frac{d\pC{s}{0}^2 f(x)}{d x^2}
(F_{1\cdot s}\circ\sigma_s(a_1))
\frac{d\pC{s}{1}^2 f(x)}{d x^2}
(F_{2\cdot s}\circ\sigma_s(a_2))
\frac{d\pC{s}{2}^2 f(x)}{dx^2}
\end{split}
\]
}

\AddEq{h(f,g):A->B}
{
\[
h(f,g):A\rightarrow B
\]
}

\AddEq{differential of product, algebra}
{
\[
\frac{df(x)g(x)}{dx}\circ dx
=\left(\frac{df(x)}{dx}\circ dx\right)\ g(x)
+f(x)\ \left(\frac{dg(x)}{dx}\circ dx\right)
\]
}

\AddEquation{derivative of product, algebra}
{
\frac{df(x)g(x)}{dx}
=\frac{df(x)}{dx}\ g(x)+f(x)\ \frac{dg(x)}{dx}
}

%% file: Preliminary.Omega.English.tex
\input{Preliminary.Omega.Eq}

\ePrints{4975-6381,1601.03259,1610.309618526,CACAA.06}
\ifx\Semafor\ValueOn
\def\Preliminary{on}
\ShowEq{SetupRefMeasure}
\Chapter{Lebesgue Integral in Abelian \texorpdfstring{$\Omega$}{Omega}-Group}
\fi

\ePrints{5148-4632,MVector,5410-9916,4975-6381,1601.03259,1610.309618526,CACAA.06}
\ifx\Semafor\ValueOn
\Section{\texorpdfstring{$\Omega$}{Omega}-Group}

\begin{definition}
\labelDefinition{additive map}
Let sum be defined
in $\Omega$\Hyph algebra $A$.
A map
\ShowEq{f:A->B}fAA
of $\Omega_1$\Hyph algebra $A$
is called
\AddIndex{additive map}{additive map}
if
\ShowEq{f(a+b)=}
\qed
\end{definition}

\ShowDefinition{polyadditive map}

\ShowDefinition{Omega group}

\ShowTheorem{operation is distributive over addition}
\begin{proof}
The theorem follows from definitions
\RefDefinition{polyadditive map},
\RefDefinition{Omega group}.
\end{proof}

\ShowEq{contents: norm on Omega group}

\ePrints{5410-9916}
\ifx\Semafor\ValueOff
\ShowDefinition{multiplicative Omega group}
\fi

\ifx\texFuture\Defined
\begin{definition}
\labelDefinition{multiplicative map}
{\it
Let $A$, $B$ be
multiplicative $\Omega$\Hyph groups.
The map
\ShowEq{f:A->B}fAB
is called
\AddIndex{multiplicative}{multiplicative map},
if
\DrawEq{f(ab)=f(a)f(b)}{}
}
\qed
\end{definition}
\fi

\ePrints{5148-4632,MVector}
\ifx\Semafor\ValueOn
\ShowDefinition{Abelian multiplicative Omega group}

\ShowDefinition{Omega ring}

\ShowTheorem{product in Omega ring is distributive over addition}
\ShowProof{product in Omega ring is distributive over addition}

\ShowDefinition{matrix over Omega ring}

\ShowEq{convention: Einstein summation}

\ShowDefinition{matrix operations}

\ShowEq{remark: pronunciation of product}

\ShowTheorem{rcstar transpose}
\ProofTheorem{\RefRepresentation}{rcstar transpose}

\ShowDefinition{biring}

\ShowTheorem{duality principle for biring}

\ShowTheorem{duality principle for biring of matrices}

\ShowEq{remark: reducible biring}
\fi
\fi

\ePrints{aaa}
\ifx\Semafor\ValueOn
\begin{definition}
\labelDefinition{orbit of left-side representation}
{\it
Let the map
\ShowEq{f:A->*B}g{A_1}{A_2}
be the left\Hyph side representation
of multiplicative $\Omega$\Hyph group $A_1$
in $\Omega_2$\Hyph algebra $A_2$.
For any
\EqParm{a in A}{=c,n=2}
we define
\AddIndex{orbit of representation}{orbit of representation}
of the multiplicative $\Omega$\Hyph group $A_1$ as set
\ShowEq{orbit of representation}
}
\qed
\end{definition}
\fi

\ePrints{5148-4632,MVector}
\ifx\Semafor\ValueOn
\Section{Representation of Multiplicative \texorpdfstring{$\Omega$}{Omega}-Group}

We assume that transformations of representation
\ShowEq{f:A->*B}g{A_1}{A_2}
of multiplicative $\Omega$\Hyph group $A_1$
in $\Omega_2$\Hyph algebra $A_2$
may act on $A_2$\Hyph numbers either on the left or on the right.
We must consider an endomorphism of $\Omega_2$\Hyph algebra $A_2$ as operator.
In such case, the following conditions must be satisfied.

\ShowDefinition{representation of group}

\ShowTheorem{left right shift}
\ProofTheorem{\RefRepresentation}{left right shift}

\EqParm{remark: notation for effective representation}{=left}

\EqParm{remark: notation for effective representation}{=right}

\ShowTheorem{duality principle, algebra representation}
\fi

\ePrints{MVector,5148-4632}
\ifx\Semafor\ValueOn

\Section{Tensor Product of Representations}

\ShowDefinition{reduced polymorphism of representations}

\ShowDefinition{tensor product of representations}

\ShowTheorem{tensor product of representations}
\ProofTheorem{\RefRepresentation}{tensor product of representations}

\ShowTheorem{representation, tensor product}
\ProofTheorem{\RefRepresentation}{representation, tensor product}

\ShowTheorem{tensor product and polymorphism}
\ProofTheorem{\RefRepresentation}{tensor product and polymorphism}

\ShowTheorem{B times->B otimes}
\ProofTheorem{\RefRepresentation}{B times->B otimes}
\fi

\ePrints{4975-6381,1601.03259,1610.309618526,CACAA.06}
\ifx\Semafor\ValueOn

\Section{Algebra of Sets}

\ShowEq{contents: Algebra of Sets}

\Section{Lebesgue Integral}

\ShowEq{contents: Lebesgue Integral}
\fi

%% file: Preliminary.Omega.Eq.tex

\ePrints{4975-6381,1601.03259}
\ifx\Semafor\ValueOn
\def\Preliminary{on}
\ShowEq{SetupRefOmegaNorm}
\ShowEq{SetupRefMeasure}
\fi

\DefEq
{
\ePrints{5410-9916,4975-6381,1610.309618526,CACAA.06,1601.03259}
\ifx\Semafor\ValueOn
\ShowDefinition{norm on Omega group}

\ShowTheorem{|a-b|>|a|-|b|}
\ePrints{CACAA.06}%
\ifx\Semafor\ValueOff%
\ProofTheorem{\RefTheoremOmegaNorm}{|a-b|>|a|-|b|}
\fi

\ShowEq{=Omega}

\ePrints{4975-6381,1601.03259}
\ifx\Semafor\ValueOff
\ShowDefinition{open ball}

\ShowDefinition{closed ball}
\fi

\ShowDefinition{open set}

\ePrints{5410-9916}
\ifx\Semafor\ValueOn
\ShowDefinition{continuous map, Omega group}

\ShowTheorem{continuous map, open set}
\ProofTheorem{\RefTheoremOmegaNorm}{continuous map, open set}

\ShowTheorem{image of interval is interval, real field}
\ProofTheorem{\RefTheoremOmegaNorm}{image of interval is interval, real field}

\ShowDefinition{norm of operation}

\ShowTheorem{|fx|<|f||x|1n}
\ProofTheorem{\RefTheoremOmegaNorm}{|fx|<|f||x|1n}

\ShowDefinition{norm of representation}

\ShowTheorem{|fab|<|f||a||b|}
\ProofTheorem{\RefTheoremOmegaNorm}{|fab|<|f||a||b|}

\ShowDefinition{compact set}

\ShowTheorem{norm of compact set is bounded}
\ProofTheorem{\RefTheoremOmegaNorm}{norm of compact set is bounded}
\fi

\ShowDefinition{limit of sequence}

\ShowDefinition{fundamental sequence}

\ePrints{5410-9916}
\ifx\Semafor\ValueOn
\ShowTheorem{lim a=lim b}
\ProofTheorem{\RefTheoremOmegaNorm}{lim a=lim b}
\fi

\ShowDefinition{complete Omega group}

\ePrints{5410-9916}
\ifx\Semafor\ValueOn
\ShowTheorem{c1+c2 in B}
\ProofTheorem{\RefTheoremOmegaNorm}{c1+c2 in B}

\ShowTheorem{set of maps to Omega group}
\ProofTheorem{\RefTheoremOmegaNorm}{set of maps to Omega group}

\ShowEq{remark: set of maps to Omega group}

\ShowDefinition{limit of sequence, map to Omega group}

\ShowDefinition{sequence converges uniformly}

\ShowTheorem{sequence converges uniformly, fn-fm}
\ProofTheorem{\RefTheoremOmegaNorm}{sequence converges uniformly, fn-fm}

\ShowTheorem{h=f+g, converges uniformly}
\ProofTheorem{\RefTheoremOmegaNorm}{h=f+g, converges uniformly}

\ShowTheorem{h=f1...fn omega, converges uniformly}
\ProofTheorem{\RefTheoremOmegaNorm}{h=f1...fn omega, converges uniformly}

\ShowTheorem{representation f generates representation fX}
\ProofTheorem{\RefTheoremOmegaNorm}{representation f generates representation fX}

\ShowTheorem{fX(g1)(g2), converges uniformly}
\ProofTheorem{\RefTheoremOmegaNorm}{fX(g1)(g2), converges uniformly}
\fi

\ePrints{5410-9916}
\ifx\Semafor\ValueOff
\ShowDefinition{series converges normally}
\fi

\fi
}
{contents: norm on Omega group}

\DefEq
{
\ShowDefinition{semiring of sets}

\ShowDefinition{algebra of sets}

\ShowEq{remark: A cup/minus B}

\ShowDefinition{sigma algebra of sets}

\ShowDefinition{Borel algebra}
}
{contents: Algebra of Sets}

\AddEq{contents: Lebesgue Integral}
{
\ShowDefinition{C_X C_Y measurable}

\ShowEq{example: measurable map}



\ShowDefinition{sigma-additive measure}

\ShowDefinition{simple map}

\ShowTheorem{measurable simple map}
\ProofTheorem{\RefMeasure}{measurable simple map}

\ShowEq{remark - measure - effective representation of real field}

\ShowDefinition{integrable map}

\ShowDefinition{Integral of Map over Set of Finite Measure}

\ShowTheorem{|int f|<int|f|}
\ProofTheorem{\RefMeasure}{|int f|<int|f|}

\ShowTheorem{int |f|<M mu}
\ProofTheorem{\RefMeasure}{int |f|<M mu}

\ShowTheorem{int f+g X}
\ProofTheorem{\RefMeasure}{int f+g X}

\ShowTheorem{|int h|<int|omega||f1n|}
\ProofTheorem{\RefMeasure}{|int h|<int|omega||f1n|}

\ShowTheorem{h=fX g1 g2 representation}
\ProofTheorem{\RefMeasure}{h=fX g1 g2 representation}
}

%% file: Derivative.16.English.tex

\input{Derivative.16.Eq}

\Chapter{Differentiable Maps}
\labelChapter{Differentiable maps}

\Section{Topological Ring}
\labelSection{Topological Ring}

\begin{definition}
Ring $D$ is called 
\AddIndex{topological ring}{topological ring}\,\footnote{
I made definition according to definition
from \citeBib{Pontryagin: Topological Group},
chapter 4}
if $D$ is topological space and the algebraic operations
defined in $D$ are continuous in the topological space $D$.
\qed
\end{definition}

According to definition, for arbitrary elements $a$, $b\in D$
and for arbitrary neighborhoods $W_{a-b}$ of the element $a-b$,
$W_{ab}$ of the element $ab$ there exists neighborhoods $W_a$
of the element $a$ and $W_b$ of the element $b$ such
that $W_a-W_b\subset W_{a-b}$, $W_aW_b\subset W_{ab}$.

\begin{definition}
\labelDefinition{norm on ring}
\AddIndex{Norm on ring}{norm on ring}
$D$ is a map\,\footnote{
I made definition according to the definition from
\citeBib{Bourbaki: General Topology: Chapter 5 - 10},
IX, \S 3.2
and the definition
\citeBib{Arnautov Glavatsky Mikhalev}-1.1.12,
p. 23.}
\[d\in D\rightarrow |d|\in R\]
which satisfies the following axioms
\begin{itemize}
\item $|a|\ge 0$
\item $|a|=0$ if, and only if, $a=0$
\item $|ab|=|a|\ |b|$
\item $|a+b|\le|a|+|b|$
\end{itemize}

Ring $D$, endowed with the structure defined by a given norm on
$D$, is called
\AddIndex{normed ring}{normed ring}.
\qed
\end{definition}

Invariant distance on additive group of ring $D$
\[d(a,b)=|a-b|\]
defines topology of metric space,
compatible with ring structure of $D$.

\begin{definition}
\labelDefinition{limit of sequence, normed ring}
Let $D$ be normed ring.
Element $a\in D$ is called 
\AddIndex{limit of a sequence}{limit of sequence}
\ShowEq{limit of sequence, normed ring}
if for every $\epsilon\in R$, $\epsilon>0$,
there exists positive integer $n_0$ depending on $\epsilon$ and such,
that
\ShowEq{an-a}
for every $n>n_0$.
\qed
\end{definition}

\begin{theorem}
\labelTheorem{limit of sequence, normed ring, product on scalar}
Let $D$ be normed ring of characteristic $0$ and let $d\in D$.
Let $a\in D$ be limit of a sequence $\{a_n\}$.
Then
\[
\lim_{n\rightarrow\infty}(a_nd)=ad
\]
\[
\lim_{n\rightarrow\infty}(da_n)=da
\]
\end{theorem}
\begin{proof}
Statement of the theorem is trivial, however I give this proof
for completeness sake. 
Since $a\in D$ is limit of the sequence $\{a_n\}$,
then according to definition
\RefDefinition{limit of sequence, normed ring}
for given $\epsilon\in R$, $\epsilon>0$,
there exists positive integer $n_0$ such, that
\ShowEq{an-a d}
for every $n>n_0$.
According to definition \RefDefinition{norm on ring}
the statement of theorem follows from inequalities
\ShowEq{an-a d1}
for any $n>n_0$.
\end{proof}

\begin{definition}
Let $D$ be normed ring.
The sequence $\{a_n\}$, $a_n\in D$ is called 
\AddIndex{fundamental}{fundamental sequence}
or \AddIndex{Cauchy sequence}{Cauchy sequence},
if for every $\epsilon\in R$, $\epsilon>0$
there exists positive integer $n_0$ depending on $\epsilon$ and such,
that $|a_p-a_q|<\epsilon$ for every $p$, $q>n_0$.
\qed
\end{definition}

\begin{definition}
Normed ring $D$ is called
\AddIndex{complete}{complete ring}
if any fundamental sequence of elements
of ring $D$ converges, i.e.
has limit in ring $D$.
\qed
\end{definition}

Later on, speaking about normed ring of characteristic $0$,
we will assume that homeomorphism of field of rational numbers $Q$
into ring $D$ is defined.

\begin{theorem}
\labelTheorem{complete ring contains real number}
Complete ring $D$ of characteristic $0$
contains as subfield an isomorphic image of the field $R$ of
real numbers.
It is customary to identify it with $R$.
\end{theorem}
\begin{proof}
Consider fundamental sequence of rational numbers $\{p_n\}$.
Let $p'$ be limit of this sequence in ring $D$.
Let $p$ be limit of this sequence in field $R$.
Since immersion of field $Q$ into division ring $D$ is homeomorphism,
then we may identify $p'\in D$ and $p\in R$.
\end{proof}

\begin{theorem}
\labelTheorem{complete ring and real number}
Let $D$ be complete ring of characteristic $0$ and let $d\in D$.
Then any real number $p\in R$ commute with $d$.
\end{theorem}
\begin{proof}
Let us represent real number $p\in R$ as
fundamental sequence of rational numbers $\{p_n\}$.
Statement of theorem follows from chain of equalities
\ShowEq{complete ring and real number}
based on statement of theorem
\RefTheorem{limit of sequence, normed ring, product on scalar}.
\end{proof}

\Section{Topological \texorpdfstring{$D$}{D}-Algebra}

\begin{definition}
\labelDefinition{norm on D module}
Let $D$ be normed commutative ring.\,\footnote{I
made definition according to definition
from \citeBib{Bourbaki: General Topology: Chapter 5 - 10},
IX, \S 3.3.
We use notation either $|a|$ or $\|a\|$ for norm.
}
\AddIndex{Norm on $D$\Hyph module}
{norm on D module} $A$
is a map
\ShowEq{norm on D module}
which satisfies the following axioms
\StartLabelItem
\begin{enumerate}
\ShowEq{norm on D module 1}
\ShowEq{norm on D module 2, 1}
if, and only if,
\ShowEq{norm on D module 2, 2}
\ShowEq{norm on D module 3}
\end{enumerate}
$D$\Hyph module $A$,
endowed with the structure defined by a given norm on
$A$, is called
\AddIndex{normed $D$\Hyph module}{normed D module}.
\qed
\end{definition}

\ShowEq{=module}

\ShowDefinition{open ball}

\ShowDefinition{closed ball}

\begin{theorem}
\labelTheorem{b in Bar}
\ShowEq{b in Bar}
iff
\ShowEq{b-a in Br}
\end{theorem}
\begin{proof}
The theorem follows from the definition
\RefDefinition{closed ball}.
\end{proof}



\ShowDefinition{limit of sequence}

\ShowTheorem{limit of sequence}
\ShowProof{limit of sequence}

\ShowDefinition{fundamental sequence}

\ShowTheorem{fundamental sequence}
\ShowProof{fundamental sequence}

\begin{definition}
Norms\,\footnote{
See also the definition
\citeBib{Shilov single 3}\Hyph 12.35.{\cyr a}
on page 53.
}
\ShowEq{|1,|2}{}
defined on $D$\Hyph module $A$
are called
\AddIndex{equivalent}{equivalent norms}
if the statement
\DrawEq{a=lim an}{}
does not depend on selected norm.
\qed
\end{definition}


In theorems
\RefTheorem{B1o in B2o},
we consider norms
\ShowEq{|1,|2}{}
defined on $D$\Hyph module $A$.
We also consider open and closed balls
\ShowEq{B1o,B2o}

\begin{theorem}
\labelTheorem{B1o in B2o}
Norms\,\footnote{
See also the lemma
\citeBib{Shilov single 3}\Hyph 12.35.{\cyr b}
on page 54.
}
\ShowEq{|1,|2}{}
defined on $D$\Hyph module $A$
are equivalent iff
there exist constants
\ShowEq{c12>0}
such that for any $\rho>0$ 
\ShowEq{cB1c in B2c}
\ShowEq{cB2c in B1c}
\end{theorem}
\begin{proof}
Let there exist constants $c_1$, $c_2$.
Let
\DrawEq{a=lim an}{}
with respect to norm $\|x\|_1$.
According to theorems
\RefTheorem{b in Bar},
\RefTheorem{limit of sequence, normed module}
and definitions
\RefDefinition{open ball},
\RefDefinition{closed ball},
for any
\ShowEq{epsilon in R}
there exists positive integer $n_0$ such, that
\ShowEq{a-an in B1e}
The statement
\ShowEq{a-an in B2}
follows from statements
\EqRef{cB1c in B2c},
\EqRef{a-an in B1e}
and the definition
\RefDefinition{closed ball}.
According to the theorem
\RefTheorem{limit of sequence, normed module},
\DrawEq{a=lim an}{}
with respect to norm $\|x\|_2$.

The similar way we prove that the statement
\DrawEq{a=lim an}{}
with respect to norm $\|x\|_1$
follows from the statement
\DrawEq{a=lim an}{}
with respect to norm $\|x\|_2$.

Let norms
\ShowEq{|1,|2}{}
be equivalent.
Let there be no constant $c_2$ such that the statement
\EqRef{cB2c in B1c}
is true.
Then, for any
\ShowEq{n,1.}
we can find closed balls
\ShowEq{B1rn B2rn}
such that
\ShowEq{cB2c not in B1c}
From the statement
\EqRef{cB2c not in B1c}
and from the definition
\RefDefinition{closed ball},
it follows that there exist $x_n$ such that
\ShowEq{xn rn <>}
Let
\ShowEq{yn=xn/n}
From the equality
\EqRef{yn=xn/n}
and from statements
\RefItem{norm on D module 3, 2},
\EqRef{xn rn <>},
it follows that
\ShowEq{yn rn <>}
From the statement
\EqRef{yn rn <>},
it follows that the statement
\ShowEq{lim y =0}
is true with respect to norm $\|x\|_2$ and
is false with respect to norm $\|x\|_1$.
This contradicts to the statement that norms are equivalent;
therefore, the constatnt $c_2$ exists.
The similar way we prove that the constatnt $c_1$ exists.
\end{proof}

\begin{theorem}
\labelTheorem{Norms are equivalent iff |1<|2<|1}
Norms\,\footnote{
See also the corollary
\citeBib{Shilov single 3}\Hyph 12.35.{\cyr v}
on page 54.
See also the definition
\citeBib{8176-4374}\Hyph 9.9
and the proposition
\citeBib{8176-4374}\Hyph 9.10
on the page 288.
}
\ShowEq{|1,|2}{}
defined on $D$\Hyph module $A$
are equivalent iff
there exist constants
\ShowEq{c12>0}
such that for any $x\in A$
\ShowEq{|1<|2<|1}
\end{theorem}
\begin{proof}
According to the theorem
\RefTheorem{B1o in B2o},
there exists the constant
$c_1$
such that
for
\ShowEq{r=x|1}
\ShowEq{B1x in B2x}
The statement
\ShowEq{|2<|1}
follows from the statement
\EqRef{B1x in B2x}.
According to the theorem
\RefTheorem{B1o in B2o},
there exists the constant
$c_2$
such that
for
\ShowEq{r=x|2}
\ShowEq{B2x in B1x}
The statement
\ShowEq{|1<|2}
follows from the statement
\EqRef{B2x in B1x}.
The statement
follows from statements
\EqRef{|2<|1},
\EqRef{|1<|2}.
\end{proof}

\begin{definition}
Normed $D$\Hyph module $A$ is called
\AddIndex{Banach $D$\Hyph module}{Banach module}
if any fundamental sequence of elements
of module $A$ converges, i.e.
has limit in module $A$.
\qed
\end{definition}

\ifx\texFuture\Defined
\begin{definition}
Given a topological commutative ring $D$ and
$D$\Hyph algebra $A$ such that $A$ has
a topology compatible with the structure of the additive
group of $A$ and maps
\ShowEq{topological D algebra}
are continuous,
then $A$ is called a
\AddIndex{topological $D$\Hyph algebra}
{topological D algebra}\,\footnote{
I made definition according to definition
from \citeBib{Bourbaki: Topological Vector Space},
p. TVS I.1}.
\qed
\end{definition}
\fi

\begin{definition}
\labelDefinition{continuous map, algebra}
A map
\ShowEq{f:A->B}f{A_1}{A_2}
of Banach $D_1$\Hyph module $A_1$ with norm $|x|_1$
into Banach $D_2$\Hyph module $A_2$ with norm $|y|_2$
is called \AddIndex{continuous}{continuous map}, if
for every as small as we please $\epsilon>0$
there exist such $\delta>0$, that
\[
|x'-x|_1<\delta
\]
implies
\[
|f(x')-f(x)|_2<\epsilon
\]
\qed
\end{definition}

\ShowTheorem{set of A->B is D module}
\begin{proof}
We define sum of maps
\ShowEq{fg in BA}
according to the equality
\ShowEq{(f+g)(a)=}
According to the definition
\RefDefinition{module over commutative ring},
the set $B^A$ is Abelian group.
We define the representation of the ring $D$
in Abelian group $B^A$
according to the equality
\ShowEq{(df)(a)=}
We can right the equality
\EqRef{(df)(a)=}
without using brackets.
Therefore, the statement
\RefItem{set of A->B is D module}
is true.

Statements
\RefItem{norm on D module 1},
\RefItem{norm on D module 2}
are evident for the value $\|f\|$.
The statement
\RefItem{norm on D module 3, 1}
follows from the inequality
\ShowEq{|(f+g)(a)|<}
and from the definition
\EqRef{norm of map, algebra}.
The statement
\RefItem{norm on D module 3, 2}
for the value $\|f\|$
follows from the equality
\EqRef{(df)(a)=},
from the statement
\RefItem{norm on D module 3, 2}
for the value $\|a\|_B$
and from the definition
\EqRef{norm of map, algebra}.
Therefore, the statement
\RefItem{norm of map}
is true.
\end{proof}

\begin{theorem}
Let
\ShowEq{f:A->B}f{A_1}{A_2}
be linear map of Banach $D$\Hyph module $A_1$ with norm $|x|_1$
into Banach $D$\Hyph module $A_2$ with norm $|y|_2$.
Then
\ShowEq{norm of linear map, algebra}
\end{theorem}
\begin{proof}
From definitions
\EqParm{ref definition linear map from A1 to A2}{=z}
and theorems
\RefTheorem{complete ring contains real number},
\RefTheorem{complete ring and real number},
it follows that
\ShowEq{linear map and real number}
From the equality \EqRef{linear map and real number}
and the definition \RefDefinition{norm on d algebra}
it follows that
\ShowEq{linear map and real number, 1}
Assuming $\displaystyle r=\frac 1{|x|_1}$, we get
\ShowEq{norm of linear map, algebra, 1}
Equality \EqRef{norm of linear map, algebra}
follows from equalities \EqRef{norm of linear map, algebra, 1}
and \EqRef{norm of map, algebra}.
\end{proof}

\begin{theorem}
Let
\ShowEq{f:A->B}f{A_1}{A_2}
be linear map of Banach $D_1$\Hyph module $A_1$ with norm $|x|_1$
into Banach $D_2$\Hyph module $A_2$ with norm $|y|_2$.
Since $\|f\|<\infty$, then map $f$ is continuous.
\end{theorem}
\begin{proof}
Since map $f$ is linear, then
according to the theorem \RefTheorem{set of A->B is D module}
\[
|f(x)-f(y)|_2=|f(x-y)|_2\le \|f\|\ |x-y|_1
\]
Let us assume arbitrary $\epsilon>0$. Let
\ShowEq{delta=epsilon/|f|}
Then
\[
|f(x)-f(y)|_2\le \|f\|\ \delta=\epsilon
\]
follows from inequality
\[
|x-y|_1<\delta
\]
According to definition \RefDefinition{continuous map, algebra}
map $f$ is continuous.
\end{proof}

\begin{theorem}
\labelTheorem{L(A;B) is Banach module}
Let $A$ be Banach $D$\Hyph module with norm $|x|_A$.
Let $B$ be Banach $D$\Hyph module with norm $|y|_B$.
$D$\Hyph module
\ShowEq{L(A;B)}DAB{}
equiped by norm
\ShowEq{norm of map to algebra}
is Banach $D$\Hyph module.
\end{theorem}
\begin{proof}
To prove the theorem, we need to prove
that limit of a sequence of linear maps
is a linear map.
The statement follows from equalities
\ShowEq{lim f(a+b)=}
\ShowEq{lim f(ca)=}
\end{proof}

\ShowDefinition{norm of polylinear map}

\ShowTheorem{|f(a)|<|f||a| 1n}
\ShowProof{|f(a)|<|f||a| 1n}

\ShowTheorem{|on|->0 ona1p->0}
\ShowProof{|on|->0 ona1p->0}

According to the definition
\RefDefinition{algebra over ring},
$D$\Hyph algebra $A$ is
$D$\Hyph module $A$
equiped by bilinear map
\ShowEq{*:AA->A}
which is called product.
In general, norm $\|*\|$ of product may be different from $1$.

\ShowEq{def H2}
\begin{example}
\labelExample{algebra of hyperbolic numbers}
{\it
Let\,\footnote{
I appreciate the remark of Nathan BeDell
where he put my attention
that there exists algebra
where norm $\|*\|$ of product may be different from $1$.
The example
\RefExample{algebra of hyperbolic numbers}
is based on the definition
\citeBib{1708.01190}\Hyph 3.1
on page 3
from the paper
\citeBib{1708.01190}
by Nathan BeDell.
See also the definition on the page
\citeBib{978-0-8176-8384-9}\Hyph 24.
}
\HH\ be algebra of hyperbolic numbers.
\HH\Hyph number has form
\ShowEq{H2 number}{}{}
where $a$, $b$ are real numbers.
Product of \HH\Hyph numbers
\ShowEq{H2 number}1,
\ShowEq{H2 number}2{}
has the following form
\ShowEq{product hyperbolic numbers}
From the equality
\EqRef{product hyperbolic numbers},
it follows that $j^2=1$.
Consider the norm\,\footnote{
Such norm was introduced in
\citeBib{978-0-8176-8384-9}
on page 25.
}
of \HH\Hyph number
\ShowEq{H2 number}{}{}
defined by the equality
\ShowEq{H2 number norm}

From equalities
\EqRef{product hyperbolic numbers},
\EqRef{H2 number norm},
it follows that norm of product has value
\ShowEq{norm hyperbolic numbers product}
Therefore, $\|*\|=1$ in algebra \HH.
}
\qed
\end{example}

\begin{example}
{\it
The norm $\|x\|$ of \HH\Hyph number
considered in example
\RefExample{algebra of hyperbolic numbers}
does not satisfy to the statement
\RefItem{norm on D module 2}.\,\footnote{
Garret Sobczyk put our attention that this norm
corresponds to geometry of event space in special relativity.
}
So we consider different norm.
If the norm of \HH\Hyph number
\ShowEq{H2 number}{}{}
is defined by the equality
\ShowEq{H2 number norm 1}
then norm of product is maximum of the map
\ShowEq{norm of product is maximum of map}
The map $f$ is symmetric with respect to its variables; so it gets its extremum
when variables are equal.
Let
\ShowEq{x1...=}
Then
\ShowEq{faaaa=}
This value is maximum because
\ShowEq{f1010=}
Therefore,
\ShowEq{|*|=v2}
}
\qed
\end{example}

According to the theorem
\RefTheorem{Norms are equivalent iff |1<|2<|1},
if we defined norm $\|x\|_1$ in algebra \HH,
then any norm $\|x\|_2$ defined by the equality
\ShowEq{|2=r|1}
is equivalent to the norm $\|x\|_1$.

\begin{theorem}
\labelTheorem{there exists equivalent norm |*|=1}
Let $A$ be $D$\Hyph algebra.
If, in $D$\Hyph module $A$, there exist norm $\|x\|_1$ such that
norm $\|*\|_1$ of product in $D$\Hyph algebra $A$ is different from $1$,
then there exists equivalent norm
\ShowEq{|2=*|1 |1}
in $D$\Hyph module $A$ such that
\ShowEq{|*|2=1}
\end{theorem}
\begin{proof}
The equality
\ShowEq{|*|2=1 1}
follows from equalities
\EqRef{norm of map, module},
\EqRef{|2=*|1 |1}.
The equality
\EqRef{|*|2=1}
follows from the equality
\EqRef{|*|2=1 1}.
\end{proof}

According to the theorem
\RefTheorem{there exists equivalent norm |*|=1},
if we introduce norm
\ShowEq{|2=v2|1}
in algebra \HH, then norm of product $\|*\|_2=1$.

\begin{definition}
\labelDefinition{norm on d algebra}
Let $D$ be normed commutative ring.
Let $A$ be $D$\Hyph algebra.
The norm\,\footnote{
I made definition according to definition
from \citeBib{Bourbaki: General Topology: Chapter 5 - 10},
IX, \S 3.3.
If $D$\Hyph algebra $A$ is division algebra,
then norm is called
\AddIndex{absolute value}{absolute value}
and we use notation $|a|$ for norm of $A$\Hyph number $a$.
See the definition
from \citeBib{Bourbaki: General Topology: Chapter 5 - 10},
IX, \S 3.2.
}
\ShowEq{norm on d algebra}
on $D$\Hyph module $A$ such that\,\footnote{
The inequality
\EqRef{norm on d algebra 3}
follows from the theorem
\RefTheorem{there exists equivalent norm |*|=1}.
Otherwise we would have to write
\ShowEq{norm on d algebra 3 1}
}
\ShowEq{norm on d algebra 3}
is called
\AddIndex{norm on $D$\Hyph algebra}
{norm on D algebra} $A$.
$D$\Hyph algebra $A$,
endowed with the structure defined by a given norm on
$A$, is called
\AddIndex{normed $D$\Hyph algebra}{normed D algebra}.
\qed
\end{definition}

\begin{definition}
Normed $D$\Hyph algebra $A$ is called
\AddIndex{Banach $D$\Hyph algebra}{Banach algebra}
if any fundamental sequence of elements
of algebra $A$ converges, i.e.
has limit in algebra $A$.
\qed
\end{definition}

\Section{The Derivative of Map of \texorpdfstring{$D$}{D}-Algebra}
\labelSection{The Derivative of Map in Algebra}

\ShowDefinition{differentiable map}DAB{algebra}

\begin{definition}
\labelDefinition{differential of map}
{\it
Since, for given $x\in A$, we consider the increment
\eqRef{derivative of map, def}{algebra}
of the map
\ShowEq{f:A->B}fAB
as function of differential $dx$ of variable $x$,
then the linear part of this function
\ShowEq{differential of independent variable}
\ShowEq{differential of map}
\DrawEq{differential of map =}{}
is called
\AddIndex{differential of map}{differential of map}
$f$.
}
\qed
\end{definition}

\ShowEq{remark: differential L(A,A)}

\ShowTheorem{derivative, representation in algebra}
\begin{proof}
From definitions
\EqParm{ref definition linear map from A1 to A2}{=z}
and the theorem
\RefTheorem{complete ring contains real number}
it follows
\ShowEq{differential is multiplicative over field R, algebra}
\ShowEq{t ne 0 in R, a ne 0 in A}
Combining equality
\EqRef{differential is multiplicative over field R, algebra}
and definition
\RefDefinition{differentiable map, algebra},
we get the definition
\EqRef{derivative, linear path}
of the derivative.
\end{proof}

\begin{corollary}
{\it
A map $f$
is called differentiable on the set $U\subset D$,
if at every point $x\in U$
the increment of the map $f$ can be represented as
\ShowEq{differential of map, t, algebra}
where
$o:R\rightarrow A$ is such continuous map that
\[
\lim_{t\rightarrow 0}\frac{|o(t)|}{|t|}=0
\]
}
\qed
\end{corollary}

\ShowTheorem{representation of derivative, algebra A->B}
\begin{proof}
The theorem
follows from the definitions
\RefDefinition{differentiable map, algebra}
and from the statement
\ShowEq{ref map f generated by basis I}
\end{proof}

\ShowDefinition{coordinates of derivative, algebra A->B}

\ShowTheorem{representation of differential, algebra A->B}
\begin{proof}
The theorem
follows from the theorem
\RefTheorem{representation of derivative, algebra A->B}
and from the definitions
\RefDefinition{coordinates of derivative, algebra A->B}.
\end{proof}

\begin{remark}
\labelRemark{linear maps generated by map I0}
If $D$\Hyph module
\ShowEq{L(A;A)}
is generated by the map $F_0\circ x=x$,
then the equality
\EqRef{df=sum dsf/dx dx A->B}
gets form
\ShowEq{df=sum dsf/dx dx}
where the expression
\ShowEq{component of derivative}
\ShowEq{component of derivative =}
is called
\AddIndex{component of derivative}
{component of derivative} of map $f(x)$.
\qed
\end{remark}

\begin{theorem}
\labelTheorem{derivative, standard form, algebra}
Let $A$ be free Banach $D$\Hyph module.
Let $B$ be free Banach $D$\Hyph algebra.
Let $\Basis F$ be the basis of left \BoxB{B}module
\ShowEq{L(A;B)}DAB.
Let $\Basis e$ be the basis of $D$\Hyph module $B$.
\AddIndex{Standard representation of derivative}
{derivative, standard representation} of map
\ShowEq{f:A->B}fAB
has form
\ShowEq{standard component of derivative}
\ShowEq{derivative, algebra, standard representation}
Expression
\ShowEq{standard component of derivative, algebra}
in equality
\EqRef{derivative, algebra, standard representation}
is called
\AddIndex{standard component of derivative}
{standard component of derivative} of the map $f$.
\end{theorem}
\begin{proof}
Statement of theorem folows from the definitions
\RefDefinition{differentiable map, algebra}
and from the statement
\ShowEq{ref standard representation of map A1 A2, associative algebra}
\end{proof}

\begin{theorem}
\labelTheorem{derivative and jacobian, algebra}
Let $A$ be free Banach $D$\Hyph module.
Let $B$ be free Banach $D$\Hyph algebra.
Let $\Basis F$ be the basis of left \BoxB{B}module
\ShowEq{L(A;B)}DAB.
Let $\Basis e_A$ be basis of the free finite dimensional
$D$\Hyph module $A$.
Let $\Basis e_B$ be basis of the free finite dimensional associative
$D$\Hyph algebra $B$.
Let
\ShowEq{structural constants, algebra}
be structural constants of algebra $B$.
Then it is possible to represent the derivative of the map
\ShowEq{f:A->B}fAB
as
\ShowEq{derivative and jacobian, algebra}
where $dx\in A$ has expansion
\ShowEq{derivative and jacobian, 1, algebra}
relative to basis $\Basis e_A$ and Jacobian matrix of map $f$ has form
\ShowEq{standard components and Jacobian, algebra A->B}
\end{theorem}
\begin{proof}
Statement of theorem follows from the theorem
\ShowEq{ref theorem coordinates of map A1 A2, algebra}
\end{proof}

\begin{remark}
If $D$\Hyph module
\ShowEq{L(A;A)}
is generated by the map $F_0\circ x=x$,
then equalities
\EqRef{derivative, algebra, standard representation},
\EqRef{standard components and Jacobian, algebra A->B}
get form
\ShowEq{standard component of derivative I0}
\ShowEq{derivative, algebra, standard representation I0}
\ShowEq{standard components and Jacobian, algebra}
where the expression
\ShowEq{standard component of derivative, algebra I0}
is called
\AddIndex{standard component of derivative}
{standard component of derivative} of the map $f$.
\qed
\end{remark}

\ShowTheorem{derivative of the sum}
\begin{proof}
According to the definition
\EqRef{derivative, linear path},
\ShowEq{derivative of the sum, 4}
The equality \EqRef{derivative of the sum, 3}
follows from the equality \EqRef{derivative of the sum, 4}.
\end{proof}

\ShowTheorem{dfa1p/dx=df/dx a1p}
\ShowProof{dfa1p/dx=df/dx a1p}

\begin{convention}
\labelConvention{bilinear map, parameter linear map}
Given bilinear map
\ShowEq{h:BxB->B}
we consider following maps
\ShowEq{h1:A3->A}
\ShowEq{h2:A3->A}
defined by equality
\ShowEq{h1()=}
\ShowEq{h2()=}
We will use letter $h$ to denote maps $h_1$, $h_2$.
\qed
\end{convention}

\ShowTheorem{bilinear map and differential}
\begin{proof}
Equality \EqRef{bilinear map and derivative, 3}
follows from chain of equalities
\ShowEq{bilinear map and derivative, 4}
based on definition
\EqRef{derivative, linear path}.
Equality \EqRef{bilinear map and derivative, 5}
follows from the equality
\EqRef{bilinear map and derivative, 3}
and from the convention
\RefConvention{bilinear map, parameter linear map}.
\end{proof}

\ShowTheorem{derivative of product, algebra}
\begin{proof}
The theorem follows from theorems
\RefTheorem{bilinear map and differential}
and the definition
\ShowEq{ref algebra over ring}
\end{proof}

\begin{theorem}
Let $A$ be free Banach $D$\Hyph module.
Let $B$ be free Banach $D$\Hyph algebra.
Let the derivative of the map
\ShowEq{f:A->B}fAB
have expansion
\DrawEq{derivative of f, algebra}{product}
Let the derivative of the map
\ShowEq{g:A->B}
have expansion
\DrawEq{derivative of g, algebra}{product}
The derivative of the map $f(x)g(x)$
has form
\ShowEq{derivative, fg, algebra}
\ShowEq{component of derivative, fg, algebra}
\end{theorem}
\begin{proof}
Let us substitute \eqRef{derivative of f, algebra}{product}
and \eqRef{derivative of g, algebra}{product}
into equality \EqRef{derivative of product, algebra}
\ShowEq{derivative of fg, division ring}
Based \EqRef{derivative of fg, division ring},
we define equalities
\EqRef{component of derivative, fg, algebra}.
\end{proof}

\begin{theorem}
\labelTheorem{derivative of representation}
Let $A$ be Banach $D$\Hyph module.
Let $B$, $C$ be Banach $D$\Hyph algebras.
Let
\ShowEq{h:B->*C}
be representation of $D$\Hyph module $B$
in $D$\Hyph module $C$.
Let $f$, $g$ be differentiable maps
\ShowEq{f,g:A->BC}
The derivative of the map
\ShowEq{h(f)(g)}
has form
\ShowEq{derivative of representation}
\end{theorem}
\begin{proof}
Since the map
\ShowEq{h:B->End C}
is homomorphism of the Abelian group
and representation
\ShowEq{h(a):C->C}
is endomorphism of the additive group,
then the map $h(b)(c)$ is bilinear map.
The theorem follows from the theorem
\RefTheorem{bilinear map and differential}.
\end{proof}

\ShowTheorem{derivative of tensor product}
\begin{proof}
The theorem follows from the theorems
\RefTheorem{bilinear map and differential},
\ShowEq{ref theorem V times->V otimes}
and the definition
\ShowEq{ref algebra over ring}
\end{proof}

\ifx\texFuture\Defined
\begin{remark}
Let
\DrawEq{derivative of f, algebra}{tensor}
\DrawEq{derivative of g, algebra}{tensor}
Then
\ShowEq{derivative of tensor product 1, algebra}
We do not write brackets, because tensor product
is associative and distributive over addition
(theorems
\ShowEq{ref tensor product is associative and distributive}
).
\qed
\end{remark}
\fi

\ShowTheorem{map is continuous, derivative}
\begin{proof}
From the theorem \RefTheorem{set of A->B is D module}
it follows
\ShowEq{derivative and continuos map, algebra, 1}
From \eqRef{derivative of map, def}{algebra},
\EqRef{derivative and continuos map, algebra, 1}
it follows
\ShowEq{derivative and continuos map, algebra, 2}
Let us assume arbitrary $\epsilon>0$. Let
\ShowEq{delta=epsilon/|df|}
Then from inequality
\ShowEq{|a|<delta}
it follows
\ShowEq{derivative and continuos map, algebra, 3}
According to definition \RefDefinition{continuous map, algebra}
map $f$ is continuous at point $x$.
\end{proof}

\begin{theorem}
\labelTheorem{derivative, 0, D algebra}
Let $A$ be free Banach $D$\Hyph module.
Let $B$ be free Banach $D$\Hyph algebra.
Let map
\ShowEq{f:A->B}fAB
be differentiable at point
$x$.
Then
\ShowEq{derivative, 0, D algebra}
\end{theorem}
\begin{proof}
The theorem follows from the definitions
\RefDefinition{differentiable map, algebra}
and from the theorem
\ShowEq{ref linear map, 0}
\end{proof}

\ShowTheorem{composite map, derivative, D algebra}
\begin{proof}
According to definition
\RefDefinition{differentiable map, algebra}
\ShowEq{derivative of map g, D algebra}
where
\ShowEq{o1:B->C}
is such continuous map that
\ShowEq{o1:B->C, 1}
According to definition
\RefDefinition{differentiable map, algebra}
\ShowEq{derivative of map f, D algebra}
where
\ShowEq{o2:A->B}
is such continuous map that
\ShowEq{o2:A->B, 1}
According to
\EqRef{derivative of map f, D algebra}
increment $a$ of value $x\in A$
leads to increment
\ShowEq{composite map, b fxa, D algebra}
of value $y$.
Using \EqRef{composite map, y fx, D algebra},
\EqRef{composite map, b fxa, D algebra}
in equality \EqRef{derivative of map g, D algebra},
we get
\ShowEq{derivative of map gf, 1, D algebra}
According to definitions
\EqParm{ref definition linear map from A1 to A2}{=c}
from equality \EqRef{derivative of map gf, 1, D algebra}
it follows
\ShowEq{derivative of map gf, 2, D algebra}
According to definition \RefDefinition{norm on d algebra}
\ShowEq{derivative of map gf, 3, D algebra}
From \EqRef{composite map, norm g, D algebra}
it follows that
\ShowEq{derivative of map gf, 4, D algebra}
From \EqRef{composite map, norm f, D algebra}
it follows that
\ShowEq{derivative of map gf, 6, D algebra}
According to the theorem \RefTheorem{derivative, 0, D algebra}
\ShowEq{derivative of map gf, 7, D algebra}
Therefore,
\ShowEq{derivative of map gf, 5, D algebra}
From equalities
\EqRef{derivative of map gf, 3, D algebra},
\EqRef{derivative of map gf, 4, D algebra},
\EqRef{derivative of map gf, 5, D algebra}
it follows
\ShowEq{derivative of map gf, 8, D algebra}
According to definition
\RefDefinition{differentiable map, algebra}
\ShowEq{derivative of map gf, 9, D algebra}
where
\ShowEq{o:A->B}oAC
is such continuous map that
\ShowEq{lim |o|/|a|}AC
Equality \EqRef{composite map, derivative, D algebra}
follows from
\EqRef{derivative of map gf, 2, D algebra},
\EqRef{derivative of map gf, 8, D algebra},
\EqRef{derivative of map gf, 9, D algebra}.

From equality \EqRef{composite map, derivative, D algebra}
and theorem
\RefTheorem{representation of derivative, algebra A->B},
it follows that
\ShowEq{composite map, derivative 2, D algebra}
\EqRef{composite map, derivative 01, D algebra}
follow from equality
\EqRef{composite map, derivative 2, D algebra}.
\end{proof}

\begin{theorem}
Let $A$ be free Banach $D$\Hyph module.
Let $B$, $C$ be free associative Banach $D$\Hyph algebras.
Let $\Basis F$ be the basis of left \BoxB{B}module
\ShowEq{L(A;B)}DAB.
Let $\Basis G$ be the basis of left \BoxB{C}module
\ShowEq{L(A;B)}DBC.
Let
\ShowEq{expansion of derivative of f with respect to basis I}
be expansion of derivative of the map
\ShowEq{f:A->B}fAB
with respect to the basis $\Basis F$.
Let
\ShowEq{expansion of derivative of g with respect to basis J}
be expansion of derivative of the map
\ShowEq{f:A->B}gBC
with respect to the basis $\Basis G$.
Then derivative of map
\DrawEq{h=g o f}{}
has expansion
\ShowEq{expansion of derivative of h with respect to basis K}
with respect to the basis
\DrawEq{JlIk}{}
where
\ShowEq{expansion of derivative of h with respect to basis K =}
\end{theorem}
\begin{proof}
The theorem follows from theorems
\RefTheorem{representation of composition of linear maps},
\RefTheorem{representation of derivative, algebra A->B}
and from the definition
\RefDefinition{differentiable map, algebra}.
\end{proof}

\begin{theorem}
Let $A$ be free associative $D$\Hyph algebra.
Let left \BoxB{A}module
\ShowEq{L(A;B)}DAA{}
is generated by the identity map $F_0=\delta$.
Let
\ShowEq{expansion of derivative of f A->A}
be expansion of derivative of the map
\ShowEq{f:A->B}fAA
Let
\ShowEq{expansion of derivative of g A->A}
be expansion of derivative of the map
\ShowEq{f:A->B}gAA
Then derivative of the map
\DrawEq{h=g o f}{}
has expansion
\ShowEq{expansion of derivative of h A->A}
where
\ShowEq{composite map, derivative 01, D algebra}
\end{theorem}
\begin{proof}
The theorem follows from the theorem
\RefTheorem{representation of composition of linear maps A->A},
from the remark
\ref{remark: linear maps generated by map I0}
and from the definition
\RefDefinition{differentiable map, algebra}.
\end{proof}

%% file: Derivative.16.Eq.tex

\input{Derivative.16.Ref}

\DefEq
{
\[
pd=\lim_{n\rightarrow\infty}(p_nd)=\lim_{n\rightarrow\infty}(dp_n)=dp
\]
}
{complete ring and real number}

\DefEq
{
\[
|a_n-a|<\epsilon
\]
}
{an-a}

\AddEq{c12>0}
{
$c_1>0$, $c_2>0$
}

\AddEq{b in Bar}
{
$b\in B_c(a,\rho)$
}

\AddEq{b-a in Br}
{
$b-a\in B_c(0,\rho)$
}

\AddEquation{cB2c in B1c}
{
B_{2c}(0,c_2\rho)\subseteq B_{1c}(0,\rho)
}

\AddEquation{cB1c in B2c}
{
B_{1c}(0,c_1\rho)\subseteq B_{2c}(0,\rho)
}

\AddEquation{a-an in B1e}
{
a-a_n\in B_{1o}(0,\frac 12c_1\epsilon)\subseteq B_{1c}(0,\frac 12c_1\epsilon)
}

\AddEq{a-an in B2}
{
\[
a-a_n\in B_{2c}(0,\frac 12\epsilon)\subseteq B_{2o}(0,\epsilon)
\]
}

\AddEq{B1rn B2rn}
{
$B_{1c}(0,\rho_n)$, $B_{2c}(0,\rho_n/n)$
}

\AddEq{lim y =0}
{
\[\lim_{n\rightarrow \infty}y_n=0\]
}

\AddEquation{cB2c not in B1c}
{
B_{2c}(0,\rho_n/n)\not\subseteq B_{1c}(0,\rho_n)
}

\AddEquation{|1<|2<|1}
{
c_2\|x\|_1\le\|x\|_2\le\frac{\|x\|_1}{c_1}
}

\AddEquation{B1x in B2x}
{
B_{1c}(0,\|x\|_1)\subseteq B_{2c}(0,\|x\|_1/c_1)
}

\AddEquation{B2x in B1x}
{
B_{2c}(0,\|x\|_2)\subseteq B_{1c}(0,\|x\|_2/c_2)
}

\AddEquation{|2<|1}
{
\|x\|_2\le\|x\|_1/c_1
}

\AddEquation{|1<|2}
{
\|x\|_1\le\|x\|_2/c_2
}

\AddEq{r=x|2}
{
$\rho=\|x\|_2/c_2$
}

\AddEq{r=x|1}
{
$\rho=\|x\|_1/c_1$
}

\AddEquation{product hyperbolic numbers}
{
z_1z_2=a_1a_2+b_1b_2+j(a_1b_2+a_2b_1)
}

\AddEquation{H2 number norm}
{
\|z\|=\sqrt{|a^2-b^2|}
}

\AddEquation{H2 number norm 1}
{
\|z\|_1=\sqrt{a^2+b^2}
}

\AddEq{f1010=}
{
\[
f(1,0,1,0)=
\frac{\displaystyle\sqrt{(1*1+0*0)^2+(1*0+1*0)^2}}
{\displaystyle\sqrt{1^2+0^2}\sqrt{1^2+0^2}}=1
\]
}

\AddEq{|*|=v2}
{
$\|*\|_1=\sqrt 2$.
}

\AddEq{norm of product is maximum of map}
{
\[
f(x_1,y_1,x_2,y_2)=
\frac{\displaystyle\sqrt{(x_1x_2+y_1y_2)^2+(x_1y_2+x_2y_1)^2}}
{\displaystyle\sqrt{x_1^2+y_1^2}\sqrt{x_2^2+y_2^2}}
\]
}

\AddEq{faaaa=}
{
\[
f(a,a,a,a)=
\frac{\displaystyle\sqrt{(a^2+a^2)^2+(a^2+a^2)^2}}
{\displaystyle\sqrt{a^2+a^2}\sqrt{a^2+a^2}}
=\frac{\sqrt{8}}{2}=\sqrt{2}
\]
}

\AddEq{x1...=}
{
$x_1=x_2=y_1=y_2=a$, $a\ne 0$.
}

\AddEq{map f has maximum when}
{
$x_1=x_2=x$, $y_1=y_2=y$, $x\ne y$.
}

\AddEq{norm hyperbolic numbers product}
{
\begin{align*}
\|*\|&=
\frac{\displaystyle\sqrt{(x_1x_2+y_1y_2)^2-(x_1y_2+x_2y_1)^2}}
{\displaystyle\sqrt{x_1^2-y_1^2}\sqrt{x_2^2-y_2^2}}\\
&=
\frac{\displaystyle\sqrt{(x_1x_2+y_1y_2-x_1y_2-x_2y_1)(x_1x_2+y_1y_2+x_1y_2+x_2y_1)}}
{\displaystyle\sqrt{(x_1-y_1)(x_1+y_1)(x_2-y_2)(x_2+y_2)}}\\
&=
\frac{\displaystyle\sqrt{(x_1(x_2-y_2)-y_1(y_1-y_2))(x_1(x_2+y_2)+y_1(x_2+y_2))}}
{\displaystyle\sqrt{(x_1-y_1)(x_1+y_1)(x_2-y_2)(x_2+y_2)}}
\end{align*}
}

\AddEq{def H2}
{
\def\HH{\ensuremath{\mathcal H_2}}
}

\AddEq [2]{H2 number}
{
$z_{#1}=a_{#1}+b_{#1}j$#2
}

\AddEq{*:AA->A}
{
\[
*:A\times A\rightarrow A
\]
}

\AddEquation{xn rn <>}
{
\begin{matrix}
x_n\le \rho_n/n&x_n>\rho_n
\end{matrix}
}

\AddEquation{yn rn <>}
{
\begin{matrix}
\|y_n\|_1>1&\|y_n\|_2\le 1/n<2/n
\end{matrix}
}

\AddEquation{yn=xn/n}
{
y_n=x_n/n
}

\AddEq{B1o,B2o}
{
\begin{align*}
B_{1o}(0,\rho)&=\{a\in A:\|x\|_1<\rho\}&
B_{1c}(0,\rho)&=\{a\in A:\|x\|_1\le\rho\}\\
B_{2o}(0,\rho)&=\{a\in A:\|x\|_2<\rho\}&
B_{2c}(0,\rho)&=\{a\in A:\|x\|_2\le\rho\}
\end{align*}
}

\AddEq[1]{|1,|2}
{
$\|x\|_1$, $\|x\|_2$#1
}

\AddEq{B(a)in B}
{
$B_o(a,R)\subseteq B$.
}

\DefEq
{
\[
|a_n-a|<\frac{\epsilon}{|d|}
\]
}
{an-a d}

\DefEq
{
\begin{align*}
|a_nd-ad|=|(a_n-a)d|&=|a_n-a||d|<\frac{\epsilon}{|d|}|d|=\epsilon
\\
|da_n-da|=|d(a_n-a)|&=|d||a_n-a|<|d|\frac{\epsilon}{|d|}=\epsilon
\end{align*}
}
{an-a d1}

\DefEq
{
\[
\frac{dg(f(x))}{df(x)}=
\left.\frac{dg(y)}{dy}\right|_{y=f(x)}
\]
}
{dg/df=}

\DefEquation
{
\begin{split}
\frac{d\pC{st}{0} h(x)}{dx}&=
\frac{d\pC{s}{0} g(f(x))}{df(x)}
\frac{d\pC{t}{0} f(x)}{dx}
\\
\frac{d\pC{st}{1} h(x)}{dx}&=
\frac{d\pC{t}{1} f(x)}{dx}
\ \frac{d\pC{s}{1} g(f(x))}{df(x)}
\end{split}
}
{composite map, derivative 01, D algebra}

\DefEquation
{
\frac{df}{dx}=\frac{d^kf}{dx}\circ F_k
}
{expansion of derivative of f with respect to basis I}

\DefEquation
{
\frac{dg}{dx}=\frac{d^lg}{dx}\circ G_l
}
{expansion of derivative of g with respect to basis J}

\DefEquation
{
\frac{dh}{dx}=\frac{d^{lk}h}{dx}\circ H_{lk}
}
{expansion of derivative of h with respect to basis K}

\DefEquation
{
\frac{df}{dx}=\frac{d_{s\cdot 0}f}{dx}\otimes \frac{d_{s\cdot 1}f}{dx}
}
{expansion of derivative of f A->A}

\DefEquation
{
\frac{dg}{dx}=\frac{d_{t\cdot 0}g}{dx}\otimes \frac{d_{t\cdot 1}g}{dx}
}
{expansion of derivative of g A->A}

\DefEquation
{
\frac{dh}{dx}=\frac{d_{ts\cdot 0}h}{dx}\otimes \frac{d_{ts\cdot 1}h}{dx}
}
{expansion of derivative of h A->A}

\DefEquation
{
\frac{d^{lk}h}{dx}=\frac{d^lg}{dx}\circ \left(G^k_m\circ \frac{d^mf}{dx}\right)
}
{expansion of derivative of h with respect to basis K =}

\DefEq
{
$o_1:B\rightarrow C$
}
{o1:B->C}

\DefEq
{
$o_2:A\rightarrow B$
}
{o2:A->B}

\DefEq
{
\[g:A\rightarrow B\]
}
{g:A->B}

\DefEq
{
\[
\lim_{b\rightarrow 0}\frac{\|o_1(b)\|_C}{\|b\|_B}=0
\]
}
{o1:B->C, 1}

\DefEq
{
\[
\lim_{a\rightarrow 0}\frac{\|o_2(a)\|_B}{\|a\|_A}=0
\]
}
{o2:A->B, 1}

\DefEquation
{
b=\frac{df(x)}{dx}\circ a+o_2(a)
}
{composite map, b fxa, D algebra}

\DefEquation
{
\begin{split}
&\,g(f(x+a))-g(f(x))
\\
=&\,g\left(f(x)+\frac{df(x)}{dx}\circ a
+o_2(a)\right)
-g(f(x))
\\
=&\,
\frac{dg(f(x))}{df(x)}\circ\left(\frac{df(x)}{dx}\circ a
+o_2(a)\right)
-o_1\left(\frac{df(x)}{dx}\circ a
+o_2(a)\right)
\end{split}
}
{derivative of map gf, 1, D algebra}

\DefEquation
{
\begin{split}
g(f(x+a))-g(f(x))
&=\frac{dg(f(x))}{df(x)}\circ\frac{df(x)}{dx}\circ a
\\&+\frac{dg(f(x))}{df(x)}\circ o_2(a)
-o_1\left(\frac{df(x)}{dx}\circ a
+o_2(a)\right)
\end{split}
}
{derivative of map gf, 2, D algebra}

\DefEquation
{
\begin{split}
&
\lim_{a\rightarrow 0}
\frac
{\displaystyle\left\|\frac{dg(f(x))}{df(x)}\circ o_2(a)
-o_1\left(\frac{df(x)}{dx}\circ a
+o_2(a)\right)\right\|_C}
{\|a\|_A}
\\
\le
&
\lim_{a\rightarrow 0}
\frac
{\displaystyle\left\|\frac{dg(f(x))}{df(x)}\circ o_2(a)\right\|_C}
{\|a\|_A}
+\lim_{a\rightarrow 0}
\frac
{\displaystyle\left\|o_1\left(\frac{df(x)}{dx}\circ a
+o_2(a)\right)\right\|_C}
{\|a\|_A}
\end{split}
}
{derivative of map gf, 3, D algebra}

\DefEquation
{
\lim_{a\rightarrow 0}
\frac
{\displaystyle\left\|\frac{dg(f(x))}{df(x)}\circ o_2(a)\right\|_C}
{|a|_A}
\le
G\lim_{a\rightarrow 0}
\frac{\|o_2(a)\|_B}{\|a\|_A}=0
}
{derivative of map gf, 4, D algebra}

\DefEq
{
\begin{align*}
&\lim_{a\rightarrow 0}
\frac
{\displaystyle\left\|o_1\left(\frac{df(x)}{dx}\circ a
+o_2(a)\right)\right\|_C}
{\|a\|_A}
\\
=&
\lim_{a\rightarrow 0}
\frac
{\displaystyle\left\|o_1\left(\frac{df(x)}{dx}\circ a
+o_2(a)\right)\right\|_C}
{\displaystyle\left\|\frac{df(x)}{dx}\circ a
+o_2(a)\right\|_B}
\lim_{a\rightarrow 0}
\frac
{\displaystyle\left\|\frac{df(x)}{dx}\circ a
+o_2(a)\right\|_B}
{\|a\|_A}
\\
\le&
\lim_{a\rightarrow 0}
\frac
{\displaystyle\left\|o_1\left(\frac{df(x)}{dx}\circ a
+o_2(a)\right)\right\|_C}
{\displaystyle\left\|\frac{df(x)}{dx}\circ a
+o_2(a)\right\|_B}
\lim_{a\rightarrow 0}
\frac
{\displaystyle\left\|\frac{df(x)}{dx}\right\|\|a\|_A
+\|o_2(a)\|_B}
{\|a\|_A}
\\
=&
\lim_{a\rightarrow 0}
\frac
{\displaystyle\left\|o_1\left(\frac{df(x)}{dx}\circ a
+o_2(a)\right)\right\|_C}
{\displaystyle\left\|\frac{df(x)}{dx}\circ a
+o_2(a)\right\|_B}
\left\|\frac{df(x)}{dx}\right\|
\end{align*}
}
{derivative of map gf, 6, D algebra}

\DefEq
{
\[
\lim_{a\rightarrow 0}
\left(\frac{df(x)}{dx}\circ a
+o_2(a)\right)=0
\]
}
{derivative of map gf, 7, D algebra}

\DefEquation
{
\lim_{a\rightarrow 0}
\frac
{\displaystyle\left\|\frac{dg(f(x))}{df(x)}\circ o_2(a)
-o_1\left(\frac{df(x)}{dx}\circ a
+o_2(a)\right)\right\|_C}
{\|a\|_A}=0
}
{derivative of map gf, 8, D algebra}

\DefEq
{
\symb{\|a\|}{norm on d algebra}1
}
{norm on d algebra}

\DefEq
{
\symb{\|a\|}{norm on D module}{}
\[a\in A\rightarrow \ShowSymbol{norm on D module}{}\in R\]
}
{norm on D module}

\DefEq
{
\item $\|a\|\ge 0$
\labelItem{norm on D module 1}
}
{norm on D module 1}

\DefEq
{
\item $\|a\|=0$
\labelItem{norm on D module 2}
}
{norm on D module 2, 1}

\DefEq
{
$a=0$
}
{norm on D module 2, 2}

\AddEquation{|*|2=1}
{
\|*\|_2=1
}

\AddEquation{|2=*|1 |1}
{
\|x\|_2=\|*\|_1\|x\|_1
}

\AddEquation{|*|2=1 1}
{
\begin{split}
\|*\|_2
&=\text{sup}\frac{\|ab\|_2}{\|a\|_2\|b\|_2}
=\text{sup}\frac{\|*\|_1\|ab\|_1}{\|*\|_1\|a\|_1\|*\|_1\|b\|_1}
=\frac 1{\|*\|_1}\text{sup}\frac{\|ab\|_1}{\|a\|_1\|b\|_1}\\
&=\frac{\|*\|_1}{\|*\|_1}
\end{split}
}

\AddEq{|2=r|1}
{
\[
\|x\|_2=r\|x\|_1\ \ \ r\in R,r>0
\]
}

\DefEquation
{
(df)(a)=d(f(a))
}
{(df)(a)=}

\DefEq
{
$f$, $g\in B^A$
}
{fg in BA}

\DefEq
{
\item $\|a+b\|\le \|a\|+\|b\|$
\labelItem{norm on D module 3, 1}
\item $\|da\|=|d|\,\|a\|$, $d\in D$, $a\in A$
\labelItem{norm on D module 3, 2}
}
{norm on D module 3}

\DefEq
{
\[\|(f+g)(a)\|_B\le\|f(a)\|_B+\|g(a)\|_B\]
}
{|(f+g)(a)|<}

\DefEq
{
\[\|f\|=
\text{sup}\frac{\|f(x)\|_B}{\|x\|_A}\]
}
{norm of map to algebra}

\DefEq
{
\[
\lim_{n\rightarrow \infty}f_n(a+b)
=\lim_{n\rightarrow\infty}f_n(a)+\lim_{n\rightarrow \infty}f_n(b)
\]
}
{lim f(a+b)=}

\DefEq
{
\[
\lim_{n\rightarrow \infty}f_n(ca)=c\lim_{n\rightarrow \infty}f_n(a)
\]
}
{lim f(ca)=}

\DefEq
{
$\mathcal L(D;A\rightarrow A)$
}
{L(A;A)}

\DefEquation
{
(g\circ f)(x+a)-(g\circ f)(x)
=\frac{dg(f(x))}{dx}\circ a
+o(a)
}
{derivative of map gf, 9, D algebra}

\DefEquation
{
\begin{array}{r@{\,}l}
&\displaystyle
\frac{d\pC{st}{0} (g\circ f)(x)}{dx}
\ a
\ \frac{d\pC{st}{1} (g\circ f)(x)}{dx}
\\[10pt]
=&\displaystyle
\frac{d\pC{s}{0} g(f(x))}{df(x)}
\ \left(\frac{df(x)}{dx}\circ a\right)
\ \frac{d\pC{s}{1} g(f(x))}{df(x)}
\\[10pt]
=&\displaystyle
\frac{d\pC{s}{0} g(f(x))}{df(x)}
\frac{d\pC{t}{0} f(x)}{dx}
\ a
\ \frac{d\pC{t}{1} f(x)}{dx}
\ \frac{d\pC{s}{1} g(f(x))}{df(x)}
\end{array}
}
{composite map, derivative 2, D algebra}

\DefEquation
{
\lim_{a\rightarrow 0}
\frac
{\displaystyle\left\|o_1\left(\frac{df(x)}{dx}\circ a
+o_2(a)\right)\right\|_C}
{\|a\|_A}=0
}
{derivative of map gf, 5, D algebra}

\DefEq
{
\[
\frac{df(x)}{dx}\circ 0=0
\]
}
{derivative, 0, D algebra}

\DefEquation
{
f(x+a)-f(x)
=\frac{df(x)}{dx}\circ a
+o_2(a)
}
{derivative of map f, D algebra}

\DefEquation
{
g(y+b)-g(y)
=\frac{dg(y)}{dy}\circ b
+o_1(b)
}
{derivative of map g, D algebra}

\DefEq
{
\frac{df(x)}{dx}
=\frac{d\pC{s}{0} f(x)}{dx}
\otimes
\frac{d\pC{s}{1} f(x)}{dx}
}
{derivative of f, algebra}

\DefEq
{
\frac{dg(x)}{dx}
=\frac{d\pC{t}{0} g(x)}{dx}
\otimes
\frac{d\pC{t}{1} g(x)}{dx}
}
{derivative of g, algebra}

\DefEquation
{
\frac{df(x)g(x)}{dx}
=\frac{d\pC{s}{0} f(x)}{dx}
\otimes
\left(
\frac{d\pC{s}{1} f(x)}{dx}g(x)
\right)
+\left(
f(x)
\frac{d\pC{t}{0} g(x)}{dx}
\right)
\otimes
\frac{d\pC{t}{1} g(x)}{dx}
}
{derivative, fg, algebra}

\DefEquation
{
\begin{array}{r@{\,}l@{\ \ \ }r@{\,}l}
\displaystyle\frac{d\pC{s}{0} f(x)g(x)}{dx}
&
\displaystyle=\frac{d\pC{s}{0} f(x)}{dx}
&
\displaystyle\frac{d\pC{t}{0} f(x)g(x)}{dx}
&
\displaystyle=f(x)\frac{d\pC{t}{0} g(x)}{dx}
\\[15pt]
\displaystyle\frac{d\pC{s}{1} f(x)g(x)}{dx}
&
\displaystyle=\frac{d\pC{s}{1} f(x)}{dx}
g(x)
&
\displaystyle\frac{d\pC{t}{1} f(x)g(x)}{dx}
&
\displaystyle=\frac{d\pC{t}{1} g(x)}{dx}
\end{array}
}
{component of derivative, fg, algebra}

\DefEquation
{
\begin{split}
\frac{df(x) g(x)}{dx}\circ a
&=\left(\frac{df(x)}{dx}\circ a\right)\ g(x)
+f(x)\ \left(\frac{dg(x)}{dx}\circ a\right)
\\
&=\frac{d\pC{s}{0} f(x)}{dx}a
\frac{d\pC{s}{1} f(x)}{dx}
g(x)
+f(x)\frac{d\pC{t}{0} g(x)}{dx}a
\frac{d\pC{t}{1} g(x)}{dx}
\end{split}
}
{derivative of fg, division ring}

\DefEq
{
\begin{align*}
\frac{dh(f(x),g(x))}{dx}\circ a
&=\lim_{t\rightarrow 0}(t^{-1}(h(f(x+ta),g(x+ta))-h(f(x),g(x))))
\\
&=\lim_{t\rightarrow 0}(t^{-1}(h(f(x+ta),g(x+ta))-h(f(x),g(x+ta))))
\\
&+\lim_{t\rightarrow 0}(t^{-1}(h(f(x),g(x+ta))-h(f(x),g(x))))
\\
&=h(\lim_{t\rightarrow 0}t^{-1}(f(x+ta)-f(x)),g(x))
\\
&+h(f(x),\lim_{t\rightarrow 0}t^{-1}(g(x+ta)-g(x)))
\end{align*}
}
{bilinear map and derivative, 4}

\DefEq
{
\begin{align*}
(a,v)\in D\times A&\rightarrow av\in A
\\
(v,w)\in A\times A&\rightarrow vw\in A
\end{align*}
}
{topological D algebra}

\AddEquation{norm on d algebra 3}
{
\|ab\|\le\|a\|\,\|b\|
}

\AddEq{norm on d algebra 3 1}
{
\[
\|ab\|\le\|*\|\|a\|\,\|b\|
\]
}

\AddEq{|2=v2|1}
{
\[
\|x\|_2=\sqrt{2}\|x\|_1
\]
}

\DefEq
{
$a\in A$, $|a|=1$,
}
{unit sphere in algebra}

\DefEq
{
$\{a_n\}$
\symb{\lim_{n\rightarrow\infty}a_n}{limit of sequence}{}
\[
a=\ShowSymbol{limit of sequence}{}
\]
}
{limit of sequence, normed ring}

\DefEquation
{
\frac{df(x)}{dx}\circ(ta)
=
t\frac{df(x)}{dx}\circ a
}
{differential is multiplicative over field R, algebra}

\DefEq
{
\[
\begin{matrix}
t\in R& t\ne 0&a\in A&a\ne 0
\end{matrix}
\]
}
{t ne 0 in R, a ne 0 in A}

\DefEquation
{
f(x+ta)-f(x)
=t\frac{df(x)}{dx}\circ a
+o(t)
}
{differential of map, t, algebra}

\DefEquation
{
\begin{split}
\frac{d(f+g)(x)}{dx}\circ a
&=
\lim_{t\rightarrow 0,\ t\in R}(t^{-1}((f+g)(x+ta)-(f+g)(x)))
\\
&=
\lim_{t\rightarrow 0,\ t\in R}(t^{-1}(f(x+ta)+g(x+ta)-f(x)-g(x)))
\\
&=
\lim_{t\rightarrow 0,\ t\in R}(t^{-1}(f(x+ta)-f(x)))
\\
&+
\lim_{t\rightarrow 0,\ t\in R}(t^{-1}(g(x+ta)-g(x)))
\\
&=
\frac{df(x)}{dx}\circ a
+\frac{dg(x)}{dx}\circ a
\end{split}
}
{derivative of the sum, 4}

\DefEquation
{
\frac{df(x)\otimes g(x)}{dx}
=\frac{d\pC{s}{0} f(x)}{dx}
\otimes
\frac{d\pC{s}{1} f(x)}{dx}
\otimes g(x)
+f(x)\otimes \frac{d\pC{t}{0} g(x)}{dx}
\otimes
\frac{d\pC{t}{1} g(x)}{dx}
}
{derivative of tensor product 1, algebra}

\DefEq
{
\[
\frac{df(x)}{dx}\circ dx
=dx^{\gii}\frac{\partial f^{\gij}}{\partial x^{\gii}}
e_{\gij}
\]
}
{derivative and jacobian, algebra}

\DefEq
{
\[
h_1:\mathcal L(D;A\rightarrow B_1)\times B_2\rightarrow \mathcal L(D;A\rightarrow B)
\]
}
{h1:A3->A}

\DefEq
{
\[
h_1(f,v)\circ u=h(f\circ u,v)
\]
}
{h1()=}

\DefEq
{
\[
h_2:B_1\times\mathcal L(D;A\rightarrow B_2)\rightarrow \mathcal L(D;A\rightarrow B)
\]
}
{h2:A3->A}

\DefEq
{
\[
h_2(u,f)\circ v=h(u,f\circ v)
\]
}
{h2()=}

\DefEquation
{
\left\|\frac{df(x)}{dx}\circ a\right\|_B
\le\left\|\frac{df(x)}{dx}\right\|\,\|a\|_A
}
{derivative and continuos map, algebra, 1}

\DefEq
{
$h(f(x))(g(x))$
}
{h(f)(g)}

\DefEquation
{
\frac{dh(f(x))(g(x))}{dx}
=h\left(\frac{df(x)}{dx}\right)(g(x))+h(f(x))\left(\frac{dg(x)}{dx}\right)
}
{derivative of representation}

\DefEq
{
\[h(a):C\rightarrow C\]
}
{h(a):C->C}

\DefEq
{
\[h:B\rightarrow \End(\{+\},C)\]
}
{h:B->End C}

\DefEq
{
\[
\xymatrix
{
h:B\ar[r]|{*}&C
}
\]
}
{h:B->*C}

\DefEquation
{
\|f(x+a)-f(x)\|_B<\left\|\frac{df(x)}{dx}\right\|\,\|a\|_A
}
{derivative and continuos map, algebra, 2}

\DefEquation
{
\frac{\partial f^{\gij}}{\partial x^{\gii}}=
\StandPartial{f(x)}{x}{kr} C_{\gi{ki}}^{\gi p}C_{\gi{pr}}^{\gij}
}
{standard components and Jacobian, algebra}

\DefEquation
{
\frac{\partial f^{\gik}}{\partial x^{\gil}}=
\frac{d^{k\cdot\gi{ij}} f(x)}{dx}
F_{k\cdot}^{}{}^{\gim}_{\gil}
C^{\gi p}_{\gi{im}}C^{\gik}_{\gi{pj}}
}
{standard components and Jacobian, algebra A->B}

\DefEq
{
\[
\begin{matrix}
dx=dx^{\gii}\ e_{\gii}&&dx^{\gii}\in D
\end{matrix}
\]
}
{derivative and jacobian, 1, algebra}

\DefEq
{
$\displaystyle\ShowSymbol{standard component of derivative}{}$
}
{standard component of derivative, algebra}

\DefEq
{
$\displaystyle\ShowSymbol{standard component of derivative}{I0}$
}
{standard component of derivative, algebra I0}

\DefEq
{
\[
\|f(x+a)-f(x)\|_B\le \left\|\frac{df(x)}{dx}\right\|\,\delta=\epsilon
\]
}
{derivative and continuos map, algebra, 3}

\DefEq
{
\[\delta=\frac\epsilon{\displaystyle\left\|\frac{df(x)}{dx}\right\|}\]
}
{delta=epsilon/|df|}

\DefEq
{
\[
\|a\|_A<\delta
\]
}
{|a|<delta}

\DefEq
{
\[\delta=\frac\epsilon{\|f\|}\]
}
{delta=epsilon/|f|}

\DefEquation
{
\frac{df(x)}{dx}
=\ShowSymbol{standard component of derivative}{}
(e_{\gii}\otimes e_{\gij})\circ F_k
}
{derivative, algebra, standard representation}

\DefEq
{
\symb{\frac{d^{k\cdot\gi{ij}} f(x)}{dx}}
{standard component of derivative}{}
}
{standard component of derivative}

\DefEquation
{
\frac{df(x)}{dx}
=\ShowSymbol{standard component of derivative}{I0}
e_{\gii}\otimes e_{\gij}
}
{derivative, algebra, standard representation I0}

\DefEq
{
\symb{\frac{d^{\gi{ij}} f(x)}{dx}}
{standard component of derivative}{I0}
}
{standard component of derivative I0}

\DefEquation
{
\frac{|f(x)|_2}{|x|_1}=\left|f\left(\frac x{|x|_1}\right)\right|_2
}
{norm of linear map, algebra, 1}

\DefEquation
{
\begin{matrix}
f(rx)=rf(x)
&
r\in R
\end{matrix}
}
{linear map and real number}

\DefEq
{
\[
\frac{|f(rx)|_2}{|rx|_1}
=\frac{|r|\ |f(x)|_2}{|r|\ |x|_1}
=\frac{|f(x)|_2}{|x|_1}
\]
}
{linear map and real number, 1}

\DefEquation
{
\|f\|=\text{sup}\{|f(x)|_2:|x|_1=1\}
}
{norm of linear map, algebra}

%% file: Derivative.16.Ref.tex

\DefEq
{
\RefDefinition{linear map from A1 to A2, commutative module},
\RefDefinition{linear map from A1 to A2, algebra},
\RefDefinition{differentiable map, algebra}\Pt
}
{ref definition linear map from A1 to A2}

\DefEq
{
\RefTheorem{Tensor product is distributive over sum},
\RefTheorem{tensor product is associative}
}
{ref tensor product is associative and distributive}

\DefEq
{
\RefItem{map f generated by basis I}
}
{ref map f generated by basis I}

\DefEq
{
\RefItem{standard representation of map A1 A2, associative algebra}
}
{ref standard representation of map A1 A2, associative algebra}

\DefEq
{
\RefTheorem{coordinates of map A1 A2, algebra}.
}
{ref theorem coordinates of map A1 A2, algebra}

\DefEq
{
\RefTheorem{linear map, 0, D algebra}.
}
{ref linear map, 0}

\DefEq
{
\RefDefinition{algebra over ring}.%
}
{ref algebra over ring}

\DefEq
{
\RefTheorem{V times->V otimes}%
}
{ref theorem V times->V otimes}

%% file: Second.Derivative.16.English.tex

\input{Second.Derivative.16.Eq}

\Chapter{Derivative of Second Order of Map of \texorpdfstring{$D$}{D}-Algebra}
\labelChapter{Derivative of Second Order of Map of Algebra}

\Section{Derivative of Second Order of Map of \texorpdfstring{$D$}{D}-Algebra}

Let $D$ be the complete commutative ring of characteristic $0$.
Let $A$, $B$ be Banach $D$\Hyph modules.
Let
\ShowEq{f:A->B}fAB
be differentiable map.
According to the remark \ref{remark: differential L(A,A)},
the derivative is map
\EqParm{x->df/dx in L(A,B)}{A=AB}
According to the theorem
\RefTheorem{L(A;B) is Banach module},
set
\ShowEq{L(A;B)}DAB{}
is Banach $D$\Hyph module.
Therefore, we may consider the question,
if map
\ShowEq{df/dx}
is differentiable.

According to the definition
\RefDefinition{differentiable map, algebra}
and to the theorem
\RefTheorem{dfa1p/dx=df/dx a1p},
\ShowEq{derivative of map df, algebra}
where
\ShowEq{derivative of map df, algebra, o2}
is such continuous map, that
\ShowEq{derivative of map df, algebra, o2 lim}
According to definition
\RefDefinition{differentiable map, algebra},
the map
\ShowEq{increment of derivative, algebra}
is linear map of variable $a_2$. From the equality
\EqRef{derivative of map df, algebra}
it follows that map
\ShowEq{increment of derivative, algebra}
is linear map of variable $a_1$.
Therefore, the map
\ShowEq{increment of derivative, algebra}
is bilinear map.

\ShowDefinition{derivative of Second Order, algebra}

\begin{remark}
According to definition
\RefDefinition{derivative of Second Order, algebra}
the derivative of second order of the map $f$ is map
\ShowEq{differential L(A,A;A)}
According to the theorem
\ShowEq{ref theorem tensor product and polylinear map}
we may consider also expression
\ShowEq{differential AA A,A}
Then
\ShowEq{differential L(AA;A)}
We use the same notation
\ShowEq{differential D2A2}
for maps
\EqRef{differential L(A,A;A)}
and
\EqRef{differential L(AA;A)}
because of the nature of the argument it is clear what kind of map we consider.
\qed
\end{remark}

\ShowTheorem{Differential of Second Order, algebra, representation}
\begin{proof}
The theorem follow from the definition
\RefDefinition{derivative of Second Order, algebra}
and from the theorem
\ShowEq{ref polylinear map}
\end{proof}

\ShowDefinition{derivative of Order n, algebra}

\Section{Taylor Series}
\labelSection{Taylor Series}

Let $D$ be the complete commutative ring of characteristic $0$.
Let $A$ be associative $D$\Hyph algebra.
Let $p_n(x)$ be the monomial of power $n$, $n>0$,
in one variable over $D$\Hyph algebra $A$.
Let $\Basis F$ be basis of algebra
\ShowEq{L(A;B)}DAA.

\begin{theorem}
\labelTheorem{derivative of f(x)x, algebra}
For any $m>0$ the following equality is true
\ShowEq{derivative of f(x)x, algebra}
where symbol $\widehat{h^i}$ means absense of variable $h^i$ in the list.
\end{theorem}
\begin{proof}
For $m=1$, this is corollary  of the equality
\EqRef{derivative of product, algebra}
and the theorem
\RefTheorem{derivative of linear map}
\ShowEq{derivative of f(x)x, 1, algebra}
Assume, \EqRef{derivative of f(x)x, algebra} is true
for $m-1$. Then
\ShowEq{derivative of f(x)x, m-1, algebra}
Since
\ShowEq{dhi=0}
then, using the equality \EqRef{derivative of product, algebra},
we get
\ShowEq{derivative of f(x)x, m, algebra}
The difference between
equalities \EqRef{derivative of f(x)x, algebra} and
\EqRef{derivative of f(x)x, m, algebra}
is only in form of presentation.
We proved the theorem.
\end{proof}

\begin{theorem}
\labelTheorem{derivative of pn = 0, algebra}
For any $n\ge 0$, the following equality is true
\ShowEq{derivative n1 of pn, algebra}
\end{theorem}
\begin{proof}
Since $p_0(x)=a_0$, $a_0\in D$, then for $n=0$ theorem is corollary
of definition
\RefDefinition{differentiable map, algebra}.
Let statement of theorem is true for $n-1$. According to theorem
\RefTheorem{derivative of f(x)x, algebra},
when $f(x)=p_{n-1}(x)$, we get
\ShowEq{d(p(n-1)xan)}
According to suggestion of induction all monomials are equal $0$.
\end{proof}

\begin{theorem}
\labelTheorem{derivative of pn = 0, m < n, algebra}
If $m<n$, then the following equality is true
\ShowEq{m derivative of polinom pn, algebra}
\end{theorem}
\begin{proof}
For $n=1$, the following equality is true
\ShowEq{d0p1(0)=0}
Assume that statement is true
for $n-1$. Then according to theorem \RefTheorem{derivative of f(x)x, algebra}
\ShowEq{dm p(n-1)xan}
First term equal $0$ because $x=0$.
Because $m-1<n-1$, then
rest terms equal $0$ according to suggestion of induction.
We proved the statement of theorem.
\end{proof}

When $h_1=...=h_n=h$, we assume
\ShowEq{dfxh}
This notation does not create ambiguity, because we can determine function
according to number of arguments.

\begin{theorem}
\labelTheorem{n derivative of polinom pn, algebra}
For any $n>0$, the following equality is true
\ShowEq{n derivative of polinom pn, algebra}
\end{theorem}
\begin{proof}
For $n=1$, the following equality is true
\ShowEq{n derivative of polinom pn, 2, algebra}
Assume the statement is true
for $n-1$. Then according to theorem \RefTheorem{derivative of f(x)x, algebra}
\ShowEq{n derivative of polinom pn, 1, algebra}
First term equal $0$ according to theorem
\RefTheorem{derivative of pn = 0, algebra}.
The rest $n$ terms equal, and according to suggestion of induction
from the equality \EqRef{n derivative of polinom pn, 1, algebra}
it follows
\ShowEq{n derivative of polinom pn, 3, algebra}
Therefore, statement of theorem is true for any $n$.
\end{proof}

Let $p(x)$ be polynomial of power $n$.\,\footnote{I consider
Taylor polynomial for polynomials by analogy with
construction of Taylor polynomial in \citeBib{Fikhtengolts: Calculus volume 1}, p. 246.}
\[
p(x)=p_0+p_{1i_1}(x)+...+p_{ni_n}(x)
\]
We assume sum by index $i_k$ which enumerates terms
of power $k$.
According to theorem \RefTheorem{derivative of pn = 0, algebra},
\RefTheorem{derivative of pn = 0, m < n, algebra},
\RefTheorem{n derivative of polinom pn, algebra}
\ShowEq{Taylor polynomial 1, algebra}
Therefore, we can write
\ShowEq{Taylor polynomial 2, algebra}
This representation of polynomial is called
\AddIndex{Taylor polynomial}{Taylor polynomial, algebra}.
If we consider substitution of variable $x=y-y_0$, then considered above construction
remain true for polynomial
\[
p(y)=p_0+p_{1i_1}(y-y_0)+...+p_{ni_n}(y-y_0)
\]
Therefore
\ShowEq{Taylor polynomial 3, algebra}

Assume that function $f(x)$ is differentiable at point $x_0$
up to any order.\,\footnote{I explore
construction of Taylor series by analogy with
construction of Taylor series in \citeBib{Fikhtengolts: Calculus volume 1}, p. 248, 249.}

\ShowTheorem{n derivative equal 0, algebra}
\begin{proof}
When $n=1$ this statement follows from the equality
\EqRef{differential of map, t, algebra}.

Let statement be true for $n-1$.
Map
\ShowEq{n derivative equal 0, 1, algebra}
satisfies to condition
\ShowEq{n derivative equal 0, 2, algebra}
According to suggestion of induction
\[
f_1(x_0+th)=o(t^{n-1})
\]
Then the equality
\EqRef{derivative, linear path}
gets form
\[
o(t^{n-1})=\lim_{t\rightarrow 0,\ t\in R}(t^{-1}f(x+th))
\]
Therefore,
\[
f(x+th)=o(t^n)
\]
\end{proof}

\ShowEq{Taylor polynomial}

\ShowEq{Taylor series}

%% file: Second.Derivative.16.Eq.tex

\input{Second.Derivative.16.Ref}

\ifx\Semafor\ValueOff

\DefEq
{
$o_2:A\rightarrow \mathcal L(D;A\rightarrow B)$
}
{derivative of map df, algebra, o2}

\DefEq
{
\[
\lim_{a_2\rightarrow 0}\frac{\|o_2(a_2)\|_B}{\|a_2\|_A}=0
\]
}
{derivative of map df, algebra, o2 lim}

\DefEq
{
\[
\frac{d^2 f(x)}{dx^2}\circ(a_1\otimes a_2)=
\frac{d^2 f(x)}{dx^2}\circ(a_1;a_2)
\]
}
{differential AA A,A}

\DefEquation
{
x\in A\rightarrow \frac{d^2 f(x)}{dx^2}\in\mathcal L(D;A\times A\rightarrow A)
}
{differential L(A,A;A)}

\DefEquation
{
\begin{split}
\frac{d^2 f(x)}{dx^2}\in&\ \mathcal L(D;A\otimes A\rightarrow A)
\\
\frac{d^2 f}{dx^2}:A\rightarrow&\ \mathcal L(D;A\otimes A\rightarrow A)
\end{split}
}
{differential L(AA;A)}

\DefEq
{
$\displaystyle\frac{d^2 f}{dx^2}$
}
{differential D2A2}

\DefEq
{
\[
\frac{df(x)(F\circ x)}{dx}\circ h_1
=\left(\frac{df(x)}{dx}\circ h_1\right)(F\circ x)+
f(x)(F\circ h_1)
\]
}
{derivative of f(x)x, 1, algebra}

\DefEq
{
\[
f_1(x_0)=
\left.\frac{d f_1(x)}{dx}\right|_{x=x_0}\circ h=...
=\left.\frac{d^{n-1} f_1(x)}{d x^{n-1}}\right|_{x=x_0}\circ h^{n-1}=0
\]
}
{n derivative equal 0, 2, algebra}

\DefEq
{
\[
f_1(x)=\frac{df(x)}{dx}\circ h
\]
}
{n derivative equal 0, 1, algebra}

\DefEq
{
\[
\frac{d^n p_n(x)}{dx^n}\circ h
=n\left(\frac{d^{n-1} p_{n-1}(x)}{dx^{n-1}}\circ h\right)ha_n
=n(n-1)!p_{n-1}(h)ha_n=n!p_n(h)
\]
}
{n derivative of polinom pn, 3, algebra}

\DefEquation
{
\begin{split}
\frac{d^n p_n(x)}{d x^n}\circ h^n
=&\left(\frac{d^n p_{n-1}(x)}{d x^n}\circ h^n\right)(F\circ x)a_n
\\+&\left(\frac{d^{n-1} p_{n-1}(x)}{dx^{n-1}}\circ h^{n-1}\right)(F\circ h)a_n
\\
+&...
+\left(\frac{d^{n-1} p_{n-1}(x)}{dx^{n-1}}\circ h^{n-1}\right)(F\circ h)a_n
\end{split}
}
{n derivative of polinom pn, 1, algebra}

\DefEq
{
\[
\frac{d p_1(x)}{dx}\circ h
=\frac{d a_0(F\circ x)a_1}{d x}\circ h=a_0(F\circ h)a_1=1!p_1(h)
\]
}
{n derivative of polinom pn, 2, algebra}

\DefEq
{
\[
\begin{split}
\frac{d^{m-1}f(x)(F\circ x)}{dx^{m-1}}\circ(h_1;...;h_{m-1})
=&\frac{d^{m-1}f(x)}{dx^{m-1}}\circ(h_1;...;h_{m-1})(F\circ x)
\\
+&\frac{d^{m-2}f(x)}{dx^{m-2}}\circ(\widehat{h_1};...;h_{m-2};h_{m-1})(F\circ h_1)
+...
\\
+&\frac{d^{m-2}f(x)}{dx^{m-2}}\circ(h_1;...;h_{m-2};\widehat{h_{m-1}})(F\circ h_{m-1})
\end{split}
\]
}
{derivative of f(x)x, m-1, algebra}

\DefEquation
{
\begin{split}
&\,\frac{d^mf(x)(F\circ x)}{dx^m}\circ(h_1;...;h_{m-1};h_m)
\\=&\,\frac{d^mf(x)}{dx^m}\circ(h_1;...;h_{m-1};h_m)(F\circ x)
\\
+&\,\frac{d^{m-1}f(x)}{dx^{m-1}}\circ(h_1;...;h_{m-2};h_{m-1})(F\circ h_m)
\\
+&\,\frac{d^{m-1}f(x)}{dx^{m-1}}\circ(\widehat{h_1};...;h_{m-2};h_{m-1};h_m)(F\circ h_1)
+...
\\
+&\,\frac{d^{m-1}f(x)}{dx^{m-1}}\circ(h_1;...;h_{m-2};\widehat{h_{m-1}};h_m)(F\circ h_{m-1})
\end{split}
}
{derivative of f(x)x, m, algebra}

\DefEq
{
\begin{align*}
\frac{d^{n+1}p_n(x)}{dx^{n+1}}\circ(h_1;...;h_{n+1})
&=\frac{d^{n+1}p_{n-1}(x)(F\circ x)a_n}{dx^{n+1}}\circ(h_1;...;h_{n+1})
\\
&=\frac{d^{n+1}p_{n-1}(x)}{dx^{n+1}}\circ (h_1;...;h_{n+1})(F\circ x)a_n
\\
&+\frac{d^np_{n-1}(x)}{dx^n}\circ(\widehat{h_1};...;h_n;h_{n+1})(F\circ h_1)a_n
+...
\\
&+\frac{d^np_{n-1}(x)}{dx^n}\circ(h_1;...;h_n;\widehat{h_{n+1}})(F\circ h_{n+1})a_n
\end{align*}
}
{d(p(n-1)xan)}

\DefEq
{
\begin{align*}
&\,\frac{d^mp_{n-1}(x)(F\circ x)a_n}{dx^m}\circ(h_1;...;h_m)
\\=&\,\frac{d^mp_{n-1}(x)}{dx^m}\circ(h_1;...;h_m)(F\circ x)a_n
\\
+&\,\frac{d^{m-1}p_{n-1}(x)}{dx^{m-1}}\circ(\widehat{h_1};...;h_{m-1};h_m)(F\circ h_1)a_n
+...
\\
+&\,\frac{d^{m-1}p_{n-1}(x)}{dx^{m-1}}\circ(h_1;...;h_{m-1};\widehat{h_m})(F\circ h_m)a_n
\end{align*}
}
{dm p(n-1)xan}

\DefEquation
{
\begin{split}
\frac{d^mf(x)(F\circ x)}{dx^m}\circ(h_1;...;h_m)
=&\frac{d^mf(x)}{dx^m}\circ(h_1;...;h_m)(F\circ x)
\\
+&\frac{d^{m-1}f(x)}{dx^{m-1}}\circ(\widehat{h_1};...;h_{m-1};h_m)(F\circ h_1)
+...
\\
+&\frac{d^{m-1}f(x)}{dx^{m-1}}\circ(h_1;...;h_{m-1};\widehat{h_m})(F\circ h_m)
\end{split}
}
{derivative of f(x)x, algebra}

\DefEq
{
\[
\frac{d^kp(0)}{dx^k}\circ x=k!p_{ki_k}(x)
\]
}
{Taylor polynomial 1, algebra}

\DefEq
{
\[
p(x)=p_0
+\frac{1}{1!}\frac{dp(0)}{dx}\circ x
+\frac{1}{2!}\frac{d^2p(0)}{dx^2}\circ x^2
+...+\frac{1}{n!}\frac{d^np(0)}{dx^n} \circ x^n
\]
}
{Taylor polynomial 2, algebra}

\DefEq
{
\[
\begin{split}
p(y)&=p_0+\frac{1}{1!}\frac{dp(y_0)}{dy}\circ(y-y_0)
+\frac{1}{2!}\frac{d^2p(y_0)}{dy^2}\circ(y-y_0)^2+...
\\
&+\frac{1}{n!}\frac{d^np(y_0)}{dy^n}\circ(y-y_0)^n
\end{split}
\]
}
{Taylor polynomial 3, algebra}

\DefEq
{
\[
\frac{d^n p_n(x)}{d x^n}\circ h^n=n!p_n(h)
\]
}
{n derivative of polinom pn, algebra}

\DefEq
{
\[
\frac{d^nf(x)}{dx^n}\circ h^n
=\frac{d^nf(x)}{d x^n}\circ(h_1;...;h_n)
\]
}
{dfxh}

\DefEq
{
\[
d^0 p_1(0)=a_0(F\circ x)a_1|_{x=0}=0
\]
}
{d0p1(0)=0}

\DefEq
{
\[
\left.\frac{d^m p_n(x)}{dx^m}\right|_{x=0}=0
\]
}
{m derivative of polinom pn, algebra}

\DefEq
{
\[
\frac{d^{n+1} p_n(x)}{dx^{n+1}}=0
\]
}
{derivative n1 of pn, algebra}

\DefEq
{
$\displaystyle\DfTwo$
}
{increment of derivative, algebra}

\DefEq
{
$\displaystyle\frac{dh_i}{dx} =0$,
}
{dhi=0}

\DefEquation
{
\frac{df(x+a_2)}{dx}\circ a_1-\frac{df(x)}{dx}\circ a_1
=\DfTwo
+o_2(a_2)
}
{derivative of map df, algebra}

%% file: Second.Derivative.16.Ref.tex

\DefEq
{
\RefTheorem{tensor product and polylinear map},
}
{ref theorem tensor product and polylinear map}

\DefEq
{
\RefTheorem{polylinear map, algebra}.
}
{ref polylinear map}

\DefEq
{
\RefTheorem{module L(A;A) is algebra}
}
{ref module L(A;A) is algebra}

%% file: Indefinite.Integral.English.tex

\input{Indefinite.Integral.Eq}

\Chapter{Method of Successive Differentiation}

\Section{Indefinite Integral}
\labelSection{indefinite integral}

\ShowDefinition{indefinite integral}

In this section we consider
integration as operation inverse to differentiation.
As a matter of fact, we consider procedure of solution of ordinary differential
equation
\ShowEq{df=g}

\ShowExample{differential equation y=xx, real number}

\ShowTheorem{int=x3 algebra}
\ShowProof{int=x3 algebra}

\ShowRemark{int=x3 algebra}

\ShowRemark{notation int=x3 algebra}

\ShowRemark{differential equation 3x^2 does not possess a solution}

\Section{Exponent}

In a field we can define exponent as solution
of differential equation
\ShowEq{exponent over field}
It is evident that we cannot write such equation for division ring.
However we can write the equation
\EqRef{exponent over field}
as follows
\ShowEq{exponent derivative over field}
This equation is closer to our goal,
however there is the question: in which order
we should multiply $y$ and $h$?
To answer this question we change format of the equation
\EqRef{exponent derivative over field}
\ShowEq{exponent derivative over division ring}
Hence, our goal is to solve differential equation
\EqRef{exponent derivative over division ring}
with initial condition \[y(0)=1\]

\begin{definition}
\labelDefinition{set of permutations in derivative}
{\it
For $n\ge 0$, let
\ShowEq{set of permutations}
be set of permutations
\EqParm{derivative, permutation}{K=n,J=}
such that each permutation $\sigma$ has following properties
\ifx\setCACAA\undefined
\StartLabelItem
\else
\StartLabelItem[definition]
\fi
\begin{enumerate}
\item If there exist $i$, $j$, $i\ne j$, such that
$\sigma(h_i)$ is situated in product
\ShowEq{sigma yh1n}
on the left side of $\sigma(h_j)$ and $\sigma(h_j)$ is situated
on the left side of $\sigma(y)$, then $i<j$.
\labelItem{derivative, permutation n, left}
\item If there exist $i$, $j$, $i\ne j$, such that
$\sigma(h_i)$ is situated in product
\ShowEq{sigma yh1n}
on the right side of $\sigma(h_j)$ and $\sigma(h_j)$ is situated
on the right side of $\sigma(y)$, then $i>j$.
\labelItem{derivative, permutation n, right}
\end{enumerate}
}
\qed
\end{definition}

\begin{lemma}
\labelLemma{permutation n+1 generated by permutation n}
{\it
For $n\ge 0$,
for any permutation
\ShowEq{t in Sn+1},
there exists unique permutation
\ShowEq{s in Sn}
such that either
\DrawEq{tau yh1...=sigma yh... 1}{->}
or
\DrawEq{tau yh1...=sigma yh... 2}{->}
}
\end{lemma}
\begin{proof}
Consider product
\DrawEq{tau yh...}{->}
Since $n+1$ is largest index then, according to conditions
\RefItem{derivative, permutation n, left},
\RefItem{derivative, permutation n, right},
$\tau(h_{n+1})$
is written either immediately before or immediately after $\tau(y)$.
Therefore, product
\eqRef{tau yh...}{->}
has either form
\eqRef{tau yh1...=sigma yh... 1}{->}
or
\eqRef{tau yh1...=sigma yh... 2}{->}.
Therefore for any permutation
\ShowEq{t in Sn+1},
there exists corresponding
permutation
\ShowEq{s in Sn}.
\end{proof}

\begin{lemma}
\labelLemma{permutation n generates permutation n+1, 1}
{\it
For $n\ge 0$,
for any permutation
\ShowEq{s in Sn}.
there exists unique permutation
\ShowEq{t in Sn+1}
such that
\DrawEq{tau yh1...=sigma yh... 1}{<-}
}
\end{lemma}
\begin{proof}
Consider permutation
\ShowEq{sigma yh... 1}
To write down permutation
\EqParm{derivative, permutation}{K=n+1,J=}
which satisfies
\eqRef{tau yh1...=sigma yh... 1}{<-},
we write down the tuple
\ShowEq{tau yh1... top}
instead of the top tuple
\ShowEq{sigma yh... 1 top}
and write down the tuple
\ShowEq{tau yh1... 1 down}
instead of the down tuple
\ShowEq{sigma yh... 1 down}
In expression
\DrawEq{tau yh...}{}
$h_{n+1}$ is written immediately before $y$.
Since $n+1$ is the largest value of index,
then permutation $\tau\in S(n+1)$.
\end{proof}

\begin{lemma}
\labelLemma{permutation n generates permutation n+1, 2}
{\it
For $n\ge 0$,
for any permutation
\ShowEq{s in Sn}.
there exists unique permutation
\ShowEq{t in Sn+1}
such that
\DrawEq{tau yh1...=sigma yh... 2}{<-}
}
\end{lemma}
\begin{proof}
Consider permutation
\ShowEq{sigma yh... 2}
To write down permutation
\EqParm{derivative, permutation}{K=n+1,J=}
which satisfies
\eqRef{tau yh1...=sigma yh... 2}{<-},
we write down the tuple
\ShowEq{tau yh1... top}
instead of the top tuple
\ShowEq{sigma yh... 2 top}
and write down the tuple
\ShowEq{tau yh1... 2 down}
instead of the down tuple
\ShowEq{sigma yh... 2 down}
In expression
\DrawEq{tau yh...}{}
$h_{n+1}$ is written immediately after $y$.
Since $n+1$ is the largest value of index,
then permutation $\tau\in S(n+1)$.
\end{proof}

\begin{theorem}
\labelTheorem{permutation sigma generates permutation tau, n=2k}
\ShowEq{tau yh1...=sigma yh... +}
\end{theorem}
\begin{proof}
The theorem follows from lemmas
\RefLemma{permutation n+1 generated by permutation n},
\RefLemma{permutation n generates permutation n+1, 1},
\RefLemma{permutation n generates permutation n+1, 2}.
\end{proof}

\begin{theorem}
\labelTheorem{exponent derivative n over division ring}
If a map $y$ is solution of differential equation
\EqRef{exponent derivative over division ring}
then the derivative of order $n$ of function $y$
has form
\ShowEq{exponent derivative n over division ring}
\end{theorem}
\begin{proof}
We prove this statement by induction.
For $n=0$, this statement is true because this is the statement $y=y$.
For $n=1$, the statement is true because this is differential equation
\EqRef{exponent derivative over division ring}.
Let the statement be true for $n=k-1$.
Hence
\ShowEq{exponent derivative n=k-1}
According to the definition
\ShowEq{ref derivative of Order n, algebra}
the derivative of order $k$ has form
\ShowEq{exponent derivative n=k, 1}
From equalities
\EqRef{exponent derivative over division ring},
\EqRef{exponent derivative n=k, 1}
it follows that
\ShowEq{exponent derivative n=k, 2}
From equalities
\EqRef{tau yh1...=sigma yh... +},
\EqRef{exponent derivative n=k, 2},
it follows that
the statement of theorem is true for $n=k$.
We proved the theorem.
\end{proof}

\begin{theorem}
\labelTheorem{exponent}
The solution of differential equation
\EqRef{exponent derivative over division ring}
with initial condition $y(0)=1$ is
\AddIndex{exponent}{exponent}
\ShowEq{exponent over division ring}
that has following Taylor series expansion
\ShowEq{exponent Taylor series}
\end{theorem}
\begin{proof}
The derivative of order $n$ has $2^n$ items.
In fact, the derivative of order $1$ has $2$ items,
and each differentiation increase number of items twice.
From initial condition $y(0)=1$ and theorem
\RefTheorem{exponent derivative n over division ring},
it follows that the derivative of order $n$
of required solution has form
\ShowEq{exponent derivative n, x=0}
Taylor series expansion
\EqRef{exponent Taylor series}
follows
from
\EqRef{exponent derivative n, x=0}.
\end{proof}

\begin{theorem}
\labelTheorem{exponent of sum}
The equality
\ShowEq{exponent of sum}
is true iff
\ShowEq{exponent of sum, 1}
\end{theorem}
\begin{proof}
To prove the theorem it is enough to consider Taylor series
\ShowEq{exponent a b a+b}
Let us multiply expressions
\EqRef{exponent a} and \EqRef{exponent b}.
The sum of monomials of order $3$ has form
\ShowEq{exponent ab 3}
and in general does not equal expression
\ShowEq{exponent a+b 3}
The proof of statement that \EqRef{exponent of sum} follows from
\EqRef{exponent of sum, 1}
is trivial.
\end{proof}

\Section{Hyperbolic Trigonometry}

To introduce hyperbolic sine
\ShowEq{y1=sh}
and hyperbolic cosine
\ShowEq{y2=ch}
we consider system of differential equations over real field
\ShowEq{d sh ch R}
and initial conditions
\DrawEq{d sh ch init}{real}
The system of differential equations
\EqRef{d sh ch R},
\eqRef{d sh ch init}{real}
is equivalent to differential equations
\ShowEq{d2 sh}
\ShowEq{d2 ch}

In Banach $D$\Hyph algebra $A$,
we write system of differential equations
\EqRef{d sh ch R},
\eqRef{d sh ch init}{real}
in the form similar to differential equation
\EqRef{exponent derivative over division ring}
\ShowEq{d sh ch A}
\DrawEq{d sh ch init}{}

\begin{theorem}
If maps $y_1$, $y_2$ are solution of system of differential equations
\EqRef{d sh ch A}
then the derivatives of order $n$ of maps $y_1$, $y_2$
have form
\ShowEq{sh ch derivative n}
\end{theorem}
\begin{proof}
We prove this statement by induction.
For $n=0$, this statement is true because this is the statement $y_1=y_1$, $y_2=y_2$.
For $n=1$ the statement is true because this is system of differential equations
\EqRef{d sh ch A}.
Let the statement be true for $n=k-1$.
Hence
\ShowEq{sh ch derivative n=k-1}
where
$i=1$, $2$ and
\ShowEq{j=i,2-i}
According to the definition
\ShowEq{ref derivative of Order n, algebra}
the derivative of order $k$ has form
\ShowEq{sh ch derivative n=k, 1}
From
\EqRef{d sh ch A},
\EqRef{sh ch derivative n=k, 1},
it follows that
\ShowEq{sh ch derivative n=k, 2}
where
\ShowEq{l=2-j}
From
\EqRef{j=i,2-i},
\EqRef{l=2-j},
it follows that
\ShowEq{l=i,2-i}
From
\EqRef{tau yh1...=sigma yh... +},
\EqRef{sh ch derivative n=k, 2},
it follows that
the statement of theorem is true for $n=k$.
We proved the theorem.
\end{proof}

\begin{theorem}
\labelTheorem{sh ch}
The solution of system of differential equations
\EqRef{d sh ch A}
with initial condition
\DrawEq{d sh ch init}{Taylor}
is the tuple of maps
\begin{itemize}
\item
\AddIndex{hyperbolic sine}{hyperbolic sine}
\ShowEq{hyperbolic sine}
that has following Taylor series expansion
\ShowEq{hyperbolic sine Taylor series}
\item
\AddIndex{hyperbolic cosine}{hyperbolic cosine}
\ShowEq{hyperbolic cosine}
that has following Taylor series expansion
\ShowEq{hyperbolic cosine Taylor series}
\end{itemize}
\end{theorem}
\begin{proof}
The derivative of order $n$ has $2^n$ items.
In fact, the derivative of order $1$ has $2$ items,
and each differentiation increase number of items twice.
From initial condition
\eqRef{d sh ch init}{Taylor}
and from the equality
\EqRef{sh ch derivative n},
it follows that the derivative of order $n$
of required solution has form
\ShowEq{sh ch derivative n, x=0}
Taylor series expansions
\EqRef{hyperbolic sine Taylor series},
\EqRef{hyperbolic cosine Taylor series}
follow
from
\EqRef{sh ch derivative n, x=0}.
\end{proof}

\begin{theorem}
The system of differential equations
\DrawEq{d sh A}{sinh}
\ShowEq{d ch A}
is equivalent to the differential equation
\ShowEq{d2 sh A}
\end{theorem}
\begin{proof}
From the differential equation
\eqRef{d sh A}{sinh},
it follows that
\DrawEq{d2 sh A 1}{sinh}
From
\EqRef{d ch A},
\eqRef{d2 sh A 1}{sinh},
it follows that
\ShowEq{d2 sh A 2}
The differential equation
\EqRef{d2 sh A}
follows from
\EqRef{d2 sh A 2}.
\end{proof}

\Section{Elliptical Trigonometry}

To introduce sine
\ShowEq{y1=sin}
and cosine
\ShowEq{y2=cos}
we consider system of differential equations over real field
\ShowEq{d sin cos R}
and initial conditions
\DrawEq{d sh ch init}{real sin}
The system of differential equations
\EqRef{d sin cos R},
\eqRef{d sh ch init}{real sin}
is equivalent to differential equations
\ShowEq{d2 sin}
\ShowEq{d2 cos}

In Banach $D$\Hyph algebra $A$,
we write system of differential equations
\EqRef{d sin cos R},
\eqRef{d sh ch init}{real sin}
in the form similar to differential equation
\EqRef{exponent derivative over division ring}
\ShowEq{d sin cos A}
\DrawEq{d sh ch init}{}

\begin{theorem}
If maps $y_1$, $y_2$ are solution of system of differential equations
\EqRef{d sin cos A}
then the derivatives of order $n$ of function $y$
have form
\begin{itemize}
\item $n=2k$
\ShowEq{sin cos derivative n=2k}
\item $n=2k+1$
\ShowEq{sin cos derivative n=2k+1}
\end{itemize}
\end{theorem}
\begin{proof}
We prove this statement by induction.
For $k=0$, this statement is true:
\begin{itemize}
\item
Since $n=2k$,
then this is the statement $y_1=y_1$, $y_2=y_2$.
\item
Since $n=2k+1$,
then the statement follows from the system of differential equations
\EqRef{d sin cos A}.
\end{itemize}
Let the statement be true for $k=l-1$.
Hence, for $n=2k+1$,
\ShowEq{sin derivative n=2k+1 k=l-1}
\ShowEq{cos derivative n=2k+1 k=l-1}
According to the definition
\ShowEq{ref derivative of Order n, algebra}
from
\EqRef{sin derivative n=2k+1 k=l-1},
\EqRef{cos derivative n=2k+1 k=l-1},
it follows that
the derivatives of order $m=n+1$ have form
\ShowEq{sin derivative m=2k+2 k=l-1, 1}
\ShowEq{cos derivative m=2k+2 k=l-1, 1}
From
\EqRef{d sin cos A},
\EqRef{sin derivative m=2k+2 k=l-1, 1},
\EqRef{cos derivative m=2k+2 k=l-1, 1},
it follows that
\ShowEq{sin derivative m=2k+2 k=l-1, 2}
\ShowEq{cos derivative m=2k+2 k=l-1, 2}
From
\EqRef{tau yh1...=sigma yh... +},
\EqRef{sin derivative m=2k+2 k=l-1, 2},
\EqRef{cos derivative m=2k+2 k=l-1, 2},
it follows that
\ShowEq{sin cos derivative m=2k+2 k=l-1, 2}
Since
\ShowEq{m=2l}
then equalities
\EqRef{sin cos derivative n=2k}
are true for $k=l$, $n=2k$.

According to the definition
\ShowEq{ref derivative of Order n, algebra}
from
\EqRef{sin cos derivative n=2k},
it follows that
the derivatives of order $m=n+1$ have form
\ShowEq{sin derivative m=2k+1 k=l, 1}
\ShowEq{cos derivative m=2k+1 k=l, 1}
From
\EqRef{d sin cos A},
\EqRef{sin derivative m=2k+1 k=l, 1},
\EqRef{cos derivative m=2k+1 k=l, 1},
it follows that
\ShowEq{sin derivative m=2k+1 k=l, 2}
\ShowEq{cos derivative m=2k+1 k=l, 2}
Since
\ShowEq{m=2k+1}
then equalities
\EqRef{sin cos derivative n=2k+1}
follow from
\EqRef{tau yh1...=sigma yh... +},
\EqRef{sin derivative m=2k+1 k=l, 2},
\EqRef{cos derivative m=2k+1 k=l, 2}.
\EqRef{sin cos derivative n=2k}
for $k=l$, $n=2k+1$.

We proved the theorem.
\end{proof}

\begin{theorem}
\labelTheorem{sin cos}
The solution of system of differential equations
\EqRef{d sin cos A}
with initial condition
\DrawEq{d sh ch init}{Taylor syn}
is the tuple of maps
\begin{itemize}
\item
\AddIndex{sine}{sine}
\ShowEq{sine}
that has following Taylor series expansion
\ShowEq{sine Taylor series}
\item
\AddIndex{cosine}{cosine}
\ShowEq{cosine}
that has following Taylor series expansion
\ShowEq{cosine Taylor series}
\end{itemize}
\end{theorem}
\begin{proof}
The derivative of order $n$ has $2^n$ items.
In fact, the derivative of order $1$ has $2$ items,
and each differentiation increase number of items twice.
From initial condition
\eqRef{d sh ch init}{Taylor syn}
and from equalities
\EqRef{sin cos derivative n=2k},
\EqRef{sin cos derivative n=2k+1},
it follows that the derivative of order $n$
of required solution has form
\begin{itemize}
\item $n=2k$
\ShowEq{sin cos derivative n=2k, x=0}
\item $n=2k+1$
\ShowEq{sin cos derivative n=2k+1, x=0}
\end{itemize}
Taylor series expansions
\EqRef{sine Taylor series},
\EqRef{cosine Taylor series}
follow
from
\EqRef{sin cos derivative n=2k, x=0},
\EqRef{sin cos derivative n=2k+1, x=0}.
\end{proof}

\begin{theorem}
The system of differential equations
\DrawEq{d sh A}{sine}
\ShowEq{d cos A}
is equivalent to the differential equation
\ShowEq{d2 sin A}
\end{theorem}
\begin{proof}
From the differential equation
\eqRef{d sh A}{sine},
it follows that
\DrawEq{d2 sh A 1}{sine}
From
\EqRef{d cos A},
\eqRef{d2 sh A 1}{sine},
it follows that
\ShowEq{d2 sin A 2}
The differential equation
\EqRef{d2 sin A}
follows from
\EqRef{d2 sin A 2}.
\end{proof}

%% file: Indefinite.Integral.Eq.tex

\DefEq%
{%
\def\S{\sigma}%
\def\K{k-1}%
}%
{K=k-1}%

\DefEq
{
\def\S{\sigma}%
\def\K{n}%
}
{K=n}

\DefEq
{
\def\S{\tau}%
\def\K{n+1}%
}
{K=n+1}

\DefEq
{
\def\S{\tau}%
\def\K{m}%
\def\J{}
}
{K=m}

\DefEq
{
\def\J{}%
}
{J=}

\DefEquation
{
\frac{d y}{d x}\circ h=yh
}
{exponent derivative over field}

\DefEq
{
\[
\begin{pmatrix}
y&h_1&...&h_n
\end{pmatrix}
\]
}
{exponent permutation tuple 1}

\DefEq
{
\[
\sigma=
\begin{pmatrix}
y&h_1&...&h_n
\\
\sigma(y)&\sigma(h_1)&...&\sigma(h_n)
\end{pmatrix}
\]
}
{exponent derivative, permutation n}

\DefEq
{
\[
\frac{d^{k-1}y}{d x^{k-1}}\circ (h_1,...,h_{k-1})=\frac 1{2^{k-1}}
\sum_{\sigma\in S(k-1)} \sigma(y)\sigma(h_1) ... \sigma(h_{k-1})
\]
}
{exponent derivative n=k-1}

\DefEq
{
\[
\frac{d^{k-1}y_i}{d x^{k-1}}\circ (h_1,...,h_{k-1})=\frac 1{2^{k-1}}
\sum_{\sigma\in S(k-1)} \sigma(y_j)\sigma(h_1) ... \sigma(h_{k-1})
\]
}
{sh ch derivative n=k-1}

\DefEq
{
\[\sigma(y)\sigma(h_1) ... \sigma(h_n)\]
}
{sigma yh1n}

\DefEq
{
\tau(y)\tau(h_1)...\tau(h_{n+1})=\sigma(h_{n+1}y)\sigma(h_1)...\sigma(h_n)
}
{tau yh1...=sigma yh... 1}

\DefEq
{
\tau(y)\tau(h_1)...\tau(h_{n+1})=\sigma(yh_{n+1})\sigma(h_1)...\sigma(h_n)
}
{tau yh1...=sigma yh... 2}

\DefEq
{
\tau(y)\tau(h_1)...\tau(h_{n+1})
}
{tau yh...}

\DefEq
{
\[
\sigma=
\begin{pmatrix}
h_{n+1}y&h_1&...&h_n
\\
\sigma(h_{n+1}y)&\sigma(h_1)&...&\sigma(h_n)
\end{pmatrix}
\in SE(n)
\]
}
{sigma yh... 1}

\DefEq
{
\[
\sigma=
\begin{pmatrix}
yh_{n+1}&h_1&...&h_n
\\
\sigma(yh_{n+1})&\sigma(h_1)&...&\sigma(h_n)
\end{pmatrix}
\in SE(n)
\]
}
{sigma yh... 2}

\DefEq
{
\[
\begin{pmatrix}
h_{n+1}y&h_1&...&h_n
\end{pmatrix}
\]
}
{sigma yh... 1 top}

\DefEq
{
\[
\begin{pmatrix}
yh_{n+1}&h_1&...&h_n
\end{pmatrix}
\]
}
{sigma yh... 2 top}

\DefEq
{
\[
\begin{pmatrix}
y&h_1&...&h_n&h_{n+1}
\end{pmatrix}
\]
}
{tau yh1... top}

\DefEq
{
\[
\begin{pmatrix}
...&h_{n+1}y&...
\end{pmatrix}
\]
}
{sigma yh... 1 down}

\DefEq
{
\[
\begin{pmatrix}
...&h_{n+1}&y&...
\end{pmatrix}
\]
}
{tau yh1... 1 down}

\DefEq
{
\[
\begin{pmatrix}
...&yh_{n+1}&...
\end{pmatrix}
\]
}
{sigma yh... 2 down}

\DefEq
{
\[
\begin{pmatrix}
...&y&h_{n+1}&...
\end{pmatrix}
\]
}
{tau yh1... 2 down}

\DefEquation
{
\begin{split}
&\,
\sum_{\sigma\in SE(n)}\sigma(h_{n+1}y)\sigma(h_1)...\sigma(h_n)
+
\sum_{\sigma\in SE(n)}\sigma(yh_{n+1})\sigma(h_1)...\sigma(h_n)
\\=&\,
\sum_{\tau\in SE(n+1)}\tau(y)\tau(h_1)...\tau(h_{n+1})
\end{split}
}
{tau yh1...=sigma yh... +}

\DefEq
{
$\tau\in SE(n+1)$
}
{t in Sn+1}

\DefEq
{
$\sigma\in SE(n)$
}
{s in Sn}

\DefEquation
{
\begin{split}
\frac{d^n y_1}{d x^n}\circ(h_1,...,h_n)&=\frac {(-1)^k}{2^n}
\sum_{\sigma\in SE(n)} \sigma(y_1)\sigma(h_1) ... \sigma(h_n)
\\
\frac{d^n y_2}{d x^n}\circ(h_1,...,h_n)&=\frac {(-1)^k}{2^n}
\sum_{\sigma\in SE(n)} \sigma(y_2)\sigma(h_1) ... \sigma(h_n)
\end{split}
}
{sin cos derivative n=2k}

\DefEq
{
\[
m=n+1=2k+1+1=2(l-1)+2=2l
\]
}
{m=2l}

\DefEq
{
\[
m=n+1=2k+1
\]
}
{m=2k+1}

\DefEquation
{
\begin{split}
\frac{d^n y_1}{dx^n}\circ(h_1,...,h_n)&=\frac {(-1)^k}{2^n}
\sum_{\sigma\in SE(n)} \sigma(y_2)\sigma(h_1) ... \sigma(h_n)
\\
\frac{d^n y_2}{d x^n}\circ(h_1,...,h_n)&=-\frac {(-1)^k}{2^n}
\sum_{\sigma\in SE(n)} \sigma(y_1)\sigma(h_1) ... \sigma(h_n)
\end{split}
}
{sin cos derivative n=2k+1}

\DefEquation
{
\frac{d^n y_1}{d x^n}\circ(h_1,...,h_n)=\frac {(-1)^{l-1}}{2^n}
\sum_{\sigma\in SE(n)} \sigma(y_2)\sigma(h_1) ... \sigma(h_n)
}
{sin derivative n=2k+1 k=l-1}

\DefEquation
{
\frac{d^n y_2}{d x^n}\circ(h_1,...,h_n)=-\frac {(-1)^{l-1}}{2^n}
\sum_{\sigma\in SE(n)} \sigma(y_1)\sigma(h_1) ... \sigma(h_n)
}
{cos derivative n=2k+1 k=l-1}

\DefEq
{
\[
\S=
\begin{pmatrix}
y_{\J}&h_1&...&h_{\K}
\\
\S(y_{\J})&\S(h_1)&...&\S(h_{\K})
\end{pmatrix}
\]
}
{derivative, permutation}

\DefEq
{
$y_{\J}$, $h_1$, ..., $h_{\K}$.
}
{derivative, permutation, vars}

\DefEquation
{
\begin{split}
\frac{d^k y}{d x^k}\circ(h_1,...,h_k)
=&\frac{d}{d x}
\left(\frac{d^{k-1}y}{d x^{k-1}}\circ(h_1,...,h_{k-1})\right)\circ h_k
\\
=&\frac 1{2^{k-1}}\frac{d}{d x}
\left(
\sum_{\sigma\in S(k-1)} \sigma(y)\sigma(h_1) ... \sigma(h_{k-1})
\right)
\circ h_k
\end{split}
}
{exponent derivative n=k, 1}

\DefEquation
{
\begin{split}
\frac{d^k y_i}{d x^k}\circ(h_1,...,h_k)
=&\frac{d}{d x}
\left(\frac{d^{k-1}y_i}{d x^{k-1}}\circ(h_1,...,h_{k-1})\right)\circ h_k
\\
=&\frac 1{2^{k-1}}\frac{d}{d x}
\left(
\sum_{\sigma\in S(k-1)} \sigma(y_j)\sigma(h_1) ... \sigma(h_{k-1})
\right)
\circ h_k
\end{split}
}
{sh ch derivative n=k, 1}

\DefEquation
{
\begin{split}
\frac{d^m y_1}{d x^m}\circ(h_1,...,h_m)
=&\frac{d}{d x}
\left(\frac{d^ny_1}{d x^n}\circ(h_1,...,h_n)\right)\circ h_m
\\
=&\frac {(-1)^{l-1}}{2^n}\frac d{d x}
\left(
\sum_{\sigma\in SE(n)} \sigma(y_2)\sigma(h_1) ... \sigma(h_n)
\right)
\circ h_m
\end{split}
}
{sin derivative m=2k+2 k=l-1, 1}

\DefEquation
{
\begin{split}
\frac{d^m y_1}{d x^m}\circ(h_1,...,h_m)
=&\frac{d}{d x}
\left(\frac{d^ny_1}{d x^n}\circ(h_1,...,h_n)\right)\circ h_m
\\
=&\frac {(-1)^k}{2^n}\frac{d}{dx}
\left(
\sum_{\sigma\in SE(n)} \sigma(y_1)\sigma(h_1) ... \sigma(h_n)
\right)
\circ h_m
\end{split}
}
{sin derivative m=2k+1 k=l, 1}

\DefEquation
{
\begin{split}
\frac{d^m y_2}{d x^m}\circ(h_1,...,h_m)
=&\frac{d}{d x}
\left(\frac{d^ny_2}{d x^n}\circ(h_1,...,h_n)\right)\circ h_m
\\
=&\frac {(-1)^k}{2^n}\frac{d}{d x}
\left(
\sum_{\sigma\in SE(n)} \sigma(y_2)\sigma(h_1) ... \sigma(h_n)
\right)
\circ h_m
\end{split}
}
{cos derivative m=2k+1 k=l, 1}

\DefEquation
{
\begin{split}
\frac{d^m y_2}{d x^m}\circ(h_1,...,h_m)
=&\frac{d}{d x}
\left(\frac{d^ny_2}{d x^n}\circ(h_1,...,h_n)\right)\circ h_m
\\
=&-\frac {(-1)^{l-1}}{2^n}\frac{d}{d x}
\left(
\sum_{\sigma\in SE(n)} \sigma(y_1)\sigma(h_1) ... \sigma(h_n)
\right)
\circ h_m
\end{split}
}
{cos derivative m=2k+2 k=l-1, 1}

\DefEquation
{
\begin{split}
\frac{d^ky}{d x^k}\circ(h_1,...,h_k)
=\frac 1{2^{k-1}}\frac 12
&\left(
\sum_{\sigma\in S(k-1)} \sigma(yh_k)\sigma(h_1) ... \sigma(h_{k-1})
\right.
\\
&+
\left.
\sum_{\sigma\in S(k-1)} \sigma(h_ky)\sigma(h_1) ... \sigma(h_{k-1})
\right)
\end{split}
}
{exponent derivative n=k, 2}

\DefEquation
{
\left\{
\begin{aligned}
&k-1=2m\Rightarrow j=i
\\
&k-1=2m-1\Rightarrow j=2-i
\end{aligned}
\right.
}
{j=i,2-i}

\DefEquation
{
\left\{
\begin{aligned}
&k=2m+1\Rightarrow l=2-i
\\
&k=2m\Rightarrow l=i
\end{aligned}
\right.
}
{l=i,2-i}

\DefEquation
{
l=2-j
}
{l=2-j}

\DefEquation
{
\begin{split}
\frac{d^ky_i}{d x^k}\circ(h_1,...,h_k)
=\frac 1{2^{k-1}}\frac 12
&\left(
\sum_{\sigma\in S(k-1)} \sigma(y_lh_k)\sigma(h_1) ... \sigma(h_{k-1})
\right.
\\
+&
\left.
\sum_{\sigma\in S(k-1)} \sigma(h_ky_l)\sigma(h_1) ... \sigma(h_{k-1})
\right)
\end{split}
}
{sh ch derivative n=k, 2}

\DefEquation
{
\begin{array}{r@{\,}l}
\displaystyle\frac{d^my_1}{d x^m}\circ(h_1,...,h_m)
&\displaystyle=-\frac {(-1)^{l-1}}{2^n}\frac 12
\\
\multicolumn{2}{l}{\displaystyle*\left(
\sum_{\sigma\in SE(n)} \sigma(y_1h_m)\sigma(h_1) ... \sigma(h_n)
+
\sum_{\sigma\in SE(n)} \sigma(h_my_1)\sigma(h_1) ... \sigma(h_n)
\right)}
\\
&\displaystyle=\frac {(-1)^l}{2^m}
\\
\multicolumn{2}{l}{\displaystyle*\left(
\sum_{\sigma\in SE(n)} \sigma(y_1h_m)\sigma(h_1) ... \sigma(h_n)
+
\sum_{\sigma\in SE(n)} \sigma(h_my_1)\sigma(h_1) ... \sigma(h_n)
\right)}
\end{array}
}
{sin derivative m=2k+2 k=l-1, 2}

\DefEquation
{
\begin{array}{r@{\,}l}
\displaystyle\frac{d^my_1}{d x^m}\circ(h_1,...,h_m)
&\displaystyle=\frac {(-1)^k}{2^n}\frac 12
\\
\multicolumn{2}{l}{\displaystyle*\left(
\sum_{\sigma\in SE(n)} \sigma(y_2h_m)\sigma(h_1) ... \sigma(h_n)
+
\sum_{\sigma\in SE(n)} \sigma(h_my_2)\sigma(h_1) ... \sigma(h_n)
\right)}
\\
&\displaystyle=\frac {(-1)^k}{2^m}
\\
\multicolumn{2}{l}{\displaystyle*\left(
\sum_{\sigma\in SE(n)} \sigma(y_2h_m)\sigma(h_1) ... \sigma(h_n)
+
\sum_{\sigma\in SE(n)} \sigma(h_my_2)\sigma(h_1) ... \sigma(h_n)
\right)}
\end{array}
}
{sin derivative m=2k+1 k=l, 2}

\DefEquation
{
\begin{array}{r@{\,}l}
\displaystyle\frac{d^my_2}{d x^m}\circ(h_1,...,h_m)
&\displaystyle=-\frac {(-1)^k}{2^n}\frac 12
\\
\multicolumn{2}{l}{\displaystyle*\left(
\sum_{\sigma\in SE(n)} \sigma(y_1h_m)\sigma(h_1) ... \sigma(h_n)
+
\sum_{\sigma\in SE(n)} \sigma(h_my_1)\sigma(h_1) ... \sigma(h_n)
\right)}
\\
&\displaystyle=-\frac {(-1)^k}{2^m}
\\
\multicolumn{2}{l}{\displaystyle*\left(
\sum_{\sigma\in SE(n)} \sigma(y_1h_m)\sigma(h_1) ... \sigma(h_n)
+
\sum_{\sigma\in SE(n)} \sigma(h_my_1)\sigma(h_1) ... \sigma(h_n)
\right)}
\end{array}
}
{cos derivative m=2k+1 k=l, 2}

\DefEquation
{
\begin{array}{r@{\,}l}
\displaystyle\frac{d^my_2}{d x^m}\circ(h_1,...,h_m)
&\displaystyle=-\frac {(-1)^{l-1}}{2^n}\frac 12
\\
\multicolumn{2}{l}{\displaystyle*\left(
\sum_{\sigma\in SE(n)} \sigma(y_2h_m)\sigma(h_1) ... \sigma(h_n)
+
\sum_{\sigma\in SE(n)} \sigma(h_my_2)\sigma(h_1) ... \sigma(h_n)
\right)}
\\
=&\displaystyle\frac {(-1)^l}{2^m}
\\
\multicolumn{2}{l}{\displaystyle*\left(
\sum_{\sigma\in SE(n)} \sigma(y_2h_m)\sigma(h_1) ... \sigma(h_n)
+
\sum_{\sigma\in SE(n)} \sigma(h_my_2)\sigma(h_1) ... \sigma(h_n)
\right)}
\end{array}
}
{cos derivative m=2k+2 k=l-1, 2}

\DefEq
{
\begin{align*}
\frac{d^my_1}{d x^m}\circ(h_1,...,h_m)
=&\frac {(-1)^l}{2^m}
\sum_{\tau\in S(m)} \tau(y_1)\tau(h_1) ... \tau(h_m)
\\
\frac{d^my_2}{d x^m}\circ(h_1,...,h_m)
=&\frac {(-1)^l}{2^m}
\sum_{\tau\in S(m)} \tau(y_2h_m)\tau(h_1) ... \tau(h_m)
\end{align*}
}
{sin cos derivative m=2k+2 k=l-1, 2}

\DefEq
{
\symb{SE(n)}{set of permutations}{1}
}
{set of permutations}

\DefEq
{
\symb{e^x}{exponent}{}
$y=\ShowSymbol{exponent}{}$
}
{exponent over division ring}

\DefEq
{
\symb{\sinh x}{hyperbolic sine}{}
$y_1=\ShowSymbol{hyperbolic sine}{}$
}
{hyperbolic sine}

\DefEq
{
\symb{\sin x}{sine}{}
$y_1=\ShowSymbol{sine}{}$
}
{sine}

\DefEq
{
\symb{\cosh x}{hyperbolic cosine}{}
$y_2=\ShowSymbol{hyperbolic cosine}{}$
}
{hyperbolic cosine}

\DefEq
{
\symb{\cos x}{cosine}{}
$y_2=\ShowSymbol{cosine}{}$
}
{cosine}

\DefEquation
{
\ShowSymbol{exponent}{}=\sum_{n=0}^{\infty}
\frac 1{n!}x^n
}
{exponent Taylor series}

\DefEquation
{
\ShowSymbol{hyperbolic sine}{}=\sum_{n=0}^{\infty}
\frac 1{(2n+1)!}x^{2n+1}
}
{hyperbolic sine Taylor series}

\DefEquation
{
\ShowSymbol{sine}{}=\sum_{n=0}^{\infty}
\frac {(-1)^n}{(2n+1)!}x^{2n+1}
}
{sine Taylor series}

\DefEquation
{
\ShowSymbol{hyperbolic cosine}{}=\sum_{n=0}^{\infty}
\frac 1{(2n)!}x^{2n}
}
{hyperbolic cosine Taylor series}

\DefEquation
{
\ShowSymbol{cosine}{}=\sum_{n=0}^{\infty}
\frac {(-1)^n}{(2n)!}x^{2n}
}
{cosine Taylor series}

\DefEquation
{
\left.\frac{d^ny}{d x^n}\right|_{x=0}\circ(h,...,h)=h^n
}
{exponent derivative n, x=0}

\DefEquation
{
\begin{matrix}
n=2k
&\displaystyle
\left.\frac{d^n y_1}{d x^n}\right|_{x=0}\circ(h_1,...,h_n)=0
&\displaystyle
\left.\frac{d^n y_2}{d x^n}\right|_{x=0}\circ(h_1,...,h_n)=h^n
\\
n=2k+1
&\displaystyle
\left.\frac{d^n y_1}{d x^n}\right|_{x=0}\circ(h_1,...,h_n)=h^n
&\displaystyle
\left.\frac{d^n y_2}{d x^n}\right|_{x=0}\circ(h_1,...,h_n)=0
\end{matrix}
}
{sh ch derivative n, x=0}

\DefEquation
{
\begin{matrix}
\displaystyle
\left.\frac{d^n y_1}{d x^n}\right|_{x=0}\circ(h_1,...,h_n)
=(-1)^kh^n
&\displaystyle
\left.\frac{d^n y_2}{d x^n}\right|_{x=0}\circ(h_1,...,h_n)=0
\end{matrix}
}
{sin cos derivative n=2k+1, x=0}

\DefEquation
{
\begin{matrix}
\displaystyle
\left.\frac{d^n y_1}{d x^n}\right|_{x=0}\circ(h_1,...,h_n)=0
&\displaystyle
\left.\frac{d^n y_2}{d x^n}\right|_{x=0}\circ(h_1,...,h_n)
=(-1)^kh^n
\end{matrix}
}
{sin cos derivative n=2k, x=0}

\DefEquation
{
ab=ba
}
{exponent of sum, 1}

\DefEq
{
\begin{align}
e^a&=\sum_{n=0}^{\infty}\frac 1{n!}a^n
\EqLabel{exponent a}
\\
e^b&=\sum_{n=0}^{\infty}\frac 1{n!}b^n
\EqLabel{exponent b}
\\
e^{a+b}&=\sum_{n=0}^{\infty}\frac 1{n!}(a+b)^n
\EqLabel{exponent a+b}
\end{align}
}
{exponent a b a+b}

\DefEquation
{
\frac 16(a+b)^3=\frac 16a^3+\frac 16a^2b+\frac 16aba+\frac 16ba^2
+\frac 16ab^2+\frac 16bab+\frac 16b^2a+\frac 16b^3
}
{exponent a+b 3}

\DefEquation
{
\frac 16a^3+\frac 12a^2b+\frac 12ab^2+\frac 16b^3
}
{exponent ab 3}

\DefEquation
{
e^{a+b}=e^ae^b
}
{exponent of sum}

\DefEquation
{
\frac{d^ny}{d x^n}\circ(h_1,...,h_n)=\frac 1{2^n}
\sum_{\sigma\in SE(n)} \sigma(y)\sigma(h_1) ... \sigma(h_n)
}
{exponent derivative n over division ring}

\DefEquation
{
\begin{matrix}
n=2k
&\displaystyle
\frac{d^n y_1}{d x^n}\circ(h_1,...,h_n)=\frac 1{2^n}
\sum_{\sigma\in SE(n)} \sigma(y_1)\sigma(h_1) ... \sigma(h_n)
\\
&\displaystyle
\frac{d^n y_2}{d x^n}\circ(h_1,...,h_n)=\frac 1{2^n}
\sum_{\sigma\in SE(n)} \sigma(y_2)\sigma(h_1) ... \sigma(h_n)
\\
n=2k+1
&\displaystyle
\frac{d^n y_1}{d x^n}\circ(h_1,...,h_n)=\frac 1{2^n}
\sum_{\sigma\in SE(n)} \sigma(y_2)\sigma(h_1) ... \sigma(h_n)
\\
&\displaystyle
\frac{d^n y_2}{d x^n}\circ(h_1,...,h_n)=\frac 1{2^n}
\sum_{\sigma\in SE(n)} \sigma(y_1)\sigma(h_1) ... \sigma(h_n)
\end{matrix}
}
{sh ch derivative n}

\DefEquation
{
\frac{d y}{d x}\circ h=\frac 12(yh+hy)
}
{exponent derivative over division ring}

\DefEquation
{
y'=y
}
{exponent over field}

\DefEq
{
\ePrints{}
\ifx\Semafor\ValueOn
\RefTheorem[1601.03259]{derivative of pn is symmetric, m < n, algebra},
\else
\RefTheorem{derivative of pn is symmetric, m < n, algebra},
\fi
}
{derivative of pn is symmetric, m < n}

\DefEq
{
\ePrints{}%
\ifx\Semafor\ValueOn%
\EqRef[1601.03259]{derivative of Order n, algebra},
\else%
\EqRef{derivative of Order n, algebra},
\fi%
}
{ref derivative of Order n, algebra}

\DefEq
{
\[
y_1=\sinh x
\]
}
{y1=sh}

\DefEq
{
\[
y_1=\sin x
\]
}
{y1=sin}

\DefEquation
{
\begin{split}
&\frac{d^2y_1}{dx^2}=y_1
\\
&x=0\ \ \ y_1=0\ \ \ \frac{dy_1}{dx}=1
\end{split}
}
{d2 sh}

\DefEquation
{
\begin{split}
&\frac{d^2y_1}{dx^2}=-y_1
\\
&x=0\ \ \ y_1=0\ \ \ \frac{dy_1}{dx}=1
\end{split}
}
{d2 sin}

\DefEquation
{
\begin{split}
&\frac{d^2y_2}{dx^2}=y_2
\\
&x=0\ \ \ y_2=1\ \ \ \frac{dy_2}{dx}=0
\end{split}
}
{d2 ch}

\DefEquation
{
\begin{split}
&\frac{d^2y_2}{dx^2}=-y_2
\\
&x=0\ \ \ y_2=1\ \ \ \frac{dy_2}{dx}=0
\end{split}
}
{d2 cos}

\DefEquation
{
\begin{split}
\frac{d y_1}{d x}&=\frac{1}{2}(y_2\otimes 1+1\otimes y_2)
\\
\frac{d y_2}{d x}&=\frac{1}{2}(y_1\otimes 1+1\otimes y_1)
\end{split}
}
{d sh ch A}

\DefEquation
{
\begin{split}
\frac{d y_1}{d x}&=\frac{1}{2}(y_2\otimes 1+1\otimes y_2)
\\
\frac{d y_2}{d x}&=-\frac{1}{2}(y_1\otimes 1+1\otimes y_1)
\end{split}
}
{d sin cos A}

\DefEquation
{
\frac{\partial^2 y}{\partial x^2}
-\frac 14(y\otimes 1\otimes 1+2\ 1\otimes y\otimes 1+\otimes 1\otimes 1\otimes y)=0
}
{d2 sh A}

\DefEquation
{
\frac{d^2 y}{d x^2}
+\frac 14(y\otimes 1\otimes 1+2\ 1\otimes y\otimes 1+\otimes 1\otimes 1\otimes y)=0
}
{d2 sin A}

\DefEq
{
\frac{d^2 y_1}{d x^2}
=\frac{1}{2}
\left(\frac{d y_2}{d x}\otimes 1+1\otimes\frac{d y_2}{d x}\right)
}
{d2 sh A 1}

\DefEquation
{
\frac{d^2 y_1}{d x^2}
=\frac{1}{2}
\left(\frac{1}{2}(y_1\otimes 1+1\otimes y_1)\otimes 1
+1\otimes\frac{1}{2}(y_1\otimes 1+1\otimes y_1)\right)
}
{d2 sh A 2}

\DefEquation
{
\frac{d^2 y_1}{d x^2}
=-\frac{1}{2}
\left(\frac{1}{2}(y_1\otimes 1+1\otimes y_1)\otimes 1
+1\otimes\frac{1}{2}(y_1\otimes 1+1\otimes y_1)\right)
}
{d2 sin A 2}

\DefEquation
{
\frac{d y_2}{d x}=\frac{1}{2}(y_1\otimes 1+1\otimes y_1)
}
{d ch A}

\DefEq
{
\frac{d y_1}{d x}=\frac{1}{2}(y_2\otimes 1+1\otimes y_2)
}
{d sh A}

\DefEquation
{
\frac{d y_2}{d x}=-\frac{1}{2}(y_1\otimes 1+1\otimes y_1)
}
{d cos A}

\DefEq
{
x=0\ \ \ y_1=0\ \ \ y_2=1
}
{d sh ch init}

\DefEq
{
\[
y_2=\cosh x
\]
}
{y2=ch}

\DefEq
{
\[
y_2=\cos x
\]
}
{y2=cos}

\DefEquation
{
\begin{split}
\frac{dy_1}{dx}&=y_2
\\
\frac{dy_2}{dx}&=y_1
\end{split}
}
{d sh ch R}

\DefEquation
{
\begin{split}
\frac{dy_1}{dx}&=y_2
\\
\frac{dy_2}{dx}&=-y_1
\end{split}
}
{d sin cos R}

%% file: Two.Integrals.English.tex

\input{Two.Integrals.Eq}

\ePrints{1702.01}
\ifx\Semafor\ValueOff
\Chapter{Lebesgue Integral}

I decided to consider theorems
\RefTheorem{Lebesgue Integral along Linear Path},
\RefTheorem{Lebesgue Integral along Path}
\ifx\setCACAA\undefined
in isolated chapter,
\else
in isolated section,
\fi
because I believe that these theorems are very important.
It is important to note that theorems
\RefTheorem{Lebesgue Integral along Linear Path},
\RefTheorem{Lebesgue Integral along Path}
provide an alternative method for solving differential equation
\ShowEq{df=g}

\Section{Lebesgue Integral along Linear Path}
\fi

\begin{theorem}
\labelTheorem{Lebesgue Integral along Linear Path}
Let $A$ be Banach $D$\Hyph module.
Let $B$ be Banach $D$\Hyph algebra.
Let
\ShowEq{g:A->BoxB}
be integrable map
\DrawEq{f=int g.}{linear}
and
$\|g\|=G<\infty$.
For any $A$\Hyph numbers $a$, $x$,
let
\ShowEq{h:[01]->A}
be linear path in $D$\Hyph algebra $A$
\ShowEq{h=a+tb}
Then
\DrawEq{int hn= 2}{linear path}
\end{theorem}
\begin{proof}
From the equality
\EqRef{h=a+tb},
it follows that
\ShowEq{dx=dt b-a}
\ShowEq{h'=b-a}
and the equality
\ShowEq{g(x)=g(h(t))}
follows from the equality
\EqRef{dx=dt b-a}.
To calculate integral
\ShowEq{int h=int 01}
we consider the partition
\DrawEq{t0...tn}{}
of segment $[0,1]$
into equal segments.
In such case, we can use simple map\,\footnote{
According to the definition
\ShowEq{ref: simple map}
and the theorem
\ShowEq{ref: measurable simple map}
a map is simple, if
range is finite or countable set
and inverse image of each value is measurable set.
}%
\ShowEq{gn(t)=}
instead of integrand.
From
\EqRef{g(x)=g(h(t))}, \EqRef{gn(t)=},
it follows that
\ShowEq{|g-gn|}
when
\ShowEq{ti<t<ti+1}
From
\EqRef{|g-gn|},
it follows that
\ShowEq{g(h)=lim gn}
and therefore
\ShowEq{g(b-a)=lim...}
The equality
\ShowEq{g(x)dx=}
follows from the equality
\eqRef{f=int g.}{linear}
and definitions
\RefDefinition{differentiable map, algebra},
\RefDefinition{indefinite integral}.
Here
$o:A\rightarrow A$
is the map such that
\ShowEq{lim |o(a)|=0}
From
\EqRef{gn(t)=}, \EqRef{g(x)dx=},
it follows that
\ShowEq{dt gn (b-a)}
where
\ShowEq{lim o1/n}
From
\EqRef{dt gn (b-a)}
it follows that
the map
$g_n(t)\circ(x-a)$
is integrable and
\ShowEq{int gn=...}
From
\EqRef{lim o1/n},
it follows that
\ShowEq{lim no=0}
From
\EqRef{g(b-a)=lim...},
\EqRef{lim no=0}
and the definition
\EqParm{ref: Integral of Map over Set of Finite Measure}{=c}
it follows that
the map
\ShowEq{g(h)o(b-a)}
is integrable on the segment $[0,1]$ and
the equality
\eqRef{int hn= 2}{linear path}
follows from
\eqRef{f=int g.}{linear}.
\end{proof}

\ePrints{1702.01}
\ifx\Semafor\ValueOff
\Section{Lebesgue Integral along Path}
\fi

\begin{theorem}
\labelTheorem{Lebesgue Integral along Path}
Let there exists indefinite integral
\DrawEq{f=int g.}{}
For any rectifiable continuous path
\ShowEq{h:[01]->A}
from $a$ to $x$ in $D$\Hyph module $A$
\DrawEq{int hn= 2}{path}
\end{theorem}
\begin{proof}
Let
\DrawEq{t0...tn}{path}
be set of points of the interval $[0,1]$
such that
\ShowEq{lim max dt = 0}
We define the path $\gamma_n$ by the equality
\ShowEq{hn=htiti1}
Sequence of maps $\gamma_n$
converges to map $h$.
We may request uniform convergence.
In theory of integral, it follows that
\ShowEq{int h=lim int hn}
According to the partition
\eqRef{t0...tn}{path}
we can represent integral
\ShowEq{int hn}
as sum
\ShowEq{int hn=sum}
For each $i$, we change variable
\ShowEq{s=...ti}
Then
\ShowEq{int hin=}
According to the theorem
\RefTheorem{Lebesgue Integral along Linear Path},
\ShowEq{int hin= 1}
From
\EqRef{int hn=sum},
\EqRef{int hin= 1},
it follows that
\ShowEq{int hn=}
\end{proof}

\Section{Solving of Differential Equation}

We will consider one time more the differential equation
\DrawEq{differential equation y=xx, 1, algebra}{path}
and initial condition
\DrawEq{differential equation, initial}{y=xx, path}
According to the definition]
\RefDefinition{indefinite integral},
the map $y$ is integral
\ShowEq{int=x3 path}
So we can apply the theorem
\RefTheorem{Lebesgue Integral along Linear Path}
to solve differential equation
\eqRef{differential equation, initial}{y=xx, path}.
We consider the linear path
\ShowEq{h:[01]->A}
in $D$\Hyph algebra $A$
\ShowEq{h=tx}
Then the integral
\eqRef{int hn= 2}{linear path}
gets the form
\ShowEq{int=x3 path tx}


How easily we solved differential equation.
However, is this procedure very easy?

To get answer on this question, we also consider differential equation
\DrawEq{differential equation y=xx, 1, a, algebra}{path}
\DrawEq{differential equation, initial}{}
If we consider the path
\EqRef{h=tx},
we will get that
\ShowEq{int 3x2=x3}
However we know that this is not true.

To solve differential equation,
we considered linear path.
However a path can be arbitrary.
Consider the path
\ShowEq{path 0ax}
We can represent the integral
\ShowEq{int h gx dx}
as sum
\ShowEq{int h gx dx sum}
where $\gamma_1$ is the path
\ShowEq{path h1}
and $\gamma_2$ is the path
\ShowEq{path h2}
We considered above the integral
\ShowEq{int h1=a3}
So we need to consider the integral
\ShowEq{int h2=int 12}
\begin{itemize}
\item
Let
\ShowEq{gx=xx+}
Then the integral
\EqRef{int h2=int 12}
gets the form
\ShowEq{int h2=int 12 xx+}
The equality
\ShowEq{int h2 xx+ =}
\ePrints{1702.01}
\ifx\Semafor\ValueOff
follows from equalities
\EqRef{int h2=int 12 xx+},
\EqRef{int h2=int 12 xx+ 5}.
\else
follows from the equalitiy
\EqRef{int h2=int 12 xx+}.
\fi
The equality
\ShowEq{int h xx+ =}
follows from equalities
\EqRef{int h gx dx sum},
\EqRef{int h1=a3},
\EqRef{int h2 xx+ =}.
\item
Let
\ShowEq{gx=3x2}
Then the integral
\EqRef{int h2=int 12}
gets the form
\ShowEq{int h2=int 12 3x2}
The equality
\ShowEq{int h2 3x2 =}
\ePrints{1702.01}
\ifx\Semafor\ValueOff
follows from equalities
\EqRef{int h2=int 12 3x2},
\EqRef{int h2=int 12 3x2 3}.
\else
follows from the equalitiy
\EqRef{int h2=int 12 3x2}.
\fi
The equality
\ShowEq{int h 3x2 =}
follows from equalities
\EqRef{int h gx dx sum},
\EqRef{int h1=a3},
\EqRef{int h2 3x2 =}.
\end{itemize}

So the integral
\ShowEq{int h 3x2 dx}
depends on choice of a path $\gamma$ from $0$ to $x$;
and the differential equation
\eqRef{differential equation y=xx, 1, a, algebra}{path}
does not have solution.
\ePrints{1702.01}
\ifx\Semafor\ValueOff
So, before we can apply Lebesgue integral
to solving of differential equation,
we need to answer one of two following questions.
How we can prove that integral
does not depend from choice of path or how we can find the path
which proves dependence of integral on path.

We suddenly found ourselves in new territory.
The name of this territory is differential forms.
\fi

%% file: Two.Integrals.Eq.tex

\input{Two.Integrals.Ref}

\DefEq
{
f(x)=\int g(x)\circ dx
}
{f=int g.}

\DefEquation
{
y=\gamma(t)=a+t(x-a)
}
{h=a+tb}

\DefEquation
{
\gamma(t)=tx
}
{h=tx}

\DefEquation
{
g_n(t):[0,1]\rightarrow A\otimes A\ \ \ t_{i-1}< t\le t_i=>g_n(t)=g(\gamma(t_i))
}
{gn(t)=}

\DefEquation
{
\begin{split}
&\,\sum_{k=0}^n
\left(1\otimes \frac {k^2}{n^2}x^2+\frac knx\otimes \frac knx
+\frac {k^2}{n^2}x^2\otimes 1\right)\circ\left(\frac 1nx\right)
\\=&
\,((1\otimes x^2+x\otimes x+x^2\otimes 1)\circ x)
\sum_{k=0}^n\frac {k^2}{n^3}
\\=&
\,3x^3
\frac {n(n+1)(2n+1)}{6n^3}
\end{split}
}
{sum x3 1}

\DefEq
{
\[
\gamma(t)=
\left\{\begin{matrix}
ta&0\le t\le 1
\\
a+(t-1)(x-a)&1\le t\le 2
\end{matrix}
\right.
\]
}
{path 0ax}

\DefEq
{
\[\int (3\otimes x^2)\circ dx=x^3+C\]
}
{int 3x2=x3}

\DefEq
{
\[g(\gamma(t))=\lim_{n\rightarrow \infty}g_n(t)\]
}
{g(h)=lim gn}

\DefEquation
{
g(\gamma(t))\circ(x-a)=\lim_{n\rightarrow \infty}g_n(t)\circ(x-a)
}
{g(b-a)=lim...}

\DefEquation
{
g(y)\circ dy=\frac{d f(y)}{dy}\circ dy=f(y+dy)-f(y)-o(dy)
}
{g(x)dx=}

\DefEq
{
\[\lim_{a\rightarrow 0}\frac{\|o(a)\|}{\|a\|}=0\]
}
{lim |o(a)|=0}

\DefEquation
{
\int_0^1dt (g(a+t(b-a))\circ(b-a))=f(b)-f(a)
}
{int gh=...}

\DefEq
{
$g(\gamma(t))\circ(x-a)$
}
{g(h)o(b-a)}

\DefEq
{
\[\int(1\otimes x^2+x\otimes x+x^2\otimes 1)\circ dx\]
}
{int=x3 path}

\DefEquation
{
\begin{split}
&\,\int_{\gamma}(1\otimes x^2+x\otimes x+x^2\otimes 1)\circ dx
\\=&\,\int_0^1dt((1\otimes (tx)^2+(tx)\otimes (tx)+(tx)^2\otimes 1)\circ x)
\\=&\,\left(\int_0^1dt\,t^2\right)((1\otimes x^2+x\otimes x+x^2\otimes 1)\circ x)
\\=&\,\frac 13(x^3+x^3+x^3)
\\=&\,x^3
\end{split}
}
{int=x3 path tx}

\DefEquation
{
\int_0^1dt (g_n(t)\circ(x-a))=f(x)-f(a)-n\,o\left(\frac 1n(x-a)\right)
}
{int gn=...}

\DefEquation
{
\lim_{n\rightarrow\infty}n\,o\left(\frac 1n(x-a)\right)=0
}
{lim no=0}

\DefEq
{
\int_{\gamma}g(y)\circ dy=
\int_0^1 dt\left(g(\gamma(t))\circ\frac{d\gamma(t)}{dt}\right)=f(x)-f(a)
}
{int hn= 2}

\DefEq
{
0=t_0<t_1<...<t_{n-1}<t_n=1
}
{t0...tn}

\DefEq
{
\[
\lim_{n\rightarrow \infty}\max(t_{i+1}-t_i)=0
\]
}
{lim max dt = 0}

\DefEq
{
$[t_i,t_{i+1}]$,
}
{[titi1]}

\DefEquation
{
\int_{\gamma} dt\left(g(\gamma(t))\circ\frac{d\gamma(t)}{dt}\right)=\lim_{n\rightarrow \infty}
\int_{\gamma_n} dt\left(g(\gamma_n(t))\circ\frac{d\gamma_n(t)}{dt}\right)
}
{int h=lim int hn}

\DefEquation
{
\int_{\gamma_n} dt\left(g(\gamma_n(t))\circ\frac{d\gamma_n(t)}{dt}\right)=\sum_{i=1}^n
\int_{t_{i-1}}^{t_i} dt\left(g(\gamma_n(t))\circ\frac{d\gamma_n(t)}{dt}\right)
}
{int hn=sum}

\DefEq
{
$\displaystyle\int_{\gamma_n} dt\left(g(\gamma_n(t))\circ\frac{d\gamma_n(t)}{dt}\right)$
}
{int hn}

\DefEq
{
$\displaystyle\int_{\gamma}g(x)\circ dx$
}
{int h gx dx}

\DefEquation
{
\int_{\gamma_1}g(x)\circ dx=a^3
}
{int h1=a3}

\DefEq
{
\[
\int_{\gamma_1}g(x)\circ dx=x^3
\]
}
{int h xx+ =}

\DefEq
{
\[
\int_{\gamma_1}g(x)\circ dx=x^3+\frac 12x^2a+\frac 12xax
-ax^2+xa^2
-\frac 12axa
-\frac 12a^2x
\]
}
{int h 3x2 =}

\DefEq
{
$\int_{\gamma}(3\otimes x^2)\circ dx$
}
{int h 3x2 dx}

\DefEq
{
\[
\gamma_1:t\in[0,1]\subset R\rightarrow ta\in A
\]
}
{path h1}

\DefEq
{
\[
\gamma_2:t\in[1,2]\subset R\rightarrow a+(t-1)(x-a)\in A
\]
}
{path h2}

\DefEq
{
$g(x)=3\otimes x^2$.
}
{gx=3x2}

\DefEq
{
$g(x)=x^2\otimes 1+x\otimes x +1\otimes x^2$.
}
{gx=xx+}

\DefEquation
{
\begin{split}
\int_{\gamma_2}g(x)\circ dx&=\int_1^2dt\,g(a+(t-1)(x-a))\circ\frac{d(a+(t-1)(x-a))}{dt}
\\&=\int_1^2dt\,g(a+(t-1)(x-a))\circ(x-a)
\\&=\int_0^1dt\,g(a+t(x-a))\circ(x-a)
\end{split}
}
{int h2=int 12}

\DefEquation
{
\begin{split}
&\,\int_{\gamma_2}g(x)\circ dx
\\=&\,\int_0^1dt\,((a+t(x-a))^2\otimes 1+(a+t(x-a))\otimes(a+t(x-a))
\\+&\,1\otimes(a+t(x-a))^2)
\circ(x-a)
\\=&\,\int_0^1dt\,((a+t(x-a))^2(x-a)+(a+t(x-a))(x-a)(a+t(x-a))
\\+&\,(x-a)(a+t(x-a))^2)
\end{split}
}
{int h2=int 12 xx+}

\DefEquation
{
\begin{split}
&\,\int_{\gamma_2}g(x)\circ dx
\\=&\,xa^2+axa+a^2x-3a^3
\\+&\,\frac 12(2x^2a+2xax+2ax^2-4xa^2-4axa-4a^2x+6a^3)
\\+&\,\frac 13(3x^3-3x^2a-3xax+3xa^2-3ax^2+3axa+3a^2x-3a^3)
\\=&\,xa^2+axa+a^2x-3a^3
\\+&\,x^2a+xax+ax^2-2xa^2-2axa-2a^2x+3a^3
\\+&\,x^3-x^2a-xax+xa^2-ax^2+axa+a^2x-a^3
\\=&\,x^3-a^3
\end{split}
}
{int h2 xx+ =}

\DefEquation
{
\begin{split}
\int_{\gamma_2}g(x)\circ dx
&=\int_0^1dt\,3\otimes(a+t(x-a))^2\circ(x-a)
\\&=3\int_0^1dt\,(x-a)(a+t(x-a))^2
\end{split}
}
{int h2=int 12 3x2}

\DefEquation
{
\begin{split}
\int_{\gamma_2}g(x)\circ dx
&=3\left(xa^2-a^3+\frac 12(x^2a+xax-2xa^2-a^2x-axa+2a^3)\right.
\\&\left.+\frac 13(x^3-x^2a-xax-ax^2+xa^2+axa+a^2x-a^3)\right)
\\&=x^3+\left(\frac 32-1\right)x^2a+\left(\frac 32-1\right)xax
-ax^2+xa^2
\\&+\left(1-\frac 32\right)axa
+\left(1-\frac 32\right)a^2x-a^3
\\&=x^3+\frac 12x^2a+\frac 12xax
-ax^2+xa^2
-\frac 12axa
-\frac 12a^2x-a^3
\end{split}
}
{int h2 3x2 =}

\DefEquation
{
\int_{\gamma}g(x)\circ dx=
\int_{\gamma_1}g(x)\circ dx+\int_{\gamma_2}g(x)\circ dx
}
{int h gx dx sum}

\DefEquation
{
s=\frac{t-t_{i-1}}{t_i-t_{i-1}}\ \ \ \ \gamma_{in}(s)=\gamma_n(t)\ \ \ \ t_{i-1}\le t\le t_i
}
{s=...ti}

\DefEquation
{
\begin{split}
&\,\int_{t_{i-1}}^{t_i} dt\left(g(\gamma_n(t))\circ\frac{d\gamma_n(t)}{dt}\right)
\\=&\,
\int_{t_{i-1}}^{t_i} dt
\\ * &\,
\left(g\left(\gamma_n(t_{i-1})
+(t-t_{i-1})\frac{\gamma_n(t_i)-\gamma_n(t_{i-1})}{t_i-t_{i-1}}\right)
\circ \frac{\gamma_n(t_i)-\gamma_n(t_{i-1})}{t_i-t_{i-1}}\right)
\\=&\,
\int_0^1 (t_i-t_{i-1})ds\left(g\left(\gamma_{in}(0)+s(\gamma_{in}(1)-\gamma_{in}(0))\right)
\circ \frac{\gamma_{in}(1)-\gamma_{in}(0)}{t_i-t_{i-1}}\right)
\\=&\,
\int_0^1 ds(g(\gamma_{in}(0)+s(\gamma_{in}(1)-\gamma_{in}(0)))
\circ (\gamma_{in}(1)-\gamma_{in}(0)))
\end{split}
}
{int hin=}

\DefEquation
{
\begin{split}
\int_{t_{i-1}}^{t_i} dt\left(g(\gamma_n(t))\circ\frac{d\gamma_n(t)}{dt}\right)
&=f(\gamma_{in}(1))-f(\gamma_{in}(0))\\&=f(\gamma_n(t_i))-f(\gamma_n(t_{i-1}))
\end{split}
}
{int hin= 1}

\DefEquation
{
dt (g_n(t_i)\circ(x-a))=f(\gamma(t_i))-f(\gamma(t_{i-1}))-o\left(\frac 1n(x-a)\right)
}
{dt gn (b-a)}

\DefEquation
{
\lim_{n\rightarrow\infty}\frac{\|o\left(\frac 1n(x-a)\right)\|}{\|\frac 1n(x-a)\|}=0
}
{lim o1/n}

\DefEquation
{
\begin{split}
\int_{\gamma_n} dt\left(g(\gamma_n(t))\circ\frac{d\gamma_n(t)}{dt}\right)&=\sum_{i=1}^n
(f(\gamma_n(t_i))-f(\gamma_n(t_{i-1})))
\\&=f(\gamma_n(1))-f(\gamma_n(0))=f(\gamma(1))-f(\gamma(0))
\end{split}
}
{int hn=}

\DefEq
{
\[
\gamma_{in}:[t_i,t_{i+1}]\rightarrow A
\ \ \  \gamma_{in}(t)=\gamma(t_i)+(t-t_i)\frac{\gamma(t_{i+1})-\gamma(t_i)}{t_{i+1}-t_i}
\]
}
{hi:01->A}

\DefEq
{
\[
\gamma_{n}(t)=\gamma(t_i)+(t-t_i)\frac{\gamma(t_{i+1})-\gamma(t_i)}{t_{i+1}-t_i}
\ \ \ t_i\le t\le t_{i+1}
\]
}
{hn=htiti1}

\DefEquation
{
\|g(\gamma(t))-g_n(t)\|<\|g\|\|\gamma(t)-\gamma(t_i)\|\le\frac{G\|x-a\|}n
}
{|g-gn|}

\DefEq
{
$t_i\le t<t_{i+1}$.
}
{ti<t<ti+1}

\DefEquation
{
\int_{\gamma}g(y)\circ dy=\int_0^1dt\left(g(\gamma(t))\circ \frac{d\gamma(t)}{dt}\right)
}
{g(x)=g(h(t))}

\DefEquation
{
dy=d\gamma(t)=dt\,\frac{d\gamma(t)}{dt}
}
{dx=dt b-a}

\DefEquation
{
\frac{d\gamma(t)}{dt}=x-a
}
{h'=b-a}

\DefEq
{
\[\gamma:[0,1]\subset R\rightarrow A\]
}
{h:[01]->A}

\DefEquation
{
\int_{\gamma}g(y)\circ dy=\int_0^1dt (g(a+t(x-a))\circ(x-a))
}
{int h=int 01}

%% file: Two.Integrals.Ref.tex

\DefEq
{
\ifx\texFuture\Defined
\ePrints{1601.03259,1610.309618526,1601.03259}
\ifx\Semafor\ValueOn
\RefDefinition[1310.5591]{simple map}
\else
\ePrints{4975-6381}
\ifx\Semafor\ValueOn
\RefDefinition[5410-9916]{simple map}
\else
\RefDefinition[CACAA.04.001]{simple map}
\fi
\fi
\fi
\RefDefinition{simple map}
}
{ref: simple map}

\DefEq
{
\ifx\texFuture\Defined
\ePrints{1601.03259,1610.309618526}
\ifx\Semafor\ValueOn
\RefTheorem[1310.5591]{measurable simple map},
\else
\ePrints{4975-6381}
\ifx\Semafor\ValueOn
\RefTheorem[5410-9916]{measurable simple map},
\else
\RefTheorem[CACAA.04.001]{measurable simple map},
\fi
\fi
\fi
\RefTheorem{measurable simple map}
}
{ref: measurable simple map}

\DefEq
{
\ePrints{1601.03259,1610.309618526}
\ifx\Semafor\ValueOn
\RefDefinition[1310.5591]{Integral of Map over Set of Finite Measure}\Pt
\else
\RefDefinition[CACAA.04.001]{Integral of Map over Set of Finite Measure}\Pt
\fi
}
{ref: Integral of Map over Set of Finite Measure}

%% file: Diff.Form.1.English.tex

\input{Diff.Form.1.Eq}

\Chapter{Differential Form}

\Section{Structure of Polylinear Map}

\begin{definition}
\labelDefinition{symmetric polylinear map}
{\it
Let $A$, $B$ be algebras over commutative ring $D$.
A polylinear map
\ShowEq{f in L(A->B)}{A^n}B{}
is called
\AddIndex{symmetric}{symmetric polylinear map},
if
\ShowEq{fa=fsa}
for any permutation $\sigma$ of the set
\ShowEq{set a1n}
}
\qed
\end{definition}

\begin{theorem}
Let
\ShowEq{f in L(A->B)}{A^n}B{}
be a polylinear map.
Then the map
\ShowEq{symmetrization of polylinear map}
\EqParm{<f>}{=z}
defined by the equality
\ShowEq{symmetrization of polylinear map =}
is symmetric polylinear map
and is called
\AddIndex{symmetrization of polylinear map}{symmetrization of polylinear map}
$f$.
\end{theorem}
\begin{proof}
The theorem follows from the equality
\ShowEq{gs=sum fss}
and the definition
\RefDefinition{symmetric polylinear map}.
\end{proof}

\begin{definition}
\labelDefinition{skew symmetric polylinear map}
{\it
Let $A$, $B$ be algebras over commutative ring $D$.
A polylinear map
\ShowEq{f in L(A->B)}{A^n}B{}
is called
\AddIndex{skew symmetric}{skew symmetric polylinear map},
if
\ShowEq{fa=sfsa}
for any permutation $\sigma$ of the set
\EqParm{set a1n}{=.}
}
\qed
\end{definition}

\begin{theorem}
Let
\ShowEq{f in L(A->B)}{A^n}B{}
be a polylinear map.
Then the map
\ShowEq{alternation of polylinear map}
\EqParm{[f]}{=z}
defined by the equality
\ShowEq{alternation of polylinear map =}
is skew symmetric polylinear map
and is called
\AddIndex{alternation of polylinear map}{alternation of polylinear map}
$f$.
\end{theorem}
\begin{proof}
The theorem follows from the equality
\ShowEq{gs=sum sfss}
and the definition
\RefDefinition{skew symmetric polylinear map}.
\end{proof}

\begin{theorem}
\labelTheorem{fx1n=0, xi=xi1}
A polylinear map
\ShowEq{f in L(A->B)}{A^n}B{}
is skew symmetric iff
\ShowEq{fx1n=0}
as soon as
$x_i=x_{i+1}$
for at list one\footnote{
In the book
\citeBib{Cartan differential form},
page 9,
Henri Cartan considered the theorem
\RefTheorem{fx1n=0, xi=xi1}
as definition of skew symmetric map.
}
\ShowEq{1<=i<n}
\end{theorem}
\begin{proof}
For given $i$, consider permutation
\ShowEq{sigma i}
It is evident that
\ShowEq{|sigma i|}

Let
\ShowEq{f in L(A->B)}{A^n}B{}
be skew symmetric map.
According the definition
\RefDefinition{skew symmetric polylinear map},
From equalities
\EqRef{sigma i},
\EqRef{|sigma i|},
it follows that
\ShowEq{fxii}

Let
\ShowEq{f in L(A->B)}{A^n}B.
Let
\ShowEq{fx1n=0}
as soon as
$x_i=x_{i+1}$
for at list one
\ShowEq{1<=i<n}
Then
\ShowEq{0=f xi+xi1}
From equalities
\EqRef{sigma i},
\EqRef{|sigma i|},
\EqRef{0=f xi+xi1},
it follows that
\ShowEq{fxii1=-fxi1i}
Since any permutation
is product of permutations $\sigma_i$,
then, from the equality
\EqRef{fxii1=-fxi1i}
and the definition
\RefDefinition{skew symmetric polylinear map},
it follows that the map $f$ is skew symmetric.
\end{proof}

\begin{theorem}
\labelTheorem{module of skew symmetric polylinear maps}
The set
\ShowEq{module of skew symmetric polylinear maps}
of skew symmetric polylinear maps
is $D$\Hyph module.
\end{theorem}
\begin{proof}
According to the theorem
\RefTheorem{L(An;B) is free D module},
linear composition of skew symmetric polylinear maps
is polylinear map.
According to the theorem
\RefTheorem{fx1n=0, xi=xi1},
linear composition of skew symmetric polylinear maps
is skew symmetric polylinear map.
\end{proof}

Without loss of generality, we assume
\ShowEq{L(A1,A0,B)=}

\begin{theorem}
Let $A$, $B$ be Banach $D$\Hyph algebras.
$D$\Hyph module $\LAnAB$ is closed in $D$\Hyph module
\ShowEq{L(A;B)}D{A^n}B.
\end{theorem}
\begin{proof}
The theorem follows from theorems
\RefTheorem{fx1n=0, xi=xi1},
\RefTheorem{module of skew symmetric polylinear maps},
since the equality
\ShowEq{fxixi=0}
holds in passage to the limit,
when sequence of maps $f_k$ converges to map $f$.
\end{proof}

\Section{Product of Skew Symmetric Polylinear Maps}

\begin{definition}
\labelDefinition{product of symmetric polylinear maps}
{\it
Let $A$, $B_2$ be free algebras over commutative ring $D$.\,\footnote{
To define product of skew symmetric polylinear maps,
I follow definition in section
\citeBib{Cartan differential form}-1.4 of chapter 1,
pages 12 - 14.
}
Let
\ShowEq{h:B1->*B2}
be left\Hyph side representation of
free associative $D$\Hyph algebra $B_1$ in $D$\Hyph module $B_2$.
The map
\ShowEq{hpq:Lpq->Lp+q}
is defined by the equality
\ShowEq{hpq(fg)=}
where, in the right side of the equality
\EqRef{hpq(fg)=},
we consider left\Hyph side product
of $B_2$\Hyph number
\EqParm{g()}{=z}
over $B_1$\Hyph number
\EqParm{f()}{=.}
}
\qed
\end{definition}

If $B_1=B_2$, then, in the right side of the equality
\EqRef{hpq(fg)=},
we consider product of $B_1$\Hyph numbers
\EqParm{f()}{=c}
\EqParm{g()}{=.}
According to the theorem
\RefTheorem{Free Algebra over Ring},
this definition is compatible with the definition
\RefDefinition{product of symmetric polylinear maps}.

Let $f$, $g$ be skew symmetric polylinear maps.
In general, the map
\ShowEq{hpq(fg)}
is not skew symmetric polylinear map.

\begin{definition}
\labelDefinition{exterior product polylinear map}
{\it
The skew symmetric polylinear map
\ShowEq{exterior product}
\ShowEq{exterior product =}
is called
\AddIndex{exterior product}{exterior product}.
}
\qed
\end{definition}

\begin{theorem}
Exterior product satisfies the following equation
\ShowEq{exterior product = skew symmetric}
\end{theorem}
\begin{proof}
The equality
\EqRef{exterior product = skew symmetric}
follows from equalities
\EqRef{alternation of polylinear map =},
\EqRef{exterior product =}.
\end{proof}

\begin{theorem}
\labelTheorem{fAf ne 0}
Let $B$ be non commutative $D$\Hyph algebra and
\ShowEq{f:A->B}fAB
be linear map. Then, in general,\,\footnote{
See also the equality
(\citeBib{Sudbery Quaternionic Analysis}\Hyph (2.25).
}
\ShowEq{fAf ne 0}
\end{theorem}
\begin{proof}
According to the equality
\EqRef{exterior product = skew symmetric},
\ShowEq{fAf=}
The expression
\EqRef{fAf=},
in general, is different from $0$.
\end{proof}

From the proof of the theorem
\RefTheorem{fAf ne 0},
it follows that
\ShowEq{fAf=0}
only when the image of the map $f$ is
a subset of the center of the algebra $B$.

\begin{convention}
\labelConvention{IJK=()}
Let
\ShowEq{IJK=()}
The order of index in sets $J$, $K$ is the same as in the set $I$.
Let
\ShowEq{RJI}
be the set of injections
\ShowEq{J->I}
For any map
\EqParm{l in RJI}{=c}
let
\ShowEq{D(l)}
be the range of the map $\lambda$ and
\ShowEq{D1=I-D}
The map
\ShowEq{K->D1}
is defined by the equality
\ShowEq{mu(l)()=}
Let
\ShowEq{S1(l)}
be the set of permutations of the set
\ShowEq{D1(l)}
\qed
\end{convention}

Let the set $I_{m.n}$ be defined by the equality
\ShowEq{Imn=}

\begin{lemma}
\labelLemma{sigma=lambda,tau}
{\it
For any permutation
\EqParm{s in S}{=c}
there exists unique map
\EqParm{l in RJI}{=z}
and unique permutation
\EqParm{t in S1l}{=z}
such that
\DrawEq{sigma=lambda,tau}{->}
}
\end{lemma}
\begin{proof}
Let
\EqParm{s in S}{=.}
We define the map
\EqParm{l in RJI}{=z}
using the table
\ShowEq{lambda=}
We define the map
\ShowEq{t:D1->D1}
by the equality
\ShowEq{t(l)=s(a)}
The equality
\ShowEq{t(m(a))=s(a)}
follows from equalities
\EqRef{mu(l)()=},
\EqRef{t(l)=s(a)}.
Since maps $\sigma$, $\mu(\lambda)$
are bijection, then the map $\tau$ is bijection.
Therefore, the map $\tau$ is permutation of the set
\ShowEq{D1(l)}
The equality
\eqRef{sigma=lambda,tau}{->}
follows from equalities
\EqRef{lambda=},
\EqRef{t(m(a))=s(a)}.
\end{proof}

\begin{lemma}
\labelLemma{lambda,tau=sigma}
{\it
For any map
\EqParm{l in RJI}{=z}
and any permutation
\EqParm{t in S1l}{=c}
there exists unique permutation
\EqParm{s in S}{=z}
such that
\DrawEq{sigma=lambda,tau}{<-}
}
\end{lemma}
\begin{proof}
According to
\EqRef{D1=I-D},
\ShowEq{D1+D=I}
Since maps
\ShowEq{l m(l) t}
are injections, then
the map $\sigma$ is permutation.
\end{proof}

\ePrints{4975-6381}
\ifx\Semafor\ValueOn
\newpage
\fi
According to lemmas
\RefLemma{sigma=lambda,tau},
\RefLemma{lambda,tau=sigma},
there exist bijection of the set $S$ into set of tuples
\ShowEq{(lambda,tau)}
where
\EqParm{l in RJI}{=c}
\EqParm{t in S1l}{=.}
We express this bijection by the equality
\ShowEq{s=(lambda,tau)}

\begin{lemma}
\labelLemma{sign of the map l}
{\it
We define sign of the map
\EqParm{l in RJI}{=z}
using the equality
\ShowEq{|l|=|l,d|}
where $\delta$ is identity permutation of the set
\ShowEq{D1(l)}
Then
\ShowEq{|lt|=|l||t|}
}
\end{lemma}
\begin{proof}
Since
\ShowEq{|d|=1}
then the theorem holds for permutation
\EqParm{(l,d)}{=.}

Let permutation $\tau$ be transposition.
Then the permutation
\ShowEq{(lambda,tau)}
is different from the permutation
\EqParm{(l,d)}{=z}
by transposition.
Therefore, in this case, the equality
\EqRef{|lt|=|l||t|}
holds.

Let permutation $\tau$ be product of $n$ transpositions.
Then the permutation
\ShowEq{(lambda,tau)}
is different from the permutation
\EqParm{(l,d)}{=z}
by product of $n$ transpositions.
Therefore, in general case, the equality
\EqRef{|lt|=|l||t|}
holds.
\end{proof}

\begin{theorem}
\labelTheorem{f1 f2 = sum l in Rp}
Let
\ShowEq{f in LAAB}1p1
\ShowEq{f in LAAB}2q2
Let $R_0$ be the set of injections
\ShowEq{I1p->I1pq}
which preserve order of index in the set $I_{1.p+q}$.
Let
\ShowEq{R1=I1p->I1pq}
Exterior product of skew symmetric polylinear maps satisfies to equalities\,\footnote{
In the book
\citeBib{Cartan differential form}, pages 12, 13,
the equality
\EqRef{f1 f2 = sum l in R0p}
is the definition of exterior product of skew symmetric polylinear maps,
since the map
\ShowEq{apq->l}
is the permutation which satisfy to the requirement (1.4.3) on the page
\citeBib{Cartan differential form}\Hyph 12.
}
\ShowEq{f1 f2 = sum l in R1p}
\ShowEq{f1 f2 = sum l in R0p}
\end{theorem}
\begin{proof}
Let $S$ be the set of permutations of the set
\ShowEq{I1pq}
According to the definition
\EqRef{exterior product = skew symmetric},
\ShowEq{f1 A f2 =}
According to lemmas
\RefLemma{sigma=lambda,tau},
\RefLemma{lambda,tau=sigma},
\RefLemma{sign of the map l},
\ShowEq{f1 A f2 = 1}
Since the map $f_2$ is skew symmetric, then, according to the definition
\RefDefinition{skew symmetric polylinear map},
\ShowEq{f2 l = |s| f2 sl}
The equality
\ShowEq{f1 A f2 = 2}
follows from equalities
\EqRef{f1 A f2 = 1},
\EqRef{f2 l = |s| f2 sl}.
The equality
\EqRef{f1 f2 = sum l in R1p}
follows from the equality
\EqRef{f1 A f2 = 2}.

For any
\EqParm{l in R0}{=c}
consider the set
\ShowEq{R2(l)=...}
The map
\ShowEq{m in R2}
differs from the map
\EqParm{l in R0}{=z}
by a permutation of the set $D(\lambda)$ and
\ShowEq{|l|f1=|m|f1}
At the same time
\ShowEq{R1=U R2}
The equality
\ShowEq{f1 A f2 = 3}
follows from equalities
\EqRef{f1 f2 = sum l in R1p},
\EqRef{|l|f1=|m|f1},
\EqRef{R1=U R2}.
The equality
\EqRef{f1 f2 = sum l in R0p}
follows from the equality
\EqRef{f1 A f2 = 3}.
\end{proof}

\begin{theorem}
Let
\ShowEq{f in LAAB}1p1
\ShowEq{f in LAAB}2q2
Let $R_0$ be the set of injections
\ShowEq{Ip1q->I1pq}
which preserve order of index in the set $I_{1.p+q}$.
Let
\ShowEq{R1=Ip1q->I1pq}
Exterior product of skew symmetric polylinear maps satisfies to equalities
\ShowEq{f1 f2 = sum l in R1q}
\ShowEq{f1 f2 = sum l in R0q}
\end{theorem}
\begin{proof}
The proof of the theorem is similar to the proof of the theorem
\RefTheorem{f1 f2 = sum l in Rp}.
\end{proof}

\begin{theorem}
Exterior product of skew symmetric polylinear maps is associative.
Let
\ShowEq{f in LAAB}1p1
\ShowEq{f in LAAB}2q1
\ShowEq{f in LAAB}3r2
Then
\ShowEq{f1(f2f3)=(f1f2)f3 wedge}
\end{theorem}
\begin{proof}
Let the representation of $D$\Hyph algebra $B_1$ in $D$\Hyph module $B_2$
is defined by the map $h$ and the product in $D$\Hyph algebra $B_1$
is defined by the map $g$.

Let $S$ be the set of permutations of the set
\ShowEq{I1pqr}

According to the definition
\RefDefinition[0912.3315]{left-side representation of group},
the following equality
\ShowEq{f1(f2f3)=(f1f2)f3}
holds for any $B_1$\Hyph numbers
\ShowEq{f1()f2()}
and for any $B_2$\Hyph number
\ShowEq{f3()}
According to the definition
\RefDefinition{product of symmetric polylinear maps},
we can write down the equality
\EqRef{f1(f2f3)=(f1f2)f3}
as
\ShowEq{f1(f2f3)=(f1f2)f3 1}
The equality
\EqRef{f1(f2f3)=(f1f2)f3 1}
remains true when we apply permutation
\ShowEq{s in SI}
So the equality
\ShowEq{f1(f2f3)=(f1f2)f3 2}
follows from the equality
\EqRef{f1(f2f3)=(f1f2)f3 1}.

Let
\ShowEq{R1=I1p->I1pqr}
According to lemmas
\RefLemma{sigma=lambda,tau},
\RefLemma{lambda,tau=sigma},
\RefLemma{sign of the map l},
we rewrite left side of the equality
\EqRef{f1(f2f3)=(f1f2)f3 2}
as
\ShowEq{f1(f2f3)}
According to the definitions
\EqRef{alternation of polylinear map =},
\EqRef{exterior product =},
the equality
\ShowEq{f1(f2f3) 2}
follows from the equality
\EqRef{f1(f2f3)}.
According to lemmas
\RefLemma{sigma=lambda,tau},
\RefLemma{lambda,tau=sigma},
\RefLemma{sign of the map l},
the equality
\ShowEq{f1(f2f3) 3}
follows from equalities
\EqRef{exterior product = skew symmetric},
\EqRef{f1(f2f3) 2}.
According to the definition
\RefDefinition{product of symmetric polylinear maps},
we can write down the equality
\EqRef{f1(f2f3) 3}
as
\ShowEq{f1(f2f3) 4}
According to definitions
\EqRef{alternation of polylinear map =},
\EqRef{exterior product =},
the equality
\ShowEq{f1(f2f3) 5}
follows from the equality
\EqRef{f1(f2f3) 4}.

Let
\ShowEq{R2=Ipqpqr->I1pqr}
According to lemmas
\RefLemma{sigma=lambda,tau},
\RefLemma{lambda,tau=sigma},
\RefLemma{sign of the map l},
we rewrite right side of the equality
\EqRef{f1(f2f3)=(f1f2)f3 2}
as
\ShowEq{(f1f2)f3}
According to the definitions
\EqRef{alternation of polylinear map =},
\EqRef{exterior product =},
the equality
\ShowEq{(f1f2)f3 2}
follows from the equality
\EqRef{(f1f2)f3}.
According to lemmas
\RefLemma{sigma=lambda,tau},
\RefLemma{lambda,tau=sigma},
\RefLemma{sign of the map l},
the equality
\ShowEq{(f1f2)f3 3}
follows from equalities
\EqRef{exterior product = skew symmetric},
\EqRef{(f1f2)f3 2}.
According to the definition
\RefDefinition{product of symmetric polylinear maps},
we can write down the equality
\EqRef{(f1f2)f3 3}
as
\ShowEq{(f1f2)f3 4}
According to definitions
\EqRef{alternation of polylinear map =},
\EqRef{exterior product =},
the equality
\ShowEq{(f1f2)f3 5}
follows from the equality
\EqRef{(f1f2)f3 4}.

The equality
\EqRef{f1(f2f3)=(f1f2)f3 wedge}
follows from equalities
\EqRef{f1(f2f3)=(f1f2)f3 1},
\EqRef{f1(f2f3) 5},
\EqRef{(f1f2)f3 5}.
\end{proof}

\Section{Differential Form}

\begin{definition}
\labelDefinition{map of class Cn}
{\it
A map is of class $C^n$,
if the map has continuous derivative of order $n$.
A map is of class $C^{\infty}$,
if the map has continuous derivative of any order.
}
\qed
\end{definition}

\begin{definition}
{\it
Let $A$, $B$ be Banach $D$\Hyph algebras.
Let
\ShowEq{U subset A}
be an open set.
The map
\ShowEq{o:U->LAAB}
is called
\AddIndex{differential form of degree $p$}{differential form of degree p}
defined in $U$ and with values in $B$
or $B$\Hyph valued
\AddIndex{differential $p$\Hyph form}{differential p form}
defined in $U$.
A differential form $\omega$ is said
to be of class
\ShowEq{class Cn}
if map $\omega$ is of class $C^n$.
The set
\ShowEq{set of differential p forms}
of $B$\Hyph valued differential $p$\Hyph forms
of class $C^n$ defined in $U$
is $D$\Hyph module.
}
\qed
\end{definition}

\begin{theorem}
Let $A$, $B_2$ be free Banach algebras over commutative ring $D$.\,\footnote{
To define product of differential forms,
I follow definition in section
\citeBib{Cartan differential form}-2.2 of chapter 1,
page 19.
}
Let
\ShowEq{h:B1->*B2}
be left\Hyph side representation of
free associative Banach $D$\Hyph algebra $B_1$ in Banach $D$\Hyph module $B_2$.
Let
\ShowEq{diff form in U}{\alpha}np{B_1},
\ShowEq{diff form in U}{\beta}nq{B_2}.
\AddIndex{Exterior product}{exterior product}
\ShowEq{exterior product, differential forms}
of differential forms $\alpha$, $\beta$
defined by the equality
\ShowEq{a wedge b x =}
is differential form of class $C^n$.
\end{theorem}
\begin{proof}
For any $x\in U$,
\EqParm{form x in LA}{form=alpha,A=AB,=c}
\EqParm{form x in LA}{form=beta,A=AB,=.}
According to definitions
\RefDefinition{product of symmetric polylinear maps},
\RefDefinition{exterior product polylinear map},
\ShowEq{a(x) wedge b(x) in}
The map
\ShowEq{x->a(x) wedge b(x)}
is a map of class $C^n$
since this map is composition of the map
\ShowEq{x->(a(x),b(x))}
of class $C^n$ and continuous bilinear map
\ShowEq{LA1xLA2->LA}
generated by representation $h$.
\end{proof}

\begin{corollary}
\labelCorollary{Exterior product generates bilinear map}
{\it
Exterior product of differential forms
generates a bilinear map
\EqParm{O(B1)xO(B2)->O(B2)}{A=B2}
}
\qed
\end{corollary}

\begin{theorem}
Let $A$, $B_2$ be free Banach algebras over commutative ring $D$.
Let
\ShowEq{h:B1->*B2}
be left\Hyph side representation of
free associative Banach $D$\Hyph algebra $B_1$ in Banach $D$\Hyph module $B_2$.
Let
\ShowEq{diff form in U}{\alpha}np{B_1},
\ShowEq{diff form in U}{\beta}nq{B_2}.
Then
\ShowEq{(a wedge b)(x)()=}
\end{theorem}
\begin{proof}
The theorem follows from equalities
\EqRef{alternation of polylinear map =},
\EqRef{hpq(fg)=},
\EqRef{exterior product =}.
\end{proof}

For a given integer $n$,
consider $D$\Hyph module
\ShowEq{OnB=...}
which is direct sum of $D$\Hyph modules
\ShowEq{OnB1}
$D$\Hyph module
\EqParm{OnB}{=z}
equipped with exterior product
\EqParm{O(B1)xO(B2)->O(B2)}{A=B1}
is graded algebra.

\begin{theorem}
\labelTheorem{1 form in finite dimensional Banach algebra}
Let $\Basis e_A$ be basis
of finite dimensional Banach algebra $A$ over commutative ring $D$.
Let $\Basis e_B$ be basis
of finite dimensional Banach algebra $B$ over commutative ring $D$.
Let
\ShowEq{o1p x}
be coordinates of polylinear map
\ShowEq{o(x) in L(An;B)}
relative to bases $\Basis e_A$, $\Basis e_B$.
Then the differential $p$\Hyph form $\omega$ has form
\ShowEq{w(a1p)=...}
where
\ShowEq{a1p=...}
\ShowEq{w1p=...}
\end{theorem}
\begin{proof}
Since
\EqParm{a1p}{n=1,=z}
are any $A$\Hyph numbers, then equalities
\EqRef{w(a1p)=...},
\EqRef{w1p=...}
follow from the equality
\EqRef{a1p=...}
and from the equality
\ShowEq{w(a1n)=a1n...}
\end{proof}

A map
\ShowEq{f:A->B}fAB
is an example of $B$\Hyph valued differential $0$\Hyph form.

A derivative of the map
\ShowEq{f:A->B}fAB
is an example of $B$\Hyph valued differential $1$\Hyph form.

\begin{theorem}
\labelTheorem{derivative in finite dimensional Banach algebra}
Let $\Basis e$ be basis
of finite dimensional Banach algebra $A$ over field $D$.
We may consider map 
\ShowEq{f:A->B}fAA
as map
\ShowEq{f(x)=f(x1n)}
where
\ShowEq{x=x1n}
Since the map $f$ is differentiable,
then there exist partial derivatives
\ShowEq{df/dxi 1n}
and the following equality is true
\DrawEq{df/dxi=}{theorem}
\end{theorem}
\begin{proof}
According to the definition
\RefDefinition{differential of map},
differential of the map $f$ has form
\ShowEq{df=df/dx o dx}
Since the map $f$ is the map
of $D$\Hyph vector space,
then differentials have form
\DrawEq{dx=dx e}{differential}
\ShowEq{df=df/dxi dxi}
The equality
\ShowEq{df/dx o dx=df/dxi dxi}
follows from equalities
\EqRef{df=df/dx o dx},
\eqRef{dx=dx e}{differential},
\EqRef{df=df/dxi dxi}.
The equality
\eqRef{df/dxi=}{theorem}
follows from the equality
\EqRef{df/dx o dx=df/dxi dxi}.
\end{proof}

It is easy to see that partial derivative
\EqParm{df/dxi}{=z}
represented as
\eqRef{df/dxi=}{theorem}
is a map of $D$\Hyph algebra $A$.
If this map is differentiable,
then we can consider partial derivative of order $2$
\ShowEq{d/dxj df/dxi=}

\begin{theorem}
\labelTheorem{second derivative is symmetric bilinear map}
Let $A$ be free Banach $D$\Hyph module.
Let the derivative of second order of the map
\ShowEq{f:A->B}fAB
be continuous.
Then the second derivative is symmetric bilinear map
\ShowEq{df/dx2 ab=ba}
\end{theorem}
\begin{proof}
Since the map
\ShowEq{d2f/dx2}
is continuous,
then the equality
\EqRef{d/dxj df/dxi=}
imply that the map
\ShowEq{d2f/dxij}
is continuous. 
According to the theorem
in \citeBib{Fikhtengolts: Calculus volume 1}
(pages 405, 406),
\ShowEq{df/dxij=df/dxji}
The equality
\ShowEq{df/dx2 ij=ji}
follows from equalities
\EqRef{d/dxj df/dxi=},
\EqRef{df/dxij=df/dxji}.
The equality
\EqRef{df/dx2 ab=ba}
follows from equalities
\ShowEq{a= b= e}
and from the equality
\EqRef{df/dx2 ij=ji}.
\end{proof}

\Section{Exterior Differentiation}

Let $A$, $B$ be Banach $D$\Hyph algebras.
Let
\ShowEq{U subset A}
be an open set.
Let the map
\ShowEq{o:U->LAAB}
be differential $p$\Hyph form
of class $C_n$, $n>0$.
According to the definition
\RefDefinition{map of class Cn},
at every point $x\in U$
the increment of map
$\omega$ can be represented as
\ShowEq{derivative of omega}
where
\ShowEq{d omega in Cn-1}
According to the definition
\EqRef{derivative of omega},
the map
\ShowEq{d omega a0p}
is skew symmetric polylinear map
with respect to variables
\EqParm{a1p}{n=1,=z}
and is a linear map with respect to variable $a_0$.
However, the map
\EqRef{d omega a0p}
is not skew symmetric polylinear map
with respect to variables
\ShowEq{a0p}

\begin{definition}
\labelDefinition{exterior differential}
{\it
Let $A$, $B$ be Banach $D$\Hyph algebras.
Let
\ShowEq{U subset A}
be an open set.
Let the map
\ShowEq{o:U->LAAB}
be differential $p$\Hyph form
of class $C_n$, $n>0$.
The map
\ShowEq{exterior differential}
\ShowEq{exterior differential def}
is called
\AddIndex{exterior differential}{exterior differential}.
}
\qed
\end{definition}

\begin{theorem}
Exterior differential holds the equality
\ShowEq{do=...sum}
\end{theorem}
\begin{proof}
The equality
\EqRef{do=...sum}
follows from equalities
\EqRef{alternation of polylinear map =},
\EqRef{exterior differential def}.
\end{proof}

\ifx\texFuture\Defined
\begin{convention}
Let, in the convention
\RefConvention{IJK=()},
\ShowEq{I=a0p,J=a0}
Then we can identify any map
\EqParm{l in RJI}{=z}
with permutation
\ShowEq{ti in S p+1}
which is defined by tuple
\ShowEq{ti a0p}
From the equality
\EqRef{ti a0p},
it follows that
\qed
\end{convention}
\fi

\begin{lemma}
\labelLemma{s->(ti mu)}
{\it
Let $S_1$ be the set of permutations of the set of $A$\Hyph numbers
\EqParm{a1p}{n=0,=.}
Let $S$ be the set of permutations of the set of $A$\Hyph numbers
\EqParm{a1p}{n=1,=.}
Let permutation
\ShowEq{ti in S1}
be defined by tuple
\ShowEq{ti a0p}
For any permutation
\ShowEq{s in S1}
there exist permutations
\ShowEq{ti in S1, s1 in S}
such that
\DrawEq{s=(ti mu)}{->}
}
\end{lemma}
\begin{proof}
The permutation $\tau_i$ is defined by the equality
\ShowEq{ti 0 = s 0}
The permutation $\mu$ is defined by the equality
\ShowEq{mu=()}
The equality
\eqRef{s=(ti mu)}{->}
follows from the equalities
\EqRef{ti 0 = s 0},
\EqRef{mu=()}.
\end{proof}

\begin{lemma}
\labelLemma{(ti mu)->s}
{\it
For any permutation $\mu$ of the set of $A$\Hyph numbers
\ShowEq{a0p - ai}
there exists the permutation $\sigma$ of the set of $A$\Hyph numbers
\EqParm{a1p}{n=0,=z}
such that
\DrawEq{s=(ti mu)}{<-}
\ShowEq{|s|=(-1)i|mu|}
}
\end{lemma}
\begin{proof}
The permutation $\sigma$ has form
\ShowEq{s=mu i ne 0}
The equality
\eqRef{s=(ti mu)}{<-}
follows from the equality
\EqRef{s=mu i ne 0}.
The equality
\ShowEq{|s|=|t mu|}
follows from the equality
\EqRef{s=mu i ne 0}.
Since the permuation $\tau_i$ is
the product of $i$ transpositions, then
\ShowEq{|ti|=(-1)i}
The equality
\EqRef{|s|=(-1)i|mu|}
follows from the equalities
\EqRef{|s|=|t mu|},
\EqRef{|ti|=(-1)i}.
\end{proof}

\begin{lemma}
\labelLemma{sum s1 =sum +-s}
{\it
Let $S_1$ be the set of permutations of the set of $A$\Hyph numbers
\EqParm{a1p}{n=0,=.}
Let $S$ be the set of permutations of the set of $A$\Hyph numbers
\EqParm{a1p}{n=1,=.}
For any polylinear map $f$
\ShowEq{sum s1 =sum +-s}
}
\end{lemma}
\begin{proof}
The lemma follows from lemmas
\RefLemma{s->(ti mu)},
\RefLemma{(ti mu)->s}.
\end{proof}

\begin{theorem}
\labelTheorem{exterior differential of differential form}
We can write exterior differential of differential $p$\Hyph form as\,\footnote{
According to the definition
\EqRef{ti a0p},
we can write the equality
\EqRef{d omega = sum tau}
as
\ShowEq{d omega = sum}
}
\ShowEq{d omega = sum tau}
\end{theorem}
\begin{proof}
The equality
\ShowEq{d omega = sum 1}
follows from equalities
\EqRef{alternation of polylinear map =},
\EqRef{exterior differential def}.
According to the lemma
\RefLemma{sum s1 =sum +-s},
the equality
\ShowEq{d omega = sum 2}
follows from the equality
\EqRef{d omega = sum 1}.
The equality
\EqRef{d omega = sum tau}
follows from the equality
\EqRef{d omega = sum 2},
since the map
\ShowEq{dx omega}
is skew symmetric with respect to variables
\EqParm{a1p}{n=1,=z}
and we can write summand corresponding to identity permutation $\mu$
instead of sum over permutations $\mu$.
\end{proof}

\begin{example}
Let
\ShowEq{diff form in U}{\omega}n2B.
Then
\ShowEq{d omega2=}
\qed
\end{example}

\begin{theorem}
Let $A$, $B_2$ be free Banach algebras over commutative ring $D$.\,\footnote{
The theorem is simiral to the theorem
\citeBib{Cartan differential form}-2.4.2,
page 22.
However I cannot use Cartan's proof
because he uses commutativity of product in real field. 
}
Let
\ShowEq{h:B1->*B2}
be left\Hyph side representation of
free associative Banach $D$\Hyph algebra $B_1$ in Banach $D$\Hyph module $B_2$.
Let
\ShowEq{diff form in U}{\alpha}np{B_1},
\ShowEq{diff form in U}{\beta}nq{B_2}.
Then
\ShowEq{d(a wedge b)=}
\end{theorem}
\begin{proof}
According to the theorem
\RefTheorem{bilinear map and differential}
and to the corollary
\RefCorollary{Exterior product generates bilinear map},
\ShowEq{d(a A b)/dx=...}
The equality
\ShowEq{d(a A b)/dx=..., 1}
follows from the equality
\EqRef{d(a A b)/dx=...}.
The equality
\ShowEq{d(a A b)/dx=..., 2}
follows from equalities
\EqRef{d omega = sum tau},
\EqRef{d(a A b)/dx=..., 1}.

Let $S_1$ be the set of permutations of the set of $A$\Hyph numbers
\EqParm{a1p+q}{n=0}
Let $S$ be the set of permutations of the set of $A$\Hyph numbers
\EqParm{a1p+q}{n=1}

According to definitions
\EqRef{alternation of polylinear map =},
\EqRef{exterior product =},
\ShowEq{d(a A b)/dx=..., 2 1.1}
According to the lemma
\RefLemma{sum s1 =sum +-s},
the equality
\ShowEq{d(a A b)/dx=..., 2 1.2}
follows from the equality
\EqRef{d(a A b)/dx=..., 2 1.1}.
Let
\ShowEq{R1=Ip+1q>I0p+q}
According to lemmas
\RefLemma{sigma=lambda,tau},
\RefLemma{lambda,tau=sigma},
\RefLemma{sign of the map l},
the equality
\ShowEq{d(a A b)/dx=..., 2 1.3}
follows from the equality
\EqRef{d(a A b)/dx=..., 2 1.2}.
According to definitions
\EqRef{alternation of polylinear map =},
\EqRef{exterior differential def},
the equality
\ShowEq{d(a A b)/dx=..., 2 1.4}
follows from equalities
\EqRef{do=...sum},
\EqRef{d(a A b)/dx=..., 2 1.3}.
The equality
\ShowEq{d(a A b)/dx=..., 2 1.5}
follows from equalities
\EqRef{f1 f2 = sum l in R1q},
\EqRef{d(a A b)/dx=..., 2 1.4}.

According to definitions
\EqRef{alternation of polylinear map =},
\EqRef{exterior product =},
\ShowEq{d(a A b)/dx=..., 2 2.1}
According to the lemma
\RefLemma{sum s1 =sum +-s},
the equality
\ShowEq{d(a A b)/dx=..., 2 2.2}
follows from the equality
\EqRef{d(a A b)/dx=..., 2 2.1}.
Let
\ShowEq{R1=I1p->I0p+q}
According to lemmas
\RefLemma{sigma=lambda,tau},
\RefLemma{lambda,tau=sigma},
\RefLemma{sign of the map l},
the equality
\ShowEq{d(a A b)/dx=..., 2 2.3}
follows from the equality
\EqRef{d(a A b)/dx=..., 2 2.2}.
According to definitions
\EqRef{alternation of polylinear map =},
\EqRef{exterior differential def},
the equality
\ShowEq{d(a A b)/dx=..., 2 2.4}
follows from equalities
\EqRef{do=...sum},
\EqRef{d(a A b)/dx=..., 2 2.3}.
The equality
\ShowEq{d(a A b)/dx=..., 2 2.5}
follows from equalities
\EqRef{f1 f2 = sum l in R1p},
\EqRef{d(a A b)/dx=..., 2 2.4}.
Since
\ShowEq{|a1p a0 apq|=}
then the equality
\ShowEq{d(a A b)/dx=..., 2 2.6}
follows from the equality
\EqRef{d(a A b)/dx=..., 2 2.5}.

The equality
\ShowEq{d(a A b)/dx=..., 3}
follows from equalities
\EqRef{d(a A b)/dx=..., 2},
\EqRef{d(a A b)/dx=..., 2 1.5},
\EqRef{d(a A b)/dx=..., 2 2.6}.
The equality
\EqRef{d(a wedge b)=}
follows from the equality
\EqRef{d(a A b)/dx=..., 3}.
\end{proof}

\begin{theorem}
\labelTheorem{d2omega=0}
Let $A$ be free Banach $D$\Hyph module.
If
\ShowEq{diff form in U}{\omega}npB,
$n\ge 2$, then\,\footnote{
See also the theorem
\citeBib{Cartan differential form}-2.5.1
on the page 23.
}
\ShowEq{d2omega=0}
\end{theorem}
\begin{proof}
According to the definition
\RefDefinition{exterior differential},
\ShowEq{d2o=d/dx do/dx}
Let $S$ be set of permutations of the set
\ShowEq{b12I1p}
According to the definition
\EqRef{alternation of polylinear map =},
the equality
\ShowEq{d2o=d/dx do/dx 1}
follows from the equality
\EqRef{d2o=d/dx do/dx}.
According to the theorem
\RefTheorem{second derivative is symmetric bilinear map},
the map
\ShowEq{d2o/dx2}
is symmetric with respect to $A$\Hyph numbers $b_1$, $b_2$.
The equality
\EqRef{d2omega=0}
follows from definitions
\RefDefinition{symmetric polylinear map},
\RefDefinition{skew symmetric polylinear map}.
\end{proof}

\Section{Poincar\'e's Theorem}

\begin{definition}
{\it
Subset $U$ of $D$\Hyph algebra $B$ is called
\AddIndex{starlike}{starlike set}
with respect to the $B$\Hyph number
\EqParm{a in U}{=c}
if,
for any $B$\Hyph number
\ShowEq{x in U}
and real number
\EqParm{0<t<1}{=c}
\EqParm{t'a+tx in U}{=.}
}
\qed
\end{definition}

\begin{definition}
\labelDefinition{integrable form}
{\it
Let $A$, $B$ be Banach algebras.
Let
\ShowEq{U subset A}
be open set.
The differential form
\ShowEq{diff form in U}{\omega}npB{}
is called
\AddIndex{integrable}{integrable form},
if there exists a differential form
\ShowEq{diff form in U}{\alpha}n{p-1}B{}
such that
\DrawEq{da=o}{}
Then we use notation
\ShowEq{indefinite integral form}
\ShowEq{a=int o}
and the differential form $\alpha$ is called
\AddIndex{indefinite integral}{indefinite integral}
of the differential form $\omega$.
}
\qed
\end{definition}

\begin{theorem}
\labelTheorem{map int continuous is continuous}
Let $A$, $B$ be Banach algebras.\,\footnote{
The statement is similar to the lemma
\citeBib{Cartan differential form}\Hyph 2.12.2
on the page 33.
}
Let
\ShowEq{[01]}
be subset of real field.
Let
\ShowEq{phi:UxI->B}
be continued map.
Then the map
\ShowEq{psi:U->B}
defined by the equality
\DrawEq{psi=int phi}{}
is continuous.
\end{theorem}
\begin{proof}
Choose
\EqParm{eps>0}{=.}
\StartLabelItem
\begin{enumerate}
\item
For each
\ShowEq{x,t in U I}
there exists
\ShowEq{eta xt}
such that the statement
\ShowEq{|dx|<,|dt|<}
implies
\ShowEq{|d phi|< 1}
\labelItem{|d phi|< 1}
\item
In particular, the statement
\ShowEq{|dt|<}
implies
\ShowEq{|d phi|< 2}
\labelItem{|d phi|< 2}
\end{enumerate}
From statements
\RefItem{|d phi|< 1},
\RefItem{|d phi|< 2},
it follows that the statement
\ShowEq{|dx|<,|dt|<}
implies
\DrawEq{|d phi|< 3}{}
Let
\ShowEq{I(t)=}
Since the set $I$ is compact,
then, for given $x$,
the finite set of intervals $I(x,t)$
which cover the set $I$,
for instance
\ShowEq{I(x,ti)}
Let $\eta(x)$ be the smallest number of $\eta(x,t_i)$.
Since, for any $t'\in I$, there exists $t_i$ such that
\ShowEq{|dti|<}
then the statement
\ShowEq{|dx|<}
implies
\DrawEq{|d phi|< 3}{1}
for any $t'\in I$.
Therefore, the map $\varphi$ is continuous at $x$
uniformly with respect to $t$.

Since
\ShowEq{d psi=}
then, according to the theorem
\RefTheorem{|int f|<int|f|},
\ShowEq{|d psi|<}
From inequalities
\eqRef{|d phi|< 3}{1},
\EqRef{|d psi|<}
and the theorem
\RefTheorem{int |f|<M mu},
the statement
\ShowEq{|dx|<}
implies
\ShowEq{|d psi|< 1}
Therefore, the map $\psi$ is continuous.
\end{proof}

\begin{lemma}
\labelLemma{o1=int o, |o1h/h|->0}
{\it
Let
\ShowEq{o:AxI->B}
be such continuous map that
\DrawEq{|oh/h|->0}{lemma}
Then
\DrawEq{o1=int o}{}
is map
\ShowEq{f:A->B}{o_1}AB
such that
\DrawEq{|o1h/h|->0}{lemma}
}
\end{lemma}
\begin{proof}
From the equality
\eqRef{|oh/h|->0}{lemma},
it follows that
we can represent the map $o$ as product
\ShowEq{o=o2 h}
where
\ShowEq{lim o2=0}
Continuity of the map $o_2$ follows from the continuity of the map $o$.
Since the set $I$ is compact, then there exists the map
\ShowEq{o3:A->R}
defined by the equality
\ShowEq{o3=max o2}

\begin{lemma}
\labelLemma{o3(h)->0}
{\it
\ShowEq{o3(h)->0}
}
\end{lemma}

{\sc Proof.}
Let
\ShowEq{o3(h) ne 0}
Then there exists
\ShowEq{eps>0}
such that
\ShowEq{o3>e}
for any
\EqParm{delta>0}{=z}
and for any
\ShowEq{h<d}
Since the set $I$ is compact, then there exists $t=t(h)$
such that
\ShowEq{o2>e}
for any
\ShowEq{h<d}
This statement contradicts to the equality
\EqRef{lim o2=0}.
Therefore, the statement
\EqRef{o3(h) ne 0}
is false.
\hfill\(\odot\)

According to the theorem
\RefTheorem{|int h|<int|omega||f1n|},
\ShowEq{|o1|<}
The statement
\eqRef{|o1h/h|->0}{lemma}
follows from the statement
\EqRef{|o1|<}
and from the lemma
\RefLemma{o3(h)->0}.
\end{proof}

\begin{theorem}
\labelTheorem{differentiation under integral}
Let $A$, $B$ be Banach algebras.\,\footnote{
The statement is similar to the lemma
\citeBib{Cartan differential form}\Hyph 2.12.2
on the page 33.
}
Let
\ShowEq{[01]}
be subset of real field.
Let
\ShowEq{phi:UxI->B}
be continued map.
Let the map
\ShowEq{psi:U->B}
be defined by the equality
\DrawEq{psi=int phi}{1}
If derivative
\ShowEq{d phi/dx}
exists at every point
\EqParm{x,t in U I}{=c}
then the map $\psi$ is differentiable and
\ShowEq{d psi/dx}
\end{theorem}
\begin{proof}
Let the derivative
\ShowEq{d phi/dx}
exist and be continuous map
\ShowEq{d phi:UxI->LAB}
According to the theorem
\RefTheorem{map int continuous is continuous},
the map
\ShowEq{l=int d phi}
is linear and continuous.
According to the theorem
\RefTheorem{h=fX g1 g2 representation},
the equality
\ShowEq{l=int d phi h}
follows from the equality
\EqRef{l=int d phi}.

According to the definition
\RefDefinition{differentiable map, algebra},
\ShowEq{phi x+h - x}
where
\DrawEq{|oh/h|->0}{integral}
According to the theorem
\RefTheorem{int f+g X},
the equality
\ShowEq{phi x+h - x int}
follows from the equality
\EqRef{phi x+h - x}.
Let
\DrawEq{o1=int o}{integral}
According to the lemma
\RefLemma{o1=int o, |o1h/h|->0},
the equality
\DrawEq{|o1h/h|->0}{integral}
follows from equalities
\eqRef{|oh/h|->0}{integral},
\eqRef{o1=int o}{integral}.
The equality
\ShowEq{phi x+h - x int 1}
follows from equalities
\ShowEq{phi x+h - x int 2}
The equality
\EqRef{d psi/dx}
follows from equalities
\EqRef{l=int d phi},
\eqRef{|o1h/h|->0}{integral},
\EqRef{phi x+h - x int 1}
and from the definition
\RefDefinition{differentiable map, algebra}.
\end{proof}

\begin{theorem}
\labelTheorem{differentiation order n under integral}
Let $A$, $B$ be Banach algebras.\,\footnote{
The statement is similar to the corollary
\citeBib{Cartan differential form}\Hyph 2.12.3
on the page 35.
}
Let
\ShowEq{[01]}
be subset of real field.
Let
\ShowEq{phi:UxI->B}
be map of class $C^n$.
Let the map
\ShowEq{psi:U->B}
be defined by the equality
\DrawEq{psi=int phi}{}
Then the map $\psi$ is of class $C^n$ and
\ShowEq{d psi/dxn}
\end{theorem}
\begin{proof}
The proof follows by induction on $n$, since the theorem follows from the theorem
\RefTheorem{differentiation under integral}
and the definition
\EqRef{derivative of Order n, algebra}
\end{proof}

\begin{theorem}[Poincar\'e's Theorem]
\labelTheorem{Poincare Theorem p=1}
Let $A$, $B$ be Banach algebras.\,\footnote{
I follow the theorem
\citeBib{Cartan differential form}\Hyph 2.12.1
on the page 33.
}
Let
\ShowEq{U subset A}
be open set, starlike with respect to the $A$\Hyph number
\EqParm{a in U}{=.}
If the differential form
\ShowEq{diff form in U}{\omega}n1B{}
satisfies to the equation
\DrawEq{do=0}{}
then the differential form $\omega$ is integrable.
\end{theorem}
\begin{proof}
Let $a=0$. In general, we can achieve this
by a translation.

Let
\ShowEq{diff form in U}{\omega}n1B.
Set
\ShowEq{f=int o(tx,x)}
The map
\ShowEq{phi(x,t)=o(tx,x)}
is the composition of following maps

\begin{itemize}
\item
the map
\ShowEq{phi1:UxI->UxA}
exists because the set $U$
is starlike with respect to $A$\Hyph number $0$
and this map is of class $C^{\infty}$;
\item
the map
\ShowEq{phi2:UxA->LxA}
is of class $C^n$;
\item
the bilinear map
\ShowEq{phi3:LxA->B}
is of class $C^{\infty}$.
\end{itemize}
Therefore, the map $\varphi$ is of class $C^n$.
According to theorems
\RefTheorem{map int continuous is continuous},
\RefTheorem{differentiation order n under integral},
the map $f$ is of class $C^n$.

According to the theorem
\RefTheorem{differentiation under integral},
\ShowEq{df/dx=int...}
According to the theorem
\RefTheorem{bilinear map and differential},
the equality
\ShowEq{d phi(x,t)/dt=}
follows from the equality
\EqRef{phi(x,t)=o(tx,x)}.
From the statement
\DrawEq{do=0}{-}
it follows that
\ShowEq{dw a1 a2}
is a bilinear symmetric map of $a_1$ and $a_2$.
Therefore,
\ShowEq{dw a1 a2=dw a2 a1}
The equality
\ShowEq{d phi(x,t)/dt= 1}
follows from equalities
\EqRef{d phi(x,t)/dt=},
\EqRef{dw a1 a2=dw a2 a1}.
The equality
\ShowEq{d phi(x,t)/dt= 2}
follows from the equality
\EqRef{d phi(x,t)/dt= 1}.
At the same time
\ShowEq{d tw/dt=}
The equality
\ShowEq{df/dx=int... 1}
follows from equalities
\EqRef{df/dx=int...},
\EqRef{d phi(x,t)/dt= 2},
\EqRef{d tw/dt=}.
Set
\ShowEq{g=t omega}
The equality
\ShowEq{df/dx=int... 2}
follows from equalities
\EqRef{df/dx=int... 1},
\EqRef{g=t omega}.
Therefore, $\omega$ is equal to the differential form $df$.
\end{proof}

\begin{theorem}[Poincar\'e's Theorem]
Let $A$, $B$ be Banach algebras.\,\footnote{
I follow the theorem
\citeBib{Cartan differential form}\Hyph 2.12.1
on the page 33.
}
Let
\ShowEq{U subset A}
be open set, starlike with respect to the $A$\Hyph number
\EqParm{a in U}{=.}
If the differential form
\ShowEq{diff form in U}{\omega}npB{}
satisfies to the equation
\DrawEq{do=0}{}
then the differential form $\omega$ is integrable.
\end{theorem}
\begin{proof}
For every $p$, we define the linear map
\ShowEq{k:Onp->Onp-1}
as follows.

\StartLabelItem
\begin{enumerate}
\item
If
\ShowEq{f in O0n}
let $k(f)=0$.
\labelItem{map k p=0}
\item
For $p\ge 1$ and
\ShowEq{diff form in U}{\omega}npB,
set
\DrawEq{k(o)=a}{item}
where
\ShowEq{alpha=int dt omega}
\end{enumerate}

\begin{lemma}
{\it
The map $k$ is well defined.
}
\end{lemma}

{\sc Proof.}
For $p=0$, the definition
\RefItem{map k p=0}
is the same as saying
that $D$\Hyph module
\ShowEq{On,-1}
is reduced to $0$.

For $p>0$,
consider skew symmetric polylinear map
\ShowEq{w(tx,x):}
Then we can write the equality
\EqRef{alpha=int dt omega}
as follows
\ShowEq{alpha=int dt omega 1}
From
\EqRef{alpha=int dt omega 1}
it follows that
\ShowEq{alpha in LA}

The map
\ShowEq{x,t->w}
is the composition of following maps

\begin{itemize}
\item
the map
\ShowEq{a1:UxI->UxIxA}
exists because the set $U$
is starlike with respect to $A$\Hyph number $0$
and this map is of class $C^{\infty}$;
\item
the map
\ShowEq{a2:UxIxA->LxA}
is of class $C^n$;
\item
the map
\ShowEq{a3:LxA->L}
is of class $C^{\infty}$.
\end{itemize}
Therefore, $\alpha$ is differential $(p-1)$\Hyph form of class $C^n$.
\hfill\(\odot\)

\begin{lemma}
{\it
If
\ShowEq{f in O0n}
$n\ge 1$, then
\ShowEq{k df=f-f0}
where $f_0$ is constant map
\ShowEq{f0:U->B}
}
\end{lemma}

{\sc Proof.}
The differential form $k(df)$ is the map
\ShowEq{x->int01}
\hfill\(\odot\)

\begin{lemma}
{\it
If
\ShowEq{diff form in U}{\omega}npB,
$n\ge 1$, $p\ge 1$, then
\ShowEq{d(k w)...}
}
\end{lemma}

{\sc Proof.}
The equality
\ShowEq{dw(tx)=}
follows from the equality
\EqRef{d omega = sum}.
From equalities
\EqRef{alpha=int dt omega},
\EqRef{dw(tx)=},
it follows that the differential form
\ShowEq{kdw=b}
is defined by the equality
\ShowEq{kdw=int}
According to the theorem
\RefTheorem{differentiation under integral},
the equality
\ShowEq{alpha'=int dt omega'}
follows from the equality
\EqRef{alpha=int dt omega}.
The equality
\ShowEq{d alpha=int dt omega'}
follows from equalities
\EqRef{d omega = sum},
\EqRef{alpha'=int dt omega'}.
Since $\omega$ is skew symmetric map,
then the following equality is true
\ShowEq{p int tw}
The equality
\ShowEq{dk+kd=}
follows from equalities
\ShowEq{dk+kd= ref}
The equality
\ShowEq{dk+kd= 1}
follows from the equality
\EqRef{dk+kd=}
At the same time
\ShowEq{d tpw/dt=}
The equality
\ShowEq{dk+kd= 2}
follows from equalities
\EqRef{dk+kd= 1},
\EqRef{d tpw/dt=}.
The equality
\EqRef{d(k w)...}
follows from the equality
\EqRef{dk+kd= 2}.
\hfill\(\odot\)

If
\DrawEq{do=0}{-},
then the equality
\ShowEq{w=d(k w)}
follows from the equality
\EqRef{d(k w)...}.
Therefore, the map $k$ furnishes the differential form
\DrawEq{k(o)=a}{-}
such that
\DrawEq{da=o}{-}.
\end{proof}

%% file: Diff.Form.1.Eq.tex

\input{Diff.Form.1.Parm}

\DefEq
{
$\displaystyle\frac{\partial f}{\partial x^{\gii}}$\Pt
}
{df/dxi}

\DefEquation
{
\frac{\partial^2 f}{\partial x^{\gij}\partial x^{\gii}}=
\frac{d}{d x}
\left(
\frac{d f}{d x}\circ e_{\gii}
\right)
\circ e_{\gij}
=
\frac{d^2 f}{d x^2}\circ (e_{\gii};e_{\gij})
}
{d/dxj df/dxi=}

\DefEq
{
\[
\omega:U\rightarrow\mathcal{LA}(D;A^p\rightarrow B)
\]
}
{o:U->LAAB}

\DefEq
{
\[
\alpha(x)\wedge\beta(x)\in\mathcal {LA}(D;A^{p+q}\rightarrow B_2)
\]
}
{a(x) wedge b(x) in}

\DefEq
{
$x \in U$
}
{x in U}

\DefEq
{
$\delta>0$\Pt
}
{delta>0}

\DefEq
{
$h$, $\|h\|_A<\delta$.
}
{h<d}

\DefEq
{
$o_3(h)\ge\epsilon$
}
{o3>e}

\DefEq
{
$\|o_2(h,t(h))\|_B\ge\epsilon$
}
{o2>e}

\DefEq
{
$I=[0,1]$
}
{[01]}

\DefEq
{
\[
\varphi:U\times I\rightarrow B
\]
}
{phi:UxI->B}

\DefEq
{
\[
\frac{\partial\varphi(x,t)}{\partial x}
:U\times I\rightarrow\mathcal L(D;A\rightarrow B)
\]
}
{d phi:UxI->LAB}

\DefEq
{
d\omega(x)=0
}
{do=0}

\DefEq
{
\symb{\int \omega(x)}{indefinite integral}{form}
}
{indefinite integral form}

\DefEq
{
\[
\alpha(x)=
\ShowSymbol{indefinite integral}{form}
\]
}
{a=int o}

\DefEq
{
$\epsilon>0$\Pt
}
{eps>0}

\DefEq
{
$\|x'-x\|_A<\eta(x,t)$, $|t'-t|<\eta(x,t)$
}
{|dx|<,|dt|<}

\DefEq
{
$\|x'-x\|_A\le\eta(x)$
}
{|dx|<}

\DefEq
{
$|t'-t|<\eta(x,t)$
}
{|dt|<}

\DefEq
{
$|t'-t_i|<\eta(x,t_i)$,
}
{|dti|<}

\DefEq
{
\[\|\varphi(x',t')-\varphi(x,t)\|_B\le\frac{\epsilon}2\]
}
{|d phi|< 1}

\DefEq
{
\[\|\varphi(x,t')-\varphi(x,t)\|_B\le\frac{\epsilon}2\]
}
{|d phi|< 2}

\DefEq
{
$\{I(x,t_i)\}$.
}
{I(x,ti)}

\DefEq
{
\[
I(x,t)=\{t'\in I:|t'-t|<\eta(x,t)\}
\]
}
{I(t)=}

\DefEq
{
\|\varphi(x',t')-\varphi(x,t')\|_B\le\epsilon
}
{|d phi|< 3}

\DefEq
{
$\eta(x,t)\in R$
}
{eta xt}

\DefEq
{
\[\psi:U\rightarrow B\]
}
{psi:U->B}

\DefEq
{
\[
\psi(x')-\psi(x)=
\int_0^1dt(\varphi(x',t)-\varphi(x,t))
\]
}
{d psi=}

\DefEquation
{
\|\psi(x')-\psi(x)\|_B=
\int_0^1dt\|\varphi(x',t)-\varphi(x,t)\|_B
}
{|d psi|<}

\DefEq
{
\[
\|\psi(x')-\psi(x)\|_B=
\int_0^1dt\epsilon=\epsilon
\]
}
{|d psi|< 1}

\DefEq
{
\[o:A\times I\rightarrow B\]
}
{o:AxI->B}

\DefEq
{
\lim_{\|h\|_A\rightarrow 0}\frac{\|o(h,t)\|_B}{\|h\|_A}=0
}
{|oh/h|->0}

\DefEq
{
\lim_{\|h\|_A\rightarrow 0}\frac{\|o_1(h)\|_B}{\|h\|_A}=0
}
{|o1h/h|->0}

\DefEquation
{
\varphi(x+h,t)-\varphi(x,t)=\frac{\partial\varphi(x,t)}{\partial x}\circ h+o(h,t)
}
{phi x+h - x}

\DefEquation
{
\int_0^1dt\,\varphi(x+h,t)-\int_0^1dt\,\varphi(x,t)
=\int_0^1dt\frac{\partial\varphi(x,t)}{\partial x}\circ h+\int_0^1dt\,o(h,t)
}
{phi x+h - x int}

\DefEquation
{
\psi(x+h)-\psi(x)=\lambda(x)\circ h+o_1(h)
}
{phi x+h - x int 1}

\DefEq
{
\eqRef{psi=int phi}{1},
\EqRef{l=int d phi},
\EqRef{phi x+h - x int},
\eqRef{o1=int o}{integral}.
}
{phi x+h - x int 2}

\DefEquation
{
\lim_{h\rightarrow 0} o_2(h,t)=0
}
{lim o2=0}

\DefEquation
{
o(h,t)=o_2(h,t)h
}
{o=o2 h}

\DefEq
{
o_1(h)=\int_0^1dt\,o(h,t)
}
{o1=int o}

\DefEquation
{
\begin{split}
\|o_1(h)\|_B=\left|\int_0^1dt\,o(h,t)\right|_B
&\le\int_0^1dt\|o_2(h,t)\|_B\|h\|_A
\\&\le\int_0^1dt\,o_3(h)\|h\|_A
=o_3(h)\|h\|_A
\end{split}
}
{|o1|<}

\DefEq
{
\[o_3:A\rightarrow R\]
}
{o3:A->R}

\AddEquation{o3(h)->0}
{
\lim_{h\rightarrow 0}o_3(h)=0
}

\AddEquation{o3(h) ne 0}
{
\lim_{h\rightarrow 0}o_3(h)\ne 0
}

\DefEq
{
\[
o_3(h)=\max(\|o_2(h,t)\|_B,t\in I)
\]
}
{o3=max o2}

\DefEq
{
$(x,t)\in U\times I$\Pt
}
{x,t in U I}

\DefEquation
{
\frac{d\psi(x)}{d x}=\int_0^1dt\frac{\partial\varphi(x,t)}{\partial x}
}
{d psi/dx}

\DefEquation
{
\frac{d^n\psi(x)}{d x^n}=\int_0^1dt\frac{\partial^n\varphi(x,t)}{\partial x^n}
}
{d psi/dxn}

\DefEquation
{
\varphi(x,t)=\omega(tx)\circ x
}
{phi(x,t)=o(tx,x)}

\DefEq
{
$(d_x\omega(x)\circ a_1)\circ a_2$
}
{dw a1 a2}

\DefEquation
{
\left(\frac{d \omega(tx)}{d tx}\circ a\right)\circ x
=
\left(\frac{d \omega(tx)}{d tx}\circ x\right)\circ a
}
{dw a1 a2=dw a2 a1}

\DefEquation
{
\begin{split}
\frac{\partial\varphi(x,t)}{\partial x}\circ a
&=\left(\frac{d \omega(tx)}{d tx}\circ\frac{d tx}{d x}\circ a\right)\circ x
+\omega(tx)\circ\left(\frac{d x}{d x}\circ a\right)
\\&=\left(\frac{d \omega(tx)}{d tx}\circ(ta)\right)\circ x
+\omega(tx)\circ a
\\&=t\left(\frac{d \omega(tx)}{d tx}\circ a\right)\circ x
+\omega(tx)\circ a
\end{split}
}
{d phi(x,t)/dt=}

\DefEquation
{
\frac{\partial \varphi(x,t)}{\partial x}\circ a
=t\left(\frac{d\omega(tx)}{dtx}\circ x\right)\circ a
+\omega(tx)\circ a
}
{d phi(x,t)/dt= 1}

\DefEquation
{
\frac{\partial \varphi(x,t)}{\partial x}
=t\left(\frac{d\omega(tx)}{dtx}\circ x\right)
+\omega(tx)
}
{d phi(x,t)/dt= 2}

\DefEquation
{
\begin{split}
\frac{dt\omega(tx)}{dt}
&=t\left(\frac{d\omega(tx)}{dtx}\circ\frac{dtx}{dt}\right)
+\omega(tx)
\\&=t\left(\frac{d\omega(tx)}{dtx}\circ x\right)
+\omega(tx)
\end{split}
}
{d tw/dt=}

\DefEquation
{
\begin{split}
\frac{d t^p\omega(tx)}{d t}
&=t^p\left(\frac{d \omega(tx)}{d tx}\circ\frac{d tx}{d t}\right)
+pt^{p-1}\omega(tx)
\\&=t^p\left(\frac{d \omega(tx)}{d tx}\circ x\right)
+pt^{p-1}\omega(tx)
\end{split}
}
{d tpw/dt=}

\DefEq
{
\[
\varphi_1:(x,t)\in U\times I\rightarrow (tx,x)\in U\times A
\]
}
{phi1:UxI->UxA}

\DefEq
{
\[
\alpha_1:(x,t)\in U\times I\rightarrow (tx,t,x)\in U\times I\times A
\]
}
{a1:UxI->UxIxA}

\DefEq
{
\[
\varphi_2:(x,a)\in U\times A\rightarrow (\omega(x),a)\in \mathcal L(D;A\rightarrow B)\times A
\]
}
{phi2:UxA->LxA}

\DefEq
{
\[
\alpha_2:(x,t,a)\in U\times I\times A\rightarrow (t^{p-1}\omega(x),a)\in \mathcal L(D;A^p\rightarrow B)\times A
\]
}
{a2:UxIxA->LxA}

\DefEq
{
\[
\varphi_3:(f,a)\in \mathcal L(D;A\rightarrow B)\times A\rightarrow f\circ a\in B
\]
}
{phi3:LxA->B}

\DefEq
{
\[
\alpha_3:(f,a)\in \mathcal L(D;A^p\rightarrow B)\times A\rightarrow f\circ a\in \mathcal L(D;A^{p-1}\rightarrow B)
\]
}
{a3:LxA->L}

\DefEq
{
\[\omega(tx)\circ x:(a_1,...,a_{p-1})\rightarrow \omega(tx)\circ(x,a_1,...,a_{p-1})\]
}
{w(tx,x):}

\DefEquation
{
\frac{d f(x)}{dx}=\int_0^1 dt\frac{\partial \varphi(x,t)}{\partial x}
}
{df/dx=int...}

\DefEquation
{
\frac{df(x)}{dx}=\int_0^1 dt\frac d{dt}(t\omega(tx))
}
{df/dx=int... 1}

\DefEquation
{
d(k(\omega))+k(d\omega)=\int_0^1 dt\frac{d t^p\omega(tx)}{d t}
}
{dk+kd= 2}

\DefEquation
{
\frac{df(x)}{dx}=g(1)-g(0)=\omega(x)
}
{df/dx=int... 2}

\DefEq
{
$\Omega^{(n)}_{-1}(U,B)$
}
{On,-1}

\DefEquation
{
\alpha(x)\circ(a_1,...,a_{p-1})=\int_0^1dt\,t^{p-1}\omega(tx)\circ(x,a_1,...,a_{p-1})
}
{alpha=int dt omega}

\DefEquation
{
\alpha(x)\circ(a_1,...,a_{p-1})=\left(\int_0^1dt\,t^{p-1}\omega(tx)\circ x\right)\circ(a_1,...,a_{p-1})
}
{alpha=int dt omega 1}

\DefEquation
{
\begin{split}
\left(\frac{d\alpha(x)}{d x}\circ a_1\right)\circ(a_2,...,a_p)
&=\int_0^1dt\,t^p\left(\frac{d\omega(tx)}{dx}\circ a_1\right)\circ(x,a_2,...,a_p)
\\&+\int_0^1dt\,t^{p-1}\omega(tx)\circ(a_1,...,a_p)
\end{split}
}
{alpha'=int dt omega'}

\DefEquation
{
\begin{split}
&\,d\alpha(x)\circ (a_1,...,a_p)
\\=&\,\sum_{i=1}^p(-1)^{i-1}
\int_0^1dt\,t^p\left(\frac{d\omega(tx)}{d x}\circ a_i\right)\circ(x,a_1,...,\overbrace{a_i},...,a_p)
\\+&\,\sum_{i=1}^p(-1)^{i-1}\int_0^1dt\,t^{p-1}\omega(tx)\circ(a_i,a_1,...,\overbrace{a_i},...,a_p)
\end{split}
}
{d alpha=int dt omega'}

\DefEq
{
$\alpha(x)\in\mathcal {LA}(D;A^{p-1}\rightarrow B)$.
}
{alpha in LA}

\DefEquation
{
\begin{split}
&\,\sum_{i=1}^p(-1)^{i-1}\int_0^1dt\,t^{p-1}\omega(tx)\circ(a_i,a_1,...,\overbrace{a_i},...,a_p)
\\=&\,p\int_0^1dt\,t^{p-1}\omega(tx)\circ(a_1,...,a_p)
\end{split}
}
{p int tw}

\DefEquation
{
\begin{split}
&\,d(k(\omega))\circ(a_1,...,a_p)+k(d\omega)\circ(a_1,...,a_p)
\\=&\,(d\alpha)\circ(a_1,...,a_p)+\beta\circ(a_1,...,a_p)
\\=&\,\int_0^1dt\,t^p
\left(\frac{d\omega(tx)}{dx}\circ x\right)\circ(a_1,...,a_p)
+p\int_0^1dt\,t^{p-1}\omega(tx)\circ(a_1,...,a_p)
\end{split}
}
{dk+kd=}

\DefEquation
{
\begin{split}
d(k(\omega))+k(d\omega)
&=\int_0^1dt\,t^p
\frac{d\omega(tx)}{dx}\circ x
+p\int_0^1dt\,t^{p-1}\omega(tx)
\\&=\int_0^1dt\left(t^p
\frac{d\omega(tx)}{dx}\circ x
+pt^{p-1}\omega(tx)\right)
\end{split}
}
{dk+kd= 1}

\DefEq
{
\eqRef{k(o)=a}{item},
\EqRef{kdw=b},
\EqRef{kdw=int},
\EqRef{d alpha=int dt omega'},
\EqRef{p int tw}.
}
{dk+kd= ref}

\DefEq
{
\[(x,t)\rightarrow t^{p-1}\omega(tx,x)\]
}
{x,t->w}

\DefEq
{
k(\omega)=\alpha
}
{k(o)=a}

\DefEquation
{
k(d\omega)=\beta
}
{kdw=b}

\DefEquation
{
\begin{split}
d\omega(tx)\circ(x,a_1,...,a_p)
&=\left(\frac{d\omega(tx)}{dx}\circ x\right)\circ(a_1,...,a_p)
\\&+\sum_{i=1}^p(-1)^p\left(\frac{d\omega(tx)}{dx}\circ a^i\right)
\circ(x,a_1,...,\overbrace{a^i},...,a^p)
\end{split}
}
{dw(tx)=}

\DefEquation
{
\begin{split}
&\,\beta(x)\circ(a_1,...,a_p)
\\=&\,\int_0^1dt\,t^p
\left(\frac{d\omega(tx)}{dx}\circ x\right)\circ(a_1,...,a_p)
\\+&\,\sum_{i=1}^p(-1)^i\int_0^1dt\,t^p\left(\frac{d\omega(tx)}{dx}\circ a^i\right)
\circ(x,a_1,...,\overbrace{a^i},...,a^p)
\end{split}
}
{kdw=int}

\DefEq
{
\[
x\rightarrow \int_0^1dt\left(\frac{df(tx)}{dx}\circ x\right)
=\int_0^1dt\frac{df(tx)}{dt}=f(x)-f(0)
\]
}
{x->int01}

\DefEq
{
$f\in\Omega_0^{(n)}$,
}
{f in O0n}

\DefEquation
{
d(k(\omega))+k(d\omega)=\omega
}
{d(k w)...}

\DefEquation
{
\omega=d(k(\omega))
}
{w=d(k w)}

\DefEquation
{
k(df)=f-f_0
}
{k df=f-f0}

\DefEq
{
\[f_0:x\in U\rightarrow f(0)\in B\]
}
{f0:U->B}

\DefEq
{
\[
k:\Omega^{(n)}_p(U,B)\rightarrow \Omega^{(n)}_{p-1}(U,B)
\]
}
{k:Onp->Onp-1}

\DefEquation
{
g(t)=t\omega(tx)
}
{g=t omega}

\DefEquation
{
f(x)=\int_0^1dt\,\omega(tx)\circ x
}
{f=int o(tx,x)}

\DefEq
{
$\partial_x\varphi(x,t)$
}
{d phi/dx}

\DefEq
{
\psi(x)=\int_0^1dt\varphi(x,t)
}
{psi=int phi}

\DefEquation
{
\lambda(x)=\int_0^1dt
\frac{\partial \varphi(x,t)}{\partial x}
}
{l=int d phi}

\DefEquation
{
\lambda(x)\circ h=\int_0^1dt
\frac{\partial \varphi(x,t)}{\partial x}\circ h
}
{l=int d phi h}

\DefEq
{
d\alpha(x)=\omega(x)
}
{da=o}

\DefEq
{
$0\le t\le 1$\Pt
}
{0<t<1}

\DefEq
{%
$(1-t)a+tx\in U$\Pt
}
{t'a+tx in U}

\AddEq [5]{diff form in U}
{
$#1\in\Omega^{(#2)}_{#3}(U,#4)$#5
}%

\DefEq
{
\[
d\omega\circ(a_1,a_2,a_3)
=\frac{d\omega}{dx}\circ(a_1,a_2,a_3)
-\frac{d\omega}{dx}\circ(a_2,a_1,a_3)
+\frac{d\omega}{dx}\circ(a_3,a_1,a_2)
\]
}
{d omega2=}

\DefEq
{
$\{b_1,b_2\}\cup I_{1.p}$.
}
{b12I1p}

\DefEquation
{
d^2\omega\circ(b_2,b_1,a_1,...,a_p)=
(p+1)(p+2)\frac{d^2\omega}{dx^2}\circ[b_2,b_1,a_1,...,a_p]
}
{d2o=d/dx do/dx}

\DefEquation
{
\begin{split}
&\,d^2\omega\circ(b_2,b_1,a_1,...,a_p)
\\=&\,
\frac{(p+2)!}{p!}\sum_{\sigma\in S}
\frac{d^2\omega}{dx^2}\circ(\sigma(b_2),\sigma(b_1),\sigma(a_1),...,\sigma(a_p))
\\=&\,
\frac{(p+2)!}{p!}\sum_{\sigma\in S}
\left(
\frac{d^2\omega}{dx^2}\circ(\sigma(b_2),\sigma(b_1))
\right)
\circ(\sigma(a_1),...,\sigma(a_p))
\end{split}
}
{d2o=d/dx do/dx 1}

\DefEq
{
$\displaystyle\frac{d^2\omega}{dx^2}$
}
{d2o/dx2}

\DefEquation
{
d^2\omega(x)=0
}
{d2omega=0}

\DefEq
{
$\Pform(x)\in\mathcal{LA}(\pD;\pA^{\Pp}\rightarrow \pB)$\Pt
}
{form x in LA}

\DefEq
{
\symb{\Omega^{(n)}_p(U,B)}{set of differential p forms}{1}
}
{set of differential p forms}

\DefEq
{
\symb{C^n}{class Cn}{1}
}
{class Cn}

\DefEq
{
$U\subseteq A$
}
{U subset A}

\DefEquation
{
\omega(x)\circ (a_1,...,a_p)
=\omega_{\gii_1,...,\gii_p}(x)a_1^{\gii_1}...a_n^{\gii_p}
}
{w(a1p)=...}

\DefEquation
{
\omega(x+a_0)-\omega(x)
=\frac{d \omega(x)}{d x}\circ a_0
+o(dx)
}
{derivative of omega}

\DefEquation
{
\left(\frac{d \omega(x)}{d x}\circ a_0\right)
\circ(a_1,...,a_p)
}
{d omega a0p}

\DefEq
{
\[
\frac{d \omega(x)}{d x}\in\mathcal {LA}(D;A^{p+1}\rightarrow B)
\ \ \ \frac{d \omega(x)}{d x}\in C_{n-1}
\]
}
{d omega in Cn-1}

\DefEq
{
$\omega(x)\in\mathcal L(D;A^n\rightarrow B)$, $x\in U$,
}
{o(x) in L(An;B)}

\DefEquation
{
\begin{split}
0&=f\circ(...,a_i+a_{i+1},a_i+a_{i+1},...)\\
&=f\circ(...,a_i,a_i+a_{i+1},...)+f\circ(...,a_{i+1},a_i+a_{i+1},...)\\
&=f\circ(...,a_i,a_i,...)+f\circ(...,a_i,a_{i+1},...)\\
&+f\circ(...,a_{i+1},a_i,...)+f\circ(...,a_{i+1},a_{i+1},...)\\
&=f\circ(...,a_i,a_{i+1},...)+f\circ(...,a_{i+1},a_i,...)
\end{split}
}
{0=f xi+xi1}

\DefEq
{
$a_0$, $a_1$, ..., $a_p$
}
{a0p}

\DefEq
{
\symb{d\omega}{exterior differential}{}
}
{exterior differential}

\DefEquation
{
\ShowSymbol{exterior differential}{}\circ(a_0,...,a_p)
=(p+1)\left[\frac{d\omega}{dx}\right]\circ(a_0,...,a_p)
}
{exterior differential def}

\DefEquation
{
d\omega\circ(a_0,...,a_p)
=\frac 1{p!}\sum_{\sigma\in S(p+1)}
\frac{d\omega}{dx}\circ(\sigma(a_0),...,\sigma(a_p))
}
{do=...sum}

\DefEquation
{
d\omega\circ(a_0,...,a_p)=\sum_{i=0}^p(-1)^i
\left(\frac{d\omega}{dx}\circ a_i\right)
\circ(a_0,...,\overbrace{a_i},...,a_p)
}
{d omega = sum}

\DefEquation
{
\begin{split}
d\omega\circ(a_0,...,a_p)&=\frac {p+1}{(p+1)!}
\sum_{\sigma\in S(p+1)}|\sigma|
\frac{d\omega(x)}{dx}
\circ(\sigma(a_0),...,\sigma(a_p))
\\
&=\frac 1{p!}
\sum_{\sigma\in S(p+1)}|\sigma|
\frac{d\omega(x)}{dx}
\circ(\sigma(a_0),...,\sigma(a_p))
\end{split}
}
{d omega = sum 1}

\DefEq
{
$\partial_x\omega(x)$
}
{dx omega}

\DefEquation
{
d\omega\circ(a_0,...,a_p)=
\sum_{i=0}^p\frac {(-1)^i}{p!}\sum_{\mu\in S(p)}
\frac{d\omega(x)}{dx}
\circ(\tau_i(a_0),\tau_i(\mu(a_1)),...,\tau_i(\mu(a_p)))
}
{d omega = sum 2}

\DefEquation
{
d\omega\circ(a_0,...,a_p)=\sum_{i=0}^p(-1)^i
\frac{d\omega}{dx}
\circ(\tau_i(a_0),...,\tau_i(a_p))
}
{d omega = sum tau}

\DefEq
{
\symb{<f>}{symmetrization of polylinear map}{}
}
{symmetrization of polylinear map}

\DefEq
{
$\ShowSymbol{symmetrization of polylinear map}{}\circ(a_1,...,a_n)$\Pt
}
{<f>}

\DefEq
{
\[
\ShowSymbol{symmetrization of polylinear map}{}\circ(a_1,...,a_n)
=\frac 1{n!}
\sum_{\sigma\in S(n)}f\circ\sigma\circ(a_1,...,a_n)
\]
}
{symmetrization of polylinear map =}

\DefEq
{
\symb{[f]}{alternation of polylinear map}{}
}
{alternation of polylinear map}

\DefEquation
{
\ShowSymbol{alternation of polylinear map}{}\circ(a_1,...,a_n)
=\frac 1{n!}\sum_{\sigma\in S(n)}|\sigma|(f\circ\sigma\circ(a_1,...,a_n))
}
{alternation of polylinear map =}

\DefEq
{
$\ShowSymbol{alternation of polylinear map}{}$\Pt
}
{[f]}

\DefEq
{
\[
f\circ(a_1, ..., a_n)=f\circ\sigma\circ(a_1,...,a_n)
\]
}
{fa=fsa}

\DefEq
{
\begin{align*}
<f>\circ\sigma\circ (a_1,...,a_n)
&=<f>\circ(\sigma(a_1),...,\sigma(a_n))
\\
&=\frac 1{n!}\sum_{\sigma_1\in S(n)}f\circ\sigma_1\circ(\sigma(a_1),...,\sigma(a_n))
\\
&=\frac 1{n!}\sum_{\sigma_2\in S(n)}f\circ(\sigma_2(a_1),...,\sigma(a_n))
&\sigma_2=\sigma_1\circ\sigma
\\
&=\frac 1{n!}\sum_{\sigma_2\in S(n)}f\circ\sigma_2\circ(a_1,...,a_n)
\\
&=<f>\circ (a_1,...,a_n)
\end{align*}
}
{gs=sum fss}

\DefEq
{
\begin{align*}
[f]\circ\sigma\circ (a_1,...,a_n)
&=[f]\circ(\sigma(a_1),...,\sigma(a_n))
\\
&=\frac 1{n!}\sum_{\sigma_1\in S(n)}|\sigma_1|f\circ\sigma_1\circ(\sigma(a_1),...,\sigma(a_n))
\\
&=\frac 1{n!}\sum_{\sigma_2\in S(n)}|\sigma_2||\sigma|f\circ(\sigma_2(a_1),...,\sigma_2(a_n))
&\sigma_2&=\sigma_1\circ\sigma
\\
&=|\sigma|\frac 1{n!}\sum_{\sigma_2\in S(n)}|\sigma_2|f\circ\sigma_2\circ(a_1,...,a_n)
&|\sigma_2|&=|\sigma_1||\sigma|
\\
&=|\sigma|[f]\circ (a_1,...,a_n)
\end{align*}
}
{gs=sum sfss}

\DefEq
{
\[
\Omega^{(n)}(U,B_1)=\bigoplus_{p\ge 0}\Omega^{(n)}_p(U,B_1)
\]
}
{OnB=...}

\DefEq
{
$\Omega^{(n)}_p(U,B_1)$, $p=0$, $1$, ... .
}
{OnB1}

\DefEq
{
\begin{align*}
&\,(\alpha\wedge\beta)(x)\circ(a_1,...,a_{p+q})
\\=&\,\frac {(p+q)!}{p!q!}\sum_{\sigma\in S(n)}|\sigma|(\alpha(x)\circ(\sigma(a_1),...,\sigma(a_p)))
(\beta(x)\circ(\sigma(a_{p+1}),...,\sigma(a_{p+q})))
\end{align*}
}
{(a wedge b)(x)()=}

\DefEquation
{
d(\alpha\wedge\beta)=d\alpha\wedge\beta+(-1)^p\alpha\wedge d\beta
}
{d(a wedge b)=}

\DefEq
{
$S(p+q)$
}
{S p+q}

\DefEquation
{
\tau_i=
\begin{pmatrix}
a_i&a_0&...&\overbrace{a_i}&...&a_p
\end{pmatrix}
}
{ti a0p}

\DefEq
{
$\sigma\in S_1$,
}
{s in S1}

\DefEq
{
\[
f\circ(a_1, ..., a_n)=|\sigma|(f\circ\sigma\circ(a_1,...,a_n))
\]
}
{fa=sfsa}

\DefEq
{
\[
\omega(x)\circ (a_1,...,a_n)
=\omega(x)\circ(a_1^{\gii_1}e_{A\cdot\gii_1},...,a_n^{\gii_n}e_{A\cdot\gii_n})
=a_1^{\gii_1}...a_n^{\gii_n}(\omega(x)\circ(e_{A\cdot\gii_1},...,e_{A\cdot\gii_n}))
\]
}
{w(a1n)=a1n...}

\DefEq
{
$\tau_i\in S_1$, $i=0$, ..., $p$,
}
{ti in S1}

\DefEq
{
$\tau_i\in S_1$, $\mu\in S$
}
{ti in S1, s1 in S}

\DefEquation
{
\sum_{\sigma_1\in S_1}|\sigma_1|f(\sigma_1(a_0),...,\sigma_1(a_p))
=\sum_{i=0}^p(-1)^i\sum_{\sigma\in S}|\sigma|f(\tau_i(a_0),\sigma(\tau_i(a_1)),
...,\sigma(\tau_i(a_p)))
}
{sum s1 =sum +-s}

\DefEq
{
\begin{pmatrix}
\sigma(a_0)&\sigma(a_1)&...&\sigma(a_p)
\end{pmatrix}
=
\begin{pmatrix}
\tau_i(a_0)&\tau_i(\mu(a_1))&...&\tau_i(\mu(a_p))
\end{pmatrix}
}
{s=(ti mu)}

\DefEquation
{
\sigma=\tau_i
\begin{pmatrix}
a_i&a^0&...&\overbrace{a_i}&...&a_p\\
a_i&\mu(a^0)&...&&...&\mu(a_p)
\end{pmatrix}
}
{s=mu i ne 0}

\DefEquation
{
|\sigma|=(-1)^i|\mu|
}
{|s|=(-1)i|mu|}

\DefEq
{
\[
\begin{matrix}
a_0&...&\overbrace{a_i}&...&a_p
\end{matrix}
\]
$i=0$, ..., $p$,
}
{a0p - ai}

\DefEquation
{
\mu=
\begin{pmatrix}
\tau_i(a_1)&....&\tau_i(a_p)\\
\tau_i(\sigma(a_1))&....&\tau_i(\sigma(a_p))
\end{pmatrix}
}
{mu=()}

\DefEquation
{
\tau_i(a_0)=\sigma(a_0)
}
{ti 0 = s 0}

\DefEquation
{
|\tau_i|=(-1)^i
}
{|ti|=(-1)i}

\DefEquation
{
|\sigma|=|\tau_i|\,|\mu|
}
{|s|=|t mu|}

\DefEq
{
$a_{\Pn}$, ..., $a_p$\Pt
}
{a1p}

\DefEq
{
$a_{\Pn}$, ..., $a_{p+q}$.
}
{a1p+q}

\DefEq
{
$a_0$, ..., $a_{p+q}$.
}
{a0p+q}

\DefEq
{
\begin{align*}
I&=\{a_1,...,a_{p+q}\}\\
J&=\{J_1,...,J_p\}\subseteq I\\
K&=\{K_1,...,K_q\}=I\setminus J
\end{align*}
}
{IJK=()}

\DefEq
{
$R(J,I)$
}
{RJI}

\DefEq
{
$\omega_{\gii_1,...,\gii_p}(x)$
}
{o1p x}

\DefEquation
{
\frac{d^2 f}{d x^2}\circ (e_{\gii};e_{\gij})
=
\frac{d^2 f}{d x^2}\circ (e_{\gij};e_{\gii})
}
{df/dx2 ij=ji}

\DefEquation
{
\frac{\partial^2 f}{\partial x^{\gij}\partial x^{\gii}}=
\frac{\partial^2 f}{\partial x^{\gii}\partial x^{\gij}}
}
{df/dxij=df/dxji}

\DefEq
{
\[
\begin{matrix}
a=a^{\gii}e_{\gii}&b=b^{\gii}e_{\gii}
\end{matrix}
\]
}
{a= b= e}

\DefEquation
{
\frac{d^2 f}{d x^2}\circ(a,b)
=
\frac{d^2 f}{d x^2}\circ(b,a)
}
{df/dx2 ab=ba}

\DefEq
{
$i$, $1\le i<n$.
}
{1<=i<n}

\DefEquation
{
|\sigma_i|=-1
}
{|sigma i|}

\DefEquation
{
\sigma_i(a_k)=\left\{
\begin{matrix}
a_{i+1}&k=i\\
a_i&k=i+1\\
a_k&k\ne i,k\ne i+1
\end{matrix}
\right.
}
{sigma i}

\DefEq
{
\begin{align*}
\mathcal{LA}(D;A\rightarrow B)&=\mathcal L(D;A\rightarrow B)\\
\mathcal{LA}(D;A^0\rightarrow B)&=B
\end{align*}
}
{L(A1,A0,B)=}

\def\LAnAB{\mathcal{LA}(D;A^n\rightarrow B)}

\DefEq
{
\symb{\LAnAB}{module of skew symmetric polylinear maps}{1}
}
{module of skew symmetric polylinear maps}

\DefEq
{
\[
fg:\mathcal L(D;A^p\rightarrow B)\times\mathcal L(D;A^q\rightarrow C)\rightarrow \mathcal L(D;A^{p+q}\rightarrow C)
\]
}
{hpq:Lpq->Lp+q}

\DefEquation
{
(fg)\circ(a_1,...,a_{p+q})=(f\circ(a_1,...,a_p))(g\circ(a_{p+1},...,a_{p+q}))
}
{hpq(fg)=}

\DefEq
{
$fg$
}
{hpq(fg)}

\AddEq [3]{f in LAAB}
{
\[f_{#1}\in\mathcal{LA}(D;A^{#2}\rightarrow B_{#3})\]
}

\DefEquation
{
\begin{split}
&\,(f_1\circ(a_1,...,a_p))((f_2\circ(a_{p+1},...a_{p+q}))(f_3\circ(a_{p+q+1},...,a_{p+q+r})))
\\=&\,
((f_1\circ(a_1,...,a_p))(f_2\circ(a_{p+1},...a_{p+q})))(f_3\circ(a_{p+q+1},...,a_{p+q+r}))
\end{split}
}
{f1(f2f3)=(f1f2)f3}

\DefEquation
{
\begin{split}
&\,\sum_{\sigma\in S}|\sigma|(f_1\circ(\sigma(a_1),...,\sigma(a_p)))
((f_2f_3)\circ(\sigma(a_{p+1}),...,\sigma(a_{p+q+r})))
\\=&\,
\sum_{\lambda\in R_1}|\lambda|(f_1\circ(\lambda(a_1),...,\lambda(a_p)))
\sum_{\tau\in S_1(\lambda)}|\tau|
((f_2f_3)\circ(\tau(\lambda_1),...,\tau(\lambda_{q+r})))
\end{split}
}
{f1(f2f3)}

\DefEquation
{
\begin{split}
&\,\sum_{\sigma\in S}|\sigma|(f_1\circ(\sigma(a_1),...,\sigma(a_p)))
((f_2f_3)\circ(\sigma(a_{p+1}),...,\sigma(a_{p+q+r})))
\\=&\,
\sum_{\lambda\in R_1}|\lambda|(q+r)!(f_1\circ(\lambda(a_1),...,\lambda(a_p)))
((f_2f_3)\circ[\lambda_1,...,\lambda_{q+r}])
\\=&\,
\sum_{\lambda\in R_1}|\lambda|(q+r)!\frac{q!r!}{(q+r)!}
\\ * &\,
(f_1\circ(\lambda(a_1),...,\lambda(a_p)))
((f_2\wedge f_3)\circ(\lambda_1,...,\lambda_{q+r}))
\end{split}
}
{f1(f2f3) 2}

\AddEquation{f1(f2f3) 3}
{
\begin{split}
&\,\sum_{\sigma\in S}|\sigma|(f_1\circ(\sigma(a_1),...,\sigma(a_p)))
((f_2f_3)\circ(\sigma(a_{p+1}),...,\sigma(a_{p+q+r})))
\\=&\,
\sum_{\lambda\in R_1}|\lambda|(f_1\circ(\lambda(a_1),...,\lambda(a_p)))
\\ * &\,\frac{q!r!}{(q+r)!}
\sum_{\tau\in S_1(\lambda)}|\tau|
((f_2\wedge f_3)\circ(\tau(\lambda_1),...,\tau(\lambda_{q+r})))
\\=&\,\frac{q!r!}{(q+r)!}
\\ * &\,
\sum_{\sigma\in S}|\sigma|(f_1\circ(\sigma(a_1),...,\sigma(a_p)))
((f_2\wedge f_3)\circ(\sigma(a_{p+1}),...,\sigma(a_{p+q+r})))
\end{split}
}

\DefEquation
{
\begin{split}
&\,\sum_{\sigma\in S}|\sigma|(f_1\circ(\sigma(a_1),...,\sigma(a_p)))
((f_2f_3)\circ(\sigma(a_{p+1}),...,\sigma(a_{p+q+r})))
\\=&\,\frac{q!r!}{(q+r)!}
\sum_{\sigma\in S}|\sigma|(f_1(f_2\wedge f_3))
\circ(\sigma(a_1),...,\sigma(a_{p+q+r}))
\end{split}
}
{f1(f2f3) 4}

\DefEquation
{
\begin{split}
&\,\sum_{\sigma\in S}|\sigma|(f_1\circ(\sigma(a_1),...,\sigma(a_p)))
((f_2f_3)\circ(\sigma(a_{p+1}),...,\sigma(a_{p+q+r})))
\\=&\,(p+q+r)!\frac{q!r!}{(q+r)!}
(f_1(f_2\wedge f_3))
\circ[a_1,...,a_{p+q+r}]
\\=&\,(p+q+r)!\frac{q!r!}{(q+r)!}\frac{p!(q+r)!}{(p+q+r)!}
(f_1\wedge(f_2\wedge f_3))
\circ(a_1,...,a_{p+q+r})
\\=&\,p!q!r!
(f_1\wedge(f_2\wedge f_3))
\circ(a_1,...,a_{p+q+r})
\end{split}
}
{f1(f2f3) 5}

\DefEquation
{
\begin{split}
&\,\sum_{\sigma\in S}|\sigma|((f_1f_2)\circ(\sigma(a_1),...,\sigma(a_{p+q})))
(f_3\circ(\sigma(a_{p+q+1}),...,\sigma(a_{p+q+r})))
\\=&\,
\sum_{\lambda\in R_2}|\lambda|
\left(\sum_{\tau\in S_1(\lambda)}|\tau|
((f_1f_2)\circ(\tau(\lambda_1),...,\tau(\lambda_{p+q})))\right)
\\ * &\,
(f_3\circ(\lambda(a_{p+q+1}),...,\lambda(a_{p+q+r})))
\end{split}
}
{(f1f2)f3}

\DefEquation
{
\begin{split}
&\,\sum_{\sigma\in S}|\sigma|((f_1f_2)\circ(\sigma(a_1),...,\sigma(a_{p+q})))
(f_3\circ(\sigma(a_{p+q+1}),...,\sigma(a_{p+q+r})))
\\=&\,
\sum_{\lambda\in R_2}|\lambda|(p+q)!
(((f_1f_2)\circ[\lambda_1,...,\lambda_{p+q}]))
\\ * &\,
(f_3\circ(\lambda(a_{p+q+1}),...,\lambda(a_{p+q+r})))
\\=&\,
\sum_{\lambda\in R_2}|\lambda|(p+q)!\frac{p!q!}{(p+q)!}
\\ * &\,
((f_1\wedge f_2)\circ(\lambda_1,...,\lambda_{p+q})))
(f_3\circ(\lambda(a_{p+q+1}),...,\lambda(a_{p+q+r})))
\end{split}
}
{(f1f2)f3 2}

\AddEquation{(f1f2)f3 3}
{
\begin{split}
&\,\sum_{\sigma\in S}|\sigma|((f_1f_2)\circ(\sigma(a_1),...,\sigma(a_{p+q})))
(f_3\circ(\sigma(a_{p+q+1}),...,\sigma(a_{p+q+r})))
\\=&\,
\sum_{\lambda\in R_2}|\lambda|\frac{p!q!}{(p+q)!}
\left(\sum_{\tau\in S_1(\lambda)}|\tau|
((f_1\wedge f_2)\circ(\tau(\lambda_1),...,\tau(\lambda_{p+q})))\right)
\\ * &\,
(f_3\circ(\lambda(a_{p+q+1}),...,\lambda(a_{p+q+r})))
\\=&\,\frac{p!q!}{(p+q)!}
\sum_{\sigma\in S}|\sigma|
((f_1\wedge f_2)\circ(\sigma(a_1),...,\sigma(a_{p+q}))))
\\ * &\,
(f_3\circ(\sigma(a_{p+q+1}),...,\sigma(a_{p+q+r})))
\end{split}
}

\DefEquation
{
\begin{split}
&\,\sum_{\sigma\in S}|\sigma|((f_1f_2)\circ(\sigma(a_1),...,\sigma(a_{p+q})))
(f_3\circ(\sigma(a_{p+q+1}),...,\sigma(a_{p+q+r})))
\\=&\,\frac{p!q!}{(p+q)!}
\sum_{\sigma\in S}|\sigma|((f_1\wedge f_2)f_3)
\circ(\sigma(a_1),...,\sigma(a_{p+q+r})))
\end{split}
}
{(f1f2)f3 4}

\DefEquation
{
\begin{split}
&\,\sum_{\sigma\in S}|\sigma|((f_1f_2)\circ(\sigma(a_1),...,\sigma(a_{p+q})))
(f_3\circ(\sigma(a_{p+q+1}),...,\sigma(a_{p+q+r})))
\\=&\,
(p+q+r)!\frac{p!q!}{(p+q)!}((f_1\wedge f_2)f_3)
\circ[a_1,...,a_{p+q+r}]
\\=&\,
(p+q+r)!\frac{p!q!}{(p+q)!}\frac{(p+q)!r!}{(p+q+r)!}
((f_1\wedge f_2)\wedge f_3)
\circ(a_1,...,a_{p+q+r})
\\=&\,
p!q!r!
((f_1\wedge f_2)\wedge f_3)
\circ(a_1,...,a_{p+q+r})
\end{split}
}
{(f1f2)f3 5}

\DefEq
{
$I=\{a_0,...,a_p\}$, $J=\{a_0\}$.
}
{I=a0p,J=a0}

\DefEquation
{
\frac{d(\alpha\wedge\beta)}{d x}\circ a_0=
\left(\frac{d\alpha}{d x}\circ a_0\right)\wedge\beta+
\alpha\wedge\left(\frac{d\beta}{d x}\circ a_0\right)
}
{d(a A b)/dx=...}

\DefEquation
{
\begin{split}
&\,
\left(\frac{d(\alpha\wedge\beta)}{d x}\circ a_0\right)
\circ(a_1,...,a_{p+q})
\\=&\,
\left(\left(\frac{d\alpha}{dd x}\circ a_0\right)\wedge\beta\right)
\circ(a_1,...,a_{p+q})
\\+&\,
\left(\alpha\wedge\left(\frac{d\beta}{d x}\circ a_0\right)\right)
\circ(a_1,...,a_{p+q})
\end{split}
}
{d(a A b)/dx=..., 1}

\DefEquation
{
\begin{split}
&\,d(\alpha\wedge\beta)\circ(a_0,...,a_{p+q})\\=&\,
\sum_{i=0}^{p+q}(-1)^i
\left(\frac{d(\alpha\wedge\beta)}{d x}\circ \tau_i(a_0)\right)
\circ(\tau_i(a_1),...,\tau_i(a_{p+q}))
\\=&\,
\sum_{i=0}^{p+q}(-1)^i
\left(\left(\frac{d\alpha}{d x}\circ \tau_i(a_0)\right)\wedge\beta\right)
\circ(\tau_i(a_1),...,\tau_i(a_{p+q}))
\\+&\,
\sum_{i=0}^{p+q}(-1)^i
\left(\alpha\wedge\left(\frac{d\beta}{d x}\circ \tau_i(a_0)\right)\right)
\circ(\tau_i(a_1),...,\tau_i(a_{p+q}))
\end{split}
}
{d(a A b)/dx=..., 2}

\DefEquation
{
\begin{split}
&\,
\sum_{i=0}^{p+q}(-1)^i
\left(\left(\frac{d\alpha}{d x}\circ \tau_i(a_0)\right)\wedge\beta\right)
\circ(\tau_i(a_1),...,\tau_i(a_{p+q}))
\\=&\,
\sum_{i=0}^{p+q}\frac{(-1)^i}{p!q!}
\sum_{\sigma\in S}|\sigma|
\left(\left(\frac{d\alpha}{d x}\circ \tau_i(a_0)\right)
\circ(\sigma(\tau_i(a_1)),...,\sigma(\tau_i(a_p)))\right)
\\ * &\,
(\beta\circ(\sigma(\tau_i(a_{p+1})),...,\sigma(\tau_i(a_{p+q}))))
\\=&\,
\sum_{i=0}^{p+q}\frac{(-1)^i}{p!q!}
\sum_{\sigma\in S}|\sigma|
\left(\frac{d\alpha}{dx}
\circ(\tau_i(a_0),\sigma(\tau_i(a_1)),...,\sigma(\tau_i(a_p)))\right)
\\ * &\,
(\beta\circ(\sigma(\tau_i(a_{p+1})),...,\sigma(\tau_i(a_{p+q}))))
\end{split}
}
{d(a A b)/dx=..., 2 1.1}

\DefEquation
{
\begin{split}
&\,
\sum_{i=0}^{p+q}(-1)^i
\left(\left(\frac{d\alpha}{dx}\circ \tau_i(a_0)\right)\wedge\beta\right)
\circ(\tau_i(a_1),...,\tau_i(a_{p+q}))
\\=&\,
\frac 1{p!q!}\sum_{\sigma_1\in S_1}|\sigma_1|
\left(\frac{d\alpha}{dx}
\circ(\sigma_1(a_0),...,\sigma_1(a_p))\right)
(\beta\circ(\sigma_1(a_{p+1}),...,\sigma_1(a_{p+q})))
\end{split}
}
{d(a A b)/dx=..., 2 1.2}

\AddEquation{d(a A b)/dx=..., 2 1.3}
{
\begin{split}
&\,
\sum_{i=0}^{p+q}(-1)^i
\left(\left(\frac{d\alpha}{dx}\circ \tau_i(a_0)\right)\wedge\beta\right)
\circ(\tau_i(a_1),...,\tau_i(a_{p+q}))
\\=&\,
\frac 1{p!q!}\sum_{\lambda\in R_1}|\lambda|
\left(
\sum_{\tau\in S_1(\lambda)}|\tau|
\frac{d\alpha}{dx}
\circ(\tau(\lambda_0),...,\tau(\lambda_p))
\right)
\\ * &\,
(\beta\circ(\lambda(a_{p+1}),...,\lambda(a_{p+q})))
\end{split}
}

\DefEquation
{
\begin{split}
&\,
\sum_{i=0}^{p+q}(-1)^i
\left(\left(\frac{d\alpha}{dx}\circ \tau_i(a_0)\right)\wedge\beta\right)
\circ(\tau_i(a_1),...,\tau_i(a_{p+q}))
\\=&\,
\frac 1{p!q!}\sum_{\lambda\in R_1}|\lambda|
p!
d\alpha
\circ(\lambda_0,...,\lambda_p)
(\beta\circ(\lambda(a_{p+1}),...,\lambda(a_{p+q})))
\\=&\,
\frac 1{q!}\sum_{\lambda\in R_1}|\lambda|
(d\alpha
\circ(\lambda_0,...,\lambda_p))
(\beta\circ(\lambda(a_{p+1}),...,\lambda(a_{p+q})))
\end{split}
}
{d(a A b)/dx=..., 2 1.4}

\DefEquation
{
\begin{split}
&\,
\sum_{i=0}^{p+q}(-1)^i
\left(\left(\frac{d\alpha}{dx}\circ \tau_i(a_0)\right)\wedge\beta\right)
\circ(\tau_i(a_1),...,\tau_i(a_{p+q}))
\\=&\,
(d\alpha\wedge\beta)
\circ(a_0,...,a_{p+q})
\end{split}
}
{d(a A b)/dx=..., 2 1.5}

\DefEquation
{
\begin{split}
&\,
\sum_{i=0}^{p+q}(-1)^i
\left(\alpha\wedge\left(\frac{d\beta}{dx}\circ \tau_i(a_0)\right)\right)
\circ(\tau_i(a_1),...,\tau_i(a_{p+q}))
\\=&\,
\sum_{i=0}^{p+q}\frac{(-1)^i}{p!q!}
\sum_{\sigma\in S}|\sigma|
(\alpha\circ(\sigma(\tau_i(a_1)),...,\sigma(\tau_i(a_p))))
\\ * &\,
\left(\left(\frac{d\beta}{dx}\circ \tau_i(a_0)\right)
\circ(\sigma(\tau_i(a_{p+1})),...,\sigma(\tau_i(a_{p+q})))\right)
\\=&\,
\sum_{i=0}^{p+q}\frac{(-1)^i}{p!q!}
\sum_{\sigma\in S}|\sigma|
(\alpha\circ(\sigma(\tau_i(a_1)),...,\sigma(\tau_i(a_p))))
\\ * &\,
\left(\frac{d\beta}{dx}
\circ(\tau_i(a_0),\sigma(\tau_i(a_{p+1})),...,\sigma(\tau_i(a_{p+q})))\right)
\end{split}
}
{d(a A b)/dx=..., 2 2.1}

\DefEquation
{
\begin{split}
&\,
\sum_{i=0}^{p+q}(-1)^i
\left(\alpha\wedge\left(\frac{d\beta}{dx}\circ \tau_i(a_0)\right)\right)
\circ(\tau_i(a_1),...,\tau_i(a_{p+q}))
\\=&\,
\frac 1{p!q!}\sum_{\sigma_1\in S_1}|\sigma_1|
(\alpha\circ(\sigma_1(a_1),...,\sigma_1(a_p)))
\\ * &\,
\left(\frac{d\beta}{dx}
\circ(\sigma_1(a_0),\sigma_1(a_{p+1}),...,\sigma_1(a_{p+q}))\right)
\end{split}
}
{d(a A b)/dx=..., 2 2.2}

\DefEquation
{
\begin{split}
&\,
\sum_{i=0}^{p+q}(-1)^i
\left(\alpha\wedge\left(\frac{d\beta}{dx}\circ \tau_i(a_0)\right)\right)
\circ(\tau_i(a_1),...,\tau_i(a_{p+q}))
\\=&\,
\frac 1{p!q!}\sum_{\lambda\in R_1}|\lambda|
(\alpha\circ(\lambda(a_1),...,\lambda(a_p))
\\ * &\,
\left(
\sum_{\tau\in S_1(\lambda)}|\tau|
\frac{d\beta}{dx}
\circ(\tau(\lambda_0),\tau(\lambda_1),...,\tau(\lambda_q))
\right)
\end{split}
}
{d(a A b)/dx=..., 2 2.3}

\DefEquation
{
\begin{split}
&\,
\sum_{i=0}^{p+q}(-1)^i
\left(\alpha\wedge\left(\frac{d\beta}{dx}\circ \tau_i(a_0)\right)\right)
\circ(\tau_i(a_1),...,\tau_i(a_{p+q}))
\\=&\,
\frac 1{p!q!}\sum_{\lambda\in R_1}|\lambda|
q!
(\alpha
\circ(\lambda(a_1),...,\lambda(a_p))
(d\beta\circ(\lambda_0,...,\lambda_q)))
\\=&\,
\frac 1{p!}\sum_{\lambda\in R_1}|\lambda|
(\alpha
\circ(\lambda(a_1),...,\lambda(a_p))
(d\beta\circ(\lambda_0,...,\lambda_q)))
\end{split}
}
{d(a A b)/dx=..., 2 2.4}

\DefEquation
{
\begin{split}
&\,
\sum_{i=0}^{p+q}(-1)^i
\left(\alpha\wedge\left(\frac{d\beta}{dx}\circ \tau_i(a_0)\right)\right)
\circ(\tau_i(a_1),...,\tau_i(a_{p+q}))
\\=&\,
(\alpha\wedge d\beta)
\circ(a_1,...,a_p,a_0,a_{p+1},...,a_{p+q})
\end{split}
}
{d(a A b)/dx=..., 2 2.5}

\DefEquation
{
\begin{split}
&\,
\sum_{i=0}^{p+q}(-1)^i
\left(\alpha\wedge\left(\frac{d\beta}{dx}\circ \tau_i(a_0)\right)\right)
\circ(\tau_i(a_1),...,\tau_i(a_{p+q}))
\\=&\,(-1)^p
(\alpha\wedge d\beta)
\circ(a_0,...,a_{p+q})
\end{split}
}
{d(a A b)/dx=..., 2 2.6}

\DefEq
{
\[
\left|
\begin{pmatrix}
a_1&a_2&...&a_p&a_0&a_{p+1}&...&a_{p+q}\\
a_0&a_1&...&a_{p-1}&a_p&a_{p+1}&...&a_{p+q}
\end{pmatrix}
\right|
=(-1)^p
\]
}
{|a1p a0 apq|=}

\DefEquation
{
\begin{split}
&\,d(\alpha\wedge\beta)\circ(a_0,...,a_{p+q})\\=&\,
(d\alpha\wedge\beta)
\circ(a_0,...,a_{p+q})
+(-1)^p
(\alpha\wedge d\beta)
\circ(a_0,...,a_{p+q})
\end{split}
}
{d(a A b)/dx=..., 3}

\DefEquation
{
(f_1\wedge f_2)\wedge f_3=f_1\wedge (f_2\wedge f_3)
}
{f1(f2f3)=(f1f2)f3 wedge}

\DefEq
{
$J\rightarrow I$.
}
{J->I}

\DefEq
{
\[I_{1.p}\rightarrow I_{1.p+q}\]
}
{I1p->I1pq}

\DefEq
{
\[I_{p+1.p+q}\rightarrow I_{1.p+q}\]
}
{Ip1q->I1pq}

\DefEq
{
$R_2=R(I_{p+q+1.p+q+r},I_{1.p+q+r})$.
}
{R2=Ipqpqr->I1pqr}

\DefEq
{
$\lambda\in R(J,I)$\Pt
}
{l in RJI}

\DefEq
{
$\tau\in S_1(\lambda)$\Pt
}
{t in S1l}

\DefEq
{
$\tau\in S_2(\nu)$\Pt
}
{t in S2n}

\DefEq
{
$S_1(\lambda)$
}
{S1(l)}

\DefEq
{
\[f\wedge f=0\]
}
{fAf=0}

\DefEq
{
\[f:B\rightarrow B]
}
{f:B->B}

\DefEquation
{
(f\wedge f)\circ(a,b)=f(a)f(b)-f(b)f(a)=[f(a),f(b)]
}
{fAf=}

\DefEquation
{
f\wedge f\ne 0
}
{fAf ne 0}

\DefEq
{
$\sigma\in S$\Pt
}
{s in S}

\DefEquation
{
\mu(\lambda)(K_k)=\lambda_k
}
{mu(l)()=}

\DefEq
{
$\displaystyle\frac{\partial^2f(x)}{\partial x^{\gii}\partial x^{\gij}}$
}
{d2f/dxij}

\DefEq
{
$\displaystyle\frac{d^2f(x)}{dx^2}$
}
{d2f/dx2}

\DefEq
{
\[
\mu(\lambda):K\rightarrow D_1(\lambda)
\]
}
{K->D1}

\DefEquation
{
\lambda=
\begin{pmatrix}
J_1&...&J_p\\
\sigma(J_1)&...&\sigma(J_p)
\end{pmatrix}
}
{lambda=}

\DefEq
{
\[\tau:D_1(\lambda)\rightarrow D_1(\lambda)\]
}
{t:D1->D1}

\DefEquation
{
\tau(\lambda_k)=\sigma(a_{K.k})
}
{t(l)=s(a)}

\DefEquation
{
\tau(\mu(\lambda)(a_{K.k}))=\sigma(a_{K.k})
}
{t(m(a))=s(a)}

\DefEq
{
\sigma=
\begin{pmatrix}
J_1&...&J_p&K_1&...&K_q\\
\lambda(J_1)&...&\lambda(J_p)&\tau(\mu(\lambda)(K_1))&...&\tau(\mu(\lambda)(K_q))
\end{pmatrix}
}
{sigma=lambda,tau}

\DefEquation
{
D_1(\lambda)=I\setminus D(\lambda)
=\{\lambda_1,...,\lambda_q\}
}
{D1=I-D}

\DefEq
{
$(\lambda,\tau)$
}
{(lambda,tau)}

\DefEq
{
\[
\sigma=(\lambda,\tau)
\]
}
{s=(lambda,tau)}

\DefEquation
{
|(\lambda,\tau)|=|\lambda||\tau|
}
{|lt|=|l||t|}

\DefEq
{
$(\lambda,\delta)$\Pt
}
{(l,d)}

\DefEq
{
$|\delta|=1$,
}
{|d|=1}

\DefEq
{
\[
|\lambda|=|(\lambda,\delta)|
\]
}
{|l|=|l,d|}

\DefEq
{
$\lambda$, $\mu(\lambda)$, $\tau$
}
{l m(l) t}

\DefEq
{
\[D_1(\lambda)\cup D(\lambda)=I\]
}
{D1+D=I}

\DefEq
{
$D_1(\lambda)$.
}
{D1(l)}

\DefEq
{
$D(\lambda)$
}
{D(l)}

\DefEq
{
\[
I_{m.n}=\{a_i:m\le i\le n\}
\]
}
{Imn=}

\DefEq
{
$I_{1.p+q+r}$.
}
{I1pqr}

\DefEq
{
$I_{1.p+q}$.
}
{I1pq}

\DefEquation
{
\begin{split}
&\,(f_1\wedge f_2)\circ(a_1,...,a_{p+q})
\\=&\,
\frac 1{p!q!}\sum_{\sigma\in S}|\sigma|
(f_1\circ(\sigma(a_1),...,\sigma(a_p)))(f_2\circ(\sigma(a_{p+1}),...,\sigma(a_{p+q})))
\end{split}
}
{f1 A f2 =}

\DefEquation
{
\begin{split}
&\,(f_1\wedge f_2)\circ(a_1,...,a_{p+q})
\\=&\,
\frac 1{p!q!}\sum_{\lambda\in R_1}|\lambda|
(f_1\circ(\lambda(a_1),...,\lambda(a_p)))
\sum_{\tau\in S_1(\lambda)}|\tau|
(f_2\circ(\tau(\lambda_1),...,\tau(\lambda_{q})))
\end{split}
}
{f1 A f2 = 1}

\DefEq
{
$\lambda\in R_0$\Pt
}
{l in R0}

\DefEq
{
$\mu\in R_2(\lambda)$
}
{m in R2}

\DefEquation
{
R_1=\bigcup_{\lambda\in R_0}R_2(\lambda)
}
{R1=U R2}

\DefEquation
{
|\lambda|(f_1\circ(\lambda(a_1),...,\lambda(a_p)))=
|\mu|(f_1\circ(\mu(a_1),...,\mu(a_p)))
}
{|l|f1=|m|f1}

\DefEquation
{
\begin{split}
&\,(f_1\wedge f_2)\circ(a_1,...,a_{p+q})
\\=&\,
\frac 1{p!q!}\sum_{\lambda\in R_1}|\lambda|
(f_1\circ(\lambda(a_1),...,\lambda(a_p)))
q!
(f_2\circ(\lambda_1,...,\lambda_{q}))
\end{split}
}
{f1 A f2 = 2}

\DefEquation
{
\begin{split}
&\,(f_1\wedge f_2)\circ(a_1,...,a_{p+q})
\\=&\,
\frac 1{p!}p!\sum_{\lambda\in R_0}|\lambda|
(f_1\circ(\lambda(a_1),...,\lambda(a_p)))
(f_2\circ(\lambda_1,...,\lambda_{q}))
\end{split}
}
{f1 A f2 = 3}

\DefEq
{
\[
R_2(\lambda)=\{\mu\in R_1:D(\mu)=D(\lambda)\}
\]
}
{R2(l)=...}

\DefEquation
{
f_2\circ(\lambda_1,...,\lambda_q)=|\tau|
(f_2\circ(\tau(\lambda_1),...,\tau(\lambda_{q})))
}
{f2 l = |s| f2 sl}

\DefEq
{
$\sigma\in S$.
}
{s in SI}

\DefEq
{
\[R_1=R(I_{1.p},I_{1.p+q+r})\]
}
{R1=I1p->I1pqr}

\DefEq
{
\[R_1=R(I_{1.p},I_{1.p+q})\]
}
{R1=I1p->I1pq}

\DefEq
{
\[R_1=R(I_{p+1.p+q},I_{1.p+q})\]
}
{R1=Ip1q->I1pq}

\DefEq
{
$R_1=R(I_{p+1.p+q},I_{0.p+q})$.
}
{R1=Ip+1q>I0p+q}

\DefEq
{
$R_1=R(I_{1.p},I_{0.p+q})$.
}
{R1=I1p->I0p+q}

\DefEquation
{
(f_1\wedge f_2)\circ(a_1,...,a_{p+q})=
\sum_{\lambda\in R_0}|\lambda|(f_1\circ(\lambda(a_1),..,\lambda(a_p)))
(f_2\circ(\lambda_1,...,\lambda_q))
}
{f1 f2 = sum l in R0p}

\DefEquation
{
(f_1\wedge f_2)\circ(a_1,...,a_{p+q})
=
\frac 1{p!}\sum_{\lambda\in R_1}|\lambda|
(f_1\circ(\lambda(a_1),...,\lambda(a_p)))
(f_2\circ(\lambda_1,...,\lambda_{q}))
}
{f1 f2 = sum l in R1p}

\DefEquation
{
\begin{split}
&\,(f_1\wedge f_2)\circ(a_1,...,a_{p+q})
\\=&\,
\sum_{\lambda\in R_0}|\lambda|(f_1\circ(\lambda_1,...,\lambda_p))
(f_2\circ(\lambda(a_{p+1}),..,\lambda(a_{p+q})))
\end{split}
}
{f1 f2 = sum l in R0q}

\DefEquation
{
\begin{split}
&\,(f_1\wedge f_2)\circ(a_1,...,a_{p+q})
\\=&\,
\frac 1{q!}\sum_{\lambda\in R_1}|\lambda|
(f_1\circ(\lambda_1,...,\lambda_p))
(f_2\circ(\lambda(a_{p+1}),...,\lambda(a_{p+q})))
\end{split}
}
{f1 f2 = sum l in R1q}

\DefEq
{
\[
\begin{pmatrix}
a_1&...&a_p&a_{p+1}&...&a_{p+q}
\\
\lambda(a_1)&...&\lambda(a_p)&\lambda_1&...&\lambda_q
\end{pmatrix}
\]
}
{apq->l}

\DefEquation
{
\begin{split}
&\,(f_1\circ(a_1,...,a_p))((f_2f_3)\circ(a_{p+1},...,a_{p+q+r}))
\\=&\,
((f_1f_2)\circ(a_1,...,a_{p+q}))(f_3\circ(a_{p+q+1},...,a_{p+q+r}))
\end{split}
}
{f1(f2f3)=(f1f2)f3 1}

\DefEquation
{
\begin{split}
&\,\sum_{\sigma\in S}|\sigma|(f_1\circ(\sigma(a_1),...,\sigma(a_p)))
((f_2f_3)\circ(\sigma(a_{p+1}),...,\sigma(a_{p+q+r})))
\\=&\,
\sum_{\sigma\in S}|\sigma|((f_1f_2)\circ(\sigma(a_1),...,\sigma(a_{p+q})))
(f_3\circ(\sigma(a_{p+q+1}),...,\sigma(a_{p+q+r})))
\end{split}
}
{f1(f2f3)=(f1f2)f3 2}

\DefEq
{
$f_3(a_{p+q+1},...,a_{p+q+r})$.
}
{f3()}

\DefEq
{
$f_1(a_1,...,a_p)$, $f_1(a_{p+1},...a_{p+q})$
}
{f1()f2()}

\DefEq
{
\symb{f\wedge g}{exterior product}{1}
}
{exterior product}

\DefEq
{
\[
x\in U\rightarrow \alpha(x)\wedge\beta(x)\in\mathcal {LA}(D;A^{p+q}\rightarrow B_2)
\]
}
{x->a(x) wedge b(x)}

\DefEq
{
\[
\mathcal{LA}(D;A^p\rightarrow B_1)\times\mathcal {LA}(D;A^q\rightarrow B_2)\rightarrow 
\mathcal {LA}(D;A^{p+q}\rightarrow B_2)
\]
}
{LA1xLA2->LA}

\DefEq
{
\[
\Omega^{(n)}_p(U,B_1)\times\Omega^{(n)}_q(U,\PA)
\rightarrow \Omega^{(n)}_{p+q}(U,\PA)
\]
}
{O(B1)xO(B2)->O(B2)}

\DefEq
{
$\Omega^{(n)}(U,B_1)$\Pt
}
{OnB}

\DefEq
{
\[
x\rightarrow (\alpha(x),\beta(x))
\]
}
{x->(a(x),b(x))}

\DefEquation
{
(\ShowSymbol{exterior product}{1})(x)=\alpha(x)\wedge\beta(x)
}
{a wedge b x =}

\DefEq
{
\symb{\alpha\wedge\beta}{exterior product}{1}
}
{exterior product, differential forms}

\DefEquation
{
(\ShowSymbol{exterior product}{1})\circ(a_1,...,a_{p+q})=\frac {(p+q)!}{p!q!}
[fg]\circ(a_1,...,a_{p+q})
}
{exterior product =}

\DefEquation
{
(f\wedge g)\circ(a_1,...,a_{p+q})=
\frac 1{p!q!}\sum_{\sigma\in S(p+q)}|\sigma|((f\wedge g)\circ\sigma\circ(a_1,...,a_{p+q}))
}
{exterior product = skew symmetric}

\DefEq
{
\[
\xymatrix
{
h:B_1\ar[r]|{*}&B_2
}
\]
}
{h:B1->*B2}

\DefEq
{
$g\circ(a_{p+1},...,a_{p+q})$\Pt
}
{g()}

\DefEq
{
$f\circ(a_1,...,a_p)$\Pt
}
{f()}

\DefEq
{
\[f_k(...,a_i,a_i,...)=0\]
}
{fxixi=0}

\DefEq
{
\[
f(...,a_i,a_i,...)=|\sigma_i|f(...,a_i,a_i,...)=-f(...,a_i,a_i,...)=0
\]
}
{fxii}

\DefEquation
{
f(...,a_i,a_{i+1},...)=-f(...,a_{i+1},a_i,...)
=|\sigma_i|f(\sigma_i(a_1),...,\sigma_i(a_n))
}
{fxii1=-fxi1i}

\DefEq
{
\[
f\circ(a_1,...,a_n)=0
\]
}
{fx1n=0}

\DefEq
{
\[x=x^{\gii}e_{\gii}\]
}
{x=x1n}

\DefEquation
{
\omega_{\gii_1,...,\gii_p}(x)=\omega(x)\circ(e_{A\cdot\gii_1},...,e_{A\cdot\gii_p})
}
{w1p=...}

\DefEquation
{
\begin{matrix}
a_1=a_1^{\gii}e_{A\cdot\gii}&...&a_p=a_p^{\gii}e_{A\cdot\gii}
\end{matrix}
}
{a1p=...}

\DefEquation
{
\frac{\partial f}{\partial x^{\gii}}dx^{\gii}
=\frac{d f}{d x}\circ (dx^{\gii}e_{\gii})
=dx^{\gii}\left(\frac{d f}{d x}\circ e_{\gii}\right)
}
{df/dx o dx=df/dxi dxi}

\DefEquation
{
df=\frac{\partial f}{\partial x^{\gii}}dx^{\gii}
}
{df=df/dxi dxi}

\DefEq
{
dx=dx^{\gii}e_{\gii} \ \ \ dx^{\gii}\in D
}
{dx=dx e}

\DefEquation
{
df=\frac{d f(x)}{d x}\circ dx
}
{df=df/dx o dx}

\DefEq
{
$\displaystyle\frac{\partial f}{\partial x^{\gii}}$,
$\gii=\gi 1$, ..., $\gi n$,
}
{df/dxi 1n}

\DefEq
{
\[f(x)=f(x^{\gi 1},...,x^{\gi n})\]
}
{f(x)=f(x1n)}

\DefEq
{
\frac{\partial f}{\partial x^{\gii}}=
\frac{d f}{d x}\circ e_{\gii}
}
{df/dxi=}

%% file: Diff.Form.1.Parm.tex

\DefEq
{%
\def\Pform{\alpha}%
\def\Pp{p}%
\def\PA{B_1}%
}
{form=alpha}

\DefEq
{%
\def\Pform{\beta}%
\def\Pp{q}%
\def\PA{B_2}%
}
{form=beta}

\DefEq
{%
\def\Pform{\omega}%
\def\Pp{p}%
\def\PA{B}%
}
{form=omega}

\DefEq
{%
\def\Pform{\omega}%
\def\Pp{2}%
\def\PA{B}%
}
{form=omega p2}

\DefEq
{%
\def\Pp{2}%
}
{p=2}

\DefEq
{
\def\PA{B_1}%
}
{A=B1}

\DefEq
{
\def\PA{B_2}%
}
{A=B2}

%% file: Diff.Form.2.English.tex

\input{Diff.Form.2.Eq}

\Chapter{Structure of Differential Form}

\Section{Polylinear Map into Associative \texorpdfstring{$D$}{D}-Algebra}

\begin{theorem}
Let $\Basis F$ be the basis of left \BoxB{B}module
\ShowEq{L(A;B)}DAB.
We identify the polylinear map
\ShowEq{f in L(A->B)}{A^n}B{}
and the tensor
\ShowEq{b in Bon}
using the equality
\ShowEq{f circ=b circ n}
\end{theorem}
\begin{proof}
The proof of the theorem is similar to the proof of the theorem
\RefTheorem{representation of algebra An in LAnA}.
\end{proof}

\begin{theorem}
Let
\ShowEq{f=(b sigma)oI s}
be polylinear map. Than
\ShowEq{[f]=(b sigma)oI s}
\end{theorem}
\begin{proof}
The equality
\ShowEq{[f]=(b sigma)oI s 1}
follows from equalities
\EqRef{alternation of polylinear map =},
\EqRef{f=(b sigma)oI s}.
The equality
\EqRef{[f]=(b sigma)oI s}
follows from the equality
\EqRef{[f]=(b sigma)oI s 1}
because, for any $s$, we can use sum over
\EqParm{s in S}{=c}
\EqParm{s=s2 o ss}{=c}
instead of sum over
\EqParm{s2 in S}{=.}
Then
\ShowEq{|s2|=|s||ss|}
\end{proof}

Tensor representation of skew symmetric polylinear map is cumbersome.
So, when it is possible, we will use expession
\ShowEq{b o I}
to represent skew symmetric polylinear map.

\ePrints{1610.309618526,CACAA.06}
\ifx\Semafor\ValueOn
\ShowDefinition{otimes -}

See also the definition
\Ref[1302.7204]{definition: otimes -}
and the theorem
\RefTheorem[1302.7204]{Am->A(n+m)}
how we use the operation
\EqParm{otimes -=}{=z}
to define product of polynomials.
\fi

\begin{theorem}
Let
\EquationParm{f in LAAB =}{nn=1}
\EquationParm{f in LAAB =}{nn=2}
Then
\ShowEq{b1 A b2}
\end{theorem}
\begin{proof}
The equality
\ShowEq{f1of2o}
follows from equalities
\EqRefParm{f in LAAB =}{nn=1},
\EqRefParm{f in LAAB =}{nn=2}.
For every $s$, $t$,
the map
\ShowEq{sigma st}
is a permutation.
The equality
\ShowEq{f1of2o 1}
follows from equalities
\EqRef{f1of2o},
\EqRef{sigma st}.
The equality
\ShowEq{f1of2o 2}
follows from equalities
\EqRef{otimes -, 1},
\EqRef{f1of2o 1}.
The equality
\EqRef{b1 A b2}
follows from equalities
\EqRef{exterior product =},
\EqRef{b o I},
\EqRef{f1of2o 2}.
\end{proof}

\begin{theorem}
Let $\Basis e$ be the basis of finite dimensional $D$\Hyph algebra $B$.
Let
\EquationParm{f in LAAB =}{nn=}
Let
\ShowEq{b=b...e}
be standard representation of the tensor $b_s$.
Then
\ShowEq{f=b...e}
\end{theorem}
\begin{proof}
The equality
\EqRef{f=b...e}
follows from equalities
\EqRefParm{f in LAAB =}{nn=},
\EqRef{b=b...e}.
\end{proof}

\Section{Differential Form with Values in Associative \texorpdfstring{$D$}{D}-Algebra}

\begin{theorem}
\labelTheorem{tensor is continuous iff}
Let $A$, $B$ be associative finite dimensional Banach $D$\Hyph algebras.
Let
\ShowEq{U in A}
be open set.
The map
\DrawEq{f:U->Bp}{continuous}
is continuous iff
the map $f$ has form
\ShowEq{f=sum f}
where maps
\ShowEq{fis:U->B}
are continuous.
\end{theorem}
\begin{proof}
Direct proof of the theorem is not simple.
For instance, let
\ShowEq{f01:U->B}
be continuous maps.
Then the tensor
\ShowEq{f=f0of1}
is continuous map
\ShowEq{f:U->B2}
However, in the representation of the tensor $f(x)$
\ShowEq{f=f0of1 1}
each summand is not
continuous map.

So, not every representation of continuous map
\eqRef{f:U->Bp}{continuous}
has continuous summands.
However we have a tool
which allows us to prove the theorem.

Let $\Basis e_B$ be the basis of finite dimensional $D$\Hyph algebra $B$.
Then standard representation of the tensor $f(x)$ has the form\,\footnote{
We assume that the basis $\Basis e_B$ does not depend on $x$.
However this statement can be relaxed,
if we consider fibered $D$\Hyph algebra $B$
over $D$\Hyph algebra $A$
(I considered the definition of fibered universal algebra in
\citeBib{0702.561}).
I considered connection
of manifold over Banach algebra
in the section
\RefSection[0906.0135]{Manifolds with D-Affine Connections}.
This construction is not important for the proof of the theorem
\RefTheorem{tensor is continuous iff}
and I will consider it in future papers.
}
\ShowEq{f(x)=...xe}

\begin{lemma}
\labelLemma{standard components of tensor are continuous maps}
{\it
The map $f$ is continuous iff
all standard components of the tensor $f(x)$
\ShowEq{fi:U->D}
are continuous maps.
}
\end{lemma}

{\sc Proof.}
Since all maps
\ShowEq{fi1p}
are continuous,
then for any
\EqParm{e>0}{=c}
there exists
\ShowEq{d>0}
such that the statement
\ShowEq{|x-x1|<d}
implies
\ShowEq{|fix-fix1|<e}
for any map
\ShowEq{fi1p}
(this statement is valid because the set of maps
\ShowEq{fi1p}
is finite).
The equality
\ShowEq{|fx-fx1|<e 1}
follows from the equality
\EqRef{|fix-fix1|<e}.

We define the norm of tensor product
\ShowEq{xei}
according to the definition
\RefDefinition{norm of polylinear map}
\ShowEq{|xei|}
Since the set of tensor products
\ShowEq{xei}
is finite,
then we set
\ShowEq{E=max|xei|}
The inequality
\ShowEq{|fx-fx1|<e 2}
where $k$ is the number of terms in the sum
\EqRef{|fx-fx1|<e 1}
follows from
\EqRef{|fx-fx1|<e 1},
\EqRef{E=max|xei|}.

From the inequality
\EqRef{|fx-fx1|<e 2},
it follows that the norm
\ShowEq{|fx-fx1|}
can be made arbitrarily small
by an appropriate choice of $\epsilon$.
Therefore the map $f$ is continuous.

If we assume that at least one of maps
\ShowEq{fi1p}
is not continuous,
then, from the arguments in the proof, it follows that
the map $f$ is not continuous.
\hfill\(\odot\)

The theorem follows from the lemma
\RefLemma{standard components of tensor are continuous maps},
since standard representation is particular case
of representation of tensor.
\end{proof}

\begin{theorem}
\labelLemma{standard components of tensor are differentiable maps}
Let $A$, $B$ be associative finite dimensional Banach $D$\Hyph algebras.
Let
\ShowEq{U in A}
be open set.
The map
\DrawEq{f:U->Bp}{standard differentiable}
is differentiable iff
standard components of tensor $f(x)$ are differentiable.
In such case\,\footnote{
I recall that the map
\ShowEq{standard df/dx}
maps $A$\Hyph number $a$ to $D$\Hyph number
\ShowEq{standard df/dx a}
}
\ShowEq{derivative of standard components of tensor}
\end{theorem}
\begin{proof}
The theorem follows from theorems
\RefTheorem{bilinear map and differential},
\RefTheorem{definition of module over commutative ring}.
\end{proof}

\begin{theorem}
\labelTheorem{components of tensor are differentiable maps}
Let $A$, $B$ be associative finite dimensional Banach $D$\Hyph algebras.
Let
\ShowEq{U in A}
be open set.
The map
\DrawEq{f:U->Bp}{differentiable}
is differentiable iff
the map $f$ has form
\ShowEq{f=sum f}
where maps
\ShowEq{fis:U->B}
are differentiable.
In such case
\ShowEq{derivative of components of tensor}
\end{theorem}
\begin{proof}
The theorem follows from the theorem
\RefTheorem{derivative of tensor product}.
\end{proof}

\begin{theorem}
Let $A$, $B$ be associative finite dimensional Banach $D$\Hyph algebras.
Let
\ShowEq{U in A}
be open map.
The form
\ShowEq{diff form in U}{\omega}npB{}
iff
the form $\omega$ has form
\ShowEq{o=sum o}
where
\ShowEq{diff form in U}{\omega_{s\cdot i}}n0B.
\end{theorem}
\begin{proof}
According to theorems
\RefTheorem[\RefCalculus]{|f(a)|<|f||a| 1n},
\RefTheorem[\RefCalculus]{|on|->0 ona1p->0},
the map
\ShowEq{os Is}
is continuous, when the map
\ShowEq{os}
is continuous.
The statement of the theorem for $n=0$ follows from theorems
\RefTheorem{exterior differential of differential form},
\RefTheorem{tensor is continuous iff}.

Let the theorem be true for $n=k$. Let
\ShowEq{diff form in U}{\omega}{k+1}pB.
Then
\ShowEq{diff form in U}{\omega}kpB{}
and we can write differential form $\omega$ as
\ShowEq{o=sum o}
where
\ShowEq{diff form in U}{\omega_{s\cdot i}}k0B.
We can write differential form
\ShowEq{diff form in U}{d^k\omega}1pB{}
as
\ShowEq{o=sum o dk}
where
\ShowEq{okis=}
where an index $s$ depends from index $t$ and $0\le j\le k$.
According to the theorem
\RefTheorem{components of tensor are differentiable maps},
the maps
\ShowEq{okis}
are differentiable
and, for any index $s$, there exists an index $t$ such that
\ShowEq{okis=k}
Therefore,
\ShowEq{diff form in U}{\omega_{s\cdot i}}{k+1}0B.
We proved the theorem for $n=k+1$.

Therefore, we proved the theorem for any value of $n$.
\end{proof}

\Section{Differential \texorpdfstring{$1$}{1}-form}

\begin{theorem}
Let $A$,$B$ be free Banach $D$\Hyph algebras.
Let
\ShowEq{U subset A}
be open starlike set.
Differential $1$\Hyph form
\ShowEq{diff form in U}{\omega}n1B{}
is integrable iff
\ShowEq{d omega=0}
\end{theorem}
\begin{proof}
If differential $1$\Hyph form $\omega$
is integrable,
then, according to the definition
\RefDefinition{integrable form},
there exist the map $f$ such that
\ShowEq{o=df}
According to the theorem
\RefTheorem{d2omega=0},
\ShowEq{do=d2f}
At the same time, if
\DrawEq{do=0}{}
then, according to the theorem
\RefTheorem{Poincare Theorem p=1},
the differential $1$\Hyph form $\omega$ is integrable.
\end{proof}

\begin{definition}
{\it
Let $A$,$B$ be free Banach $D$\Hyph algebras.
Let
\ShowEq{U subset A}
be open set.
Let
\ShowEq{g:ab->U}
be a path of class $C^1$ in $U$.
We define
\AddIndex{the integral of the differential $1$\Hyph form $\omega$ along the path}
{integral of differential 1 form along path}
$\gamma$ by the equality
\ShowEq{integral of differential 1 form along path}
\ShowEq{integral of differential 1 form along path =}
}
\qed
\end{definition}

\begin{theorem}
\labelTheorem{int o=sum int a=...<t<...=b}
Let
\ShowEq{a=...<t<...=b}
be subdivision of segment $[a,b]$.
Let $\gamma_i$ denote the path corresponding to  segment
\ShowEq{[ti]}
Then
\ShowEq{int o=sum int}
\end{theorem}
\begin{proof}
The equality
\EqRef{int o=sum int}
follows from the equality
\ShowEq{int ac=int ab+int bc}
and induction over number of terms.
\end{proof}

\begin{theorem}
\labelTheorem{int df=fb-fa}
Let $A$,$B$ be free Banach $D$\Hyph algebras.
Let
\ShowEq{U subset A}
be open set.
Let
\ShowEq{g:ab->U}
be a path of class $C^1$ in $U$.
Let
\ShowEq{f:A->B}fAB
be differentiable map.
Then
\ShowEq{int df=fb-fa}
\end{theorem}
\begin{proof}
The theorem follows from the theorem
\RefTheorem{Lebesgue Integral along Path}.\,\footnote{
See also the theorem
\citeBib{Cartan differential form}\Hyph 3.4.1
on the page 43.
}
\end{proof}

\begin{theorem}
\labelTheorem{integral depends only on origin and extremity}
Let $A$,$B$ be free Banach $D$\Hyph algebras.
Let
\ShowEq{U subset A}
be open set.
Let
\ShowEq{f:A->B}fAB
be differentiable map.
The integral
\ShowEq{int df}
depends only on origin $\gamma(a)$ and extremity $\gamma(b)$
of piecewise path $\gamma$ of class $C^1$.
\end{theorem}
\begin{proof}
The theorems follows from theorems
\RefTheorem{int o=sum int a=...<t<...=b},
\RefTheorem{int df=fb-fa}.
\end{proof}

\begin{theorem}
Let $A$, $B$ be free Banach $D$\Hyph algebras.
Let
\ShowEq{U subset A}
be open connected set.
The following properties of differential form
\ShowEq{diff form in U}{\omega}11B{}
are equivalent
\StartLabelItem
\begin{enumerate}
\item
$\omega$ is integrable in $U$
\labelItem{omega is integrable}
\item
\DrawEq{int g o=0}{-}
for any loop $\gamma$, piecewise of class $C^1$, contatined in $U$
\labelItem{int g o=0}
\end{enumerate}
\end{theorem}
\begin{proof}
The statement
\RefItem{int g o=0}
follows from the statement
\RefItem{omega is integrable}
according to theorems
\RefTheorem{int o=sum int a=...<t<...=b},
\RefTheorem{integral depends only on origin and extremity}.

Let the statement
\RefItem{int g o=0}
be true.
Let
\ShowEq{x0 x in U}
Since the set $U$ is connected,
then there exist piecewise path of class $C^1$
\ShowEq{g1:[ab]->U}
and piecewise path of class $C^1$
\ShowEq{g2:[bc]->U}
The piecewise path of class $C^1$
\ShowEq{g:[ac]->U}
defined by the equality
\ShowEq{g:[ac]->U =}
is a loop.
According to the statement
\RefItem{int g o=0},
\DrawEq{int g o=0}{proof}
According to the theorem
\RefTheorem{int o=sum int a=...<t<...=b},
\ShowEq{int g o=1+2}
From equalities
\eqRef{int g o=0}{proof},
\EqRef{int g o=1+2},
it follows that the integral
\ShowEq{int g1 o=int g o-int g2 o}
does not depends on the path.
Let
\ShowEq{fx=int g1 o}

Since the set $U$ is open, then there exist
\ShowEq{eta1>0}
such that linear path
\ShowEq{g:[01]->U x+th}
is contained in $U$ when
\ShowEq{|h|<eta1}
According to the theorem
\RefTheorem{int o=sum int a=...<t<...=b},
the equality
\ShowEq{fx+h=fx+int o}
follows from the equality
\EqRef{fx=int g1 o}.
The equality
\ShowEq{fx+h-fx=int o}
follows from equalities
\EqRef{g:[01]->U x+th},
\EqRef{fx+h=fx+int o}.
The equality
\ShowEq{fx+h-fx-oh=int o}
follows from equality
\EqRef{fx+h-fx=int o}.

Since the map
\ShowEq{o:U->LAB}
is continuous, then for any
\EqParm{eps>0}{=z}
there exist
\ShowEq{0<eta<eta1}
such that
\ShowEq{oh+t-oh<e}
The statement
\ShowEq{(oh+t-oh)h<eh}
follows from the statement
\EqRef{oh+t-oh<e}.
The statement
\ShowEq{|fx+h-fx-oh|<eh}
follows from the statement
\EqRef{(oh+t-oh)h<eh}
and from the equality
\EqRef{fx+h-fx-oh=int o}.
The statement
\ShowEq{|fx+h-fx-oh|<o(h)}
follows from the statement
\EqRef{|fx+h-fx-oh|<eh}.
From the statement
\EqRef{|fx+h-fx-oh|<o(h)}
and the definition
\RefDefinition{differentiable map, algebra},
it follows that differential form $\omega$
is derivative of the map $f$.
Therefore, the statement
\RefItem{omega is integrable}
follows from the statement
\RefItem{int g o=0}.
\end{proof}

\begin{definition}
{\it
Let $A$, $B$ be free Banach $D$\Hyph algebras.
Let
\ShowEq{U subset A}
be open connected set.
Let differential form
\ShowEq{diff form in U}{\omega}n1B{}
be integrable.
For any $A$\Hyph numbers $a$, $b$,
we define
\AddIndex{definite integral}{definite integral}
by the equality
\ShowEq{definite integral}
\ShowEq{definite integral ab}
for any path $\gamma$ from $a$ to $b$.
}
\qed
\end{definition}



\Section{Complex field}

We may consider complex field
either as $C$\Hyph algebra or as $R$\Hyph algebra.

\begin{itemize}
\item
$C$\Hyph vector space
\ShowEq{L(A;B)}CCC{}
is generated by the map $I_0=E$.
If derivative of the map $f$ of complex field
belongs to the set $CE$,
then the map $f$ is called
\AddIndex{holomorphic}{holomorphic map}.
According to theorems
\RefTheorem{L(CCC)=CE},
\RefTheorem{Cauchy Riemann equations linear},
holomorphic map satisfies to the differential equation
(the Cauchy\Hyph Riemann equations)
\ShowEq{Cauchy Riemann equations, complex field, 1}
\item
According to the theorem
\RefTheorem{linear map of complex field},
$C$\Hyph vector space
\ShowEq{L(A;B)}RCC{}
has the basis consisting of
identity map $E$
and conjugation $I$.
If derivative of the map $f$ of complex field
belongs to the set $CI$,
then the map $f$ is called
\AddIndex{conjugate holomorphic}{conjugate holomorphic map}.
According to theorems
\RefTheorem{L(CCC)=CE},
\RefTheorem{matrix of linear map f in CI},
conjugate holomorphic map satisfies to the differential equation
\ShowEq{conjugate holomorphic map}
In this section we consider how this affects the derivative of a map
\ShowEq{f:A->B}fCC
\end{itemize}

\begin{example}
Consider the derivative of the map
\ShowEq{y=x*n}
The equality
\ShowEq{dx*n}
follows from the equality
\EqRef{y=x*n}.
According to the definition
\RefDefinition{differentiable map, algebra},
the equality
\ShowEq{dx*n/dx}
follows from the equality
\EqRef{dx*n}.
Therefore, Taylor series for the map
\EqRef{y=x*n}
has the following form
\ShowEq{y=(Ix)n}
\qed
\end{example}

Now we are ready to consider more advanced example.

\begin{example}
Consider the derivative of the map
\ShowEq{y=xx*2}
The equality
\ShowEq{xx*2 h1}
follows from the equality
\EqRef{y=xx*2}.
According to the definition
\RefDefinition{differentiable map, algebra},
the equality
\ShowEq{dxx*2/dx}
follows from the equality
\EqRef{xx*2 h1}.
The equality
\ShowEq{xx*2 h1h2}
follows from the equality
\EqRef{dxx*2/dx}.
According to the definition
\RefDefinition{derivative of Second Order, algebra},
the equality\,\footnote{
At the first moment when I had seen the symmetry in the equality
\EqRef{d2xx*2/dx2},
it was a surprise for me.
The initial intent of this example was the expectation
of a differential form similar to differential form
\ShowEq{o=3 dx x2}
considered in the remark
\ref{remark: differential equation 3x^2 does not possess a solution}.
But then I realized that the equality
\EqRef{d2xx*2/dx2}
follows from the theorem
\RefTheorem{d2omega=0}.
}
\ShowEq{d2xx*2/dx2}
follows from the equality
\EqRef{xx*2 h1h2}.
The equality
\ShowEq{d3xx*2/dx3}
\ePrints{1610.309618526,CACAA.06}
\ifx\Semafor\ValueOn
follows from the definition
\RefDefinition{derivative of Order n, algebra}.
\else
follows from the definition
\RefDefinition{derivative of Order n, algebra}
and from the theorem
\RefTheorem{derivative, fx=axb, algebra}.
\fi
Therefore, Taylor series for the map
\EqRef{y=xx*2}
has the following form
\ShowEq{y=(Ex)(Ix)2}
\qed
\end{example}

Immediately following question arises.
If the equality
\ShowEq{dx*/dx}
is true,
then what does this have to do with the theory of functions of complex variable
(\citeBib{Shabat: Complex Analysis})?
There is simple answer on this question.
The theory of functions of complex variable considers
$C$\Hyph vector space
\ShowEq{L(A;B)}CCC.
In this case, there is no linear map
approximating the map $y=x^*$.
We express this statement by the equality
\ShowEq{dx*/dx=0}
which means independence of $x$ and $x^*$.

Since $C$\Hyph vector space $C$
has dimension $1$,
then there is no
skew symmetric polylinear map
\EqParm{f in LA}{f=f,p=2,A=CC,=.}
However $R$\Hyph vector space $C$
has dimension $2$.
So there exists
skew symmetric polylinear map
\EqParm{f in LA}{f=f,p=2,A=C,=.}

For instance, the map
\ShowEq{fa12=a11 a22-a12 a21}
is skew symmetric polylinear map.
Since
\ShowEq{a11= a12=}
where
\ShowEq{Ia=a*}
is map of conjugation,
then the equality
\ShowEq{fa12=re im-im re}
follows from the equality
\EqRef{fa12=a11 a22-a12 a21}.
We can write polylinear map
\EqRef{fa12=re im-im re}
as
\ShowEq{f a1 a2=(I,E)}

\begin{theorem}
\labelTheorem{df/dx= complex field}
If we consider the map of complex field
\ShowEq{f:A->B}fCC
as function of
\ShowEq{x=x0+ix1}
then
\ShowEq{df/dx= complex field}
\end{theorem}
\begin{proof}
According to the theorem
\RefTheorem{linear map of complex field},
derivative of the map $f$ has form
\ShowEq{df/dx=a0E+a1I}
where
\ShowEq{a01= df/dx}
The equality
\EqRef{df/dx= complex field}
follows from equalities
\EqRef{df/dx=a0E+a1I},
\EqRef{a01= df/dx}.
\end{proof}

According to the theorem
\RefTheorem{linear map of complex field},
the differential form
\EqParm{form x in LA}{form=omega p2,A=C,=z}
has the following form
\DrawEq{omega=aE+bI}{}
where $a$, $b$ are complex valued maps.
For any $C$\Hyph number $c_1$
\ShowEq{omega c1}

\begin{theorem}
\labelTheorem{omega=aE+bI is integrable}
The differential form
\DrawEq{omega=aE+bI}{theorem}
is integrable iff
\DrawEq{d omega=0 condition}{theorem}
\end{theorem}
\begin{proof}
Derivative of the differential form
\eqRef{omega=aE+bI}{theorem}
is
\ShowEq{do/dx=da/dx+db/dx}
According to the theorem
\RefTheorem{df/dx= complex field},
\ShowEq{da db /dx}
The equality
\ShowEq{d omega/dx c1 c2}
follows from equalities
\EqRef{do/dx=da/dx+db/dx},
\EqRef{d omega/dx c1 c2}.
The exterior differential of the differential form $\omega$ has the following form
\ShowEq{d omega a12}
According to the theorem
\RefTheorem{Poincare Theorem p=1},
the equality
\eqRef{d omega=0 condition}{theorem}
follows from the equality
\EqRef{d omega a12}.
\end{proof}

Let
\ShowEq{f=int o}
We have two different cases.

\begin{itemize}
\item
If the equality
\ShowEq{df dx01 =0}
is true, then we can represent the map $f$
as sum of the holomorphic map
and the conjugate holomorphic map.
\item
If the equality
\EqRef{df dx01 =0}
is not true, then we cannot represent the map $f$
as sum of the holomorphic map
and the conjugate holomorphic map.
In the section
\RefSection{differential form in complex field},
we consider such form.
\end{itemize}

%% file: Diff.Form.2.Eq.tex

\input{Diff.Form.2.Parm}

\DefEq
{
$b\in B^{\otimes(n+1)}$
}
{b in Bon}

\DefEquation
{
y=(x^*)^n
}
{y=x*n}

\DefEquation
{
\begin{split}
\frac{\partial f^{\giA}}{\partial x^{\gi 0}}
&=\hphantom{-\,}\frac{\partial f^{\gi 0}}{\partial x^{\giA}}
\\[5pt]
\frac{\partial f^{\gi 0}}{\partial x^{\gi 0}}
&=-\frac{\partial f^{\giA}}{\partial x^{\giA}}
\end{split}
}
{conjugate holomorphic map}

\DefEquation
{
\begin{split}
((x+h_1)^*)^n-(x^*)^n
&=(x^*+h_1^*)^n-(x^*)^n
\\&=(x^*)^n-n(x^*)^{n-1}h_1^*+o(h_1)-(x^*)^n
\\&=n(x^*)^{n-1}h_1^*+o(h_1)
\end{split}
}
{dx*n}

\DefEq
{
$\omega(x)=3\,dx\,x^2$
}
{o=3 dx x2}

\DefEquation
{
\frac{d(x^*)^n}{dx}=n(x^*)^{n-1}I
}
{dx*n/dx}

\DefEq
{
\[
y=\frac 1{n!}n!(I\circ x)^n=(I\circ x)^n
\]
}
{y=(Ix)n}

\DefEq
{
\begin{align*}
y&=\frac 16
2((E\circ x)(I\circ x)^2+(I\circ x)(E\circ x)(I\circ x)+(I\circ x)^2(E\circ x))
\\&=(E\circ x)(I\circ x)^2
\end{align*}
}
{y=(Ex)(Ix)2}

\DefEq
{
\[
\frac{dx^*}{dx}=0
\]
}
{dx*/dx=0}

\DefEq
{
\[
\frac{dx^*}{dx}\circ h=h^*
\]
}
{dx*/dx}

\DefEquation
{
y=x(x^*)^2
}
{y=xx*2}

\DefEquation
{
\begin{split}
&\,(x+h_1)((x+h_1)^*)^2-x(x^*)^2
\\=&\,(x+h_1)(x^*+h^*_1)^2-x(x^*)^2
\\=&\,(x+h_1)((x^*)^2+2x^*h^*_1+o(h_1))-x(x^*)^2
\\=&\,(x+h_1)(x^*)^2+2(x+h_1)x^*h^*_1+o(h_1)-x(x^*)^2
\\=&\,x(x^*)^2+h_1(x^*)^2+2xx^*h^*_1+o(h_1)-x(x^*)^2
\\=&\,(x^*)^2h_1+2xx^*h^*_1+o(h_1)
\end{split}
}
{xx*2 h1}

\DefEquation
{
\begin{split}
&\,((x+h_2)^*)^2h_1+2(x+h_2)(x+h_2)^*h^*_1-(x^*)^2h_1-2xx^*h^*_1
\\=&\,(x^*+h_2^*)^2h_1+2(x+h_2)(x^*+h_2^*)h^*_1-(x^*)^2h_1-2xx^*h^*_1
\\=&\,((x^*)^2+2x^*h_2^*+o(h_2))h_1+2((x+h_2)x^*+(x+h_2)h_2^*)h^*_1
\\-&\,(x^*)^2h_1-2xx^*h^*_1
\\=&\,(x^*)^2h_1+2x^*h_1h_2^*+o(h_2)+2(xx^*+x^*h_2+xh_2^*+o(h_2))h^*_1
\\-&\,(x^*)^2h_1-2xx^*h^*_1
\\=&\,2x^*h_1h_2^*+o(h_2)+2xx^*h^*_1+2x^*h^*_1h_2
\\+&\,2xh^*_1h_2^*-2xx^*h^*_1
\\=&\,2x^*h_1h_2^*+2x^*h^*_1h_2+2xh^*_1h_2^*+o(h_2)
\end{split}
}
{xx*2 h1h2}

\DefEquation
{
\frac {d^2x(x^*)^2}{dx^2}=2x^*(E,I)+2x^*(I,E)+2x(I,I)
}
{d2xx*2/dx2}

\DefEquation
{
\frac {dx(x^*)^2}{dx}=(x^*)^2E+2xx^*I
}
{dxx*2/dx}

\DefEquation
{
\frac {d^3x(x^*)^2}{dx^3}=2(E,I,I)+2(I,E,I)+2(I,I,E)
}
{d3xx*2/dx3}

\DefEquation
{
a_0=-\frac 12(f_{\gi 0}+f_{\gi 1}i+f_{\gi 2}j+f_{\gi 3}k)
}
{a0=}

\DefEquation
{
a_1=\frac 12(f_{\gi 0}+f_{\gi 1}i)
}
{a1=}

\DefEquation
{
a_2=\frac 12(f_{\gi 0}+f_{\gi 2}j)
}
{a2=}

\DefEquation
{
a_3=\frac 12(f_{\gi 0}+f_{\gi 3}k)
}
{a3=}

\DefEq
{
\omega(x)=a(x)\circ E+b(x)\circ I
}
{omega=aE+bI}

\DefEq
{
\[
\omega(x)\circ c_1=a(x)\circ E\circ c_1+b(x)\circ I\circ c_1
=a(x)c_1+b(x)c^*_1
\]
}
{omega c1}

\DefEquation
{
\frac{df(x)}{dx}\circ c_1=a_0(x)\circ E\circ c_1+a_1(x)\circ I\circ c_1
}
{df/dx=a0E+a1I}

\DefEquation
{
\begin{split}
\begin{pmatrix}
a_0\\a_1
\end{pmatrix}
&=
\begin{pmatrix}
a_0^{\gi 0}&a_0^{\gi 1}
\\
a_1^{\gi 0}&a_1^{\gi 1}
\end{pmatrix}
\begin{pmatrix}
1\\i
\end{pmatrix}
=\frac 12
\begin{pmatrix}
1&1
\\
1&-1
\end{pmatrix}
\begin{pmatrix}
\displaystyle\frac{\partial f^{\gi 0}}{\partial x^{\gi 0}}&
\displaystyle\frac{\partial f^{\gi 1}}{\partial x^{\gi 0}}
\\[10pt]
\displaystyle\frac{\partial f^{\gi 1}}{\partial x^{\gi 1}}&
-\displaystyle\frac{\partial f^{\gi 0}}{\partial x^{\gi 1}}
\end{pmatrix}
\begin{pmatrix}
1\\i
\end{pmatrix}
\\
&=\frac 12
\begin{pmatrix}
1&1
\\
1&-1
\end{pmatrix}
\begin{pmatrix}
\displaystyle
\frac{\partial f^{\gi 0}}{\partial x^{\gi 0}}+
\frac{\partial f^{\gi 1}}{\partial x^{\gi 0}}i
\\[10pt]
\displaystyle
\frac{\partial f^{\gi 1}}{\partial x^{\gi 1}}
-\frac{\partial f^{\gi 0}}{\partial x^{\gi 1}}i
\end{pmatrix}
\\
&=\frac 12
\begin{pmatrix}
1&1
\\
1&-1
\end{pmatrix}
\begin{pmatrix}
\displaystyle
\frac{\partial f^{\gi 0}}{\partial x^{\gi 0}}+
\frac{\partial f^{\gi 1}}{\partial x^{\gi 0}}i
\\[10pt]
\displaystyle
-i\left(
\frac{\partial f^{\gi 0}}{\partial x^{\gi 1}}
+\frac{\partial f^{\gi 1}}{\partial x^{\gi 1}}i
\right)
\end{pmatrix}
=\frac 12
\begin{pmatrix}
1&1
\\
1&-1
\end{pmatrix}
\begin{pmatrix}
\displaystyle
\frac{\partial f}{\partial x^{\gi 0}}
\\[10pt]
\displaystyle
-i\frac{\partial f}{\partial x^{\gi 1}}
\end{pmatrix}
\\
&=\frac 12
\begin{pmatrix}
\displaystyle
\frac{\partial f}{\partial x^{\gi 0}}
-i\frac{\partial f}{\partial x^{\gi 1}}
\\[10pt]
\displaystyle
\frac{\partial f}{\partial x^{\gi 0}}
+i\frac{\partial f}{\partial x^{\gi 1}}
\end{pmatrix}
\end{split}
}
{a01= df/dx}

\DefEquation
{
\frac{d\omega(x)\circ c_1}{dx}=
\frac{da(x)}{dx}c_1+
\frac{db(x)}{dx}c_1^*
}
{do/dx=da/dx+db/dx}

\DefEquation
{
\begin{split}
\frac{d\omega(x)\circ c_1}{dx}\circ c_2
&=
\left(\frac{da(x)}{dx}\circ c_2\right)c_1+
\left(\frac{db(x)}{dx}\circ c_2\right)c_1^*
\\&=
\left(
\frac 12\left(
\frac{\partial a}{\partial x^{\gi 0}}
-i\frac{\partial a}{\partial x^{\gi 1}}
\right)c_2
+\frac 12\left(
\frac{\partial a}{\partial x^{\gi 0}}
+i\frac{\partial a}{\partial x^{\gi 1}}
\right)c_2^*
\right)c_1
\\&+
\left(
\frac 12\left(
\frac{\partial b}{\partial x^{\gi 0}}
-i\frac{\partial b}{\partial x^{\gi 1}}
\right)c_2
+\frac 12\left(
\frac{\partial b}{\partial x^{\gi 0}}
+i\frac{\partial b}{\partial x^{\gi 1}}
\right)c_2^*
\right)c_1^*
\end{split}
}
{d omega/dx c1 c2}

\DefEquation
{
\begin{split}
\frac{da(x)}{dx}
&=
\frac 12\left(
\frac{\partial a}{\partial x^{\gi 0}}
-i\frac{\partial a}{\partial x^{\gi 1}}
\right)E
+\frac 12\left(
\frac{\partial a}{\partial x^{\gi 0}}
+i\frac{\partial a}{\partial x^{\gi 1}}
\right)I
\\
\frac{db(x)}{dx}
&=
\frac 12\left(
\frac{\partial b}{\partial x^{\gi 0}}
-i\frac{\partial b}{\partial x^{\gi 1}}
\right)E
+\frac 12\left(
\frac{\partial b}{\partial x^{\gi 0}}
+i\frac{\partial b}{\partial x^{\gi 1}}
\right)I
\end{split}
}
{da db /dx}

\DefEq
{
$x^{\gi 0}$, $x^{\gi 1}$,
$x=x^{\gi 0}+x^{\gi 1}i$,
}
{x=x0+ix1}

\DefEquation
{
\frac{df(x)}{dx}=
\frac 12\left(
\frac{\partial f}{\partial x^{\gi 0}}
-i\frac{\partial f}{\partial x^{\gi 1}}
\right)E
+\frac 12\left(
\frac{\partial f}{\partial x^{\gi 0}}
+i\frac{\partial f}{\partial x^{\gi 1}}
\right)I
}
{df/dx= complex field}

\DefEquation
{
\begin{split}
d\omega(x)\circ (c_1,c_2)
&=\frac 12\left(
\frac{\partial a}{\partial x^{\gi 0}}
+i\frac{\partial a}{\partial x^{\gi 1}}
\right)(c_1c_2^*-c_1^*c_2)
\\
&+\frac 12\left(
\frac{\partial b}{\partial x^{\gi 0}}
-i\frac{\partial b}{\partial x^{\gi 1}}
\right)(c_1^*c_2-c_1c_2^*)
\\
&=\frac 12\left(
\frac{\partial a}{\partial x^{\gi 0}}
+i\frac{\partial a}{\partial x^{\gi 1}}
-\frac{\partial b}{\partial x^{\gi 0}}
+i\frac{\partial b}{\partial x^{\gi 1}}
\right)(c_1c_2^*-c_1^*c_2)
\end{split}
}
{d omega a12}

\DefEq
{
\frac{\partial a}{\partial x^{\gi 0}}
+i\frac{\partial a}{\partial x^{\gi 1}}
-\frac{\partial b}{\partial x^{\gi 0}}
+i\frac{\partial b}{\partial x^{\gi 1}}
=0
}
{d omega=0 condition}

\DefEquation
{
\frac{\partial a}{\partial x^{\gi 0}}
+i\frac{\partial a}{\partial x^{\gi 1}}
=0
\ \ \ \ \,
-\frac{\partial b}{\partial x^{\gi 0}}
+i\frac{\partial b}{\partial x^{\gi 1}}
=0
}
{df dx01 =0}

\DefEq
{
\[f=\int\omega\]
}
{f=int o}

\DefEquation
{
d\omega(x)=0
}
{d omega=0}

\DefEq
{
\[\omega(x)=df(x)\]
}
{o=df}

\DefEq
{
\[d\omega(x)=d^2f(x)=0\]
}
{do=d2f}

\DefEq
{
\[
a=t_0<...<t_i<...<t_n=b
\]
}
{a=...<t<...=b}

\DefEquation
{
\int_{\gamma}\omega=\sum_{i=1}^n\int_{\gamma_i}\omega
}
{int o=sum int}

\DefEq
{
$x_0$, $x\in U$.
}
{x0 x in U}

\DefEq
{
\[
\gamma_1:[a,b]\subset R\rightarrow U\ \ \ \gamma_1(a)=x_0\ \ \ \gamma_1(b)=x
\]
}
{g1:[ab]->U}

\DefEq
{
\[
\gamma_2:[b,c]\subset R\rightarrow U\ \ \ \gamma_2(b)=x\ \ \ \gamma_2(c)=x_0
\]
}
{g2:[bc]->U}

\DefEq
{
\[
\gamma:[a,c]\subset R\rightarrow U
\]
}
{g:[ac]->U}

\DefEq
{
\[
\gamma(x)=\left\{
\begin{matrix}
\gamma_1(x)&x\in[a,b]\\
\gamma_2(x)&x\in[b,c]
\end{matrix}
\right.
\]
}
{g:[ac]->U =}

\DefEq
{
\displaystyle\int_{\gamma}\omega=0
}
{int g o=0}

\DefEq
{
$[t_{i-1},t_i]$, $i=1$, ..., $n$.
}
{[ti]}

\DefEquation
{
\int_{\gamma}\omega=\int_{\gamma_1}\omega+\int_{\gamma_2}\omega
}
{int g o=1+2}

\DefEq
{
\[
\int_{\gamma_1}\omega=-\int_{\gamma_2}\omega
\]
}
{int g1 o=int g o-int g2 o}

\DefEquation
{
f(x)=\int_{\gamma_1}\omega
}
{fx=int g1 o}

\DefEquation
{
\gamma:t\in[0,1]\subset R\rightarrow x+th\in A
}
{g:[01]->U x+th}

\DefEquation
{
f(x+h)=f(x)+\int_x^{x+h}\omega
}
{fx+h=fx+int o}

\DefEquation
{
f(x+h)-f(x)=\int_0^1dt\,\omega(x+th)\circ h
}
{fx+h-fx=int o}

\DefEquation
{
f(x+h)-f(x)-\omega(x)\circ h=\int_0^1dt(\omega(x+th)\circ h-\omega(x)\circ h)
}
{fx+h-fx-oh=int o}

\DefEquation
{
\|f(x+h)-f(x)-\omega(x)\circ h\|\le\epsilon\|h\|
}
{|fx+h-fx-oh|<eh}

\DefEquation
{
\|f(x+h)-f(x)-\omega(x)\circ h\|=o(h)
}
{|fx+h-fx-oh|<o(h)}

\DefEq
{
\[
\omega:U\rightarrow\mathcal{L}(D;A\rightarrow B)
\]
}
{o:U->LAB}

\DefEquation
{
\|\omega(x+th)-\omega(x)\|<\epsilon
\ \ \ \ \|h\|<\eta\ \ \ \ t\in [0,1]
}
{oh+t-oh<e}

\DefEquation
{
\|\omega(x+th)\circ h-\omega(x)\circ h\|<\epsilon\|h\|
\ \ \ \ \|h\|<\eta\ \ \ \ t\in [0,1]
}
{(oh+t-oh)h<eh}

\DefEq
{
$\eta$, $0<\eta<\eta_1$,
}
{0<eta<eta1}

\DefEq
{
\begin{align*}
\int_{\gamma}\omega&=\int_{t_0}^{t_2}dt\omega(\gamma(t))\frac{d\gamma(t)}{dt}
\\&=\int_{t_0}^{t_1}dt\omega(\gamma(t))\frac{d\gamma(t)}{dt}
+\int_{t_1}^{t_2}dt\omega(\gamma(t))\frac{d\gamma(t)}{dt}
\\&=\int_{\gamma_1}\omega+\int_{\gamma_2}\omega
\end{align*}
}
{int ac=int ab+int bc}

\DefEq
{
\symb{\int_a^b\omega}{definite integral}{}
}
{definite integral}

\DefEq
{
\[
\ShowSymbol{definite integral}{}=\int_{\gamma}\omega
\]
}
{definite integral ab}

\DefEq
{
\[
\omega(x)=a_0(x)\otimes a_1(x)\otimes a_2(x)
\]
}
{omega=a012}

\DefEquation
{
\begin{split}
\frac{d\omega(x)}{dx}
&=\frac{d_{s_0\cdot 0}a_0(x)}{dx}\otimes\frac{d_{s_0\cdot 1}a_0(x)}{dx}\otimes a_1(x)\otimes a_2(x)
\\&+a_0(x)\otimes \frac{d_{s_1\cdot 0}a_1(x)}{dx}\otimes\frac{d_{s_1\cdot 1}a_1(x)}{dx}\otimes a_2(x)
\\&+a_0(x)\otimes a_1(x)\otimes \frac{d_{s_2\cdot 0}a_2(x)}{dx}\otimes\frac{d_{s_2\cdot 1}a_2(x)}{dx}
\end{split}
}
{d omega/dx=a012}

\DefEquation
{
\begin{split}
\frac{d\omega(x)}{dx}\circ (b_1,b_2,b_3)
&=\frac{d_{s_0\cdot 0}a_0(x)}{dx}b_1\frac{d_{s_0\cdot 1}a_0(x)}{dx}b_2 a_1(x)b_3 a_2(x)
\\&+a_0(x)b_2 \frac{d_{s_1\cdot 0}a_1(x)}{dx}b_1\frac{d_{s_1\cdot 1}a_1(x)}{dx}b_3 a_2(x)
\\&+a_0(x)b_2 a_1(x)b_3 \frac{d_{s_2\cdot 0}a_2(x)}{dx}b_1\frac{d_{s_2\cdot 1}a_2(x)}{dx}
\end{split}
}
{d omega/dx=a012o}

\DefEq
{
$\displaystyle\int_{\gamma}\omega$
}
{int h omega}

\DefEq
{
$\eta_1>0$
}
{eta1>0}

\DefEq
{
$\|h\|<\eta_1$.
}
{|h|<eta1}

\DefEq
{
$\displaystyle\int_{\gamma}df(x)$
}
{int df}

\DefEquation
{
\int_{\gamma}df=f(g(b))-f(g(a))
}
{int df=fb-fa}

\DefEq
{
\symb{\int_{\gamma}\omega}{integral of differential 1 form along path}{}
}
{integral of differential 1 form along path}

\DefEq
{
\[
\ShowSymbol{integral of differential 1 form along path}{}
=\int_a^bdt\omega(\gamma(t))\frac{d\gamma(t)}{dt}
\]
}
{integral of differential 1 form along path =}

\DefEq
{
\[
\gamma:[a,b]\rightarrow U
\]
}
{g:ab->U}

\DefEquation
{
\begin{split}
f\circ(a_1,...,a_n)
&=
(b_s\circ F_{[s]}\circ\sigma_s)\circ(a_1,...,a_n)
\\&=
b_{0\cdot s}(F_{[s](1)}\circ\sigma_s( a_1))b_{1\cdot s}...
(F_{[s](n)}\circ \sigma_s(a_n))b_{n\cdot s}
\end{split}
}
{f circ=b circ n}

\DefEq
{
\[
f\circ(a_1,a_2)=(1\otimes 1\otimes 1)\circ(I,E)\circ(a_1,a_2)
-(1\otimes 1\otimes 1)\circ(I,E)\circ(a_2,a_1)
\]
}
{f a1 a2=(I,E)}

\DefEquation
{
f=b_s\circ F_{[s]}\circ\sigma_s
}
{f=(b sigma)oI s}

\DefEquation
{
[f]=\frac 1{n!}\sum_{\sigma\in S(n)}|\sigma||\sigma_s|
b_s\circ F_{[s]}\circ\sigma
}
{[f]=(b sigma)oI s}

\DefEquation
{
f=[b_s\circ F_{[s]}]
}
{b o I}

\DefEquation
{
\begin{split}
&\,(f\circ(a_1,...,a_p))(g\circ(a_{p+1},...,a_{p+q}))
\\=&\,
(b_s\circ(F_{[s](1)},...,F_{[s](p)})\circ\sigma_s\circ(a_1,...,a_p))
\\ *&\,
(c_t\circ(F_{[t](1)},...,F_{[t](q)})\circ\tau_t\circ(a_{p+1},...,a_{p+q}))
\\=&\,
(b_{s\cdot 0}(F_{[s](1)}\circ(\sigma_s(a_1)))b_{s\cdot 1}
...(F_{[s](p)}\circ(\sigma_s(a_p)))b_{s\cdot p})
\\ *&\,
(c_{t\cdot 0}(F_{[t](1)}\circ(\tau_t(a_{p+1})))c_{t\cdot 1}
...(F_{[t](q)}\circ(\tau_t(a_{p+q})))c_{t\cdot q})
\end{split}
}
{f1of2o}

\DefEquation
{
\begin{split}
&\,(f\circ(a_1,...,a_p))(g\circ(a_{p+1},...,a_{p+q}))
\\=&\,
b_{s\cdot 0}(F_{[s](1)}\circ(\sigma_{st}(a_1))b_{s\cdot 1}
...(F_{[s](p)}\circ(\sigma_{st}(a_p)))b_{s\cdot p}
\\ *&\,
c_{t\cdot 0}(F_{[t](1)}\circ(\sigma_{st}(a_{p+1}))c_{t\cdot 1}
...(F_{[t](q)}\circ(\sigma_{st}(a_{p+q})))c_{t\cdot q}
\end{split}
}
{f1of2o 1}

\DefEquation
{
\begin{split}
&\,(f_1\circ(a_1,...,a_p))(f_2\circ(a_{p+1},...,a_{p+q}))
\\=&\,
(b_s\underline{\otimes}c_t)
\circ(F_{[s](1)},...,F_{[s](p)},F_{[t](1)},...,F_{[t](q)})
\circ\sigma_{st}\circ(a_1,...,a_{p+q})
\end{split}
}
{f1of2o 2}

\DefEquation
{
\sigma_{st}=
\begin{pmatrix}
a_1&...&a_p&a_{p+1}&...&a_{p+q}\\
\sigma_s(a_1)&...&\sigma_s(a_p)&\tau_t(a_{p+1})&...&\tau_t(a_{p+q})
\end{pmatrix}
}
{sigma st}

\DefEq
{
\Pf=[\Pa_{\Pp}\circ(F_{[\Pp](1)},...,F_{[\Pp](\Pb)})]
\in\mathcal{LA}(D;A^{\Pb}\rightarrow B)
}
{f in LAAB =}

\DefEq
{
$\Pf\in\mathcal{LA}(\pD;\pA^{\Pp}\rightarrow \pB)$\Pt
}
{f in LA}

\DefEquation
{
b_s=b_s^{\gi{i_0...i_p}}e_{\gi{i_0}}\otimes...\otimes e_{\gi{i_p}}
}
{b=b...e}

\DefEquation
{
f=b_s^{\gi{i_0...i_p}}[(e_{\gi{i_0}}\otimes...\otimes e_{\gi{i_p}})
\circ(F_{[s](1)},...,F_{[s](p)})]
}
{f=b...e}

\DefEquation
{
f_1\wedge f_2=
[(b_s\underline{\otimes}c_t)\circ(F_{[s](1)},...,F_{[s](p)},F_{[t](1)},...,F_{[t](q)})]
}
{b1 A b2}

\DefEq
{
$U\subseteq A$
}
{U in A}

\AddEq{fis:U->B}
{
\[
f_{s\cdot i}:U\rightarrow B
\]
}

\DefEq
{
$f=f_0\otimes f_1$
}
{f=f0of1}

\DefEq
{
$\epsilon>0$\Pt
}
{e>0}

\DefEq
{
$\delta>0$
}
{d>0}

\DefEq
{
\[
\|x-x_1\|_A<\delta
\]
}
{|x-x1|<d}

\DefEquation
{
\|f^{\gi{i_0...i_p}}(x)-f^{\gi{i_0...i_p}}(x_1)\|<\epsilon
}
{|fix-fix1|<e}

\DefEquation
{
\begin{split}
&\,\|f(x)-f(x_1)\|
\\=&\,\|f^{\gii_0....\gii_p}(x)e_{B\cdot\gii_0}\otimes...\otimes e_{B\cdot\gii_p}
-f^{\gii_0....\gii_p}(x_1)e_{B\cdot\gii_0}\otimes...\otimes e_{B\cdot\gii_p}\|
\\=&\,\|(f^{\gii_0....\gii_p}(x)
-f^{\gii_0....\gii_p}(x_1))e_{B\cdot\gii_0}\otimes...\otimes e_{B\cdot\gii_p}\|
\\=&\,\|(f^{\gii_0....\gii_p}(x)
-f^{\gii_0....\gii_p}(x_1))\|\|e_{B\cdot\gii_0}\otimes...\otimes e_{B\cdot\gii_p}\|
\\ \le &\,
\epsilon\sum_{(\gii_0...\gii_p)}
\|e_{B\cdot\gii_0}\otimes...\otimes e_{B\cdot\gii_p}\|
\end{split}
}
{|fx-fx1|<e 1}

\DefEquation
{
\|f(x)-f(x_1)\|
\le
\epsilon k E
}
{|fx-fx1|<e 2}

\DefEq
{
$\|f(x)-f(x_1)\|$
}
{|fx-fx1|}

\DefEq
{
$e_{B\cdot\gii_0}\otimes...\otimes e_{B\cdot\gii_p}$
}
{xei}

\DefEquation
{
\|e_{B\cdot\gii_0}\otimes...\otimes e_{B\cdot\gii_p}\|=
\max\frac{\|(e_{B\cdot\gii_0}\otimes...\otimes e_{B\cdot\gii_p})\circ(b_1,...,b_p)\|}
{\|b_1\|...\|b_p\|}
}
{|xei|}

\DefEquation
{
E=\max\|e_{B\cdot\gii_0}\otimes...\otimes e_{B\cdot\gii_p}\|
}
{E=max|xei|}

\DefEquation
{
f(x)=f^{\gii_0....\gii_p}(x)e_{\gii_0}\otimes...\otimes e_{\gii_p}
}
{f(x)=...xe}

\DefEq
{
\[
f^{\gii_0....\gii_p}:U\rightarrow D
\]
}
{fi:U->D}

\DefEq
{
$f^{\gii_0...\gii_p}$
}
{fi1p}

\DefEq
{
\[
f(x)=((f_0(x)+sign(x^{\gi 1}))\otimes f_2(x))+((f_0(x)-sign(x^{\gi 1}))\otimes f_2(x))
\]
}
{f=f0of1 1}

\DefEq
{
\begin{align*}
f_0:&U\rightarrow B\\
f_1:&U\rightarrow B
\end{align*}
}
{f01:U->B}

\AddEq{f=sum f}
{
\[f=f_{s\cdot 0}\otimes...\otimes f_{s\cdot p}\]
}

\AddEq{o=sum o}
{
\[\omega=[(\omega_{s\cdot 0}\otimes...\otimes \omega_{s\cdot p})\circ F_{[s]}]\]
}

\AddEq{o=sum o dk}
{
\[\omega=[(\omega^k_{t\cdot 0}\otimes...\otimes \omega^k_{t\cdot p})\circ F_{[t]}]\]
}

\AddEq{okis=}
{
\[
\omega^k_{t\cdot i}(x)=\frac{d^j \omega_{s\cdot i}(x)}{dx^j}
\]
}

\AddEq{okis=k}
{
\[
\omega^k_{t\cdot i}(x)=\frac{d^k \omega_{s\cdot i}(x)}{d x^k}
\]
}

\AddEq{okis}
{
$\omega^k_{t\cdot i}$
}

\AddEq{os Is}
{
\[(\omega_{s\cdot 0}\otimes...\otimes \omega_{s\cdot p})\circ F_{[s]}\]
}

\DefEq
{
\[(\omega_{s\cdot 0}\otimes...\otimes \omega_{s\cdot p})\]
}
{os}

\DefEquation
{
\frac{df(x)}{dx}\circ a=
\left(\frac{df^{\gii_0....\gii_p}(x)}{dx}\circ a\right)
e_{\gii_0}\otimes...\otimes e_{\gii_p}
}
{derivative of standard components of tensor}

\DefEq
{
$\displaystyle\frac{df^{\gii_0....\gii_p}(x)}{dx}$
}
{standard df/dx}

\DefEq
{
$\displaystyle\frac{df^{\gii_0....\gii_p}(x)}{dx}\circ a$.
}
{standard df/dx a}

\DefEquation
{
\begin{split}
\frac{df(x)}{dx}\circ a&=
\left(\frac{df_{0\cdot s}(x)}{dx}\circ a\right)
\otimes f_{1\cdot s}\otimes...\otimes f_{p\cdot s}
\\&+
f_{0\cdot s}\otimes
\left(\frac{df_{1\cdot s}(x)}{dx}\circ a\right)
\otimes ...\otimes f_{p\cdot s}+...
\\&+f_{0\cdot s}\otimes\otimes ...\otimes
\left(\frac{df_{p\cdot s}(x)}{dx}\circ a\right)
\end{split}
}
{derivative of components of tensor}

\DefEq
{
f:U\rightarrow B^{\otimes p}
}
{f:U->Bp}

\DefEquation
{
f(a_1,a_2)=a_1^{\gi 0}a_2^{\gi 1}-a_1^{\gi 1}a_2^{\gi 0}
}
{fa12=a11 a22-a12 a21}

\DefEq
{
\[
\begin{matrix}
a_1^{\gi 0}=\displaystyle\frac 12(a_1+I\circ a_1)&a_1^{\gi 1}=\displaystyle\frac 12(a_1-I\circ a_1)
\\[10pt]
a_2^{\gi 0}=\displaystyle\frac 12(a_2+I\circ a_2)&a_2^{\gi 1}=\displaystyle\frac 12(a_2-I\circ a_2)
\end{matrix}
\]
}
{a11= a12=}

\DefEq
{
\[I\circ a=a^*\]
}
{Ia=a*}

\DefEquation
{
\begin{split}
f(a_1,a_2)&=\frac 14(a_1+I\circ a_1)(a_2-I\circ a_2)-\frac 14(a_1-I\circ a_1)(a_2+I\circ a_2)
\\
&=\frac 14(a_1(a_2-I\circ a_2)+(I\circ a_1)(a_2-I\circ a_2)
\\
&-a_1(a_2+I\circ a_2)+(I\circ a_1)(a_2+I\circ a_2))
\\
&=\frac 14(a_1a_2-a_1(I\circ a_2)+(I\circ a_1)a_2-(I\circ a_1)(I\circ a_2)
\\
&-a_1a_2-a_1(I\circ a_2)+(I\circ a_1)a_2+(I\circ a_1)(I\circ a_2))
\\
&=\frac 12((I\circ a_1)a_2-a_1(I\circ a_2))
\end{split}
}
{fa12=re im-im re}

\DefEq
{
\[
f:U\rightarrow B^{\otimes 2}
\]
}
{f:U->B2}

\DefEq
{
$|\sigma_2|=|\sigma||\sigma_s|$.
}
{|s2|=|s||ss|}

\DefEquation
{
[f]=\frac 1{n!}\sum_{\sigma_2\in S(n)}|\sigma_2|
b_s\circ F_{[s]}\circ\sigma_2\circ\sigma_s
}
{[f]=(b sigma)oI s 1}

\DefEq
{
$\underline{\otimes}$\Pt
}
{otimes -=}

\DefEq
{
$\sigma=\sigma_2\circ\sigma_s$\Pt
}
{s=s2 o ss}

\DefEq
{
$\sigma_2\in S(n)$\Pt
}
{s2 in S}

%% file: Diff.Form.2.Parm.tex

\DefEq%
{%
\def\Pa{b}%
\def\Pf{f}%
\def\Pp{s}%
\def\Pb{p}%
}%
{nn=}%

\DefEq%
{%
\def\Pa{b}%
\def\Pf{f}%
\def\Pp{s}%
\def\Pb{p}%
}%
{nn=1}%

\DefEq%
{%
\def\Pa{c}%
\def\Pf{g}%
\def\Pp{t}%
\def\Pb{q}%
}%
{nn=2}%

%% file: Appendix.16.Calculate.English.tex

\input{Appendix.16.Calculate.Eq}

\Chapter{Supporting Calculations}

In this appendix, I put calculations
which I need to design this book.
But I did not want that these calculations distracted from the basic logic of the text.

\Section{Calculations to Estimate Integral (\ref{eq: int h2=int 12 xx+})}

We reduce integrand of the integral
\EqRef{int h2=int 12 xx+}
\ShowEq{int h2=int 12 xx+ 1}
The equality
\ShowEq{int h2=int 12 xx+ 2}
follows from the equality
\EqRef{int h2=int 12 xx+ 1}.
The equality
\ShowEq{int h2=int 12 xx+ 21}
follows from the equality
\EqRef{int h2=int 12 xx+ 2}.
The equality
\ShowEq{int h2=int 12 xx+ 3}
follows from the equality
\EqRef{int h2=int 12 xx+ 21}.
The equality
follows from the equality
\ShowEq{int h2=int 12 xx+ 4}
follows from the equality
\EqRef{int h2=int 12 xx+ 3}.
The equality
\ShowEq{int h2=int 12 xx+ 5}
follows from the equality
\EqRef{int h2=int 12 xx+ 4}.

\Section{Calculations to Estimate Integral (\ref{eq: int h2=int 12 3x2})}

We reduce integrand of the integral
\EqRef{int h2=int 12 3x2}
\ShowEq{int h2=int 12 3x2 1}
The equality
\ShowEq{int h2=int 12 3x2 2}
follows from the equality
\EqRef{int h2=int 12 3x2 1}.
The equality
\ShowEq{int h2=int 12 3x2 3}
follows from the equality
\EqRef{int h2=int 12 3x2 2}.

\ePrints{1601.03259}
\ifx\Semafor\ValueOn
\Section{Example of Differential Form}
\else
\Section{Differential Form in Complex Field}
\fi
\labelSection{differential form in complex field}

According to the theorem
\RefTheorem{omega=aE+bI is integrable},
differential form
\DrawEq{omega=aE+bI}{example}
in complex field
is integrable iff
\DrawEq{d omega=0 condition}{example}
Let
\ShowEq{a=x2}
The equality
\ShowEq{da/d}
follows from the equality
\EqRef{a=x2}.
The equality
\ShowEq{da dx0+dx1}
follows from the equality
\EqRef{da/d}.
The equality
\ShowEq{db -dx0+dx1}
follows from equalities
\eqRef{d omega=0 condition}{example},
\EqRef{da dx0+dx1}.
Let
\ShowEq{db/dx0=0}
The equation
\ShowEq{db/dx1=}
follows from equalities
\EqRef{db -dx0+dx1},
\EqRef{db/dx0=0}.
The equality
\ShowEq{b=x1x1}
follows from the equation
\EqRef{db/dx1=}.
The equality
\ShowEq{omega(x)= example}
follows from equalities
\eqRef{omega=aE+bI}{example},
\EqRef{a=x2},
\EqRef{b=x1x1}.

Let the map
\ShowEq{f:A->B}fCC
be integral of the differential form
\EqRef{omega(x)= example}
\ShowEq{f=int o}
Equations
\ShowEq{df dx0 dx1}
follow from equalities
\EqRef{df/dx= complex field},
\EqRef{omega(x)= example}.
Equations
\ShowEq{df a b}
follow from equations
\EqRef{df dx0 dx1 a},
\EqRef{df dx0 dx1 b}.
The equality
\ShowEq{f(x0)+C(x1)}
follows from the equation
\EqRef{df dx0 a b}.
The equation
\ShowEq{dC/dx1}
follows from equalities
\EqRef{df dx1 a b},
\EqRef{f(x0)+C(x1)}.
The equality
\ShowEq{C(x1)=}
follows from the equation
\EqRef{dC/dx1}.
The equality
\ShowEq{f(x)=}
follow from equalities
\EqRef{f(x0)+C(x1)},
\EqRef{C(x1)=}.

%% file: Appendix.16.Calculate.Eq.tex

\DefEquation
{
\begin{split}
&\,(x-a)(a+t(x-a))^2
\\=&\,(x-a)(a^2+at(x-a)+t(x-a)a+t^2(x-a)^2
\\=&\,(x-a)a^2+(x-a)at(x-a)
\\+&\,(x-a)t(x-a)a+(x-a)t^2(x-a)^2
\\=&\,(x-a)a^2+t((x-a)a(x-a)+(x-a)(x-a)a)
\\+&\,t^2(x-a)(x^2-ax-xa+a^2)
\end{split}
}
{int h2=int 12 3x2 1}

\DefEquation
{
a=3(x^{\gi 0})^2+6x^{\gi 0}x^{\gi 1}i
}
{a=x2}

\DefEquation
{
\begin{split}
\frac{\partial a}{\partial x^{\gi 0}}&=6x^{\gi 0}+6x^{\gi 1}i
\\
\frac{\partial a}{\partial x^{\gi 1}}&=6x^{\gi 0}i
\end{split}
}
{da/d}

\DefEquation
{
\frac{\partial a}{\partial x^{\gi 0}}+i\frac{\partial a}{\partial x^{\gi 1}}
=6x^{\gi 0}+6x^{\gi 1}i-6x^{\gi 0}
=6x^{\gi 1}i
}
{da dx0+dx1}

\DefEquation
{
6x^{\gi 1}i
-\frac{\partial b}{\partial x^{\gi 0}}
+i\frac{\partial b}{\partial x^{\gi 1}}
=0
}
{db -dx0+dx1}

\DefEquation
{
\frac{\partial b}{\partial x^{\gi 0}}=0
}
{db/dx0=0}

\DefEquation
{
\frac{\partial b}{\partial x^{\gi 1}}=-6x^{\gi 1}
}
{db/dx1=}

\DefEquation
{
\begin{split}
\omega(x)&=(3(x^{\gi 0})^2+6x^{\gi 0}x^{\gi 1}i)\circ E+(-3(x^{\gi 1})^2)\circ I
\\&=\frac 34(4 x^2-(x-x^*)^2)\circ E+\frac 34(x-x^*)^2\circ I
\end{split}
}
{omega(x)= example}

\DefEq
{
\begin{align}
\frac{\partial f}{\partial x^{\gi 0}}
&=3(x^{\gi 0})^2-3(x^{\gi 1})^2+6x^{\gi 0}x^{\gi 1}i
\label{eq: df dx0 a b}
\\
i\frac{\partial f}{\partial x^{\gi 1}}
&=-3(x^{\gi 0})^2-3(x^{\gi 1})^2-6x^{\gi 0}x^{\gi 1}i
\nonumber
\\
\frac{\partial f}{\partial x^{\gi 1}}
&=-6x^{\gi 0}x^{\gi 1}+3(x^{\gi 0})^2i+3(x^{\gi 1})^2i
\label{eq: df dx1 a b}
\end{align}
}
{df a b}

\DefEq
{
\begin{align*}
f(x)
&=(x^{\gi 0})^3+3(x^{\gi 0})^2x^{\gi 1}i-3x^{\gi 0}(x^{\gi 1})^2+(x^{\gi 1})^3i+C
\\
&=x^3-\frac 18(x-x^*)^3+C
\end{align*}
}
{f(x)=}

\DefEquation
{
C_1(x^{\gi 1})
=(x^{\gi 1})^3i+C
}
{C(x1)=}

\DefEquation
{
\frac{\partial f}{\partial x^{\gi 1}}
=-6x^{\gi 0}x^{\gi 1}+3(x^{\gi 0})^2i+3(x^{\gi 1})^2i
=-6x^{\gi 0}x^{\gi 1}+3(x^{\gi 0})^2i+\frac{dC_1}{dx^{\gi 1}}
}
{dC/dx1}

\DefEquation
{
f(x)
=(x^{\gi 0})^3-3x^{\gi 0}(x^{\gi 1})^2+3(x^{\gi 0})^2x^{\gi 1}i+C_1(x^{\gi 1})
}
{f(x0)+C(x1)}

\DefEquation
{
b=-3(x^{\gi 1})^2
}
{b=x1x1}

\DefEq
{
\begin{align}
\frac 12\left(
\frac{\partial f}{\partial x^{\gi 0}}
-i\frac{\partial f}{\partial x^{\gi 1}}
\right)&=3(x^{\gi 0})^2+6x^{\gi 0}x^{\gi 1}i
\label{eq: df dx0 dx1 a}
\\
\frac 12\left(
\frac{\partial f}{\partial x^{\gi 0}}
+i\frac{\partial f}{\partial x^{\gi 1}}
\right)&=-3(x^{\gi 1})^2
\label{eq: df dx0 dx1 b}
\end{align}
}
{df dx0 dx1}

\DefEquation
{
\begin{split}
&\,(x-a)(a+t(x-a))^2
\\=&\,(xa^2-a^3+t((x-a)(ax-a^2)+(x-a)(xa-a^2)
\\+&\,t^2(x^3-xax-x^2a+xa^2-ax^2+a^2x+axa-a^3)
\\=&\,xa^2-a^3+t((x-a)ax-(x-a)a^2+(x-a)xa-(x-a)a^2)
\\+&\,t^2(x^3-xax-x^2a+xa^2-ax^2+a^2x+axa-a^3)
\end{split}
}
{int h2=int 12 3x2 2}

\DefEquation
{
\begin{split}
&\,(x-a)(a+t(x-a))^2
\\=&\,xa^2-a^3+t(xax-a^2x-xa^2+a^3+x^2a-axa-xa^2+a^3)
\\+&\,t^2(x^3-x^2a-xax-ax^2+xa^2+axa+a^2x-a^3)
\\=&\,xa^2-a^3+t(x^2a+xax-2xa^2-a^2x-axa+2a^3)
\\+&\,t^2(x^3-x^2a-xax-ax^2+xa^2+axa+a^2x-a^3)
\end{split}
}
{int h2=int 12 3x2 3}

\DefEquation
{
\begin{split}
&\,(a+t(x-a))^2(x-a)+(a+t(x-a))(x-a)(a+t(x-a))
\\+&\,(x-a)(a+t(x-a))^2
\\=&\,(a^2+ta(x-a)+t(x-a)a+t^2(x-a)^2)(x-a)+(a(x-a)
\\+&\,t(x-a)^2)(a+t(x-a))
\\+&\,(x-a)(a^2+ta(x-a)+t(x-a)a+t^2(x-a)^2)
\\=&\,(a^2+t(ax-a^2+xa-a^2)+t^2(x^2-xa-ax+a^2))(x-a)
\\+&\,(ax-a^2)(a+t(x-a))+t(x^2-xa-ax+a^2)(a+t(x-a))
\\+&\,(x-a)(a^2+t(ax-a^2+xa-a^2)+t^2(x^2-xa-ax+a^2))
\end{split}
}
{int h2=int 12 xx+ 1}

\DefEquation
{
\begin{split}
&\,(a+t(x-a))^2(x-a)+(a+t(x-a))(x-a)(a+t(x-a))
\\+&\,(x-a)(a+t(x-a))^2
\\=&\,a^2(x-a)+t(ax-a^2+xa-a^2)(x-a)
\\+&\,t^2(x^2-xa-ax+a^2)(x-a)+(ax-a^2)a
\\+&\,t(ax-a^2)(x-a)+t(x^2(a+t(x-a))-xa(a+t(x-a))
\\-&\,ax(a+t(x-a))+a^2(a+t(x-a)))+(x-a)a^2
\\+&\,t(x-a)(ax-a^2+xa-a^2)+t^2(x-a)(x^2-xa-ax+a^2)
\end{split}
}
{int h2=int 12 xx+ 2}

\DefEquation
{
\begin{split}
&\,(a+t(x-a))^2(x-a)+(a+t(x-a))(x-a)(a+t(x-a))
\\+&\,(x-a)(a+t(x-a))^2
\\=&\,a^2x-a^3+t(ax(x-a)-a^2(x-a)+xa(x-a)-a^2(x-a))
\\+&\,t^2(x^2(x-a)-xa(x-a)-ax(x-a)+a^2(x-a))
\\+&\,axa-a^3+t(ax(x-a)-a^2(x-a))+t(x^2a+tx^2(x-a)
\\-&\,(xa^2+txa(x-a))
\\-&\,axa-tax(x-a)+a^3+ta^2(x-a))
\\+&\,xa^2-a^3+t((x-a)ax-(x-a)a^2+(x-a)xa-(x-a)a^2)
\\+&\,t^2((x-a)x^2-(x-a)xa-(x-a)ax+(x-a)a^2)
\end{split}
}
{int h2=int 12 xx+ 21}

\DefEquation
{
\begin{split}
&\,(a+t(x-a))^2(x-a)+(a+t(x-a))(x-a)(a+t(x-a))
\\+&\,(x-a)(a+t(x-a))^2
\\=&\,a^2x-a^3+t(ax^2-axa-a^2x+a^3+xax-xa^2-a^2x+a^3)
\\+&\,t^2(x^3-x^2a-xax+xa^2-ax^2+axa+a^2x-a^3)
\\+&\,axa-a^3+t(ax^2-axa-a^2x+a^3+x^2a+t(x^3-x^2a)-xa^2
\\-&\,t(xax-xa^2)
-axa-t(ax^2-axa)+a^3+t(a^2x-a^3))
\\+&\,xa^2-a^3+t(xax-a^2x-xa^2+a^3+x^2a-axa-xa^2+a^3)
\\+&\,t^2(x^3-ax^2-xxa+axa-xax+a^2x+xa^2-a^3)
\end{split}
}
{int h2=int 12 xx+ 3}

\DefEquation
{
\begin{split}
&\,(a+t(x-a))^2(x-a)+(a+t(x-a))(x-a)(a+t(x-a))
\\+&\,(x-a)(a+t(x-a))^2
\\=&\,a^2x-a^3+xa^2-a^3+axa-a^3
\\+&\,t(ax^2-axa-a^2x+a^3+xax-xa^2-a^2x+a^3)
\\+&\,t(xax-a^2x-xa^2+a^3+x^2a-axa-xa^2+a^3)
\\+&\,t(ax^2-axa-a^2x+a^3+x^2a-xa^2-axa+a^3)
\\+&\,t^2(x^3-x^2a-xax+xa^2-ax^2+axa+a^2x-a^3)
\\+&\,t^2(x^3-x^2a-xax+xa^2-ax^2+axa+a^2x-a^3)
\\+&\,t^2(x^3-ax^2-x^2a+axa-xax+a^2x+xa^2-a^3)
\end{split}
}
{int h2=int 12 xx+ 4}

\DefEquation
{
\begin{split}
&\,(a+t(x-a))^2(x-a)+(a+t(x-a))(x-a)(a+t(x-a))
\\+&\,(x-a)(a+t(x-a))^2
\\=&\,xa^2+axa+a^2x-3a^3
\\+&\,t(2x^2a+2xax+2ax^2-4xa^2-4axa-4a^2x+6a^3)
\\+&\,t^2(3x^3-3x^2a-3xax+3xa^2-3ax^2+3axa+3a^2x-3a^3)
\end{split}
}
{int h2=int 12 xx+ 5}

%% file: Summary.Calculus.English.tex

\input{Summary.Calculus.Eq}

\chapter{Summary of Statements}

Let $D$ be the complete commutative ring of characteristic $0$.

\section{Table of Derivatives}

\begin{theorem}
\labelTheorem{derivative of const, algebra}
For any $b\in A$
\ShowEq{derivative of const, algebra}
\end{theorem}
\begin{proof}
Immediate corollary of definition
\RefDefinition{differentiable map, algebra}.
\end{proof}

\begin{theorem}
\labelTheorem{dx/dx=1}
\ShowEq{dx/dx=1}
\end{theorem}
\begin{proof}
According to the definition
\ShowEq{ref derivative, representation in algebra}
\ShowEq{dx/dx=}
The equality
\EqRef{dx/dx=1}
follows from
\EqRef{dx/dx=}.
\end{proof}

\begin{theorem}
\labelTheorem{derivative, product over constant, algebra}
For any $b$, $c\in A$
\ShowEq{derivative, product over constant, algebra}
\end{theorem}
\begin{proof}
Immediate corollary of equalities
\ShowEq{ref derivative of product}
and the theorem
\RefTheorem{derivative of const, algebra}.
\end{proof}

\begin{theorem}
\labelTheorem{d(f+g)=df+dg}
Let
\ShowEq{f:A->B}fAB
\ShowEq{f:A->B}gAB
be maps of Banach $D$\Hyph module $A$
into associative Banach $D$\Hyph algebra $A$.
Since there exist
the derivatives
\ShowEq{df,dg},
then there exists
the derivative
\ShowEq{d(f+g)}
\ShowEq{d(f+g)dx=df+dg}
\ShowEq{d(f+g)=df+dg}
Since
\ShowEq{df=+}
\ShowEq{dg=+}
then
\ShowEq{d(f+g)=+}
\end{theorem}
\begin{proof}
Since there exist
the derivatives
\ShowEq{df,dg},
then according to the definition
\ShowEq{ref derivative, representation in algebra}
\ShowEq{df=}
\ShowEq{dg=}
According to the definition
\ShowEq{ref derivative, representation in algebra}
\ShowEq{d(f+g)=}
The equality
\EqRef{d(f+g)dx=df+dg}
follows from
\EqRef{df=},
\EqRef{dg=},
\EqRef{d(f+g)=}.
The equality
\EqRef{d(f+g)=df+dg}
follows from
\ShowEq{ref sum of maps}
\EqRef{d(f+g)dx=df+dg}.
The equality
\EqRef{d(f+g)=+}
follows from
\EqRef{df=+},
\EqRef{dg=+},
\EqRef{d(f+g)=df+dg}.
\end{proof}

\begin{theorem}
\labelTheorem{derivative, fx=axb, algebra}
For any $b$, $c\in A$
\ShowEq{derivative, fx=axb, algebra}
\end{theorem}
\begin{proof}
Corollary of theorems
\RefTheorem{dx/dx=1},
\RefTheorem{derivative, product over constant, algebra},
when $f(x)=x$.
\end{proof}

\begin{theorem}
\labelTheorem{derivative of linear map}
Let $f$ be linear map
\ShowEq{derivative linear map associative algebra}
Then
\ShowEq{derivative linear map associative algebra, 1}
\end{theorem}
\begin{proof}
Corollary of theorems
\ShowEq{ref derivative of the sum}
\end{proof}

\begin{corollary}
For any $b\in A$
\ShowEq{derivative, xb-bx, algebra}
\qed
\end{corollary}

\begin{convention}
\labelConvention{set of permutations SO}
For $n\ge 0$, let
\ShowEq{set of permutations SO}
be set of permutations
\EqParm{s in SO(kn)}{K=n}
such that each permutation $\sigma$
preserves the order of variables $x_i$:
since $i<j$, then in the tuple
\ShowEq{s(y1kxn)}
$x_i$
precedes
$x_j$.
\qed
\end{convention}

\begin{lemma}
\labelLemma{enumerate set of permutations S(1,n)}
{\it
We can enumerate
the set of permutations $S(1,n)$ by index
\ShowEq{i 1n}
such way that
\StartLabelItem
\begin{enumerate}
\item
\ShowEq{s1(y)=y}
\item
Since $i>1$, then
\ShowEq{si(xi)=y}
\end{enumerate}
}
\end{lemma}
\begin{proof}
Since the order of variables
\ShowEq{x2...n}
in permutation
\ShowEq{s in S1n}
does not depend on permutation, permutation
\ShowEq{s in S1n}
are different by position which variable $y$ has.
Accordingly, we can enumerate
the set of permutations $S(1,n)$ by index
whose value corresponds to the number of position of variable $y$.
\end{proof}

The lemma
\RefLemma{enumerate set of permutations S(1,n)}
has simple interpretation.
Let $n-1$ white balls and $1$ black ball be in narrow box.
The black ball is the most left ball; white balls
are numbered from $2$ to $n$ in the order as we put them into the box.
The essence of the permutation $\sigma_k$ is that we take out
the black ball from the box and then we put it into cell with number $k$.
At the same time, white balls with the number not exceeding $k$ shift to the left.

\begin{lemma}
\labelLemma{S(1,n+1)=S(1,n)+}
{\it
For $n>0$, let
\ShowEq{S+(1,n)=}
Then
\ShowEq{S(1,n+1)=S(1,n)+}
}
\end{lemma}
\begin{proof}
Let
\ShowEq{sigma in S+}
According to the definition
\EqRef{S+(1,n)=},
there exists permutation
\ShowEq{tau in S1n}
such that
\ShowEq{tau x1n n+1}
According to the convention
\RefConvention{set of permutations SO},
the statement
\ShowEq{i<j<n+1}
implyes that in the tuple
\ShowEq{tau x1n n+1}
the variable $x_j$ is located between variables $x_i$ and $x_{n+1}$.
According to the convention
\RefConvention{set of permutations SO},
\ShowEq{sigma in S1n}
Therefore
\ShowEq{S+ subset S}
According to the lemma
\RefLemma{enumerate set of permutations S(1,n)},
the set
\ShowEq{S+(1,n)}
has $n$ permutations.

Let
\ShowEq{sigma=x2n y}
According to the convention
\RefConvention{set of permutations SO},
\ShowEq{sigma in S}
According to the definition
\EqRef{S+(1,n)=},
\ShowEq{sigma not in S+}

Therefore, we have listed $n+1$ elements of the set $S(1,n+1)$.
According to the lemma
\RefLemma{enumerate set of permutations S(1,n)},
the statement
\EqRef{S(1,n+1)=S(1,n)+}
follows from statements
\EqRef{S+ subset S},
\EqRef{sigma in S}.
\end{proof}

\begin{theorem}
\labelTheorem{dpn dx=+SO}
For any monomial
\ShowEq{pn(x)=a xn}
derivative has form
\ShowEq{dpn dx=+SO}
\end{theorem}
\begin{proof}
For $n=1$,
the map
\ShowEq{p1(x)=}
is linear map.
According to the theorem
\RefTheorem{derivative, fx=axb, algebra}
and convention
\RefConvention{set of permutations SO}
\ShowEq{dp1(x)=}

Let the statement be true for $n-1$
\ShowEq{dpn-1}
Since
\ShowEq{pn(x)=}
then according to the theorem
\ShowEq{ref derivative of product, theorem}
and the definition
\EqRef{pn(x)=}
\ShowEq{dpn(x)=}
The equality
\ShowEq{dpn(x)= 1}
follows from
\EqRef{dx/dx=1},
\EqRef{derivative, product over constant, algebra},
\EqRef{dpn-1},
\EqRef{dpn(x)=}.
The equality
\ShowEq{dpn(x)= 2}
follows from
\EqRef{dpn(x)= 1}
and multiplication rule of monomials (definitions
\RefDefinition{otimes -},
\RefDefinition{product of homogeneous polynomials}).
According to the lemma
\RefLemma{S(1,n+1)=S(1,n)+},
the equality
\EqRef{dpn dx=+SO}
follows from the equality
\EqRef{dpn(x)= 2}.
\end{proof}

\begin{lemma}
\labelLemma{sigma->(mu,nu)}
{\it
Let
\ShowEq{k<n}
For any permutation
\ShowEq{s in S(k+1,n)},
there exists unique pair of permutations
\ShowEq{mn in S(k,n)}
such that
\ShowEq{s=nu(mu)}
}
\end{lemma}
\begin{proof}
\end{proof}

\begin{lemma}
\labelLemma{(mu,nu)->sigma}
{\it
Let
\ShowEq{k<n}
For any pair of permutations
\ShowEq{mn in S(k,n)}
there exists unique permutation
\ShowEq{s in S(k+1,n)},
such that
\ShowEq{s=nu(mu)}
}
\end{lemma}
\begin{proof}
\end{proof}

\begin{theorem}
\labelTheorem{dkpn dx=+SO}
For any monomial
\ShowEq{pn(x)=a xn}
derivative of order $k$ has form
\ShowEq{dkpn dx=+SO}
\end{theorem}
\begin{proof}
For $k=1$, the statement of the theorem
is the statement of the theorem
\RefTheorem{dpn dx=+SO}.
Let the theorem be true for $k-1$. Then
\ShowEq{dk-1pn dx=+SO}
\end{proof}

\ifx\texFuture\Defined
\begin{theorem}
\labelTheorem{derivative of pn is symmetric, m < n, algebra}
The derivative
\ShowEq{derivative of pn is symmetric}
is symmetric polynomial with respect to variables $h_1$, ..., $h_m$.
\end{theorem}
\begin{proof}
To prove the theorem we consider algebraic properties
of the derivative and give equivalent definition.
We start from construction of monomial.

\begin{lemma}
{\it
We define the derivative of power $k$ according to rule
\ShowEq{derivative of monomial, algebraic definition, algebra}
}
\end{lemma}
According to construction, polynomial $r_n(h_1,...,h_k,x_{k+1},...,x_n)$
is symmetric with respect to variables $h_1$, ..., $h_k$, $x_{k+1}$, ..., $x_n$.
Therefore, polynomial
\EqRef{derivative of monomial, algebraic definition, algebra}
is symmetric with respect to variables $h_1$, ..., $h_k$.

For $k=1$, we will prove that definition
\EqRef{derivative of monomial, algebraic definition, algebra}
of the derivative coincides with definition
\ShowEq{ref derivative of map}

Let us prove now that definition
\EqRef{derivative of monomial, algebraic definition, algebra}
of the derivative coincides with definition
\ShowEq{ref derivative of Order n}
for $k>1$.

Let equality
\EqRef{derivative of monomial, algebraic definition, algebra}
be true for $k-1$.
Consider arbitrary monomial of polynomial $r_n(h_1,...,h_{k-1},x_k,...,x_n)$.
Identifying variables $h_1$, ..., $h_{k-1}$
with elements of division ring $D$, we consider polynomial
\begin{equation}
R_{n-k}(x_k,...,x_n)=r_n(h_1,...,h_{k-1},x_k,...,x_n)
\EqLabel{reduced polynom, algebra}
\end{equation}
Assume $P_{n-k}(x)=R_{n-k}(x_k,...,x_n)$, $x_k=...=x_n=x$.
Therefore
\ShowEq{derivative of pn is symmetric, 2}
According to definition
\ShowEq{ref derivative of Order n}
of the derivative
\ShowEq{derivative of Order n, 1, algebra}
According to definition
\EqRef{derivative of monomial, algebraic definition, algebra}
of the derivative
\ShowEq{derivative of monomial, algebraic definition, 1, algebra}
According to definition \EqRef{reduced polynom, algebra},
from the equality
\EqRef{derivative of monomial, algebraic definition, 1, algebra},
it follows that
\begin{equation}
\partial P_{n-k}(x)(h_k)=r_n(h_1,...,h_k,x_{k+1},...,x_n)\ \ \ x_{k+1}=x_n=x
\EqLabel{derivative of monomial, algebraic definition, 2, algebra}
\end{equation}
From comparison of \EqRef{derivative of Order n, 1, algebra}
and \EqRef{derivative of monomial, algebraic definition, 2, algebra}
it follows that
\[
\partial^k p_n(x)(h_1;...;h_k)=r_n(h_1,...,h_k,x_{k+1},...,x_n)\ \ \ x_{k+1}=x_n=x
\]
Therefore the equality
\EqRef{derivative of monomial, algebraic definition, algebra}
is true for any $k$ and $n$.

We proved the statement of theorem.
\end{proof}
\fi

\ShowTheorem{derivative x2}
\ShowProof{derivative x2}

\begin{remark}
The theorem
\RefTheorem{derivative x2}
also follows from the theorem
\RefTheorem{dpn dx=+SO}
since
\ShowEq{x2=1x1x1}
\qed
\end{remark}

\begin{theorem}
\labelTheorem{derivative x3}
Let $D$ be the complete commutative ring of characteristic $0$.
Let $A$ be associative Banach $D$\Hyph algebra.
Then
\ShowEq{derivative x3, algebra}
\ShowEq{derivative x3 differential, algebra}
\end{theorem}
\begin{proof}
According to the theorem
\RefTheorem{derivative of product, algebra},
\ShowEq{derivative x3 =}
The equality
\EqRef{derivative x3, algebra}
follows from the equality
\EqRef{derivative x3 =}.
The equality
\EqRef{derivative x3 differential, algebra}
follows from the equality
\EqRef{derivative x3, algebra}
and the definition
\RefDefinition{differential of map}.
\end{proof}

\begin{theorem}
Let $D$ be the complete commutative ring of characteristic $0$.
Let $A$ be associative division $D$\Hyph algebra.
Then\,\footnote{
The statement of the theorem is similar to example IX,
\citeBib{Hamilton Elements of Quaternions 1}, p. 451.
If product is commutative, then the equality
\EqRef{derivative x power -1, algebra}
gets form
\ShowEq{derivative x power -1, field}}
\ShowEq{derivative x power -1, algebra}
\ShowEq{differential x power -1, algebra}
\ShowEq{derivative x power -1 component, algebra}
\end{theorem}
\begin{proof}
Let us substitute $f(x)=x^{-1}$ in definition
\ShowEq{ref derivative, representation in algebra}
\ShowEq{derivative x power -1, algebra, 1}
The equality \EqRef{derivative x power -1, algebra}
follows from chain of equalities
\EqRef{derivative x power -1, algebra, 1}.
The equality
\EqRef{differential x power -1, algebra}
follows from the equality
\EqRef{derivative x power -1, algebra}
and the definition
\RefDefinition{differential of map}.
The equality
\EqRef{derivative x power -1 component, algebra}
follows from the equality
\EqRef{derivative x power -1, algebra}.
\end{proof}

\begin{theorem}
Let $D$ be the complete commutative ring of characteristic $0$.
Let $A$ be associative division $D$\Hyph algebra.
Then\,\footnote{
If product is commutative, then
\[
y=xax^{-1}=a
\]
Accordingly, the derivative is $0$.}
\ShowEq{derivative xax power -1, algebra}
\ShowEq{differential xax power -1, algebra}
\ShowEq{derivative xax power -1 component, algebra}
\end{theorem}
\begin{proof}
The equality \EqRef{derivative xax power -1, algebra}
is corollary of equalities
\ShowEq{ref derivative, fx=axb}
The equality
\EqRef{differential xax power -1, algebra}
follows from the equality
\EqRef{derivative xax power -1, algebra}
and the definition
\RefDefinition{differential of map}.
The equality
\EqRef{derivative xax power -1 component, algebra}
follows from the equality
\EqRef{derivative xax power -1, algebra}.
\end{proof}

\begin{theorem}
\labelTheorem{dex=ex}
\ShowEq{dex=ex}
\end{theorem}
\begin{proof}
The theorem follows from the theorem
\RefTheorem{exponent}.
\end{proof}

\begin{theorem}
\labelTheorem{dsh dch}
\ShowEq{dsh=ch}
\ShowEq{dch=sh}
\end{theorem}
\begin{proof}
The theorem follows from the theorem
\RefTheorem{sh ch}.
\end{proof}

\begin{theorem}
\labelTheorem{dsin dcos}
\ShowEq{dsin=cos}
\ShowEq{dcos=sin}
\end{theorem}
\begin{proof}
The theorem follows from the theorem
\RefTheorem{sin cos}.
\end{proof}

\section{Table of Integrals}

\begin{theorem}
Let the map
\ShowEq{f:A->B}fAB
be differentiable map.
Then
\ShowEq{int df=f}
\end{theorem}
\begin{proof}
The theorem follows from the definition
\RefDefinition{indefinite integral}.
\end{proof}

\begin{theorem}
\ShowEq{int=axb}
\ShowEq{int=axb=}
\end{theorem}
\begin{proof}
According to the definition
\RefDefinition{indefinite integral},
the map $y$ is integral
\EqRef{int=axb},
when the map $y$ satisfies to differential equation
\ShowEq{differential equation, additive function, algebra}
and initial condition
\DrawEq{differential equation, initial}{y=axb, algebra}
According to the theorem
\RefTheorem{derivative of const, algebra},
from the equality
\EqRef{differential equation, additive function, algebra},
it follows that the derivative of second order has form
\ShowEq{differential equation, additive function, 1, algebra}
The expansion into Taylor series
\ShowEq{differential equation, additive function, solution, algebra}
follows from
\ShowEq{differential equation y=axb, ref, algebra}
The equality
\EqRef{int=axb}
follows from equalities
\EqRef{differential equation, additive function, algebra},
\eqRef{differential equation, initial}{y=axb, algebra},
\EqRef{differential equation, additive function, solution, algebra}.
\end{proof}

\begin{remark}
According to the definition
\EqRef{a ox b c=},
we can present integral
\EqRef{int=axb}
following way
\ShowEq{int=axb 1}
\qed
\end{remark}

\begin{theorem}
\ShowEq{int ex=ex+C}
\end{theorem}
\begin{proof}
The theorem follows from the theorem
\RefTheorem{dex=ex}
and from the definition
\RefDefinition{indefinite integral}.
\end{proof}

\begin{remark}
According to the definition
\EqRef{a ox b c=},
we can present integral
\EqRef{int ex=ex+C}
following way
\ShowEq{int ex=ex+C 1}
\qed
\end{remark}

\begin{theorem}
\ShowEq{int sh=ch+C}
\ShowEq{int ch=sh+C}
\end{theorem}
\begin{proof}
The theorem follows from the theorem
\RefTheorem{dsh dch}
and from the definition
\RefDefinition{indefinite integral}.
\end{proof}

\begin{remark}
According to the definition
\EqRef{a ox b c=},
we can present integrals
\EqRef{int sh=ch+C},
\EqRef{int ch=sh+C}
following way
\ShowEq{int sh=ch+C 1}
\ShowEq{int ch=sh+C 1}
\qed
\end{remark}

\begin{theorem}
\ShowEq{int sin=cos+C}
\ShowEq{int cos=sin+C}
\end{theorem}
\begin{proof}
The theorem follows from the theorem
\RefTheorem{dsin dcos}
and from the definition
\RefDefinition{indefinite integral}.
\end{proof}

\begin{remark}
According to the definition
\EqRef{a ox b c=},
we can present integrals
\EqRef{int sin=cos+C},
\EqRef{int cos=sin+C}
following way
\ShowEq{int sin=cos+C 1}
\ShowEq{int cos=sin+C 1}
\qed
\end{remark}

%% file: Summary.Calculus.Eq.tex

\input{Summary.Calculus.Ref}

\DefEq
{
\[
\frac{db}{dx}=0
\]
}
{derivative of const, algebra}

\DefEquation
{
\frac{d^2 y}{d x^2}=0
}
{differential equation, additive function, 1, algebra}

\DefEq
{
$\displaystyle\frac{\partial^m p_n(x)}{\partial x^m}\circ(h_1;...;h_m)$
}
{derivative of pn is symmetric}

\DefEquation
{
\begin{split}
\frac{dp_n(x)}{dx}\circ dx
&=\sum_{\sigma\in SO(1,n)}(a_0\otimes...\otimes a_n)\circ
(\sigma(dx),\sigma(x_2),...,\sigma(x_n))
\\
x_2&=...=x_n=x
\end{split}
}
{dpn dx=+SO}

\DefEquation
{
\begin{split}
&\frac{d^k p_n(x)}{d x^k}\circ (dx^1; ... ;dx^k)
\\=&\sum_{\sigma\in SO(k,n)}(a_0\otimes...\otimes a_n,\sigma)\circ
(dx_1;...;dx_k;x_{k+1};...;x_n)
\\
x_{k+1}&=...=x_n=x
\end{split}
}
{dkpn dx=+SO}

\DefEquation
{
\begin{split}
&\frac{d^{k-1} p_n(x)}{d x^{k-1}}\circ (dx^1; ... ;dx^{k-1})
\\=&\sum_{\mu\in SO(k-1,n)}(a_0\otimes...\otimes a_n,\mu)\circ
(dx_1;...;dx_{k-1};x_k;...;x_n)
\\
x_k&=...=x_n=x
\end{split}
}
{dk-1pn dx=+SO}

\DefEq
{
\[
p_n(x)=(a_0\otimes...\otimes a_n)\circ x^n
\]
}
{pn(x)=a xn}

\DefEquation
{
\begin{split}
\frac{d p_{n-1}(x)}{dx}\circ dx
&=\sum_{\sigma\in S(1,n-1)}(a_0\otimes...\otimes a_{n-1})\circ \sigma(dx,x_2,...,x_{n-1})
\\
x_2&=...=x_{n-1}=x
\end{split}
}
{dpn-1}

\DefEquation
{
\frac{dp_n(x)}{dx}\circ dx
=\left(\frac{dp_{n-1}(x)}{dx}\circ dx\right)\ xa_n
+p_{n-1}(x)\ \left(\frac{dxa_n}{dx}\circ dx\right)
}
{dpn(x)=}

\DefEquation
{
\begin{split}
&\,\frac{dp_n(x)}{dx}\circ dx
\\=&\,\sum_{\sigma\in S(1,n-1)}(a_0\otimes...\otimes a_{n-1})
\circ (\sigma(dx),\sigma(x_2),...,\sigma(x_{n-1}))xa_n
\\+&\,((a_0\otimes...\otimes a_{n-1})\circ(x_2,...,x_n))dxa_n
\\x_2&\,=...=x_n=x
\end{split}
}
{dpn(x)= 1}

\DefEquation
{
\begin{split}
\frac{d p_n(x)}{dx}\circ dx
&=\sum_{\sigma\in S(1,n-1)}(a_0\otimes...\otimes a_n)
\circ (\sigma(dx),\sigma(x_2),...,\sigma(x_{n-1}),x_n))
\\&+(a_0\otimes...\otimes a_n)\circ(x_2,...,x_n,dx)
\\x_2&=...=x_n=x
\end{split}
}
{dpn(x)= 2}

\DefEquation
{
\frac{d x^3}{dx}=x^2\otimes 1+x\otimes x+1\otimes x^2
}
{derivative x3, algebra}

\DefEquation
{
dx^3=x^2\,dx+x\,dx\,x+dx\,x^2
}
{derivative x3 differential, algebra}

\DefEquation
{
\frac{dx^3}{dx}=\frac{dx^2x}{dx}
=\frac{dx^2}{dx}x+x^2\frac{dx}{dx}
=(x\otimes 1+1\otimes x)x+x^2(1\otimes 1)
}
{derivative x3 =}

\DefEq
{
\[
\S=
\begin{pmatrix}
y_1&...&y_k&x_{k+1}&...&x_{\K}
\\
\S(y_1)&...&\S(y_k)&\S(x_{k+1})&...&\S(x_{\K})
\end{pmatrix}
\]
}
{s in SO(kn)}

\DefEq
{
\[(\sigma(y_1),...,\sigma(y_k),\sigma(x_{k+1}),...,\sigma(x_n))\]
}
{s(y1kxn)}

\DefEq
{
\symb{SO(k,n)}{set of permutations}{1}
}
{set of permutations SO}

\DefEq
{
$\sigma_i(x_i)=y$.
\labelItem{si(xi)=y}
}
{si(xi)=y}

\DefEq
{
$\sigma_1(y)=y$.
\labelItem{s1(y)=y}
}
{s1(y)=y}

\DefEq
{
$i$, $1\le i\le n$,
}
{i 1n}

\DefEquation
{
p_n(x)=p_{n-1}(x)xa_n
}
{pn(x)=}

\DefEquation
{
\begin{matrix}
\displaystyle\frac{\partial^k p_n(x)}{\partial x^k}\circ(h_1;...;h_k)
=r_n(h_1,...,h_k,x_{k+1},...,x_n)
&
x_{k+1}=...=x_n=x
\end{matrix}
}
{derivative of monomial, algebraic definition, algebra}

\DefEq
{
$1<k<n$.
}
{k<n}

\DefEq
{
$\sigma\in SO(k+1,n)$
}
{s in S(k+1,n)}

\DefEq
{
$\mu\in SO(k,n)$, $\nu\in SO(1,n-k)$
}
{mn in S(k,n)}

\DefEq
{
\[
\sigma=(\mu(y_1),...,\mu(y_k),\nu(\mu(y_{k+1})),\nu(\mu(x_{k+2})),...,\nu(\mu(x_n)))
\]
}
{s=nu(mu)}

\DefEquation
{
\frac{dp_1(x)}{dx}\circ dx=a_0dxa_1=(a_0\otimes a_1)\circ dx
=\sum_{\sigma\in S(1,1)}(a_0\otimes a_1)\circ \sigma(dx)
}
{dp1(x)=}

\DefEq
{
\[p_1(x)=(a_0\otimes a_1)\circ x=a_0xa_1\]
}
{p1(x)=}

\DefEquation
{
S(1,n+1)=S^+(1,n)\cup\{(x_2,...,x_{n+1},y)\}
}
{S(1,n+1)=S(1,n)+}

\DefEq
{
\[x^2=(1\otimes 1\otimes 1)\circ x^2\]
}
{x2=1x1x1}

\DefEquation
{
S^+(1,n)\subseteq S(1,n+1)
}
{S+ subset S}

\DefEq
{
$\sigma\not\in S^+(1,n)$.
}
{sigma not in S+}

\DefEq
{
$\sigma=(x_2,...,x_{n+1},y)$.
}
{sigma=x2n y}

\DefEquation
{
(x_2,...,x_{n+1},y)\in S(1,n+1)
}
{sigma in S}

\DefEquation
{
S^+(1,n)=\{\sigma:\sigma=(\tau(y),\tau(x_2),...,\tau(x_n),x_{n+1}),\tau\in S(1,n)\}
}
{S+(1,n)=}

\DefEq
{
$S^+(1,n)$
}
{S+(1,n)}

\DefEq
{
$\sigma\in S^+(1,n)$.
}
{sigma in S+}

\DefEq
{
\[(\sigma(y),\sigma(x_1),...,\sigma(x_{n+1}))=
(\tau(y),\tau(x_2),...,\tau(x_n),x_{n+1})\]
}
{tau x1n n+1}

\DefEq
{
$i<j<n+1$
}
{i<j<n+1}

\DefEq
{
$\tau\in S(1,n)$
}
{tau in S1n}

\DefEq
{
$\sigma\in S(1,n+1)$.
}
{sigma in S1n}

\DefEq
{
$\sigma\in S(1,n)$
}
{s in S1n}

\DefEq
{
$x_2$, ..., $x_n$
}
{x2...n}

\DefEq
{
\[
P_{n-k}(x)=\frac{\partial^{k-1} p_n(x)}{\partial x^{k-1}}\circ(h_1;...;h_{k-1})
\]
}
{derivative of pn is symmetric, 2}

\DefEquation
{
\begin{split}
\frac{\partial P_{n-k}(x)}{\partial x}\circ h_k
=&\frac{\partial}{\partial x}
\left(\frac{\partial^{k-1} p_n(x)}{\partial x^{k-1}}
\circ(h_1;...;h_{k-1})\right)\circ h_k
\\
=&\frac{\partial^k p_n(x)}{\partial x^k}\circ(h_1;...;h_{k-1}; h_k)
\end{split}
}
{derivative of Order n, 1, algebra}

\DefEquation
{
\frac{\partial P_{n-k}(x)}{\partial x}\circ h_k
=R_{n-k}(h_k,x_{k+1},...,x_n)\ \ \ x_{k+1}=...=x_n=x
}
{derivative of monomial, algebraic definition, 1, algebra}

\DefEq
{
\begin{align*}
\frac{\partial f\circ x}{\partial x}&=f
\\
\frac{\partial f\circ x}{\partial x}\circ dx&=f\circ dx
\end{align*}
}
{derivative linear map associative algebra, 1}

\DefEq
{
\[
f\circ x=(a_{s\cdot 0}\otimes a_{s\cdot 1})\circ x
=a_{s\cdot 0}\ x\  a_{s\cdot 1}
\]
}
{derivative linear map associative algebra}

\DefEq
{
\[
\left\{
\begin{array}{r@{\ }lr@{\ }l}
\multicolumn{4}{c}
{
\begin{array}{r@{\ }l}
\displaystyle\frac{d (xb-bx)}{dx}
&=1\otimes b-b\otimes 1
\\
\displaystyle\frac{d (xb-bx)}{dx}\circ dx
&=dx\,b-b\,dx
\end{array}
}
\\
\displaystyle\VirtFrac
\frac{d\pC{1}{0} (xb-bx)}{dx}
&=1
&\displaystyle
\frac{d\pC{1}{1} (xb-bx)}{dx}
&=b
\\
\displaystyle\VirtFrac
\frac{d\pC{2}{0} (xb-bx)}{dx}
&=-b
&\displaystyle
\frac{d\pC{2}{1} (xb-bx)}{dx}
&=1
\end{array}
\right.
\]
}
{derivative, xb-bx, algebra}

\DefEquation
{
\frac{de^x}{dx}=\frac 12(e^x\otimes 1+1\otimes e^x)
}
{dex=ex}

\DefEquation
{
\frac{dx}{dx}\circ dx=dx
\ \ \ \ \frac{dx}{dx}=1\otimes 1
}
{dx/dx=1}

\DefEquation
{
\frac{d(f(x)+g(x))}{dx}\hphantom{\circ dx}
=\frac{df(x)}{dx}\hphantom{\circ dx}
+\frac{dg(x)}{dx}\hphantom{\circ dx}
}
{d(f+g)=df+dg}

\DefEquation
{
\frac{df(x)}{dx}=
\frac{d\pC{s}{0} f(x)}{dx}\otimes
\frac{d\pC{s}{1} f(x)}{dx}
}
{df=+}

\DefEquation
{
\frac{dg(x)}{dx}=
\frac{d\pC{t}{0} g(x)}{dx}\otimes
\frac{d\pC{t}{1} g(x)}{dx}
}
{dg=+}

\DefEquation
{
\frac{d(f(x)+g(x))}{dx}=
\frac{d\pC{s}{0} f(x)}{dx}\otimes
\frac{d\pC{s}{1} f(x)}{dx}
+\frac{d\pC{t}{0} g(x)}{dx}\otimes
\frac{d\pC{t}{1} g(x)}{dx}
}
{d(f+g)=+}

\DefEquation
{
\frac{d(f(x)+g(x))}{dx}\circ dx
=\frac{df(x)}{dx}\circ dx
+\frac{dg(x)}{dx}\circ dx
}
{d(f+g)dx=df+dg}

\DefEq
{
$\displaystyle\frac{df(x)}{dx}$, $\displaystyle\frac{dg(x)}{dx}$
}
{df,dg}

\DefEq
{
$\displaystyle\frac{d(f(x)+g(x))}{dx}$
}
{d(f+g)}

\DefEquation
{
\frac{dx}{dx}\circ dx
=\lim_{t\rightarrow 0,\ t\in R}(t^{-1}(x+tdx-x))
=dx
}
{dx/dx=}

\DefEquation
{
\frac{df(x)}{dx}\circ dx
=\lim_{t\rightarrow 0,\ t\in R}(t^{-1}(f(x+tdx)-f(x)))
}
{df=}

\DefEquation
{
\frac{dg(x)}{dx}\circ dx
=\lim_{t\rightarrow 0,\ t\in R}(t^{-1}(g(x+tdx)-g(x)))
}
{dg=}

\DefEquation
{
\begin{split}
&\,\frac{d(f(x)+g(x))}{dx}\circ dx
\\=&\,\lim_{t\rightarrow 0,\ t\in R}(t^{-1}(f(x+tdx)+g(x+tdx)-f(x)-g(x)))
\\=&\,\lim_{t\rightarrow 0,\ t\in R}(t^{-1}(f(x+tdx)-f(x)))
\\+&\,\lim_{t\rightarrow 0,\ t\in R}(t^{-1}(g(x+tdx)-g(x)))
\end{split}
}
{d(f+g)=}

\DefEquation
{
\frac{d\sinh x}{dx}=\frac 12(\cosh x\otimes 1+1\otimes \cosh x)
}
{dsh=ch}

\DefEquation
{
\frac{d\cosh x}{dx}=\frac 12(\sinh x\otimes 1+1\otimes \sinh x)
}
{dch=sh}

\DefEquation
{
\int(\sinh x\otimes 1+1\otimes\sinh x)\circ dx=2\cosh x+C
}
{int sh=ch+C}

\DefEquation
{
\int\sinh x\, dx+dx\,\sinh x=2\cosh x+C
}
{int sh=ch+C 1}

\DefEquation
{
\int(\cosh x\otimes 1+1\otimes\cosh x)\circ dx=2\sinh x+C
}
{int ch=sh+C}

\DefEquation
{
\int\cosh x\,dx+dx\,\cosh x=2\sinh x+C
}
{int ch=sh+C 1}

\DefEquation
{
\int(\sin x\otimes 1+1\otimes\sin x)\circ dx=-2\cos x+C
}
{int sin=cos+C}

\DefEquation
{
\int\sin x\, dx+dx\,\sin x=-2\cos x+C
}
{int sin=cos+C 1}

\DefEquation
{
\int(\cos x\otimes 1+1\otimes\cos x)\circ dx=2\sin x+C
}
{int cos=sin+C}

\DefEquation
{
\int\cos x\,dx+dx\,\cos x=2\sin x+C
}
{int cos=sin+C 1}

\DefEquation
{
\int(e^x\otimes 1+1\otimes e^x)\circ dx=2e^x+C
}
{int ex=ex+C}

\DefEquation
{
\int e^x\,dx+dx\,e^x=2e^x+C
}
{int ex=ex+C 1}

\DefEquation
{
\frac{d\sin x}{dx}=\frac 12(\cos x\otimes 1+1\otimes \cos x)
}
{dsin=cos}

\DefEquation
{
\frac{d\cos x}{dx}=-\frac 12(\sin x\otimes 1+1\otimes \sin x)
}
{dcos=sin}

\DefEq
{
\begin{align*}
dx^{-1}&=-x^{-2}dx
\\
\frac{dx^{-1}}{dx}&=-x^{-2}
\end{align*}
}
{derivative x power -1, field}

\DefEquation
{
\frac{dx^{-1}}{dx}=-x^{-1}\otimes x^{-1}
}
{derivative x power -1, algebra}

\DefEquation
{
d x^{-1}=-x^{-1}\ dx\ x^{-1}
}
{differential x power -1, algebra}

\DefEquation
{
\begin{matrix}
\displaystyle\frac{d\pC{1}{0} x^{-1}}{dx}=-x^{-1}
&
\displaystyle\frac{d\pC{1}{1} x^{-1}}{dx}=x^{-1}
\end{matrix}
}
{derivative x power -1 component, algebra}

\DefEquation
{
\begin{split}
\frac{dx^{-1}}{dx}\circ h
&=\lim_{t\rightarrow 0,\ t\in R}(t^{-1}((x+th)^{-1}-x^{-1}))
\\
&=\lim_{t\rightarrow 0,\ t\in R}(t^{-1}((x+th)^{-1}-x^{-1}(x+th)(x+th)^{-1}))
\\
&=\lim_{t\rightarrow 0,\ t\in R}(t^{-1}(1-x^{-1}(x+th))(x+th)^{-1})
\\
&=\lim_{t\rightarrow 0,\ t\in R}(t^{-1}(1-1-x^{-1}th)(x+th)^{-1})
\\
&=\lim_{t\rightarrow 0,\ t\in R}(-x^{-1}h(x+th)^{-1})
\end{split}
}
{derivative x power -1, algebra, 1}

\DefEquation
{
\frac{(xax^{-1})}{dx}=1\otimes ax^{-1}-xax^{-1}\otimes x^{-1}
}
{derivative xax power -1, algebra}

\DefEquation
{
d (xax^{-1})=dx\ ax^{-1}-xax^{-1}\ dx\ x^{-1}
}
{differential xax power -1, algebra}

\DefEquation
{
\left\{
\begin{matrix}
\displaystyle\frac{d\pC{1}{0} x^{-1}}{dx}=1
&
\displaystyle\frac{d\pC{1}{1} x^{-1}}{dx}=ax^{-1}
\\
\displaystyle\frac{d\pC{2}{0} x^{-1}}{dx}=-xax^{-1}
&
\displaystyle\frac{d\pC{2}{1} x^{-1}}{dx}=x^{-1}
\end{matrix}
\right.
}
{derivative xax power -1 component, algebra}

\DefEquation
{
\left\{
\begin{array}{r@{\ }lr@{\ }l}
\displaystyle\frac{dbxc}{dx}&=b\otimes c
&
\displaystyle\frac{dbxc}{dx}\circ dx&=b\,dx\,c
\\[10pt]
\displaystyle
\frac{d\pC{1}{0} bxc}{dx}&=b
&
\displaystyle
\frac{d\pC{1}{1} bxc}{dx}&=c
\end{array}
\right.
}
{derivative, fx=axb, algebra}

\DefEquation
{
\left\{
\begin{aligned}
\frac{d bf(x)c}{d x}
&=(b\otimes c)\circ\frac{d f(x)}{d x}
\\
\frac{dbf(x)c}{dx}\circ dx
&=b\left(\frac{df(x)}{dx}\circ dx\right)c
\\
\VirtFrac
\frac{d\pC{s}{0} bf(x)c}{dx}
&
=b\frac{d\pC{s}{0} f(x)}{dx}
\\
\VirtFrac
\frac{d\pC{s}{1} bf(x)c}{dx}
&
=\frac{d\pC{s}{1} f(x)}{dx}c
\end{aligned}
\right.
}
{derivative, product over constant, algebra}

\DefEquation
{
y=f\pC{s}{0}\ x\ f\pC{s}{1}+C
}
{differential equation, additive function, solution, algebra}

\DefEquation
{
\int f\pC{s}{0}\,dx\,f\pC{s}{1} =f\pC{s}{0}\,x\,f\pC{s}{1}+C
}
{int=axb 1}

\DefEquation
{
\int\frac{df(x)}{dx}\circ dx=f(x)+C
}
{int df=f}

\DefEquation
{
\int(f\pC{s}{0}\otimes f\pC{s}{1})\circ dx
=(f\pC{s}{0}\otimes f\pC{s}{1})\circ x+C
}
{int=axb}

\DefEq
{
\[
f\pC{s}{0}\in A\ \ \ f\pC{s}{1}\in A
\]
}
{int=axb=}

\DefEquation
{
\frac{dy}{dx}=f\pC{s}{0}\otimes f\pC{s}{1}
}
{differential equation, additive function, algebra}

%% file: Summary.Calculus.Ref.tex

\DefEq
{
\ePrints{}%
\ifx\Semafor\ValueOn%
\EqRef[MVector]{sum of maps, D module},
\else%
\EqRef{sum of maps, D module},
\fi%
}
{ref sum of maps}

\DefEq
{
\RefTheorem{d(f+g)=df+dg},
\RefTheorem{derivative, fx=axb, algebra},
\ePrints{}%
\ifx\Semafor\ValueOn%
\RefTheorem[1601.03259]{derivative of the sum}.
\else%
\RefTheorem{derivative of the sum}.
\fi%
}
{ref derivative of the sum}

\DefEq
{
\ePrints{}%
\ifx\Semafor\ValueOn%
\EqRef[1601.03259]{derivative of product, algebra},
\EqRef[1601.03259]{component of derivative, fg, algebra},
\else%
\EqRef{derivative of product, algebra},
\EqRef{component of derivative, fg, algebra},
\fi%
}
{ref derivative of product}

\DefEq
{
\ePrints{}%
\ifx\Semafor\ValueOn%
\RefTheorem[1601.03259]{derivative of product, algebra}
\else%
\RefTheorem{derivative of product, algebra}
\fi%
}
{ref derivative of product, theorem}

\DefEq
{
\ePrints{}%
\ifx\Semafor\ValueOn%
\EqRef[1601.03259]{derivative of Order n, algebra}
\else%
\EqRef{derivative of Order n, algebra}
\fi%
}
{ref derivative of Order n}

\DefEq
{
\ePrints{}%
\ifx\Semafor\ValueOn%
\EqRef[1601.03259]{derivative of map, algebra}.
\else%
\EqRef{derivative of map, algebra}.
\fi%
}
{ref derivative of map}

\DefEq
{
\ePrints{}%
\ifx\Semafor\ValueOn%
\EqRef[1601.03259]{derivative of product, algebra},
\else%
\EqRef{derivative of product, algebra},
\fi%
\EqRef{derivative, fx=axb, algebra}.
}
{ref derivative, fx=axb}

\DefEq
{
\EqRef{differential equation, additive function, algebra},
\eqRef{differential equation, initial}{y=axb, algebra},
\EqRef{differential equation, additive function, 1, algebra}.
}
{differential equation y=axb, ref, algebra}

%% file: Biblio.English.tex
\OpenBiblio


\BiblioItem{Doctor Ouch}
{
Kornei Chukovsky. Doctor Ouch.
\\
Translator and illustrator Jan Seabaugh.
\\
Viveca Smith Publishing, 2004, ISBN-10: 0974055107.
}%

\BiblioItem{Einstein: Electrodynamics of Moving Bodies}
{
Albert Einstein,
On the Electrodynamics of Moving Bodies, 1905,
\\
The Principle of Relativity: A Collection of Original
Memoirs on the Special and General Theory of Relativity , 37 - 65,
\\
Courier Dover Publications, 1952; ISBN-13: 978-0486600819
\\
Zur Elektrodynamik der bewegter K\"orper. Ann. Phys., 1905, 17, 891-921. 
}%

\BiblioItem{Einstein: On the Relativity Principle}
{
Albert Einstein,
On the Relativity Principle and the Conclusions Drawn from It, 1907,
\\
The Collected Papers of Albert Einstein, Volume 2:
The Swiss Years: Writings, 1900-1909. English translation. 252 - 311.
\\
Anna Beck, translator; Peter Havas, consultant.
Princeton University Press, 1989; ISBN-13: 9780691085494
\\
\"Uber das Relativit\"atsprinzip und die aus demselben gezogenen Folgerungen. 
Jahrb. d. Radioaktivit\"at u. Elektronik, 1907, 4, 411-462. 
}%

\BiblioItem{Einstein: Foundations of general relativity}
{
Albert Einstein,
Die Grundlage der allgemeinen Relativit\"atstheorie,
Ann. Phys., 1916, {\bf 49}, 769 - 822,\\
Einstein's Annalen Papers: The Complete Collection 1901-1922,
edited by J\"urgen Renn, 517 - 571,\\
Wiley-VCH Verlag GmbH \& Co. KGaA, 2005
}%

\BiblioItem{Einstein: Geometry and Experience}
{
Albert Einstein, Geometry and Experience, (1921)\\
Albert Einstein, Sidelights on Relativity, 25 - 56,\\
Courier Dover Publications, 1983
}%

\BiblioItem{Einstein: Main problems of general relativity}
{
Albert Einstein,
Grundgedanken und Probleme der Relativit\"atstheorie, (1923),\\
Nobelstiftelsen, Les Prix Nobel en 1921 - 1922,
Imprimerie Royale, Stockholm, 1923
}%

\BiblioItem{Einstein: Noneuclidean Geometry and Physics}
{
Albert Einstein,
Nichtenklidische Geometrie in der Physik Neue Rundschan, (1925)
Berlin, S. 16 - 20
}%

\BiblioItem{Einstein: Isaak Newton}
{
Albert Einstein,
Isaak Newton, 1927,
Out of My Later Years, 
Citadel Press, 1995, 219 - 223
}%

\BiblioItem{Einstein: On Science}
{
Albert Einstein,
On Science, 
Cosmic Religion, with Other Opinions and Aphorisms,142 - 146,
New York, 1931, 97 - 103
}%

\BiblioItem{Einstein: Autobiographical Notes}
{
Albert Einstein,
Autobiographical Notes, 1949,\\
Paul A. Schilpp, editor, Albert Einstein: Philosopher-Scientist,
Evanston, 
Illinois, The Library of Living Philosophers, 1949, 1 - 95
}%

\BiblioItem{Feynman 1}
{
Richard Phillips Feynman, Robert B. Leighton, Matthew Linzee Sands.
The Feynman lectures on physics: Volume 1.
Mainly Mechanics, Radiation, and Heat.
Addison\Hyph Wesley, 1965.
}%

\BiblioItem{0538731877}
{
James Shipman, Jerry D. Wilson and Aaron Todd.
Introduction to Physical Science.
Cengage Learning, 2009; ISBN 0538731877.
}%

\BiblioItem{Cite: 104}
{
Cite 104, Source unknown
}%

\BiblioItem{Ghez}
{
Ghez et al.,
The First Measurement of Spectral Lines in a Short-Period Star Bound to the Galaxy's Central Black Hole: A Paradox of Youth,
\href{http://www.journals.uchicago.edu/ApJ/journal/issues/ApJL/v586n2/16990/brief/16990.abstract.html}{ApJL, 586, L127} (2003),
eprint \href{http://arxiv.org/abs/astro-ph/0302299}{arXiv:astro-ph/0302299} (2003)
}%

\BiblioItem{Schodel}
{
R. Sch\"odel et al.,
A star in a 15.2-year orbit around the supermassive black hole at the centre of the Milky Way,
\href{http://www.nature.com/cgi-taf/DynaPage.taf?file=/nature/journal/v419/n6908/abs/nature01121_fs.html}{Nature 419, 694} (2002)
}%

\BiblioItem{Mielke}
{
Eckehard W. Mielke, Affine generalization of the Komar complex of general relativity,
\href{http://prola.aps.org/searchabstract/PRD/v63/i4/e044018}{Phys. Rev. D 63, 044018} (2001)
}%

\BiblioItem{Obukhov}
{
Yu. N. Obukhov and J. G. Pereira, Metric\hyph affine approach to teleparallel gravity,
\href{http://scitation.aip.org/getabs/servlet/GetabsServlet?prog=normal&id=PRVDAQ000067000004044016000001&idtype=cvips&gifs=Yes}
{Phys. Rev. D 67, 044016} (2003),
eprint \href{http://arxiv.org/abs/gr-qc/0212080}{arXiv:gr-qc/0212080} (2002)
}%

\BiblioItem{Sardanashvily}
{
Giovanni Giachetta, Gennadi Sardanashvily, Dirac Equation in Gauge and Affine-Metric Gravitation Theories,
eprint \href{http://arxiv.org/abs/gr-qc/9511035}{arXiv:gr-qc/9511035} (1995)
}%

\BiblioItem{Gauge}
{
Frank Gronwald and Friedrich W. Hehl, On the Gauge Aspects of Gravity, eprint
\href{http://arxiv.org/abs/gr-qc/9602013}{arXiv:gr-qc/9602013} (1996)
}%

\BiblioItem{Neeman}
{
Yuval Neeman, Friedrich W. Hehl, Test Matter in a Spacetime with Nonmetricity, eprint
\href{http://arxiv.org/abs/gr-qc/9604047}{arXiv:gr-qc/9604047} (1996)
}%

\BiblioItem{torsion}
{
F. W. Hehl, P. von der Heyde, G. D. Kerlick, and J. M. Nester,
General relativity with spin and torsion: Foundations and prospects,\\
\href{http://prola.aps.org/abstract/RMP/v48/i3/p393_1}{Rev. Mod. Phys. 48, 393} (1976)
}%

\BiblioItem{Megged}
{
O. Megged, Post-Riemannian Merger of Yang-Mills Interactions with Gravity,
eprint \href{http://arxiv.org/abs/hep-th/0008135}{arXiv:hep-th/0008135} (2001)
}%


\BiblioItem{gr-qc-9604027}
{
Yu.N. Obukhov, E.J. Vlachynsky, W. Esser, R. Tresguerres and F.W. Hehl,
An exact solution of the metric\hyph affine gauge theory with dilation, shear, and spin charges,
eprint \href{http://arxiv.org/abs/gr-qc/9604027}{arXiv:gr-qc/9604027} (1996)
}%

\BiblioItem{4419-7514}
{
Mari\'an Fabian, Petr Habala, Petr H\'ajek, Vicente Montesinos, V\'aclav Zizler.
Banach Space Theory: The Basis for Linear and Nonlinear Analysis.
\\
Springer; New York, 2010; ISBN-13: 978-1441975140
}%

\BiblioItem{Weinberg I}
{
Steven Weinberg.
The Quantum Theory of Fields. Volume I. Foundations.
Cambridge university press, 1995
}%

\BiblioItem{Weinberg II}
{
Steven Weinberg.
The Quantum Theory of Fields. Volume II. Modern applications.
Cambridge university press, 1996
}%

\BiblioItem{Reinhardt}
{
Walter Greiner, Joachim Reinhardt. Field Quantization. Springer.
}%

\BiblioItem{978-3540875604}
{
Walter Greiner, Joachim Reinhardt. Quantum Electrodynamics. Springer, 2009.
}%

\BiblioItem{978-1898563020}
{
H. Robert Mills. Practical Astronomy. Woodhead Publishing, 1994. ISBN-13: 978-1898563020.
}%

\BiblioItem{Landau I}
{
L. D. Landau, E. M. Lifshich.
Course of theoretical physics, volume 1.
Mechanics.
\\
Translated from the Russian by J. B. Sykes and J. S. Bell.
Pergamon Press, 1969
}%

\BiblioItem{Landau}
{
L. D. Landau, E. M. Lifshich, The classical theory of fields.
\\
Translated from the Russian by Morton Hamermesh.
Pergamon Press, 1971
}%

\BiblioItem{Landau III}
{
L. D. Landau, E. M. Lifshich,
Course of Theoretical Physics, Volume 3.
Quantum Mechanics Non-Relativistic Theory, Third Edition.
\\
Translated from the Russian by J. B. Sykes and J. S. Bell.
Butterworth-Heinemann, 1981, ISBN 978-0750635394.
}%

\BiblioItem{Wheeler}
{
Ignazio Ciufolini, John Wheeler. Gravitation and Inertia.
Princeton university press.
}%

\BiblioItem{Gravitation MTW}
{
Charles W. Misner, Kip S. Thorne, John Archibald Wheeler.
Gravitation.
W. H. Freeman and Company, San Francisco, 1973.
}%

\BiblioItem{Anderson98}
{
J. D. Anderson, P. A. Laing, E. L. Lau, A. S. Liu, M. M. Nieto, and S. G. Turyshev,
Indication, from Pioneer 10/11, Galileo, and Ulysses Data, of an Apparent Anomalous, Weak, Long-Range Acceleration,
\href{http://prola.aps.org/abstract/PRL/v81/i14/p2858_1}{Phys. Rev. Lett. 81, 2858}, (1998),
eprint \href{http://arxiv.org/abs/gr-qc/9808081}{arXiv:gr-qc/9808081} (1998)
}%

\BiblioItem{Anderson02}
{
J. D. Anderson, P. A. Laing, E. L. Lau, A. S. Liu, M. M. Nieto, and S. G. Turyshev,
Study of the anomalous acceleration of Pioneer 10 and 11,
\href{http://prola.aps.org/searchabstract/PRD/v65/i8/e082004}{Phys. Rev. D 65, 082004, 50 pp.}, (2002),
eprint \href{http://arxiv.org/abs/gr-qc/0104064}{arXiv:gr-qc/0104064} (2001)
}%


\BiblioItem{H. Aslaksen}
{
H. Aslaksen.  Quaternionic determinants \textit{Math.
Intelligencer} {\bf 18}(3), pp.57-65, (1996).
}%

\BiblioItem{L. Chen: Definition of determinant}
{
L. Chen, Definition of determinant and Cramer solutions over
quaternion field, \textit{Acta Math. Sinica (N.S.)} {\bf 7},
pp.171-180, (1991).
}%

\BiblioItem{L. Chen: Inverse matrix}
{
L. Chen,
Inverse matrix and properties of double determinant over quaternion
field, \textit{Sci. China, Ser. A} {\bf 34}, pp.528-540, (1991).
}%

\BiblioItem{N. Cohen S. De Leo}
{
N. Cohen, S. De Leo, The quaternionic determinant, \textit{The Electronic Journal Linear
Algebra} {\bf 7}, pp.100-111, (2000).
}%

\BiblioItem{Dyson: Quaternion determinants}
{
F. J. Dyson, Quaternion determinants, \textit{Helvetica Phys.
Acta} {\bf 45}, pp. 289-302, (1972).
}%

\BiblioItem{Melvin Hausner}
{
Melvin Hausner,
A Vector Space Approach to Geometry,
Dover Publications, 1998
}%

\BiblioItem{Serge Lang}
{
Serge Lang,
Algebra, Springer, 2002
}%

\BiblioItem{9780534423230}
{
Charles Lanski.
Concepts In Abstract Algebra.
American Mathematical Soc., 2005, ISBN 978-0534423230
}%

\BiblioItem{Burris Sankappanavar}
{
S. Burris, H.P. Sankappanavar,
A Course in Universal Algebra, Springer-Verlag (March, 1982),
\\eprint
\href{http://www.math.uwaterloo.ca/~snburris/htdocs/ualg.html}
{http://www.math.uwaterloo.ca/~snburris/htdocs/ualg.html}
\\(The Millennium Edition)
}%

\BiblioItem{Shilov single 12}
{
G. E. Shilov,
Calculus, Single Variable Functions, Parts 1 - 2,
Moscow, Nauka, 1969
}%

\BiblioItem{Shilov single 3}
{
G. E. Shilov,
Calculus, Single Variable Functions, Part 3,
Moscow, Nauka, 1970
}%

\BiblioItem{Shilov}
{
G. E. Shilov,
Calculus, Multivariable Functions,
Moscow, Nauka, 1972
}%

\BiblioItem{Kolmogorov Fomin}
{
A. N. Kolmogorov and S. V. Fomin.
Introductory Real Analysis.
\\
Translated and edited by Richard A. Silverman.
\\
Dover Publication, 1975, ISBN-13: 978-0486612263
}%

\BiblioItem{Lebedev Vorovich}
{
L. P. Lebedev, I. I. Vorovich,
Functional Analysis in Mechanics,
Springer, 2002
}%

\BiblioItem{8176-4374}
{
Mariano Giaquinta, Giuseppe Modica,
Mathematical Analysis: Linear and Metric Structures and Continuity.
\\
Springer, 2007, ISBN-13: 978-0-8176-4374-4
}%

\BiblioItem
{Rashevsky}
{
P. K. Rashevsky, Riemann Geometry and Tensor Calculus,\\
Moscow, Nauka, 1967
}%

\BiblioItem
{Kurosh: High Algebra}
{
A. G. Kurosh, Higher Algebra,
\\
George Yankovsky translator,
\\
Mir Publishers, 1988, ISBN: 978-5030001319
}%

\BiblioItem
{Kurosh: General Algebra}
{
A. G. Kurosh, Lectures on General Algebra,
Chelsea Pub Co, 1965 
}%

\BiblioItem{Sabinin: Smooth Quasigroups}
{
Lev V. Sabinin, Smooth Quasigroups and Loops,
Kluwer Academic Publisher, 1999 
}%

\BiblioItem{978-0-8176-8384-9}
{
Garret Sobczyk, New Foundations in Mathematics: The Geometric Concept of Number,
\\
Springer, 2013, ISBN: 978-0-8176-8384-9
}%

\BiblioItem{Dubrovin Fomenko Novikov part 1}
{
B. A. Dubrovin, A. T. Fomenko, S. P. Novikov,
Modern Geometry - Methods and Applications,\\
Part I, The Geometry of Surfaces, Transformation Groups, and Fields,\\
Translated by Robert G. Burns,\\
Springer - New York, 1992
}%

\BiblioItem{Dubrovin Fomenko Novikov part 2}
{
B. A. Dubrovin, A. T. Fomenko, S. P. Novikov,
Modern Geometry - Methods and Applications,
Part II: The Geometry and Topology of Manifolds,\\
Translated by Robert G. Burns,\\
Springer - New York, 1985
}%

\BiblioItem{Kobayashi Nomizu vol 1}
{
Kobayashi S, Nomizu K,
Foundations of Differential Geometry, volume I,\\
Interscience Publishers, 1963
}%

\BiblioItem{Lichnerowicz}
{
Andre Lichnerowicz,
Global Theory of Connections and Holonomy Groups,\\
Kluwer Academic Publishers, 1976, ISBN-13: 978-9028604964
}%

\BiblioItem{Korn}
{
Granino A. Korn, Theresa M. Korn,
Mathematical Handbook for Scientists and Engineer,
McGraw-Hill Book Company, New York, San Francisco,
Toronto, London, Sydney, 1968
}%

\BiblioItem{Hocking Young Topology}
{
John G. Hocking, Gail S. Young,
Topology,\\
Courier Dover Publications, 1988
}%

\BiblioItem{Olver: Lie groups to differential equations}
{
Peter J. Olver,
Applications of Lie groups to differential equations,\\
Springer, 2000
}%

\BiblioItem{1708.01190}
{
Nathan BeDell,
Doing Algebra over an Associative Algebra,
\\
eprint \href{https://arxiv.org/abs/1708.01190}{arXiv:1708.01190} (2017)
}%

\BiblioItem{Tartaglia}
{
Angelo Tartaglia and Matteo Luca Ruggiero,
Angular Momentum Effects in Michelson\Hyph Morley Type Experiments,
Gen.Rel.Grav. 34, 1371-1382 (2002),\\
eprint \href{http://arxiv.org/abs/gr-qc/0110015}{arXiv:gr-qc/0110015} (2001)
}%

\BiblioItem{Tomozawa}
{
Yukio Tomozawa, Speed of Light in Gravitational Fields, eprint
\href{http://arxiv.org/abs/astro-ph/0303047}{arXiv:astro-ph/0303047} (2004)
}%

\BiblioItem{Magueijo}
{
Joao Magueijo,
Covariant and locally Lorentz-invariant varying speed of light theories,
\href{http://prola.aps.org/abstract/PRD/v62/i10/e103521}{Phys. Rev. D 62, 103521} (2000),
eprint \href{http://arxiv.org/abs/gr-qc/0007036}{arXiv:gr-qc/0007036} (2000)
}%

\BiblioItem{Bassett}
{
Bruce A. Bassett, Stefano Liberati, Carmen Molina-Paris, and Matt Visser,
Geometrodynamics of variable-speed-of-light cosmologies,
\href{http://prola.aps.org/abstract/PRD/v62/i10/e103518}{Phys. Rev. D 62}, 103518 (2000),
eprint \href{http://arxiv.org/abs/astro-ph/0001441}{arXiv:astro-ph/0001441} (2000)
}%

\BiblioItem{C.A. Deavours The Quaternion Calculus}
{
C.A. Deavours, The Quaternion Calculus, 
American Mathematical Monthly, {\bf 80} (1973), pp. 995 - 1008
}%

\BiblioItem{Straumann}
{
Lochlainn O'Raifeartaigh and Norbert Straumann,
Gauge theory: Historical origins and some modern developments,
\href{http://prola.aps.org/abstract/RMP/v72/i1/p1_1}{Rev. Mod. Phys. 72, 1} (2000)
}%

\BiblioItem{Lammerzahl}
{
Claus L\"ammerzahl, Mark P. Haugan,
On the interpretation of Michelson\Hyph Morley experiments,
{Phys. Lett. A282 223-229} (2001),\\
eprint \href{http://arxiv.org/abs/gr-qc/0103052}{arXiv:gr-qc/0103052} (2001)
}%

\BiblioItem{0305117}
{
Holger Mueller, Sven Herrmann, Claus Braxmaier, Stephan Schiller, Achim Peters.
Modern Michelson-Morley Experiment using Cryogenic Optical Resonators.
eprint \href{http://arxiv.org/abs/physics/0305117}{arXiv:physics/0305117} (2003)
\\
Phys. Rev. Lett. 91:020401, 2003
}%

\BiblioItem{0706.2031}
{
Holger Mueller, Paul Louis Stanwix, Michael Edmund Tobar,
Eugene Ivanov, Peter Wolf, Sven Herrmann, Alexander Senger,
Evgeny Kovalchuk, Achim Peters.
Relativity tests by complementary rotating Michelson-Morley experiments.
eprint \href{http://arxiv.org/abs/0706.2031}{arXiv:0706.2031 [physics.class-ph]} (2006)
\\
Phys. Rev. Lett. 99:050401, 2007
}%

\BiblioItem{1008.1205}
{
M. Nagel, K. M\"ohle, K. D\"oringshoff, S. Herrmann, A. Senger, E.V. Kovalchuk, A. Peters.
Testing Lorentz Invariance by Comparing Light Propagation in Vacuum and Matter.
eprint \href{http://arxiv.org/abs/1008.1205}{arXiv:1008.1205 [physics.ins-det]} (2010)
}%

\BiblioItem{1109.4897}
{
The OPERA Collaboration.
Measurement of the neutrino velocity with the OPERA detector in the CNGS beam.
eprint \href{http://arxiv.org/abs/1109.4897}{arXiv:1109.4897 [hep-ex]} (2011)
}%

\BiblioItem{Ranada}
{
Antonio F. Ranada,
Pioneer acceleration and variation of light speed: experimental situation,
eprint \href{http://arxiv.org/abs/gr-qc/0402120}{arXiv:gr-qc/0402120} (2004)
}%

\BiblioItem{Gelfand Minlos: rotation and Lorentz groups}
{
Izrail Moiseevich Gelfand, Robert Adolfovich Minlos,
Representations of the rotation and Lorentz groups and their applications;\\
Engl. transl. ed. H. K. Farahat; Transl. by G. Cummins and T. Boddongton;\\
Pergamon Press, 1963
}%

\BiblioItem{math.QA-0208146}
{
I. Gelfand, S. Gelfand, V. Retakh, R. Wilson,
Quasideterminants,\\
eprint \href{http://arxiv.org/abs/math.QA/0208146}{arXiv:math.QA/0208146} (2002)
}%

\BiblioItem{q-alg-9705026}
{
I. Gelfand, V. Retakh,
Quasideterminants, I,\\
eprint \href{http://arxiv.org/abs/q-alg/9705026}{arXiv:q-alg/9705026} (1997)
}%

\BiblioItem{Gelfand Retakh 1991}
{
I. Gelfand and V. Retakh, Determinants of Matrices over Noncommutative Rings, Funct.
Anal. Appl. 25 (1991), no. 2, 91-102
}%

\BiblioItem{Gelfand Retakh 1992}
{
I. Gelfand and V. Retakh, A Theory of Noncommutative Determinants and Characteristic
Functions of Graphs, Funct. Anal. Appl. 26 (1992), no. 4, 1-20
}%

\BiblioItem{hep-th-9407124}
{
I. M. Gelfand, D. Krob, A. Lascoux, B. Leclerc, V.S. Retakh and J.-Y. Thibon,
Noncommutative symmetric functions,\\
eprint \href{http://arxiv.org/abs/hep-th/9407124}{arXiv:hep-th/9407124} (1994)
}%

\BiblioItem{0911.4454}
{
Vladimir Retakh,
From factorizations of noncommutative polynomials to combinatorial topology,\\
eprint \href{http://arxiv.org/abs/0911.4454}{arXiv:0911.4454} (2009)
}%

\BiblioItem{Naimark Shtern: Theory of group representations}
{
Mark Aronovich Naimark, Aleksandr Isaakovich Shtern,
Theory of group representations;\\
Heidelberg, 1982
}%

\BiblioItem{Barut Raczka: Theory of group representations}
{
Asim Orhan Barut; Ryszard R\c{a}czka;
Theory of group representations and applications;\\
World Scientific Publishing Co. Pre. Ltd., 1986
}%

\BiblioItem{Mihalev Pilz: concise handbook of algebra}
{
Aleksandr Vasilevich Mikhalev; G\"{u}nter Pilz;
The concise handbook of algebra;\\
Kluwer Academic Publishers, 2002
}%

\BiblioItem{McCrimmon: Jordan Algebras}
{
Kevin McCrimmon;
A Taste of Jordan Algebras;\\
Springer, 2004
}%

\BiblioItem{Shafarevich: Basic notions of algebra}
{
I. R. Shafarevich,
Basic notions of algebra,\\
Translated from the Russian by M. Reid,\\
Springer, 2005
}%

\BiblioItem{Coppel: Number Theory}
{
W.A. Coppel,
Number Theory: An Introduction to Mathematics,\\
Springer, 2009
}%

\BiblioItem{978-0486497952}
{
Michael J. Field,
Differential Calculus and Its Applications,\\
Dover Publications, 2012; ISBN-13: 978-0486497952
}%

\BiblioItem{Elsgolts: Differential Equations}
{
Lev Elsgolts,
Differential Equations and the Calculus of Variations,\\
University Press of the Pacific, 2003 
}%

\BiblioItem{Baez Huerta: algebra of grand unified theories}
{
John Baez; John Huerta;
The algebra of grand unified theories;\\
Bull. Amer. Math. Soc. {\bf 47} (2010), 483-552
}%

\BiblioItem{J. Fan: Determinants}
{
J. Fan, Determinants and multiplicative functionals
on quaternion matrices, \textit{Linear Algebra and Its
Applications} {\bf 369}, pp. 193-201, (2003).
}%

\BiblioItem{Carl Faith 1}
{
Carl Faith, Algebra: Rings, Modules and Categories I,
Springer - Verlag, Berlin - Heidelberg - New York, 1973
}%

\BiblioItem{Gilson Nimmo Ohta}
{
 C.R.Gilson, J.J.C.Nimmo, Y.Ohta, Quasideterminant solutions of a non-Abelian Hirota-Miwa
 equation, \textit{Journal of Physics A: Mathematical and Theoretical} {\bf 40}(42), pp.
 12607-12617,(2007).
}%

\BiblioItem{Haider Hassan}
{
B. Haider, M. Hassan, Quasideterminant solutions of an integrable chiral model in two
 dimensions, \textit{Journal of Physics A: Mathematical and Theoretical} {\bf 42} (35), art. no.
 355211, (2009).
}%



\BiblioItem{0702447}
{
I.I. Kyrchei, Cramer's rule for quaternion systems of linear equations,
\textit{Journal of Mathematical Sciences} {\bf 155}(6), 839-858, (2008).
 Translated from  \textit{Fundamental and Appl. Math.}
 {\bf 13}(4), pp.67-94, (2007). (in Russian)\\
eprint
\href{http://arxiv.org/abs/math/0702447}{arXiv:math.RA/0702447}
(2007)
}%

\BiblioItem{1004.4380}
{
I.I. Kyrchei, Cramer's rule for some quaternion matrix
    equations,  \textit{Applied Mathematics and Computation} {\bf 217}(5), pp.2024-2030, (2010).\\eprint
\href{http://arxiv.org/abs/1004.4380
}{arXiv:math.RA/arXiv:1004.4380 } (2010)
}%

\BiblioItem{1005.0736}
{
I.I. Kyrchei,Determinantal representations of the Moore-Penrose inverse
 over the quaternion skew field and corresponding Cramer's rules,
 \\
eprint
\href{http://arxiv.org/abs/1005.0736}{arXiv:math.RA/1005.0736}
(2010)
}%

\BiblioItem{0412.391}
{
Aleks Kleyn,
Basis Manifold,
eprint \href{http://arxiv.org/abs/math.DG/0412391}{arXiv:math.DG/0412391} (2007)
}%

\BiblioItem{0405.027}
{
Aleks Kleyn,
Reference Frame in General Relativity,\\
eprint \href{http://arxiv.org/abs/gr-qc/0405027}{arXiv:gr-qc/0405027} (2008)
}%

\BiblioItem{0405.028}
{
Aleks Kleyn, Metric\hyph Affine Manifold,\\
eprint \href{http://arxiv.org/abs/gr-qc/0405028}{arXiv:gr-qc/0405028} (2008)
}%

\BiblioItem{0612.111}
{
Aleks Kleyn,
Biring of Matrices,\\
eprint \href{http://arxiv.org/abs/math.OA/0612111}{arXiv:math.OA/0612111} (2007)
}%

\BiblioItem{0701.238}
{
Aleks Kleyn,
Lectures on Linear Algebra over Division Ring,\\
eprint \href{http://arxiv.org/abs/math.GM/0701238}{arXiv:math.GM/0701238} (2010)
}%

\BiblioItem{0702.561}
{
Aleks Kleyn,
Fibered Universal Algebra,\\
eprint \href{http://arxiv.org/abs/math.DG/0702561}{arXiv:math.DG/0702561} (2007)
}%

\BiblioItem{math.RA-0501237}
{
Aleks Kleyn,
Vector Space Over Division Ring,\\
eprint \href{http://arxiv.org/abs/math.RA/0412391}{arXiv:math.RA/0501237} (2007)
}%

\BiblioItem{math.RA-0501237v1}
{
Aleks Kleyn,
Module Over Division Ring, version 1,\\
eprint \href{http://arxiv.org/abs/math/0501237v1}{arXiv:math.RA/0501237v1} (2005)
}%

\BiblioItem{0707.2246}
{
Aleks Kleyn,
Fibered Correspondence,\\
eprint \href{http://arxiv.org/abs/0707.2246}{arXiv:0707.2246} (2007)
}%

\BiblioItem{0803.2620}
{
Aleks Kleyn,
Morphism of \Ts{T}Representations,\\
eprint \href{http://arxiv.org/abs/0803.2620}{arXiv:0803.2620} (2008)
}%

\BiblioItem{0803.3276}
{
Aleks Kleyn,
Lorentz Transformation and General Covariance Principle,\\
eprint \href{http://arxiv.org/abs/0803.3276}{arXiv:0803.3276} (2009)
}%

\BiblioItem{0812.4763}
{
Aleks Kleyn,
Introduction into Calculus over Division Ring,\\
eprint \href{http://arxiv.org/abs/0812.4763}{arXiv:0812.4763} (2010)
}%

\BiblioItem{0906.0135}
{
Aleks Kleyn,
Introduction into Geometry over Division Ring,\\
eprint \href{http://arxiv.org/abs/0906.0135}{arXiv:0906.0135} (2010)
}%

\BiblioItem{0909.0855}
{
Aleks Kleyn,
Quaternion Rhapsody,\\
eprint \href{http://arxiv.org/abs/0909.0855}{arXiv:0909.0855} (2010)
}%

\BiblioItem{0912.3315}
{
Aleks Kleyn,
Representation of Universal Algebra,\\
eprint \href{http://arxiv.org/abs/0912.3315}{arXiv:0912.3315} (2009)
}%

\BiblioItem{0912.4061}
{
Aleks Kleyn,
Linear Equation in Finite Dimensional Algebra,\\
eprint \href{http://arxiv.org/abs/0912.4061}{arXiv:0912.4061} (2010)
}%

\BiblioItem{1001.4852}
{
Aleks Kleyn,
The Matrix of Linear Maps,\\
eprint \href{http://arxiv.org/abs/1001.4852}{arXiv:1001.4852} (2010)
}%

\BiblioItem{1003.1544}
{
Aleks Kleyn,
Linear Maps of Free Algebra,\\
eprint \href{http://arxiv.org/abs/1003.1544}{arXiv:1003.1544} (2010)
}%

\BiblioItem{1006.2597}
{
Aleks Kleyn,
The G\^ateaux Derivative and Integral over Banach Algebra,\\
eprint \href{http://arxiv.org/abs/1006.2597}{arXiv:1006.2597} (2010)
}%

\BiblioItem{1011.3102}
{
Aleks Kleyn,
Polylinear Map of Free Algebra,\\
eprint \href{http://arxiv.org/abs/1011.3102}{arXiv:1011.3102} (2010)
}%

\BiblioItem{1104.5197}
{
Aleks Kleyn,
$C^*$-Rhapsody,\\
eprint \href{http://arxiv.org/abs/1104.5197}{arXiv:1104.5197} (2011)
}%

\BiblioItem{1105.4307}
{
Aleks Kleyn,
Algebra with Conjugation,\\
eprint \href{http://arxiv.org/abs/1105.4307}{arXiv:1105.4307} (2011)
}%

\BiblioItem{1107.1139}
{
Aleks Kleyn,
Linear Maps of Quaternion Algebra,\\
eprint \href{http://arxiv.org/abs/1107.1139}{arXiv:1107.1139} (2011)
}%

\BiblioItem{1107.5037}
{
Aleks Kleyn,
Orthogonal Basis and Motion in Finsler Geometry,\\
eprint \href{http://arxiv.org/abs/1107.5037}{arXiv:1107.5037} (2011)
}%

\BiblioItem{1111.6035}
{
Aleks Kleyn,
Basis of Representation of Universal Algebra,\\
eprint \href{http://arxiv.org/abs/1111.6035}{arXiv:1111.6035} (2011)
}%

\BiblioItem{1201.4158}
{
Aleks Kleyn, Alexandre Laugier,
Orthonormal Basis in Minkowski Space,\\
eprint \href{http://arxiv.org/abs/1201.4158}{arXiv:1201.4158} (2012)
}%

\BiblioItem{1202.6021}
{
Aleks Kleyn,
Maps of Conjugation of Quaternion Algebra,\\
eprint \href{http://arxiv.org/abs/1202.6021}{arXiv:1202.6021} (2012)
}%

\BiblioItem{1206.0200}
{
Aleks Kleyn,
Algebra of Fractions of Algebra with Conjugation,\\
eprint \href{http://arxiv.org/abs/1206.0200}{arXiv:1206.0200} (2012)
}%

\BiblioItem{1211.6965}
{
Aleks Kleyn,
Free Algebra with Countable Basis,\\
eprint \href{http://arxiv.org/abs/1211.6965}{arXiv:1211.6965} (2012)
}%

\BiblioItem{1302.7204}
{
Aleks Kleyn,
Polynomial over Associative $D$-Algebra,\\
eprint \href{http://arxiv.org/abs/1302.7204}{arXiv:1302.7204} (2013)
}%

\BiblioItem{1305.4547}
{
Aleks Kleyn,
Normed $\Omega$-Group,\\
eprint \href{http://arxiv.org/abs/1305.4547}{arXiv:1305.4547} (2013)
}%

\BiblioItem{1310.5591}
{
Aleks Kleyn,
Integral of Map into Abelian $\Omega$\Hyph group,\\
eprint \href{http://arxiv.org/abs/1310.5591}{arXiv:1310.5591} (2013)
}%

\BiblioItem{1502.04063}
{
Aleks Kleyn,
Linear Map of $D$\Hyph Algebra,\\
eprint \href{http://arxiv.org/abs/1502.04063}{arXiv:1502.04063} (2015)
}%

\BiblioItem{1505.03625}
{
Aleks Kleyn,
Derivative of Map of Banach algebra,\\
eprint \href{http://arxiv.org/abs/1505.03625}{arXiv:1505.03625} (2015)
}%

\BiblioItem{1601.03259}
{
Aleks Kleyn,
Introduction into Calculus over Banach Algebra,\\
eprint \href{http://arxiv.org/abs/1601.03259}{arXiv:1601.03259} (2016)
}%

\BiblioItem{MRepro}
{
Aleks Kleyn,
Representation Theory of Universal Algebra,\\
eprint \href{http://arxiv.org/abs/MRepro}{arXiv:MRepro} (2015)
}%

\BiblioItem{MVector}
{
Aleks Kleyn,
Linear Map of D-Algebra,\\
eprint \href{http://arxiv.org/abs/MVector}{arXiv:MVector} (2017)
}%

\BiblioItem{8433-5163}
{
Aleks Kleyn,
Linear Maps of Free Algebra: First Steps in Noncommutative Linear Algebra,\\
Lambert Academic Publishing, 2010
}%

\BiblioItem{8443-0072}
{
Aleks Kleyn,
Representation Theory: Representation of Universal Algebra,\\
Lambert Academic Publishing, 2011
}%

\BiblioItem{4776-3181}
{
Aleks Kleyn.\\
Linear Algebra over Division Ring: System of Linear Equations.\\
CreateSpace Independent Publishing Platform, 2012;\\
ISBN-13: 978-1477631812
}%

\BiblioItem{4993-2400}
{
Aleks Kleyn.\\
Linear Algebra over Division Ring: Vector Space.\\
CreateSpace Independent Publishing Platform, 2014;\\
ISBN-13: 978-1499324006
}%

\BiblioItem{5059-9176}
{
Aleks Kleyn.\\
Normed \(\Omega\)-Group.\\
CreateSpace Independent Publishing Platform, 2015;\\
ISBN-13: 978-1505991765
}%

\BiblioItem{5114-6019}
{
Aleks Kleyn.\\
Representation of Universal Algebra: Polymorphism.\\
CreateSpace Independent Publishing Platform, 2015;\\
ISBN-13: 978-1511460194
}%

\BiblioItem{5148-4632}
{
Aleks Kleyn.\\
Representation of Universal Algebra: Polymorphism.\\
CreateSpace Independent Publishing Platform, 2015;\\
ISBN-13: 978-1511460194
}%

\BiblioItem{5410-9916}
{
Aleks Kleyn.\\
Lebesgue Integral in Abelian $\Omega$-Group.\\
CreateSpace Independent Publishing Platform, 2016;\\
ISBN-13: 978-1541099166
}%

\BiblioItem{BRepro}
{
Aleks Kleyn,
Representation Theory of Universal Algebra,\\
CreateSpace Independent Publishing Platform, 2015;\\
ISBN-13: 
}%

\BiblioItem{CACAA.01.291}
{
Aleks Kleyn,
Introduction into Calculus over Division Ring.\\
Clifford Analysis, Clifford Algebras and their applications,
volume 1, Issue 4, pages 291 - 355, 2012
}%

\BiblioItem{CACAA.02.097}
{
Aleks Kleyn,
Polynomial over Associative $D$-Algebra.\\
Clifford Analysis, Clifford Algebras and their applications,
volume 2, Issue 2, pages 97 - 115, 2013
}%

\BiblioItem{CACAA.04.001}
{
Aleks Kleyn,
Integral of Map into Abelian $\Omega$-group.\\
Clifford Analysis, Clifford Algebras and their applications,
volume 4, Issue 1, pages 1 - 68, 2013
}%

\BiblioItem{CACAA.05.001}
{
Aleks Kleyn,
Introduction into Calculus over Division Ring.\\
Clifford Analysis, Clifford Algebras and their applications,
volume 5, issue 1, pages 1 - 68, 2016 
}%

\BiblioItem{GJSFRA.13.1.39}
{
Aleks Kleyn,
Reference frame and Lorentz transformation,\\
Global Journals of Science Frontier Research A,
volume 13, issue 1, pages 39 - 55, 2013 
}%

\BiblioItem{1506.05848}
{
Rida T. Farouki, Graziano Gentili, Carlotta Giannelli, Alessandra Sestini,
Caterina Stoppato,\\
Solution of a quadratic quaternion equation with mixed coefficients,\\
eprint \href{http://arxiv.org/abs/1506.05848}{arXiv:1506.05848} (2015)
}%

\BiblioItem{Lauve: Quantum coordinates}
{
A. Lauve, Quantum- and quasi-Plucker coordinates,
\textit{Journal of Algebra} {\bf 296}(2), pp.440-461,
(2006).
}%

\BiblioItem{Lewis D. W. Quaternion algebras}
{
Lewis D. W. Quaternion algebras and the algebraic legacy
of Hamilton's quaternions, \textit{Irish Math. Soc. Bulletin} {\bf
57}, pp. 41-64, (2006).
}%

\BiblioItem{0812.2865}
{
Jos\'e Miguel Figueroa-O'Farrill,
Three lectures on 3-algebras,
eprint \href{http://arxiv.org/abs/0812.2865}{arXiv:0812.2865} (2008)
}%

\BiblioItem{1202.0951}
{
Daniel Edward Clark,
Deconvolution of point processes,
eprint \href{http://arxiv.org/abs/1202.0951}{arXiv:1202.0951} (2012)
}%

\BiblioItem{1202.4546}
{
Ming-Liang Hu,
Disentanglement, Bell-nonlocality violation
and teleportation capacity of the decaying tripartite states,
eprint \href{http://arxiv.org/abs/1202.4546}{arXiv:1202.4546} (2012)
}%

\BiblioItem{1203.1629}
{
Borivoje Dakic, Yannick Ole Lipp, Xiaosong Ma, Martin Ringbauer,
Sebastian Kropatschek, Stefanie Barz, Tomasz Paterek, Vlatko Vedral,
Anton Zeilinger, Caslav Brukner, Philip Walther,
Quantum Discord as Optimal Resource for Quantum Communication,
eprint \href{http://arxiv.org/abs/1203.1629}{arXiv:1203.1629} (2012)
}%

\BiblioItem{Li Nimmo: Darboux transformations}
{
C.X.Li, J.J.C. Nimmo, Darboux transformations for a twisted
derivation and quasideterminant solutions to the super KdV
equation, \textit{Proceedings of the Royal Society A:
Mathematical, Physical and Engineering Sciences} {\bf 466} (2120),
pp. 2471-2493, (2010).
}%

\BiblioItem{Schiebold: Cauchy-type determinants}
{
C. Schiebold, Cauchy-type determinants and integrable
systems, \textit{Linear Algebra and Its Applications} {\bf 433}
(2), pp. 447-475, (2010)
}%

\BiblioItem{Suzuki: Noncommutative spectral decomposition}
{
T. Suzuki, Noncommutative
spectral decomposition with qua\-si\-de\-ter\-mi\-nant, \textit{Advances in
Mathematics} {\bf 217}(5), pp. 2141-2158, (2008).
}%

\BiblioItem{1105.3456}
{
C. W. F. Everitt, D. B. DeBra, B. W. Parkinson, J. P. Turneaure, J. W. Conklin,
M. I. Heifetz, G. M. Keiser, A. S. Silbergleit, T. Holmes, J. Kolodziejczak,
M. Al-Meshari, J. C. Mester, B. Muhlfelder, V. Solomonik, K. Stahl, P. Worden,
W. Bencze, S. Buchman, B. Clarke, A. Al-Jadaan, H. Al-Jibreen, J. Li, J. A. Lipa,
J. M. Lockhart, B. Al-Suwaidan, M. Taber, S. Wang,\\
Gravity Probe B: Final Results of a Space Experiment to Test General Relativity,\\
eprint \href{http://arxiv.org/abs/1105.3456}{arXiv:1105.3456[gr-qc]} (2011)
}%

\BiblioItem{0009305}
{
G. S. Asanov.
Can Neutrinos and High-Energy Particles Test Finsler Metric of Space-Time?\\
eprint \href{http://arxiv.org/abs/hep-ph/0009305}{arXiv:hep-ph/0009305} (2000)
}%

\BiblioItem{Asanov 2004}
{
G. S. Asanov.
Finsleroid - space supplemented by angle and scalar product.\\
Hypercomplex Numbers in Geometry and Physics, {\bf 1}, 2004, p. 40 - 62
}%

\BiblioItem{1004.3007}
{
Sergiu I. Vacaru,
Principles of Einstein-Finsler Gravity and Perspectives in Modern Cosmology,\\
eprint \href{http://arxiv.org/abs/1004.3007}{arXiv:1004.3007[math-ph]} (2010)
}%

\BiblioItem{1012.4148}
{
Sergiu I. Vacaru.
Principles of Einstein-Finsler Gravity and Cosmology.\\
eprint \href{http://arxiv.org/abs/1012.4148}{arXiv:1012.4148[physics.gen-ph]} (2010)
}%

\BiblioItem{1112.5641}
{
Christian Pfeifer, Mattias N.R. Wohlfarth.
Finsler geometric extension of Einstein gravity.\\
eprint \href{http://arxiv.org/abs/1112.5641}{arXiv:1112.5641[gr-qc]} (2011)
}%

\BiblioItem{0711.0056}
{
Zhe Chang, Xin Li.
Lorentz Invariance Violation and Symmetry in Randers\Hyph Finsler Spaces.\\
eprint \href{http://arxiv.org/abs/0711.0056}{arXiv:0711.0056[hep-th]} (2011)
}%

\BiblioItem{Zharinov Kursy NOC, 9}
{
В. В. Жаринов,
Алгебро-геометрические основы математической физики.
\\
Лекц. курсы НОЦ, 9, МИАН, М., 2008, 3–209
}%

\BiblioItem{Rund Finsler geometry}
{
Hanno Rund,
The differential geometry of Finsler spaces.
\\
Springer - Verlag, Berlin - G\"ottingen - Heidelberg, 1959
}%

\BiblioItem{Smirnov vol 1}
{
V. I. Smirnov,
A Course of Higher Mathematics, volume I.
\\
Translated by D. E. Brown.
\\
Translation, edited and additions made by I. N. Sneddon.
\\
Pergamon Press, Addison-Wesley Publishing Company, 1964
}%

\BiblioItem{Beem Dostoglou Ehrlich}
{
John K. Beem, Stamatis A. Dostoglou, Paul E. Ehrlich,
Advances in differential geometry and general relativity.
\\
American Mathematical Society, 2004
}%

\BiblioItem{978-0719033414}
{
Malcolm Pemberton, Nicholas Rau,
Mathematics for economists: an introductory textbook.
\\
Manchester University Press, November 2001; ISBN-13: 978-0719033414
}%

\BiblioItem{0 521 59180 5}
{
Cyrus D. Cantrell,
Modern mathematical methods for physicists and engineers.
\\
Cambridge University Press, 2000
}%

\BiblioItem{Arveson spectral theory}
{
William Arveson,
A short course on spectral theory.
\\
Springer - Verlag, New York, 2002
}%

\BiblioItem{Robert Hermann}
{
Robert Hermann,
Topics in the mathematics of quantum mechanics.
\\
Math Sci Press, 1973
}%

\BiblioItem{9705.009}
{
John C. Baez,
An Introduction to n-Categories,\\
eprint \href{http://arxiv.org/abs/q-alg/9705009}{arXiv:q-alg/9705009} (1997)
}%

\BiblioItem{0105.155}
{
John C. Baez,
The Octonions,\\
eprint \href{http://arxiv.org/abs/math.RA/0105155}{arXiv:math.RA/0105155} (2002)
}%

\BiblioItem{John Baez: Math Blogs}
{
John C. Baez,
What do mathematicians need to know about blogging?,\\
Notices of the American Mathematical Society,
(2010), 3, {\bf 57}, 333,\\
\url{http://www.ams.org/notices/201003/rtx100300333p.pdf}
}%

\BiblioItem{Tolstoi about Anna Karenina}
{
Tolstoi about Anna Karenina,
in book A Karenina Companion, by C. J. G. Turner,
published by Wilfrid Laurier University Press (August 1993)
}%

\BiblioItem
{Cohn: Universal Algebra}
{
Paul M. Cohn,
Universal Algebra,
Springer, 1981
}%

\BiblioItem
{Cohn: Algebra 1}
{
Paul M. Cohn,
Algebra, Volume 1,
John Wiley \& Sons, 1982
}%

\BiblioItem
{Cohn: Algebra 3}
{
Paul M. Cohn,
Algebra, Volume 3,
John Wiley \& Sons, 1991
}%

\BiblioItem
{Cohn: Skew Fields}
{
Paul M. Cohn,
Skew Fields,
Cambridge University Press, 1995
}%

\BiblioItem
{Lam: Noncommutative Rings}
{
T. Y. Lam,
A First Course in
Noncommutative Rings,
Springer-Verlag, 1991
}%

\BiblioItem
{Maunder: Algebraic Topology}
{
C. R. F. Maunder,
Algebraic Topology,
Dover Publications, Inc, Mineola, New York, 1996
}%

\BiblioItem{Pommaret: Partial Differential Equations}
{
J.-F. Pommaret,
Partial Differential Equations and Group Theory,
Springer, 1994
}%

\BiblioItem{Bourbaki: Set Theory}
{
N. Bourbaki,
Theory of sets,
Springer, 2004
}%

\BiblioItem{Bourbaki: Algebra 1}
{
N. Bourbaki,
Algebra 1,
Springer, 2004
}%

\BiblioItem{Bourbaki: Algebra 2}
{
N. Bourbaki,
Algebra II, Chapters 4 - 7,//
Translated by P. M. Cohn & J. Howie,//
Springer, 2004
}%

\BiblioItem
{Bourbaki: General Topology 1}
{
N. Bourbaki,
General Topology, Chapters 1 - 4,
Springer, 1989
}

\BiblioItem{Bourbaki: General Topology: Chapter 5 - 10}
{
N. Bourbaki,
General Topology, Chapters 5 - 10,
Springer, 1989
}

\BiblioItem{Bourbaki: Topological Vector Space}
{
N. Bourbaki,
Topological Vector Spaces, Chapters 1 - 5,
Transl. by H. G. Eggleston $\&$ S. Madan,
Springer, 2003
}

\BiblioItem{Bourbaki: Coxeter Group Lie}
{
N. Bourbaki,
Lie Groups and Lie Algebras, Chapters 4 - 6,
Translator Andrew Pressley,
Springer, 2002
}

\BiblioItem{Bourbaki: Real Group Lie}
{
N. Bourbaki,
Lie Groups and Lie Algebras, Chapters 7 - 9,
Translator Andrew Pressley,
Springer, 2005
}

\BiblioItem{Shabat: Complex Analysis}
{
Shabat B. V.,
Introduction to Complex Analysis,
Moscow, Nauka, 1969
}

\BiblioItem{Pontryagin: Topological Group}
{
L. S. Pontryagin,
Selected Works, Volume Two, Topological Groups,
Gordon and Breach Science Publishers, 1986
}

\BiblioItem
{Eisenhart: Riemannian Geometry}
{
Eisenhart,
Riemannian Geometry,
Princeton University Press, Princeton, 1949
}

\BiblioItem
{Eisenhart: Continuous Groups of Transformations}
{
Eisenhart,
Continuous Groups of Transformations,
Dover Publications, New York, 1961
}

\BiblioItem
{Condon Odabasi}
{
Edward Uhler Condon, Halis Odabasi,
Atomic Structure,
CUP Archive, 1980
}

\BiblioItem{Postnikov: Differential Geometry}
{
Postnikov M. M.,
Geometry IV: Differential geometry,
Moscow, Nauka, 1983
}

\BiblioItem{Fikhtengolts: Calculus volume 1}
{
Fikhtengolts G. M.,
Differential and Integral Calculus Course, volume 1,
Moscow, Nauka, 1969
}

\BiblioItem{Fikhtengolts: Calculus volume 2}
{
Fikhtengolts G. M.,
Differential and Integral Calculus Course, volume 2,
Moscow, Nauka, 1969
}

\BiblioItem{Fikhtengolts: Calculus volume 3}
{
Fikhtengolts G. M.,
Differential and Integral Calculus Course, volume 3,
Moscow, Nauka, 1969
}

\BiblioItem{Hatcher: Algebraic Topology}
{
Allen Hatcher,
Algebraic Topology,
Cambridge University Press, 2002
}

\BiblioItem{geometry of differential equations}
{
Krasil'shchik I. S., Lychagin V. V., Vinogradov A. M.,
Geometry of Jet Spaces and Nonlinear Partial Differential Equations,
\\
Translated from the Russian by A. B. Sosinskii,
\\
Gordon and Breach Science Publishers, 1985
}

\BiblioItem{Basic Concepts of Differential Geometry}
{
Alekseyevskii D. V., Vinogradov A. M., Lychagin V. V.,
Basic Concepts of Differential Geometry
\\
VINITI Summary 28
\\
Moscow. VINITI, 1988
}

\BiblioItem{cohomological analysis}
{
A. M. Vinogradov,
Cohomological Analysis of Partial Differential Equations
and Secondary Calculus,
American Mathematical Society, 2001
}

\BiblioItem{0801.1734}
{
Brandon S. DiNunno, Richard A. Matzner,
The Volume Inside a Black Hole,\\
eprint \href{http://arxiv.org/abs/0801.1734v1}{arXiv:0801.1734v1} (2008)
}

\BiblioItem{0702.447}
{
Ivan Kyrchei,
Cramer's rule for some quaternion matrix equations,\\
eprint \href{http://arxiv.org/abs/math/0702447}{arXiv:math.RA/0702447} (2007)
}

\BiblioItem{Izrail M. Gelfand: Quaternion Groups}
{
I. M. Gelfand, M. I. Graev,
Representation of Quaternion Groups over Localy Compact and
Functional Fields,\\
Funct. Anal. Appl. {\bf 2} (1968) 19 - 33;\\
Izrail Moiseevich Gelfand, Semen Grigorevich Gindikin,\\
Izrail M. Gelfand: Collected Papers, volume II, 435 - 449,\\
Springer, 1989
}

\BiblioItem{Richard D. Schafer}
{
Richard D. Schafer,
An Introduction to Nonassociative Algebras,
Dover Publications, Inc., New York, 1995
}

\BiblioItem{Bamberg Sternberg}
{
Paul Bamberg, Shlomo Sternberg,
A course in mathematics for students of physics,
Cambridge University Press, 1991
}

\BiblioItem{Conway Smith}
{
John Horton Conway, Derek Alan Smith,
On quaternions and octonions: their geometry, arithmetic, and symmetry,
A K Peters, Natick, Massachussets, 2003
}

\BiblioItem{Fueter}
{
Fueter, R.
Die Funktionentheorie der Differentialgleichungen $\Delta u = 0$ und
$\Delta \Delta u = 0$ mit vier reellen Variablen.
Comment. Math. Helv. {\bf 7} (1935), 307-330
}

\BiblioItem{Sudbery Quaternionic Analysis}
{
A. Sudbery,
Quaternionic Analysis,
Math. Proc. Camb. Phil. Soc. (1979), {\bf 85}, 199 - 225
}

\BiblioItem{0902.4771}
{
Fabrizio Colombo, Graziano Gentili, Irene Sabadini,
A Cauchy kernel for slice regular functions,\\
eprint \href{http://arxiv.org/abs/0902.4771v1}{arXiv:0902.4771v1} (2009)
}

\BiblioItem{Vadim Komkov}
{
Vadim Komkov,
Variational Principles of Continuum Mechanics with Engineering Applications: Critical Points Theory,\\
Springer, 1986
}

\BiblioItem{Alain Connes 1994}
{
Alain Connes,
Noncommutative Geometry,\\
Academic Press, 1994
}

\BiblioItem{Hamilton papers 3}
{
Sir William Rowan Hamilton,
The Mathematical Papers, Vol. III, Algebra,\\
Cambridge at the University Press, 1967
}

\BiblioItem{Hamilton Elements of Quaternions 1}
{
Sir William Rowan Hamilton,
Elements of Quaternions, Volume I,\\
Longmans, Green, and Co., London, New York, and Bombay, 1899
}

\BiblioItem{Cartan geometry in reper}
{
Elie Cartan, Vladislav V. Goldberg, Serge\u{i} Pavlovich Finikov,\\
Riemannian geometry in an orthogonal frame:
from lectures delivered by Elie Cartan at the Sorbonne in 1926-1927,\\
translated by Vladislav V. Goldberg,\\
World Scientific, 2001
}

\BiblioItem{Cartan differential form}
{
Henri Cartan.
Differential forms.\\
Kershaw Publishing Company Limited, London, 1971
}

\BiblioItem{Arnautov Glavatsky Mikhalev}
{
V. I. Arnautov, S. T. Glavatsky, A. V. Mikhalev,\\
Introduction to the theory of topological rings and modules,
Volume 1995,\\
Marcel Dekker, Inc, 1996
}

\BiblioItem{Moore Yaqub}
{
Hal G. Moore, Adil Yaqub,
A first course in linear algebra with applications,
Edition 3, Academic Press, 1998 
}

\BiblioItem{math.CV-0405471}
{
S. V. Ludkovsky,
Differentiable functions of Cayley-Dickson numbers,\\
eprint \href{http://arxiv.org/abs/math.CV/0405471}{arXiv:math.CV/0405471} (2004)
}%

\BiblioItem{W.Bertram H.Glockner K.Neeb}
{
W.Bertram, H.Glockner, K.Neeb,
Differential Calculus over General Base Fields and Rings,
Expositiones Mathematicae (2004), Volume 22, Issue 3, Pages 213-282
}

\CloseBiblio

%% file: Index.English.tex
\OpenIndex
\SetIndexSpace%
\Index
   {$*A$\Hyph vector space}%
   {*A-vector space}%
\Index
   {\CR product}%
   {cr-product}%
\Index
   {\RC product}%
   {rc-product}%
\SetIndexSpace%
\Index
   {$1$\Hyph form}%
   {1-form}%
\SetIndexSpace%
\Index
   {$2$\Hyph ary fibered relation}%
   {2 ary fibered relation}%
\SetIndexSpace%
\Index
   {$A$\Hyph algebra of polynomials over $D$\Hyph algebra $A$}%
   {algebra of polynomials over algebra}%
\Index
   {$A$\Hyph number}%
   {A number}%
\Index
   {$\mathcal A(A)$\Hyph map}%
   {A(A) map}%
\Index
   {$A*$\Hyph module}%
   {A*-module}%
\Index
   {$A*$\Hyph vector space}%
   {A*-vector space}%
\Index
   {$A$\Hyph module}%
   {module over algebra}%
\Index
   {$A$\Hyph valued function}%
   {A valued function}%
\Index
   {$A$\Hyph representation in $\Omega$\Hyph algebra}%
   {A representation of algebra}%
\Index
   {Abelian multiplicative $\Omega$\Hyph group}%
   {Abelian multiplicative Omega group}%
\Index
   {Abelian $\Omega$\Hyph group}%
   {Abelian Omega group}%
\Index
   {Abelian semigroup}%
   {Abelian semigroup}%
\Index
   {absolute value}%
   {absolute value}%
\Index
   {active \sT{G}representation}%
   {active representation, vector space}%
\Index
   {active representation}%
   {active representation}%
\Index
   {active representation of group $G(f)$ in basis manifold of representation}%
   {active representation in basis manifold}%
\Index
   {active representation of group $G(\Vector f)$ in basis manifold of tower of representations}%
   {active representation in basis manifold, tower of representations}%
\Index
   {active transformation of basis manifold of representation}%
   {active transformation of basis, representation}%
\Index
   {active transformation of basis manifold of tower of representations}%
   {active transformation of basis, tower of representations}%
\Index
   {active transformation on basis manifold}%
   {active transformation}%
\Index
   {active transformation on the set of \rcd bases}%
   {active transformation, vector space}%
\Index
   {additive map}%
   {additive map}%
\Index
   {affine basis}%
   {Affine Basis}%
\Index
   {affine functional}%
   {affine functional}%
\Index
   {affine representation of Lie group}%
   {affine representation of Lie group}%
\Index
   {affine space}%
   {affine space}%
\Index
   {affine structure on set}%
   {affine structure on set}%
\Index
   {affine transformation}%
   {affine transformation}%
\Index
   {affine transformation group}%
   {affine transformation group}%
\Index
   {affine transformation group}%
   {affine transformation group}%
\Index
   {affine transformation on basis manifold}%
   {affine transformation}%
\Index
   {algebra of fractions of algebra with conjugation}%
   {algebra of fractions of algebra with conjugation}%
\Index
   {algebra of polynomials over $D$\Hyph algebra}%
   {algebra of polynomials over D algebra}%
\Index
   {algebra of rational mappings of algebra}%
   {algebra of rational mappings of algebra}%
\Index
   {algebra of sets}%
   {algebra of sets}%
\Index
   {algebra over ring}%
   {algebra over ring}%
\Index
   {algebra with conjugation}%
   {algebra with conjugation}%
\Index
   {alternation of polylinear map}%
   {alternation of polylinear map}%
\Index
   {alternative representation of matrix}%
   {Alternative representation}%
\Index
   {anholonomic coordinate}%
   {anholonomic coordinate}%
\Index
   {anholonomic coordinates of connection}%
   {anholonomic coordinates of connection}%
\Index
   {anholonomic coordinates of vector}%
   {vector anholonomic coordinates}%
\Index
   {anholonomic coordinates on manifold}%
   {anholonomic coordinates on manifold}%
\Index
   {anholonomity object}%
   {anholonomity object}%
\Index
   {antilinear map}%
   {antilinear map}%
\Index
   {antisymmetric $2$\Hyph ary fibered relation}%
   {antisymmetric 2 ary fibered relation}%
\Index
   {$A\RCstar$\Hyph basis for vector space}%
   {Arc basis, vector space}%
\Index
   {arity}%
   {arity}%
\Index
   {arity of operation}%
   {arity of operation}%
\Index
   {associative $D$\Hyph algebra}%
   {associative D algebra}%
\Index
   {associative law}%
   {associative law}%
\Index
   {associative $\Omega$\Hyph group}%
   {associative Omega group}%
\Index
   {associative operation}%
   {associative operation}%
\Index
   {associator of $D$\Hyph algebra}%
   {associator of algebra}%
\Index
   {auto parallel line}%
   {auto parallel line}%
\Index
   {automorphism}%
   {automorphism}%
\Index
   {automorphism of representation of $\Omega$\Hyph algebra}%
   {automorphism of representation}%
\Index
   {automorphism of tower of representations}%
   {automorphism of tower of representations}%
\Index
   {automorphism of vector space}%
   {automorphism of vector space}%
\Index
   {$(^j_i)$\hyph \CR quasideterminant}%
   {j i cr-quasideterminant}%
\Index
   {norm of quaternion}%
   {norm of quaternion}%
\SetIndexSpace%
\Index
   {$B$\Hyph set}%
   {B set}%
\Index
   {Banach $D$\Hyph algebra}%
   {Banach algebra}%
\Index
   {Banach $D$\Hyph module}%
   {Banach module}%
\Index
   {base of fibered correspondence}%
   {base of fibered correspondence}%
\Index
   {base of mapping}%
   {base of map}%
\Index
   {basis}%
   {Basis}%
\Index
   {basis dual to basis}%
   {basis dual to basis}%
\Index
   {basis dual to basis}%
   {dual basis}%
\Index
   {basis for \crd vector space}%
   {basis, crd vector space}%
\Index
   {basis for \dcr vector space}%
   {basis, dcr vector space}%
\Index
   {basis for \drc vector space}%
   {basis, drc vector space}%
\Index
   {basis for module}%
   {basis, module}%
\Index
   {basis for \rcd vector space}%
   {basis, rcd vector space}%
\Index
   {basis for vector space}%
   {basis, vector space}%
\Index
   {basis manifold of affine space}%
   {Basis Manifold, Affine Space}%
\Index
   {basis manifold of central affine space}%
   {Basis Manifold, Central Affine Space}%
\Index
   {basis manifold of Euclid space}%
   {Basis Manifold, Euclid Space}%
\Index
   {basis manifold of Euclid space}%
   {Basis Manifold, Euclid Space, division ring}%
\Index
   {basis manifold of \rcd vector space}%
   {basis manifold of rcd vector space}%
\Index
   {basis manifold of representation}%
   {basis manifold representation F algebra}%
\Index
   {basis manifold of tower of representations}%
   {basis manifold tower of representations}%
\Index
   {basis manifold of vector space}%
   {basis manifold of vector space}%
\Index
   {basis of Abelian group}%
   {basis of Abelian group}%
\Index
   {basis of algebra $\mathcal L(A;A)$}%
   {basis of algebra L(A,A)}%
\Index
   {basis of representation}%
   {basis of representation}%
\Index
   {basis of tower of representations}%
   {basis of tower of representations}%
\Index
   {basis vector of representation of Lie group over algebra $A$}%
   {basis vector of representation of Lie group over algebra A}%
\Index
   {biring}%
   {biring}%
\Index
   {Borel algebra}%
   {Borel algebra}%
\Index
   {Borel set}%
   {Borel set}%
\Index
   {Borel\Hyph measurable map}%
   {Borel-measurable map}%
\Index
   {bundle of level $2$}%
   {bundle of level 2}%
\Index
   {bundle of level $n$}%
   {bundle of level n}%
\SetIndexSpace%
\Index
   {\subs row of matrix}%
   {c row}%
\Index
   {$c$\hyph row of matrix}%
   {c-row}%
\Index
   {can be embeded}%
   {can be embeded}%
\Index
   {canonical remainder of the division}%
   {canonical remainder of the division}%
\Index
   {canonical representation of division with remainder}%
   {canonical representation of division with remainder}%
\Index
   {carrier of $\Omega$\Hyph algebra}%
   {carrier of Omega-algebra}%
\Index
   {Cartan connection}%
   {Cartan connection}%
\Index
   {Cartan curvature}%
   {Cartan curvature}%
\Index
   {Cartan derivative}%
   {Cartan derivative}%
\Index
   {Cartan symbol}%
   {Cartan symbol}%
\Index
   {Cartan transport}%
   {Cartan transport}%
\Index
   {Cartesian power}%
   {Cartesian power}%
\Index
   {Cartesian power $\Bundle A$ of bundle $\Bundle B$}%
   {Cartesian power A of bundle B}%
\Index
   {Cartesian power $A$ of set $B$}%
   {Cartesian power of set}%
\Index
   {Cartesian power $n$ of bundle $\Bundle E$}%
   {Cartesian power n of bundle E}%
\Index
   {Cartesian power $n$ of $\mathfrak{H}$\Hyph algebra}%
   {Cartesian power of algebra}%
\Index
   {Cartesian power of systems of subsets}%
   {Cartesian power of systems of subsets}%
\Index
   {Cartesian product of groups}%
   {Cartesian product of groups}%
\Index
   {Cartesian product of measures}%
   {Cartesian product of measures}%
\Index
   {Cartesian product of \(\Omega\)\Hyph algebras}%
   {Cartesian product of Omega algebras}%
\Index
   {Cartesian product of systems of subsets}%
   {Cartesian product of systems of subsets}%
\Index
   {category of \drc vector spaces}%
   {category of drc vector spaces}%
\Index
   {category of fibered correspondences over diagonal}%
   {category of fibered correspondences over diagonal}%
\Index
   {category of left-side representations}%
   {category of left-side representations}%
\Index
   {category of left-side representations of $\Omega_1$\Hyph algebra $A$}%
   {category of left-side representations of Omega1 algebra}%
\Index
   {category of reduced fibered correspondences}%
   {category of reduced fibered correspondences}%
\Index
   {category of representations}%
   {category of representations}%
\Index
   {Cauchy sequence}%
   {Cauchy sequence}%
\Index
   {center of $D$\Hyph algebra $A$}%
   {center of algebra}%
\Index
   {center of ring $D$}%
   {center of ring}%
\Index
   {central affine basis}%
   {Central Affine Basis}%
\Index
   {closed ball}%
   {closed ball}%
\Index
   {closure of set}%
   {closure of set}%
\Index
   {coefficient of polynomial}%
   {coefficient of polynomial}%
\Index
   {column $D*$\Hyph vector}%
   {column D* vector}%
\Index
   {column determinant}%
   {column determinant}%
\Index
   {column vector}%
   {column vector}%
\Index
   {common factor}%
   {common factor}%
\Index
   {commutative $D$\Hyph algebra}%
   {commutative D algebra}%
\Index
   {commutative diagram of correspondences}%
   {commutative diagram of correspondences}%
\Index
   {commutative operation}%
   {commutative operation}%
\Index
   {commutativity of representations}%
   {commutativity of representations}%
\Index
   {commutator of $D$\Hyph algebra}%
   {commutator of algebra}%
\Index
   {compact set}%
   {compact set}%
\Index
   {compact\hyph open topology}%
   {compact open topology}%
\Index
   {complete division ring}%
   {complete division ring}%
\Index
   {complete measure}%
   {complete measure}%
\Index
   {complete normed $\Omega$\Hyph group}%
   {complete Omega group}%
\Index
   {complete ring}%
   {complete ring}%
\Index
   {complete system of linear partial differential equations}%
   {Complete System of Linear Partial Differential Equations}%
\Index
   {completely integrable system}%
   {completely integrable system}%
\Index
   {completion of normed $\Omega$\Hyph group}%
   {completion of normed Omega group}%
\Index
   {completion of representation}%
   {completion of representation}%
\Index
   {component of derivative}%
   {component of derivative}%
\Index
   {component of derivative of Second Order}%
   {component of derivative of Second Order}%
\Index
   {component of linear map}%
   {component of linear map}%
\Index
   {component of polylinear map}%
   {component of polylinear map}%
\Index
   {component of the G\^ateaux derivative}%
   {component of Gateaux derivative}%
\Index
   {component of the G\^ateaux derivative of second order}%
   {component of Gateaux derivative of Second Order}%
\Index
   {composition of fibered correspondences}%
   {composition of fibered correspondences}%
\Index
   {composition of reduced fibered correspondences}%
   {composition of reduced fibered correspondences}%
\Index
   {condition of reducibility of products}%
   {condition of reducibility of products}%
\Index
   {congruence}%
   {congruence}%
\Index
   {conjugate of quaternion $x$}%
   {conjugate of quaternion}%
\Index
   {conjugated affine space}%
   {conjugated affine space}%
\Index
   {conjugated $D$\Hyph  module}%
   {conjugated D module}%
\Index
   {conjugated vector space}%
   {conjugated vector space}%
\Index
   {conjugation in algebra}%
   {conjugation in algebra}%
\Index
   {conjugation in ring}%
   {conjugation in ring}%
\Index
   {conjugation transformation}%
   {conjugation transformation}%
\Index
   {connected set}%
   {connected set}%
\Index
   {connection coefficients in affine space}%
   {connection coefficients, affine space}%
\Index
   {connection in principal fibre bundle}%
   {connection in principal bundle}%
\Index
   {contact point of set}%
   {contact point of set}%
\Index
   {continues basis}%
   {continues basis}%
\Index
   {continuous correspondence}%
   {continuous correspondence}%
\Index
   {continuous map}%
   {continuous map}%
\Index
   {continuous multivariable map}%
   {continuous multivariable map}%
\Index
   {convex set}%
   {convex set}%
\Index
   {coordinate isomorphism}%
   {coordinate isomorphism}%
\Index
   {coordinate matrix of set of vectors}%
   {coordinate matrix of set of vectors}%
\Index
   {coordinate matrix of vector}%
   {coordinate matrix of vector}%
\Index
   {coordinate matrix of vector field in \rcD basis}%
   {coordinate matrix of vector field in drc basis}%
\Index
   {coordinate \rcd vector space}%
   {coordinate rcd vector space}%
\Index
   {coordinate reference frame}%
   {coordinate reference frame}%
\Index
   {coordinate representation in $\Omega_2$\Hyph algebra}%
   {coordinate representation, Omega_2 algebra}%
\Index
   {coordinate representation in \rcd vector space}%
   {coordinate representation, rcd vector space}%
\Index
   {coordinate representation in tuple of $\VX\Omega$\Hyph algebras}%
   {coordinate tower of representations, Omega algebra}%
\Index
   {coordinate representation of group in vector space}%
   {coordinate representation, vector space}%
\Index
   {coordinate representation of vector}%
   {coordinate representation of vector}%
\Index
   {coordinate vector bundle}%
   {coordinate vector bundle}%
\Index
   {coordinate vector space}%
   {coordinate vector space}%
\Index
   {coordinates}%
   {coordinates}%
\Index
   {coordinates of a geometric object in $\Omega_2$\Hyph algebra $M$}%
   {coordinates of geometric object, representation g}%
\Index
   {coordinates of a geometric object in tuple of $\VX\Omega$\Hyph algebras}%
   {coordinates of geometric object, tower of representations g}%
\Index
   {coordinates of basis}%
   {coordinates of basis}%
\Index
   {coordinates of basis of representation}%
   {coordinates of basis relative to basis, representation}%
\Index
   {coordinates of element $m$ of representation $f$ relative to set $X$}%
   {coordinates of element relative to set, representation}%
\Index
   {coordinates of endomorphism of representation}%
   {coordinates of endomorphism, representation}%
\Index
   {coordinates of endomorphism of tower of representations}%
   {coordinates of endomorphism, tower of representations}%
\Index
   {coordinates of geometric object}%
   {coordinates of geometric object, vector space}%
\Index
   {coordinates of geometric object in coordinate \rcd vector space}%
   {coordinates of geometric object, coordinate rcd vector space}%
\Index
   {coordinates of geometric object in coordinate representation}%
   {coordinates of geometric object, coordinate vector space}%
\Index
   {coordinates of geometric object in coordinate space of representation}%
   {coordinates of geometric object, coordinate representation}%
\Index
   {coordinates of geometric object in coordinate space of tower of representations}%
   {coordinates of geometric object, coordinate tower of representations}%
\Index
   {coordinates of geometric object in \rcd vector space}%
   {coordinates of geometric object, rcd vector space}%
\Index
   {coordinates of point $A$ of affine space $\overset{\circ}{A}$ relative to basis $(O,\Basis e)$}%
   {coordinates in affine space}%
\Index
   {coordinates of representation}%
   {coordinates of representation, drc vector space}%
\Index
   {coordinates of representation}%
   {coordinates of representation}%
\Index
   {coordinates of set of vectors}%
   {coordinates of set of vectors}%
\Index
   {coordinates of vector}%
   {coordinates of vector}%
\Index
   {coordinates of vector field in \Drc basis}%
   {coordinates of vector field in drc basis}%
\Index
   {coordinates of vector relative to Hamel basis}%
   {coordinates of vector, Hamel basis}%
\Index
   {coordinates of vector relative to Schauder basis}%
   {coordinates of vector, Schauder basis}%
\Index
   {coproduct of objects in category}%
   {coproduct in category}%
\Index
   {correspondence continuous on the set}%
   {correspondence continuous on the set}%
\Index
   {correspondence of homomorphism}%
   {correspondence of homomorphism}%
\Index
   {cosine}%
   {cosine}%
\Index
   {\CR inverse element of biring}%
   {cr-inverse element}%
\Index
   {\CR matrix group}%
   {cr-matrix group}%
\Index
   {\CR power}%
   {cr power}%
\Index
   {$\CRcirc$\Hyph product of matrices of maps}%
   {cr product of matrices of maps}%
\Index
   {\crd vector}%
   {crd vector}%
\Index
   {\crd vector space}%
   {crd vector space}%
\Index
   {$C^*$\Hyph algebra}%
   {Cstar-algebra}%
\Index
   {curvilinear coordinates of point in affine space}%
   {curvilinear coordinates of point in affine space}%
\SetIndexSpace%
\Index
   {$D$\Hyph linear functional}%
   {D linear functional}%
\Index
   {$D*$\hyph matrices vector space}%
   {matrices vector space}%
\Index
   {$D*$\hyph  vector space}%
   {D* vector space}%
\Index
   {$D*$\Hyph module}%
   {D*-module}%
\Index
   {$D$\Hyph affine connection on manifold with affine connections}%
   {D affine connection, affine manifold}%
\Index
   {$D$\Hyph algebra}%
   {D algebra}%
\Index
   {$D$\Hyph module}%
   {D-module}%
\Index
   {$D$\Hyph module}%
   {D module}%
\Index
   {$D$\Hyph valued variable}%
   {D valued variable}%
\Index
   {$D$\Hyph vector function}%
   {d vector function}%
\Index
   {$D$\Hyph affine connection coefficients on manifold}%
   {D affine connection coefficients, manifold}%
\Index
   {$D$\hyph vector space}%
   {D vector space}%
\Index
   {\dcr vector}%
   {dcr vector}%
\Index
   {\dcr vector space}%
   {dcr vector space}%
\Index
   {definite integral}%
   {definite integral}%
\Index
   {derivative of map}%
   {derivative of map}%
\Index
   {derivative of order $n$}%
   {derivative of Order n}%
\Index
   {derivative of second order}%
   {derivative of Second Order}%
\Index
   {determinant of matrix}%
   {determinant}%
\Index
   {deviation of trajectories}%
   {deviation of trajectories}%
\Index
   {diagonal in bundle}%
   {diagonal in bundle}%
\Index
   {diagram of correspondences}%
   {diagram of correspondences}%
\Index
   {diagram of representations}%
   {diagram of representations}%
\Index
   {differentiable map}%
   {differentiable map}%
\Index
   {differential form of degree $p$}%
   {differential form of degree p}%
\Index
   {differential of independent variable}%
   {differential of independent variable}%
\Index
   {differential of map}%
   {differential of map}%
\Index
   {differential $p$\Hyph form}%
   {differential p form}%
\Index
   {dimension of \rcd vector space}%
   {dimension of vector space}%
\Index
   {direct product of bundles}%
   {Cartesian product of bundles}%
\Index
   {direct product of $D$\Hyph vector spaces}%
   {direct product of D vector spaces}%
\Index
   {direct product of division rings}%
   {direct product of division rings}%
\Index
   {direct product of \Ts{G}representations}%
   {direct product of G* representations}%
\Index
   {direct product of \(\Omega\)\Hyph algebras}%
   {direct product of Omega algebras}%
\Index
   {direct product of \rcd vector spaces}%
   {direct product, rcd vector space}%
\Index
   {direct product of representations of fibered group}%
   {direct product of representations of fibered group}%
\Index
   {direct product of representations of group}%
   {direct product of representations of group}%
\Index
   {direct product of total spaces}%
   {Cartesian product of total spaces}%
\Index
   {direct sum}%
   {direct sum}%
\Index
   {direct sum of representations}%
   {direct sum of representations}%
\Index
   {direction over commutative ring}%
   {direction over commutative ring}%
\Index
   {distributive law}%
   {distributive law}%
\Index
   {division algebra}%
   {division algebra}%
\Index
   {division with remainder}%
   {division with remainder}%
\Index
   {division without remainder}%
   {division without remainder}%
\Index
   {divisor of polynomial}%
   {divisor of polynomial}%
\Index
   {double determinant}%
   {double determinant}%
\Index
   {\Drc linear map of vector bundles}%
   {drc linear map of vector bundles}%
\Index
   {\drc vector}%
   {drc vector}%
\Index
   {\drc vector space}%
   {drc vector space}%
\Index
   {$D\star$\Hyph antilinear homomorphism}%
   {Dstar antilinear homomorphism}%
\Index
   {$\mathcal D\star$\Hyph vector bundle}%
   {Dstar vector bundle}%
\Index
   {$\mathcal D\star$\Hyph vector field}%
   {Dstar vector field}%
\Index
   {$\mathcal D\star$\hyph linear composition of vector fields}%
   {linear composition of vector fields}%
\Index
   {$\mathcal D\star$\hyph product of vector field over scalar}%
   {Dstar product of vector field over scalar, vector space}%
\Index
   {dual space of \rcd vector space}%
   {dual space of rcd vector space}%
\Index
   {duality principle for biring}%
   {duality principle for biring}%
\Index
   {duality principle for biring of matrices}%
   {duality principle for biring of matrices}%
\SetIndexSpace%
\Index
   {effective \Ts{G}representation}%
   {effective G* representation}%
\Index
   {effective representation}%
   {effective representation}%
\Index
   {effective representation of division ring}%
   {effective representation of division ring}%
\Index
   {effective representation of fibered $\Omega$\Hyph algebra}%
   {effective representation of fibered Omega-algebra}%
\Index
   {effective representation of group}%
   {effective representation of group}%
\Index
   {effective representation of ring}%
   {effective representation of ring}%
\Index
   {effective \Ts representation of fibered division ring}%
   {effective representation of fibered division ring}%
\Index
   {effective \Ts representation of fibered group}%
   {effective representation of fibered group}%
\Index
   {endomorphism}%
   {endomorphism}%
\Index
   {endomorphism of representation of $\Omega$\Hyph algebra}%
   {endomorphism of representation}%
\Index
   {endomorphism of representation regular on generating set $X$}%
   {endomorphism of representation, regular on set}%
\Index
   {endomorphism of representation singular on generating set $X$}%
   {endomorphism of representation, singular on set}%
\Index
   {endomorphism of tower of representations}%
   {endomorphism of tower of representations}%
\Index
   {endomorphism of tower of representations regular on tuple of generating sets}%
   {endomorphism of representation, regular on tuple}%
\Index
   {endomorphism of tower of representations singular on tuple of generating sets}%
   {endomorphism of representation, singular on tuple}%
\Index
   {enhanced Lie group}%
   {enhanced Lie group}%
\Index
   {epimorphism}%
   {epimorphism}%
\Index
   {equivalence generated by representation $f$}%
   {equivalence of representation}%
\Index
   {equivalent norms}%
   {equivalent norms}%
\Index
   {essential parameters in a set of functions}%
   {essential parameters}%
\Index
   {Euclidean metric on division ring}%
   {Euclidean metric on division ring}%
\Index
   {Euclidean scalar product in $D$\Hyph vector space}%
   {Euclidean scalar product, vector space}%
\Index
   {Euclidean scalar product on division ring}%
   {Euclidean scalar product on division ring}%
\Index
   {everywhere dense subset}%
   {everywhere dense subset}%
\Index
   {expansion of vector relative to basis converges}%
   {expansion converges}%
\Index
   {expansion of vector relative to basis converges normally}%
   {expansion converges normally}%
\Index
   {exponent}%
   {exponent}%
\Index
   {extended matrix of \drc linear equations}%
   {extended matrix, system of drc linear equations}%
\Index
   {extended matrix of \rcd linear equations}%
   {extended matrix, system of rcd linear equations}%
\Index
   {extension of correspondence}%
   {extension of correspondence}%
\Index
   {extension of measure}%
   {extension of measure}%
\Index
   {exterior differential}%
   {exterior differential}%
\Index
   {exterior product}%
   {exterior product}%
\Index
   {extreme line}%
   {extreme line}%
\SetIndexSpace%
\Index
   {fibered coordinate isomorphism}%
   {fibered coordinate isomorphism}%
\Index
   {fibered correspondence from $\Bundle A$ to $\Bundle B$}%
   {fibered correspondence from A to B}%
\Index
   {fibered correspondence in $\Bundle{A}$}%
   {fibered correspondence in A}%
\Index
   {fibered correspondence of homomorphism}%
   {fibered correspondence of homomorphism}%
\Index
   {fibered equivalence}%
   {fibered equivalence}%
\Index
   {fibered group}%
   {fibered group}%
\Index
   {fibered identification morphism}%
   {fibered identification morphism}%
\Index
   {fibered little group}%
   {fibered little group}%
\Index
   {fibered morphism from bundle $\Bundle A$ into $\Bundle B$}%
   {fibered morphism from A into B}%
\Index
   {fibered natural morphism}%
   {fibered natural morphism}%
\Index
   {fibered $\Omega$\Hyph algebra}%
   {fibered Omega-algebra}%
\Index
   {fibered $\Omega$\Hyph subalgebra}%
   {fibered Omega-subalgebra}%
\Index
   {fibered ordering}%
   {fibered ordering}%
\Index
   {fibered preordering}%
   {fibered preordering}%
\Index
   {fibered ring}%
   {fibered ring}%
\Index
   {fibered stability group}%
   {fibered stability group}%
\Index
   {fibered subset}%
   {fibered subset}%
\Index
   {field-strength tensor}%
   {field-strength tensor}%
\Index
   {filter $\mathfrak{F}$ converges to $A$}%
   {filter converges}%
\Index
   {finite expansion of set}%
   {finite expansion of set}%
\Index
   {Finsler metric}%
   {Finsler metric}%
\Index
   {Finsler space}%
   {Finsler space}%
\Index
   {Finsler structure}%
   {Finsler structure}%
\Index
   {first Newton law}%
   {First Newton law}%
\Index
   {frame\Hyph dragging effect}%
   {frame dragging effect}%
\Index
   {free $A$\Hyph module}%
   {free A module}%
\Index
   {free Abelian group}%
   {free Abelian group}%
\Index
   {free algebra over ring}%
   {free algebra over ring}%
\Index
   {free module}%
   {free module}%
\Index
   {free representation}%
   {free representation}%
\Index
   {free representation of group}%
   {free representation of group}%
\Index
   {free \Ts representation of fibered group}%
   {free representation of fibered group}%
\Index
   {Frenet transport}%
   {Frenet transport}%
\Index
   {function homogeneous of degree $k$}%
   {function homogeneous}%
\Index
   {function of division ring \Ds differentiable in the Fr\'echet sense}%
   {function Dstar differentiable in Frechet sense, division ring}%
\Index
   {fundamental sequence}%
   {fundamental sequence}%
\SetIndexSpace%
\Index
   {$G$\Hyph reference frame}%
   {G reference frame}%
\Index
   {$G$\Hyph basis of vector space}%
   {G-basis}%
\Index
   {$G$\Hyph coordinates of basis}%
   {G-coordinates}%
\Index
   {$G$\Hyph space}%
   {GSpace}%
\Index
   {the G\^ateaux \dcr derivative of map $f$ of $D$\Hyph vector space $V$ to $D$\Hyph vector space $W$}%
   {Gateaux dcr derivative of map, D vector space}%
\Index
   {the G\^ateaux derivative of map}%
   {Gateaux derivative of map}%
\Index
   {the G\^ateaux derivative of order $n$}%
   {Gateaux derivative of Order n}%
\Index
   {the G\^ateaux derivative of second order}%
   {Gateaux derivative of Second Order}%
\Index
   {the G\^ateaux \Ds derivative of map $f$ of division ring $D$}%
   {Gateaux Dstar derivative of map, division ring}%
\Index
   {the G\^ateaux mixed partial derivative}%
   {Gateaux partial derivative of Second Order}%
\Index
   {the G\^ateaux partial \dcr derivative of map $f^{\gi b}$ with respect to variable $x^{\gi a}$}%
   {Gateaux partial dcr derivative of map with respect to variable, D vector space}%
\Index
   {the G\^ateaux partial derivative}%
   {Gateaux partial derivative}%
\Index
   {the G\^ateaux partial \rcd derivative of map $f^{\gi b}$ with respect to variable $x^{\gi a}$}%
   {Gateaux partial rcd derivative of map with respect to variable, D vector space}%
\Index
   {the G\^ateaux \rcd derivative of map $f$ of $D$\hyph vector space $V$ to $D$\hyph vector space $W$}%
   {Gateaux rcd derivative of map, D vector space}%
\Index
   {the G\^ateaux \sD derivative of map $f$ of division ring $D$}%
   {Gateaux starD derivative of map, division ring}%
\Index
   {generating set}%
   {generating set}%
\Index
   {generator of linear map}%
   {generator of linear map}%
\Index
   {geodetic effect}%
   {geodetic effect}%
\Index
   {geometric object defined in $\Omega_2$\Hyph algebra $M$}%
   {geometric object, representation g}%
\Index
   {geometric object defined in \rcd vector space}%
   {geometric object, rcd vector space}%
\Index
   {geometric object defined in tuple of $\VX\Omega$\Hyph algebras $\VX A$}%
   {geometric object, tower of representations g}%
\Index
   {geometric object in coordinate representation}%
   {geometric object, coordinate vector space}%
\Index
   {geometric object in coordinate representation defined in $\Omega_2$\Hyph algebra $M$}%
   {geometric object, coordinate representation g}%
\Index
   {geometric object in coordinate representation defined in \rcd vector space}%
   {geometric object, coordinate rcd vector space}%
\Index
   {geometric object in coordinate representation defined in tuple of $\VX\Omega$\Hyph algebras $\VX A$}%
   {geometric object, coordinate tower of representations g}%
\Index
   {geometric object in vector space}%
   {geometric object, vector space}%
\Index
   {geometric object of type $H$}%
   {geometric object of type H, representation g}%
\Index
   {geometric object of type $A$ in vector space}%
   {geometric object of type A, vector space}%
\Index
   {group algebra}%
   {group algebra}%
\Index
   {group of automorphisms of representation}%
   {group of automorphisms of representation}%
\SetIndexSpace%
\Index
   {Hadamard inverse of matrix}%
   {Hadamard inverse of matrix}%
\Index
   {Hamel basis}%
   {Hamel basis}%
\Index
   {hermitian conjugated vector}%
   {hermitian conjugated vector}%
\Index
   {hermitian conjugation in division ring}%
   {hermitian conjugation, division ring}%
\Index
   {hermitian matrix}%
   {hermitian matrix}%
\Index
   {hermitian metric on division ring}%
   {hermitian metric on division ring}%
\Index
   {hermitian scalar product in $D$\Hyph vector space}%
   {hermitian scalar product, vector space}%
\Index
   {hermitian scalar product on division ring}%
   {hermitian scalar product on division ring}%
\Index
   {highest common factor}%
   {highest common factor}%
\Index
   {holomorphic map}%
   {holomorphic map}%
\Index
   {holonomic coordinates of connection}%
   {holonomic coordinates of connection}%
\Index
   {holonomic coordinates of vector}%
   {vector holonomic coordinates}%
\Index
   {homogeneous bundle of fibered group}%
   {homogeneous bundle of fibered group}%
\Index
   {homogeneous linear geometric object}%
   {homogeneous linear geometric object}%
\Index
   {homogeneous map of degree $k$ over field $F$}%
   {homogeneous map of degree over field, D vector space}%
\Index
   {homogeneous polynomial of power $k$}%
   {homogeneous polynomial of power}%
\Index
   {homogeneous space}%
   {homogeneous space}%
\Index
   {homomorphic image}%
   {homomorphic image}%
\Index
   {homomorphism}%
   {homomorphism}%
\Index
   {homomorphism of fibered groups}%
   {homomorphism of fibered groups}%
\Index
   {homomorphism of fibered universal algebras}%
   {homomorphism of fibered universal algebras}%
\Index
   {horizontal component of vector}%
   {horizontal component of vector}%
\Index
   {horizontal subspace}%
   {horizontal subspace}%
\Index
   {horizontal vector}%
   {horizontal vector}%
\Index
   {hyperbolic cosine}%
   {hyperbolic cosine}%
\Index
   {hyperbolic sine}%
   {hyperbolic sine}%
\SetIndexSpace%
\Index
   {ideal of algebra}%
   {ideal of algebra}%
\Index
   {indefinite integral}%
   {indefinite integral}%
\Index
   {independent points}%
   {independent points}%
\Index
   {infinitesimal generator of representation}%
   {infinitesimal generator}%
\Index
   {infinitesimal generators of group Lie}%
   {infinitesimal generators of group Lie}%
\Index
   {integrable form}%
   {integrable form}%
\Index
   {integrable map}%
   {integrable map}%
\Index
   {integral of differential $1$\Hyph form along path}%
   {integral of differential 1 form along path}%
\Index
   {invariance principle in \drc vector space}%
   {invariance principle}%
\Index
   {invariance principle in representation of universal algebra}%
   {invariance principle, representation g}%
\Index
   {invariance principle in tower of representations of universal algebras}%
   {invariance principle, tower of representations g}%
\Index
   {invariance principle in vector space}%
   {invariance principle, vector space}%
\Index
   {inverse fibered correspondence}%
   {inverse fibered correspondence}%
\Index
   {inverse reduced fibered correspondence}%
   {inverse reduced fibered correspondence}%
\Index
   {involution in quaternion algebra}%
   {involution, quaternion algebra}%
\Index
   {isomorphism}%
   {isomorphism}%
\Index
   {isomorphism of fibered $\Omega$\Hyph algebras}%
   {isomorphism of fibered Omega-algebras}%
\Index
   {isomorphism of repesentations of $\Omega$\Hyph algebra}%
   {isomorphism of repesentations of Omega algebra}%
\Index
   {isomorphism of vector spaces}%
   {isomorphism of vector spaces}%
\Index
   {isotropic vector}%
   {isotropic vector}%
\Index
   {Lebesgue integral}%
   {Lebesgue integral}%
\SetIndexSpace%
\Index
   {$(^j_i)$\hyph $\RCcirc$\Hyph quasideterminant}%
   {j i RCcirc-quasideterminant}%
\Index
   {the Jacobi matrix of map}%
   {Jacobi matrix of map}%
\Index
   {Jacobian complete system of differential equations}%
   {Jacobian complete system of differential equations}%
\Index
   {Jacobian complete system of \drv differential equations}%
   {Jacobian complete system of drc differential equations}%
\Index
   {$(ji)$\hyph quasideterminant}%
   {j i quasideterminant}%
\Index
   {the Jacobi\Hyph G\^ateaux matrix of map}%
   {Jacobi Gateaux matrix of map}%
\SetIndexSpace%
\Index
   {kernel of homomorphism}%
   {kernel of homomorphism}%
\Index
   {kernel of inefficiency of \Ts{G}representation}%
   {kernel of inefficiency of G* representation}%
\Index
   {kernel of inefficiency of representation of fibered group}%
   {kernel of inefficiency of representation of fibered group}%
\Index
   {kernel of inefficiency of representation of group}%
   {kernel of inefficiency of representation of group}%
\Index
   {kernel of linear map}%
   {kernel of linear map}%
\Index
   {Killing equation}%
   {Killing equation}%
\Index
   {Killing equation of second type}%
   {Killing equation second type}%
\Index
   {Killing vector of second type}%
   {Killing vector second type}%
\Index
   {Kronecker symbol}%
   {Kronecker symbol}%
\SetIndexSpace%
\Index
   {latitude}%
   {latitude}%
\Index
   {leading coefficient of polynomial}%
   {leading coefficient of polynomial}%
\Index
   {Lebesgue extension of measure}%
   {Lebesgue extension of measure}%
\Index
   {Lebesgue measurable set}%
   {Lebesgue measurable}%
\Index
   {Lebesgue measure}%
   {Lebesgue measure}%
\Index
   {left $A$\Hyph module}%
   {left A module}%
\Index
   {left $A$\Hyph vector space}%
   {left A vector space}%
\Index
   {left cofactor of entry of matrix}%
   {left cofactor, matrix}%
\Index
   {left $D$\hyph vector space of columns}%
   {left vector space of columns}%
\Index
   {left $D$\hyph vector space of rows}%
   {left vector space of rows}%
\Index
   {left defined Lie algebra of Lie group}%
   {left defined Lie algebra}%
\Index
   {left double cofactor of entry of matrix}%
   {left double cofactor}%
\Index
   {left fraction}%
   {left fraction}%
\Index
   {left ideal of algebra}%
   {left ideal of algebra}%
\Index
   {left invariant vector field}%
   {left invariant vector}%
\Index
   {left module}%
   {left module}%
\Index
   {left principal ideal}%
   {left principal ideal}%
\Index
   {left shift of module}%
   {left shift of module}%
\Index
   {left shift on fibered group}%
   {left shift, fibered group}%
\Index
   {left shift on group}%
   {left shift}%
\Index
   {left shift on group}%
   {left shift, group}%
\Index
   {left structural constant of Lie algebra}%
   {left structural constant of Lie algebra}%
\Index
   {left vector space}%
   {left vector space}%
\Index
   {left zero divisor}%
   {left zero divisor}%
\Index
   {left-ordered cycle notation of permutation}%
   {left-ordered cycle notation of permutation}%
\Index
   {left\Hyph side $A_1$\Hyph representation}%
   {left-side A representation}%
\Index
   {left\Hyph side product}%
   {left-side product}%
\Index
   {left-side product of map over scalar}%
   {left-side product of map over scalar}%
\Index
   {left\Hyph side product of vector over scalar}%
   {left-side product of vector over scalar}%
\Index
   {left-side representation}%
   {left-side representation}%
\Index
   {left-side representation of fibered $\Omega$\Hyph algebra}%
   {left-side representation of fibered Omega-algebra}%
\Index
   {left-side representation of $\Omega_1$\Hyph algebra $A$ in $\Omega_2$\Hyph algebra $M$}%
   {left-side representation of algebra}%
\Index
   {left-side transformation}%
   {left-side transformation}%
\Index
   {left-side transformation on bundle}%
   {left-side transformation of bundle}%
\Index
   {Lie algebra of Lie group}%
   {algebra Lie group Lie}%
\Index
   {Lie derivative}%
   {Lie derivative}%
\Index
   {Lie derivative of connection}%
   {Lie derivative of connection}%
\Index
   {Lie derivative of metric}%
   {Lie derivative of metric}%
\Index
   {Lie group basic operators}%
   {Lie group basic operators}%
\Index
   {lift of correspondence}%
   {lift of correspondence}%
\Index
   {lift of mapping}%
   {lift of map}%
\Index
   {limit of correspondence with respect to the filter}%
   {limit of correspondence with respect to the filter}%
\Index
   {limit of filter}%
   {limit of filter}%
\Index
   {limit of sequence}%
   {limit of sequence}%
\Index
   {limit set of filter}%
   {limit set of filter}%
\Index
   {linear combination}%
   {linear combination}%
\Index
   {linear functional}%
   {linear functional}%
\Index
   {linear \Ts{G}representation}%
   {linear G* representation}%
\Index
   {linear geometric object}%
   {linear geometric object}%
\Index
   {linear homomorphism}%
   {linear homomorphism}%
\Index
   {linear map}%
   {linear map}%
\Index
   {linear map generated by map}%
   {linear map generated by map}%
\Index
   {linear map of division ring}%
   {linear map of division ring}%
\Index
   {linear representation of group}%
   {linear representation of group}%
\Index
   {linear representation of Lie group}%
   {linear representation of Lie group}%
\Index
   {linear span in vector space}%
   {linear span, vector space}%
\Index
   {linear transformation group}%
   {linear transformation group}%
\Index
   {linear transformation of affine space}%
   {linear transformation, affine space}%
\Index
   {linearly dependent}%
   {linearly dependent}%
\Index
   {linearly dependent set}%
   {linearly dependent set}%
\Index
   {linearly dependent vector fields}%
   {linearly dependent vector fields}%
\Index
   {linearly independent set}%
   {linearly independent set}%
\Index
   {little group}%
   {little group}%
\Index
   {local reference frame}%
   {local reference frame}%
\Index
   {locally compact at point $p$ space}%
   {locally compact at point space}%
\Index
   {locally compact space}%
   {locally compact space}%
\Index
   {longitude}%
   {longitude}%
\Index
   {Lorentz transformation}%
   {Lorentz transformation}%
\SetIndexSpace%
\Index
   {$m$\Hyph dimensional parallelepiped}%
   {m dimensional parallelepiped}%
\Index
   {$m$\Hyph vector}%
   {m-vector}%
\Index
   {manifold with $D$\Hyph affine connections}%
   {manifold with D- affine connections}%
\Index
   {map continuous with respect to set of arguments}%
   {map continuous with respect to set of arguments}%
\Index
   {map differentiable in the G\^ateaux sense}%
   {map differentiable in Gateaux sense}%
\Index
   {map is compatible with operation}%
   {map is compatible with operation}%
\Index
   {map of conjugation}%
   {map of conjugation}%
\Index
   {map of $\gi n$ $D$\Hyph valued variables}%
   {map of n D valued variables}%
\Index
   {map of type $G$ on manifold}%
   {map of type G on manifold}%
\Index
   {map polylinear over finite dimensional algebras}%
   {map polylinear over finite dimensional algebras}%
\Index
   {map projective over commutative ring}%
   {map projective over commutative ring}%
\Index
   {mapping of rings polylinear over commutative ring}%
   {map polylinear over commutative ring, ring}%
\Index
   {mapping space}%
   {mapping space}%
\Index
   {matrix}%
   {matrix}%
\Index
   {matrix of antilinear homomorphism}%
   {matrix of antilinear homomorphism}%
\Index
   {matrix of bilinear function}%
   {matrix of bilinear function}%
\Index
   {matrix of endomorphisms of $\Omega$\Hyph algebra}%
   {matrix of endomorphisms of Omega algebra}%
\Index
   {matrix of fibered \Drc linear map}%
   {matrix of fibered drc linear map}%
\Index
   {matrix of linear homomorphism}%
   {matrix of linear homomorphism}%
\Index
   {matrix of linear map}%
   {matrix of linear map}%
\Index
   {matrix of linear maps}%
   {matrix of linear maps}%
\Index
   {matrix of maps}%
   {matrix of maps}%
\Index
   {matrix of quadratic map}%
   {matrix of quadratic map, division ring}%
\Index
   {measurable map}%
   {measurable map}%
\Index
   {measure}%
   {measure}%
\Index
   {method of successive differentiation}%
   {method of successive differentiation}%
\Index
   {metric tensor in Minkowski space}%
   {metric tensor, Minkowski space}%
\Index
   {metric-affine manifold}%
   {metric-affine manifold}%
\Index
   {Minkowski space}%
   {Minkowski space, Finsler}%
\Index
   {minor matrix}%
   {minor matrix}%
\Index
   {module over ring}%
   {module over ring}%
\Index
   {monomial of power $k$}%
   {monomial of power}%
\Index
   {monomorphism}%
   {monomorphism}%
\Index
   {morphism from tower of representations into tower of representations}%
   {morphism from tower of representations into tower of representations}%
\Index
   {morphism of fibered \Ts representations from $\Bundle F$ into $\Bundle G$}%
   {morphism of fibered representations from f into g}%
\Index
   {morphism of representation $f$}%
   {morphism of representation f}%
\Index
   {morphism of representations from $f$ into $g$}%
   {morphism of representations from f into g}%
\Index
   {morphism of representations of $\Omega_1$\Hyph algebra in $\Omega_2$\Hyph algebra}%
   {morphism of representations of Omega1 algebra in Omega2 algebra}%
\Index
   {morphism of \Ts representations of fibered $\Omega$\Hyph algebra}%
   {morphism of representations of fibered Omega algebra}%
\Index
   {motion of Minkowski space}%
   {motion, Minkowski space}%
\Index
   {movement on basis manifold}%
   {movement transformation}%
\Index
   {multiplicative map}%
   {multiplicative map}%
\Index
   {multiplicative $\Omega$\Hyph group}%
   {multiplicative Omega group}%
\SetIndexSpace%
\Index
   {$n$\Hyph ary fibered relation}%
   {fibered relation}%
\Index
   {$n$\Hyph ary operation on set}%
   {n-ary operation on set}%
\Index
   {natural homomorphism}%
   {natural homomorphism}%
\Index
   {neutral element of operation}%
   {neutral element of operation}%
\Index
   {nonmetricity}%
   {nonmetricity}%
\Index
   {nonsingular bilinear function}%
   {nonsingular bilinear function}%
\Index
   {nonsingular system of \rcd linear equations}%
   {nonsingular system of linear equations}%
\Index
   {nonsingular tensor}%
   {nonsingular tensor}%
\Index
   {nonsingular transformation}%
   {nonsingular transformation}%
\Index
   {norm in quaternion algebra}%
   {norm, quaternion algebra}%
\Index
   {norm of functional}%
   {norm of functional}%
\Index
   {norm of map}%
   {norm of map}%
\Index
   {norm of operation}%
   {norm of operation}%
\Index
   {norm of polylinear map}%
   {norm of polymap}%
\Index
   {norm of representation}%
   {norm of representation}%
\Index
   {norm on $D$\Hyph algebra}%
   {norm on D algebra}%
\Index
   {norm on $D$\Hyph vector space}%
   {norm on D vector space}%
\Index
   {norm on $D$\Hyph module}%
   {norm on D module}%
\Index
   {norm on $\Omega$\Hyph group}%
   {norm on Omega group}%
\Index
   {norm on ring}%
   {norm on ring}%
\Index
   {normal basis}%
   {normal basis}%
\Index
   {normed $D$\Hyph algebra}%
   {normed D algebra}%
\Index
   {normed $D$\Hyph module}%
   {normed D module}%
\Index
   {normed $D$\Hyph vector space}%
   {normed D vector space}%
\Index
   {normed $\Omega$\Hyph group}%
   {normed Omega group}%
\Index
   {normed ring}%
   {normed ring}%
\Index
   {not complete group}%
   {not complete group}%
\Index
   {not complete $\Omega$\Hyph algebra}%
   {not complete Omega algebra}%
\Index
   {nucleus of $D$\Hyph algebra $A$}%
   {nucleus of algebra}%
\SetIndexSpace%
\Index
   {octonion algebra}%
   {octonion algebra}%
\Index
   {open ball}%
   {open ball}%
\Index
   {open set}%
   {open set}%
\Index
   {operation on bundle}%
   {operation on bundle}%
\Index
   {operation on set}%
   {operation on set}%
\Index
   {operator domain}%
   {operator domain}%
\Index
   {opposite algebra to algebra $P$}%
   {opposite algebra}%
\Index
   {opposite fibered preordering}%
   {opposite fibered preordering}%
\Index
   {orbit of linear map}%
   {orbit of linear map}%
\Index
   {orbit of representation}%
   {orbit of representation}%
\Index
   {orbit of representation of fibered group}%
   {orbit of representation of fibered group}%
\Index
   {orbit of representation of group}%
   {orbit of representation of group}%
\Index
   {origin of coordinate system of affine space}%
   {origin of coordinate system of affine space}%
\Index
   {orthogonal basis in Minkowski space}%
   {orthogonal basis, Minkowski space}%
\Index
   {orthogonality in Minkowski space}%
   {Minkowski orthogonality}%
\Index
   {orthonormal basis}%
   {Orthonormal Basis, division ring}%
\Index
   {orthonormal basis in Minkowski space}%
   {orthonormal basis, Minkowski space}%
\Index
   {orthonornal basis}%
   {Orthonornal Basis}%
\Index
   {outer measure}%
   {outer measure}%
\SetIndexSpace%
\Index
   {passive representation of group $G(f)$ in basis manifold of representation}%
   {passive representation in basis manifold}%
\Index
   {parallel shift of affine space}%
   {parallel shift, affine space}%
\Index
   {parallelogram}%
   {parallelogram}%
\Index
   {parity of permutation}%
   {parity of permutation}%
\Index
   {partial derivative}%
   {partial derivative}%
\Index
   {partial derivative of second order}%
   {partial derivative of second order}%
\Index
   {partial linear map}%
   {partial linear map}%
\Index
   {passive \sT{G}representation}%
   {passive *G representation}%
\Index
   {passive representation}%
   {passive representation}%
\Index
   {passive representation of group $G(\Vector f)$ in basis manifold of tower of representations}%
   {passive representation in basis manifold, tower of representations}%
\Index
   {passive transformation of the basis manifold of representation}%
   {passive transformation of basis, representation}%
\Index
   {passive transformation of the basis manifold of tower of representations}%
   {passive transformation of basis, tower of representations}%
\Index
   {passive transformation on basis manifold}%
   {passive transformation}%
\Index
   {passive transformation on the set of \rcd bases}%
   {passive transformation, vector space}%
\Index
   {permutability property of trace}%
   {permutability property of trace}%
\Index
   {permutation}%
   {permutation}%
\Index
   {pfaffian derivative}%
   {pfaffian derivative}%
\Index
   {polyadditive map}%
   {polyadditive map}%
\Index
   {polylinear map}%
   {polylinear map}%
\Index
   {polylinear skew symmetric map}%
   {polylinear map skew symmetric}%
\Index
   {polylinear symmetric map}%
   {polylinear map symmetric}%
\Index
   {polymorphism of representations}%
   {polymorphism of representations}%
\Index
   {polynomial}%
   {polynomial}%
\Index
   {polyvector}%
   {polyvector}%
\Index
   {potential energy}%
   {potential energy}%
\Index
   {power of measure}%
   {power of measure}%
\Index
   {prime $A$\Hyph number}%
   {prime A number}%
\Index
   {principal ideal}%
   {principal ideal}%
\Index
   {product in category}%
   {product in category}%
\Index
   {product of geometric object and constant}%
   {product of geometric object and constant}%
\Index
   {product of geometric object and constant in vector space}%
   {product of geometric object and constant, vector space}%
\Index
   {product of measures}%
   {product of measures}%
\Index
   {product of morphisms of representations of universal algebra}%
   {product of morphisms of representations of universal algebra}%
\Index
   {product of morphisms of tower of representations}%
   {product of morphisms of tower of representations}%
\Index
   {product of morphisms of \Ts representations of fibered $\Omega$\Hyph algebra}%
   {product of morphisms of representations of fibered Omega algebra}%
\Index
   {product of polynomials}%
   {product of polynomials}%
\Index
   {product of rings of sets}%
   {product of rings of sets}%
\Index
   {projection of bundle $\Bundle E$ along fiber $E$}%
   {projection of bundle along fiber}%
\Index
   {projective map is continuous in direction over field}%
   {projective map is continuous in direction over field}%
\Index
   {pseudo\Hyph Euclidean metric on division ring}%
   {pseudo-Euclidean metric on division ring}%
\Index
   {pseudo\Hyph Euclidean scalar product in $D$\Hyph vector space}%
   {pseudo-Euclidean scalar product, vector space}%
\Index
   {pseudo-Euclidean scalar product on division ring}%
   {pseudo-Euclidean scalar product on division ring}%
\SetIndexSpace%
\Index
   {quadratic equation}%
   {quadratic equation}%
\Index
   {quadratic form in division ring}%
   {quadratic form, division ring}%
\Index
   {quadratic map of division ring}%
   {Quadratic Map of Division Ring}%
\Index
   {quasi affine transformation on basis manifold}%
   {quasi affine transformation}%
\Index
   {quasi affine transformation on basis manifold}%
   {quasi affine drc transformation}%
\Index
   {quasi movement on basis manifold}%
   {quasi movement, division ring}%
\Index
   {quasi movement on basis manifold}%
   {quasi movement}%
\Index
   {quasiclosed ring of maps}%
   {quasiclosed ring of maps}%
\Index
   {quasideterminant}%
   {quasideterminant definition}%
\Index
   {quasimotion of Minkowski space}%
   {Quasimotion, Minkowski space}%
\Index
   {quaternion algebra}%
   {quaternion algebra}%
\Index
   {quaternion algebra $E$ over the field $F$}%
   {quaternion algebra over the field}%
\Index
   {quotient}%
   {quotient divided by}%
\Index
   {quotient bundle}%
   {quotient bundle}%
\SetIndexSpace%
\Index
   {$(^j_i)$\hyph \RC quasideterminant}%
   {j i RC-quasideterminant}%
\Index
   {\sups row of matrix}%
   {r row}%
\Index
   {$R$\Hyph module}%
   {R- module}%
\Index
   {$r$\hyph row of matrix}%
   {r-row}%
\Index
   {rank of Hermitian matrix by principal minors}%
   {rank of Hermitian matrix by principal minors}%
\Index
   {rank of quadratic map of division ring}%
   {rank of quadratic map, division ring}%
\Index
   {\RC inverse element of biring}%
   {rc-inverse element}%
\Index
   {\RC major minor matrix}%
   {RC-major minor}%
\Index
   {\RC matrix group}%
   {rc-matrix group}%
\Index
   {\RC nonsingular matrix}%
   {RC nonsingular matrix}%
\Index
   {\RC power}%
   {rc power}%
\Index
   {$\RCcirc$\Hyph product of matrices of maps}%
   {rc product of matrices of maps}%
\Index
   {\RC quasideterminant}%
   {RC-quasideterminant}%
\Index
   {\RC rank of matrix}%
   {rc-rank of matrix}%
\Index
   {\RC singular matrix}%
   {RC singular matrix}%
\Index
   {$\RCcirc$\Hyph nonsingular matrix}%
   {RCcirc nonsingular matrix}%
\Index
   {$\RCcirc$\Hyph nonsingular system of additive equations}%
   {RCcirc nonsingular system of additive equations}%
\Index
   {$\RCcirc$\Hyph quasideterminant}%
   {RCcirc-quasideterminant definition}%
\Index
   {$\RCcirc$\Hyph singular matrix}%
   {RCcirc singular matrix}%
\Index
   {\rcd affine plane}%
   {rcd affine plane}%
\Index
   {\rcd affine space}%
   {rcd affine space}%
\Index
   {\rcd vector}%
   {rcd vector}%
\Index
   {\rcd vector space}%
   {rcd vector space}%
\Index
   {reduced Cartesian product of bundles}%
   {reduced Cartesian product of bundles}%
\Index
   {reduced Cartesian product of total spaces}%
   {reduced Cartesian product of total spaces}%
\Index
   {reduced fibered correspondence from $\Bundle{A}$ to $\Bundle B$}%
   {reduced fibered correspondence from A to B}%
\Index
   {reduced fibered correspondence in $\Bundle{A}$}%
   {reduced fibered correspondence in A}%
\Index
   {reduced morphism of representations}%
   {reduced morphism of representations}%
\Index
   {reduced polymorphism of representations}%
   {reduced polymorphism of representations}%
\Index
   {reduced quadratic equation}%
   {reduced quadratic equation}%
\Index
   {reducible biring}%
   {reducible biring}%
\Index
   {reference frame in event space}%
   {reference frame in event space}%
\Index
   {reference frame manifold}%
   {reference frame manifold}%
\Index
   {reflexive $2$\Hyph ary fibered relation}%
   {reflexive 2 ary fibered relation}%
\Index
   {regular endomorphism of representation}%
   {regular endomorphism of representation}%
\Index
   {regular endomorphism of tower of representations}%
   {regular endomorphism of tower of representations}%
\Index
   {regular quadratic map in division ring}%
   {regular quadratic map, division ring}%
\Index
   {relatively prime $A$\Hyph numbers}%
   {relatively prime A numbers}%
\Index
   {remainder of the division}%
   {remainder of the division}%
\Index
   {representation conjugated to representation}%
   {representation conjugated to representation}%
\Index
   {\Ts{A}representation in $\Omega_2$\Hyph algebra}%
   {A* representation of algebra}%
\Index
   {representation of group}%
   {representation of group}%
\Index
   {representation of $\Omega$\Hyph algebra in representation}%
   {representation of Omega algebra in representation}%
\Index
   {representation of $\Omega$\Hyph algebra in tower of representations}%
   {representation of Omega algebra in tower of representations}%
\Index
   {representation of $\Omega$\Hyph algebra $A$ in category $\mathcal B$}%
   {representation of Omega algebra in category}%
\Index
   {\sT{A}representation of $\Omega_1$\Hyph algebra $A$ in $\Omega_2$\Hyph algebra}%
   {*A representation of algebra}%
\Index
   {representation of $\Omega_1$\Hyph algebra $A$ in $\Omega_2$\Hyph algebra $M$}%
   {representation of algebra}%
\Index
   {representative of geometric object in \drc vector space}%
   {representative of geometric object, drc vector space}%
\Index
   {representative of geometric object in $\Omega_2$\Hyph algebra}%
   {representative of geometric object, representation g}%
\Index
   {representative of geometric object in \rcd vector space}%
   {representative of geometric object, rcd vector space}%
\Index
   {representative of geometric object in tuple of $\VX\Omega$\Hyph algebras}%
   {representative of geometric object, tower of representations g}%
\Index
   {representative of geometric object in vector space}%
   {representative of geometric object, vector space}%
\Index
   {restriction of correspondence $\Phi$ to set $C$}%
   {restriction of correspondence}%
\Index
   {right $A$\Hyph vector space}%
   {right A vector space}%
\Index
   {right cofactor of entry of matrix}%
   {right cofactor, matrix}%
\Index
   {right $D$\Hyph module}%
   {right D module}%
\Index
   {right $D$\hyph vector space of columns}%
   {right vector space of columns}%
\Index
   {right $D$\hyph vector space of rows}%
   {right vector space of rows}%
\Index
   {right defined Lie algebra of Lie group}%
   {right defined Lie algebra}%
\Index
   {right double cofactor of entry of matrix}%
   {right double cofactor}%
\Index
   {right fraction}%
   {right fraction}%
\Index
   {right ideal of algebra}%
   {right ideal of algebra}%
\Index
   {right invariant vector field}%
   {right invariant vector}%
\Index
   {right module}%
   {right module}%
\Index
   {right module over $D$\Hyph algebra $A$}%
   {right module over algebra}%
\Index
   {right principal ideal}%
   {right principal ideal}%
\Index
   {right shift on group}%
   {right shift}%
\Index
   {right shift on group}%
   {right shift, group}%
\Index
   {right structural constant of Lie algebra}%
   {right structural constant of Lie algebra}%
\Index
   {right vector space}%
   {right vector space}%
\Index
   {right zero divisor}%
   {right zero divisor}%
\Index
   {right-ordered cycle notation of permutation}%
   {right-ordered cycle notation of permutation}%
\Index
   {right\Hyph side $A_1$\Hyph representation}%
   {right-side A representation}%
\Index
   {right\Hyph side product}%
   {right-side product}%
\Index
   {right\Hyph side product of vector over scalar}%
   {right-side product of vector over scalar}%
\Index
   {right-side representation}%
   {right-side representation}%
\Index
   {right-side representation of fibered $\Omega$\Hyph algebra}%
   {right-side representation of fibered Omega-algebra}%
\Index
   {right-side representation of $\Omega_1$\Hyph algebra $A$ in $\Omega_2$\Hyph algebra $M$}%
   {right-side representation of algebra}%
\Index
   {right-side transformation}%
   {right-side transformation}%
\Index
   {ring has characteristic $0$}%
   {ring has characteristic 0}%
\Index
   {ring has characteristic $p$}%
   {ring has characteristic p}%
\Index
   {ring of sets}%
   {ring of sets}%
\Index
   {ring of sets generated by semiring of sets}%
   {ring of sets generated by semiring}%
\Index
   {ring with conjugation}%
   {ring with conjugation}%
\Index
   {root of polynomial}%
   {root of polynomial}%
\Index
   {row $*D$\Hyph vector}%
   {row *D vector}%
\Index
   {row $D*$\Hyph vector}%
   {row D* vector}%
\Index
   {row determinant}%
   {row determinant}%
\Index
   {row vector}%
   {row vector}%
\SetIndexSpace%
\Index
   {$\star A$\Hyph module}%
   {starA-module}%
\Index
   {scalar algebra of algebra}%
   {scalar algebra of algebra}%
\Index
   {scalar algebra of ring}%
   {scalar algebra of ring}%
\Index
   {scalar of element of algebra}%
   {scalar of algebra}%
\Index
   {scalar of element of ring}%
   {scalar of ring}%
\Index
   {scalar potential}%
   {scalar potential}%
\Index
   {Schauder basis}%
   {Schauder basis}%
\Index
   {second axiom of countability}%
   {second axiom of countability}%
\Index
   {second Newton law}%
   {Second Newton law}%
\Index
   {section of bundle}%
   {section of bundle}%
\Index
   {semigroup}%
   {semigroup}%
\Index
   {semiring of sets}%
   {semiring of sets}%
\Index
   {sequence converges}%
   {sequence converges}%
\Index
   {sequence converges almost everywhere}%
   {converges almost everywhere}%
\Index
   {sequence converges uniformly}%
   {sequence converges uniformly}%
\Index
   {series converges normally}%
   {series converges normally}%
\Index
   {set admits operation}%
   {set admits operation}%
\Index
   {set is closed with respect to operation}%
   {set is closed with respect to operation}%
\Index
   {set is dense in set}%
   {dense in set}%
\Index
   {set of coordinates of representation}%
   {coordinate set of representation}%
\Index
   {set of invertible elements of algebra}%
   {set of invertible elements of algebra}%
\Index
   {set of $\Omega_2$\Hyph words of representation}%
   {word set of representation}%
\Index
   {set of tuples of coordinates of tower of representations}%
   {coordinate set of tower of representations}%
\Index
   {set of tuples of $\Vector\Omega$\Hyph words of tower of representations}%
   {word set of tower of representations}%
\Index
   {set of zeros of algebra}%
   {set of zeros of algebra}%
\Index
   {simple map}%
   {simple map}%
\Index
   {simple polyvector}%
   {simple polyvector}%
\Index
   {simplex}%
   {simplex}%
\Index
   {sine}%
   {sine}%
\Index
   {single transitive representation of fibered $\Omega$\Hyph algebra}%
   {single transitive representation of fibered Omega-algebra}%
\Index
   {single transitive representation of group}%
   {single transitive representation of group}%
\Index
   {single transitive representation of $\Omega$\Hyph algebra $A$}%
   {single transitive representation of algebra}%
\Index
   {singular linear map}%
   {singular linear map}%
\Index
   {skew product of vectors}%
   {skew product of vectors}%
\Index
   {skew symmetric polylinear map}%
   {skew symmetric polylinear map}%
\Index
   {space of orbits of \Ts{G}representation}%
   {space of orbits of G* representation}%
\Index
   {space of orbits of left\Hyph side representation}%
   {space of orbits of left side representation}%
\Index
   {spacelike vector}%
   {spacelike vector}%
\Index
   {speed of deviation}%
   {speed of deviation}%
\Index
   {spherical coordinates}%
   {spherical coordinates}%
\Index
   {square root}%
   {square root}%
\Index
   {$(\mathcal S\RCstar,\mathcal T\RCstar)$\Hyph linear map of vector bundles}%
   {src trc linear map of vector bundles}%
\Index
   {($S\star$, $\star T$)\hyph bimodule}%
   {(Sstar,starT)-bimodule}%
\Index
   {stability group}%
   {stability group}%
\Index
   {stable set of representation}%
   {stable set of representation}%
\Index
   {standard component of derivative}%
   {standard component of derivative}%
\Index
   {standard component of the G\^ateaux derivative}%
   {standard component of Gateaux derivative}%
\Index
   {standard component of linear map}%
   {standard component of linear map}%
\Index
   {standard component of polylinear map}%
   {standard component of polylinear map}%
\Index
   {standard component of tensor}%
   {standard component of tensor}%
\Index
   {standard component over field $F$ of bilitnear map $f$}%
   {standard component of bilinear map, division ring}%
\Index
   {standard coordinates of basis}%
   {standard coordinates of basis}%
\Index
   {standard coordinates of basis}%
   {standard coordinates of basis}%
\Index
   {standard representation of the derivative}%
   {derivative, standard representation}%
\Index
   {standard representation of the G\^ateaux derivative}%
   {Gateaux derivative, standard representation}%
\Index
   {standard representation of linear map}%
   {linear map, standard representation}%
\Index
   {standard representation of matrix}%
   {Standard representation}%
\Index
   {standard representation of polylinear map}%
   {polylinear map, standard representation}%
\Index
   {standard representation of quadratic map of division ring over field $F$}%
   {quadratic map, standard representation, division ring}%
\Index
   {standard representation over field $F$ of bilinear map of division ring}%
   {bilinear map, standard representation, division ring}%
\Index
   {$\star R$\hyph module}%
   {starR-module}%
\Index
   {$\star D$\hyph product of vector over scalar}%
   {starD product of vector over scalar, vector space}%
\Index
   {starlike set}%
   {starlike set}%
\Index
   {\sT representation of fibered group}%
   {starT representation of fibered group}%
\Index
   {\sT representation of fibered group}%
   {starT representation of fibered group}%
\Index
   {\sT representation of fibered $\Omega$\Hyph algebra}%
   {starT representation of fibered Omega-algebra}%
\Index
   {\sT shift on fibered group}%
   {starT shift, fibered group}%
\Index
   {\sT transformation on bundle}%
   {starT transformation of bundle}%
\Index
   {structural constants}%
   {structural constants}%
\Index
   {subalgebra of $\Omega$\Hyph algebra}%
   {subalgebra of Omega-algebra}%
\Index
   {subbundle}%
   {subbundle}%
\Index
   {subbundle of $\mathcal D\star$\hyph vector space}%
   {subbundle of Dstar vector bundle}%
\Index
   {subgroup of $\Omega$\Hyph group}%
   {subgroup of Omega group}%
\Index
   {submodule}%
   {submodule}%
\Index
   {submodule generated by set}%
   {submodule generated by set}%
\Index
   {subrepresentation generated by set $X$}%
   {subrepresentation generated by set}%
\Index
   {subrepresentation of representation}%
   {subrepresentation of representation}%
\Index
   {sum of geometric objects in vector space}%
   {sum of geometric objects, vector space}%
\Index
   {sum of geometric objects}%
   {sum of geometric objects}%
\Index
   {sum of maps}%
   {sum of maps}%
\Index
   {sum of polynomials}%
   {sum of polynomials}%
\Index
   {superposition of coordinates of the representation $f$ and the element $m$}%
   {superposition of coordinates, representation}%
\Index
   {superposition of coordinates of the tower of representations $\Vector f$ and the element $\VX a$}%
   {superposition of coordinates, tower of representations}%
\Index
   {symmetric $2$\Hyph ary fibered relation}%
   {symmetric 2 ary fibered relation}%
\Index
   {symmetric bilinear map of $D$\Hyph vector space to division ring}%
   {symmetric bilinear map, vector space to division ring}%
\Index
   {symmetric polylinear map}%
   {symmetric polylinear map}%
\Index
   {symmetric polylinear mapping into associative algebra}%
   {polylinear map symmetric, associative algebra}%
\Index
   {symmetrization of polylinear map}%
   {symmetrization of polylinear map}%
\Index
   {symmetry group}%
   {symmetry group}%
\Index
   {symmetry group}%
   {SymmetryGroup}%
\Index
   {synchronization of reference frame}%
   {synchronization of reference frame}%
\Index
   {system of additive equations}%
   {system of additive equations}%
\Index
   {system of \drc linear equations}%
   {system of drc linear equations}%
\Index
   {system of linear equations}%
   {system of linear equations}%
\Index
   {system of \rcd linear equations}%
   {system of rcd linear equations}%
\SetIndexSpace%
\Index
   {$T_1$\Hyph space}%
   {T1 space}%
\Index
   {Taylor polynomial}%
   {Taylor polynomial, division ring}%
\Index
   {Taylor series}%
   {Taylor series, division ring}%
\Index
   {tensor inverse to tensor}%
   {inverse tensor}%
\Index
   {tensor power}%
   {tensor power}%
\Index
   {tensor product}%
   {tensor product}%
\Index
   {the Fr\'echet \Ds derivative of map $f$ of division ring $D$ at point $x$}%
   {Frechet Dstar derivative of map, division ring}%
\Index
   {timelike vector}%
   {timelike vector}%
\Index
   {topological $D$\Hyph vector space}%
   {topological D vector space}%
\Index
   {topological $D$\Hyph algebra}%
   {topological D algebra}%
\Index
   {topological division ring}%
   {topological division ring}%
\Index
   {topological ring}%
   {topological ring}%
\Index
   {torsion form}%
   {torsion form}%
\Index
   {torsion tensor}%
   {torsion tensor}%
\Index
   {tower of bundles}%
   {tower of bundles}%
\Index
   {tower of effective representations}%
   {tower of effective representations}%
\Index
   {tower of representations of $\Vector{\Omega}$\Hyph algebras}%
   {tower of representations of algebras}%
\Index
   {tower of subrepresentations}%
   {tower of subrepresentations}%
\Index
   {tower of subrepresentations of tower of representations $\Vector f$ generated by tuple of sets $\VX X$}%
   {subrepresentation generated by tuple of sets}%
\Index
   {trace of quaternion}%
   {trace, quaternion algebra}%
\Index
   {transformation coordinated with equivalence}%
   {transformation coordinated with equivalence}%
\Index
   {transformation of universal algebra}%
   {transformation of universal algebra}%
\Index
   {transformation on bundle}%
   {transformation of bundle}%
\Index
   {transitive $2$\Hyph ary fibered relation}%
   {transitive 2 ary fibered relation}%
\Index
   {transitive representation of fibered $\Omega$\Hyph algebra}%
   {transitive representation of fibered Omega-algebra}%
\Index
   {transitive representation of group}%
   {transitive representation of group}%
\Index
   {transitive representation of $\Omega$\Hyph algebra $A$}%
   {transitive representation of algebra}%
\Index
   {\Ts representation of fibered group}%
   {Tstar representation of fibered group}%
\Index
   {\Ts representation of fibered $\Omega$\Hyph algebra}%
   {Tstar representation of fibered Omega-algebra}%
\Index
   {tuple of coordinates of element $\Vector a$ relative to tuple of sets $\VX X$}%
   {coordinates of element, tower of representations}%
\Index
   {tuple of equivalence generated by tower of representations $\Vector f$}%
   {tuple of equivalence of tower of representations}%
\Index
   {tuple of generating sets of tower of representations}%
   {tuple of generating sets of tower of representations}%
\Index
   {tuple of generating sets of tower subrepresentations}%
   {tuple of generating sets of tower subrepresentations}%
\Index
   {tuple of $\Vector{\Omega}$\Hyph words of element of tower of representations relative to tuple of generating sets}%
   {tuple of words relative to tuple of generating sets, tower of representations}%
\Index
   {tuple of stable sets of tower of representation}%
   {tuple of stable sets of tower of representations}%
\Index
   {twin representations}%
   {twin representations}%
\Index
   {twin representations of division ring}%
   {twin representations of division ring}%
\Index
   {twin representations of fibered group}%
   {twin representations of fibered group}%
\Index
   {twin representations of group}%
   {twin representations of group}%
\SetIndexSpace%
\Index
   {unit interval}%
   {unit interval}%
\Index
   {unit of ring of sets}%
   {unit of ring of sets}%
\Index
   {unit sphere in $D$\Hyph algebra}%
   {unit sphere in algebra}%
\Index
   {unit sphere in division ring}%
   {unit sphere in division ring}%
\Index
   {unit vector}%
   {unit vector}%
\Index
   {unital extension}%
   {unital extension}%
\Index
   {unital ring}%
   {unital ring}%
\Index
   {unitarity law}%
   {unitarity law}%
\Index
   {universal algebra}%
   {universal algebra}%
\Index
   {universally attracting object of category}%
   {universally attracting}%
\Index
   {universally repelling  object of category}%
   {universally repelling}%
\SetIndexSpace%
\Index
   {basis for vector  bundle}%
   {basis, vector bundle}%
\Index
   {valued division ring}%
   {valued division ring}%
\Index
   {vector}%
   {vector}%
\Index
   {vector bundle}%
   {vector bundle}%
\Index
   {vector module of algebra}%
   {vector module of algebra}%
\Index
   {vector module of ring}%
   {vector module of ring}%
\Index
   {vector of element of algebra}%
   {vector of algebra}%
\Index
   {vector of element of ring}%
   {vector of ring}%
\Index
   {vector potential}%
   {vector potential}%
\Index
   {vector space}%
   {vector space}%
\Index
   {vector space type}%
   {vector space type}%
\Index
   {vertical component of vector}%
   {vertical component of vector}%
\Index
   {vertical subspace}%
   {vertical subspace}%
\Index
   {vertical vector}%
   {vertical vector}%
\SetIndexSpace%
\Index
   {zero divisor}%
   {zero divisor}%
\SetIndexSpace%
\Index
   {$\mu$\Hyph measurable map}%
   {mu measurable map}%
\SetIndexSpace%
\Index
   {$\Omega$\Hyph group}%
   {Omega group}%
\Index
   {$\Omega$\Hyph linear mapping}%
   {Omega linear map}%
\Index
   {\(\Omega\)\Hyph ring}%
   {Omega ring}%
\Index
   {$\Omega$\Hyph algebra}%
   {Omega-algebra}%
\Index
   {$\Omega_2$\Hyph word of element of representation relative to generating set}%
   {word of element relative to generating set, representation}%
\SetIndexSpace%
\Index
   {$\sigma$\Hyph algebra of sets}%
   {sigma algebra of sets}%
\Index
   {$\sigma$\Hyph ring of sets}%
   {sigma ring of sets}%
\Index
   {\(\sigma\)\Hyph additive measure}%
   {sigma-additive measure}%

\CloseIndex

%% file: Symbol.English.tex
\def\indexname{Special Symbols and Notations}
\OpenIndex

\SetIndexSpace
\Symb%
   {direct sum}%
   {direct sum}%
   {0}{0}%
\Symb%
   {unit interval}%
   {unit interval}%
   {0}{0}%

\SetIndexSpace
\Symb%
   {minor matrix}%
   {A from b a}%
   {A}{0}%
\Symb%
   {minor matrix}%
   {A from columns T}%
   {A}{0}%
\Symb%
   {minor matrix}%
   {A from rows S}%
   {A}{0}%
\Symb%
   {set of vectors whose expansion relative to the basis $\Basis e$ converges normally}%
   {A plus Schauder}%
   {A}{0}%
\Symb%
   {minor matrix}%
   {A without column a}%
   {A}{0}%
\Symb%
   {minor matrix}%
   {A without columns T}%
   {A}{0}%
\Symb%
   {minor matrix}%
   {A without row b}%
   {A}{0}%
\Symb%
   {minor matrix}%
   {A without rows S}%
   {A}{0}%
\Symb%
   {active representation of group $G(f)$ in basis manifold $\mathcal B(f)$}%
   {active representation in basis manifold}%
   {A}{0}%
\Symb%
   {active representation of group $G(\Vector f)$ in basis manifold $\mathcal B(\Vector f)$}%
   {active representation in basis manifold, tower of representations}%
   {A}{0}%
\Symb%
   {$A$\Hyph algebra of polynomials over $D$\Hyph algebra $A$}%
   {algebra of polynomials over algebra}%
   {A}{0}%
\Symb%
   {algebra of polynomials over $D$\Hyph algebra $A$}%
   {algebra of polynomials over D algebra}%
   {A}{0}%
\Symb%
   {algebra of rational mappings of algebra $A$}%
   {algebra of rational mappings of algebra}%
   {A}{0}%
\Symb%
   {affine space}%
   {An}%
   {A}{0}%
\Symb%
   {associator of $D$\Hyph algebra}%
   {associator of algebra}%
   {A}{0}%
\Symb%
   {\subs row ($c$\hyph row) of matrix}%
   {c row}%
   {A}{0}%
\Symb%
   {category of left-side representations of $\Omega_1$\Hyph algebra $A$}%
   {category of left-side representations of Omega1 algebra}%
   {A}{0}%
\Symb%
   {category of representations}%
   {category of representations}%
   {A}{0}%
\Symb%
   {commutator of $D$\Hyph algebra}%
   {commutator of algebra}%
   {A}{0}%
\Symb%
   {component of linear map}%
   {component of linear map, vector}%
   {A}{0}%
\Symb%
   {component $p$ of polylinear mapping $\Vector A$}%
   {component of polyadditive map, D vector space}%
   {A}{0}%
\Symb%
   {component of polylinear map}%
   {component of polylinear map, vector}%
   {A}{0}%
\Symb%
   {conjugated $D$\Hyph  module}%
   {conjugated D module}%
   {A}{0}%
\Symb%
   {\CR power of element $A$ of biring}%
   {cr power}%
   {A}{0}%
\Symb%
   {\crd vector}%
   {crd vector}%
   {A}{0}%
\Symb%
   {\CR inverse element of biring}%
   {cr-inverse element}%
   {A}{0}%
\Symb%
   {\CR product}%
   {cr-product}%
   {A}{0}%
\Symb%
   {\dcr vector}%
   {dcr vector}%
   {A}{0}%
\Symb%
   {derivative of left shift}%
   {derivative of left shift}%
   {A}{0}%
\Symb%
   {derivative of left shift in $1$\Hyph parameter Lie group}%
   {derivative of left shift, 1-Parameter Group}%
   {A}{0}%
\Symb%
   {derivative of left shift in $1$\Hyph parameter Lie D group}%
   {derivative of left shift, 1-Parameter Group, algebra}%
   {A}{0}%
\Symb%
   {derivative of right shift}%
   {derivative of right shift}%
   {A}{0}%
\Symb%
   {derivative of right shift in $1$\Hyph parameter Lie group}%
   {derivative of right shift, 1-Parameter Group}%
   {A}{0}%
\Symb%
   {derivative of right shift in $1$\Hyph parameter Lie D group}%
   {derivative of right shift, 1-Parameter Group, algebra}%
   {A}{0}%
\Symb%
   {derivative of left shift}%
   {derivative of Tstar shift}%
   {A}{0}%
\Symb%
   {\drc vector}%
   {drc vector}%
   {A}{0}%
\Symb%
   {coordinates of vector $a$ relative to Hamel basis}%
   {Hamel basis, coordinates}%
   {A}{0}%
\Symb%
   {hermitian conjugation in division ring}%
   {hermitian conjugation, division ring}%
   {A}{0}%
\Symb%
   {tensor inverse to tensor $a$}%
   {inverse tensor}%
   {A}{0}%
\Symb%
   {isomorphic}%
   {isomorphic}%
   {A}{0}%
\Symb%
   {$(^j_i)$\hyph\CR quasideterminant}%
   {j i CR quasideterminant definition}%
   {A}{0}%
\Symb%
   {$(ji)$\hyph quasideterminant of matrix $\bfA$}%
   {j i quasideterminant definition}%
   {A}{0}%
\Symb%
   {$(^j_i)$\hyph $\RCcirc$\Hyph quasideterminant}%
   {j i RCcirc-quasideterminant definition}%
   {A}{0}%
\Symb%
   {$(^j_i)$\hyph \RC quasideterminant}%
   {j i RC-quasideterminant definition}%
   {A}{0}%
\Symb%
   {left fraction}%
   {left fraction}%
   {A}{0}%
\Symb%
   {left principal ideal}%
   {left principal ideal}%
   {A}{0}%
\Symb%
   {left shift in $D$\Hyph algebra}%
   {left shift, D algebra}%
   {A}{0}%
\Symb%
   {linear combination}%
   {linear combination}%
   {A}{0}%
\Symb%
   {little group}%
   {little group}%
   {A}{0}%
\Symb%
   {transformation of matrix}%
   {matrix, replacing its column}%
   {A}{0}%
\Symb%
   {transformation of matrix}%
   {matrix, replacing its row}%
   {A}{0}%
\Symb%
   {$A$\Hyph module of homogeneous polynomials over $D$\Hyph algebra $A$}%
   {module of homogeneous polynomials over algebra}%
   {A}{0}%
\Symb%
   {norm on $D$\Hyph module}%
   {norm on D module}%
   {A}{0}%
\Symb%
   {$\Omega$\Hyph algebra}%
   {Omega-algebra}%
   {A}{0}%
\Symb%
   {opposite algebra to algebra $A$}%
   {opposite algebra}%
   {A}{0}%
\Symb%
   {orbit of effective left\Hyph side representation}%
   {orbit of effective left-side representation}%
   {A}{0}%
\Symb%
   {orbit of effective right\Hyph side representation}%
   {orbit of effective right-side representation}%
   {A}{0}%
\Symb%
   {orbit of linear map}%
   {orbit of linear map}%
   {A}{0}%
\Symb%
   {derivative}%
   {overline nabla_l, definition 2}%
   {A}{0}%
\Symb%
   {partial linear map}%
   {partial linear map}%
   {A}{0}%
\Symb%
   {principal ideal}%
   {principal ideal}%
   {A}{0}%
\Symb%
   {quasideterminant of matrix $\bfA$}%
   {quasideterminant definition}%
   {A}{0}%
\Symb%
   {\sups row ($r$\hyph row) of matrix}%
   {r row}%
   {A}{0}%
\Symb%
   {\RC power of element $A$ of biring}%
   {rc power}%
   {A}{0}%
\Symb%
   {$\RCcirc$\Hyph quasideterminant}%
   {RCcirc-quasideterminant definition}%
   {A}{0}%
\Symb%
   {\rcd vector}%
   {rcd vector}%
   {A}{0}%
\Symb%
   {\RC inverse element of biring}%
   {rc-inverse element}%
   {A}{0}%
\Symb%
   {\RC product}%
   {rc-product}%
   {A}{0}%
\Symb%
   {\RC quasideterminant}%
   {RC-quasideterminant definition}%
   {A}{0}%
\Symb%
   {right principal ideal}%
   {right principal ideal}%
   {A}{0}%
\Symb%
   {right shift in $D$\Hyph algebra}%
   {right shift, D algebra}%
   {A}{0}%
\Symb%
   {coordinates of vector $a$ relative to Schauder basis}%
   {Schauder basis, coordinates}%
   {A}{0}%
\Symb%
   {set of additive maps}%
   {set additive maps}%
   {A}{0}%
\Symb%
   {set of invertible elements of algebra $A$}%
   {set of invertible elements of algebra}%
   {A}{0}%
\Symb%
   {set of zeros of algebra $A$}%
   {set of zeros of algebra}%
   {A}{0}%
\Symb%
   {set of polylinear maps of rings $R_1$, ..., $R_n$ into module $S$}%
   {set polylinear maps, ring}%
   {A}{0}%
\Symb%
   {simplex}%
   {simplex}%
   {A}{0}%
\Symb%
   {skew product of vectors $\Vector a_1$, ..., $\Vector a_m$}%
   {skew product of vectors}%
   {A}{0}%
\Symb%
   {space of orbits of effective left\Hyph side representation}%
   {space of orbits of effective left-side representation}%
   {A}{0}%
\Symb%
   {space of orbits of effective right\Hyph side representation}%
   {space of orbits of effective right-side representation}%
   {A}{0}%
\Symb%
   {space of orbits of left\Hyph side representation}%
   {space of orbits of left side representation}%
   {A}{0}%
\Symb%
   {square root}%
   {square root}%
   {A}{0}%
\Symb%
   {stability group}%
   {stability group}%
   {A}{0}%
\Symb%
   {\sT shift}%
   {starT shift, fibered group}%
   {A}{0}%
\Symb%
   {tensor power of algebra $A$}%
   {tensor power of algebra}%
   {A}{0}%
\Symb%
   {anholonomic coordinates of vector}%
   {vector anholonomic coordinates}%
   {A}{0}%
\Symb%
   {holonomic coordinates of vector}%
   {vector holonomic coordinates}%
   {A}{0}%

\SetIndexSpace
\Symb%
   {basis manifold of \rcd vector space $\Vector V$}%
   {basis manifold of rcd vector space}%
   {B}{0}%
\Symb%
   {basis manifold of vector space}%
   {basis manifold of vector space}%
   {B}{0}%
\Symb%
   {basis manifold of representation $f$}%
   {basis manifold representation F algebra}%
   {B}{0}%
\Symb%
   {basis manifold of tower of representations $\Vector f$}%
   {basis manifold tower of representations}%
   {B}{0}%
\Symb%
   {basis manifold of affine space}%
   {Basis Manifold, Affine Space}%
   {B}{0}%
\Symb%
   {basis manifold of central affine space}%
   {BCAn}%
   {B}{0}%
\Symb%
   {basis manifold of Euclid space}%
   {BEn}%
   {B}{0}%
\Symb%
   {Borel algebra}%
   {Borel algebra}%
   {B}{0}%
\Symb%
   {canonical remainder of the division}%
   {canonical remainder of the division}%
   {B}{0}%
\Symb%
   {Cartesian power}%
   {Cartesian power}%
   {B}{0}%
\Symb%
   {Cartesian power $\Bundle A$ of bundle $\Bundle B$}%
   {Cartesian power A of bundle B}%
   {B}{0}%
\Symb%
   {Cartesian power $A$ of set $B$}%
   {Cartesian power of set}%
   {B}{0}%
\Symb%
   {closed ball}%
   {closed ball}%
   {B}{0}%
\Symb%
   {closure of set}%
   {closure of set}%
   {B}{0}%
\Symb%
   {coproduct in category}%
   {coproduct in category}%
   {B}{0}%
\Symb%
   {basis manifold of central affine space}%
   {FCAn}%
   {B}{0}%
\Symb%
   {basis manifold of Euclid space}%
   {FEn}%
   {B}{0}%
\Symb%
   {lattice of subrepresentations of representation $f$}%
   {lattice of subrepresentations}%
   {B}{0}%
\Symb%
   {lattice of towers of subrepresentations of tower of representations $\Vector f$}%
   {lattice of subrepresentations, tower of representations}%
   {B}{0}%
\Symb%
   {open ball}%
   {open ball}%
   {B}{0}%
\Symb%
   {product in category}%
   {product in category}%
   {B}{0}%
\Symb%
   {right fraction}%
   {right fraction}%
   {B}{0}%
\Symb%
   {tensor power of representation}%
   {tensor power of representation}%
   {B}{0}%

\SetIndexSpace
\Symb%
   {$\sigma$\Hyph algebra of sets measurable with respect to measure $\mu$}%
   {algebra of sets measurable with respect to measure}%
   {C}{0}%
\Symb%
   {central affine space}%
   {CAn}%
   {C}{0}%
\Symb%
   {central affine space}%
   {central affine space}%
   {C}{0}%
\Symb%
   {continuity class}%
   {class Cn}%
   {C}{0}%
\Symb%
   {$j$th column determinant of matrix $\bfA$}%
   {column determinant}%
   {C}{0}%
\Symb%
   {cosine}%
   {cosine}%
   {C}{0}%
\Symb%
   {$\CRcirc$\Hyph product of matrices of maps}%
   {cr product of matrices of maps}%
   {C}{0}%
\Symb%
   {hyperbolic cosine}%
   {hyperbolic cosine}%
   {C}{0}%
\Symb%
   {left structural constant of Lie algebra}%
   {left structural constant of Lie algebra}%
   {C}{0}%
\Symb%
   {right structural constant of Lie algebra}%
   {right structural constant of Lie algebra}%
   {C}{0}%
\Symb%
   {set of continuous multivariable maps}%
   {set continuous multivariable maps}%
   {C}{0}%
\Symb%
   {structural constants}%
   {structural constants}%
   {C}{0}%

\SetIndexSpace
\Symb%
   {basis vector of representation of Lie group over algebra $A$}%
   {basis vector of representation of Lie group over algebra A}%
   {D}{0}%
\Symb%
   {coordinates of basis vector of representation of Lie group over algebra $A$}%
   {basis vector of representation of Lie group over algebra A, coordinates}%
   {D}{0}%
\Symb%
   {component of derivative of map $f(x)$}%
   {component of derivative}%
   {D}{0}%
\Symb%
   {component of derivative of second order of map $f(x)$}%
   {component of derivative of Second Order}%
   {D}{0}%
\Symb%
   {component of the G\^ateaux derivative of map $f(x)$}%
   {component of Gateaux derivative}%
   {D}{0}%
\Symb%
   {component of the G\^ateaux derivative of map $f(x)$}%
   {component of Gateaux derivative of map, D vector space, short}%
   {D}{0}%
\Symb%
   {component of the G\^ateaux derivative of second order of map $f(x)$}%
   {component of Gateaux derivative of Second Order}%
   {D}{0}%
\Symb%
   {component of the G\^ateaux derivative of second order of map $f(x)$}%
   {component of Gateaux derivative of Second Order, D vector space}%
   {D}{0}%
\Symb%
   {component of the G\^ateaux derivative of map $f(x)$}%
   {component of Gateaux derivative, vector space}%
   {D}{0}%
\Symb%
   {conjugation in algebra}%
   {conjugation in algebra}%
   {D}{0}%
\Symb%
   {conjugation in ring}%
   {conjugation in ring}%
   {D}{0}%
\Symb%
   {coordinate \rcd vector space}%
   {coordinate rcd vector space}%
   {D}{0}%
\Symb%
   {coordinate reference frame}%
   {coordinate reference frame, extensive definition}%
   {D}{0}%
\Symb%
   {coordinate vector bundle}%
   {coordinate vector bundle}%
   {D}{0}%
\Symb%
   {derivative of map $f$}%
   {derivative of map}%
   {D}{0}%
\Symb%
   {derivative of map $f$}%
   {derivative of map inline}%
   {D}{0}%
\Symb%
   {derivative of order $n$}%
   {derivative of Order n}%
   {D}{0}%
\Symb%
   {derivative of order $n$}%
   {derivative of Order n inline}%
   {D}{0}%
\Symb%
   {derivative of second order}%
   {derivative of Second Order}%
   {D}{0}%
\Symb%
   {derivative of second order}%
   {derivative of Second Order inline}%
   {D}{0}%
\Symb%
   {diagonal in bundle $\Bundle A$}%
   {diagonal in bundle, 1}%
   {D}{0}%
\Symb%
   {differential of independent variable}%
   {differential of independent variable}%
   {D}{0}%
\Symb%
   {differential of map $f$}%
   {differential of map}%
   {D}{0}%
\Symb%
   {direct product of division rings $D_1$, ..., $D_n$}%
   {direct product of division rings, 1 n}%
   {D}{0}%
\Symb%
   {double determinant of matrix $\bfA$}%
   {double determinant}%
   {D}{0}%
\Symb%
   {exterior differential}%
   {exterior differential}%
   {D}{0}%
\Symb%
   {the Fr\'echet \Ds derivative of map $f$ of division ring}%
   {Frechet Dstar derivative of map, division ring}%
   {D}{0}%
\Symb%
   {the G\^ateaux \dcr derivative of map $f$ of $D$\Hyph vector space $V$ to $D$\Hyph vector space $W$}%
   {Gateaux dcr derivative of map, D vector space}%
   {D}{0}%
\Symb%
   {the G\^ateaux derivative of map $f$}%
   {Gateaux derivative of map}%
   {D}{0}%
\Symb%
   {the G\^ateaux derivative of map $f$}%
   {Gateaux derivative of map, fraction}%
   {D}{0}%
\Symb%
   {the G\^ateaux derivative of order $n$}%
   {Gateaux derivative of Order n}%
   {D}{0}%
\Symb%
   {the G\^ateaux derivative of order $n$ of map $f$ of division ring}%
   {Gateaux derivative of Order n, division ring}%
   {D}{0}%
\Symb%
   {the G\^ateaux derivative of order $n$ of map $f$ of algebra}%
   {Gateaux derivative of Order n, fraction, algebra}%
   {D}{0}%
\Symb%
   {the G\^ateaux derivative of order $n$ of map $f$ of division ring}%
   {Gateaux derivative of Order n, fraction, division ring}%
   {D}{0}%
\Symb%
   {the G\^ateaux derivative of second order}%
   {Gateaux derivative of Second Order}%
   {D}{0}%
\Symb%
   {the G\^ateaux derivative of second order of mapping $f$ of algebra}%
   {Gateaux derivative of Second Order, fraction, algebra}%
   {D}{0}%
\Symb%
   {the G\^ateaux derivative of second order of map $f$ of division ring}%
   {Gateaux derivative of Second Order, fraction, division ring}%
   {D}{0}%
\Symb%
   {the G\^ateaux differential of map $f$}%
   {Gateaux differential of map, vector}%
   {D}{0}%
\Symb%
   {the G\^ateaux \Ds derivative of map $f$ of division ring $D$}%
   {Gateaux Dstar derivative of map, division ring}%
   {D}{0}%
\Symb%
   {the G\^ateaux Jacobian of map of $D$\Hyph vector space}%
   {Gateaux Jacobian of map, D vector space}%
   {D}{0}%
\Symb%
   {the G\^ateaux partial \dcr derivative of map $f^{\gi b}$ with respect to variable $v^{\gi a}$}%
   {Gateaux partial dcr derivative of map, 1, D vector space}%
   {D}{0}%
\Symb%
   {the G\^ateaux partial \dcr derivative of map $f^{\gi b}$ with respect to variable $v^{\gi a}$}%
   {Gateaux partial dcr derivative of map, 2, D vector space}%
   {D}{0}%
\Symb%
   {the G\^ateaux partial \dcr derivative of map $f^{\gi b}$ with respect to variable $v^{\gi a}$}%
   {Gateaux partial dcr derivative of map, 3, D vector space}%
   {D}{0}%
\Symb%
   {the G\^ateaux partial derivative}%
   {Gateaux partial derivative}%
   {D}{0}%
\Symb%
   {the G\^ateaux mixed partial derivative}%
   {Gateaux partial derivative of Second Order}%
   {D}{0}%
\Symb%
   {the G\^ateaux partial \rcd derivative of map $f^{\gi b}$ with respect to variable $x^{\gi a}$}%
   {Gateaux partial rcd derivative of map, 1, D vector space}%
   {D}{0}%
\Symb%
   {the G\^ateaux partial \rcd derivative of map $f^{\gi b}$ with respect to variable $x^{\gi a}$}%
   {Gateaux partial rcd derivative of map, 2, D vector space}%
   {D}{0}%
\Symb%
   {the G\^ateaux partial \rcd derivative of map $f^{\gi b}$ with respect to variable $x^{\gi a}$}%
   {Gateaux partial rcd derivative of map, 3, D vector space}%
   {D}{0}%
\Symb%
   {the G\^ateaux \rcd derivative of map $f$ of $D$\hyph vector space $V$ to $D$\hyph vector space $W$}%
   {Gateaux rcd derivative of map, D vector space}%
   {D}{0}%
\Symb%
   {the G\^ateaux \sD derivative of map $f$ of division ring $D$}%
   {Gateaux starD derivative of map, division ring}%
   {D}{0}%
\Symb%
   {matrices vector space}%
   {matrices vector space}%
   {D}{0}%
\Symb%
   {Cartan derivative}%
   {overbrace D}%
   {D}{0}%
\Symb%
   {derivative}%
   {overline D}%
   {D}{0}%
\Symb%
   {partial derivative}%
   {partial derivative}%
   {D}{0}%
\Symb%
   {partial derivative of second order}%
   {partial derivative of second order}%
   {D}{0}%
\Symb%
   {derivative $e_{(k)}$}%
   {partial(k)}%
   {D}{0}%
\Symb%
   {product of map over scalar}%
   {product of map over scalar}%
   {D}{0}%
\Symb%
   {speed of deviation}%
   {speed of deviation}%
   {D}{0}%
\Symb%
   {standard component of derivative}%
   {standard component of derivative}%
   {D}{0}%
\Symb%
   {standard component of the G\^ateaux derivative}%
   {standard component of Gateaux derivative}%
   {D}{0}%
\Symb%
   {vector space type}%
   {vector space type}%
   {D}{0}%

\SetIndexSpace
\Symb%
   {Jacobian matrix of left shift}%
   {aE, quaternion, Jacobian matrix}%
   {E}{0}%
\Symb%
   {affine basis}%
   {Affine Basis}%
   {E}{0}%
\Symb%
   {basis of vector space}%
   {Basis e}%
   {E}{0}%
\Symb%
   {basis for module}%
   {basis for module}%
   {E}{0}%
\Symb%
   {basis in vector space $\Vector V$}%
   {basis in V}%
   {E}{0}%
\Symb%
   {basis of $D$\Hyph module $\mathcal L(D;A_1;A_2)$}%
   {basis L(A1,A2)}%
   {E}{0}%
\Symb%
   {basis for \crd vector space}%
   {basis, crd vector space}%
   {E}{0}%
\Symb%
   {basis for $D$\Hyph vector space}%
   {basis, D vector space}%
   {E}{0}%
\Symb%
   {basis for \dcr vector space}%
   {basis, dcr vector space}%
   {E}{0}%
\Symb%
   {basis for \drc vector space}%
   {basis, drc vector space}%
   {E}{0}%
\Symb%
   {basis for \rcd vector space}%
   {basis, rcd vector space}%
   {E}{0}%
\Symb%
   {basis for vector bundle}%
   {basis, vector bundle}%
   {E}{0}%
\Symb%
   {basis of $(n)$\hyph vector space}%
   {basis,n vector space}%
   {E}{0}%
\Symb%
   {Cartesian power of total spaces}%
   {Cartesian power of total spaces}%
   {E}{0}%
\Symb%
   {Cartesian product of total spaces}%
   {Cartesian product of total spaces, definition 1}%
   {E}{0}%
\Symb%
   {central affine basis}%
   {Central Affine Basis}%
   {E}{0}%
\Symb%
   {form of reference frame}%
   {dual forms, reference frame}%
   {E}{0}%
\Symb%
   {Euclid space}%
   {Euclid space}%
   {E}{0}%
\Symb%
   {Euclid space}%
   {Euclid space, division ring}%
   {E}{0}%
\Symb%
   {exponent}%
   {exponent}%
   {E}{0}%
\Symb%
   {Hamel basis}%
   {Hamel basis}%
   {E}{0}%
\Symb%
   {identical transformation of bundle}%
   {identical transformation of bundle}%
   {E}{0}%
\Symb%
   {linear automorphism of quaternioin algebra}%
   {mapping E, quaternion}%
   {E}{0}%
\Symb%
   {linear automorphism of quaternioin algebra}%
   {mapping E_1, quaternion}%
   {E}{0}%
\Symb%
   {linear automorphism of quaternioin algebra}%
   {mapping E_2, quaternion}%
   {E}{0}%
\Symb%
   {orthonornal basis}%
   {Orthonornal Basis}%
   {E}{0}%
\Symb%
   {pseudo Euclid space}%
   {pseudo Euclid space}%
   {E}{0}%
\Symb%
   {pseudo Euclid space}%
   {pseudo Euclid space, division ring}%
   {E}{0}%
\Symb%
   {quaternion algebra over the field $F$}%
   {quaternion algebra over the field}%
   {E}{0}%
\Symb%
   {quaternion division algebra over the field}%
   {quaternion division algebra over the fieldL}%
   {E}{0}%
\Symb%
   {reduced Cartesian product of total spaces}%
   {reduced Cartesian product of total spaces, definition 1}%
   {E}{0}%
\Symb%
   {Schauder basis}%
   {Schauder basis}%
   {E}{0}%
\Symb%
   {set of endomorphisms}%
   {set of endomorphisms}%
   {E}{0}%
\Symb%
   {set of nonsingular \sT transformations of bundle $\Bundle E$}%
   {set of starT nonsingular transformations of bundle}%
   {E}{0}%
\Symb%
   {set of transformations of universal algebra}%
   {set of transformations}%
   {E}{0}%
\Symb%
   {set of nonsingular \Ts transformations of bundle $\Bundle E$}%
   {set of Tstar nonsingular transformations of bundle}%
   {E}{0}%
\Symb%
   {standard coordinates of basis}%
   {standard coordinates of basis}%
   {E}{0}%
\Symb%
   {standard coordinates of reference frame}%
   {standard coordinates of reference frame}%
   {E}{0}%
\Symb%
   {vector field of reference frame}%
   {vector field of reference frame}%
   {E}{0}%
\Symb%
   {vector of basis}%
   {vector of basis}%
   {E}{0}%

\SetIndexSpace
\Symb%
   {alternation of polylinear map}%
   {alternation of polylinear map}%
   {F}{0}%
\Symb%
   {component of linear map $f$ of division ring}%
   {component of linear map, division ring}%
   {F}{0}%
\Symb%
   {component of polylinear map}%
   {component of polylinear map}%
   {F}{0}%
\Symb%
   {conjugation transformation}%
   {conjugation transformation}%
   {F}{0}%
\Symb%
   {exterior product}%
   {exterior product}%
   {F}{0}%
\Symb%
   {fibered morphism from bundle $\Bundle A$ into $\Bundle B$}%
   {fibered morphism from A into B}%
   {F}{0}%
\Symb%
   {filter $\mathfrak{F}$ converges to set $A$}%
   {filter converges}%
   {F}{0}%
\Symb%
   {homomorphism of fibered universal algebras}%
   {homomorphism of fibered universal algebras}%
   {F}{0}%
\Symb%
   {inverse fibered correspondence}%
   {inverse fibered correspondence, 1}%
   {F}{0}%
\Symb%
   {inverse reduced fibered correspondence}%
   {inverse reduced fibered correspondence, 1}%
   {F}{0}%
\Symb%
   {map to Cartesian product}%
   {map to Cartesian product}%
   {F}{0}%
\Symb%
   {norm of functional}%
   {norm of functional}%
   {F}{0}%
\Symb%
   {norm of map}%
   {norm of map}%
   {F}{0}%
\Symb%
   {norm of polylinear map}%
   {norm of polymap}%
   {F}{0}%
\Symb%
   {norm of representation}%
   {norm of representation}%
   {F}{0}%
\Symb%
   {orbit of representation}%
   {orbit of representation}%
   {F}{0}%
\Symb%
   {orthonormal basis}%
   {Orthonormal Basis, division ring}%
   {F}{0}%
\Symb%
   {quaternion algebra  over field ${\rm {\mathbb{F}}}$}%
   {quaternion algebra F a b}%
   {F}{0}%
\Symb%
   {reference frame}%
   {reference frame}%
   {F}{0}%
\Symb%
   {reference frame, extensive definition}%
   {reference frame, extensive definition}%
   {F}{0}%
\Symb%
   {standard component of biadditive map $f$ over field $F$}%
   {standard component of biadditive map, division ring}%
   {F}{0}%
\Symb%
   {standard component of linear map}%
   {standard component of linear map, G}%
   {F}{0}%
\Symb%
   {standard component of polylinear map}%
   {standard component of polylinear map}%
   {F}{0}%
\Symb%
   {standard component of quadratic map $f$ over field $F$}%
   {standard component of quadratic map, division ring}%
   {F}{0}%
\Symb%
   {standard component of tensor}%
   {standard component of tensor}%
   {F}{0}%
\Symb%
   {sum of maps}%
   {sum of maps}%
   {F}{0}%
\Symb%
   {symmetrization of polylinear map}%
   {symmetrization of polylinear map}%
   {F}{0}%

\SetIndexSpace
\Symb%
   {affine transformation group}%
   {affine transformation group}%
   {G}{0}%
\Symb%
   {affine transformation group}%
   {affine transformation group}%
   {G}{0}%
\Symb%
   {Cartesian product of groups $G_1$, ..., $G_n$}%
   {Cartesian product of groups, 1 n}%
   {G}{0}%
\Symb%
   {\CR matrix group}%
   {cr-matrix group}%
   {G}{0}%
\Symb%
   {fibered little group of section $h$}%
   {fibered little group}%
   {G}{0}%
\Symb%
   {fibered stability group of section $h$}%
   {fibered stability group}%
   {G}{0}%
\Symb%
   {group of automorphisms of representation $f$}%
   {group of automorphisms of representation}%
   {G}{0}%
\Symb%
   {group of homomorphisms of vector space $\Vector V$}%
   {GV}%
   {G}{0}%
\Symb%
   {indefinite integral}%
   {indefinite integral}%
   {G}{0}%
\Symb%
   {left defined Lie algebra of Lie group}%
   {left defined Lie algebra of Lie group}%
   {G}{0}%
\Symb%
   {Lie algebra of Lie group}%
   {Lie algebra of Lie group}%
   {G}{0}%
\Symb%
   {linear transformation group}%
   {linear transformation group}%
   {G}{0}%
\Symb%
   {little group}%
   {little group}%
   {G}{0}%
\Symb%
   {orbit of effective Ts representation of group}%
   {orbit of effective starT representation of fibered group}%
   {G}{0}%
\Symb%
   {orbit of effective \Ts representation of fibered group}%
   {orbit of effective Tstar representation of fibered group}%
   {G}{0}%
\Symb%
   {\RC matrix group}%
   {rc-matrix group}%
   {G}{0}%
\Symb%
   {right defined Lie algebra of Lie group}%
   {right defined Lie algebra}%
   {G}{0}%
\Symb%
   {stability group}%
   {stability group}%
   {G}{0}%

\SetIndexSpace
\Symb%
   {Hadamard inverse of matrix}%
   {Hadamard inverse of matrix}%
   {H}{0}%
\Symb%
   {horizontal component of vector}%
   {horizontal component of vector}%
   {H}{0}%
\Symb%
   {horizontal subspace}%
   {horizontal subspace}%
   {H}{0}%
\Symb%
   {quaternion algebra}%
   {quaternion algebra}%
   {H}{0}%
\Symb%
   {quaternion algebra}%
   {quaternion algebra H a b}%
   {H}{0}%
\Symb%
   {set of homomorphisms}%
   {set of homomorphisms}%
   {H}{0}%

\SetIndexSpace
\Symb%
   {infinitesimal generator of representation}%
   {infinitesimal generator i of representation}%
   {I}{0}%
\Symb%
   {infinitesimal generator of representation}%
   {infinitesimal generator of representation}%
   {I}{0}%
\Symb%
   {Lie group infinitesimal generators}%
   {Lie group infinitesimal generators}%
   {I}{0}%
\Symb%
   {vector module of algebra $A$}%
   {vector module of algebra}%
   {I}{0}%
\Symb%
   {vector module of ring $D$}%
   {vector module of ring}%
   {I}{0}%
\Symb%
   {vector of element $d$ of algebra}%
   {vector of algebra}%
   {I}{0}%
\Symb%
   {vector of element $d$ of ring}%
   {vector of ring}%
   {I}{0}%

\SetIndexSpace
\Symb%
   {closure operator of representation $f$}%
   {closure operator, representation}%
   {J}{0}%
\Symb%
   {closure operator of tower of representations $\Vector f$}%
   {closure operator, tower of representations}%
   {J}{0}%
\Symb%
   {Jacobian matrix of right shift}%
   {Ea, quaternion, Jacobian matrix}%
   {J}{0}%
\Symb%
   {subrepresentation generated by set $X$}%
   {subrepresentation generated by set}%
   {J}{0}%
\Symb%
   {tower of subrepresentations of tower of representations $\Vector f$ generated by tuple of sets $\VX X$}%
   {subrepresentation generated by tuple of sets}%
   {J}{0}%

\SetIndexSpace
\Symb%
   {kernel of homomorphism}%
   {kernel of homomorphism}%
   {K}{0}%
\Symb%
   {kernel of linear map}%
   {kernel of linear map}%
   {K}{0}%

\SetIndexSpace
\Symb%
   {Cartesian power of systems of subsets}%
   {Cartesian power of systems of subsets}%
   {L}{0}%
\Symb%
   {Cartesian product of systems of subsets}%
   {Cartesian product of systems of subsets}%
   {L}{0}%
\Symb%
   {left $ij$th cofactor of entry of matrix}%
   {left cofactor, matrix}%
   {L}{0}%
\Symb%
   {left double $ij$th cofactor of entry of matrix}%
   {left double cofactor}%
   {L}{0}%
\Symb%
   {left shift}%
   {left shift}%
   {L}{0}%
\Symb%
   {Lie derivative of connection}%
   {Lie derivative of connection}%
   {L}{0}%
\Symb%
   {Lie derivative of metric}%
   {Lie derivative of metric}%
   {L}{0}%
\Symb%
   {limit of correspondence $\Phi$ with respect to the filter $\mathfrak{F}$}%
   {limit of correspondence with respect to the filter}%
   {L}{0}%
\Symb%
   {limit of sequence}%
   {limit of sequence}%
   {L}{0}%
\Symb%
   {linear combination}%
   {linear combination}%
   {L}{0}%
\Symb%
   {module of skew symmetric polylinear maps}%
   {module of skew symmetric polylinear maps}%
   {L}{0}%
\Symb%
   {passive transformation}%
   {passive transformation}%
   {L}{0}%
\Symb%
   {$D$\Hyph module of continuous linear mappings of normed $D$\Hyph module $A_1$ into normed $D$\Hyph module $A_2$}%
   {set continuous linear mappings, module}%
   {L}{0}%
\Symb%
   {set of continuous linear maps}%
   {set continuous linear maps, vector}%
   {L}{0}%
\Symb%
   {set of continuous polylinear maps}%
   {set continuous polylinear maps}%
   {L}{0}%
\Symb%
   {set of linear maps}%
   {set linear maps}%
   {L}{0}%
\Symb%
   {set of left-side nonsingular transformations of universal algebra $M$}%
   {set of left-side nonsingular transformations}%
   {L}{0}%
\Symb%
   {set of polylinear maps}%
   {set polylinear maps}%
   {L}{0}%
\Symb%
   {set of $n$\hyph linear maps}%
   {set polylinear maps An}%
   {L}{0}%
\Symb%
   {set of polylinear maps}%
   {set polylinear maps, D vector space}%
   {L}{0}%
\Symb%
   {set of polylinear maps of algebras $A_1$, ..., $A_n$ into algebra $A$}%
   {set polylinear maps, finite dimensional algebra}%
   {L}{0}%

\SetIndexSpace
\Symb%
   {set of left-side transformations of the universal algebra $M$}%
   {set of left-side transformations}%
   {M}{0}%
\Symb%
   {set of maps to $\Omega$\Hyph group $A$}%
   {set of maps to Omega group}%
   {M}{0}%
\Symb%
   {set of right-side transformations of universal algebra $M$}%
   {set of right-side transformations}%
   {M}{0}%
\Symb%
   {space of orbits of \Ts{G}representation}%
   {space of orbits of G* representation}%
   {M}{0}%

\SetIndexSpace
\Symb%
   {norm of quaternion $x$}%
   {norm, quaternion algebra}%
   {N}{0}%
\Symb%
   {nucleus of $D$\Hyph algebra $A$}%
   {nucleus of algebra}%
   {N}{0}%

\SetIndexSpace
\Symb%
   {geometric object in coordinate representation defined in \rcd vector space}%
   {geometric object, coordinate rcd vector space}%
   {O}{0}%
\Symb%
   {geometric object in coordinate representation}%
   {geometric object, coordinate vector space}%
   {O}{0}%
\Symb%
   {geometric object defined in \rcd vector space}%
   {geometric object, rcd vector space}%
   {O}{0}%
\Symb%
   {octonion algebra}%
   {octonion algebra}%
   {O}{0}%
\Symb%
   {orbit of representation of fibered group $\Bundle G$}%
   {orbit of representation of fibered group}%
   {O}{0}%
\Symb%
   {orbit of \Ts{G}representation of the group}%
   {orbit of representation of group}%
   {O}{0}%
\Symb%
   {tensor product}%
   {tensor product}%
   {O}{0}%

\SetIndexSpace
\Symb%
   {bundle}%
   {bundle}%
   {P}{0}%
\Symb%
   {bundle of level $2$}%
   {bundle of level 2}%
   {P}{0}%
\Symb%
   {bundle of level $n$}%
   {bundle of level n}%
   {P}{0}%
\Symb%
   {Cartesian power $n$ of bundle $\bundle{}{p}{E}{}$}%
   {Cartesian power of bundle}%
   {P}{0}%
\Symb%
   {Cartesian product of bundles}%
   {Cartesian product of bundles, definition 1}%
   {P}{0}%
\Symb%
   {passive representation of group $G(f)$ in basis manifold $\mathcal B(f)$}%
   {passive representation in basis manifold}%
   {P}{0}%
\Symb%
   {passive representation of group $G(\Vector f)$ in basis manifold $\mathcal B(\Vector f)$}%
   {passive representation in basis manifold, tower of representations}%
   {P}{0}%
\Symb%
   {reduced Cartesian product of bundles}%
   {reduced Cartesian product of bundles, definition 1}%
   {P}{0}%
\Symb%
   {set of nonsingular \sT transformations of bundle $\bundle{}pE{}$}%
   {set of starT nonsingular transformations of bundle, projection}%
   {P}{0}%
\Symb%
   {set of nonsingular \Ts transformations of bundle $\bundle{}pE{}$}%
   {set of Tstar nonsingular transformations of bundle, projection}%
   {P}{0}%

\SetIndexSpace
\Symb%
   {active transformation}%
   {active transformation}%
   {R}{0}%
\Symb%
   {Cartan curvature}%
   {Cartan curvature}%
   {R}{0}%
\Symb%
   {\CR rank of matrix}%
   {cr-rank of matrix}%
   {R}{0}%
\Symb%
   {diagonal in bundle  $\bundle{}pA{}$}%
   {diagonal in bundle, 2}%
   {R}{0}%
\Symb%
   {diagonal in bundle $\Bundle A$}%
   {diagonal in reduced bundle, 2}%
   {R}{0}%
\Symb%
   {image of $m$ under endomorphism $R$ of effective representation}%
   {endomorphism image, effective representation}%
   {R}{0}%
\Symb%
   {image of tuple $\VX a$ under endomorphism $\VX r$ of tower of effective representations}%
   {endomorphism image, tower of effective representations}%
   {R}{0}%
\Symb%
   {curvature}%
   {GLn curvature_overline}%
   {R}{0}%
\Symb%
   {product of rings of sets}%
   {product of rings of sets}%
   {R}{0}%
\Symb%
   {$\RCcirc$\Hyph product of matrices of maps}%
   {rc product of matrices of maps}%
   {R}{0}%
\Symb%
   {\RC rank of matrix}%
   {rc-rank of matrix}%
   {R}{0}%
\Symb%
   {right $ij$th cofactor of entry of matrix}%
   {right cofactor, matrix}%
   {R}{0}%
\Symb%
   {right double $ij$th cofactor of entry of matrix}%
   {right double cofactor}%
   {R}{0}%
\Symb%
   {right shift}%
   {right shift}%
   {R}{0}%
\Symb%
   {$i$th row determinant of matrix $\bfA$}%
   {row determinant}%
   {R}{0}%
\Symb%
   {scalar algebra of algebra $A$}%
   {scalar algebra of algebra}%
   {R}{0}%
\Symb%
   {scalar algebra of ring $D$}%
   {scalar algebra of ring}%
   {R}{0}%
\Symb%
   {scalar of element $d$ of algebra}%
   {scalar of algebra}%
   {R}{0}%
\Symb%
   {scalar of element $d$ of ring}%
   {scalar of ring}%
   {R}{0}%
\Symb%
   {set of right-side nonsingular transformations of universal algebra $M$}%
   {set of right-side nonsingular transformations}%
   {R}{0}%
\Symb%
   {spherical coordinates}%
   {spherical coordinates}%
   {R}{0}%

\SetIndexSpace
\Symb%
   {composition of fibered correspondences}%
   {composition of fibered correspondences}%
   {S}{0}%
\Symb%
   {hyperbolic sine}%
   {hyperbolic sine}%
   {S}{0}%
\Symb%
   {inverse fibered correspondence}%
   {inverse fibered correspondence, 2}%
   {S}{0}%
\Symb%
   {inverse reduced fibered correspondence}%
   {inverse reduced fibered correspondence, 2}%
   {S}{0}%
\Symb%
   {Lebesgue integral}%
   {Lebesgue integral}%
   {S}{0}%
\Symb%
   {linear span in vector space}%
   {linear span, vector space}%
   {S}{0}%
\Symb%
   {image of basis $X$ under passive transformation $S$}%
   {passive transformation of basis, representation}%
   {S}{0}%
\Symb%
   {image of basis $\VX  X$ under passive transformation $\VX s$}%
   {passive transformation of basis, tower of representations}%
   {S}{0}%
\Symb%
   {set of permutations}%
   {set of permutations}%
   {S}{0}%
\Symb%
   {set of transpositions}%
   {set of transpositions}%
   {S}{0}%
\Symb%
   {sine}%
   {sine}%
   {S}{0}%
\Symb%
   {symmetric group}%
   {symmetric group}%
   {S}{0}%

\SetIndexSpace
\Symb%
   {category of left-side representations}%
   {category of left-side representations}%
   {T}{0}%
\Symb%
   {tangent plane to Lie group $G$}%
   {tangent plane to Lie group}%
   {T}{0}%
\Symb%
   {trace of quaternion $x$}%
   {trace, quaternion algebra}%
   {T}{0}%

\SetIndexSpace
\Symb%
   {affine space}%
   {affine space}%
   {V}{0}%
\Symb%
   {conjugated affine space}%
   {conjugated affine space}%
   {V}{0}%
\Symb%
   {conjugated vector space}%
   {conjugated vector space}%
   {V}{0}%
\Symb%
   {coordinate vector space}%
   {coordinate vector space}%
   {V}{0}%
\Symb%
   {coordinates in vector space}%
   {coordinates in vector space}%
   {V}{0}%
\Symb%
   {direct product of $\RCstar D_i$\hyph vector spaces $\Vector V_1$, ..., $\Vector V_n$}%
   {direct product, rcd vector space, 1 n}%
   {V}{0}%
\Symb%
   {dual space of \rcd vector space $\Vector V$}%
   {dual space of rcd vector space}%
   {V}{0}%
\Symb%
   {hermitian conjugated vector}%
   {hermitian conjugated vector}%
   {V}{0}%
\Symb%
   {linear composition of vectors}%
   {linear composition of vectors}%
   {V}{0}%
\Symb%
   {vector space}%
   {V}%
   {V}{0}%
\Symb%
   {vertical component of vector}%
   {vertical component of vector}%
   {V}{0}%
\Symb%
   {vertical subspace}%
   {vertical subspace}%
   {V}{0}%

\SetIndexSpace
\Symb%
   {set of coordinates of representation $J(f,X)$}%
   {coordinate set of representation}%
   {W}{0}%
\Symb%
   {set of tuples of coordinates of tower of representations $\Vector J(\Vector f,\VX X)$}%
   {coordinate set of tower of representations}%
   {W}{0}%
\Symb%
   {coordinates of basis $X'$ relative to basis $X$ of representation}%
   {coordinates of basis relative to basis, representation}%
   {W}{0}%
\Symb%
   {coordinates of element $m$ of representation $f$ relative to set $X$}%
   {coordinates of element relative to generating set, representation}%
   {W}{0}%
\Symb%
   {coordinates of element $m$ relative to set $X$}%
   {coordinates of element relative to set, representation}%
   {W}{0}%
\Symb%
   {tuple of coordinates of element $\Vector a*$ relative to tuple of sets $\VX X$}%
   {coordinates of element, tower of representations}%
   {W}{0}%
\Symb%
   {geometric object in coordinate representation defined in $\Omega_2$\Hyph algebra $M$}%
   {geometric object, coordinate representation g}%
   {W}{0}%
\Symb%
   {geometric object in coordinate representation defined in tuple of $\VX\Omega$\Hyph algebras $\VX A$}%
   {geometric object, coordinate tower of representations g}%
   {W}{0}%
\Symb%
   {geometric object defined in $\Omega_2$\Hyph algebra $M$}%
   {geometric object, representation g}%
   {W}{0}%
\Symb%
   {geometric object defined in tuple of $\VX\Omega$\Hyph algebras $\VX A$}%
   {geometric object, tower of representations g}%
   {W}{0}%
\Symb%
   {geometric object in vector space}%
   {geometric object, vector space}%
   {W}{0}%
\Symb%
   {set of coordinates of set $B\subset J(f,X)$}%
   {subset of coordinates of representation}%
   {W}{0}%
\Symb%
   {coordinates of tuple of sets $\VX B$ relative to tuple of sets $\VX X$}%
   {subset of coordinates of tower of representations}%
   {W}{0}%
\Symb%
   {coordinates of set $B_k$ relative to tuple of sets $\VX X$}%
   {subset of coordinates of tower of representations, k}%
   {W}{0}%
\Symb%
   {set of $\Omega_2$\Hyph words representing set $B\subset J(f,X)$}%
   {subset of words of representation}%
   {W}{0}%
\Symb%
   {superposition of coordinates of the representation $f$ and the element $m$}%
   {superposition of coordinates, representation}%
   {W}{0}%
\Symb%
   {superposition of coordinates of the tower of representations $\Vector f$ and the element $\VX a$}%
   {superposition of coordinates, tower of representations}%
   {W}{0}%
\Symb%
   {$\Omega_2$\Hyph word representing element $m\in J(f,X)$}%
   {word of element relative to generating set, representation}%
   {W}{0}%
\Symb%
   {set of $\Omega_2$\Hyph words of representation $J(f,X)$}%
   {word set of representation}%
   {W}{0}%
\Symb%
   {set of tuples of $\VX{\Omega}$\Hyph words of tower of representations $\Vector J(\Vector f,\VX X)$}%
   {word set of tower of representations}%
   {W}{0}%
\Symb%
   {tuple of words of element $\Vector a*$ relative to tuple of sets $\VX X$}%
   {words of element, tower of representations}%
   {W}{0}%

\SetIndexSpace
\Symb%
   {conjugate of quaternion $x$}%
   {conjugate of quaternion}%
   {X}{0}%
\Symb%
   {local basis of affine space}%
   {local basis of affine space}%
   {X}{0}%
\Symb%
   {anholonomic coordinate}%
   {x(k)}%
   {X}{0}%

\SetIndexSpace
\Symb%
   {center of $D$\Hyph algebra $A$}%
   {center of algebra}%
   {Z}{0}%
\Symb%
   {center of ring $D$}%
   {center of ring}%
   {Z}{0}%

\SetIndexSpace
\Symb%
   {deviation of trajectories}%
   {deviation of trajectories}%
   {Delta}{1}%
\Symb%
   {identical transformation}%
   {identical transformation}%
   {Delta}{1}%
\Symb%
   {image of vector $\Vector e_k\in\Basis e$ under isomorphism to coordinate vector space}%
   {image of vector e_k, coordinate vector space}%
   {Delta}{1}%
\Symb%
   {Kronecker symbol}%
   {Kronecker symbol}%
   {Delta}{1}%

\SetIndexSpace
\Symb%
   {anholonomic coordinates of connection}%
   {anholonomic coordinates of connection}%
   {Gamma}{1}%
\Symb%
   {Cartan symbol}%
   {Cartan symbol}%
   {Gamma}{1}%
\Symb%
   {connection}%
   {conection overline}%
   {Gamma}{1}%
\Symb%
   {connection coefficients in $D$\Hyph affine space}%
   {connection coefficients, D affine space}%
   {Gamma}{1}%
\Symb%
   {connection in $D$\Hyph affine manifold}%
   {connection, affine manifold}%
   {Gamma}{1}%
\Symb%
   {$D$\Hyph affine connection coefficients on manifold}%
   {D affine connection coefficients, manifold}%
   {Gamma}{1}%
\Symb%
   {holonomic coordinates of connection}%
   {holonomic coordinates of connection}%
   {Gamma}{1}%
\Symb%
   {Cartan connection}%
   {overbrace Gamma i kl}%
   {Gamma}{1}%
\Symb%
   {set of sections of bundle}%
   {set of sections of bundle}%
   {Gamma}{1}%

\SetIndexSpace
\Symb%
   {inverse operator to operator $\psi_l$}%
   {inverse operator to operator psi l}%
   {Lambda}{1}%
\Symb%
   {inverse operator to operator $\psi_r$}%
   {inverse operator to operator psi r}%
   {Lambda}{1}%

\SetIndexSpace
\Symb%
   {Cartesian product of measures}%
   {Cartesian product of measures}%
   {Mu}{1}%
\Symb%
   {power of measure}%
   {power of measure}%
   {Mu}{1}%
\Symb%
   {product of measures}%
   {product of measures}%
   {Mu}{1}%
\Symb%
   {product of measures}%
   {product of measures, otimes}%
   {Mu}{1}%

\SetIndexSpace
\Symb%
   {anholonomity object}%
   {anholonomity object}%
   {Omega}{1}%
\Symb%
   {definite integral}%
   {definite integral}%
   {Omega}{1}%
\Symb%
   {integral of differential $1$\Hyph form along path}%
   {integral of differential 1 form along path}%
   {Omega}{1}%
\Symb%
   {norm of operation}%
   {norm of operation}%
   {Omega}{1}%
\Symb%
   {operator domain}%
   {operator domain}%
   {Omega}{1}%
\Symb%
   {set of differential $p$\Hyph forms}%
   {set of differential p forms}%
   {Omega}{1}%
\Symb%
   {set of $n$\Hyph ary operations of $\Omega$\Hyph algebra}%
   {set of n-ary operations}%
   {Omega}{1}%
\Symb%
   {set of $n$\Hyph ary operators}%
   {set of n-ary operators}%
   {Omega}{1}%

\SetIndexSpace
\Symb%
   {left basic operator of Lie group over algebra $A$}%
   {L basic operator of Lie group over algebra A}%
   {Psi}{1}%
\Symb%
   {left basic operator of group Lie}%
   {Lie Basic Operator L}%
   {Psi}{1}%
\Symb%
   {left basic operator of Lie 1-parameter group}%
   {Lie Basic Operator L, 1-Parameter Group}%
   {Psi}{1}%
\Symb%
   {left basic operator of Lie 1-parameter group over algebra $A$}%
   {Lie Basic Operator L, 1-Parameter Group, algebra}%
   {Psi}{1}%
\Symb%
   {right basic operator of group Lie}%
   {Lie Basic Operator R}%
   {Psi}{1}%
\Symb%
   {right basic operator of Lie 1-parameter group}%
   {Lie Basic Operator R, 1-Parameter Group}%
   {Psi}{1}%
\Symb%
   {right basic operator of Lie 1-parameter group over algebra $A$}%
   {Lie Basic Operator R, 1-Parameter Group, algebra}%
   {Psi}{1}%
\Symb%
   {right basic operator of Lie group over algebra $A$}%
   {R basic operator of Lie group over algebra A}%
   {Psi}{1}%

\SetIndexSpace
\Symb%
   {fibered subset}%
   {fibered subset}%
   {Sigma}{1}%
\Symb%
   {parity of permutation}%
   {parity of permutation}%
   {Sigma}{1}%
\Symb%
   {subbundle}%
   {subbundle}%
   {Sigma}{1}%

\SetIndexSpace
\Symb%
   {Cartan derivative}%
   {overbrace nabla_l}%
   {Nabla}{2}%
\Symb%
   {derivative}%
   {overline nabla_l, definition 1}%
   {Nabla}{2}%

\SetIndexSpace
\Symb%
   {Lie group composition law}%
   {Lie group composition law}%
   {Phi}{1}%
\Symb%
   {restriction of correspondence $\Phi$ to set $C$}%
   {restriction of correspondence}%
   {Phi}{1}%

\SetIndexSpace
\Symb%
   {Cartesian product of bundles}%
   {Cartesian product of bundles, definition 2}%
   {Pi}{1}%
\Symb%
   {Cartesian product of groups $G_i$, $i\in I$}%
   {Cartesian product of groups}%
   {Pi}{1}%
\Symb%
   {Cartesian product of groups $G_1$, ..., $G_n$}%
   {Cartesian product of groups, i 1 n}%
   {Pi}{1}%
\Symb%
   {Cartesian product of total spaces}%
   {Cartesian product of total spaces, definition 2}%
   {Pi}{1}%
\Symb%
   {coproduct in category}%
   {coproduct in category}%
   {Pi}{1}%
\Symb%
   {direct product of division rings $D_i$, $i\in I$}%
   {direct product of division rings}%
   {Pi}{1}%
\Symb%
   {direct product of division rings $D_1$, ..., $D_n$}%
   {direct product of division rings, i 1 n}%
   {Pi}{1}%
\Symb%
   {direct product of $\RCstar D_i$\hyph vector spaces $\Vector V_i$, $i\in I$}%
   {direct product, rcd vector space}%
   {Pi}{1}%
\Symb%
   {direct product of $\RCstar D_i$\hyph vector spaces}%
   {direct product, rcd vector space, i 1 n}%
   {Pi}{1}%
\Symb%
   {product in category}%
   {product in category}%
   {Pi}{1}%
\Symb%
   {reduced Cartesian product of bundles}%
   {reduced Cartesian product of bundles, definition 2}%
   {Pi}{1}%
\Symb%
   {reduced Cartesian product of total spaces}%
   {reduced Cartesian product of total spaces, definition 2}%
   {Pi}{1}%

\CloseIndex